\documentclass[fleqn]{article}

\usepackage{amsmath}
\usepackage{amssymb}
\usepackage{amsthm}
\usepackage{graphicx}

\usepackage{cite}

\newtheorem{theorem}{Theorem}
\newtheorem{lemma}[theorem]{Lemma}

\title{
  Where the Links--Gould invariant first fails to distinguish nonmutant prime knots
}

\author{
  David De Wit
  and
  Jon Links%
  \footnote{%
    Department of Mathematics, The University of Queensland,
    4072, Brisbane, Australia.
    \newline
    \texttt{Dr\_David\_De\_Wit@yahoo.com.au},
    \texttt{jrl@maths.uq.edu.au}
  }
}

\begin{document}

\maketitle

\begin{abstract}
  \noindent%
  It is known that the first two-variable Links--Gould quantum link invariant
  $LG\equiv LG^{2,1}$ is more powerful than the HOMFLYPT and Kauffman
  polynomials, in that it distinguishes all prime knots (including
  reflections) of up to $10$ crossings. Here we report investigations which
  greatly expand the set of evaluations of $LG$ for prime knots. Through them,
  we show that the invariant is complete, modulo mutation, for all prime knots
  (including reflections) of up to $11$ crossings, but fails to distinguish
  some nonmutant pairs of $12$-crossing prime knots.  As a byproduct, we
  classify the mutants within the prime knots of $11$ and $12$ crossings. In
  parallel, we learn that $LG$ distinguishes the chirality of all chiral prime
  knots of at most $12$ crossings. We then demonstrate that every
  mutation-insensitive link invariant fails to distinguish the chirality of a
  number of $14$-crossing prime knots. This provides $14$-crossing examples of
  chiral prime knots whose chirality is undistinguished by $LG$.
\end{abstract}


\section{The Links--Gould invariant}

For any positive integers $m$ and $n$, the Links--Gould quantum link invariant
$LG^{m,n}$ is a two-variable invariant of oriented links. We here describe it
using the variables $q$ and $p\equiv q^{\alpha+(m-n)/2}$ which are inherited
from its definition via the $2^{mn}$-dimensional $\alpha$-parametric
representation of highest weight $(\dot{0}_m\,|\,\dot{\alpha}_n)$ of the
quantum superalgebra $U_q[gl(m|n)]$. Its construction was originally described
in~\cite{GouldLinksZhang:1996b,LinksGould:1992b}, and some of its properties
together with some concrete evaluations are available
in~\cite{DeWit:2001b,DeWit:2002a,DeWitIshiiLinks:2004,DeWitKauffmanLinks:1999a}.

The case $LG^{1,1}$, in which the two variables degenerate to a single
variable ($p$), is actually the Alexander--Conway polynomial, and the next
simplest case, $LG^{2,1}$, is the first truly two-variable invariant of the
family. Apart from $LG^{1,1}$, more is known about $LG^{2,1}$ than about any
of its
kin~\cite{DeWit:2000,Ishii:2004b,Ishii:2004e,Ishii:2004a,Ishii:2003,Ishii:2004d,IshiiKanenobu:2004b,IshiiKanenobu:2004a},
and we generally refer to $LG^{2,1}$ as \emph{the}\/ Links--Gould invariant
$LG$. In particular, it has been proven to be polynomial~\cite{Ishii:2004d};
this is still only surmised for general $LG^{m,n}$.

As with all such quantum link invariants, $LG^{m,n}$ is insensitive to the
inversion of links, and this manifests in its invariance under the inversion
of $p$ (modulo a sign --- this is a well-known property of the
Alexander--Conway polynomial). That is, for any oriented link $L$, the
(polynomial) evaluation $LG^{m,n}_L$ is necessarily palindromic in $p$, meaning
that $LG^{m,n}_{-L}(q,p)\dot{=}LG^{m,n}_L(q,p^{-1})$, where by $-L$ we intend
the inversion of $L$, and by $\dot{=}$ we intend equality up to a sign $\pm 1$.
Also, as for most of the more familiar link invariants, $LG^{m,n}$ fails to
distinguish mutant links.

As we know these limitations, it is of interest to determine where $LG$ first
fails to distinguish nonmutant prime knots. Significantly, $LG$ is known to be
complete for all prime knots (including reflections) of up to $10$
crossings~\cite{DeWit:2000}. In this, it is known to be more powerful than the
HOMFLYPT and Kauffman polynomials. Unlike the situation for those polynomials,
however, a set of skein relations sufficient to determine $LG$ for any
arbitrary link is not yet known, and most evaluations to date have been
performed using a computationally expensive state model method.

Here, we report on evaluations of $LG$ for a range of the prime knots from the
Hoste--Thistlethwaite--Weeks (HTW) tables of prime knots of up to $16$
crossings~\cite{HosteThistlethwaiteWeeks:1998}. The present list of evaluations expands the
previous published list~\cite{DeWit:2000} by a factor of $200$. Of significance, we determine that
$LG$ is complete, modulo mutation, for all prime knots (including reflections) of up to $11$
crossings, but that it fails to distinguish several nonmutant $12$-crossing prime knots.

A byproduct of these investigations is a listing of all mutant cliques within
the prime knots of $11$ and $12$ crossings, which should serve as a useful
test-bed for future experiments with invariants that are not
mutation-insensitive. We also mention a number of open basic questions about
mutation.

It is also known that $LG$ is not always able to distinguish between the
reflections of chiral links, as Ishii and Kanenobu have constructed an
infinite family of chiral knots whose chirality is undetected by
$LG$~\cite{IshiiKanenobu:2004b}. At least some of these links are prime knots,
although of unknown crossing numbers. Herein, we demonstrate that $LG$
distinguishes the chirality of all chiral prime knots of up to $12$ crossings.
We also identify several $14$-crossing chiral prime knots whose chirality is
not distinguished by $LG$, or indeed by any mutation-insensitive invariant.
That said, we have as yet no example of a chiral link whose chirality is not
detected by $LG$ yet is detected by either the HOMFLYPT or the Kauffman
polynomial.


\section{Evaluations of $LG$ for the HTW prime knots}

Our evaluations of $LG$ are with respect to the HTW tables of prime knots of
up to $16$ crossings~\cite{HosteThistlethwaiteWeeks:1998}.
Dowker--Thistlethwaite (DT) codes and many other data associated with the
knots in these tables may be accessed via the program \textsc{Knotscape}
(version 1.01), as mentioned in~\cite{HosteThistlethwaiteWeeks:1998}. We let
$c^A_i$ (respectively $c^N_i$) denote the $i$th alternating (respectively
nonalternating) prime knot (type) of $c$ crossings in the HTW tables. The
presence of the superscript decorating the crossing number identifies the knot
as from the HTW tables, and facilitates the concurrent use of the classical
undecorated Alexander--Briggs notation; for instance $5^A_1$ and $5_2$ denote
the same knot type. In~\cite{DeWit:2000} is found a listing of evaluations of
$LG$ for all prime knots of up to $10$ crossings, listed with respect to the
Alexander--Briggs ordering. Herein, we abandon that ordering, and instead
adopt the HTW ordering.

In the classical tables, as reproduced for instance in~\cite{Kawauchi:1996},
canonical representatives of knot types are described graphically, however
these diagrams omit an orientation, so only describe \emph{unoriented} knots.
In contrast, the HTW tables describe knot types via DT codes, and as any DT
code necessarily encodes both a link and its reflection, the tables do not
prescribe chiralities (orientations of space), although they do prescribe
string orientations. (As quantum link invariants are insensitive to inversion,
this is unimportant to our considerations.) An evaluation of a link invariant
for a link described in terms of a DT code thus should be regarded as unique
up to a symmetry transformation corresponding to that of reflection of the
link. For $LG^{m,n}$ this transformation corresponds to the inversion of the
variable $q$. Using the HTW tables to define knots thus means that $LG$
polynomials should be regarded as attached to knot \emph{types}, rather than
knots per se. Herein, the statement that $LG$ cannot distinguish distinct
links $L_1$ and $L_2$ means that the polynomials $LG_{L_1}$ and $LG_{L_2}$ are
either equal or related by inversion of $q$. In this situation, we describe the
links as \emph{$LG$-equivalent}. Thus the $LG$-polynomial of an achiral link
is necessarily palindromic in $q$.

Our evaluation of some $LG^{m,n}$ for a particular link $L$ is generally
obtained by the application of a state model method to a braid $\beta$, where
the closure $\widehat{\beta}$ corresponds to $L$. We have implemented this in
a \textsc{Mathematica} package called \textsc{Links--Gould Explorer}; a
description of the algorithms used appears in~\cite{DeWit:BraidedStreams}.
Evaluations of $LG^{2,1}$ are currently generally feasible for braids of
string index at most $5$ (\emph{sometimes} $6$ and even $7$) and a
`reasonable' number of crossings (say less than $30$). By `feasible', we
intend computations which demand at most about $1$Gb of memory and run-times
of at most a few CPU-hours on current commonly-available hardware. Memory
requirements increase explosively with braid width.

Braids for the HTW knots are obtainable via the implementation of Vogel's
Algorithm~\cite{Vogel:1990} in the \textsc{Mathematica} package \textsc{K2K}
(version 1.3.3) by Imafuji and Ochiai~\cite{ImafujiOchiai:2002}. (\textsc{K2K}
is an interface to an underlying \textsc{C} program
\texttt{KnotTheoryByComputer} by Ochiai.) Unfortunately, a bug in the
implementation of Vogel's Algorithm in \textsc{K2K} means that it fails to
terminate for exactly $11$ of the $1,701,936$ HTW knots. (Otherwise, it yields
\emph{correct} braids!) In~\cite{DeWit:BraidsforPretzels} we explain the
application of other methods to determine braids for these $11$ knots; the net
result is that we do have braids for all the HTW knots.

Recall that a \emph{minimal} braid is one whose width is minimal over all candidate braids for a
given knot, and this width is called the \emph{string index} of the knot. Sadly, the
automatically-generated \textsc{K2K}-braids are generally far from minimal, and where possible, we
have used thinner ones. Specifically, we have to hand minimal braids for all the HTW knots of up to
$12$ crossings. This data is obtained from Livingston's website \textsc{Table of Knot
Invariants}~\cite{Livingston:TableofKnotInvariants} together with that of Stoimenow, called
\textsc{Knot data tables}~\cite{Stoimenow:Knotdatatables}. The braids for knots of up to $10$
crossings were constructed by Gittings~\cite{Gittings:2004} and those of $11$ and $12$ crossings
by Stoimenow. These braids are also (generally) of minimal length for their minimal width, and of
a structure which facilitates their use in state model algorithms. For the knots of $11$
(respectively $12$) crossings, they are of maximum string (that is, braid) index $6$ (respectively
$7$).

We refer to our full collection of braids from these diverse sources as the
\textsc{K2K}$'$-braids. To describe them, we introduce some convenient
notation. Where $P$ is either $A$ or $N$, denote by $\mathbb{K}^P_c$ the set
of all prime knots $c^P_i$ of $c$ crossings, and set
$\mathbb{K}_c\triangleq\mathbb{K}^A_c\cup\mathbb{K}^N_c$. Similarly, denote by
$\mathbb{K}^P_{c,s}$ the set of all prime knots $c^P_i$ of $c$ crossings whose
\textsc{K2K}$'$-braids have $s$ strings, and set
$\mathbb{K}_{c,s}\triangleq\mathbb{K}^A_{c,s}\cup\mathbb{K}^N_{c,s}$.

\renewcommand{\arraystretch}{1.5}
\renewcommand{\tabcolsep}{4pt}

\begin{table}[htbp]
  \tiny
  \begin{centering}
  \begin{tabular}{r@{\hspace{0pt}}r|*{15}{r}|r}
    \multicolumn{2}{c}{} & \multicolumn{15}{c}{\small $c$} & \\
    &   & 0& 3& 4& 5& 6& 7& 8&  9& 10&  11&  12&   13&   14&    15&      16&totals \\
    \cline{2-18}
    &  1& 1& .& .& .& .& .& .&  .&  .&   .&   .&    .&    .&     .&      . &1\\
    &  2& .& 1& .& 1& .& 1& .&  1&  .&   1&   .&    1&    .&     1&      . &7\\
    &  3& .& .& 1& 1& 2& 2& 8&  4& 23&   6&  71&    .&  225&     .&    746 &1089\\
    &  4& .& .& .& .& 1& 4& 7& 23& 48& 149& 290& 1017&    .&  8018&      . &9557\\
    \cline{17-17}
    &  5& .& .& .& .& .& .& 3& 13& 47& 164& 652&    .& 6905&\multicolumn{1}{r|}{.}&  80342 &88126\\
    \cline{14-16}
    &  6& .& .& .& .& .& .& .&  .&  5&  47&\multicolumn{1}{r|}{258}& 2570&    .& 40388&      . &43268\\
    {\small $s$}
    &    7& .& .& .& .& .& .& .&  .&  .&   .&\multicolumn{1}{r|}{17}&    .& 9479&     .& 190642 &200138\\
    \cline{3-13}
    &  8& .& .& .& .& .& .& .&  .&  .&   .&   .& 1145&    .& 30684&      . &31829\\
    &  9& .& .& .& .& .& .& .&  .&  .&   .&   .&    .& 2728&     .&  94924 &97652\\
    & 10& .& .& .& .& .& .& .&  .&  .&   .&   .&  140&    .&  5838&      . &5978\\
    & 11& .& .& .& .& .& .& .&  .&  .&   .&   .&    .&  196&     .&  12697 &12893\\
    & 12& .& .& .& .& .& .& .&  .&  .&   .&   .&    5&    .&   334&      . &339\\
    & 13& .& .& .& .& .& .& .&  .&  .&   .&   .&    .&    3&     .&    448 &451\\
    \cline{2-18}
    & totals&1&1&1& 2& 3& 7& 18& 41& 123& 367& 1288& 4878&19536&85263&379799&491328
  \end{tabular}

  \vspace{1\baselineskip}
  {\small Alternating}
  \vspace{2\baselineskip}

  \begin{tabular}{r@{\hspace{0pt}}r|*{15}{r}|r}
    \multicolumn{2}{c}{} & \multicolumn{15}{c}{\small $c$} & \\
    &   & 0& 3& 4& 5& 6& 7& 8&  9& 10&  11&  12&   13&   14&    15&      16&totals \\
    \cline{2-18}
    & 1&.& .& .& .& .& .& .& .&  .&  .&  .&    .&     .&     .&      .&0\\
    & 2&.& .& .& .& .& .& .& .&  .&  .&  .&    .&     .&     .&      .&0\\
    & 3&.& .& .& .& .& .& 3& .& 14&  .& 50&    .&   121&     .&    478&666\\
    & 4&.& .& .& .& .& .& .& 8& 26&130&333& 1247&     .& 15559&      .&17303\\
    \cline{17-17}
    & 5&.& .& .& .& .& .& .& .&  2& 55&486&    .& 11103&\multicolumn{1}{r|}{.}& 244352&255998\\
    \cline{14-16}
    & 6&.& .& .& .& .& .& .& .&  .&.&\multicolumn{1}{r|}{19}& 3067&     .& 92169&      .&95255\\
    {\small $s$}
    & 7&.& .& .& .& .& .& .& .&  .&  .&\multicolumn{1}{r|}{.}&    .& 13589&     .& 546065&559654\\
    \cline{3-13}
    & 8&.& .& .& .& .& .& .& .&  .&  .&  .&  768&     .& 55215&      .&55983\\
    & 9&.& .& .& .& .& .& .& .&  .&  .&  .&    .&  2578&     .& 203041&205639\\
    &10&.& .& .& .& .& .& .& .&  .&  .&  .&   28&     .&  5021&      .&5049\\
    &11&.& .& .& .& .& .& .& .&  .&  .&  .&    .&    45&     .&  14856&14901\\
    &12&.& .& .& .& .& .& .& .&  .&  .&  .&    .&     .&    66&      .&66\\
    &13&.& .& .& .& .& .& .& .&  .&  .&  .&    .&     .&     .&    114&114\\
    \cline{2-18}
    & totals&0&0&0& 0& 0& 0& 3& 8& 42& 185& 888& 5110&27436&168030&1008906&1210608
  \end{tabular}

  \vspace{1\baselineskip}
  {\small Nonalternating}
  \vspace{2\baselineskip}

  \begin{tabular}{r@{\hspace{0pt}}r|*{15}{r}|r}
    \multicolumn{2}{c}{} & \multicolumn{15}{c}{\small $c$} & \\
    &   & 0& 3& 4& 5& 6& 7& 8&  9& 10&  11&  12&   13&   14&    15&      16&totals \\
    \cline{2-18}
    & 1& 1& .& .& .& .& .&  .& .&   .&  .&   .&   .&    .&     .&     . &1\\
    & 2& .& 1& .& 1& .& 1&  .&  1&  .&  1&   .&   1&    .&     1&     . &7\\
    & 3& .& .& 1& 1& 2& 2& 11&  4& 37&  6& 121&   .&  346&     .&  1224 &1755\\
    & 4& .& .& .& .& 1& 4&  7& 31& 74&279& 623&2264&    .& 23577&     . &26860\\
    \cline{17-17}
    & 5& .& .& .& .& .& .&  3& 13& 49&219&1138&   .&18008&\multicolumn{1}{r|}{.}&324694 &344124\\
    \cline{14-16}
    & 6& .& .& .& .& .& .&  .&  .&  5& 47&\multicolumn{1}{r|}{277}&5637&    .&132557&     . &138523\\
    {\small $s$}
    & 7& .& .& .& .& .& .&  .&  .&  .&  .&\multicolumn{1}{r|}{17}&   .&23068&     .&736707 &759792\\
    \cline{3-13}
    & 8& .& .& .& .& .& .&  .&  .&  .&  .&   .&1913&    .& 85899&     . &87812\\
    & 9& .& .& .& .& .& .&  .&  .&  .&  .&   .&   .& 5306&     .&297965 &303271\\
    &10& .& .& .& .& .& .&  .&  .&  .&  .&   .& 168&    .& 10859&     . &11027\\
    &11& .& .& .& .& .& .&  .&  .&  .&  .&   .&   .&  241&     .& 27553 &27794\\
    &12& .& .& .& .& .& .&  .&  .&  .&  .&   .&   5&    .&   400&     . &405\\
    &13& .& .& .& .& .& .&  .&  .&  .&  .&   .&   .&    3&     .&   562 &565\\
    \cline{2-18}
    & totals&1&1&1& 2& 3& 7& 21& 49& 165& 552& 2176& 9988&46972&253293&1388705&1701936
  \end{tabular}

  \vspace{1\baselineskip}
  {\small All}
  \vspace{1\baselineskip}

  \caption{%
    The sizes of $\mathbb{K}^A_{c,s}$, $\mathbb{K}^N_{c,s}$ and $\mathbb{K}_{c,s}$.
    Above the zigzags lie sets for which $LG$ has
    been evaluated using a state model method.
  }%
  \label{table:stringindicesofprimeknots}
  \end{centering}
\end{table}

Table~\ref{table:stringindicesofprimeknots} lists the sizes of the
$\mathbb{K}^P_{c,s}$. An artefact of the algorithm used to construct the
\textsc{K2K}-braids is that a \textsc{K2K}-braid corresponding to an HTW prime
knot of an odd (respectively even) number of crossings has an even
(respectively odd) string index.

We have evaluated $LG$ for a total of $53,418$ of these knot types (not
counting reflections), including all complete sets $\mathbb{K}^P_{c,s}$ lying
above the zigzags in Table~\ref{table:stringindicesofprimeknots}.

\begin{itemize}
\item
  The state model method has been applied to the \textsc{K2K}$'$-braids
  of a total of $48,399$ prime knots, namely those from:
  \begin{itemize}
  \item
    all $\mathbb{K}_c$ for $c\leqslant 12$, and
  \item
    all $\mathbb{K}_{c,s}$ for $c\leqslant 15$ and $s\leqslant 5$, and $\mathbb{K}_{16,3}$.%
    \footnote{%
      The $LG$-feasible $\mathbb{K}_{16,5}$ is omitted due to its
      size, as insufficient CPU-years were available.
    } %
  \end{itemize}

\item
  A formula of Ishii~\cite{Ishii:2004b} for $LG$ for $2$-bridge knots has been
  applied to the $5546$ (necessarily alternating) $2$-bridge knots within the
  HTW tables (as determined in~\cite{DeWit:TwoBridgeKnots}).
\end{itemize}

These classes have an overlap of only $527$ knots as the
\textsc{K2K}$'$-braids for the $2$-bridge knots are generally on more than $5$
strings. Thus, we have evaluations of $LG$ for a net $48,399+5546-527=53,418$
prime knots. Table~\ref{table:LG21knownnumbers} provides a breakdown of these
totals with crossing numbers.

\renewcommand{\arraystretch}{1.5}
\renewcommand{\tabcolsep}{4pt}

\begin{table}[htbp]
  \tiny
  \begin{centering}
  \begin{tabular}{r|*{15}{r}|r}
   \multicolumn{1}{r}{} & \multicolumn{15}{c}{$c$} \\
               & 0 & 3 & 4 & 5 & 6 & 7 &  8 &  9 &  10 &  11 &   12 &   13 &    14 &    15 &   16 & totals \\
   \hline
   state model & 1 & 1 & 1 & 2 & 3 & 7 & 21 & 49 & 165 & 552 & 2176 & 2265 & 18354 & 23578 & 1224 & 48399 \\
   $2$-bridge  & . & 1 & 1 & 2 & 3 & 7 & 12 & 24 &  45 &  91 &  176 &  352 &   693 &  1387 & 2752 &  5546 \\
   overlap     & . & 1 & 1 & 2 & 3 & 7 & 12 & 24 &  45 &  91 &  176 &   26 &    95 &    37 &    7 &   527 \\
   \hline
   net         & 1 & 1 & 1 & 2 & 3 & 7 & 21 & 49 & 165 & 552 & 2176 & 2591 & 18952 & 24928 & 3969 & 53418
   \end{tabular}
  \caption{%
    Numbers of $c$-crossing prime knots for which $LG$ has been evaluated.
  }%
  \label{table:LG21knownnumbers}
  \end{centering}
\end{table}

For comparison with Table~\ref{table:stringindicesofprimeknots},
Table~\ref{table:stringindicesoftwobridgeknots} presents the numbers of
$c$-crossing $2$-bridge knots with $s$-string \textsc{K2K}$'$-braids.

\renewcommand{\arraystretch}{1.5}
\renewcommand{\tabcolsep}{4pt}

\begin{table}[hbt]
  \tiny
  \begin{centering}
  \begin{tabular}{r@{\hspace{0pt}}r|*{14}{r}|r}
    \multicolumn{2}{c}{} & \multicolumn{14}{c}{\small $c$} \\
    & & 3& 4& 5& 6& 7& 8&  9& 10& 11& 12& 13& 14& 15&   16 &totals\\
    \cline{2-17}
    & 2& 1& .& 1& .& 1& .&  1&  .&  1&  .&  1&  .&  1&   .&    7\\
    & 3& .& 1& 1& 2& 2& 3&  3&  4&  4&  5&  .&  6&  .&   7&   38\\
    & 4& .& .& .& 1& 4& 6& 12& 15& 24& 28& 25&  .& 36&   .&  151\\
    & 5& .& .& .& .& .& 3&  8& 22& 40& 73&  .& 89&  .& 151&  386\\
    & 6& .& .& .& .& .& .&  .&  4& 22& 60&123&  .&276&   .&  485\\
    {\small $s$}
    & 7& .& .& .& .& .& .&  .&  .&  .& 10&  .&280&  .& 736& 1026\\
    & 8& .& .& .& .& .& .&  .&  .&  .&  .&151&  .&592&   .&  743\\
    & 9& .& .& .& .& .& .&  .&  .&  .&  .&  .&251&  .&1145& 1396\\
    &10& .& .& .& .& .& .&  .&  .&  .&  .& 49&  .&400&   .&  449\\
    &11& .& .& .& .& .& .&  .&  .&  .&  .&  .& 64&  .& 610&  674\\
    &12& .& .& .& .& .& .&  .&  .&  .&  .&  3&  .& 82&   .&   85\\
    &13& .& .& .& .& .& .&  .&  .&  .&  .&  .&  3&  .& 103&  106\\
    \cline{2-17}
    & totals& 1& 1& 2& 3&7& 12& 24&45 & 91&176&352&693&1387&2752&5546 \\
 \end{tabular}
  \caption{%
    The numbers of $c$-crossing $2$-bridge knots with $s$-string
    \textsc{K2K}$'$-braids.
  }%
  \label{table:stringindicesoftwobridgeknots}
  \end{centering}
\end{table}

\begin{figure}[ht]
  \begin{centering}
  \includegraphics[width=120pt]{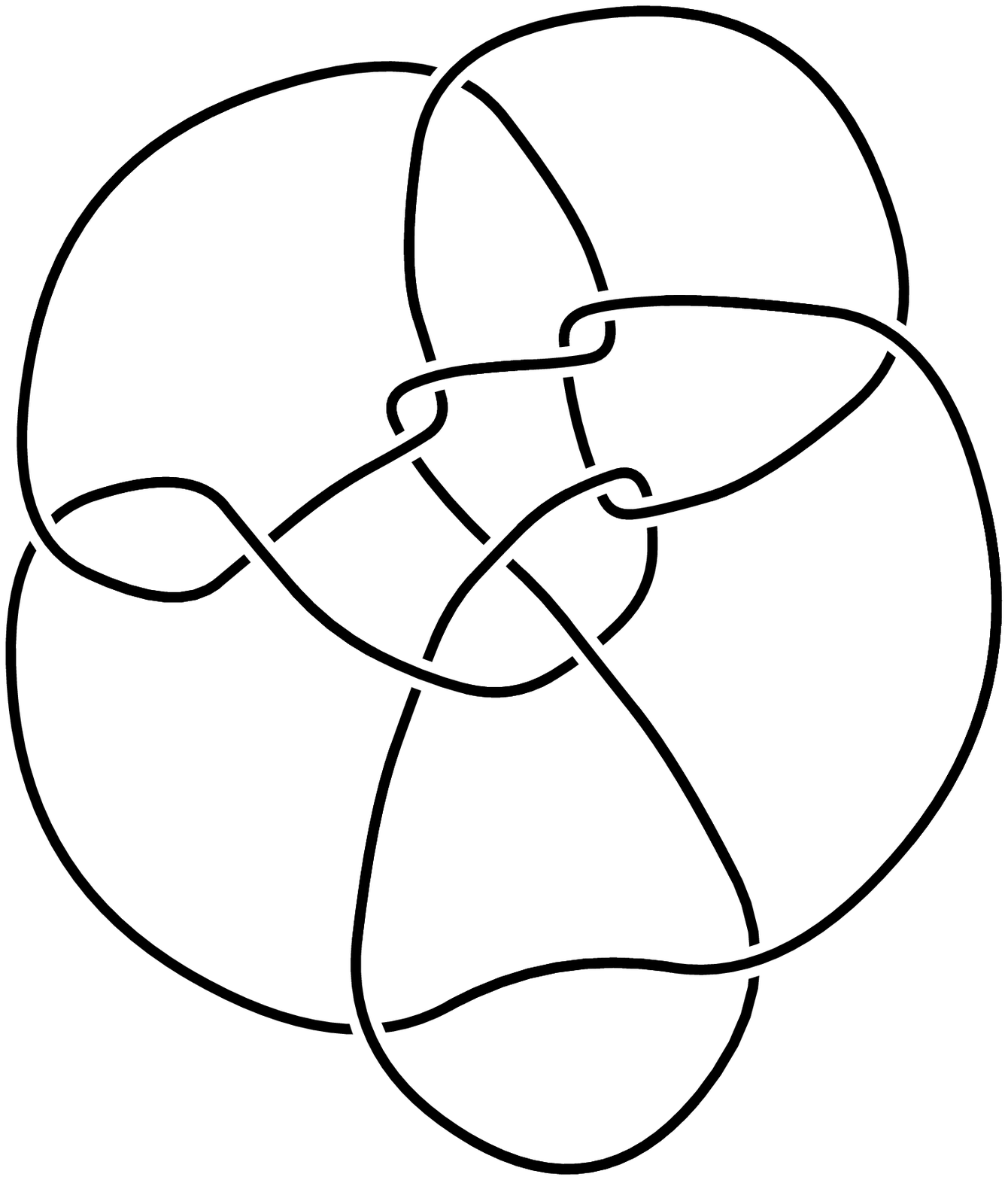}
  \caption{%
    $15^N_{139717}$: The first achiral knot of odd crossing number.
  }
  \label{figure:15N139717}
  \end{centering}
\end{figure}

We illustrate an $LG$-polynomial using $15^N_{139717}$, the only achiral knot
of odd crossing number in the HTW tables, depicted in
Figure~\ref{figure:15N139717}.%
\footnote{
 This knot is hyperbolic, and has symmetry group $D_1$
 (see~\S\ref{section:MutantCliques}
 and~\S\ref{section:WhereLGfirstfailstodistinguishreflections}).
} %
Although its \textsc{K2K}-braid is on $6$ strings, it is still $LG$-feasible.
We obtain the following polynomial, where we have written
$\overline{q}\triangleq q^{-1}$ and $\overline{p}\triangleq p^{-1}$. The
polynomial is palindromic in $p$, as are all $LG$-polynomials. It is also
palindromic in $q$, as the knot is achiral.

\tiny

\begin{eqnarray*}
  \begin{array}{@{\hspace{0pt}}r@{\hspace{5pt}}r@{\hspace{2pt}}l}
      &                              & (2\overline{q}^6 - 10\overline{q}^4 + 72\overline{q}^2 + 181 + 72 q^2 - 10 q^4 + 2 q^6)
    \\
    - &       (\overline{p}^2 + p^2) & (\overline{q}^7 - 4\overline{q}^5 + 8\overline{q}^3 + 118\overline{q} + 118 q + 8 q^3 - 4 q^5 + q^7)
    \\
    + &       (\overline{p}^4 + p^4) & (- 4\overline{q}^4 + 32\overline{q}^2 + 80 + 32 q^2 - 4 q^4)
    \\
    - &       (\overline{p}^6 + p^6) & (7\overline{q}^3 + 30\overline{q} + 30 q + 7 q^3)
    \\
    + &       (\overline{p}^8 + p^8) & (2\overline{q}^4 + 11\overline{q}^2 + 20 + 11 q^2 + 2 q^4)
    \\
    - & (\overline{p}^{10} + p^{10}) & (-\overline{q}^5 + \overline{q}^3 + 10\overline{q} + 10 q + q^3 - q^5)
    \\
    + & (\overline{p}^{12} + p^{12}) & (- 2\overline{q}^4 + \overline{q}^2 + 7 + q^2 - 2 q^4)
    \\
    - & (\overline{p}^{14} + p^{14}) & (-\overline{q}^3 + 2\overline{q} + 2 q - q^3)
    \\
    + & (\overline{p}^{16} + p^{16}) & (1)
  \end{array}
\end{eqnarray*}

\normalsize

Note that $LG$-polynomials often have a regular sign pattern, in that the
coefficients within the $q$-polynomial corresponding to any given power of $p$
are all of the same sign, and this sign alternates with rising (even) powers of
$p$. Examples of this pattern are seen in $LG$-polynomials presented
in~\S\ref{section:WhereLGfirstfailstodistinguishprimeknots}. It is not known
what property of a link causes this pattern. Observe that $LG_{15^N_{139717}}$
is a (slight) exception to this rule. (The prime knots of up to $10$ crossings
with such exceptional $LG$-polynomials are identified in~\cite{DeWit:2000}.) It
does, however, satisfy another general (and still unexplained) property of
$LG$-polynomials: Each of its $q$-polynomials contains either only odd or only
even powers of $q$.

We conclude this section with a few remarks about the computability of
$LG^{3,1}$, the next available Links--Gould invariant. The
$LG^{3,1}$-polynomials are generally of higher degree than those of
$LG^{2,1}$, and are correspondingly more expensive to evaluate. Using our
state model method we have only been able to evaluate $LG^{3,1}$ for prime
knots in the following complete sets $\mathbb{K}_{c,s}$, which together
contain only $380$ knot types:
\begin{itemize}
\item
  all $\mathbb{K}_c$ for $c\leqslant 10$, and
\item
  all $\mathbb{K}_{c,s}$ for $c\leqslant 12$ and $s\leqslant 3$, and
\item
  $\mathbb{K}_{13,2}\equiv\{13^A_{4878}\}$ and  $\mathbb{K}_{15,2}\equiv\{15^A_{85263}\}$.
\end{itemize}
We find that $LG^{3,1}$ is complete for these knots (including reflections).


\section{Mutant cliques within the HTW tables}
\label{section:MutantCliques}

Recall that a (Conway) mutation of a $(2,2)$-tangle (and thence of a link in
which the tangle is embedded) amounts to a rotation of $\pi$ about one of three
orthogonal axes of the tangle. More generally, the result of a sequence of
mutations is regarded as a mutation in itself, and as mutation is reflexive,
symmetric and transitive, we may discuss equivalence classes of knots under
mutation.

Many link invariants fail to distinguish mutants: the list includes each
$LG^{m,n}$ together with all the well-known polynomial link invariants
(see~\S\ref{section:WhereLGfirstfailstodistinguishprimeknots}), and sundry
others including the hyperbolic volume (see below) and Khovanov
homology~\cite{Wehrli:2003}. With a view to identifying where $LG$ fails to
identify nonmutant prime knots, we here describe a partial classification of
the mutants within the HTW tables. Specifically, we classify those within
$\mathbb{K}_{11}$ and $\mathbb{K}_{12}$.  Of course, we need only describe the
equivalence classes of size greater than $1$, and these we refer to as
\emph{mutant cliques} (within the HTW tables). Although the identification
turns out to be relatively straightforward to obtain, this information
apparently does not appear in the existing literature. This situation reflects
a remarkable lack of research into mutation, the notable exception being the
search for link invariants which are sensitive to it, for
example~\cite{MortonCromwell:1996}.

We determine mutant cliques by first using mutation-insensitive invariants to
filter the HTW knots into candidate mutant cliques, and then inspect the
diagrams to determine if mutations are visible between their elements. The
inspection process is nonalgorithmic, and may fail to identify some genuine
mutants, and cannot demonstrate when candidates are \emph{not} mutants.
Nevertheless, in this manner, we obtain a complete classification of the
mutant cliques within $\mathbb{K}_{11}$ and $\mathbb{K}_{12}$. Moreover, the
filtering allows us to make some qualitative statements about the remaining
undetermined mutant cliques of higher crossing numbers.

\renewcommand{\baselinestretch}{1.2}

The main invariant we use to determine candidate mutant cliques is the
hyperbolic volume, so we digress to review it. To that end, if the complement
$S^3 - L$ of a link $L$ can be assigned a hyperbolic metric, it can be
decomposed with respect to that metric, and $L$ is called \emph{hyperbolic}.%
\footnote{%
  Nonhyperbolic knots are rare in the HTW tables. Including the unknot, there
  are only $33$ of them (out of $1,701,936$). They include $13$ torus knots,
  which always form only a small proportion of all prime knots of any given crossing
  number, and $20$ satellite knots (it appears that these will also always be
  a small number of the total of any given crossing number). They are:

  \begin{itemize}
  \item
    $13$ torus knots, which include the unknot $0^A_1$ (which may be regarded as the $(2,1)$
    torus knot);
    the closed odd $2$-braids $3^A_1$, $5^A_2$, $7^A_7$, $11^A_{367}$, $13^A_{4878}$, and
    $15^A_{85263}$; and sundries $8^N_3$, $9^A_{41}$, $10^N_{21}$, $14^N_{21881}$,
    $15^N_{41185}$, $16^N_{783154}$; together with
  \item
    $20$ satellites of the trefoil: $13^N_{4587}$, $13^N_{4639}$,
    $14^N_{22180}$, $14^N_{26039}$, $15^N_{40211}$, $15^N_{59184}$,
    $15^N_{115646}$, $15^N_{124802}$, $15^N_{142188}$, $15^N_{156076}$,
    $16^N_{253502}$, $16^N_{400459}$, $16^N_{696530}$, $16^N_{697612}$,
    $16^N_{703714}$, $16^N_{703716}$, $16^N_{800356}$, $16^N_{800378}$,
    $16^N_{958969}$, $16^N_{958982}$.
  \end{itemize}
  Importantly, note that there are no $12$-crossing nonhyperbolic knots, and
  that the only $11$-crossing example is easily identified.
} %
In such a case, the complement can be assigned a \emph{hyperbolic volume},
which is an (algebraic) positive real number. This quantity may be regarded as
a property of $L$ itself, and is an invariant of $L$, for the following reason.

\renewcommand{\baselinestretch}{1}

The Gordon--Luecke theorem~\cite{GordonLuecke:1989} tells us that topologically
equivalent (true) knots have homeomorphic knot complements and vice-versa.
(Unfortunately, this statement does not hold more generally for multicomponent
links.) It follows that the hyperbolic decomposition of the knot complement
distinguishes all hyperbolic knots. In particular, the hyperbolic
decomposition distinguishes \emph{mutant} hyperbolic knots. Unfortunately
however, the hyperbolic \emph{volume} alone does \emph{not} distinguish mutant
hyperbolic knots. As composite knots are nonhyperbolic, a mutation of a
hyperbolic prime knot is necessarily a (hyperbolic) prime knot, rather than a
composite knot. It appears that we do not know whether a mutation of a
nonhyperbolic prime knot is a nonhyperbolic prime knot (as opposed to a
composite knot). The software available for determining hyperbolic
decompositions, and thus hyperbolic volumes, is the \textsc{SnapPea} program by
Weeks~\cite{HosteThistlethwaiteWeeks:1998}; we use it via its embedding within
\textsc{Knotscape}.

The hyperbolic volume is a powerful invariant, and generally distinguishes many
more knots than polynomial invariants, however, distinct knots with the same
hyperbolic volume are not \emph{necessarily} mutants. Indeed, it sometimes
fails to distinguish between prime knots of different crossing number, for
example $(5^A_1,12^N_{242})$, which are distinguished by both the HOMFLYPT and
the Kauffman polynomials. Moreover \textsc{SnapPea} only determines a
hyperbolic volume as a decimal approximation (to $10$ digits of accuracy and
precision) to its true value. Thus, as for polynomial invariants, whilst it can
be used to prove that given knots are not mutants, it can only be used as an
\emph{indicator} of possible mutants.

We continue the main thread of our discussion with some other salient facts
about mutation.

\begin{itemize}
\item
  A mutation of a true knot is a true knot, rather than a multicomponent link.
  No mutation of a nontrivial knot can create the
  unknot~\cite[pp50--51]{Adams:1994}.

\item
  Although a mutation of an alternating projection of an (alternating) knot is
  always an alternating knot~\cite[pp50--51]{Adams:1994}, it is apparently not
  known whether a mutation of a nonalternating projection of an alternating
  knot is always an alternating knot. However, at least for prime knots of low
  crossing numbers, mutant cliques appear to only contain knots of common
  `alternatingness'.

\item
  Mutation certainly preserves the numbers of crossings of link
  \emph{projections}, and it can of course leave a link unchanged. It appears
  that the question of whether mutation can convert a chiral link into its
  reflection is still unanswered; in any case we know no explicit examples of
  this phenomenon, and as we are only considering knot \emph{types}, we
  may ignore this issue.

  The first point at which mutation changes a prime knot into a distinct knot
  is within $\mathbb{K}_{11}$~\cite[p44]{Kawauchi:1996} (recall the
  Kinoshita--Terasaka pair $(11^N_{34},11^N_{42})$). In fact, at least for
  prime knots of low crossing numbers, mutation often preserves crossing
  number. This may not more generally be the case, as crossing number is not
  an invariant defined by a symmetry.
\end{itemize}

%
%
%
%
%

To determine candidate mutant cliques within the HTW tables, we first reduce
the full set of HTW knots into equivalence classes of common (apparent)
hyperbolic volume, and then further reduce these equivalence classes with the
HOMFLYPT and Kauffman polynomials, only retaining sets of size greater than
$1$. (For good measure, we further filter the candidate sets with $LG$ where
possible.) By inspection, each such candidate mutant clique lies within a
particular class $\mathbb{K}^P_c$, so every genuine mutant clique also lies
within a class $\mathbb{K}^P_c$. For the prime knots of up to $16$ crossings,
mutation preserves crossing number and alternatingness.

\begin{figure}[htbp]
  \begin{centering}
  \begin{tabular}{cc@{\hspace{10pt}}|@{\hspace{10pt}}cc}
    \includegraphics[width=75pt]{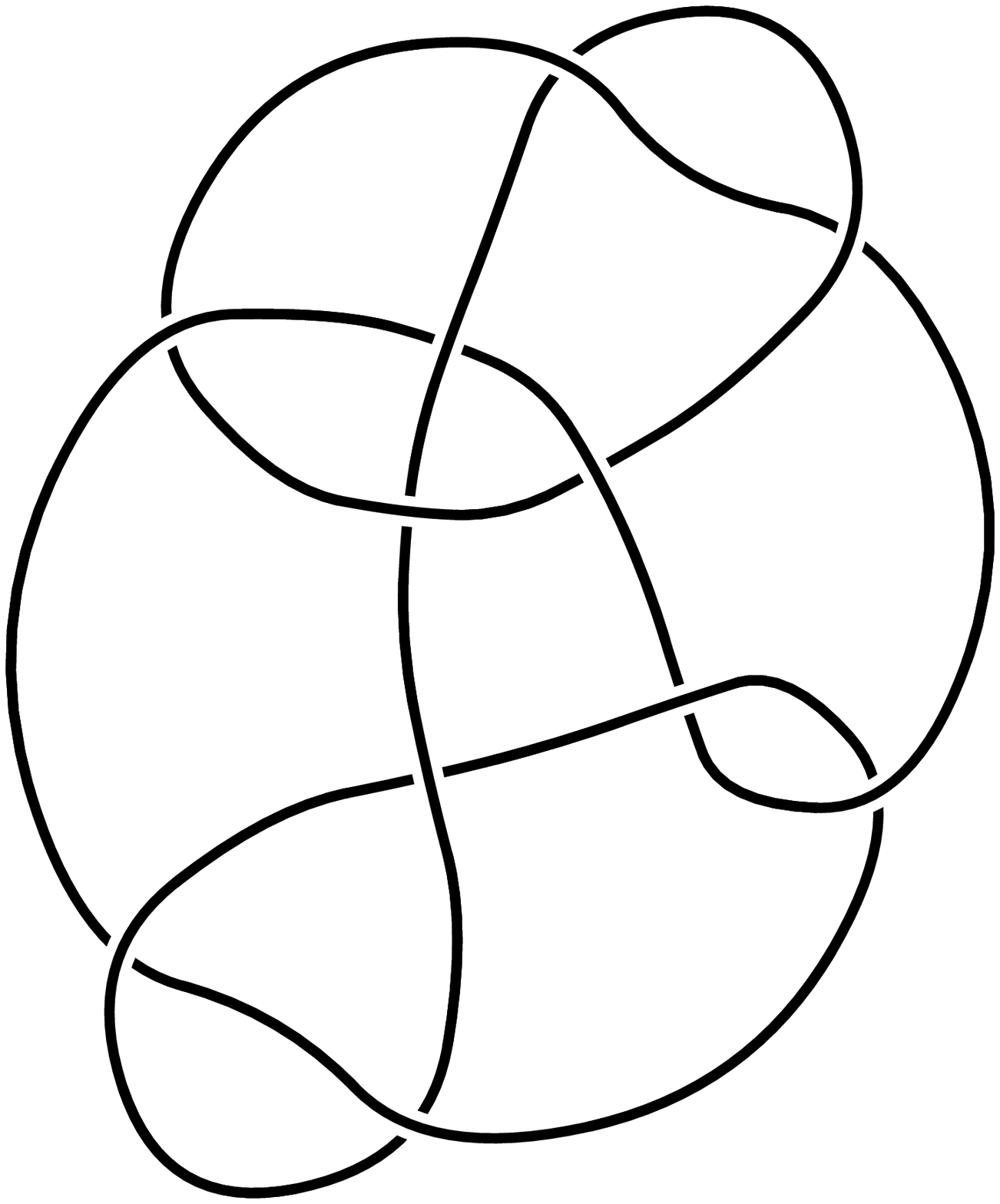}
    &
    \includegraphics[width=75pt]{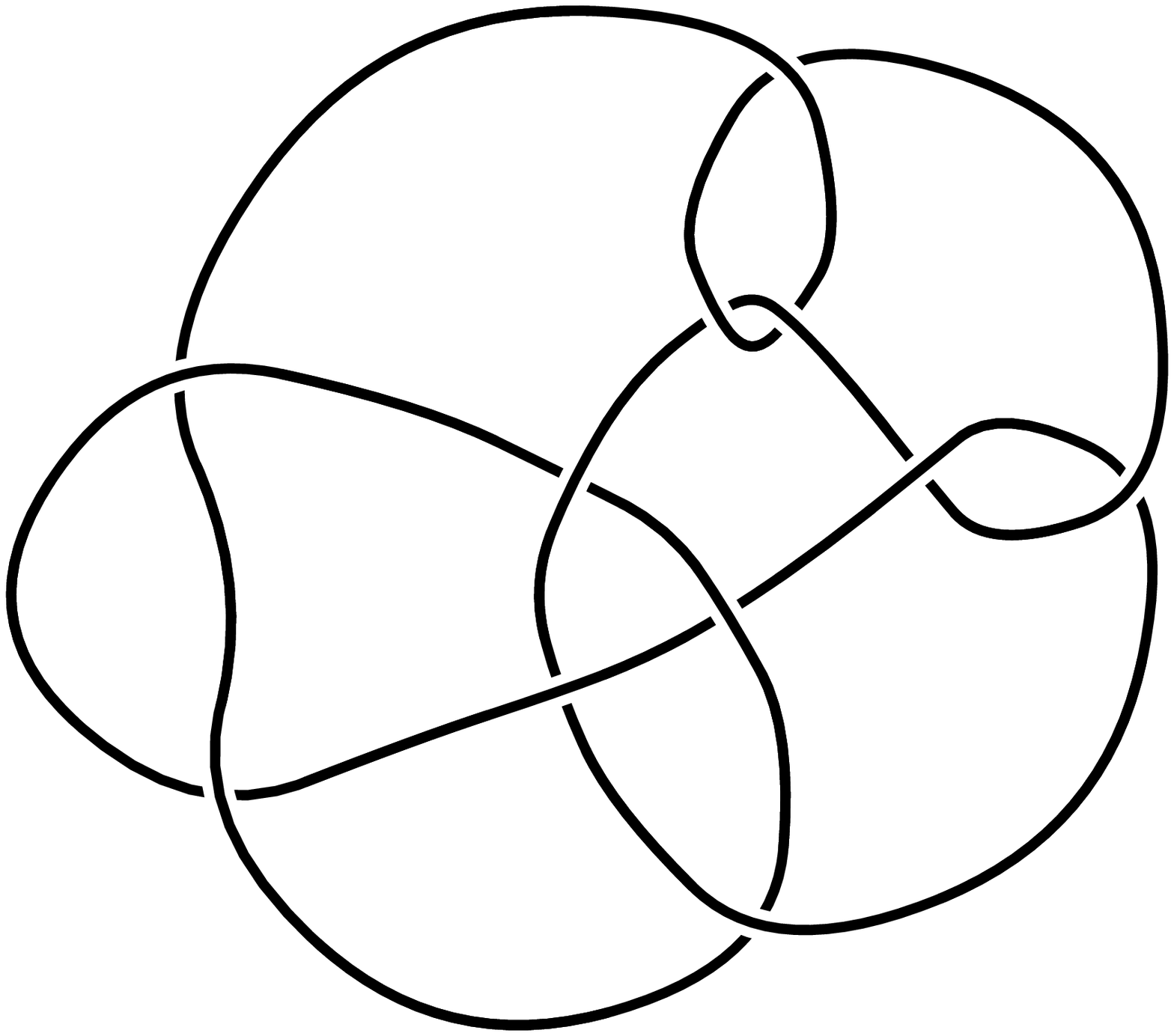}
    &
    \includegraphics[width=75pt]{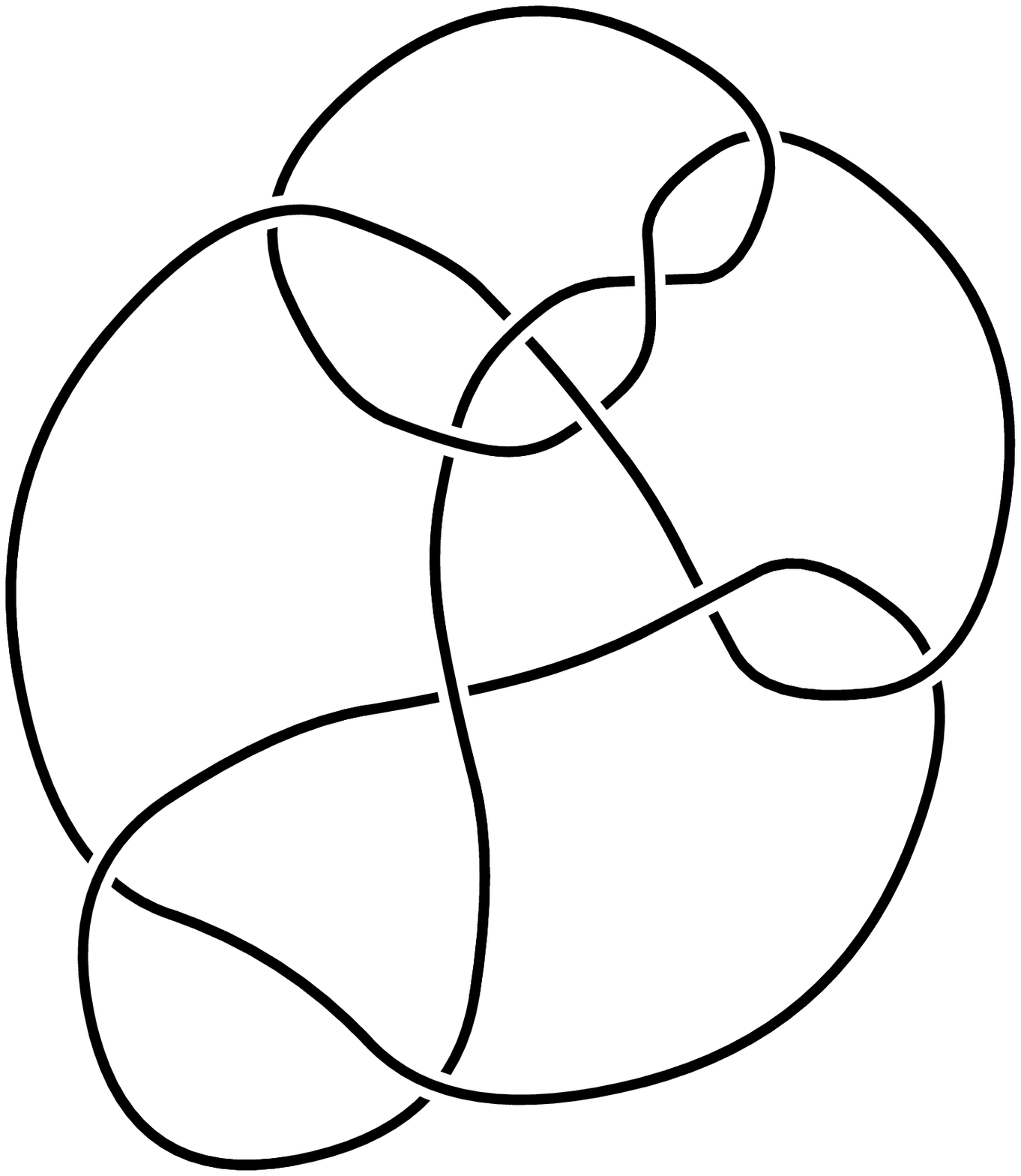}
    &
    \includegraphics[width=75pt]{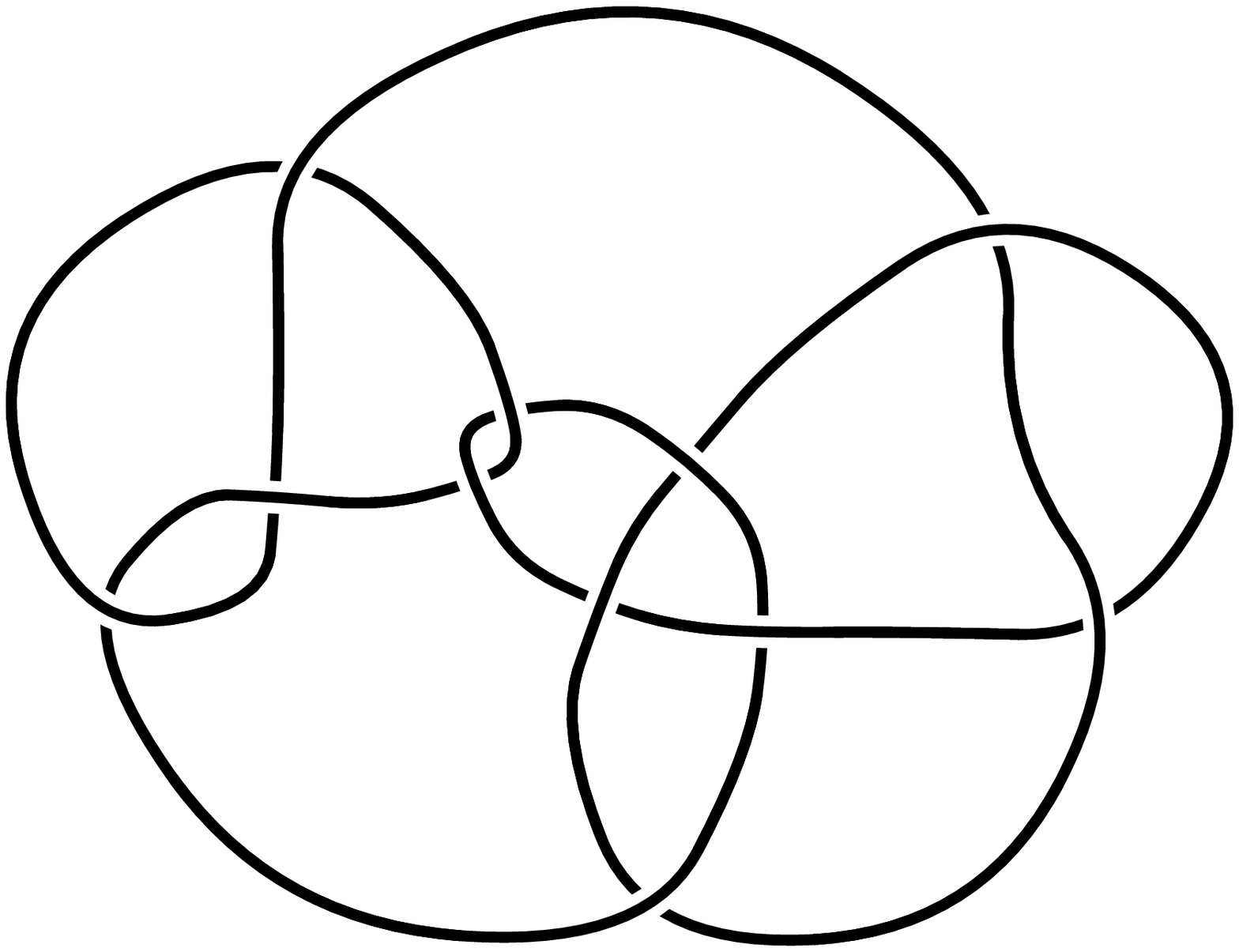}
    \\[-10pt]
    $11^A_{19}$ & $11^A_{25}$ & $11^A_{24}$ & $11^A_{26}$
    \\[10pt]
    \hline
    &&&\\[-10pt]
    \includegraphics[width=75pt]{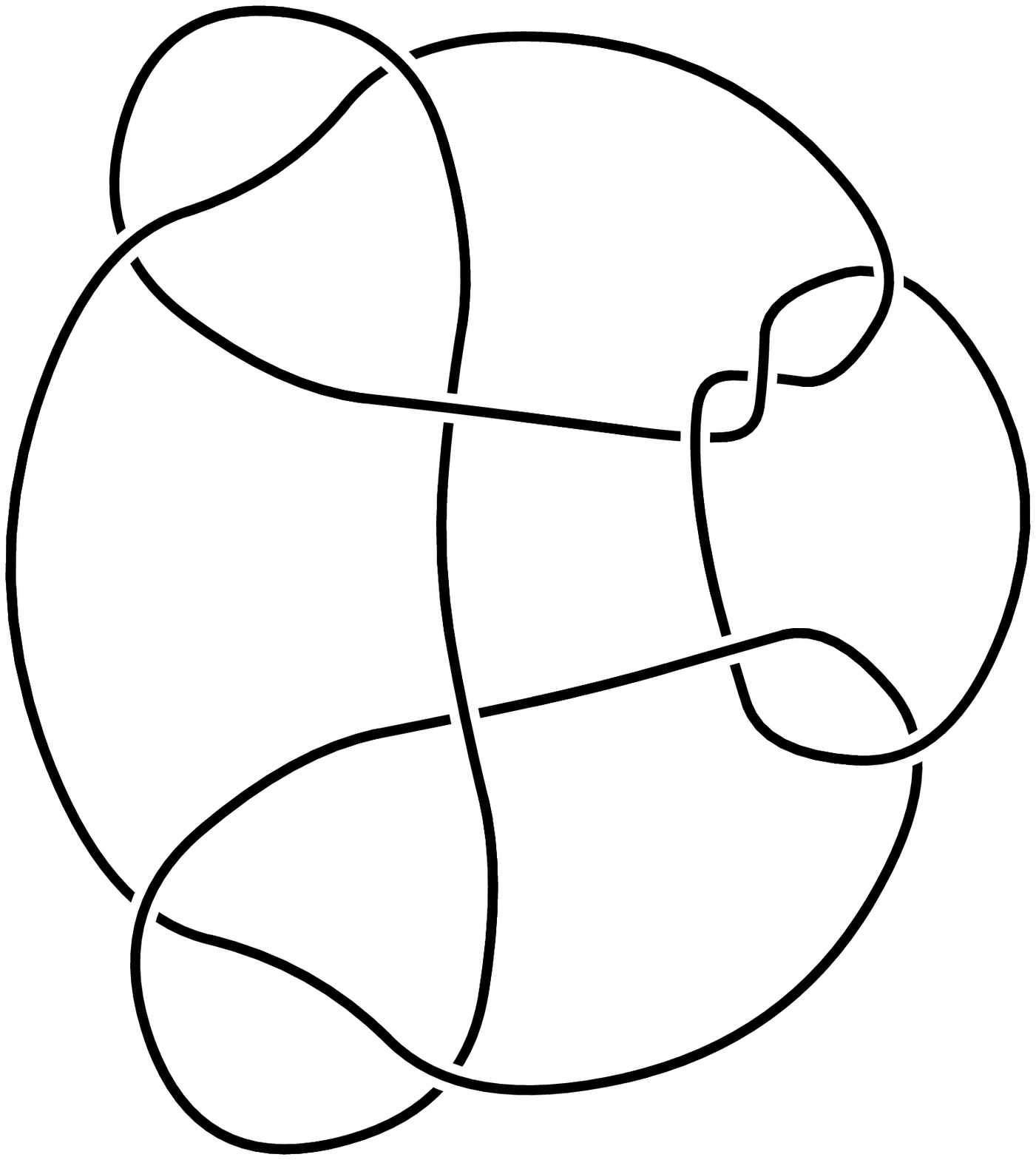}
    &
    \includegraphics[width=75pt]{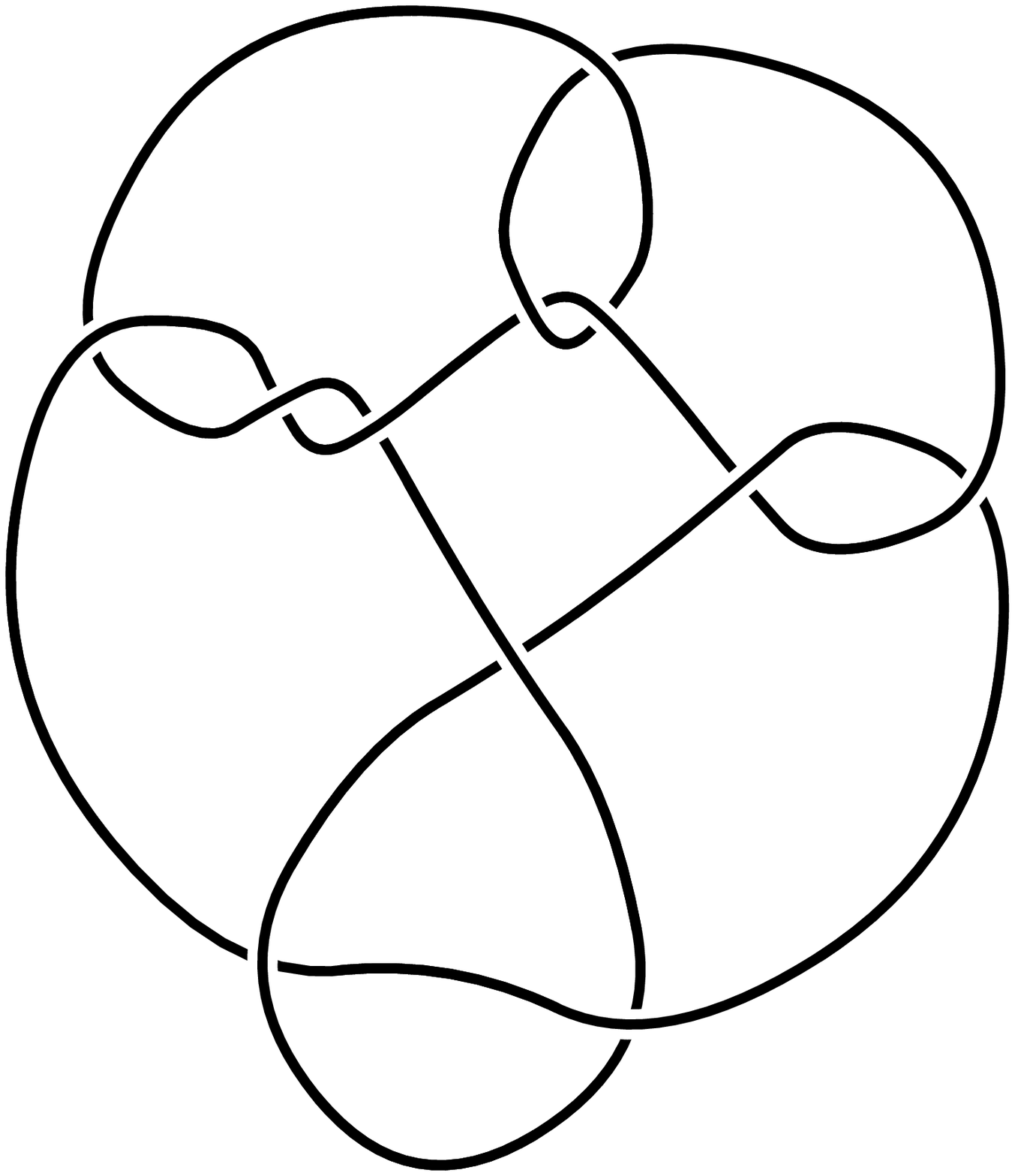}
    &
    \includegraphics[width=75pt]{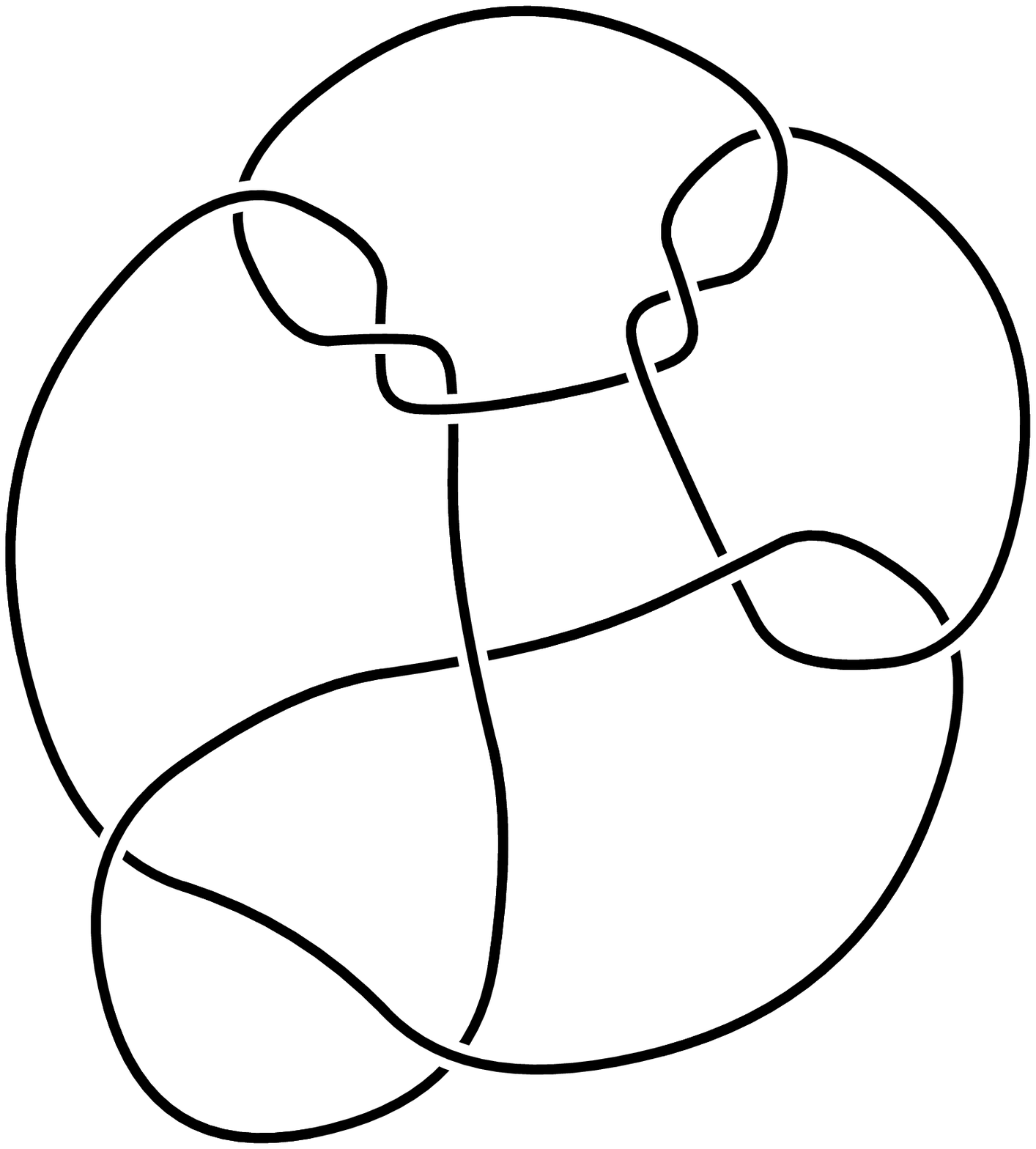}
    &
    \includegraphics[width=75pt]{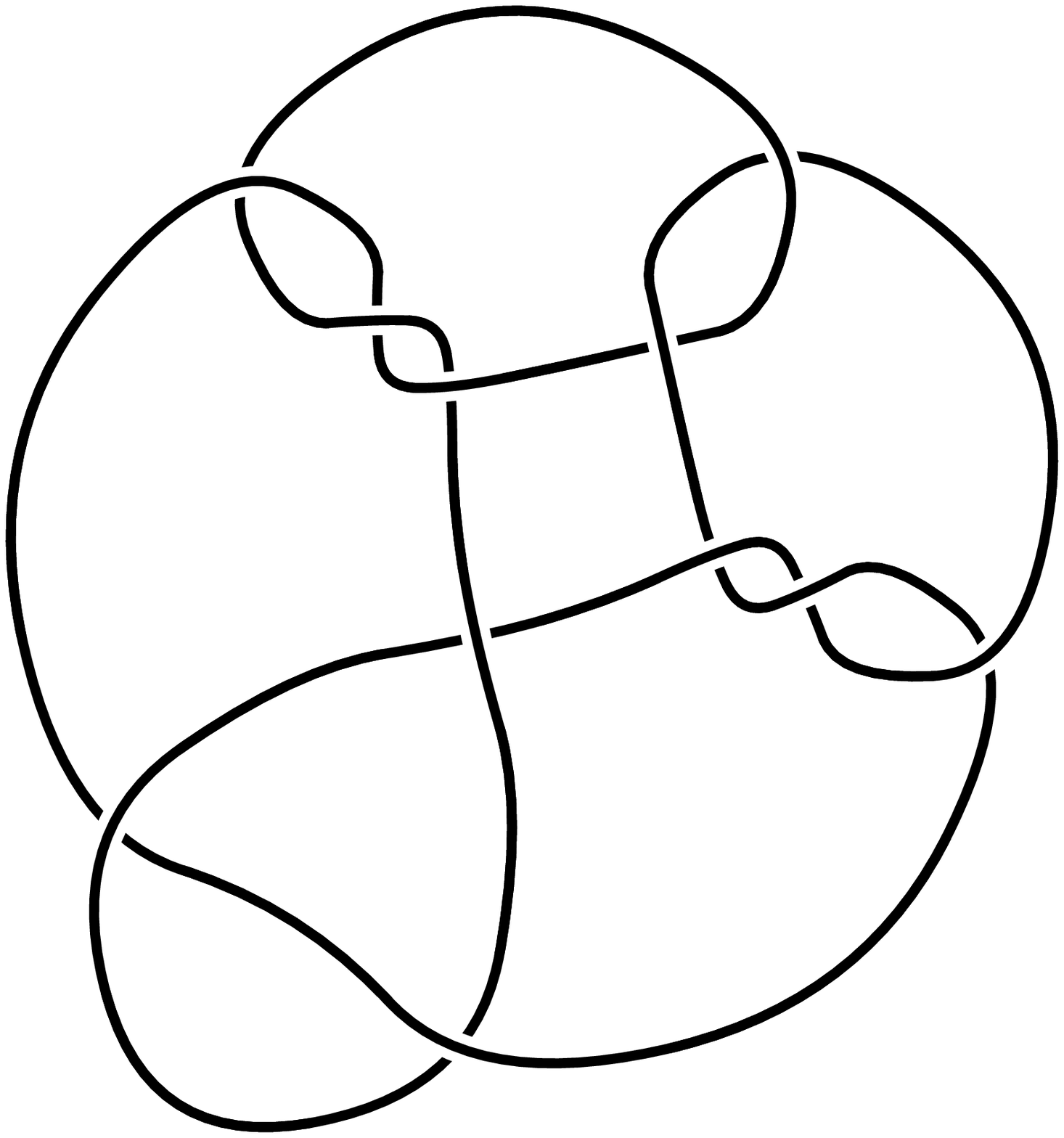}
    \\[-10pt]
    $11^A_{44}$ & $11^A_{47}$ & $11^A_{57}$ & $11^A_{231}$
    \\[10pt]
    \hline
    &&&\\[-10pt]
    \includegraphics[width=75pt]{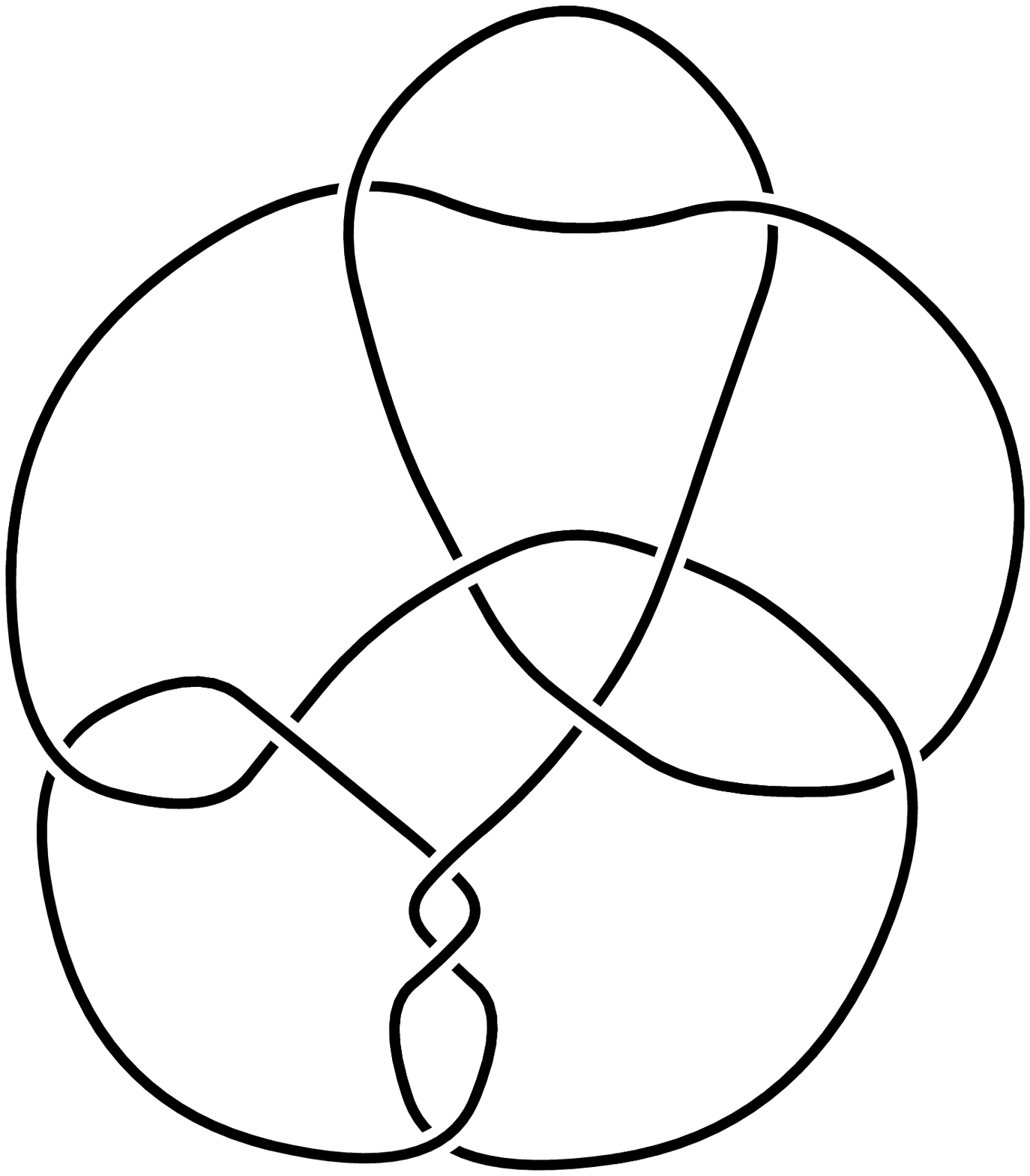}
    &
    \includegraphics[width=75pt,angle=90]{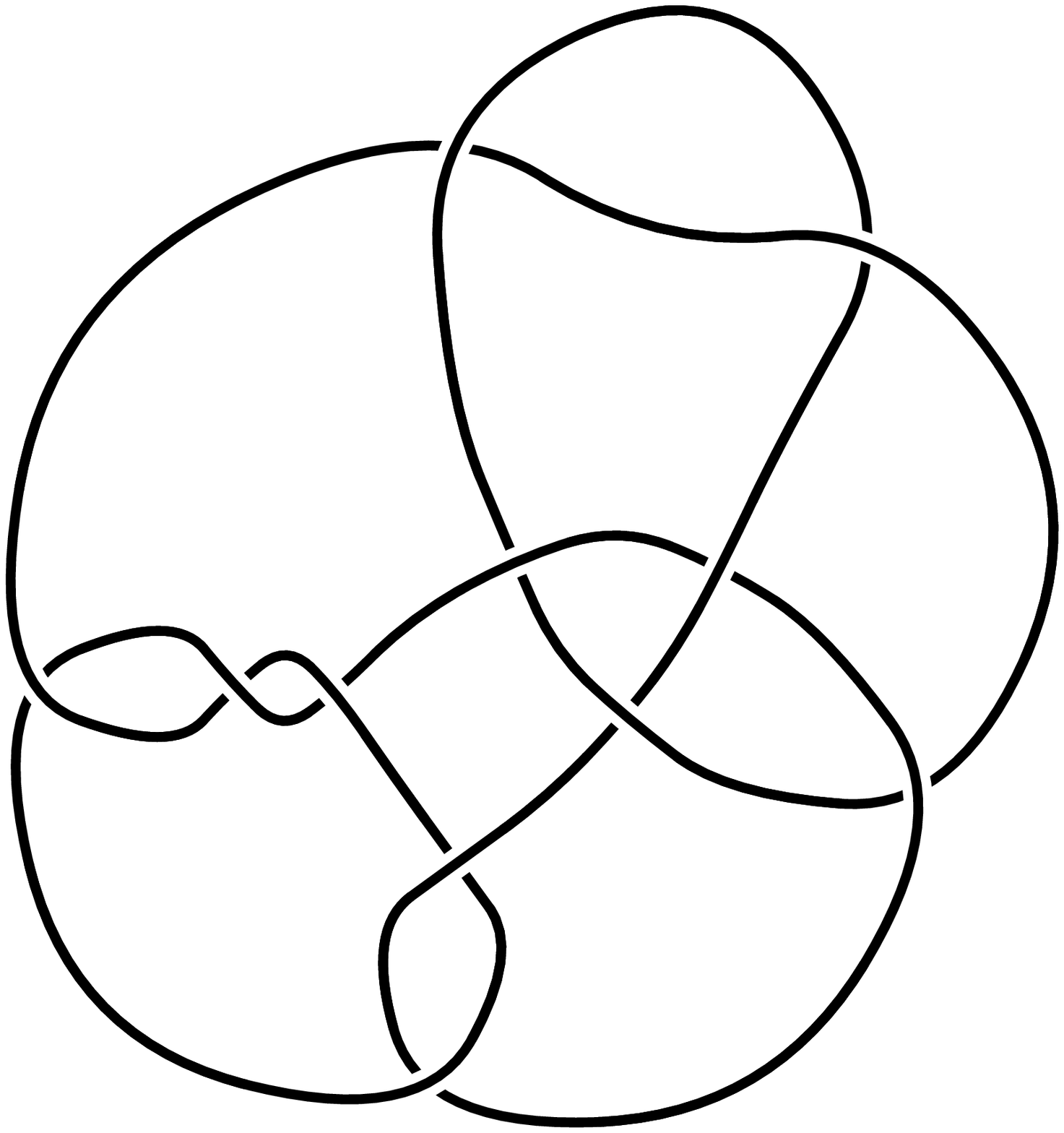}
    &
    \includegraphics[width=75pt]{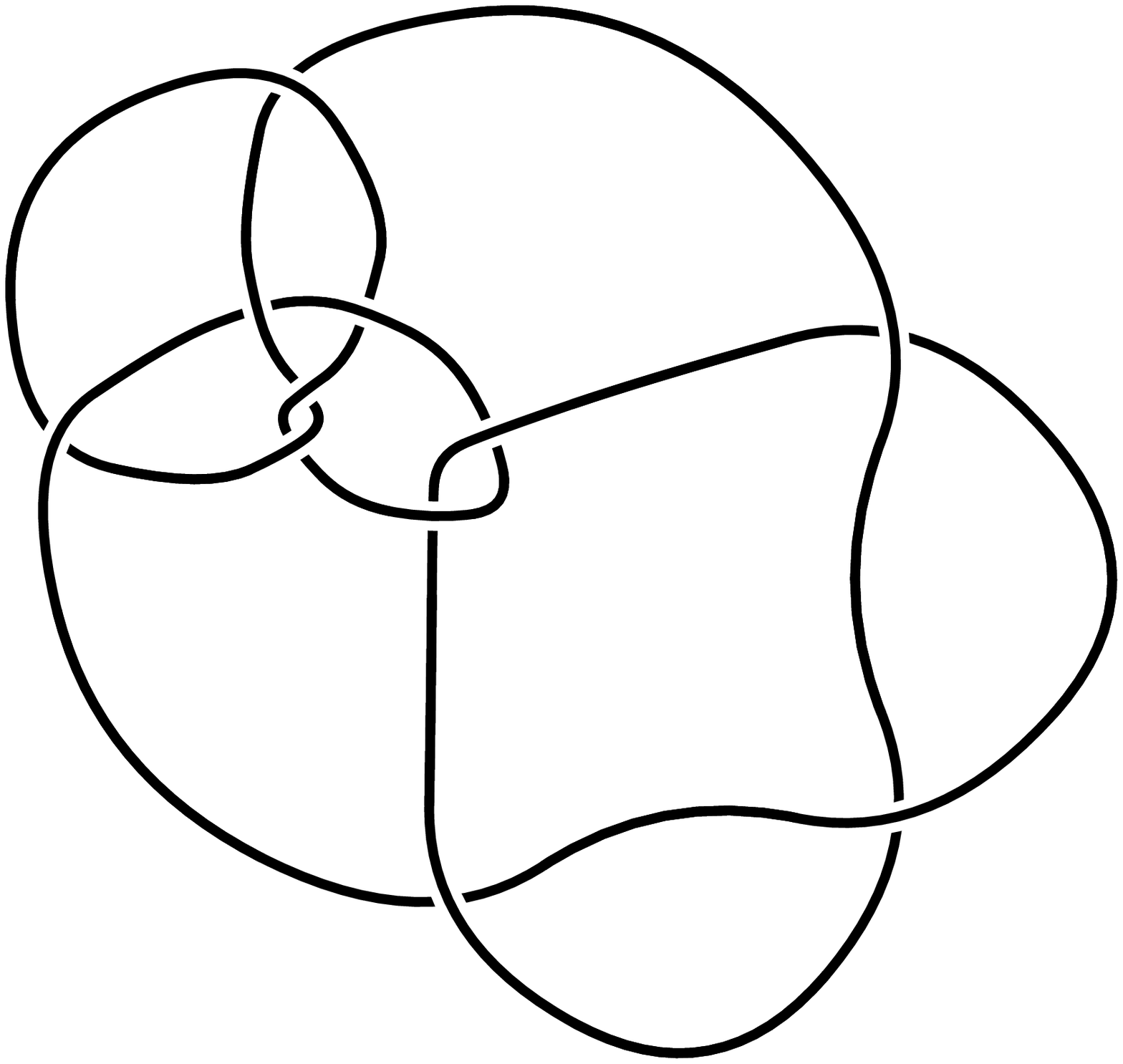}
    &
    \includegraphics[width=75pt]{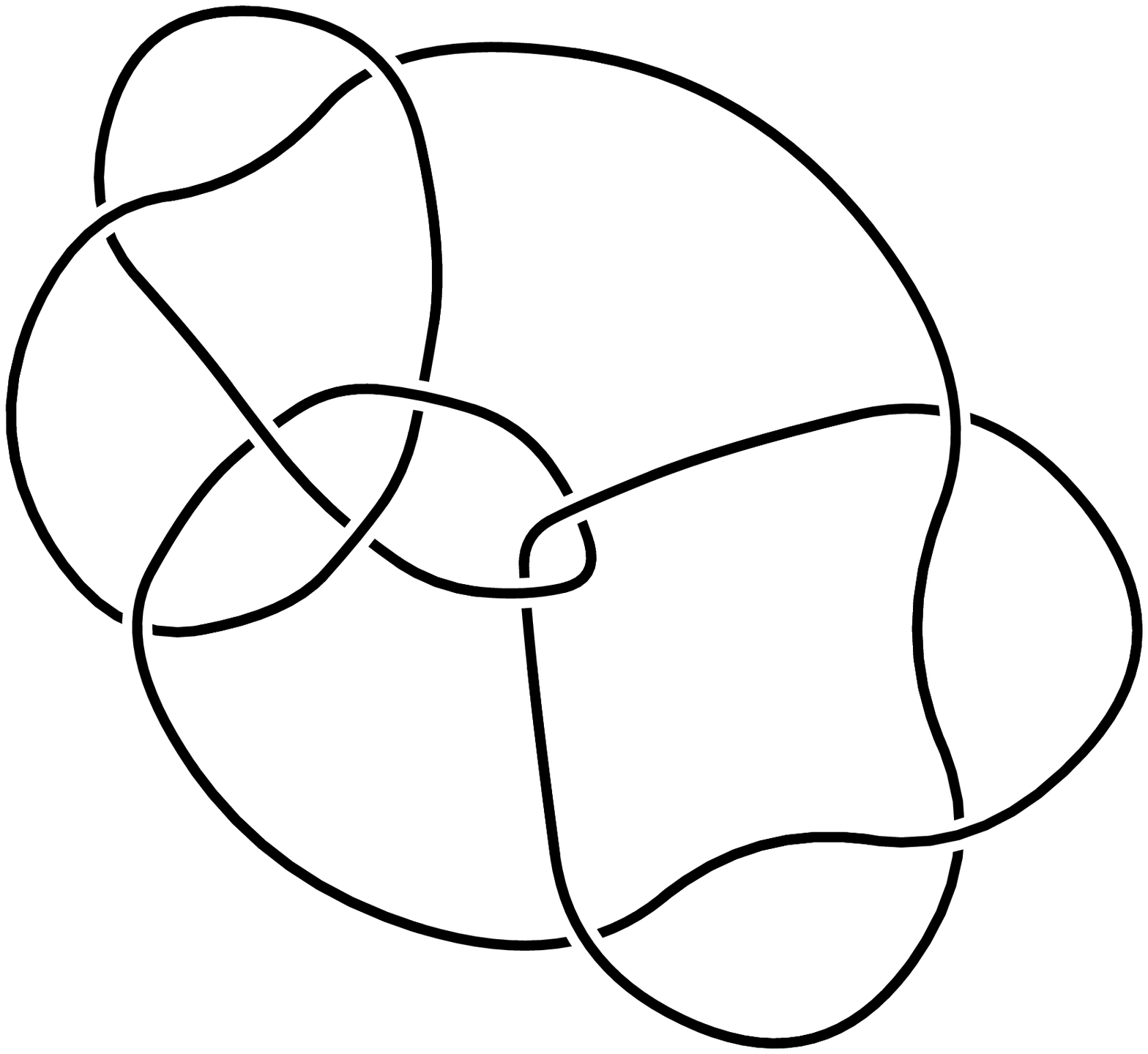}
    \\[-10pt]
    $11^A_{251}$ & $11^A_{253}$ & $11^A_{252}$ & $11^A_{254}$
  \end{tabular}
  \caption{Mutant cliques within $\mathbb{K}^A_{11}$.}
  \label{figure:Alternating11crossingmutantcliques}
  \end{centering}
\end{figure}

\begin{figure}[htbp]
  \begin{centering}
  \begin{tabular}{cc@{\hspace{10pt}}|@{\hspace{10pt}}cc}
    \includegraphics[width=75pt]{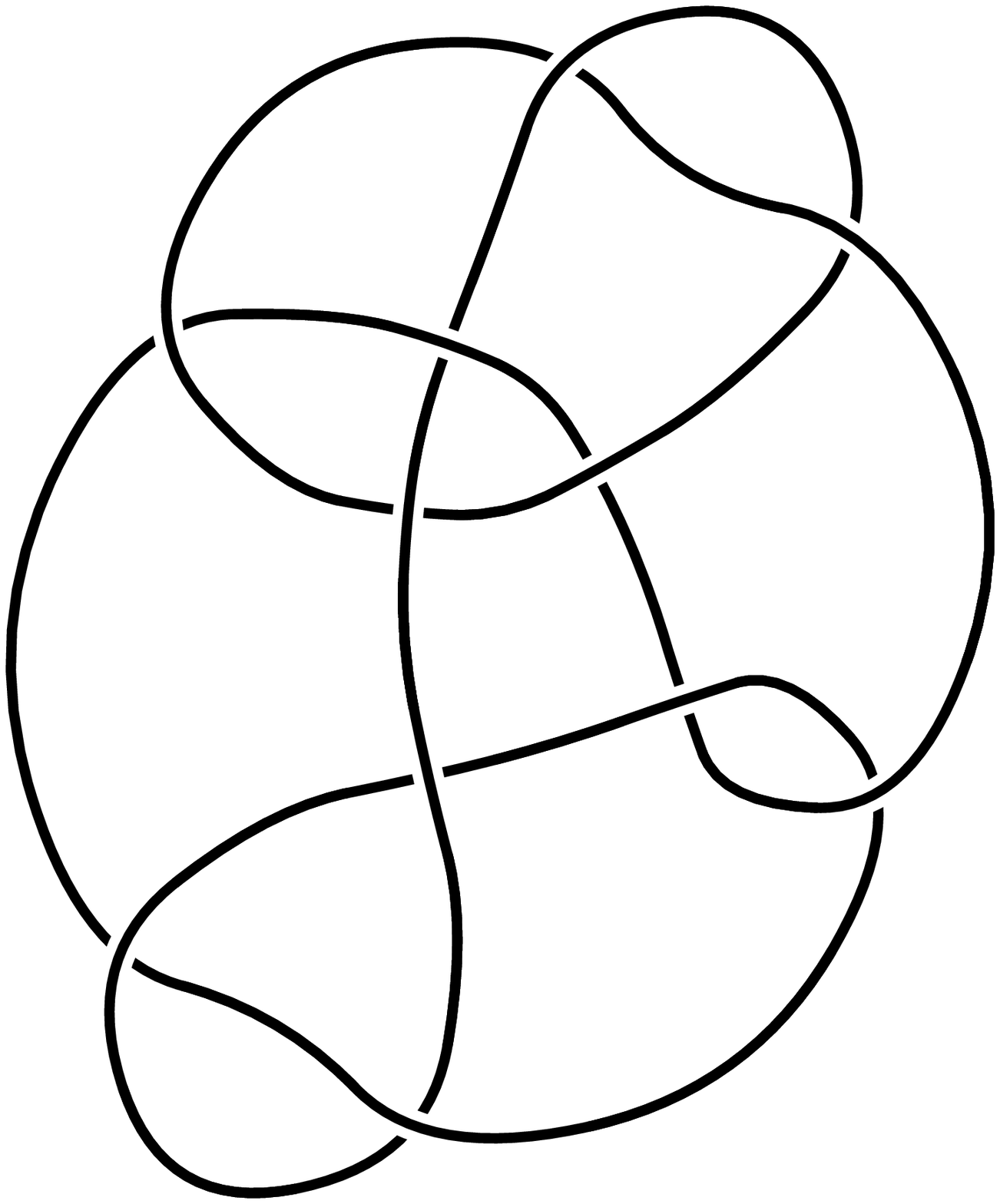}
    &
    \includegraphics[width=75pt]{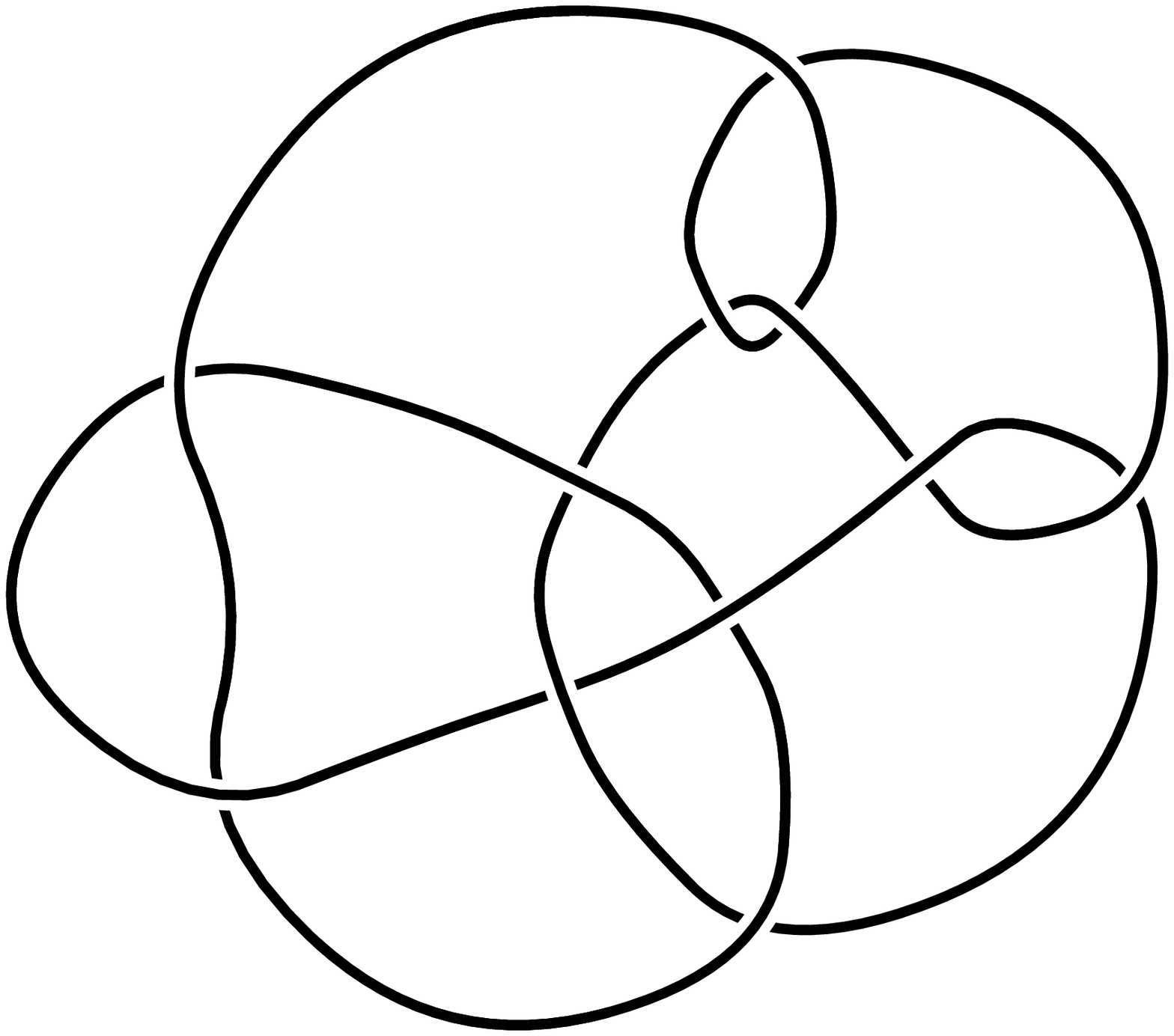}
    &
    \includegraphics[width=75pt]{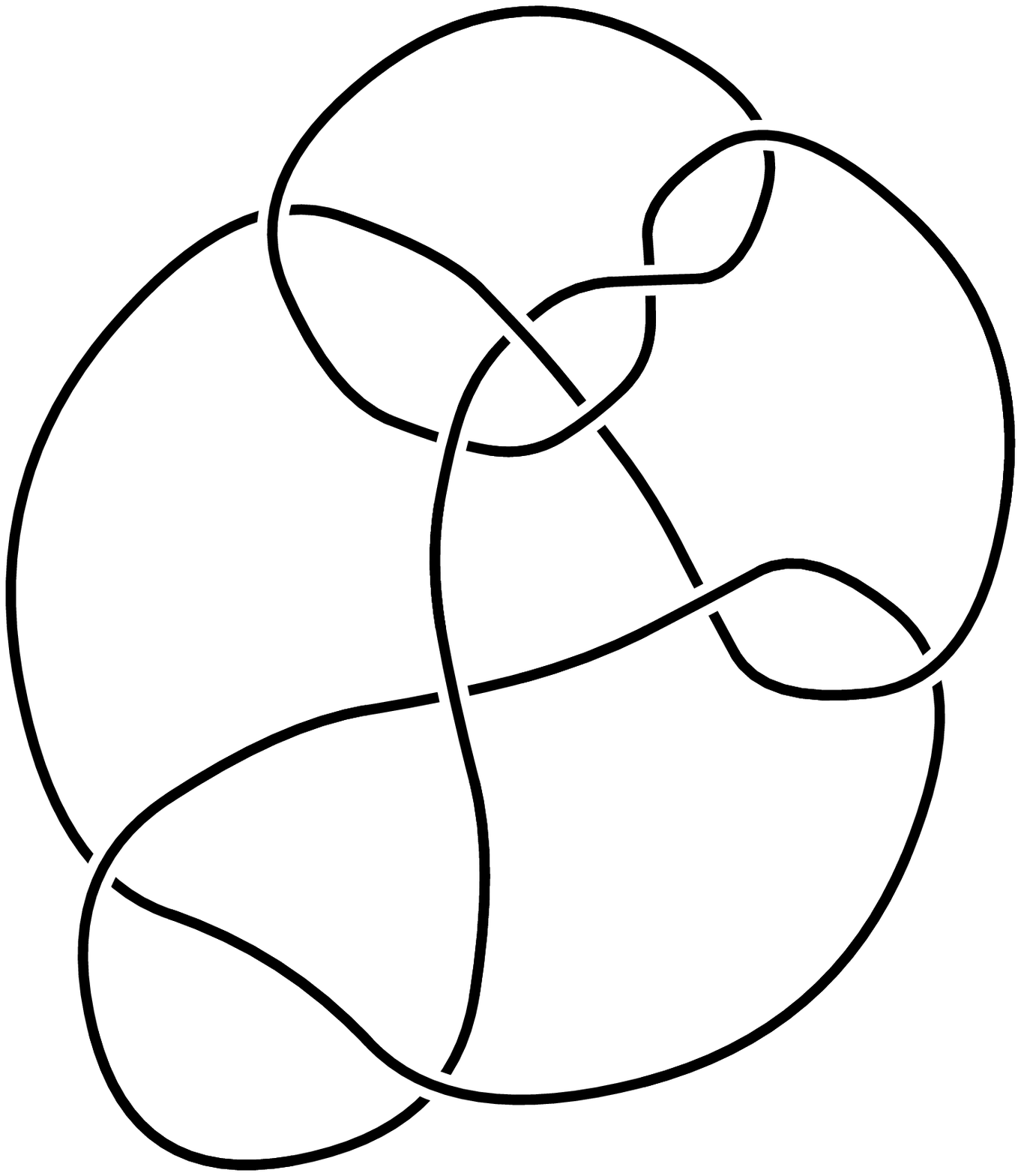}
    &
    \includegraphics[width=75pt]{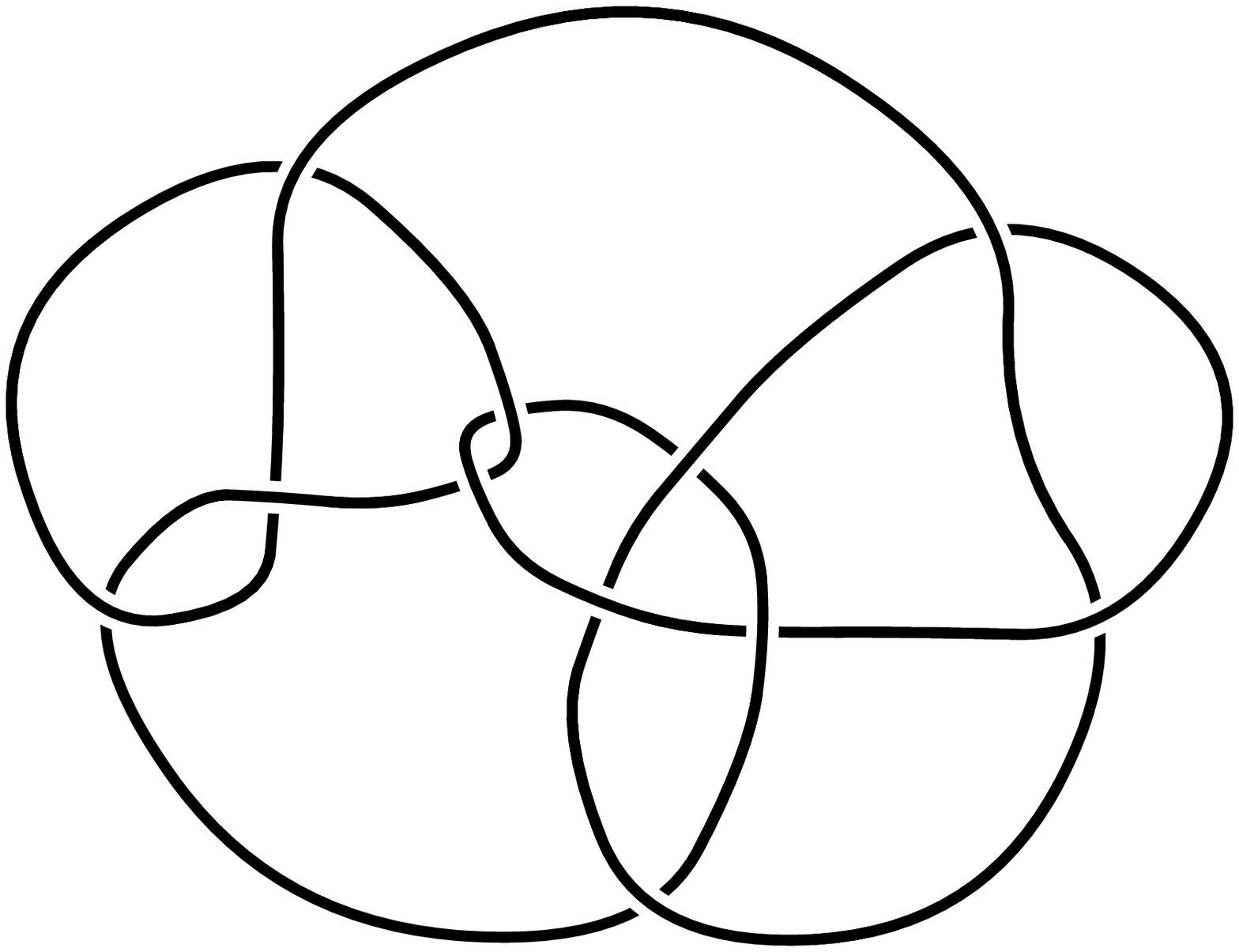}
    \\[-10pt]
    $11^N_{34}$ & $11^N_{42}$ & $11^N_{39}$ & $11^N_{45}$
    \\[10pt]
    \hline
    &&&\\[-10pt]
    \includegraphics[width=75pt]{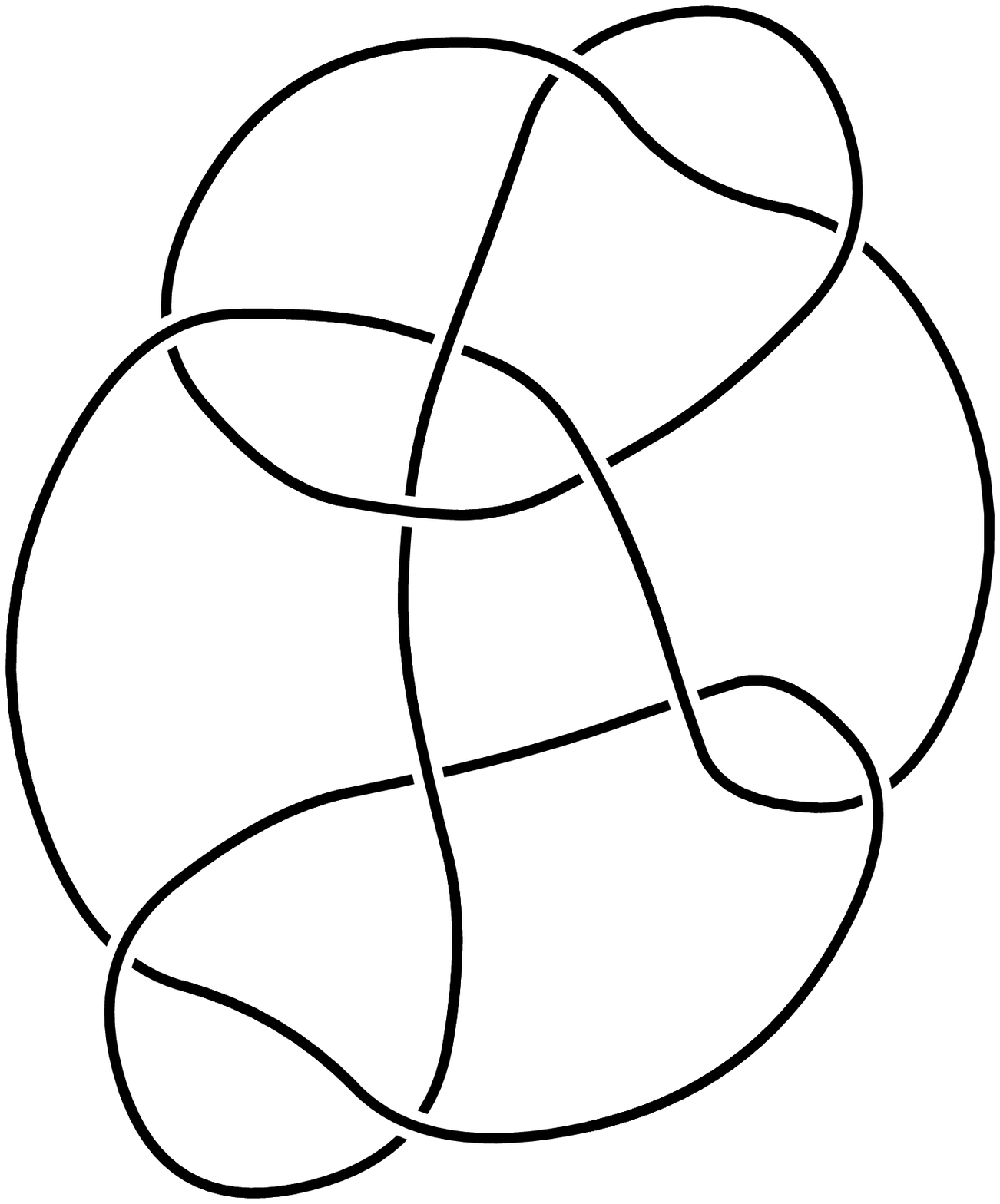}
    &
    \includegraphics[width=75pt]{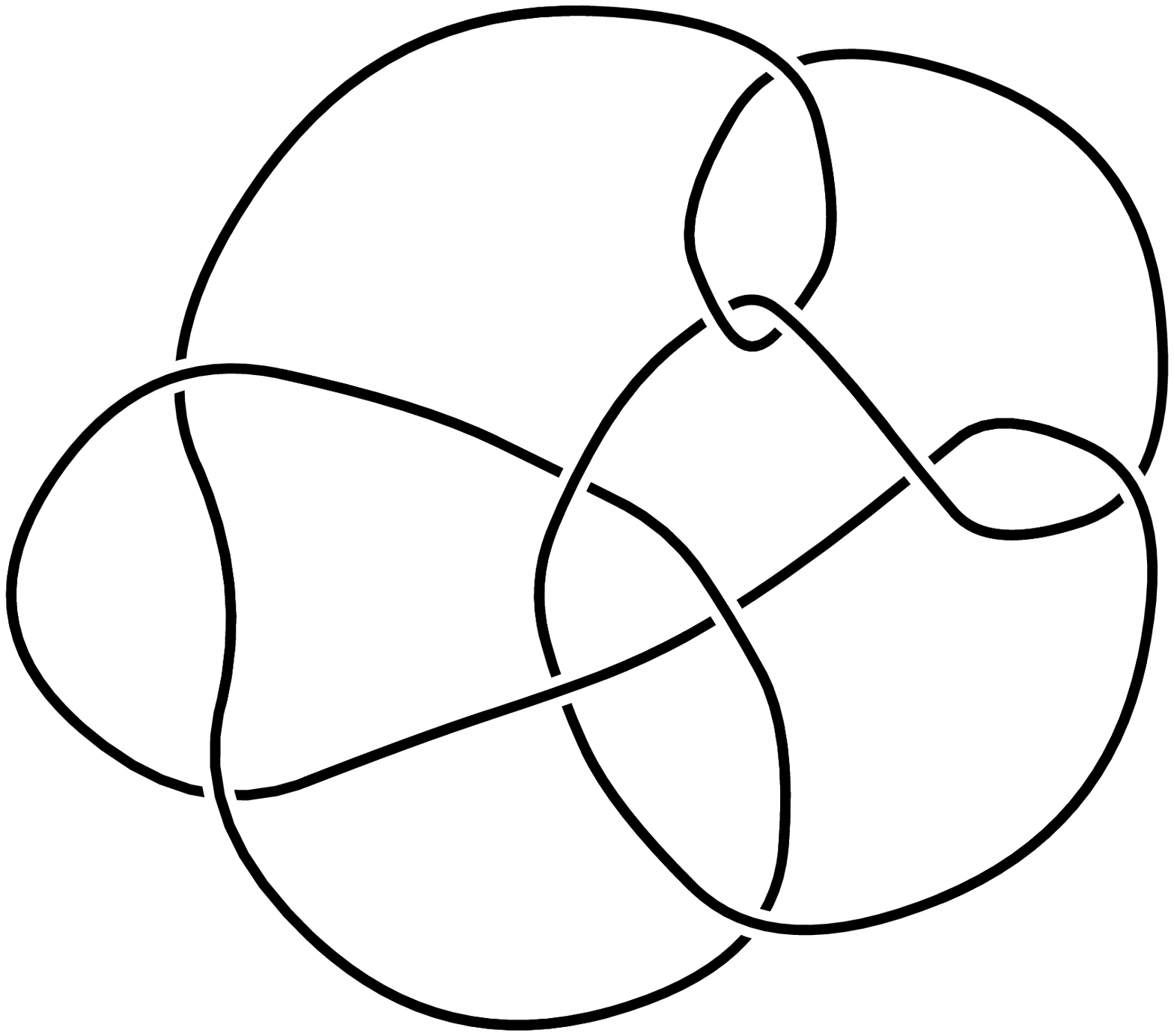}
    &
    \includegraphics[width=75pt]{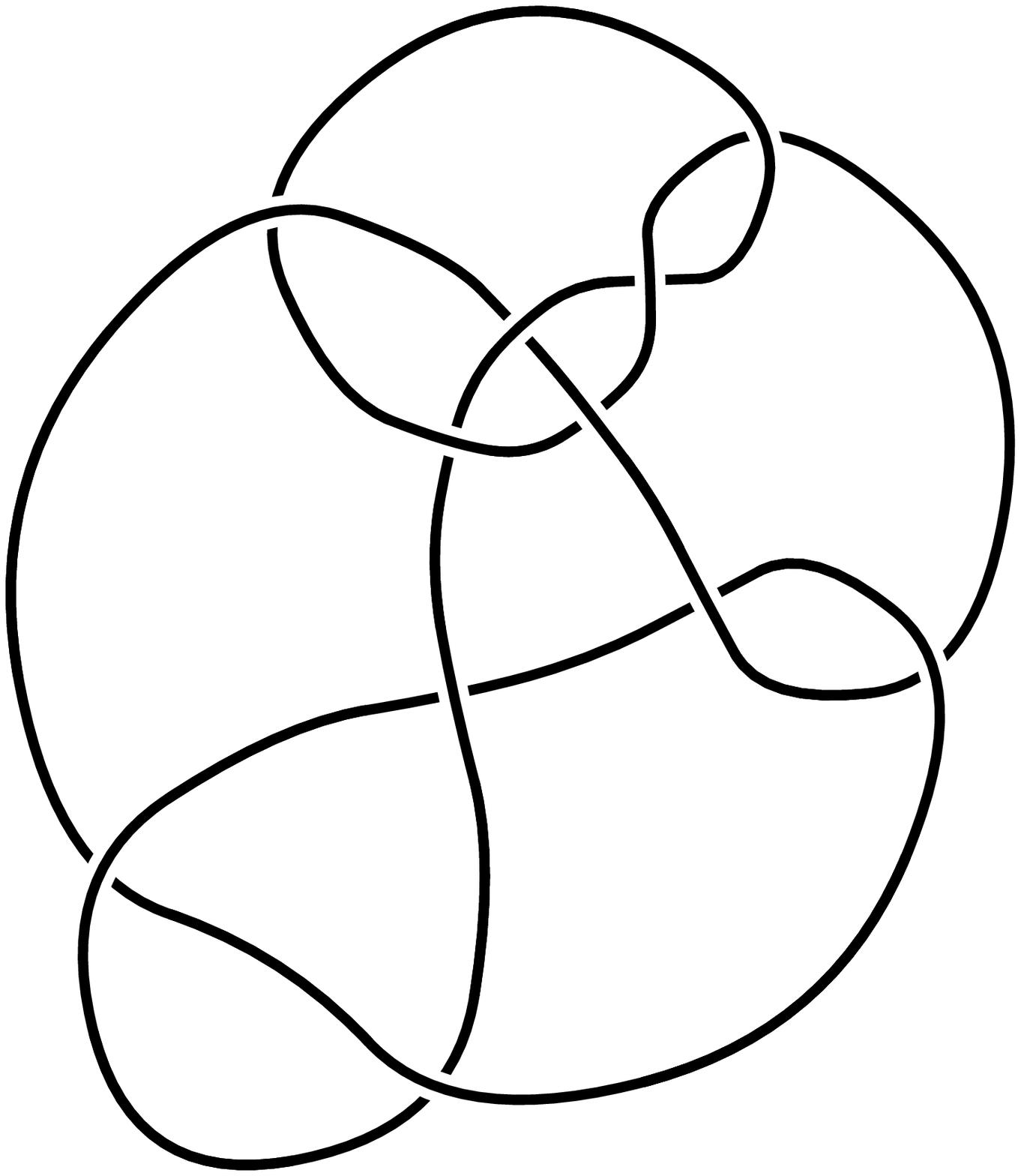}
    &
    \includegraphics[width=75pt]{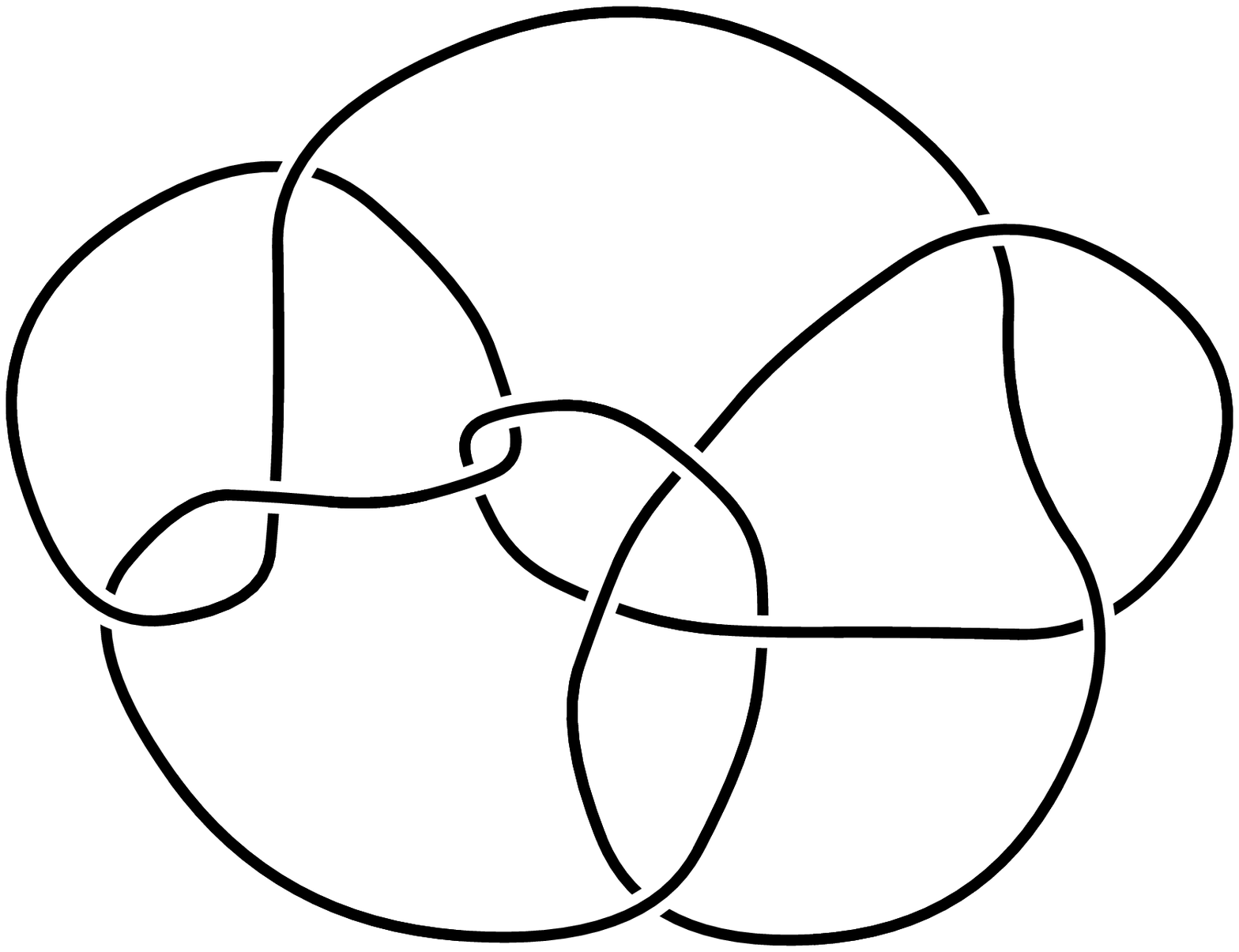}
    \\[-10pt]
    $11^N_{35}$ & $11^N_{43}$ & $11^N_{40}$ & $11^N_{46}$
    \\[10pt]
    \hline
    &&&\\[-10pt]
    \includegraphics[width=75pt]{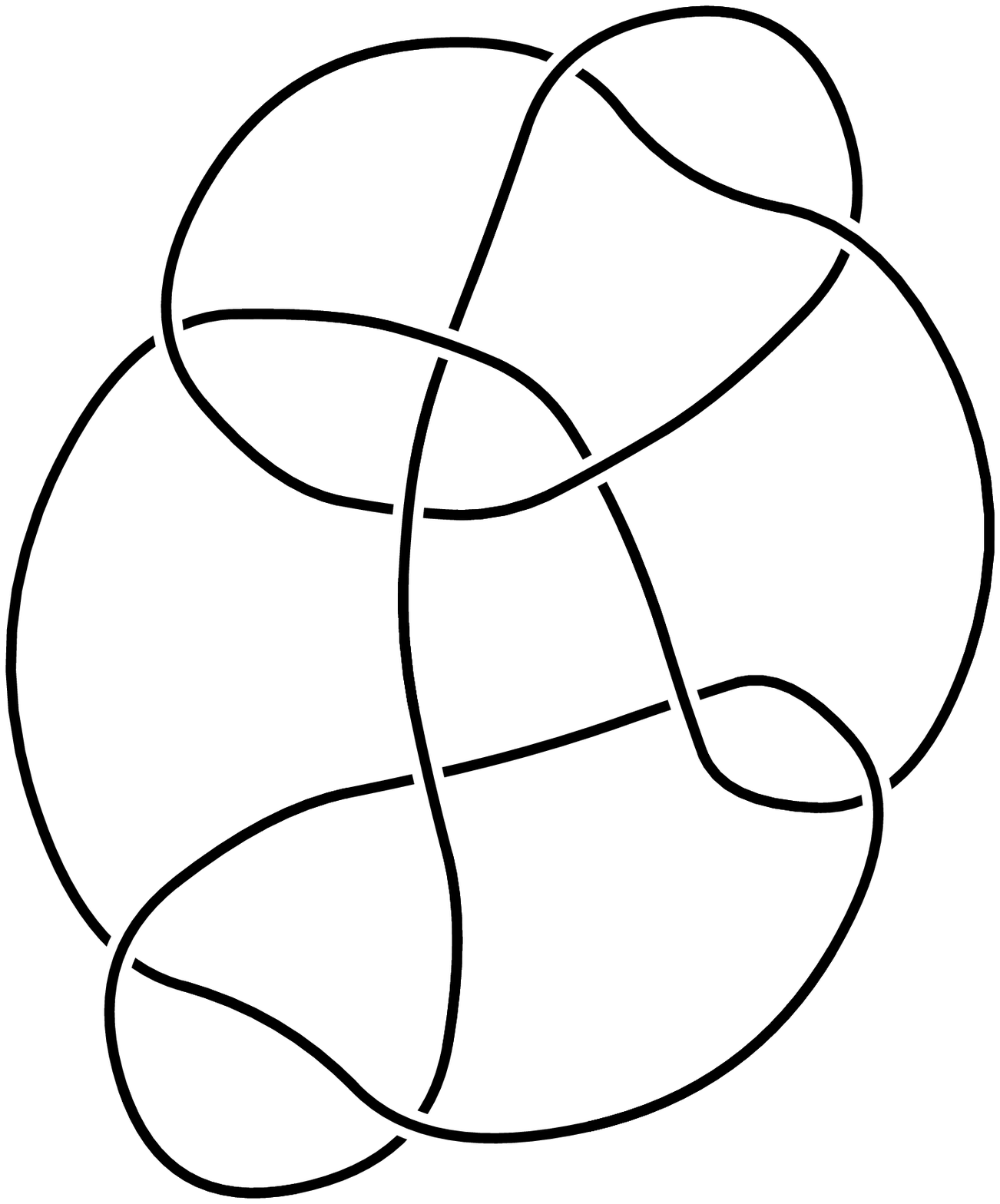}
    &
    \includegraphics[width=75pt]{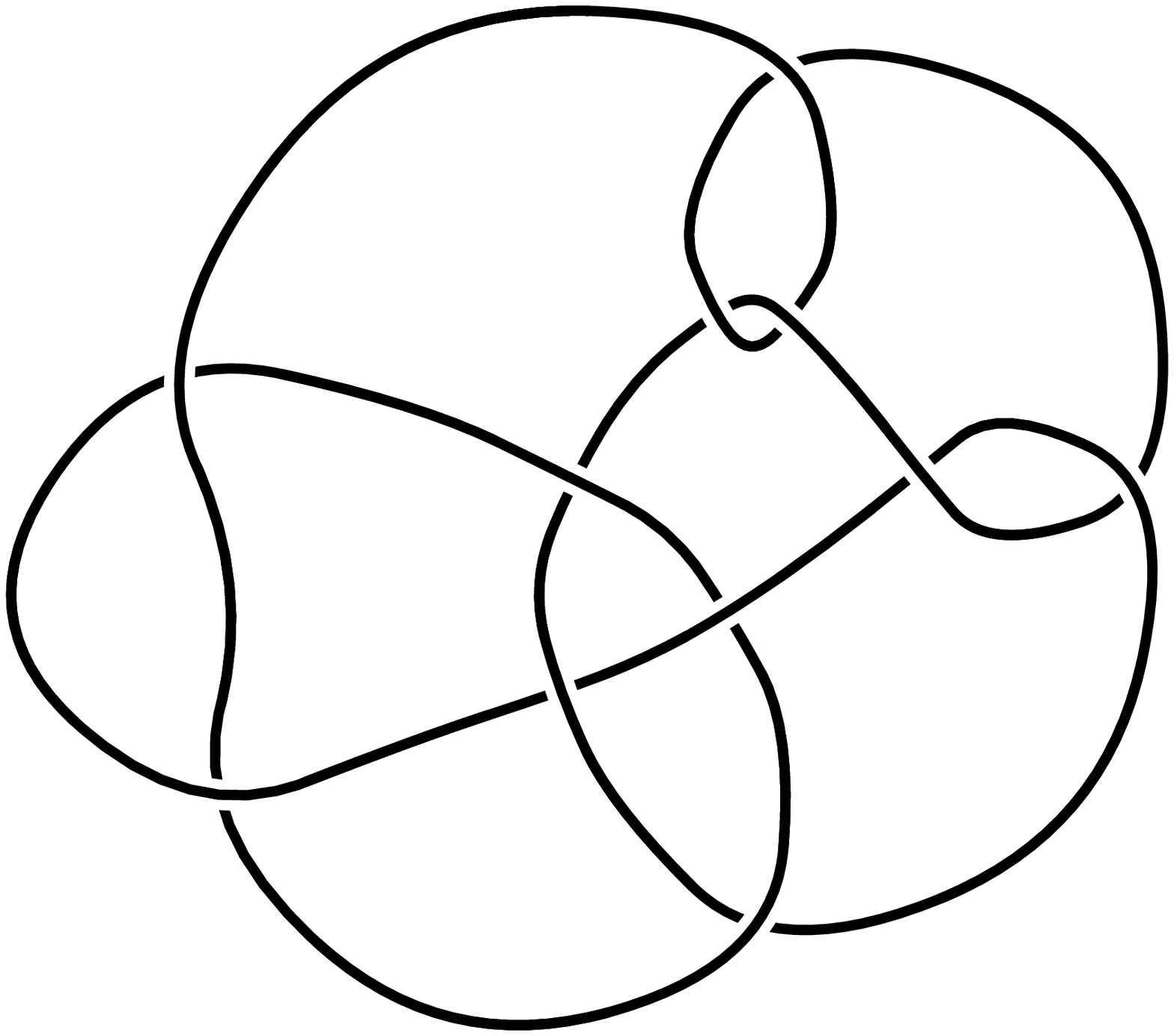}
    &
    \includegraphics[width=75pt]{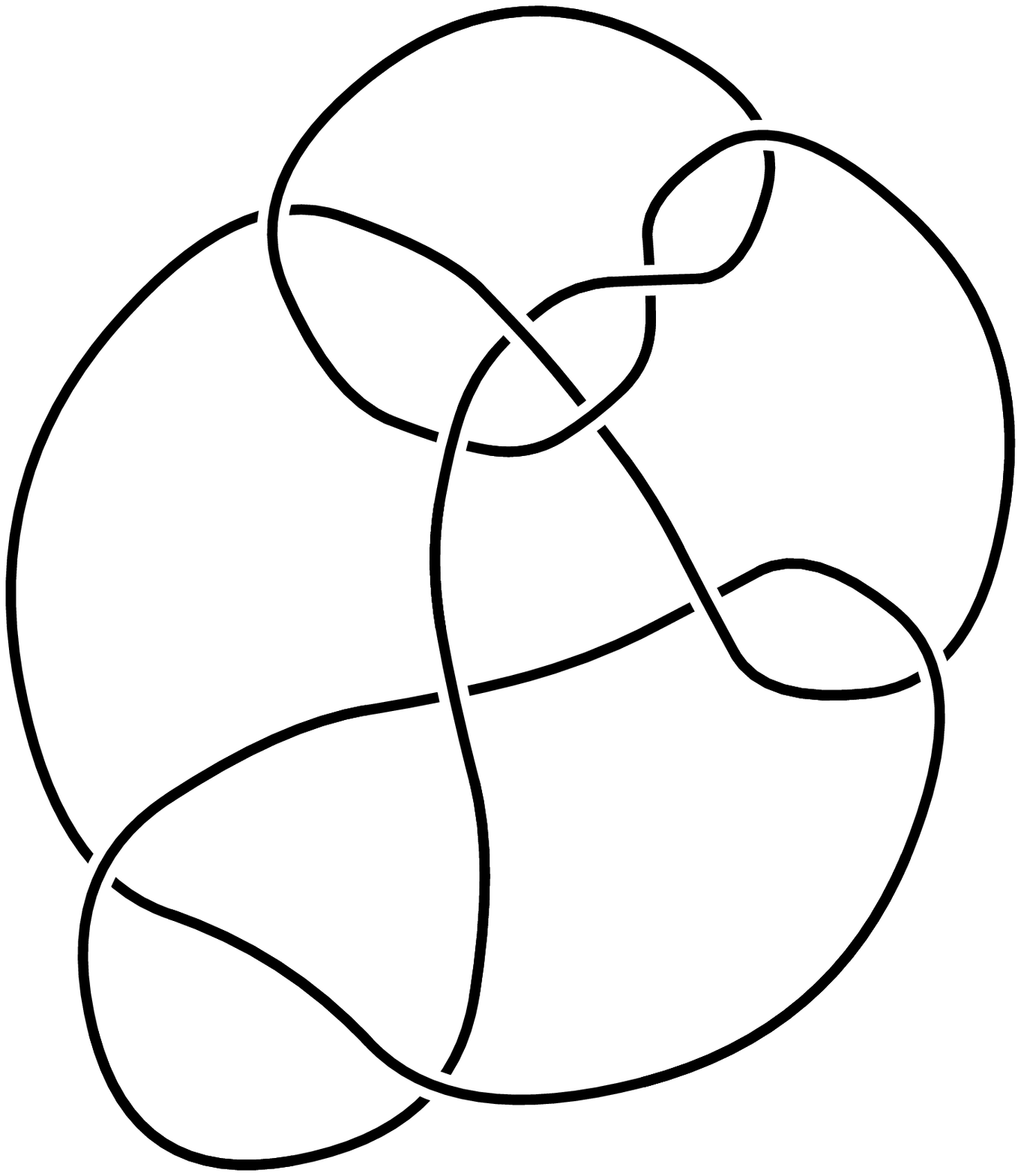}
    &
    \includegraphics[width=75pt]{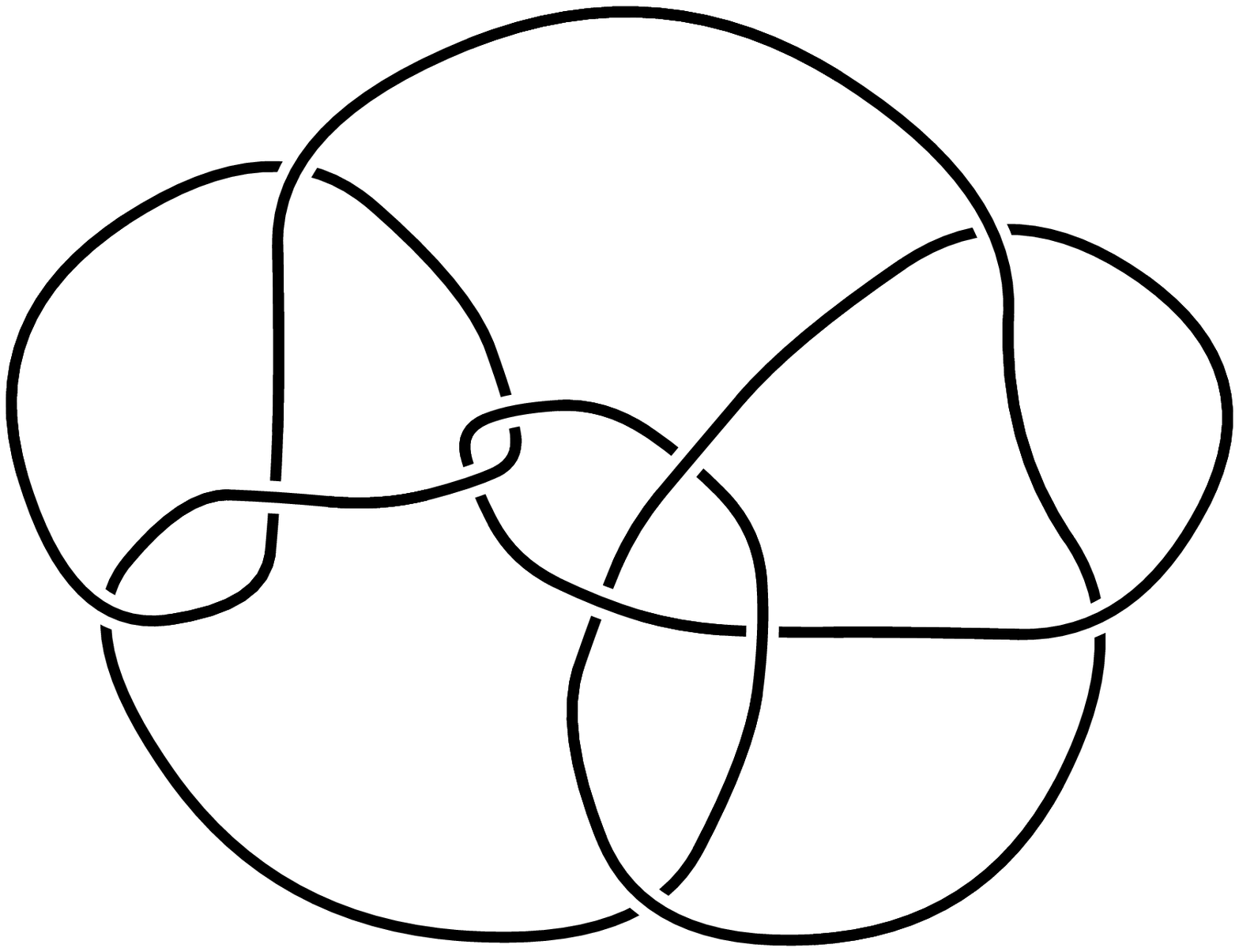}
    \\[-10pt]
    $11^N_{36}$ & $11^N_{44}$ & $11^N_{41}$ & $11^N_{47}$
    \\[10pt]
    \hline
    &&&\\[-10pt]
    \includegraphics[width=75pt]{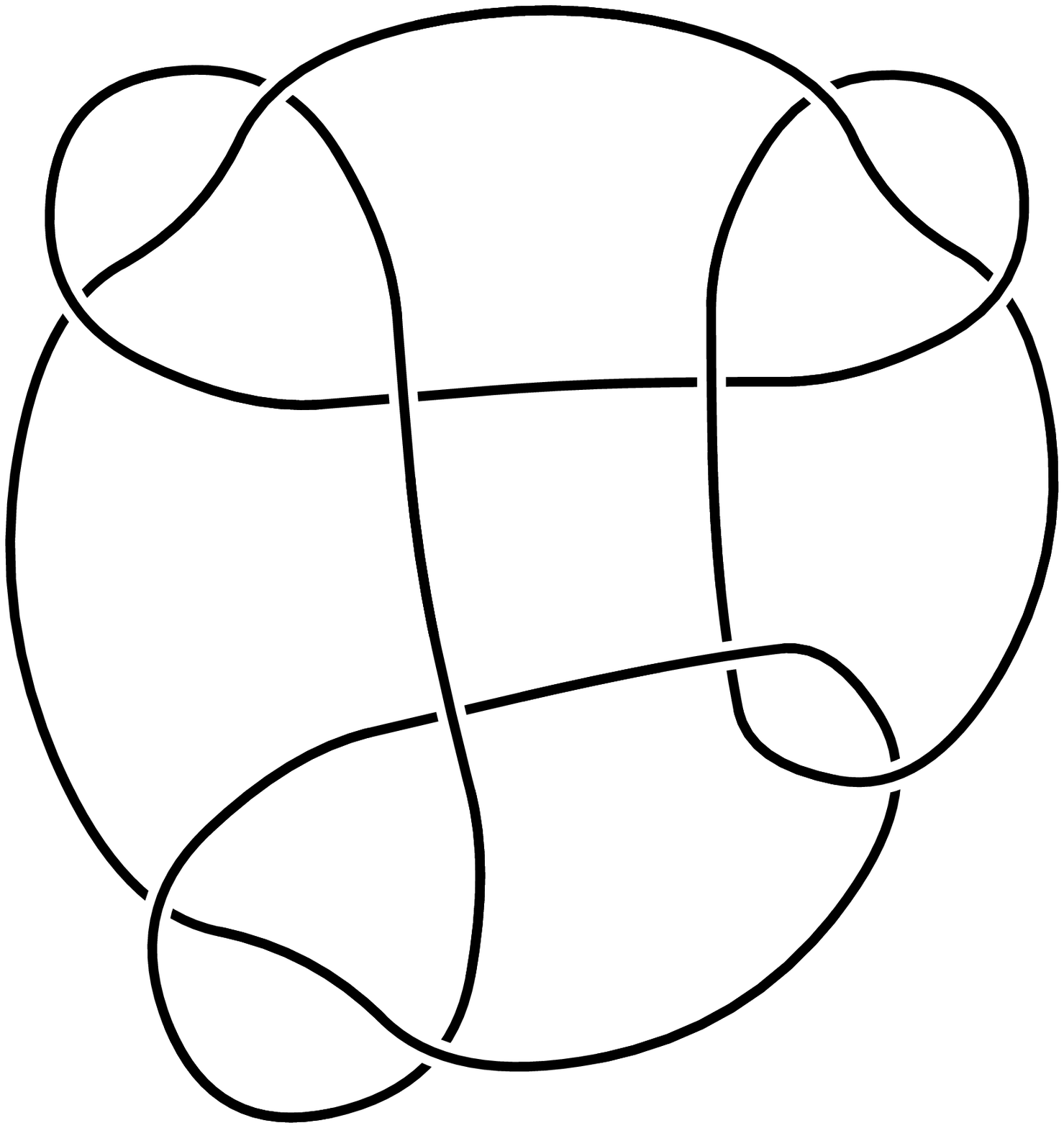}
    &
    \includegraphics[width=75pt]{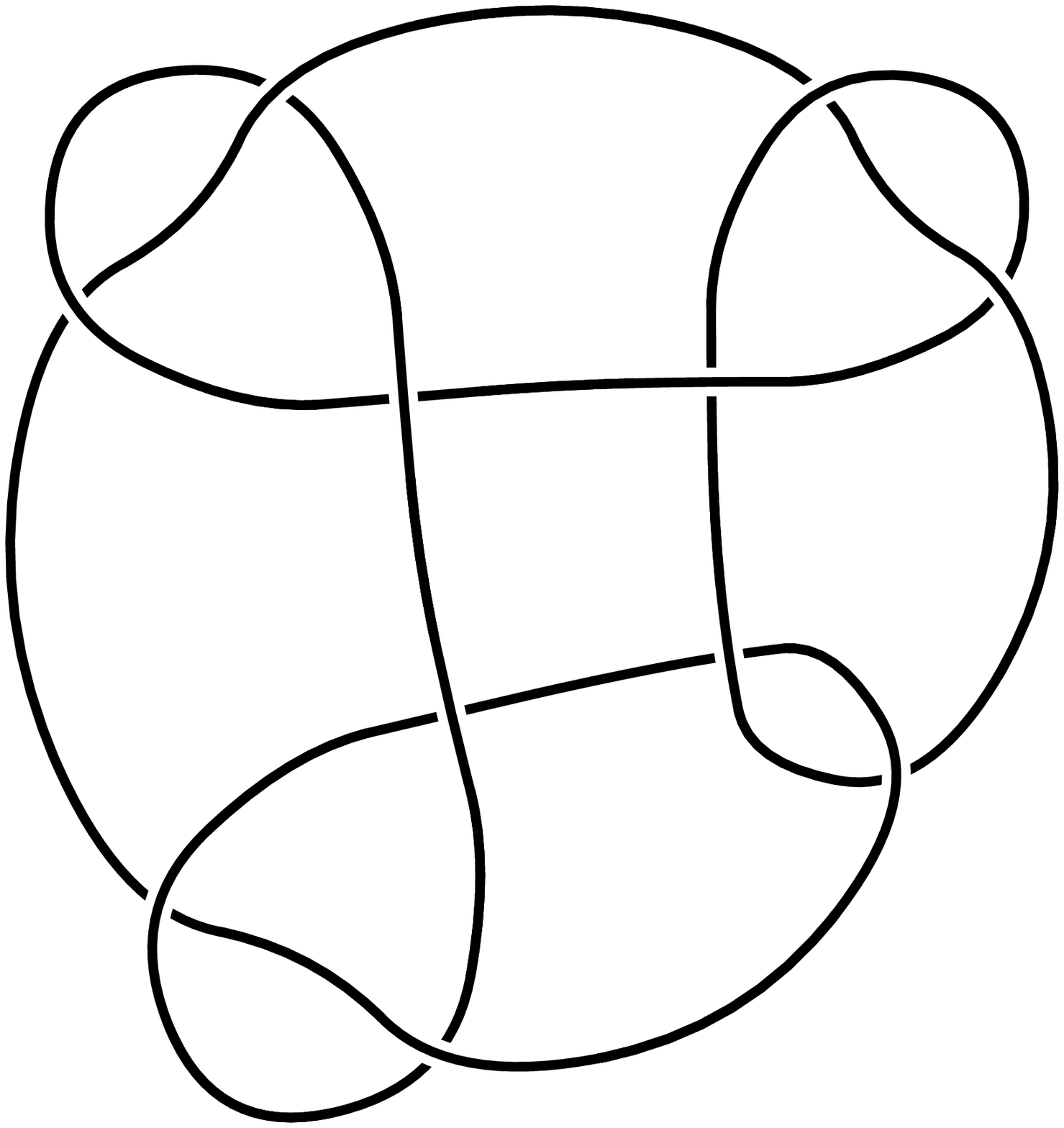}
    &
    \includegraphics[width=75pt]{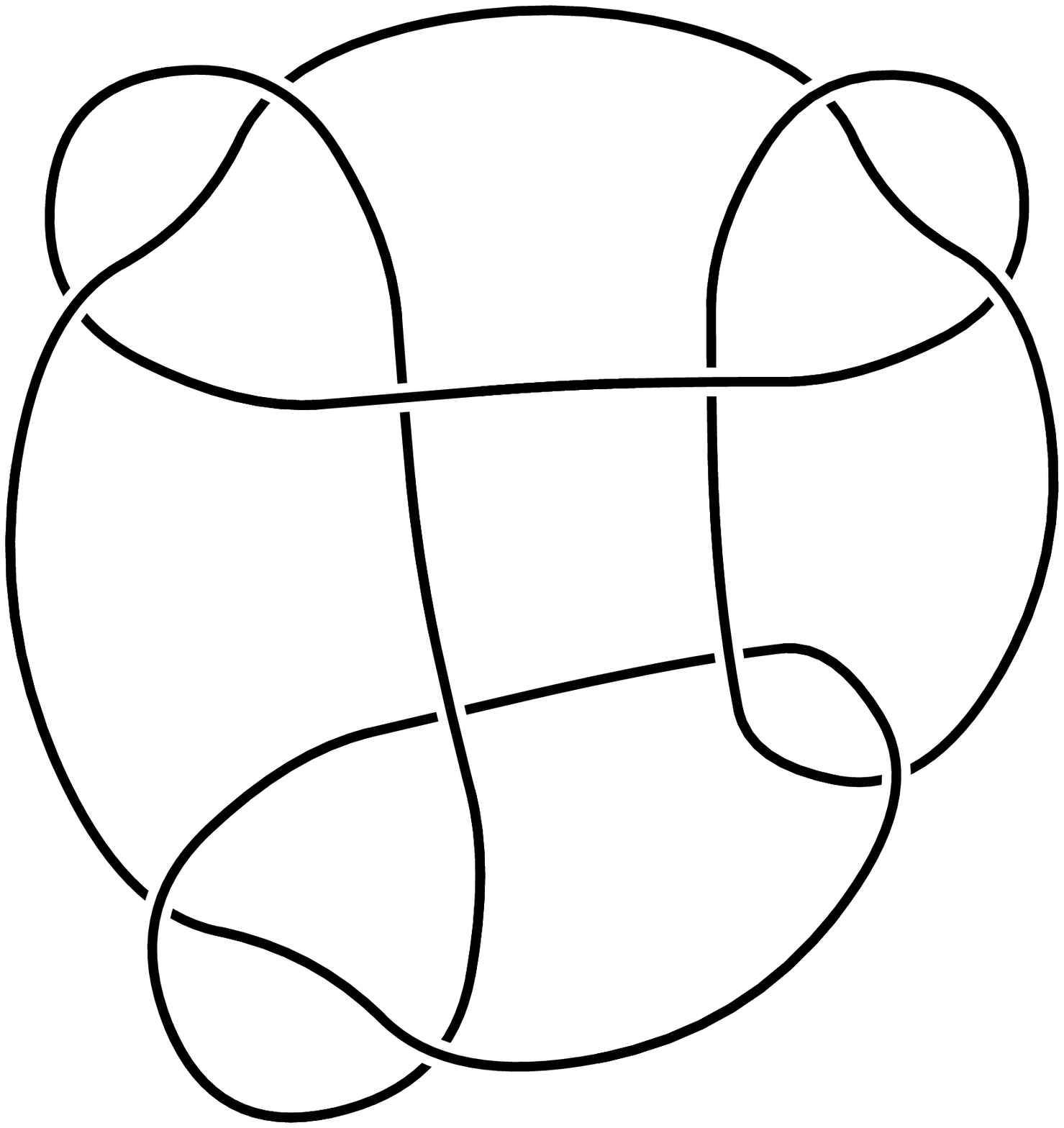}
    &
    \includegraphics[width=75pt]{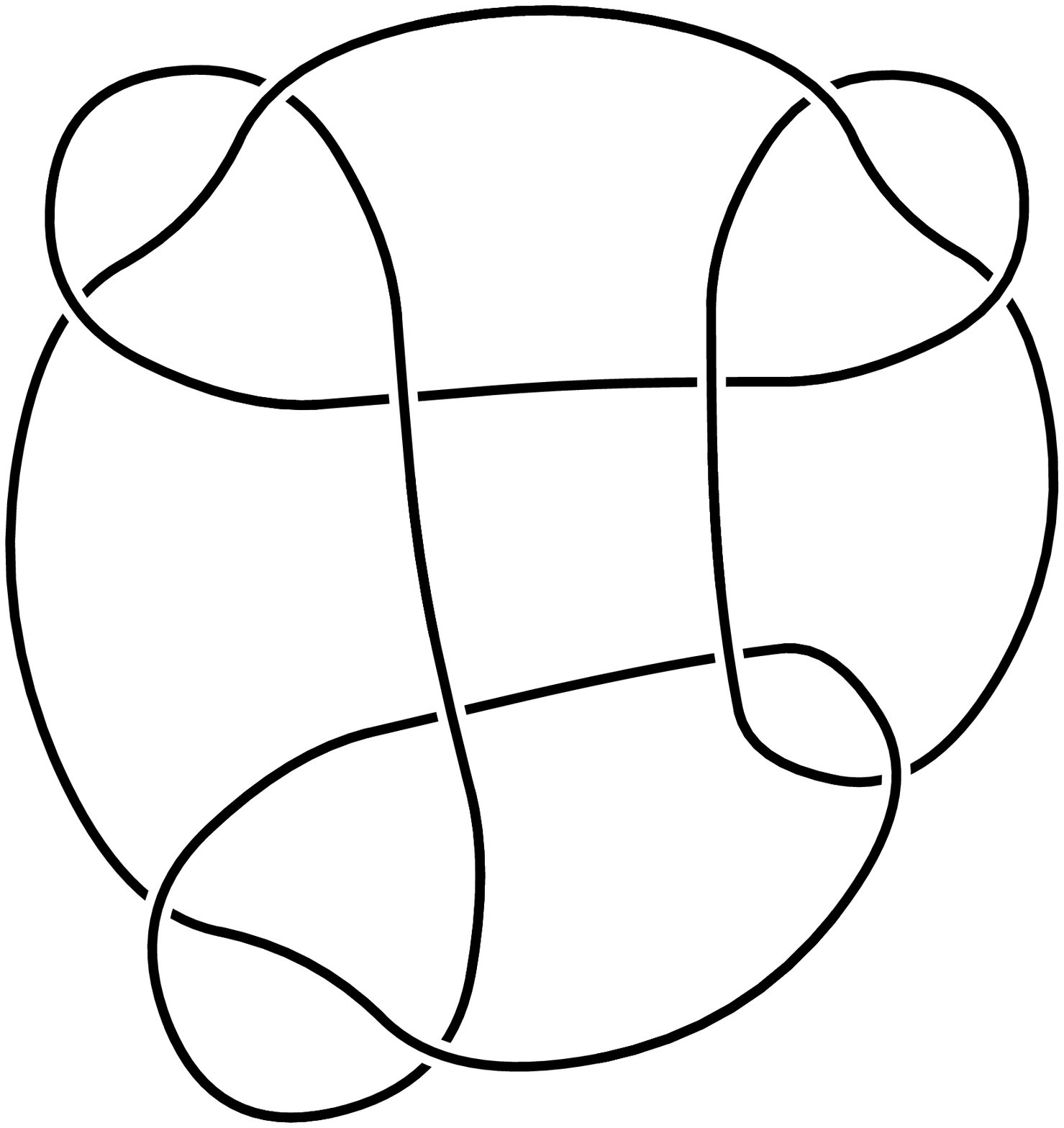}
    \\[-10pt]
    $11^N_{71}$ & $11^N_{75}$ & $11^N_{73}$ & $11^N_{74}$
    \\[10pt]
    \hline
    &&&\\[-10pt]
    \includegraphics[width=75pt]{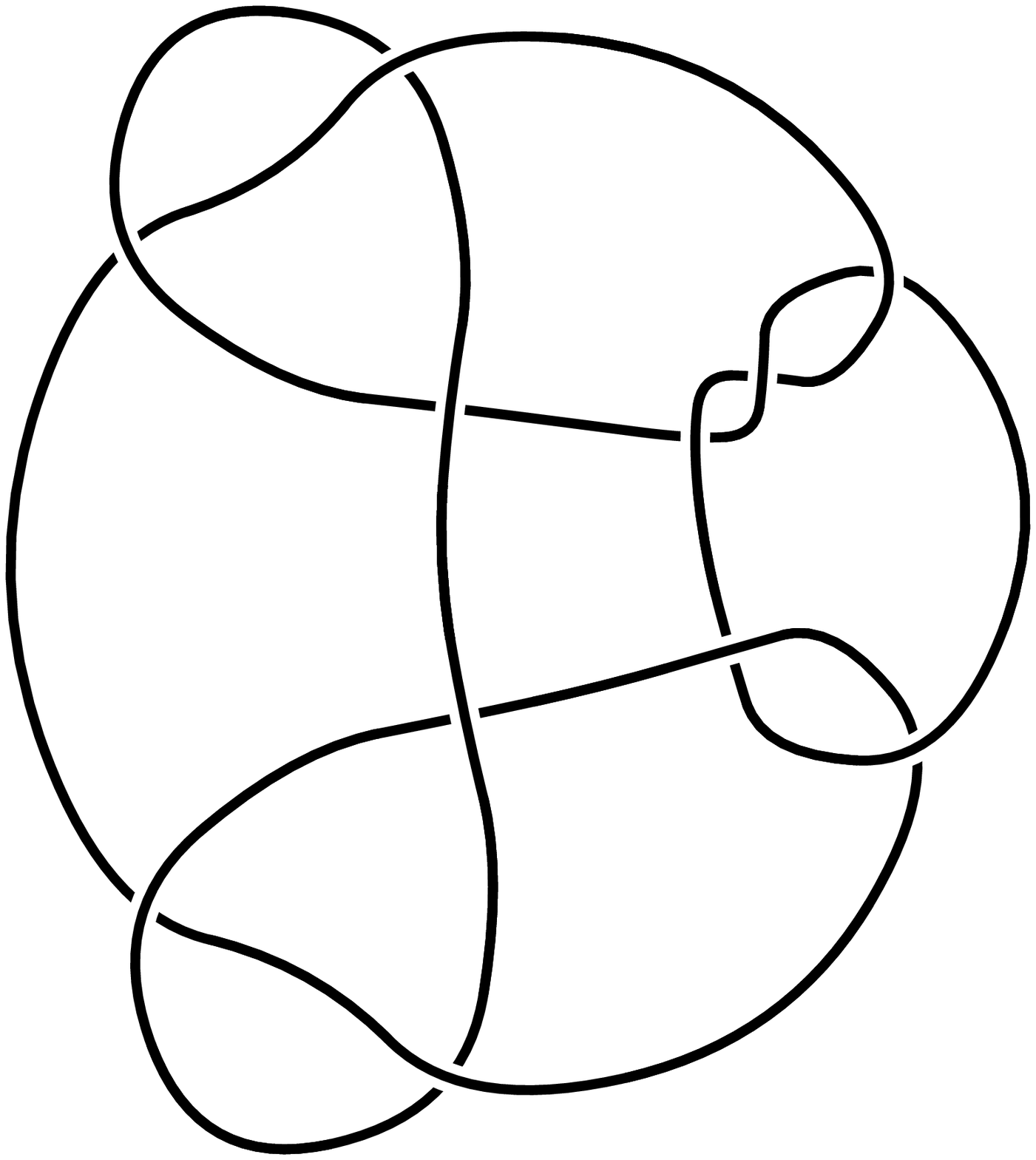}
    &
    \includegraphics[width=75pt]{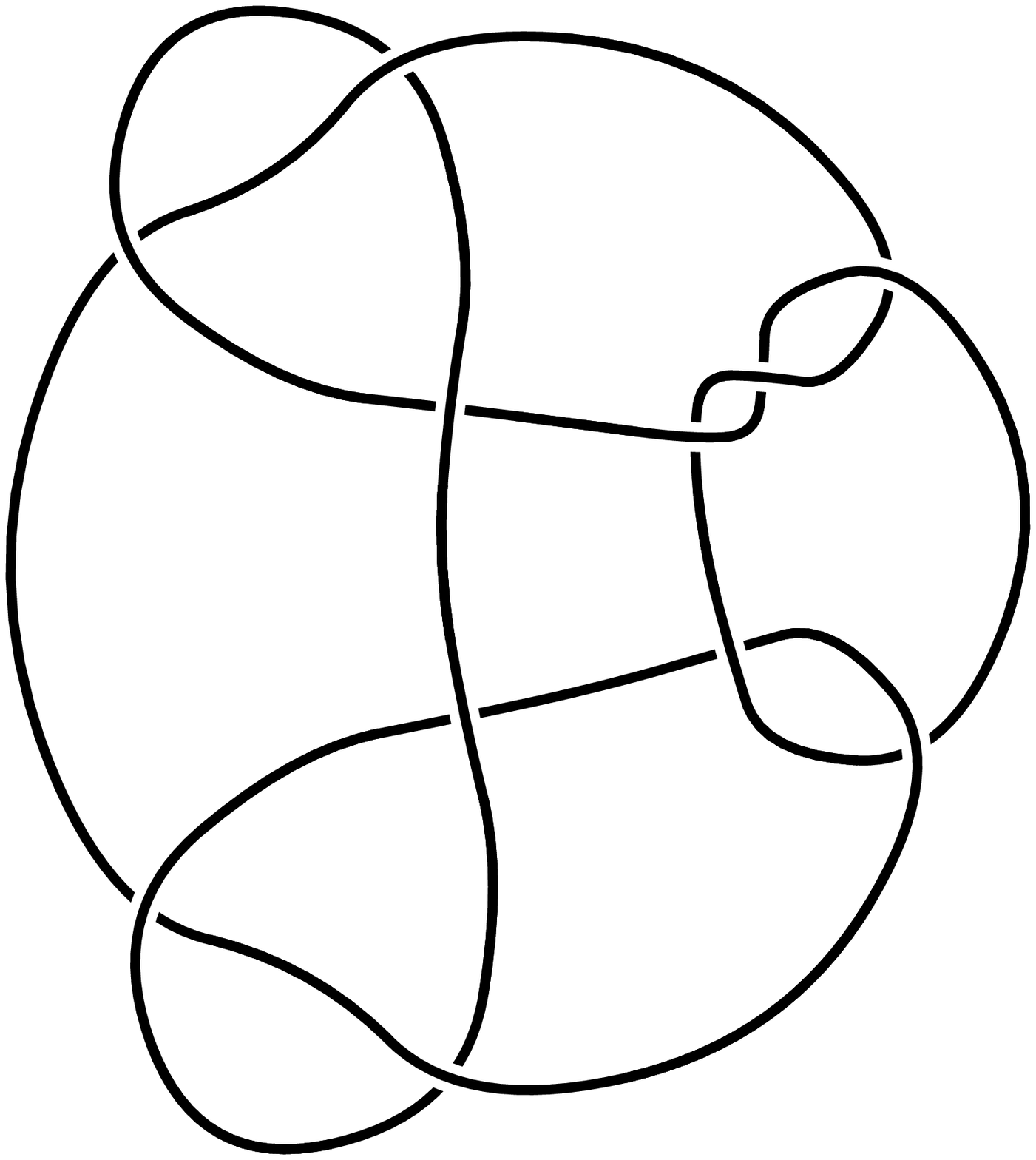}
    &
    \includegraphics[width=75pt]{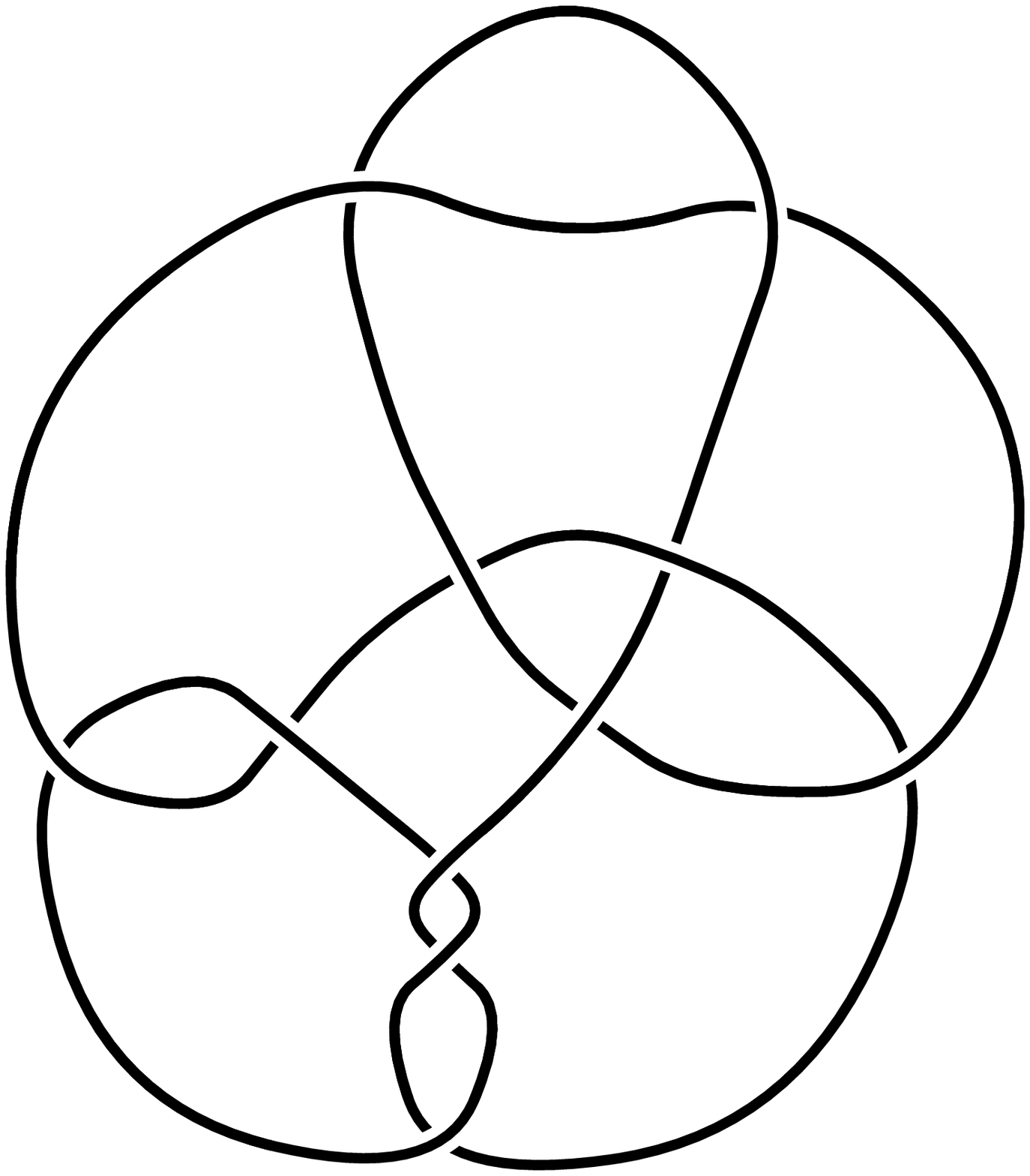}
    &
    \includegraphics[width=75pt,angle=90]{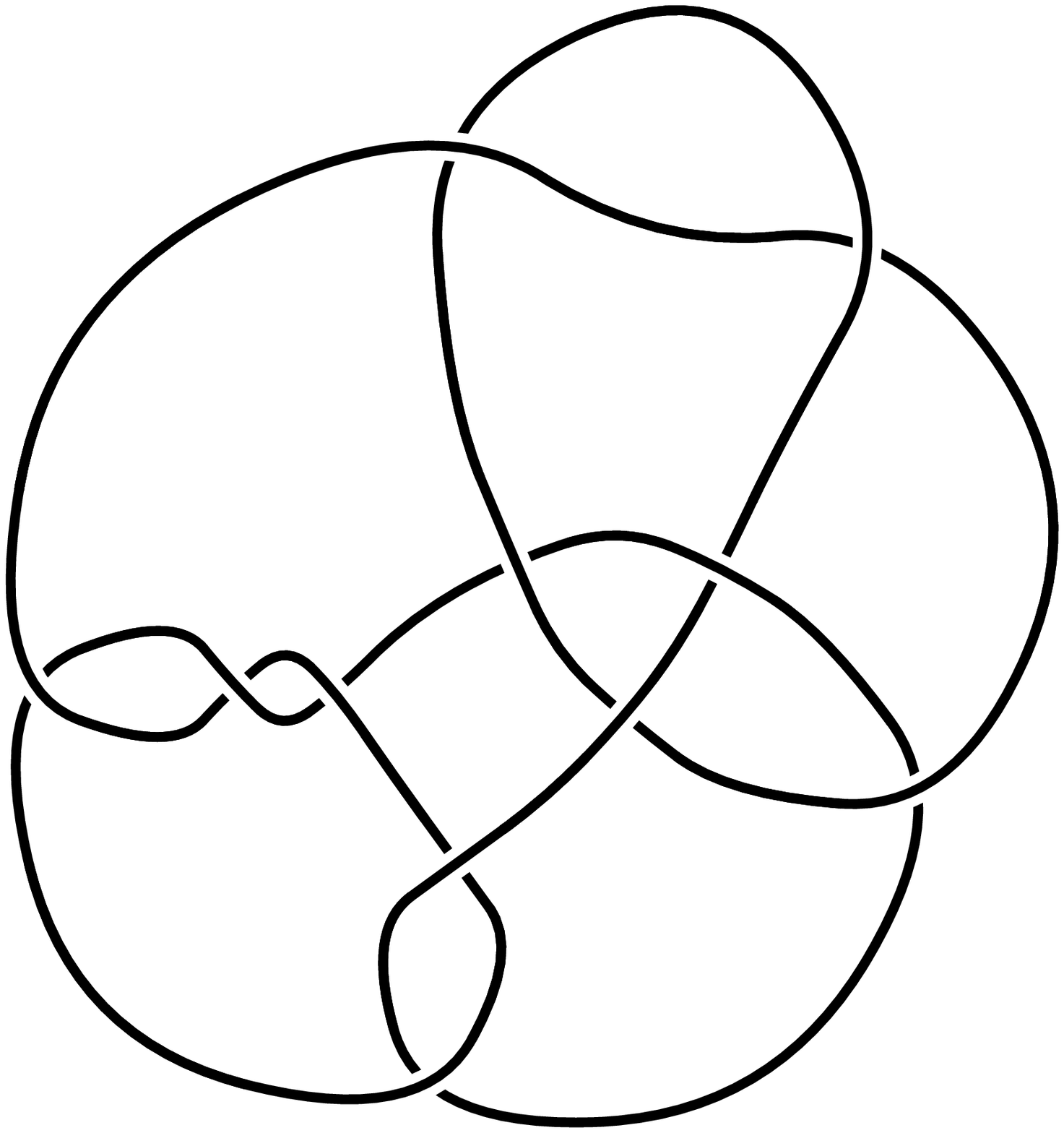}
    \\[-10pt]
    $11^N_{76}$ & $11^N_{78}$ & $11^N_{151}$ & $11^N_{152}$
  \end{tabular}
  \caption{Mutant cliques within $\mathbb{K}^N_{11}$.}
  \label{figure:Nonalternating11crossingmutantcliques}
  \end{centering}
\end{figure}

In this manner, within $\mathbb{K}_{11}$ we determine $16$ mutant candidate
pairs, all of which actually are mutants. These are depicted in
Figures~\ref{figure:Alternating11crossingmutantcliques}
and~\ref{figure:Nonalternating11crossingmutantcliques}. In many of the
pictures, a mutation can be seen by immediate inspection; others require a
little more work. We mention that we describe the mutant cliques of knot
types, ignoring reflections, so the identification of a mutation relating
pairs of diagrams may include imposing a reflection on one of them. The
inspection process is in some cases made easier by the fact that
nonalternating mutant cliques sometimes share the same shadows
(crossing-information-oblitterated diagrams) as alternating mutant cliques ---
compare $(11^A_{24},11^A_{26})$ with $(11^N_{39},11^N_{45})$,
$(11^N_{40},11^N_{46})$ and $(11^N_{41},11^N_{47})$ --- and this in turn
reflects the way in which the HTW tables are compiled.

\tiny
\begin{tabular}{c@{\hspace{0pt}}c}
  \\
  {\normalsize $\begin{array}{c}11^A\end{array}$}
  &
  $
  \begin{array}{p{275pt}}
    (19,25), (24,26), (44,47), (57,231), (251,253), (252,254)
  \end{array}
  $
  \\
  {\normalsize $\begin{array}{c}11^N\end{array}$}
  &
  $
  \begin{array}{p{275pt}}
  (34,42), (35,43), (36,44), (39,45), (40,46), (41,47), (71,75), (73,74),
  (76,78), (151,152)
  \end{array}
  $
  \\
  \\
\end{tabular}
\normalsize

Similarly, within $\mathbb{K}_{12}$ we determine via the invariants the following $75$ candidate
mutant cliques, all of which turn out to be genuine mutant cliques. These are depicted in the
Appendix as
Figures~\ref{figure:Alternating12crossingmutantcliques1of4}--\ref{figure:Nonalternating12crossingmutantcliques6of6}.
We mention that our list of mutant cliques is in agreement with that of
Stoimenow~\cite{Stoimenow:Knotdatatables}.

\tiny
\begin{tabular}{c@{\hspace{0pt}}c}
  \\
  {\normalsize $\begin{array}{c}12^A\end{array}$}
  &
  $
  \begin{array}{p{275pt}}
    (7,14), (13,15), (29,113), (36,694), (44,64), (45,65), (48,60), (59,63),
    (67,136), (91,111), (101,115), (102,107), (108,120), (114,117), (126,132),
    (131,133), (134,188), (154,162), (164,166), (167,692), (195,693), (639,680),
    (675,688), (811,817), (829,832), (830,831), (844,846), (30,33,157),
    (116,122,182)
  \end{array}
  $
  \\
  \\
  {\normalsize $\begin{array}{c}12^N\end{array}$}
  &
  $
  \begin{array}{p{275pt}}
    (21,29), (22,30), (23,31), (26,32), (27,33), (28,34), (55,223), (58,222),
    (59,220), (63,225), (64,261), (67,229), (85,130), (86,131), (87,132), (88,133),
    (89,134), (90,135), (91,136), (92,137), (93,138), (98,125), (99,126),
    (122,127), (123,128), (124,129), (205,226), (206,227), (207,228), (208,212),
    (209,213), (210,214), (231,232), (252,262), (255,263), (256,264), (364,365),
    (421,422), (553,556), (670,681), (671,682), (691,692), (693,696), (56,57,221),
    (60,61,219), (62,66,224)
  \end{array}
  $
  \\
  \\
\end{tabular}
\normalsize

\clearpage

\begin{figure}[ht]
  \begin{centering}
  \includegraphics[width=120pt]{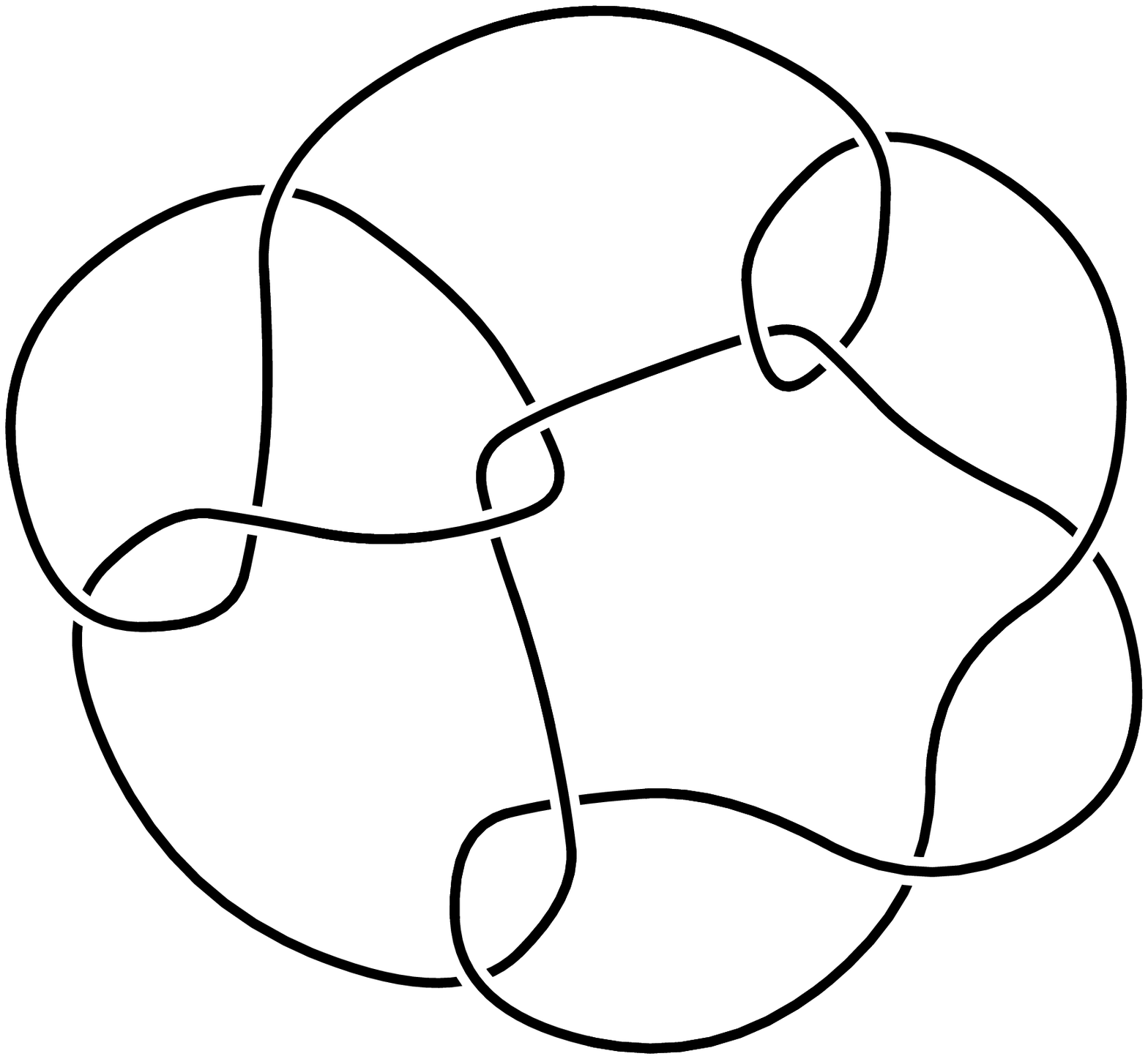}
  \hspace{20pt}
  \includegraphics[width=120pt]{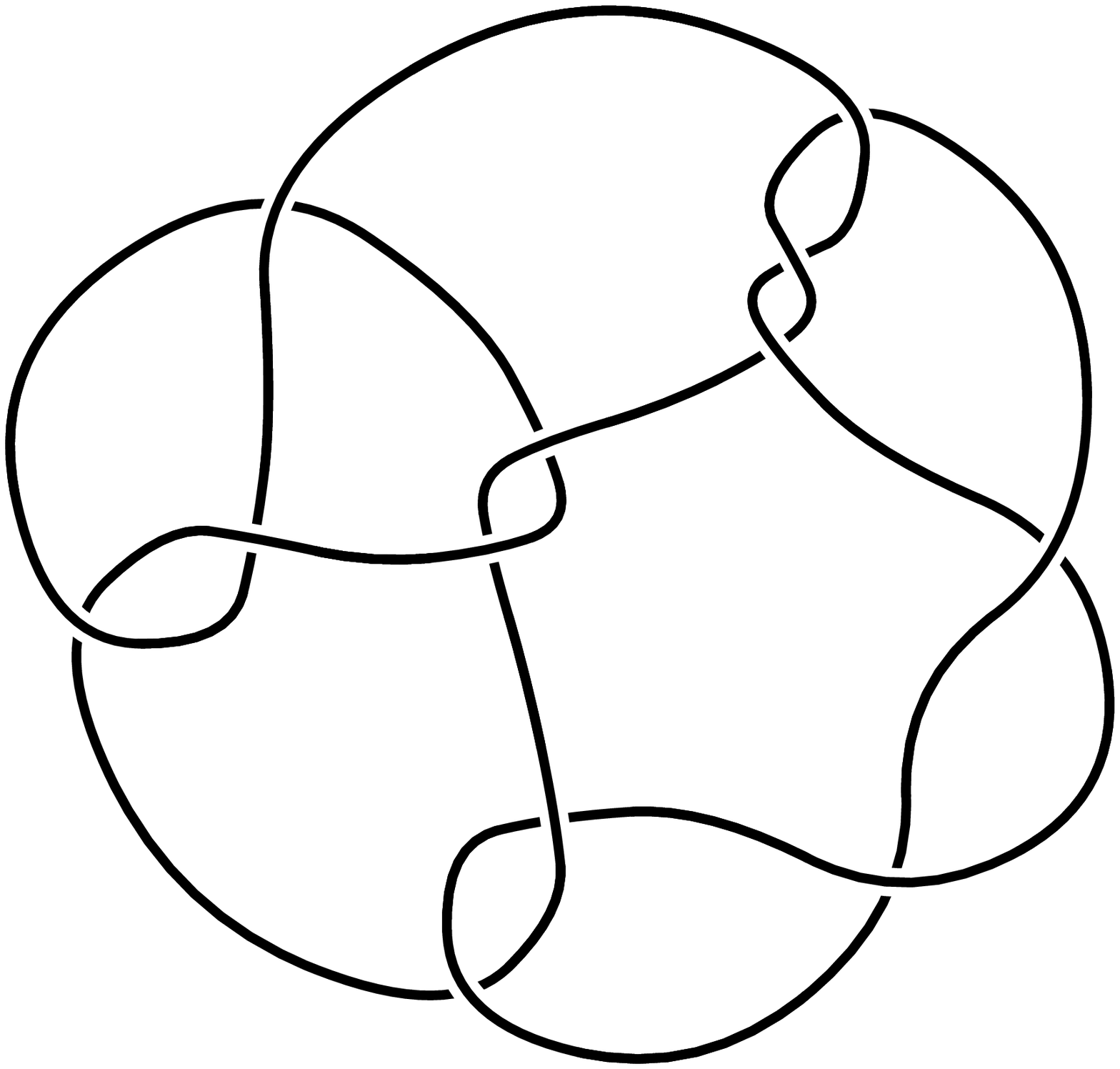}
  \caption{%
    $12^N_{119}$ and $12^N_{120}$ are \emph{not} mutants, but have
    the same hyperbolic volume.
  }
  \label{figure:12N119and12N120}
  \end{centering}
\end{figure}

Of interest is that a single pair of knots $(12^N_{119},12^N_{120})$ are
distinguished by all three polynomials, but not by the (apparent) hyperbolic
volume (to $10$ places, $9.8759424404$). This pair is depicted in
Figure~\ref{figure:12N119and12N120}; they are related by the exchange of a $2$
tangle with a $1/2$ tangle. For knots of greater crossing numbers, the same
phenomenon appears more frequently.

The maximum size of mutant candidate cliques, and presumably also of true
mutant cliques, rises with crossing number. Amongst the knots of $11$ to $14$
crossings, the maximum candidate clique sizes are respectively $2,3,6$ and $8$.

With respect to $\mathbb{K}_{11}$ and $\mathbb{K}_{12}$, we also make the
following observations.

\begin{itemize}
\item
  All elements of each (candidate) mutant clique are of the same braid index: $4$
  for all those in $\mathbb{K}_{11}$, and $4$, $5$ or $6$ for those in
  $\mathbb{K}_{12}$. (There is no reason to suppose that mutation preserves
  braid index more generally.)

\item
  None of the mutant cliques and mutant clique candidates are $2$-bridge knots.
  This reflects the fact that $2$-bridge links are known to be invariant under
  mutation~\cite[p48]{Kawauchi:1996}.
\end{itemize}


\section{Where $LG$ fails to distinguish prime knots}
\label{section:WhereLGfirstfailstodistinguishprimeknots}

We now wish to use these new evaluations to identify the point at which $LG$
first fails to distinguish pairs of distinct prime knot types. Recall
from~\cite{DeWit:2000} that $LG$ is known to be complete for all prime knots
(including reflections) of at most $10$ crossings. Of course $LG$, in common
with many well-known quantum link invariants, fails to distinguish
mutant links.%
\footnote{%
  A quantum link invariant based on a tensor product of representations of a quantum
  (super)algebra will fail to distinguish mutant links if the decomposition of the
  tensor product into irreducible submodules contains no multiplicities.

  Some quantum link invariants for which the decomposition contains
  multiplicities distinguish at least \emph{some} mutants. For instance each
  invariant associated with the $SU_q(N)$ representation with Young diagram
  $[2,1]$ distinguishes the Kinoshita--Terasaka
  mutants~\cite{MortonCromwell:1996}. An example closer to $LG$ is the
  invariant associated with the $8$-dimensional $U_q[gl(2|1)]$ adjoint
  representation~\cite{GouldLinksZhang:1996b}.
} %
To illustrate this, it is demonstrated in~\cite{DeWitKauffmanLinks:1999a} that
$LG$ assigns the same polynomials to the Kinoshita--Terasaka mutant pair of
$11$-crossing prime knots. Note that there are no mutant pairs of prime knots
of less than $11$ crossings~\cite[p44]{Kawauchi:1996}.

So, we seek pairs of distinct \emph{nonmutant} prime knots undistinguished by
$LG$. Our starting point is the following fact. Recent work by Ishii and
Kanenobu~\cite{IshiiKanenobu:2004b} has determined a class of nonmutant pairs
of prime knots undistinguished by $LG$; the smallest example of these is a
pair of $14$-crossing knots $(14^A_{6955},14^A_{7393})$ described by
$16$-crossing projections.%
\footnote{%
  In their notation, these knots are respectively $K[\beta;\pm 2,\mp 2]$, where
  $\beta=\sigma_2^2 \sigma_1^{-2}$.
} %
Thus, we need only search $\mathbb{K}_{11}$--$\mathbb{K}_{13}$ to identify
where $LG$ \emph{first} fails to distinguish nonmutant prime knots.

We find that, modulo mutation, $LG$ is indeed complete over $\mathbb{K}_{11}$.
Moreover, all $11$-crossing prime knots are chiral, so we know that $LG$ is
complete, modulo mutation, for all prime knots including reflections, of up to
$11$ crossings.

\begin{figure}[ht]
  \begin{centering}
  \begin{tabular}{cc@{\hspace{10pt}}|@{\hspace{10pt}}cc}
    \includegraphics[width=75pt]{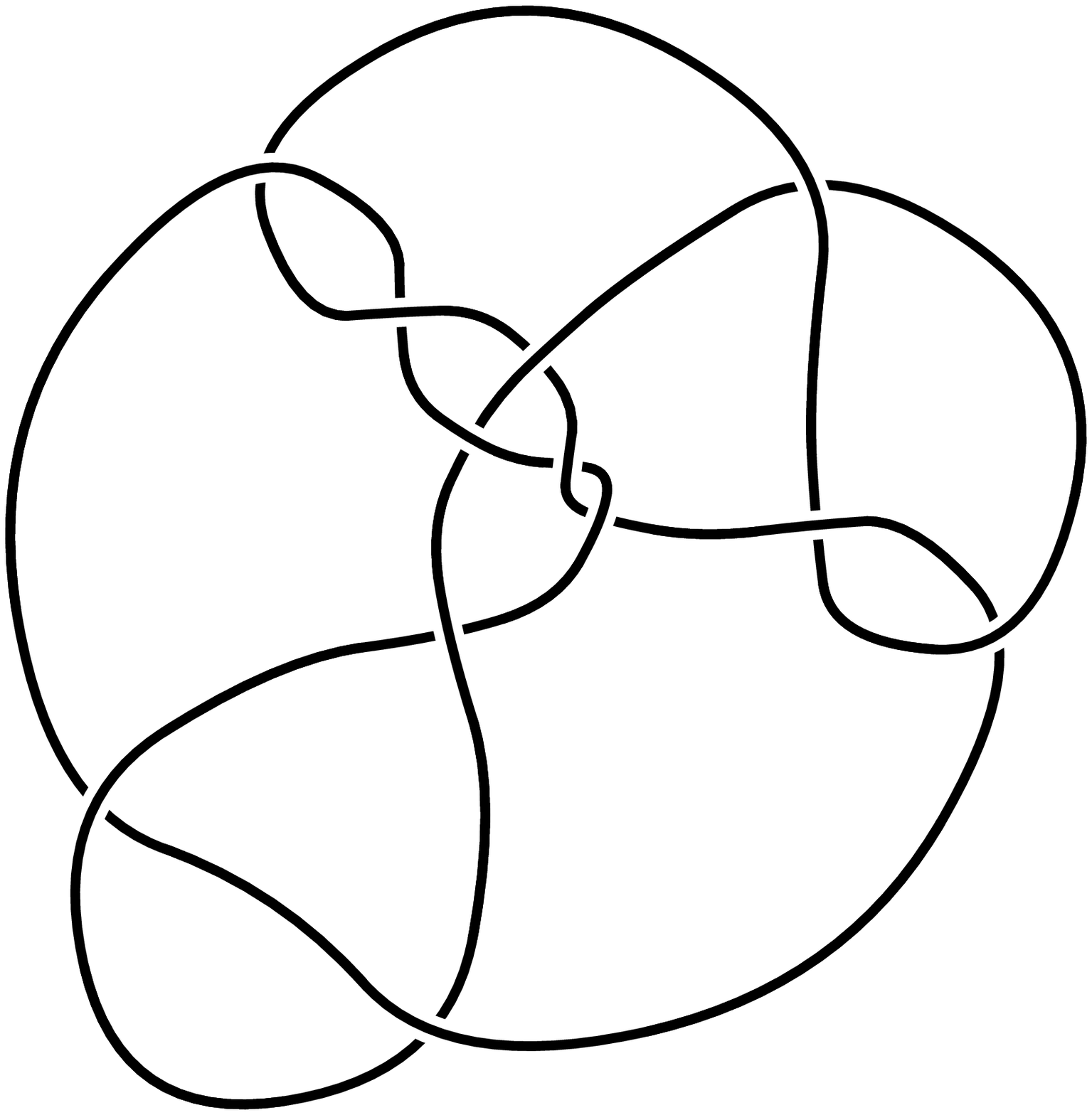}
    &
    \includegraphics[width=75pt]{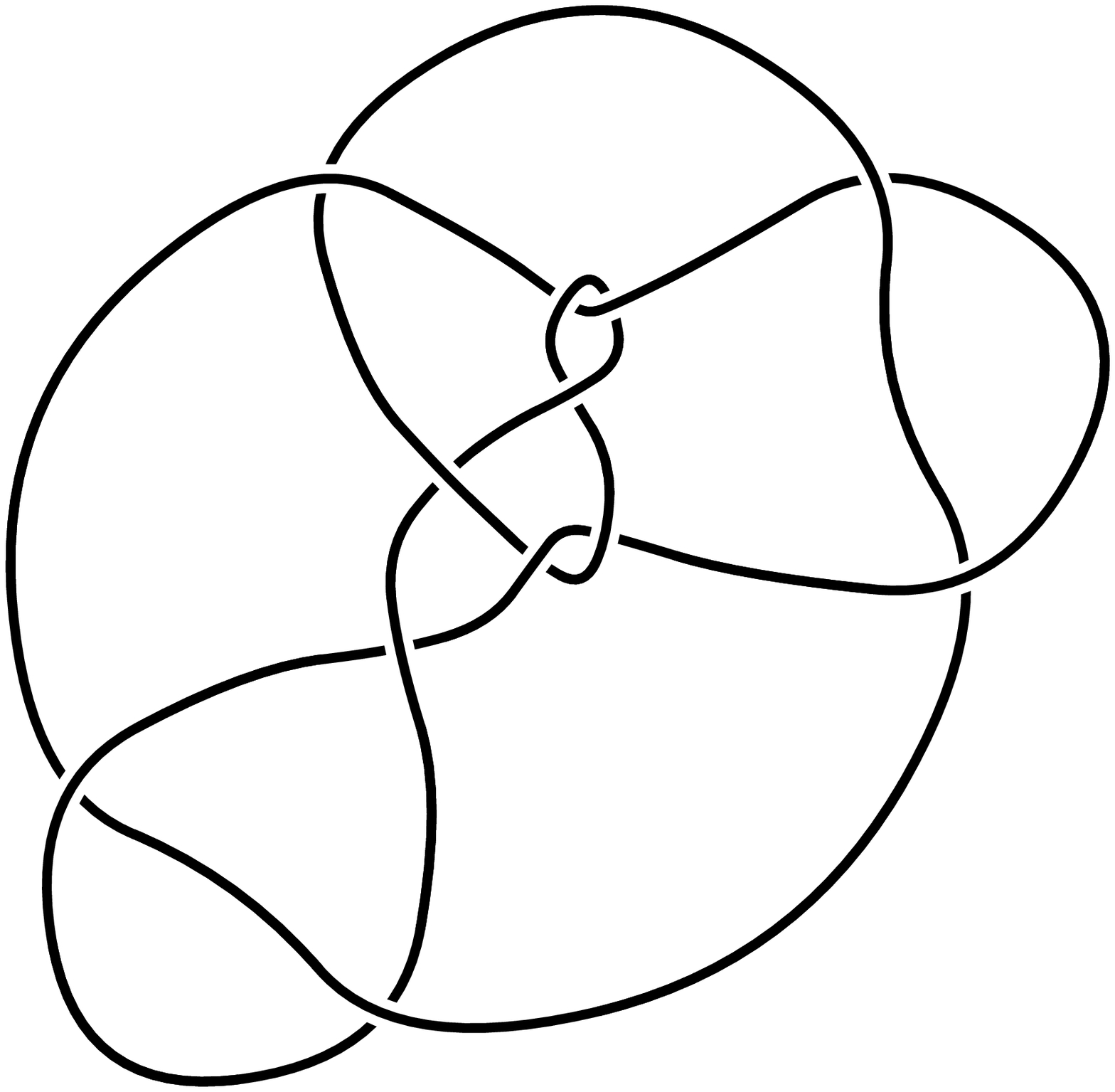}
    &
    \includegraphics[width=75pt]{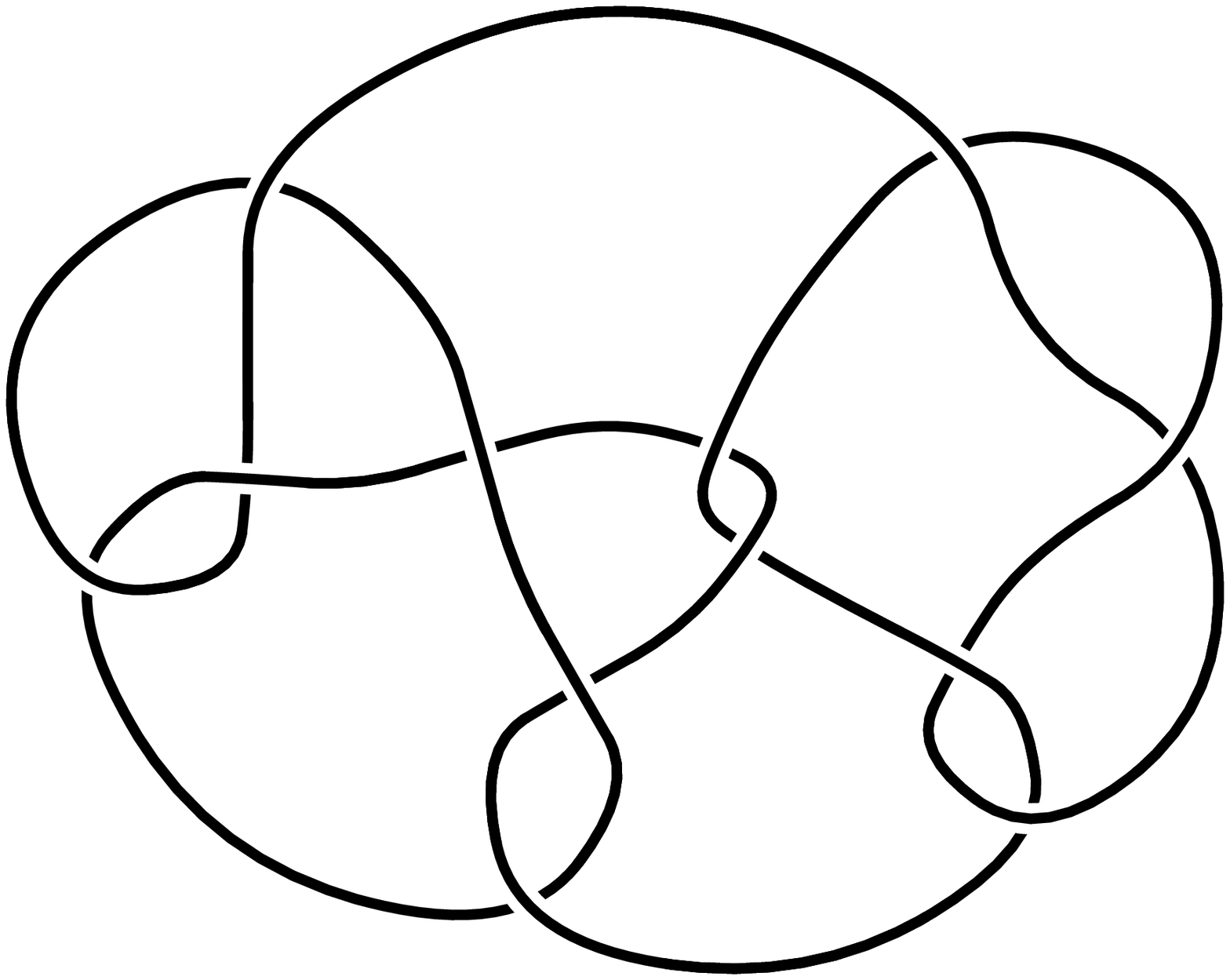}
    &
    \includegraphics[width=75pt]{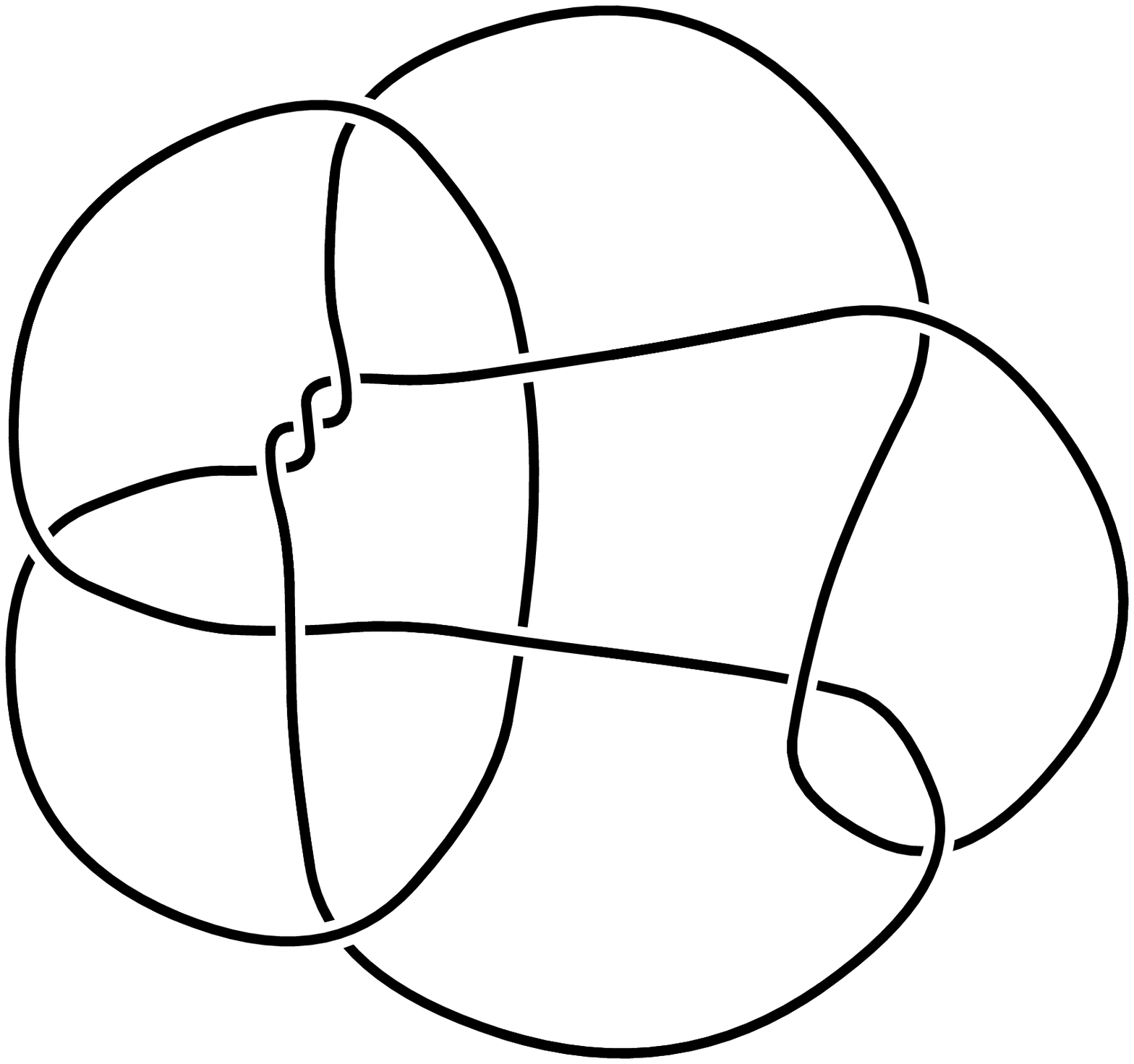}
    \\[-10pt]
    $12^A_{341}$ & $12^A_{627}$ & $12^N_{17}$ & $12^N_{584}$
    \\[10pt]
    \hline
    &&&\\[-10pt]
    \includegraphics[width=75pt]{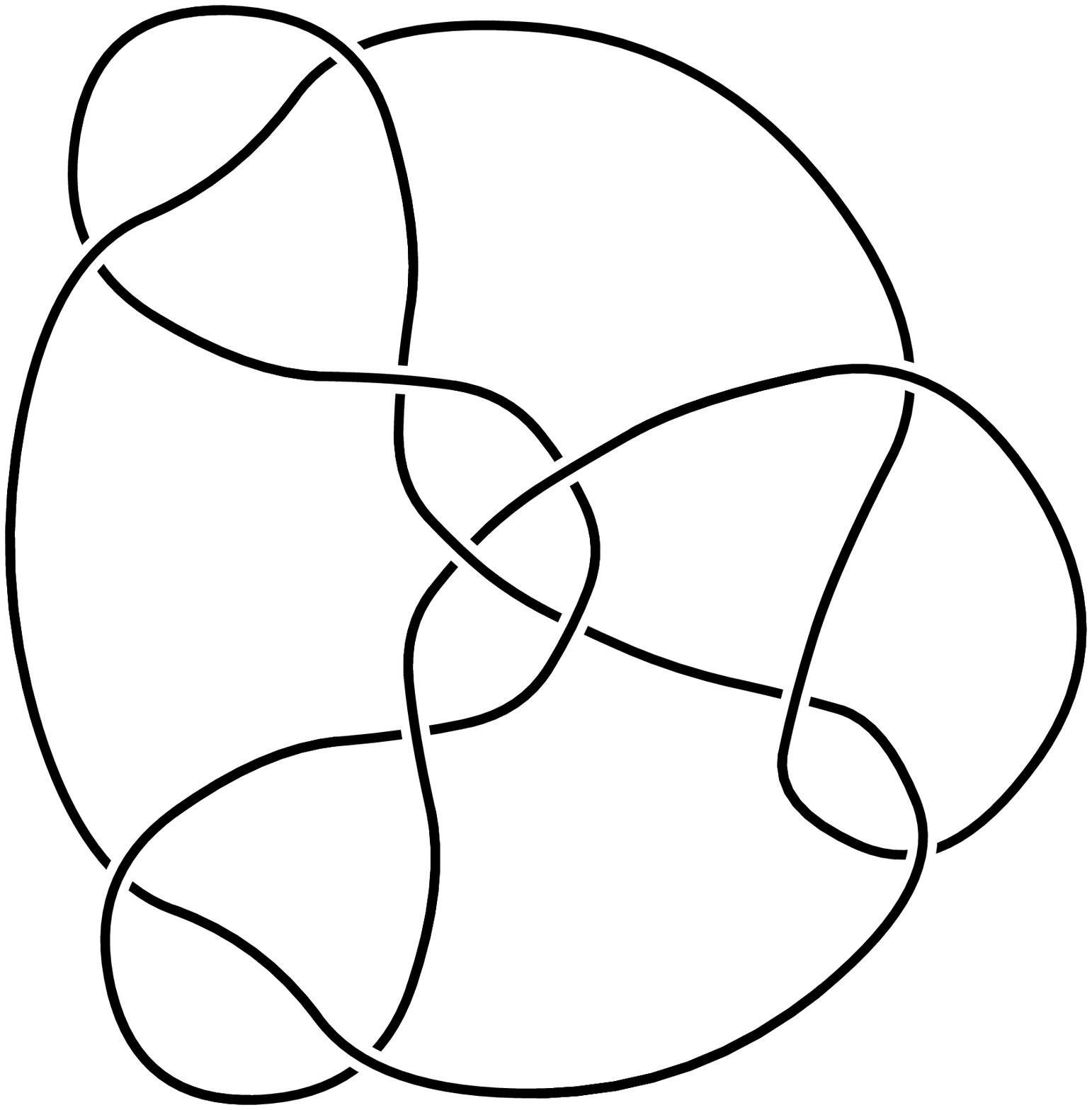}
    &
    \includegraphics[width=75pt]{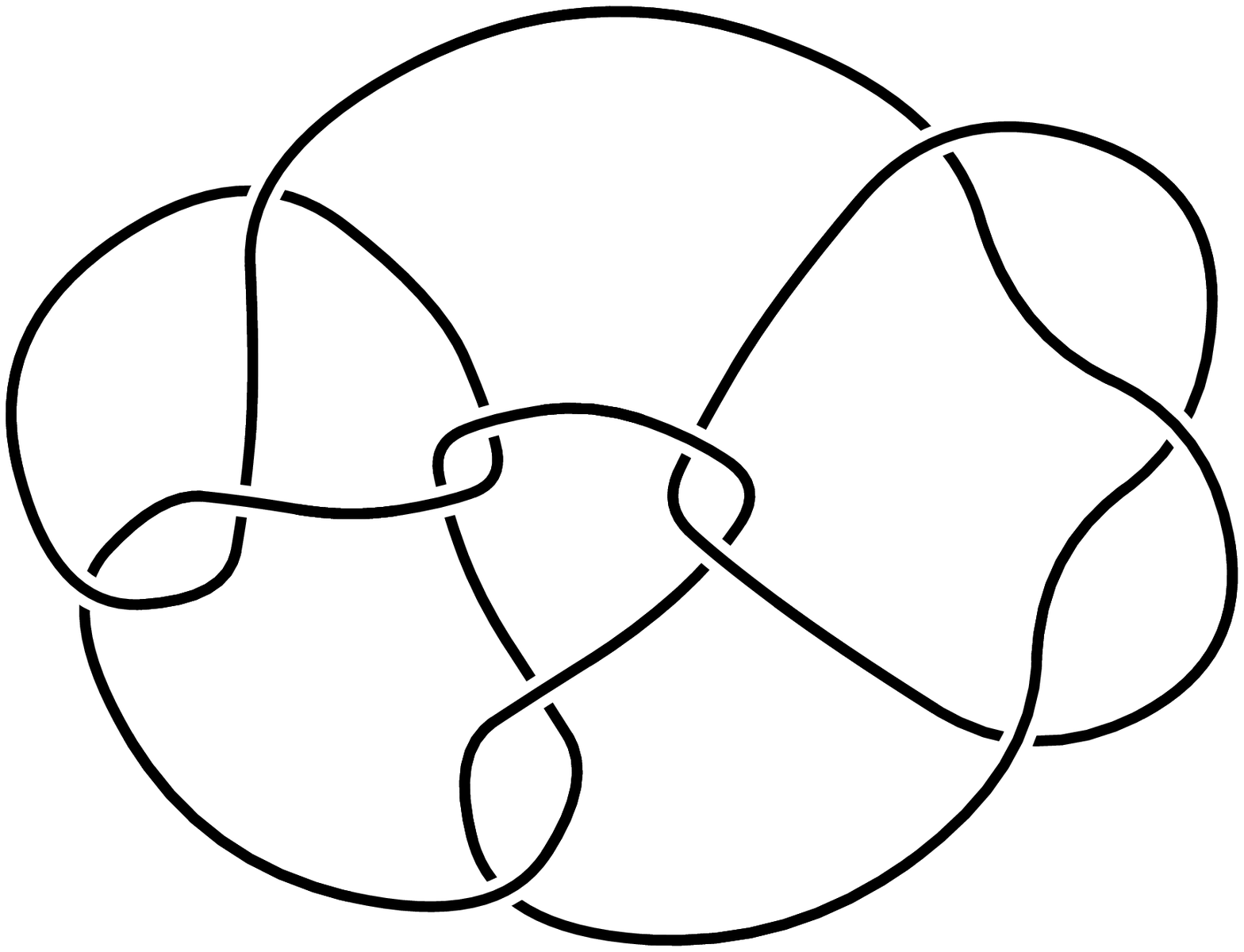}
    &
    \includegraphics[width=75pt]{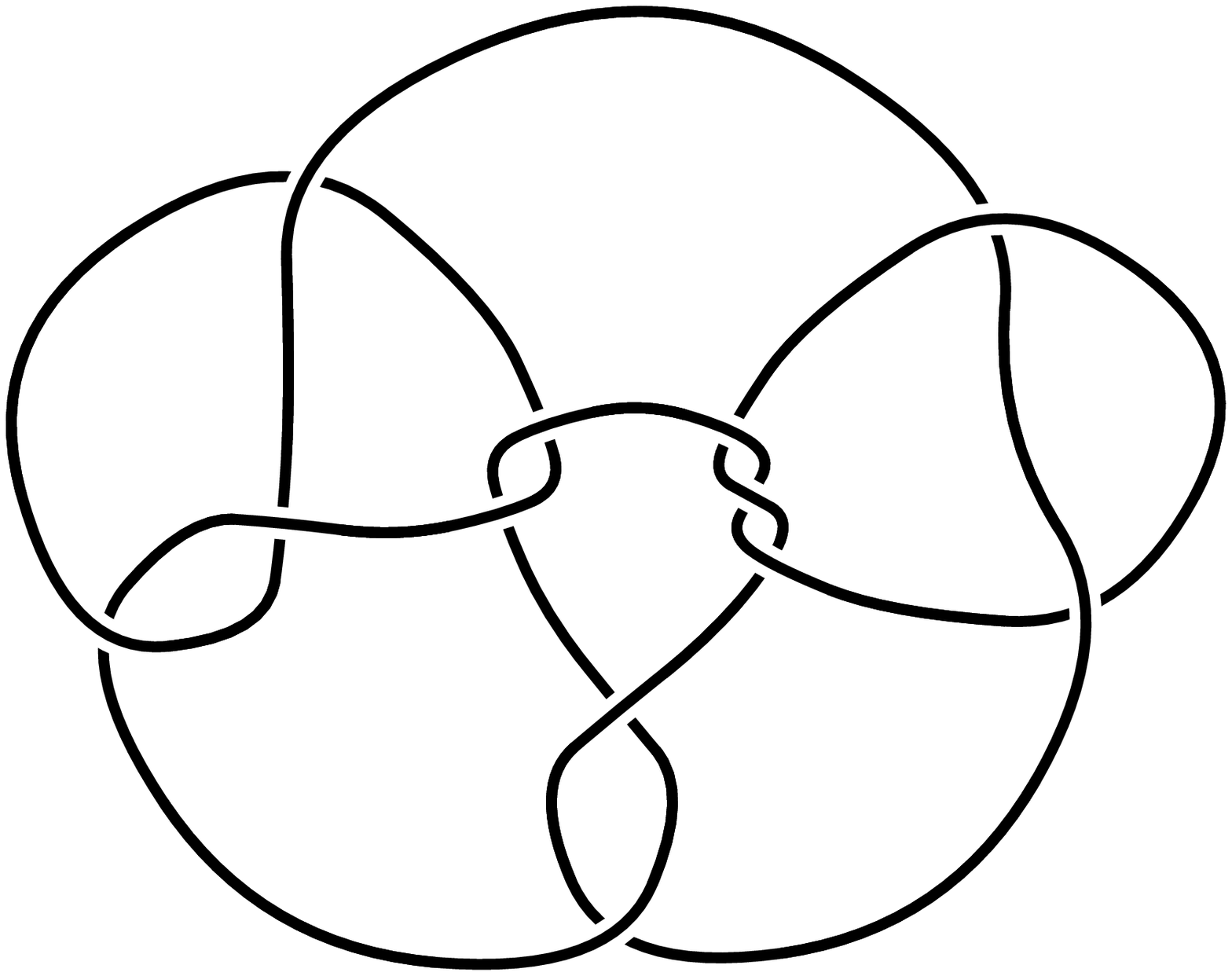}
    &
    \includegraphics[width=75pt]{Graphics/Nonmutants/12/N/12N_416}
    \\[-10pt]
    $12^N_{416}$ & $12^N_{90}$ & $12^N_{135}$ & $12^N_{416}$
  \end{tabular}
  \caption{%
    The nonmutant $LG$-pairs within $\mathbb{K}^N_{12}$.
    A mutation relating $12^N_{90}$ and $12^N_{135}$ is clearly visible.
  }
  \label{figure:NonmutantLGpairsamongstthe12crossingprimeknots}
  \end{centering}
\end{figure}

Within $\mathbb{K}_{12}$, we discover only the following nonmutant $LG$-pairs:
\begin{eqnarray*}
  (12^A_{341},12^A_{627}),
  \quad
  (12^N_{ 17},12^N_{584}),
  \quad
  (12^N_{ 90},12^N_{416})
  \quad
  \mathrm{and}
  \quad
  (12^N_{135},12^N_{416}),
\end{eqnarray*}
which we depict in
Figure~\ref{figure:NonmutantLGpairsamongstthe12crossingprimeknots}. More
precisely, $LG$ cannot distinguish the triple
$(12^N_{90},12^N_{135},12^N_{416})$, but $12^N_{90}$ and $12^N_{135}$ are
mutants.

The respective polynomials are listed below. Note that the first polynomial is
palindromic in $q$, as expected as $12^A_{341}$ and $12^A_{627}$ are achiral.
(Also note that neither of these are $2$-bridge knots.)

\vspace{\baselineskip}

\tiny

\begin{tabular}{p{40pt}|@{\hspace{20pt}}p{240pt}}
  \\
  {\normalsize
  $
  \begin{array}{c}
    12^A_{341}
    \\
    12^A_{627}
  \end{array}
  $}
  &
  $
  \begin{array}{@{\hspace{0pt}}r@{\hspace{5pt}}r@{\hspace{2pt}}l}
      &                        & (12\overline{q}^{6} + 476\overline{q}^{4} + 2474\overline{q}^{2} + 4081 + 2474 q^{2} + 476 q^{4} + 12 q^{6})
    \\
    - & (\overline{p}^2 + p^2) & (103\overline{q}^{5} + 1147\overline{q}^{3} + 3229\overline{q} + 3229 q + 1147 q^{3} + 103 q^{5})
    \\
    + & (\overline{p}^4 + p^4) & (10\overline{q}^{6} + 303\overline{q}^{4} + 1580\overline{q}^{2} + 2623 + 1580 q^{2} + 303 q^{4} + 10 q^{6})
    \\
    - & (\overline{p}^6 + p^6) & (37\overline{q}^{5} + 449\overline{q}^{3} + 1325\overline{q} + 1325 q + 449 q^{3} + 37 q^{5})
    \\
    + & (\overline{p}^8 + p^8) & (58\overline{q}^{4} + 389\overline{q}^{2} + 688 + 389 q^{2} + 58 q^{4})
    \\
    - & (\overline{p}^{10} + p^{10}) & (51\overline{q}^{3} + 206\overline{q} + 206 q + 51 q^{3})
    \\
    + & (\overline{p}^{12} + p^{12}) & (27\overline{q}^{2} + 62 + 27 q^{2})
    \\
    - & (\overline{p}^{14} + p^{14}) & (8\overline{q} + 8 q)
    \\
    + & (\overline{p}^{16} + p^{16}) & (1)
  \end{array}
  $
  \\
  \\
  \hline
  \\
  {\normalsize
  $
  \begin{array}{c}
    12^N_{17}
    \\
    12^N_{584}
  \end{array}
  $}
  &
  $
  \begin{array}{@{\hspace{0pt}}r@{\hspace{5pt}}r@{\hspace{2pt}}l}
      &                        & (12\overline{q}^{4} + 188\overline{q}^{2} +  647 + 774 q^{2} + 296 q^{4} + 22 q^{6})
    \\
    - & (\overline{p}^2 + p^2) & (44\overline{q}^{3} + 328\overline{q} + 685 q + 492 q^{3} + 93 q^{5} + 2 q^{7})
    \\
    + & (\overline{p}^4 + p^4) & (67\overline{q}^{2} + 313 + 427 q^{2} + 172 q^{4} + 13 q^{6})
    \\
    - & (\overline{p}^6 + p^6) & (56\overline{q} + 179 q + 149 q^{3} + 26 q^{5})
    \\
    + & (\overline{p}^8 + p^8) & (28 + 58 q^{2} + 22 q^{4})
    \\
    - & (\overline{p}^{10} + p^{10}) & (8 q + 8 q^{3})
    \\
    + & (\overline{p}^{12} + p^{12}) & (q^{2})
  \end{array}
  $
  \\
  \\
  \hline
  \\
  {\normalsize
  $
  \begin{array}{c}
    12^N_{90}
    \\
    12^N_{135}
    \\
    12^N_{416}
  \end{array}
  $}
  &
  $
  \begin{array}{@{\hspace{0pt}}r@{\hspace{5pt}}r@{\hspace{2pt}}l}
      &                        & (2\overline{q}^{8} + 72\overline{q}^{6} + 292\overline{q}^{4} + 201 + 374\overline{q}^{2} + 32 q^{2})
    \\
    - & (\overline{p}^2 + p^2) & (18\overline{q}^{7} + 168\overline{q}^{5} + 348\overline{q}^{3} + 277\overline{q} + 83 q + 4 q^{3})
    \\
    + & (\overline{p}^4 + p^4) & (2\overline{q}^{8} + 64\overline{q}^{6} + 240\overline{q}^{4} + 268\overline{q}^{2} + 118 + 13 q^{2})
    \\
    - & (\overline{p}^6 + p^6) & (14\overline{q}^{7} + 118\overline{q}^{5} + 199\overline{q}^{3} + 115\overline{q} + 20 q)
    \\
    + & (\overline{p}^8 + p^8) & (1\overline{q}^{8} + 36\overline{q}^{6} + 111\overline{q}^{4} + 86\overline{q}^{2} + 2)
    \\
    - & (\overline{p}^{10} + p^{10}) & (5\overline{q}^{7} + 40\overline{q}^{5} + 50\overline{q}^{3} + 15\overline{q})
    \\
    + & (\overline{p}^{12} + p^{12}) & (7\overline{q}^{6} + 20\overline{q}^{4} + 9\overline{q}^{2})
    \\
    - & (\overline{p}^{14} + p^{14}) & (4\overline{q}^{5} + 4\overline{q}^{3})
    \\
    + & (\overline{p}^{16} + p^{16}) & (\overline{q}^{4})
  \end{array}
  $
  \\
  \\
\end{tabular}

\normalsize

\vspace{\baselineskip}

Thus, we learn that the lowest crossing number of a nonmutant pair of prime
knots indistinguishable by $LG$ is $12$. We remark that there are no examples
of $LG$-pairs of knots of \emph{different} crossing numbers or different
alternatingness amongst the currently $LG$-evaluated prime knots. We also
confirm that $LG$ is not unity for any of these but the unknot.

Note that none of the four nonmutant $LG$-pairs within $\mathbb{K}_{12}$ are distinguished by
either the HOMFLYPT or Kauffman polynomials (they are distinguished by the hyperbolic volume,
which demonstrates that they are not mutants). More generally, it remains an open question whether
$LG$ ever fails to distinguish a pair of links distinguished by \emph{either} the HOMFLYPT or the
Kauffman polynomial (that is, Problem~6.4 of \cite{IshiiKanenobu:2004b}).


\section{Where $LG$ first fails to distinguish reflections}
\label{section:WhereLGfirstfailstodistinguishreflections}

\begin{figure}[ht]
  \begin{centering}
   \includegraphics[width=120pt]{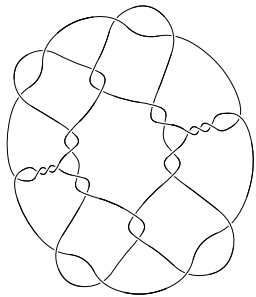}
  \caption{%
    A $28$-crossing projection of a chiral prime knot whose chirality is undetected by $LG$.
  }
  \label{figure:IK28}
  \end{centering}
\end{figure}

It is known that $LG$ is not always able to distinguish the reflections of
chiral links. Specifically, Ishii and Kanenobu~\cite{IshiiKanenobu:2004b} have
constructed an infinite family of chiral prime knots whose chirality is
undetected by $LG$. The smallest example they provide is a prime knot with a
$28$-crossing projection.%
\footnote{%
  In their notation, this knot is
  $L[\sigma_2^{-2}\sigma_1^{2};\frac{1}{2},\frac{1}{4},-\frac{1}{2},-\frac{1}{4}]$.
} %
\textsc{Knotscape} cannot reduce this projection to one of less than $28$
crossings, but does reduce it to the elegant diagram of
Figure~\ref{figure:IK28}.

We next wish to discover where $LG$ \emph{first} fails to distinguish the
chirality of a chiral prime knot, and to that end we search within
$\mathbb{K}_{11}$--$\mathbb{K}_{16}$ as we know that $LG$ distinguishes the
chiralities of all chiral prime knots of up to $10$
crossings~\cite{DeWit:2000}. Clearly, we require foreknowledge of \emph{which}
knots are chiral, and fortunately \textsc{Knotscape}, together
with~\cite{HosteThistlethwaiteWeeks:1998}, provides us with information of the
symmetry classes (chiral or not, invertible or not) of the HTW knots.

Firstly, \textsc{Knotscape} uses the embedded \textsc{SnapPea} to determine the
symmetry classes of the hyperbolic HTW knots. \textsc{SnapPea} determines the
symmetry class of a given hyperbolic knot $K$ via construction of its symmetry
group $\mathrm{Sym}(S^3,K)$, that is, the group of homotopy-equivalence
classes of homeomorphisms of $(S^3,K)$. The symmetry group is unlikely to be
discernible by casual inspection of an arbitrary projection; indeed, for any
given link, it remains an open question as to whether there exists a
(maximally-symmetric) projection in which the symmetry group is visible. It
turns out that the symmetry groups of hyperbolic knots are all finite groups,
either cyclic or dihedral~\cite[Chapter~10]{Kawauchi:1996}, as per those of
the regular plane figures. Rather than by inspection of diagrams,
\textsc{SnapPea} actually computes $\mathrm{Sym}(S^3,K)$ from the hyperbolic
decomposition of the knot complement $S^3-K$. In contrast, the symmetry groups
of nonhyperbolic knots are more general. Those of the nonhyperbolic HTW knots
are deduced manually in~\cite{HosteThistlethwaiteWeeks:1998}; they are all
dihedral groups, some finite (in particular each torus knot has group $D_1$),
and some infinite. The nonhyperbolic HTW knots are all chiral and invertible.

So, amongst the chiral knots within each set $\mathbb{K}^P_c$, we search for
knots whose $LG$ polynomials are palindromic, indicating that they fail to
distinguish chirality. In this manner, we quickly discover that $LG$
distinguishes the chirality of all chiral prime knots of up to $12$ crossings.
However, we do not have evaluations of $LG$ for any complete set
$\mathbb{K}^P_c$ beyond that point, and we uncover no failures amongst the
knots within $\mathbb{K}_{13}$ (all of which are chiral) for which we have
evaluations of $LG$.

Instead, we adopt a different approach to the problem. We begin by observing
that mutation can change a chiral knot into an achiral knot; an illustration
of such a pair of $16$-crossing mutant projections appears
in~\cite{Kauffman:1990} (this example is originally due to
Kanenobu~\cite{Kanenobu:1981}). We here confirm that the pair are
$16$-crossing prime knots: the chiral is $16^A_{259088}$ and the achiral is
$16^A_{259984}$. This example illustrates a more general situation.

\begin{lemma}
  Mutation can change the symmetry group (and thus the symmetry class) of a
  link.
\end{lemma}

The first example of this phenomenon amongst the prime knots is the sole
instance of it within $\mathbb{K}_{11}$, in the pair $(11^N_{39},11^N_{45})$.
Although both are chiral, the former is noninvertible (with symmetry group
$\mathbb{Z}_2$) and the latter invertible (with symmetry group $D_1$).

Now, say that we have a mutant pair of links, one chiral, the other achiral.
Any mutation-insensitive link invariant will necessarily assign the same value
to the achiral link and to both the chiralities of the chiral link, and thus
fail to distinguish the chirality of the latter. Thus, the fact that mutation
can change the symmetry group of a link means that a link invariant which
cannot distinguish mutants will also fail to distinguish some chiral links
from their reflections. Incompleteness for mutation implies incompleteness for
chirality. We formalise this as the following lemma.

\begin{lemma}
  Any mutation-insensitive link invariant is necessarily also incomplete for
  chirality.
\end{lemma}

So, we seek an example of such a mutant pair of minimal crossing number amongst
the prime knots. To that end, we may at first examine \emph{candidate} mutant
cliques (that is, cliques of common (apparent) hyperbolic volume, HOMFLYPT and
Kauffman polynomial), and we don't need to prove that a given chiral knot is a
member of a mutant clique before we examine its $LG$-polynomial. Thus, we
search the HTW knots for candidate mutant cliques containing at least one
chiral and at least one achiral element, and of them, we need only the chiral
elements. The first point at which we find such cliques is at $14$ crossings.
(There are no achiral $13$-crossing prime knots.) Curiously, all such
$14$-crossing candidate mutant cliques manifest in triples, each containing
exactly one chiral element. Inspection demonstrates that all candidates
actually are mutants. Moreover, as oriented knots, all the elements of the
triples are noninvertible, and so the achiral are actually all negative
achiral, meaning that they are equal to their inverse reflections, but not
equal to their reflections per se. We list the $13$ triples found in
Table~\ref{table:interesting14-crossingmutantcliques}; their diagrams are
included in the draft of this paper lodged with the LANL arXiv.

\renewcommand{\arraystretch}{1.3}
\renewcommand{\tabcolsep}{4pt}

\begin{table}[hbt]
  \small
  \begin{centering}
  \begin{tabular}{l|l@{\hspace{5pt}}l}
    \multicolumn{1}{c}{chiral} & \multicolumn{2}{c}{achiral mutants} \\
    \hline
     $14^A_{  506}$ & $14^A_{  486}$ & $14^A_{  731}$ \\
     $14^A_{  680}$ & $14^A_{  509}$ & $14^A_{  585}$ \\
     $14^A_{12813}$ & $14^A_{12807}$ & $14^A_{12875}$ \\
     $14^A_{12858}$ & $14^A_{12815}$ & $14^A_{12830}$ \\
     $14^A_{13107}$ & $14^A_{13109}$ & $14^A_{13489}$ \\
     $14^A_{13262}$ & $14^A_{13269}$ & $14^A_{13506}$ \\
     $14^A_{14042}$ & $14^A_{14043}$ & $14^A_{14671}$ \\
     $14^A_{17268}$ & $14^A_{17265}$ & $14^A_{17275}$ \\
     $14^A_{17533}$ & $14^A_{17531}$ & $14^A_{17680}$ \\
     $14^N_{ 1309}$ & $14^N_{ 1327}$ & $14^N_{ 1497}$ \\
     $14^N_{ 1641}$ & $14^N_{ 1552}$ & $14^N_{ 2132}$ \\
     $14^N_{ 1644}$ & $14^N_{ 1555}$ & $14^N_{ 1671}$ \\
     $14^N_{ 2164}$ & $14^N_{ 1669}$ & $14^N_{ 1925}$
  \end{tabular}
  \caption{%
    Mutant cliques within $\mathbb{K}_{14}$ containing both chiral and achiral
    elements.
  }%
  \label{table:interesting14-crossingmutantcliques}
  \end{centering}
\end{table}

So, the first point at which mutation changes a chiral prime knot into an
achiral prime knot is at $14$ crossings. Thus, we have the following.

\begin{lemma}
  Any mutation-insensitive link invariant will fail to distinguish the
  chirality of some $14$-crossing chiral prime knots.
\end{lemma}

\begin{figure}[htbp]
  \begin{centering}
  \includegraphics[width=120pt]{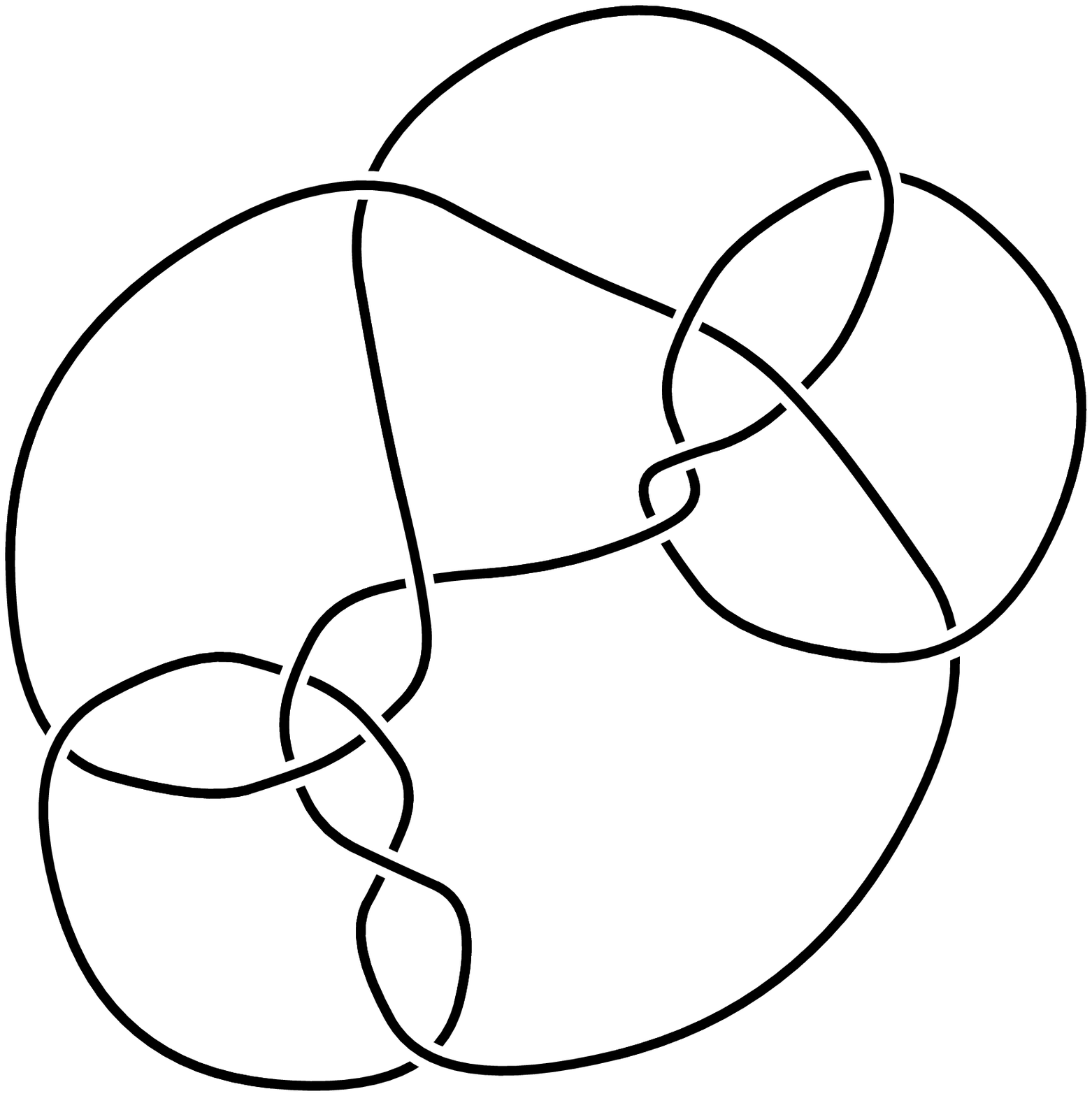}
  \hspace{20pt}
  \includegraphics[width=120pt]{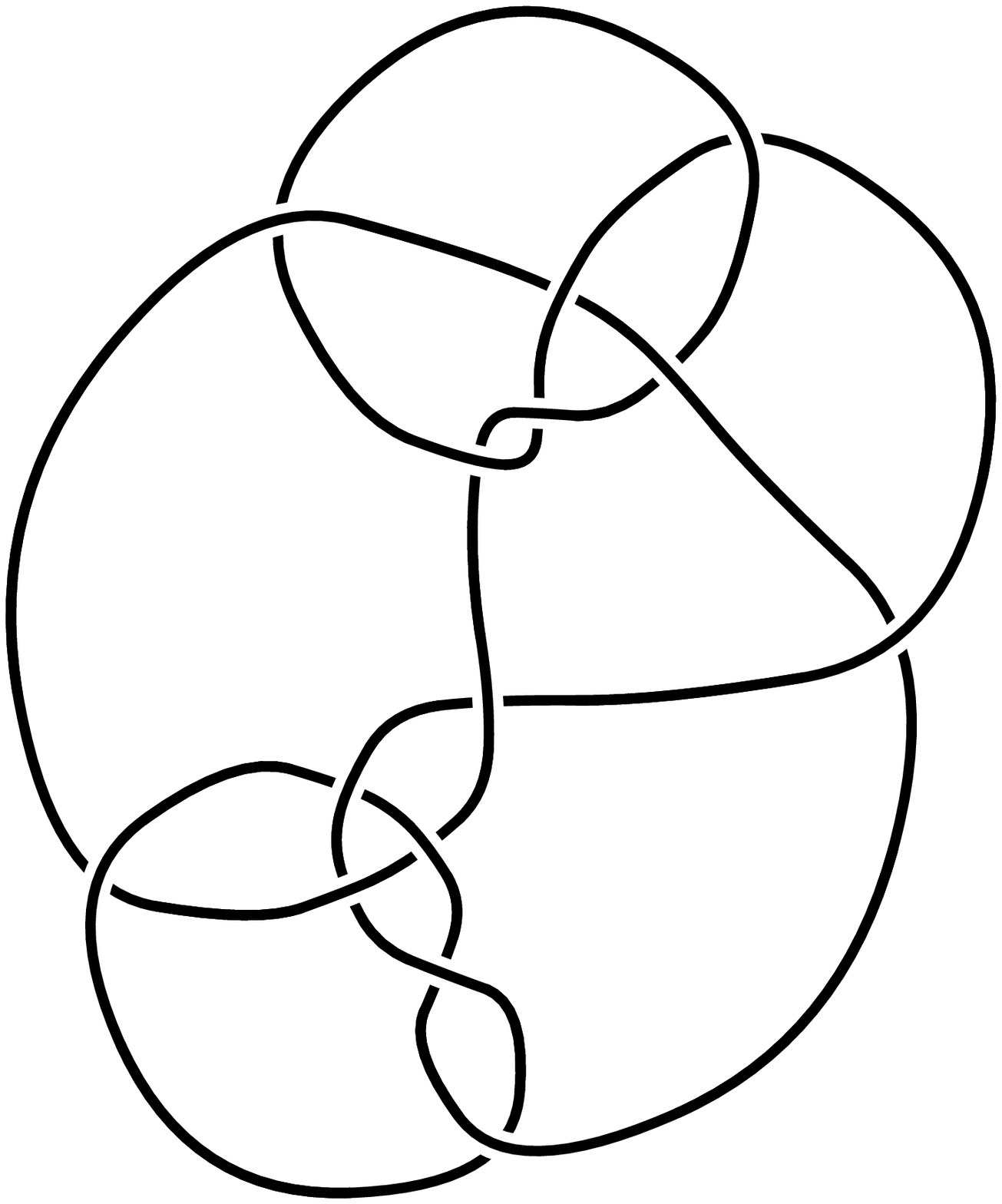}
  \caption{%
    The chiral $14^A_{13107}$ and the achiral $14^A_{13109}$ are mutants.
  }
  \label{figure:14A13107and14A13109}
  \end{centering}
\end{figure}

Such an invariant will also fail for other prime knots, some of which may be of
lower crossing number. The best that can be hoped for is that it is able to
distinguish the chirality of all chiral prime knots of up to $13$ crossings.
We illustrate using the chiral $14^A_{13107}$ and the achiral $14^A_{13109}$.
A mutation relating them is visible in Figure~\ref{figure:14A13107and14A13109}.
Indeed, $LG$ cannot distinguish this pair, and is insensitive to the chirality
of $14^A_{13107}$, yielding the following palindromic polynomial.

\tiny
\begin{eqnarray*}
  \begin{array}{@{\hspace{0pt}}r@{\hspace{5pt}}r@{\hspace{2pt}}l}
      &                              & (348\overline{q}^6 + 5106\overline{q}^4 + 19778\overline{q}^2 + 30201 + 19778 q^2 + 5106 q^4 + 348 q^6) \\
    - & (\overline{p}^{ 2} + p^{ 2}) & (33\overline{q}^7 + 1525\overline{q}^5 + 10516\overline{q}^3 + 25021\overline{q} + 25021 q + 10516 q^3 + 1525 q^5 + 33 q^7) \\
    + & (\overline{p}^{ 4} + p^{ 4}) & (248\overline{q}^6 + 3615\overline{q}^4 + 14119\overline{q}^2 + 21641 + 14119 q^2 + 3615 q^4 + 248 q^6) \\
    - & (\overline{p}^{ 6} + p^{ 6}) & (15\overline{q}^7 + 724\overline{q}^5 + 5241\overline{q}^3 + 12779\overline{q} + 12779 q + 5241 q^3 + 724 q^5 + 15 q^7) \\
    + & (\overline{p}^{ 8} + p^{ 4}) & (62\overline{q}^6 + 1156\overline{q}^4 + 4982\overline{q}^2 + 7861 + 4982 q^2 + 1156 q^4 + 62 q^6) \\
    - & (\overline{p}^{10} + p^{10}) & (113\overline{q}^5 + 1159\overline{q}^3 + 3182\overline{q} + 3182 q + 1159 q^3 + 113 q^5) \\
    + & (\overline{p}^{12} + p^{12}) & (121\overline{q}^4 + 766\overline{q}^2 + 1326 + 766 q^2 + 121 q^4) \\
    - & (\overline{p}^{14} + p^{14}) & (83\overline{q}^3 + 326\overline{q} + 326 q + 83 q^3) \\
    + & (\overline{p}^{16} + p^{16}) & (36\overline{q}^2 + 81 + 36 q^2) \\
    - & (\overline{p}^{18} + p^{18}) & (9\overline{q} + 9 q) \\
    + & (\overline{p}^{20} + p^{20}) & (1)
   \end{array}
\end{eqnarray*}
\normalsize

So, as each $LG^{m,n}$ is insensitive to mutation, it will certainly fail to
detect the chirality of some $14$-crossing prime knots. In particular, we now
have $14$ as an upper bound on the point at which $LG^{2,1}$ first fails to
distinguish the chirality of a chiral prime knot. Note that whilst the
mutation-changing-chirality-to-achirality phenomenon doesn't occur within
$\mathbb{K}_{13}$ (all the knots of which are chiral), there may still be a
$13$-crossing prime knot whose chirality is undistinguished by $LG$. Our
evaluations of $LG$ for 25\% of these knots haven't uncovered such an example;
indeed a wider search within $\mathbb{K}_{13}$ of the $LG$-feasible knots
whose chirality is undetected by the HOMFLYPT or Kauffman polynomials also
uncovers no such examples. Moreover, of the $53,418$ HTW knots for which we
have evaluations of $LG$ (most of which are chiral), the \emph{only} examples
of chiral knots whose chirality is undetected by $LG$ are $8$ already contained
within Table~\ref{table:interesting14-crossingmutantcliques}:
\begin{eqnarray*}
  14^A_{  680},
  14^A_{12813},
  14^A_{12858},
  14^A_{13107},
  14^A_{13262},
  14^A_{17268},
  14^N_{ 1309},
  14^N_{ 2164}.
\end{eqnarray*}

In contrast to the statements about $LG$, the first example of a chiral prime
knot whose chirality is undetected by both the HOMFLYPT and Kauffman
polynomials is the well-known case $9^N_4 \equiv 9_{42}$. (This is the only
such case where \emph{either} polynomial fails within $\mathbb{K}_9$.) We know
that the chirality of $9^N_4$ is detected by $LG$, and that it has no mutants
amongst the HTW prime knots. We then conclude with some open questions.

\begin{enumerate}
\item
  Does there exist a chiral prime knot whose chirality is undetected by $LG$
  yet is detected by the HOMFLYPT or Kauffman polynomials?

\item
  Can we find a chiral link whose chirality is undetected by $LG$ yet is
  \emph{not} a mutant of an achiral link? Recalling that we don't know whether
  mutation preserves crossing number or alternatingness, an easier version of
  this problem is to find a chiral prime knot whose chirality is
  undetected by $LG$ and has no mutants of the same alternatingness and
  crossing number.

\item
  Recall that the question of whether there exists a chiral link which is
  mutation-equivalent to its reflection appears to remain unanswered. Can we
  find such a link? As every mutation-insensitive link invariant will be
  unable to detect the chirality of such a link, and as $LG$ detects the
  chirality of all prime knots of up to $12$ crossings, there can be no such
  links amongst the prime knots of up to $12$ crossings. Candidates for such
  links might be found amongst the chiral links whose chirality is
  undistinguished by all mutation-insensitive link invariants (for example
  $14^A_{13107}$). Then again, are such links \emph{necessarily} mutants of
  their reflections?
\end{enumerate}


\subsection*{Acknowledgements}

The knot diagrams were made by \textsc{Knotscape}. Jon Links thanks the Australian Research
Council for support through an Australian Research Fellowship. The authors thank Taizo Kanenobu
for the keen observation of an error in an earlier draft of this paper.


\bibliographystyle{plain}
\bibliography{Where_LG_Fails}

\vspace{2\baselineskip}

\section*{Appendix}

\begin{figure}[htbp]
  \begin{centering}
  \begin{tabular}{cc@{\hspace{10pt}}|@{\hspace{10pt}}cc}
    \includegraphics[width=75pt]{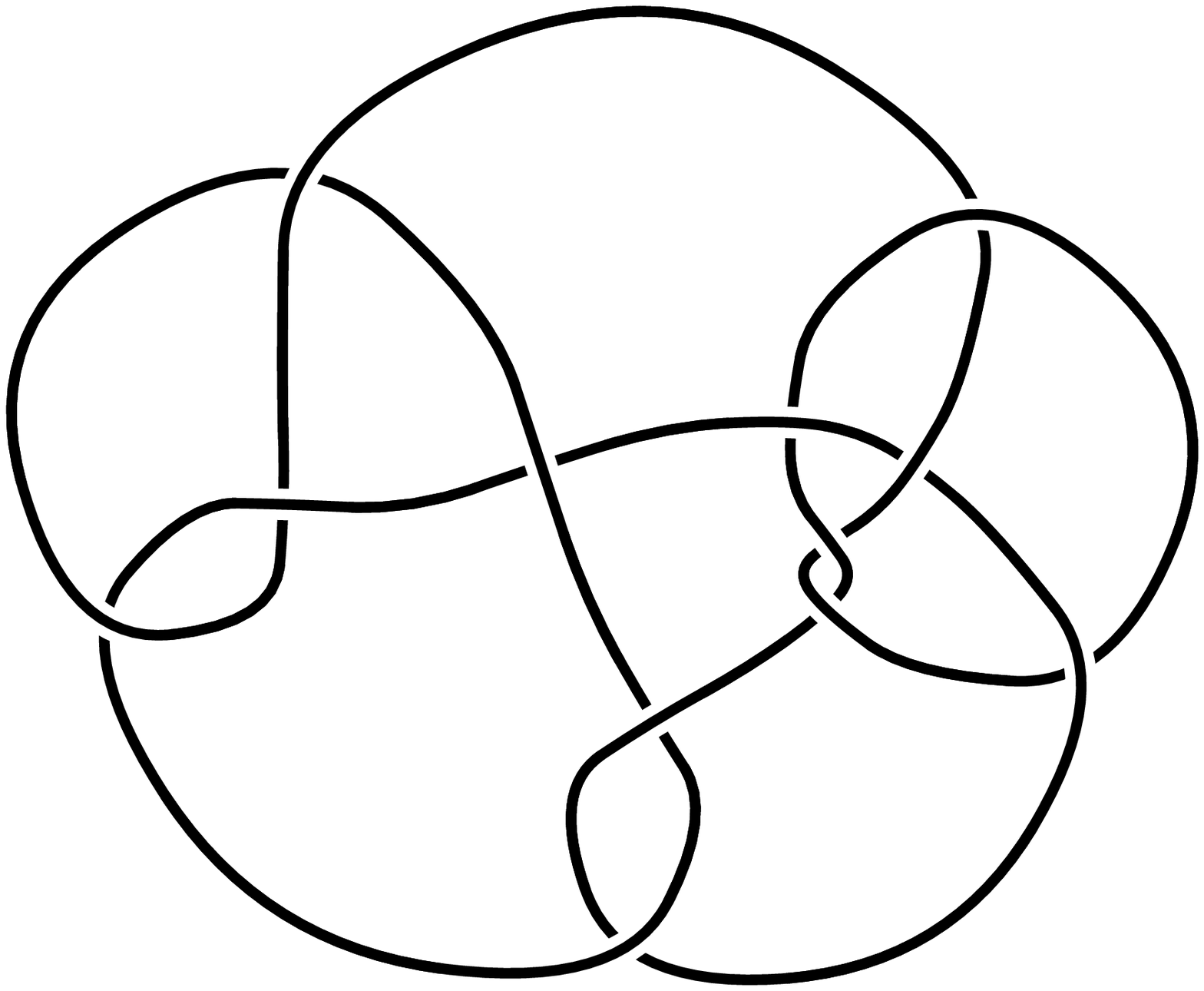}
    &
    \includegraphics[width=75pt]{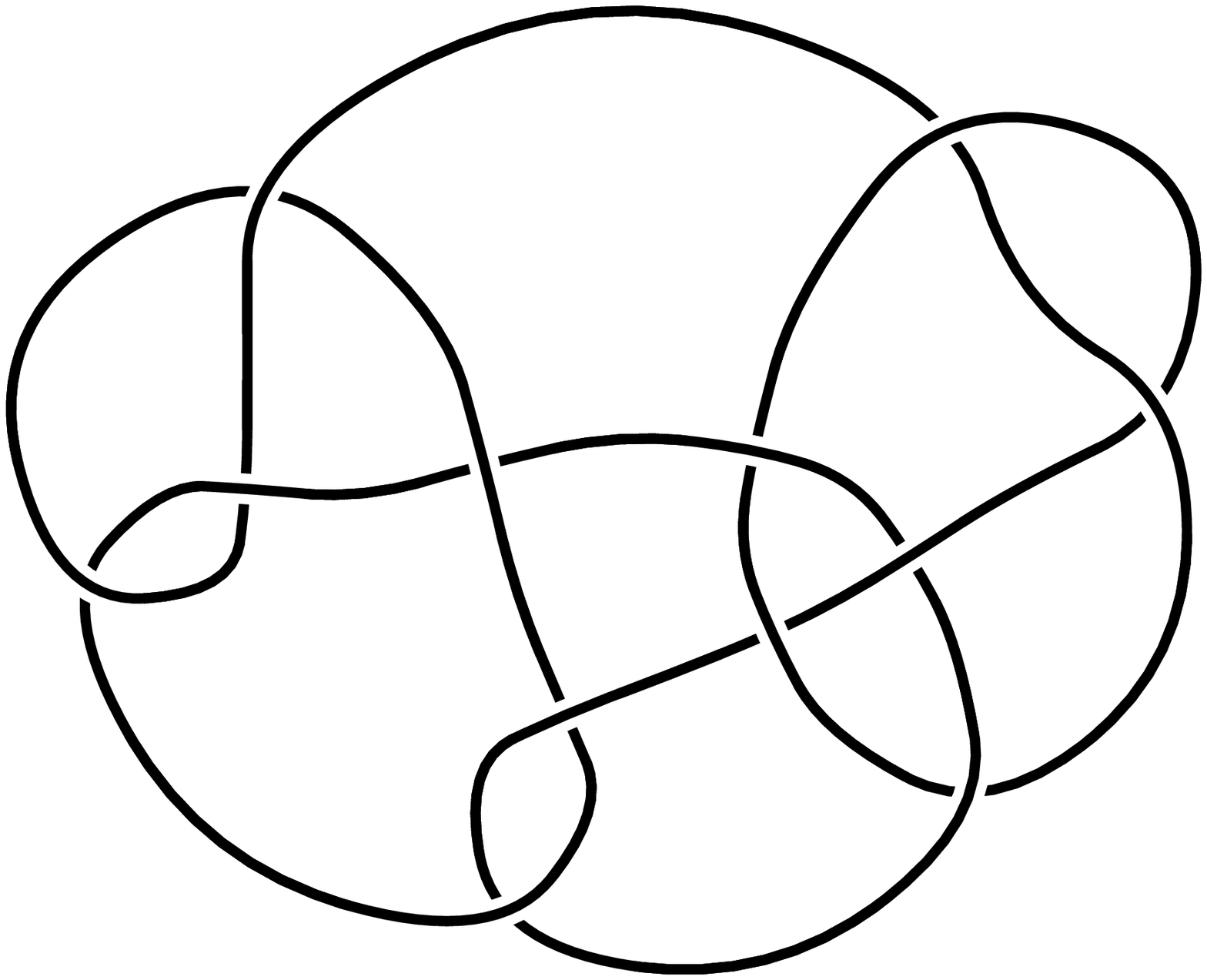}
    &
    \includegraphics[width=75pt]{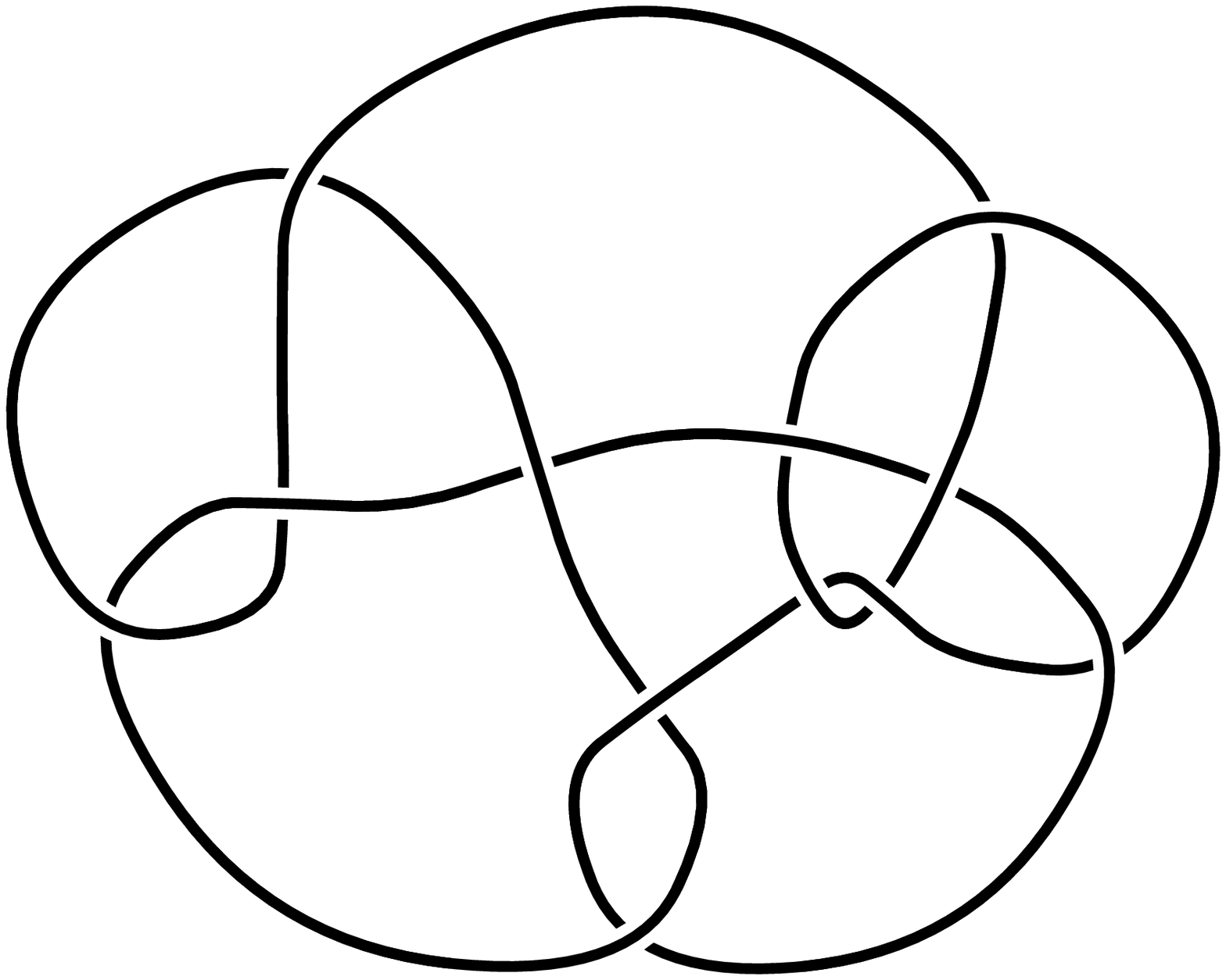}
    &
    \includegraphics[width=75pt]{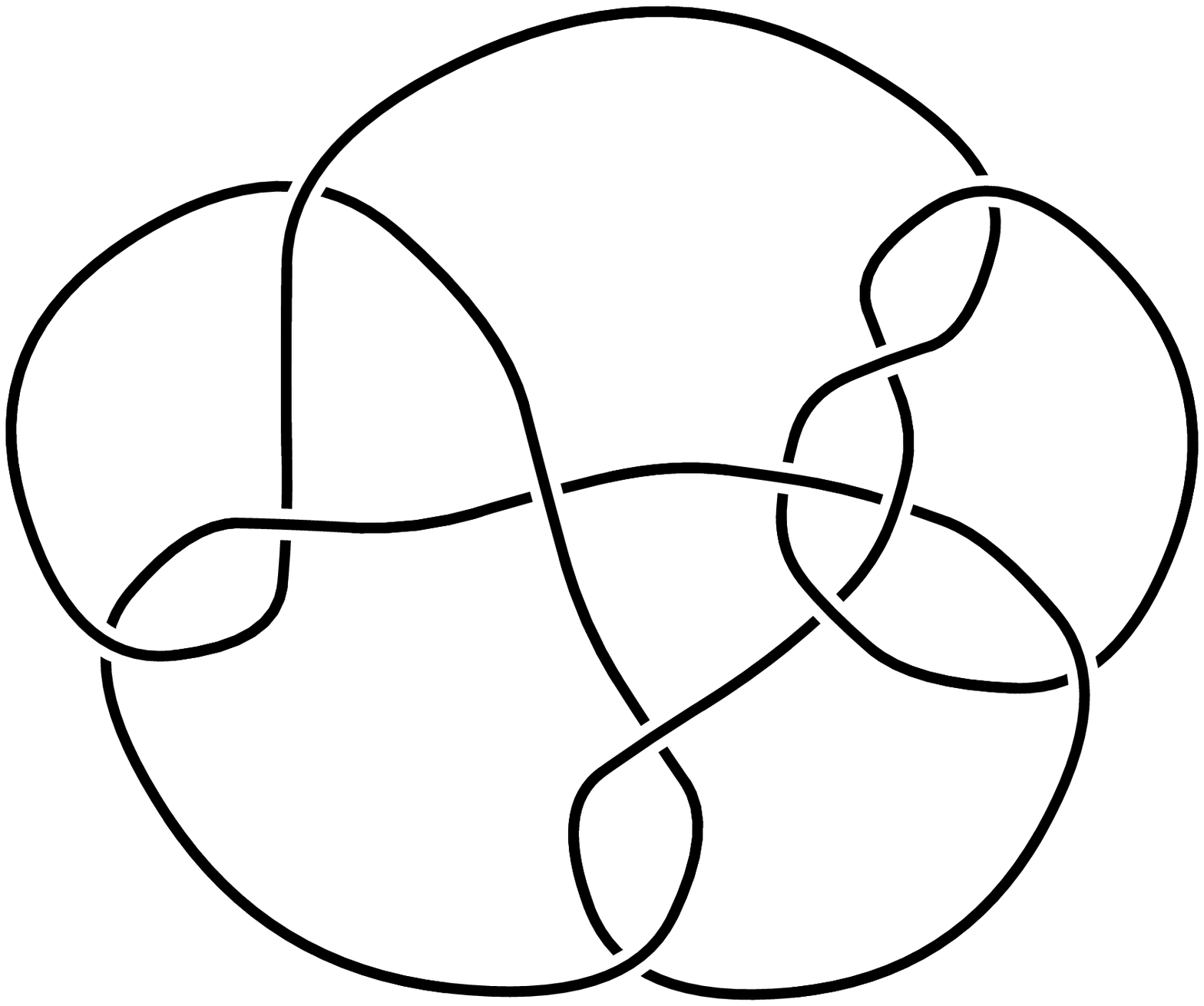}
    \\[-10pt]
    $12^A_{7}$ & $12^A_{14}$ & $12^A_{13}$ & $12^A_{15}$
    \\[10pt]
    \hline
    &&&\\[-10pt]
    \includegraphics[width=75pt]{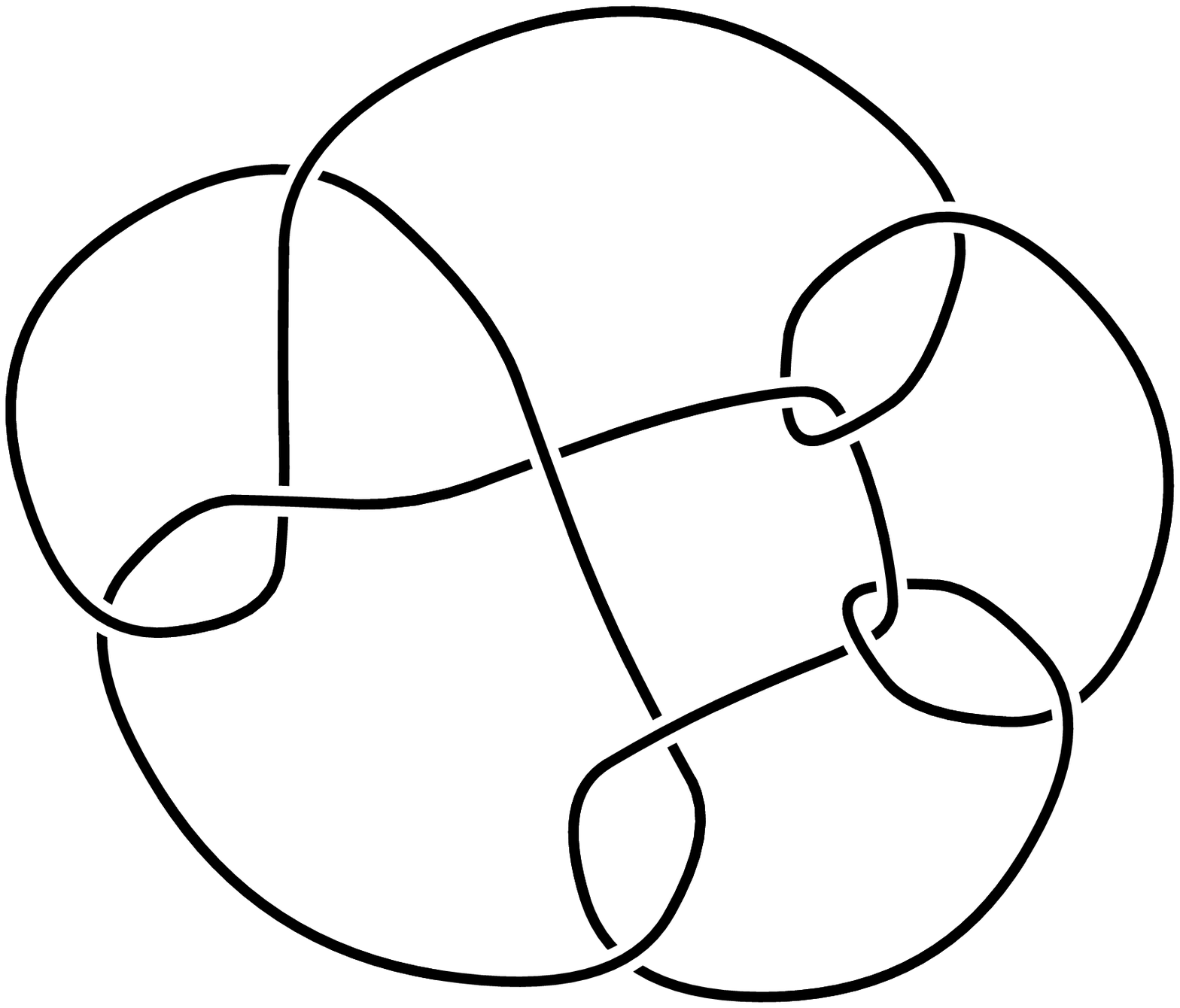}
    &
    \includegraphics[width=75pt]{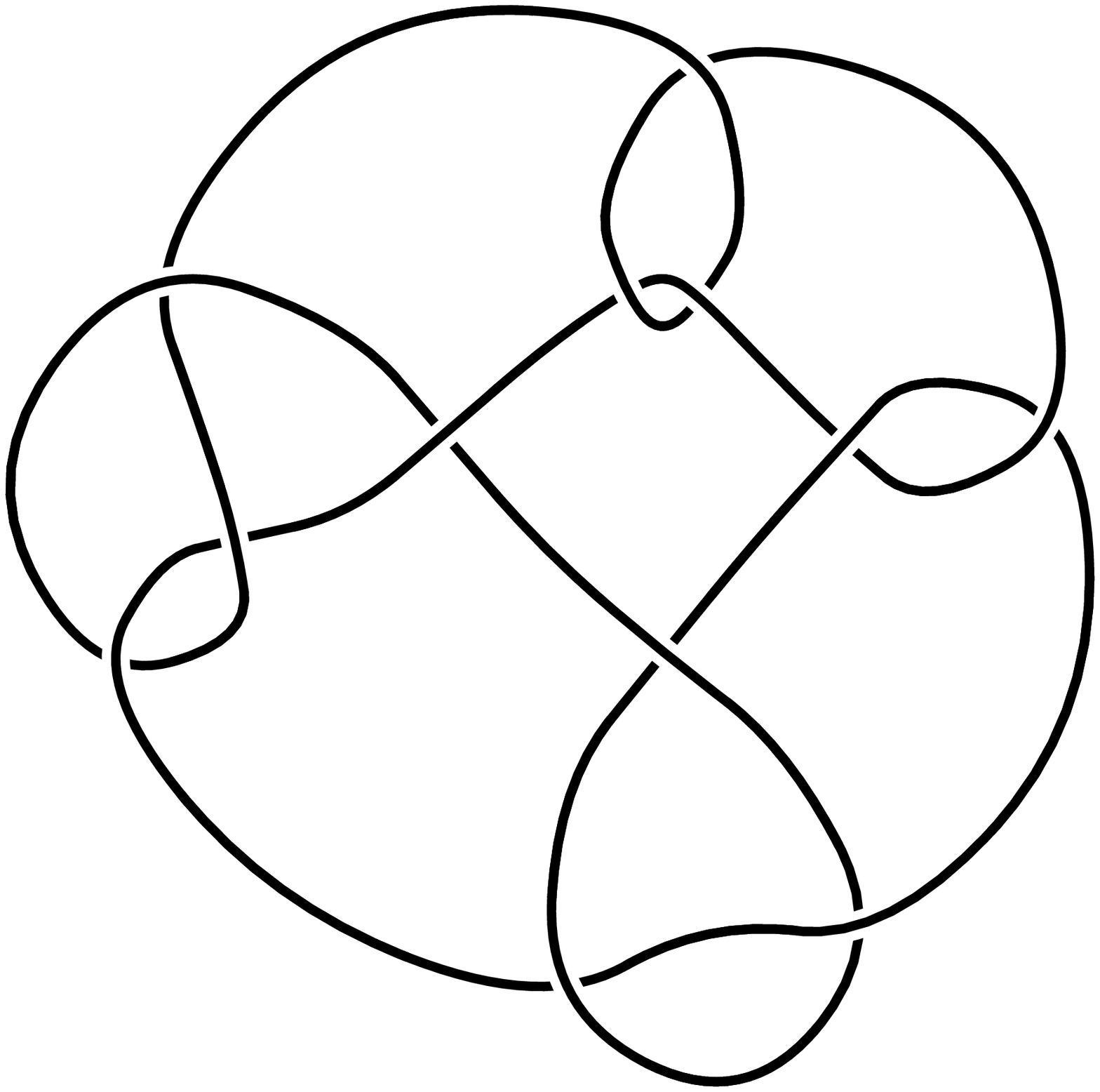}
    &
    \includegraphics[width=75pt]{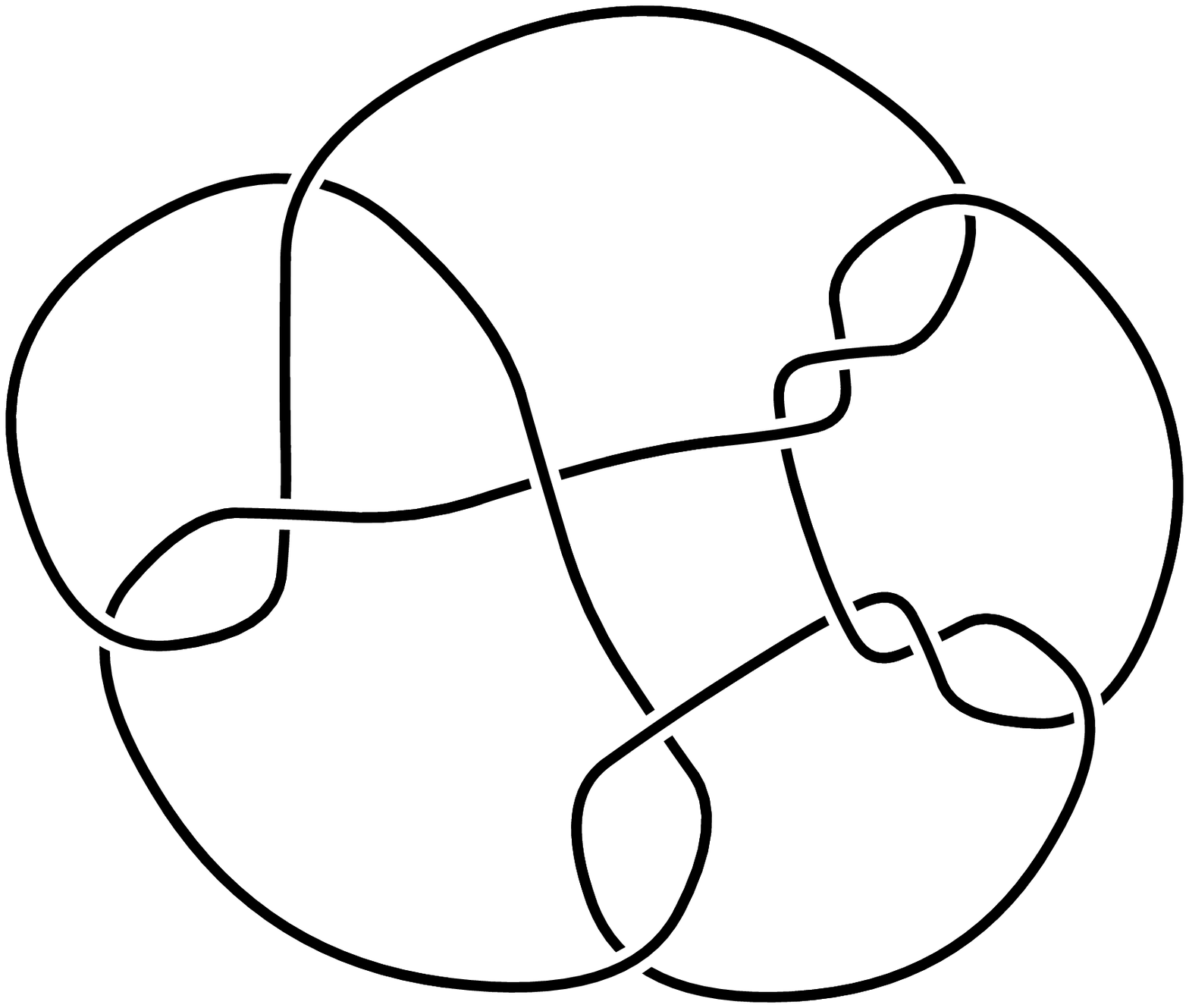}
    &
    \includegraphics[width=75pt]{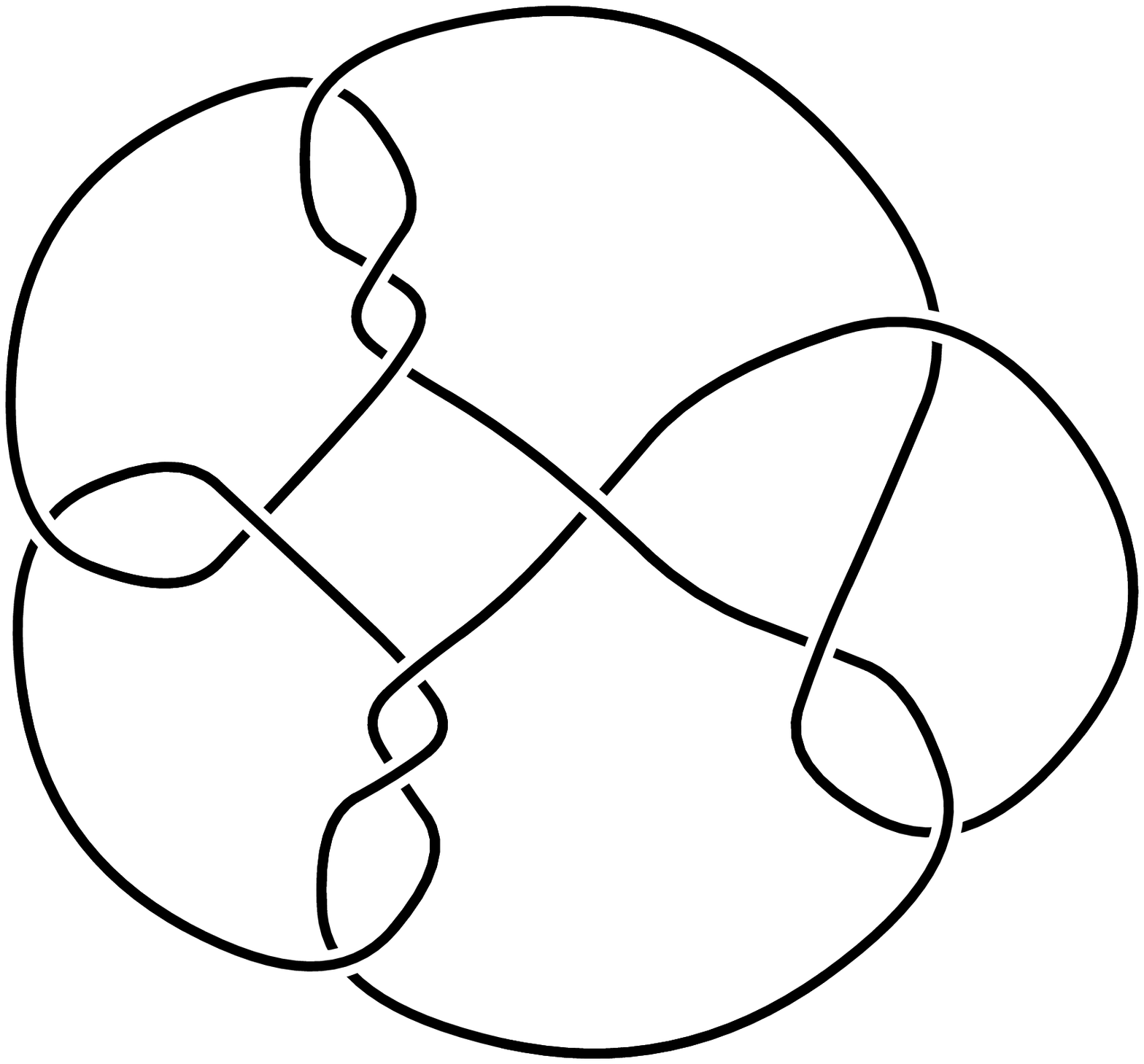}
    \\[-10pt]
    $12^A_{29}$ & $12^A_{113}$ & $12^A_{36}$ & $12^A_{694}$
    \\[10pt]
    \hline
    &&&\\[-10pt]
    \includegraphics[width=75pt]{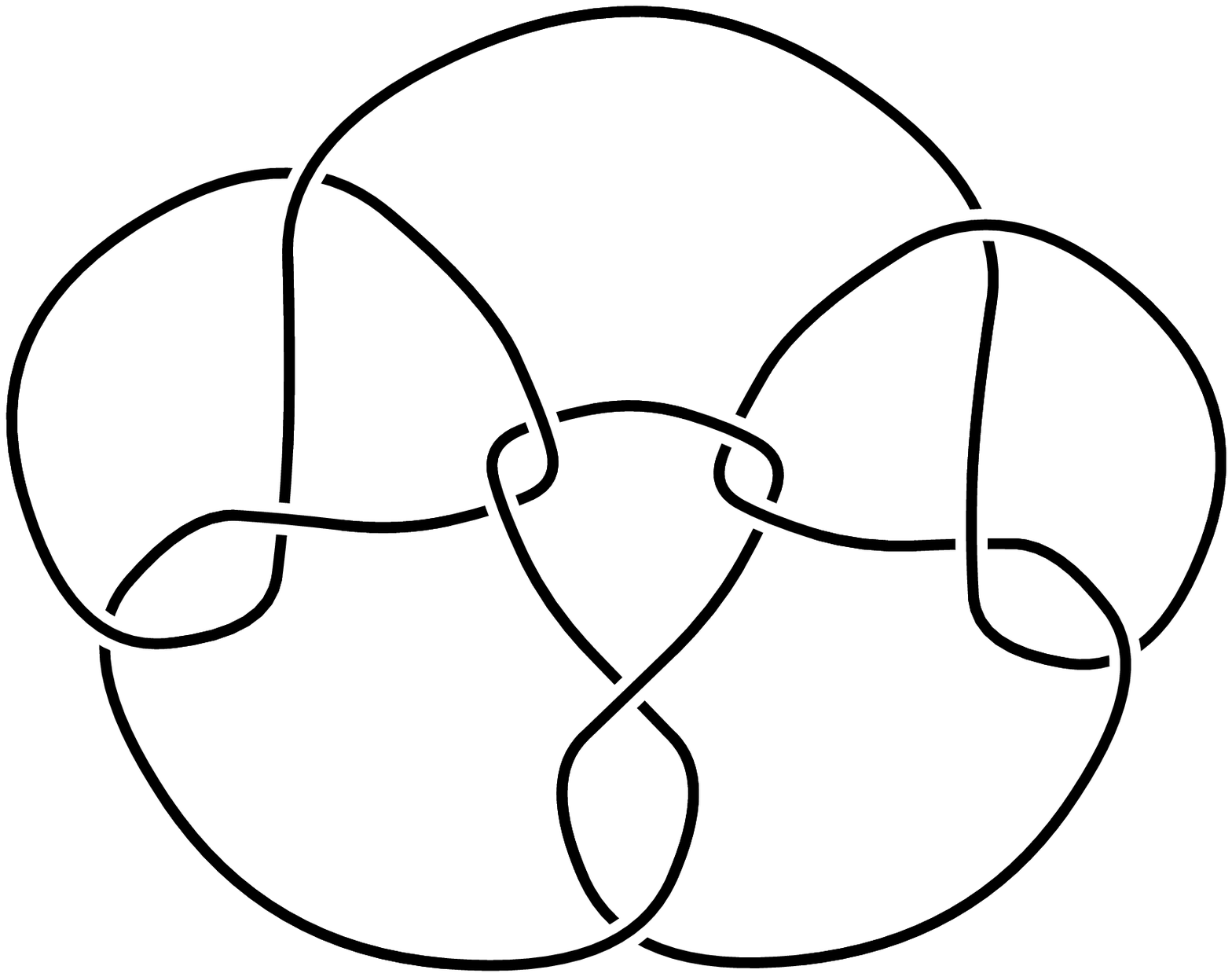}
    &
    \includegraphics[width=75pt]{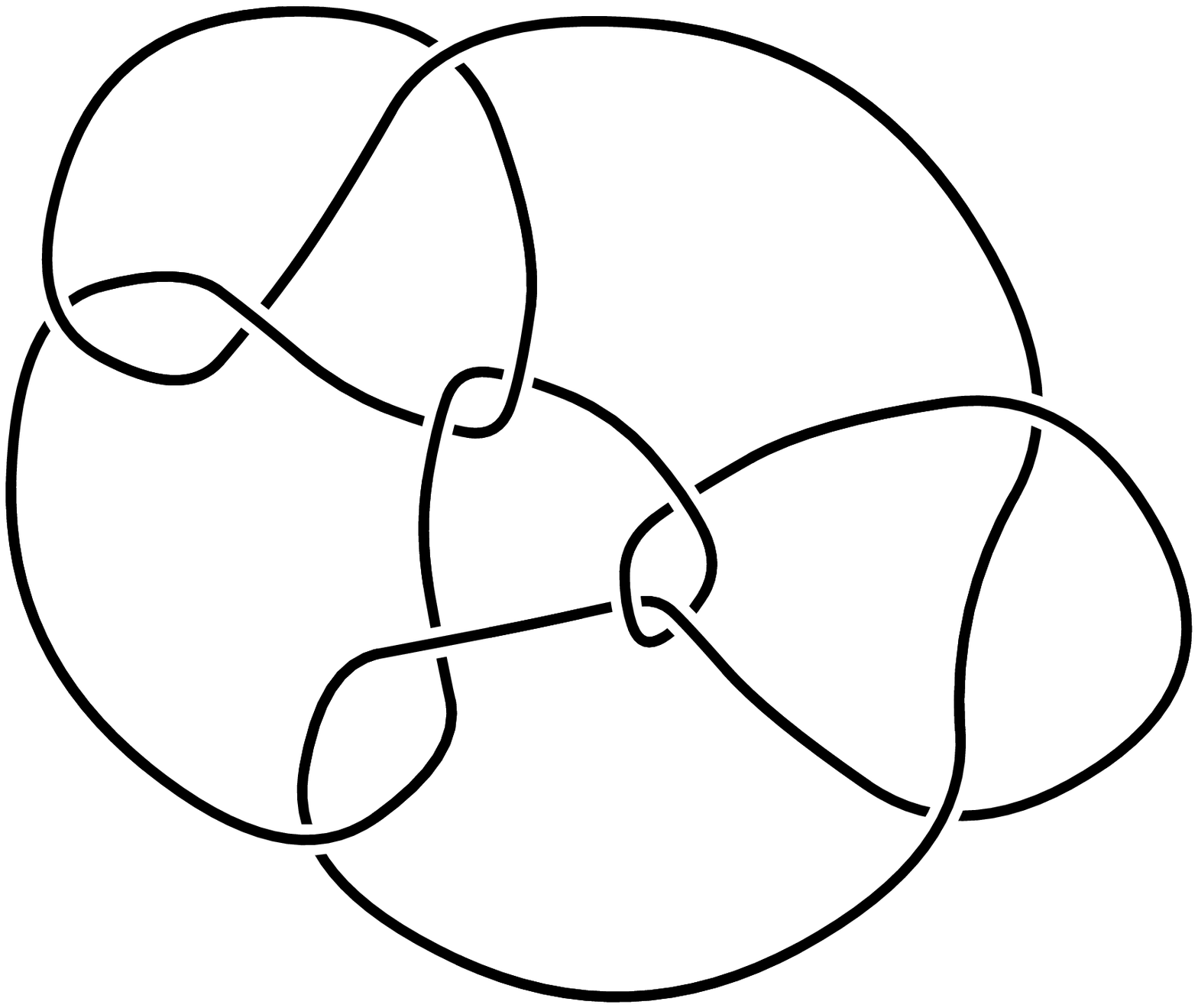}
    &
    \includegraphics[width=75pt]{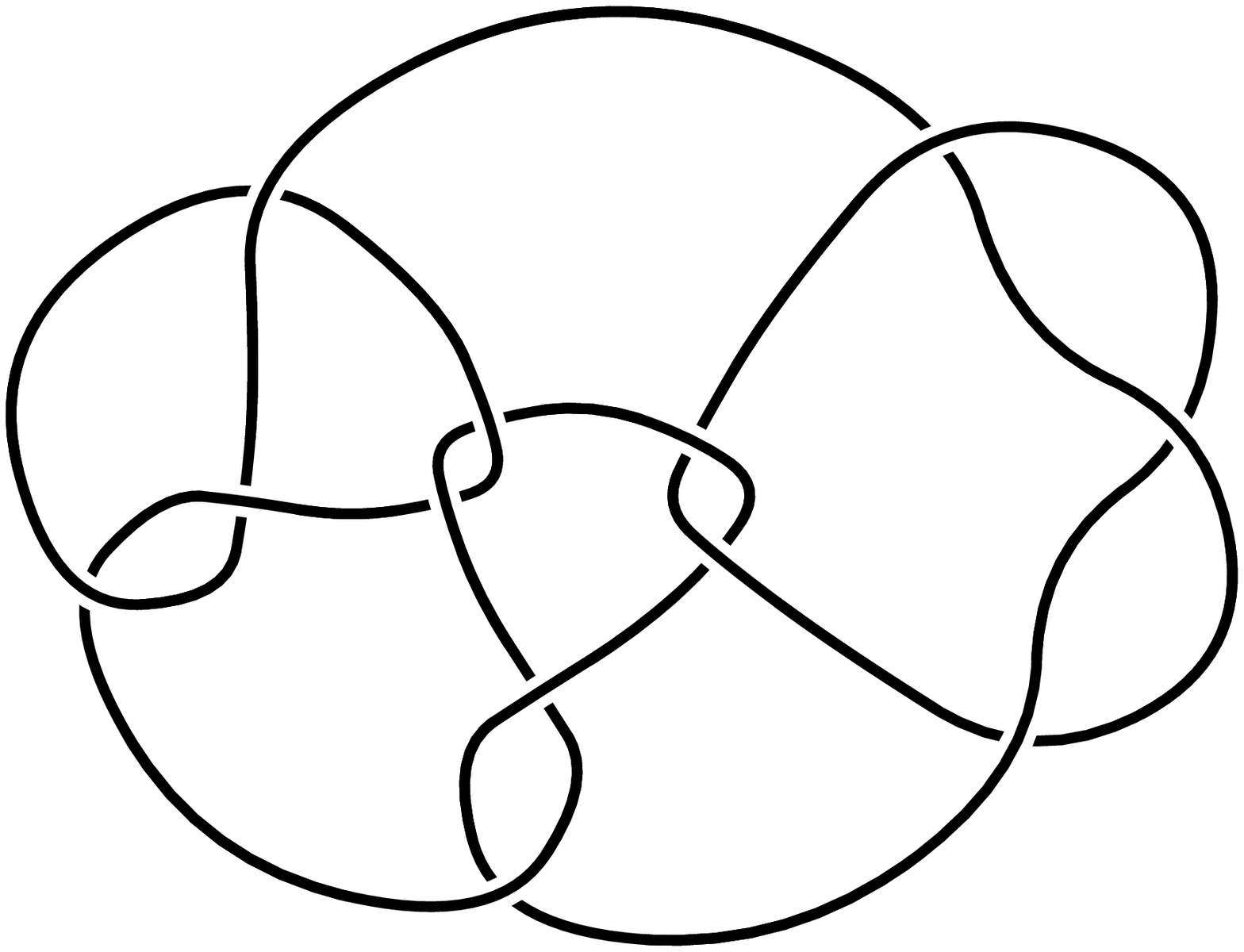}
    &
    \includegraphics[width=75pt]{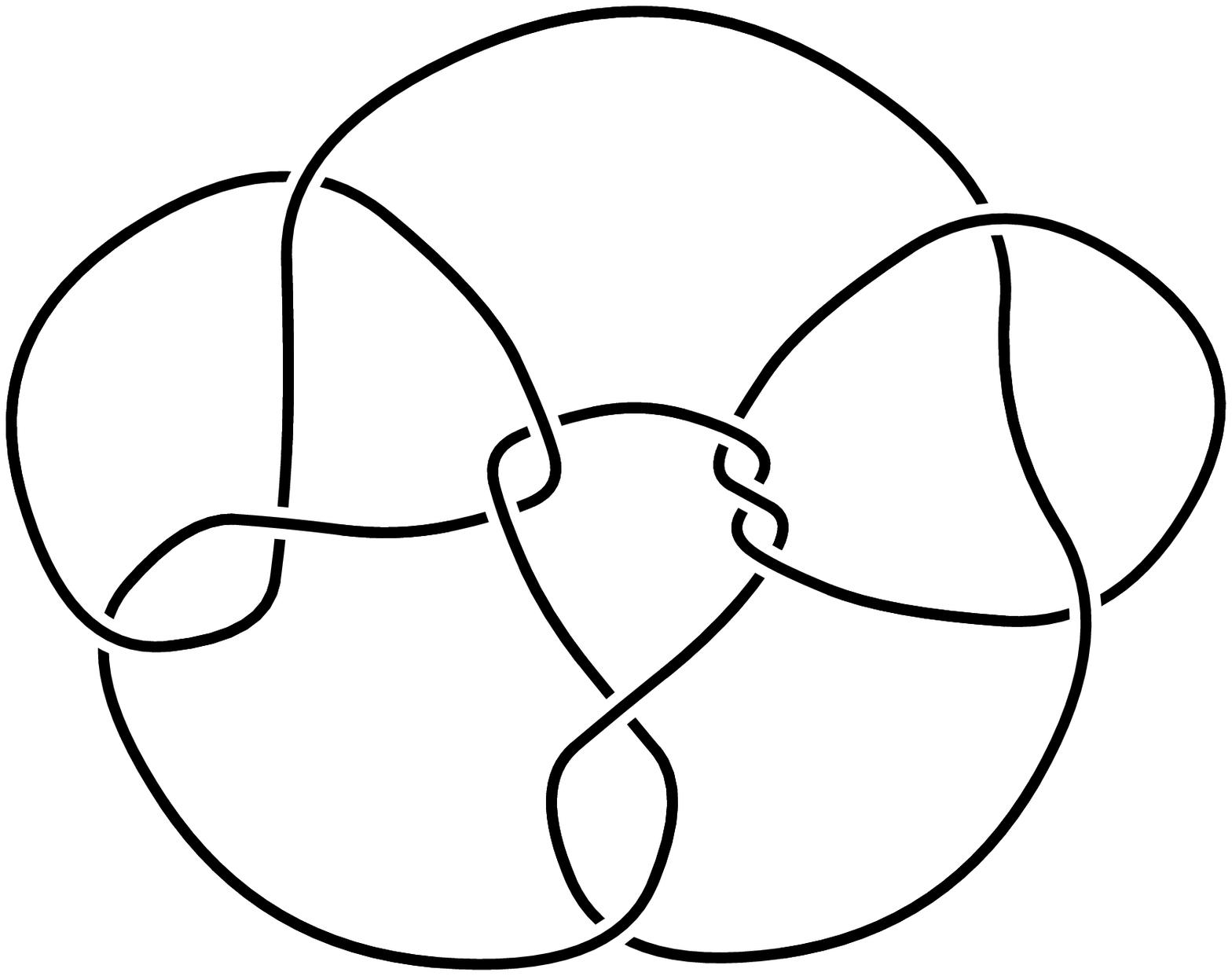}
    \\[-10pt]
    $12^A_{44}$ & $12^A_{64}$ & $12^A_{45}$ & $12^A_{65}$
    \\[10pt]
    \hline
    &&&\\[-10pt]
    \includegraphics[width=75pt]{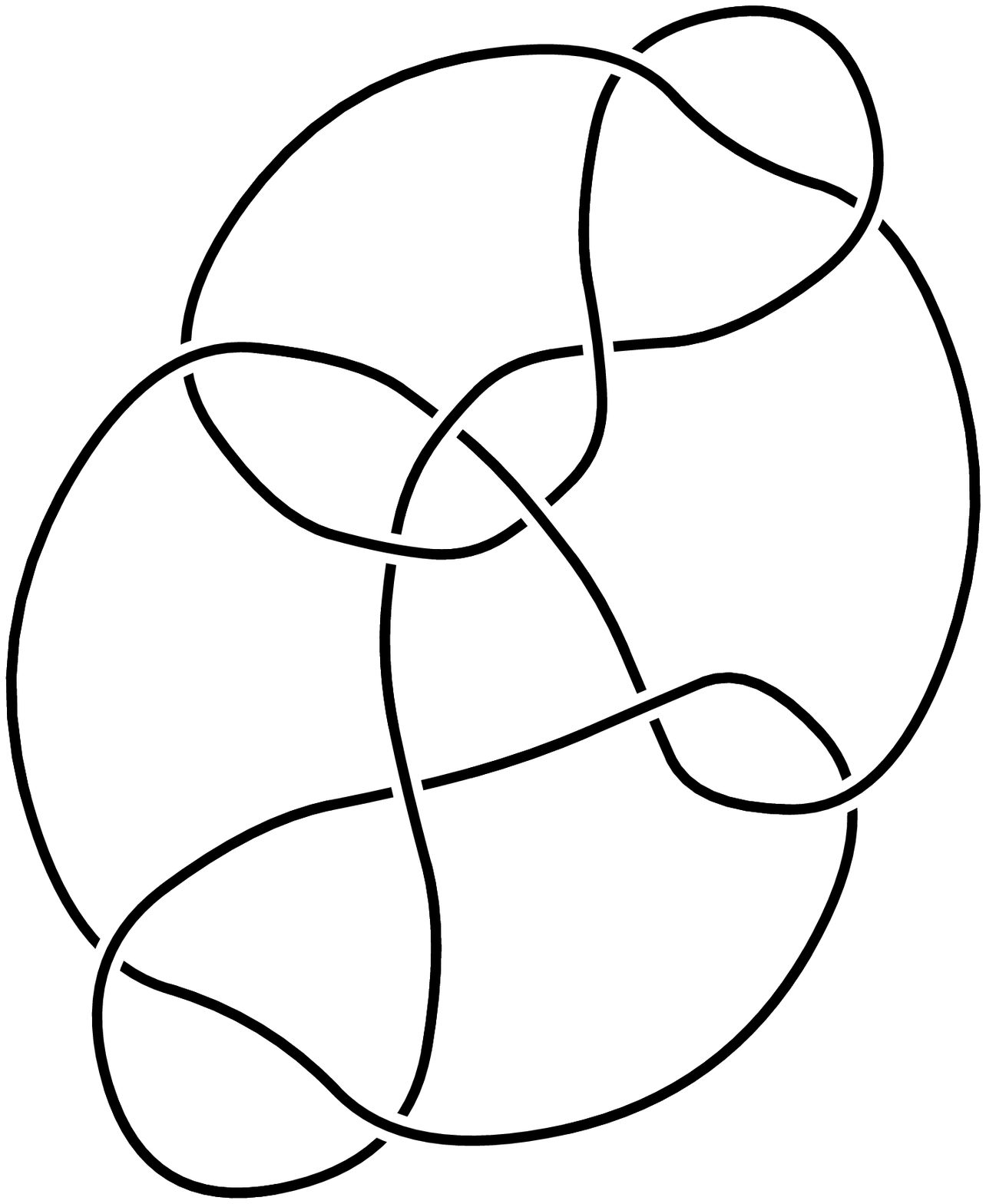}
    &
    \includegraphics[width=75pt]{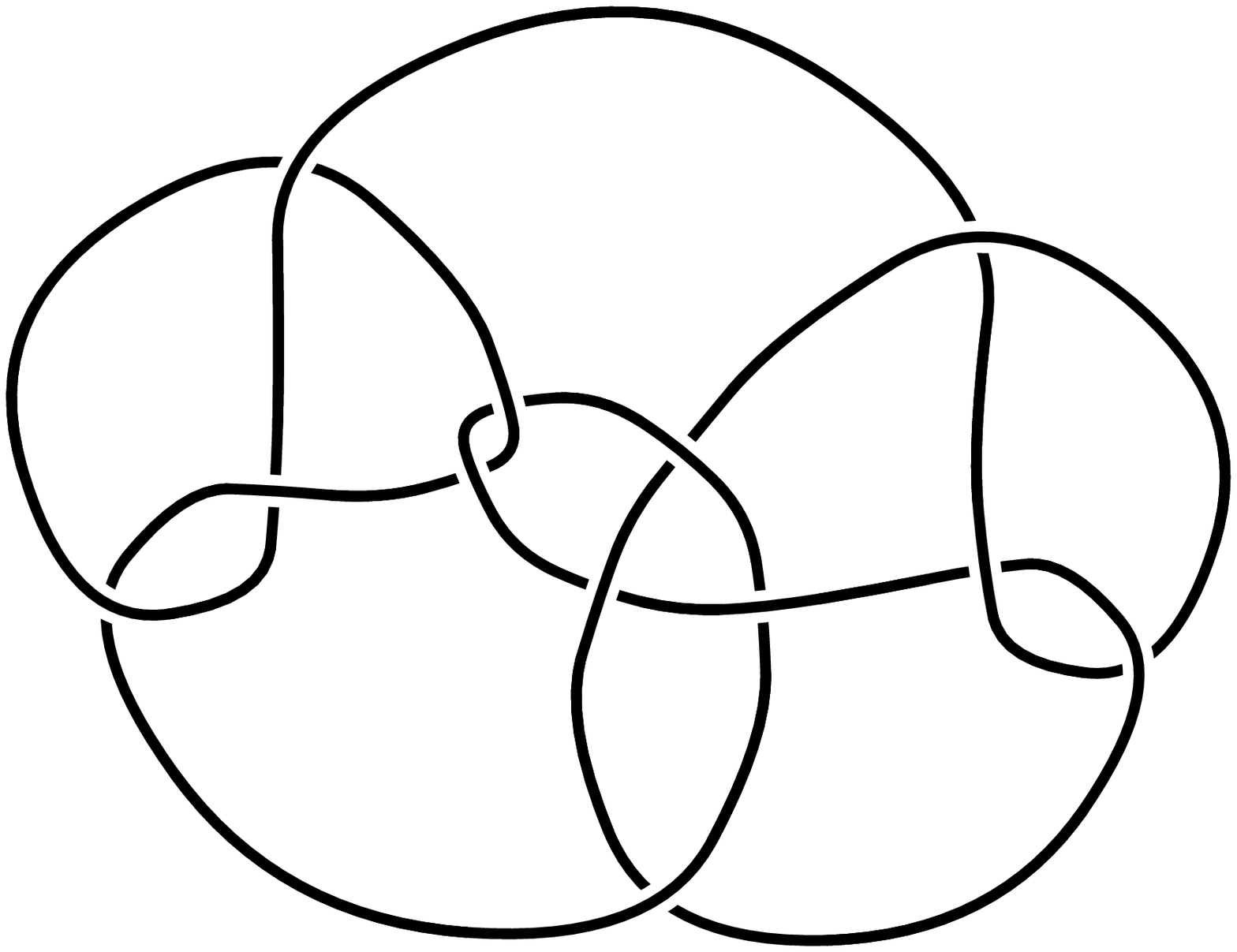}
    &
    \includegraphics[width=75pt]{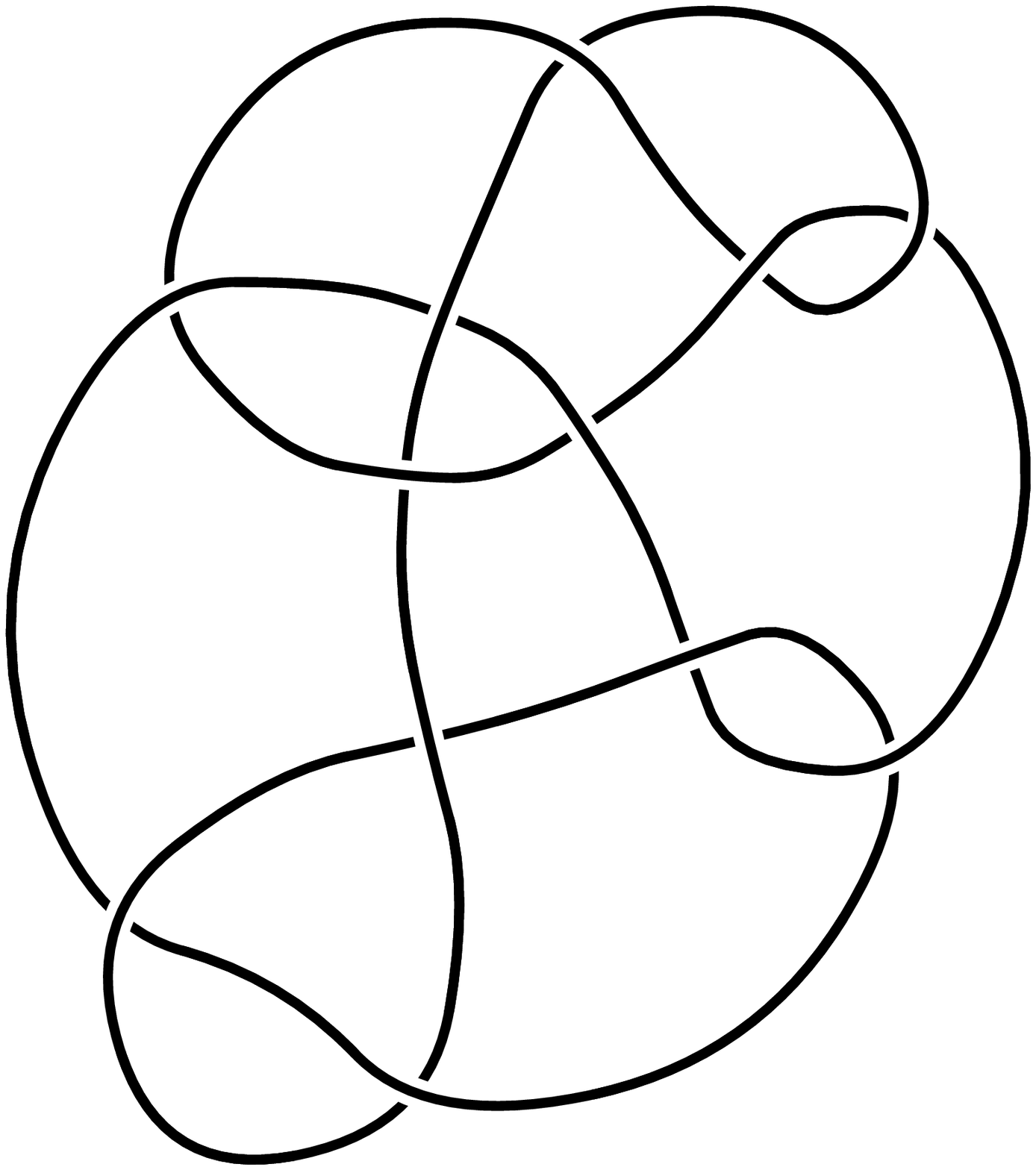}
    &
    \includegraphics[width=75pt]{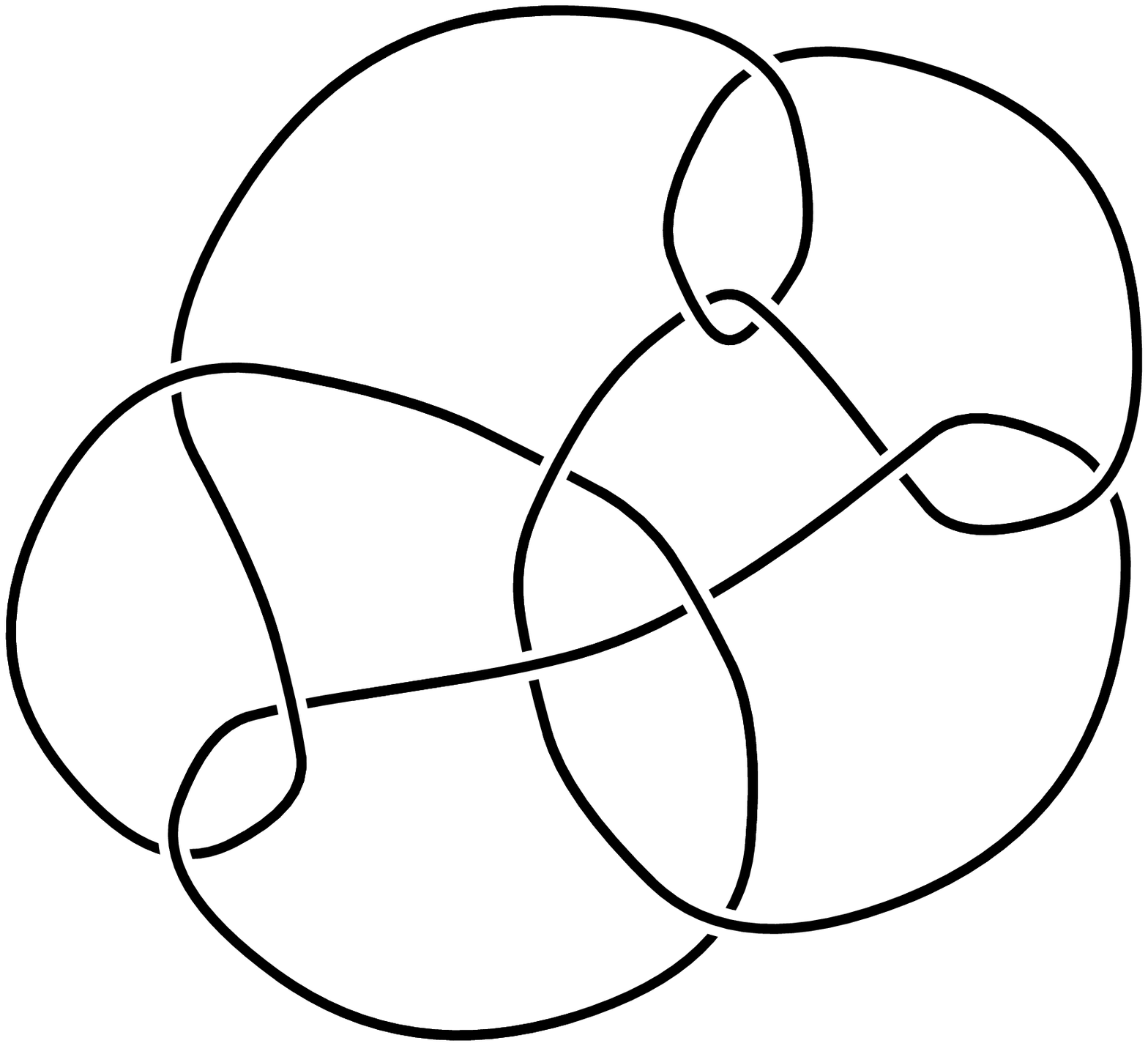}
    \\[-10pt]
    $12^A_{48}$ & $12^A_{60}$ & $12^A_{59}$ & $12^A_{63}$
    \\[10pt]
    \hline
    &&&\\[-10pt]
    \includegraphics[width=75pt]{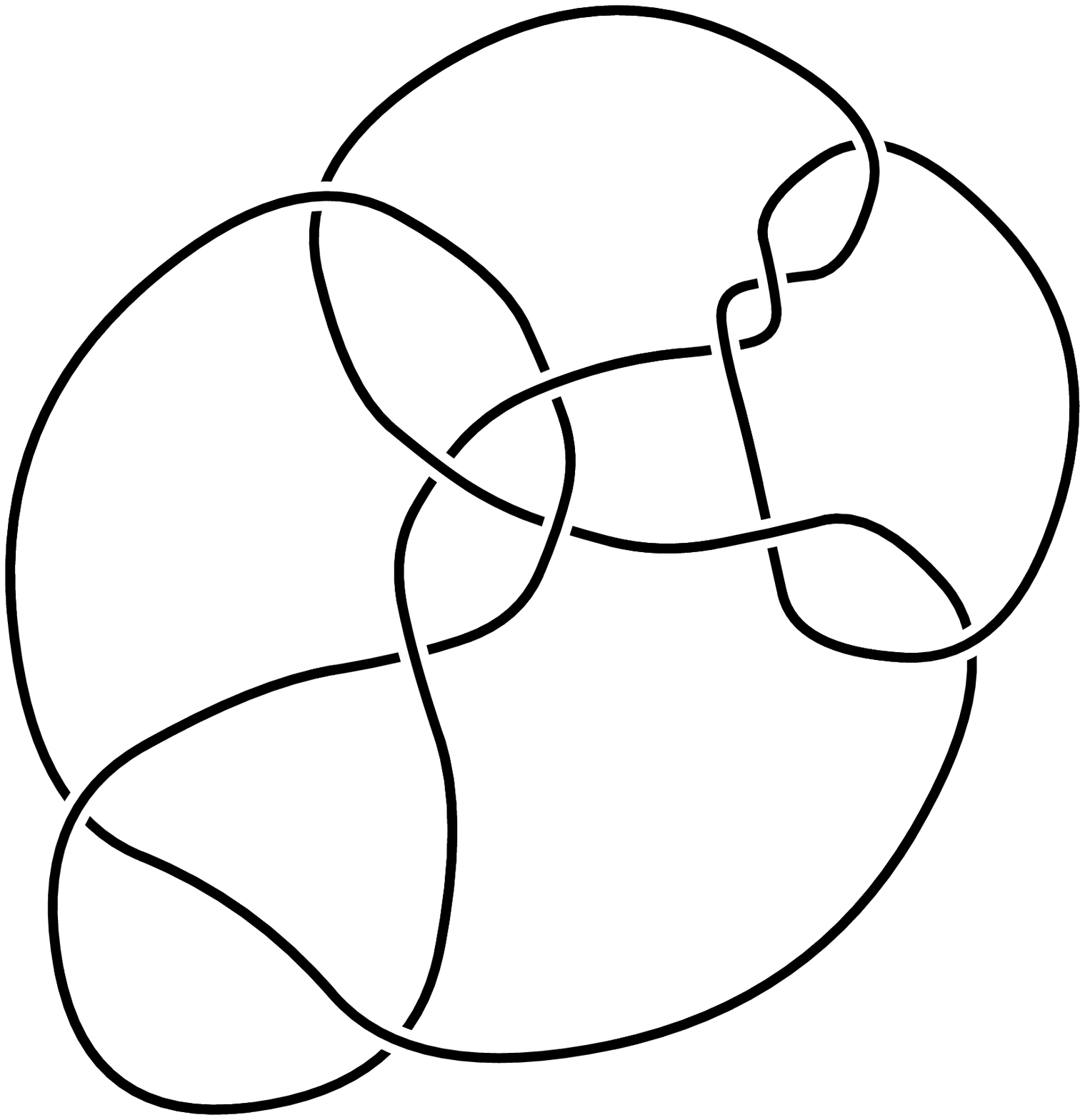}
    &
    \includegraphics[width=75pt]{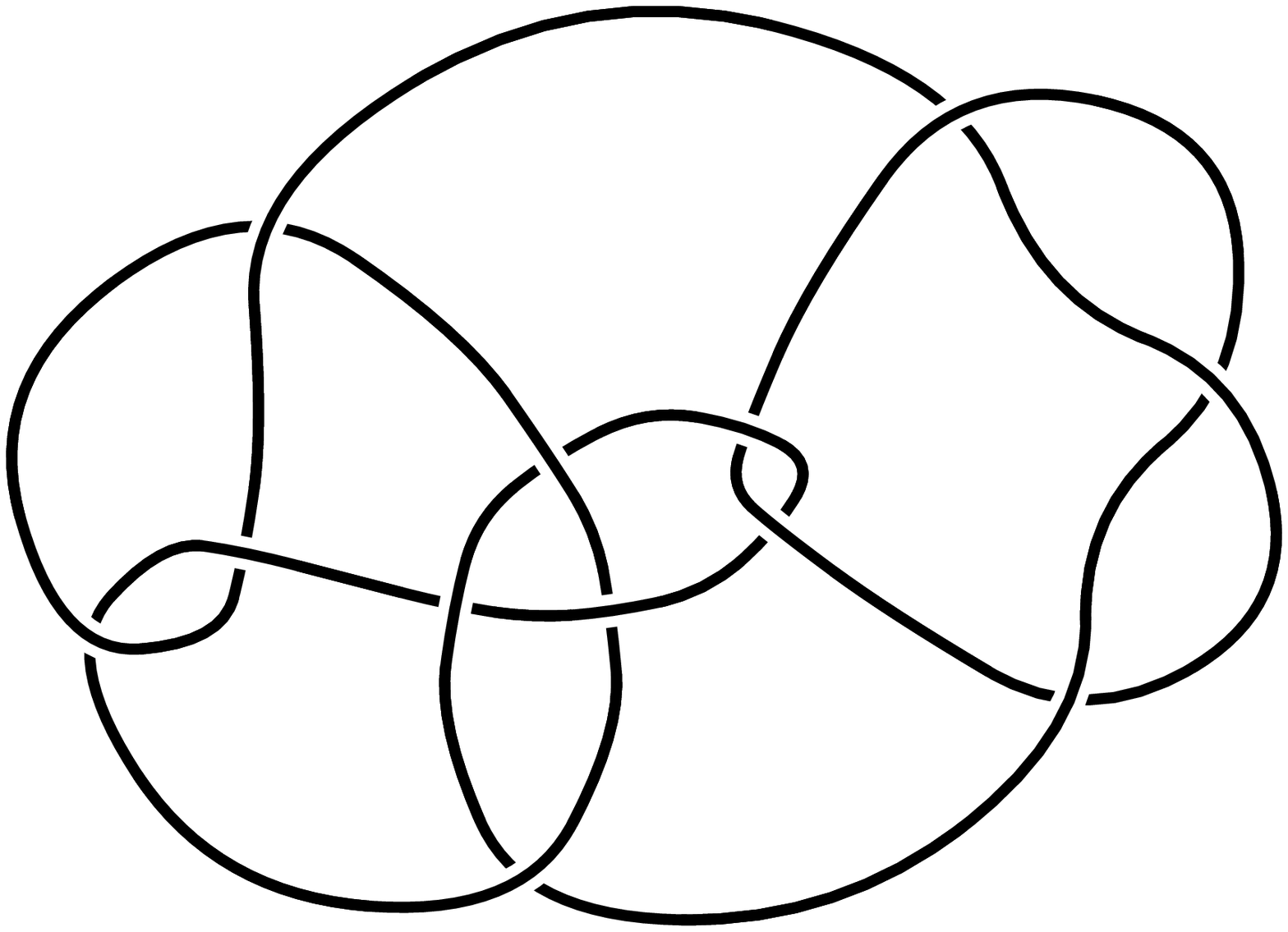}
    &
    \includegraphics[width=75pt]{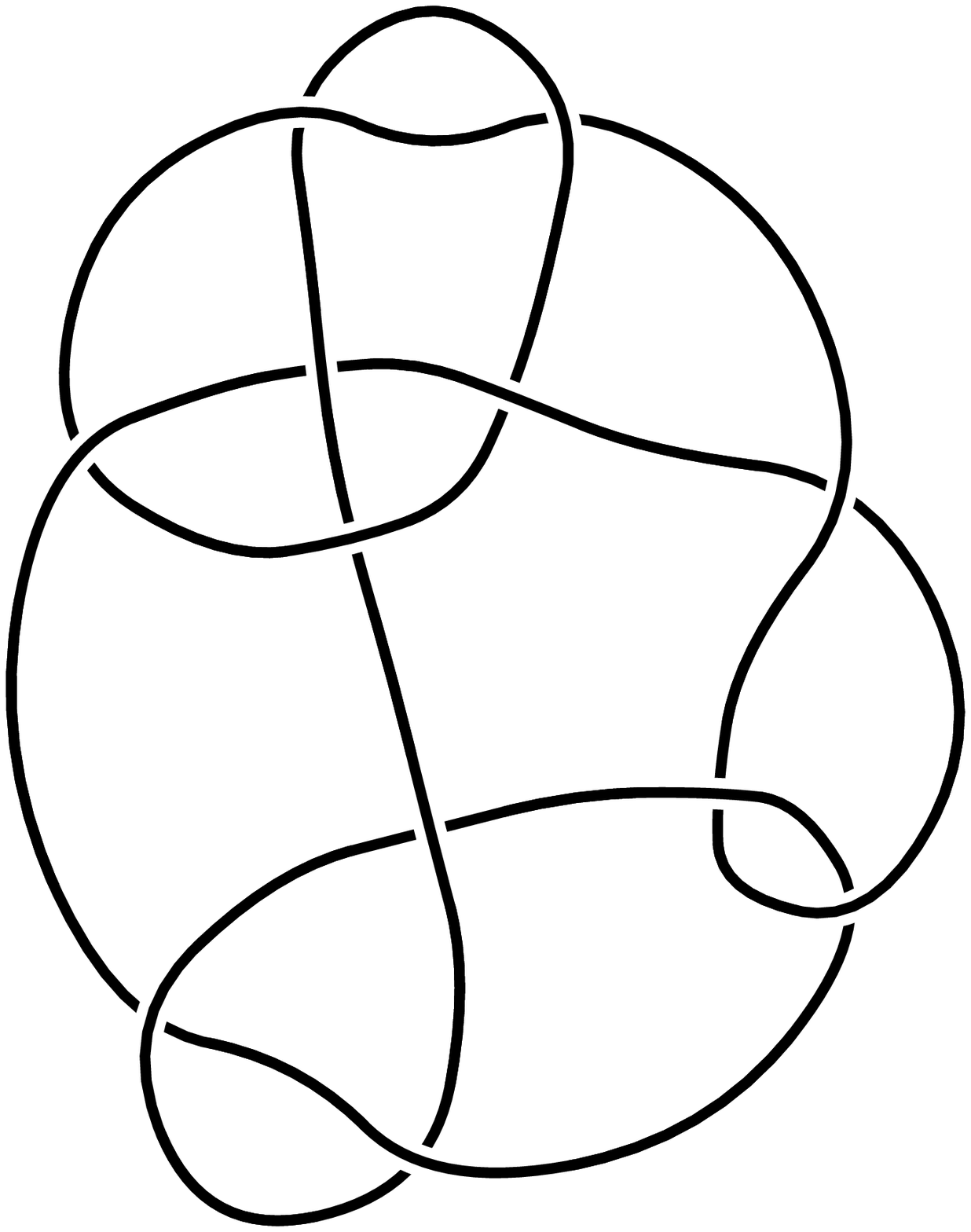}
    &
    \includegraphics[width=75pt]{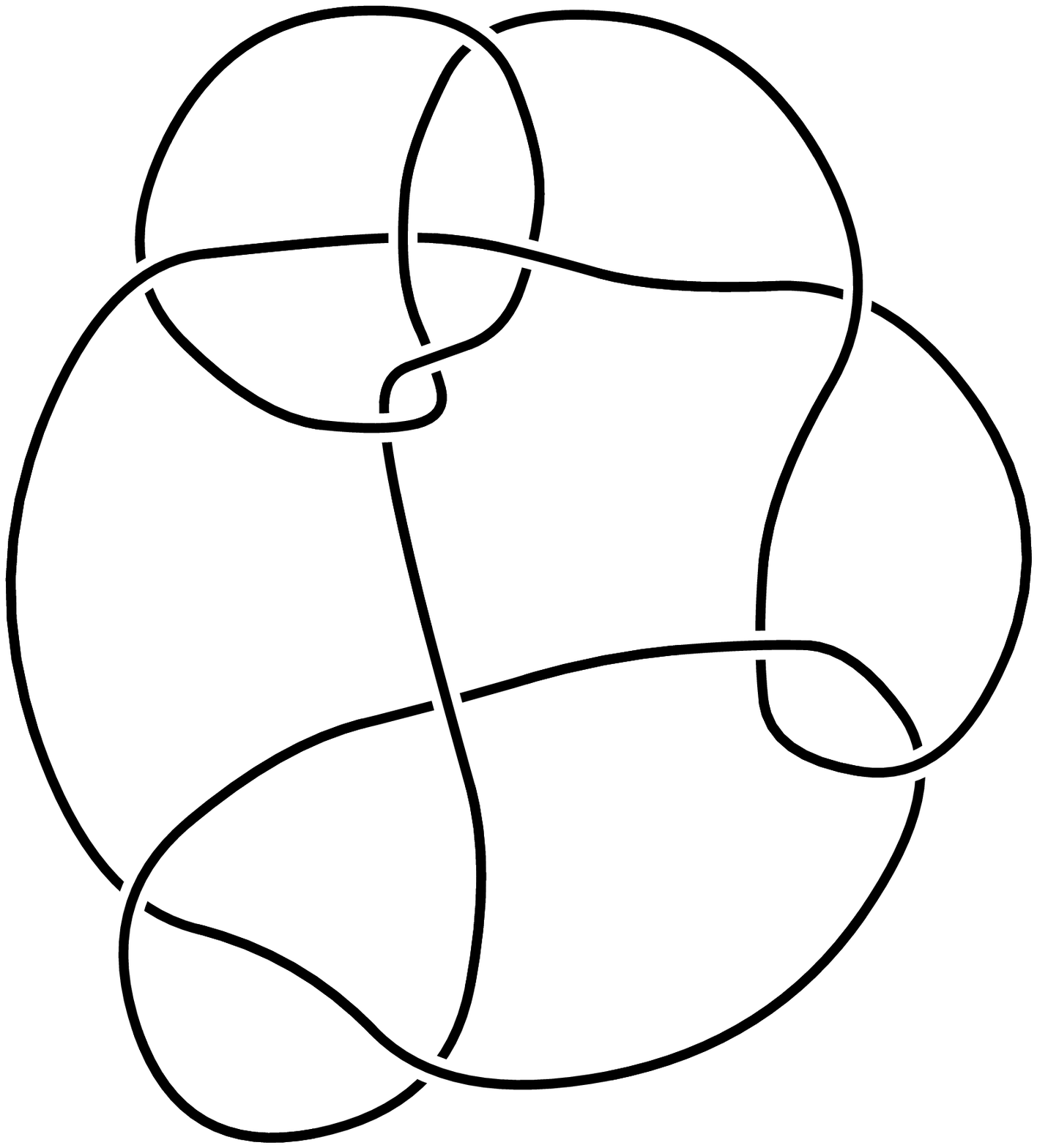}
    \\[-10pt]
    $12^A_{67}$ & $12^A_{136}$ & $12^A_{91}$ & $12^A_{111}$
  \end{tabular}
  \caption{Alternating $12$-crossing mutant cliques 1/4}
  \label{figure:Alternating12crossingmutantcliques1of4}
  \end{centering}
\end{figure}

\begin{figure}[htbp]
  \begin{centering}
  \begin{tabular}{ccc}
    \includegraphics[width=75pt]{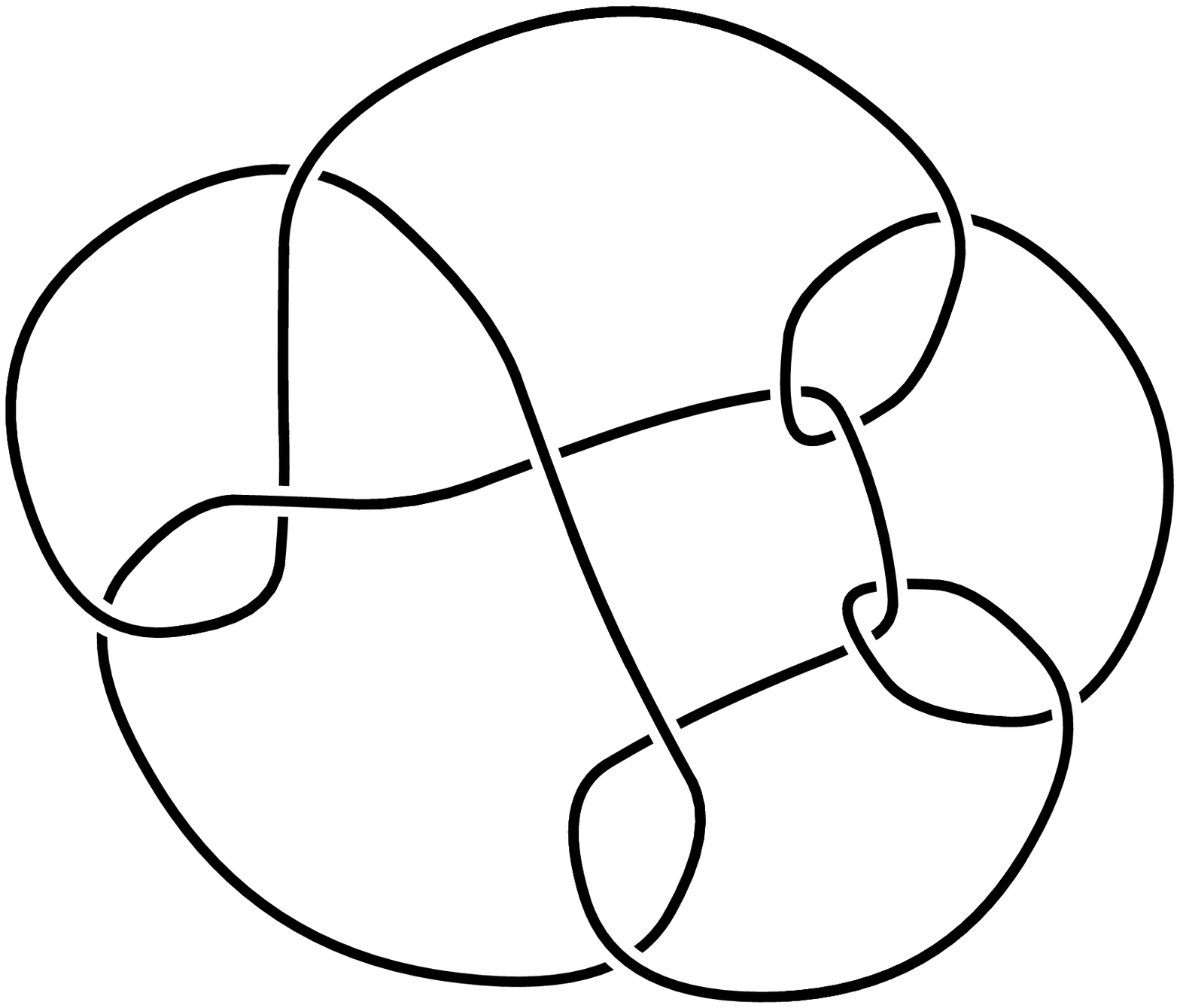}
    &
    \includegraphics[width=75pt]{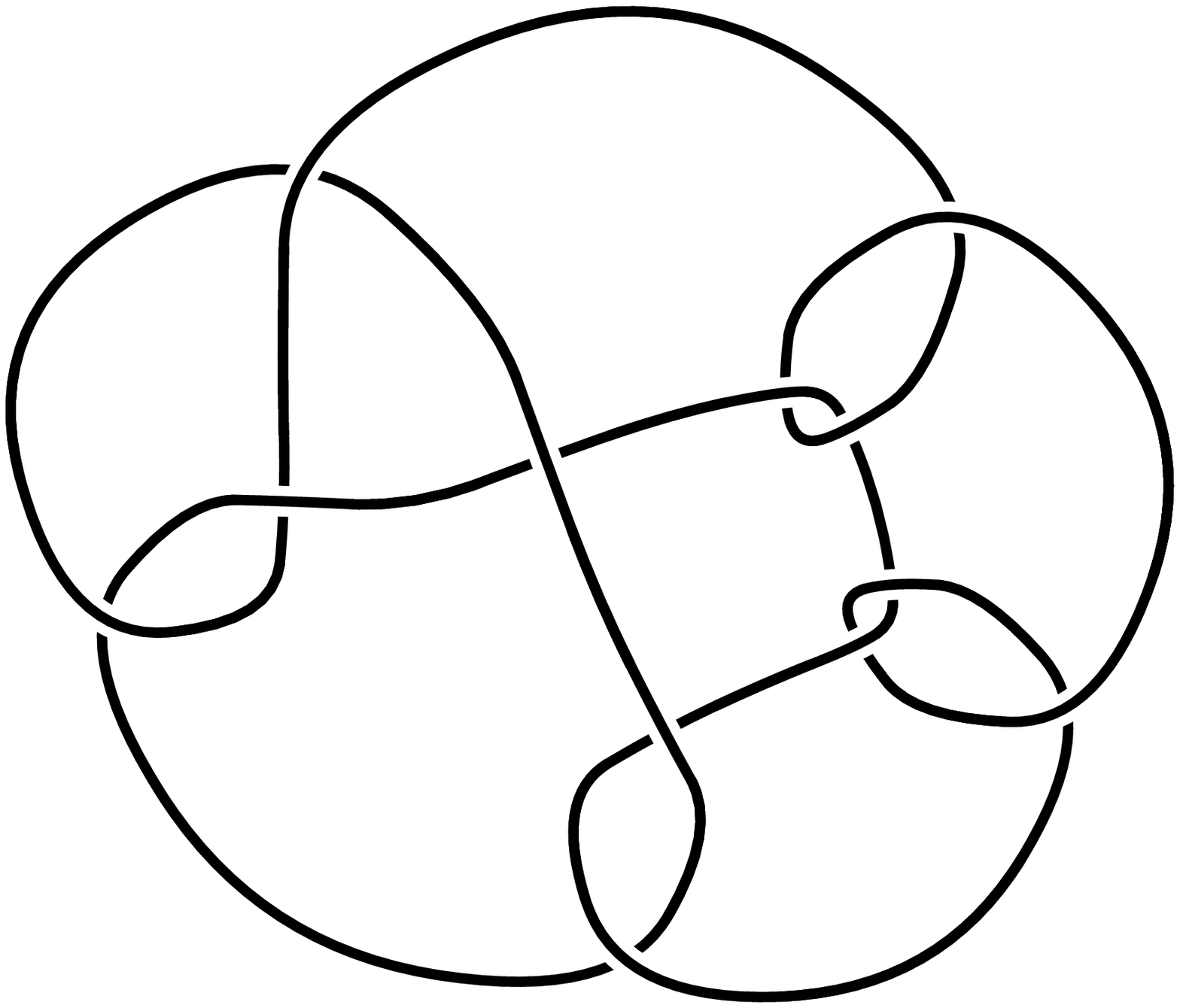}
    &
    \includegraphics[width=75pt]{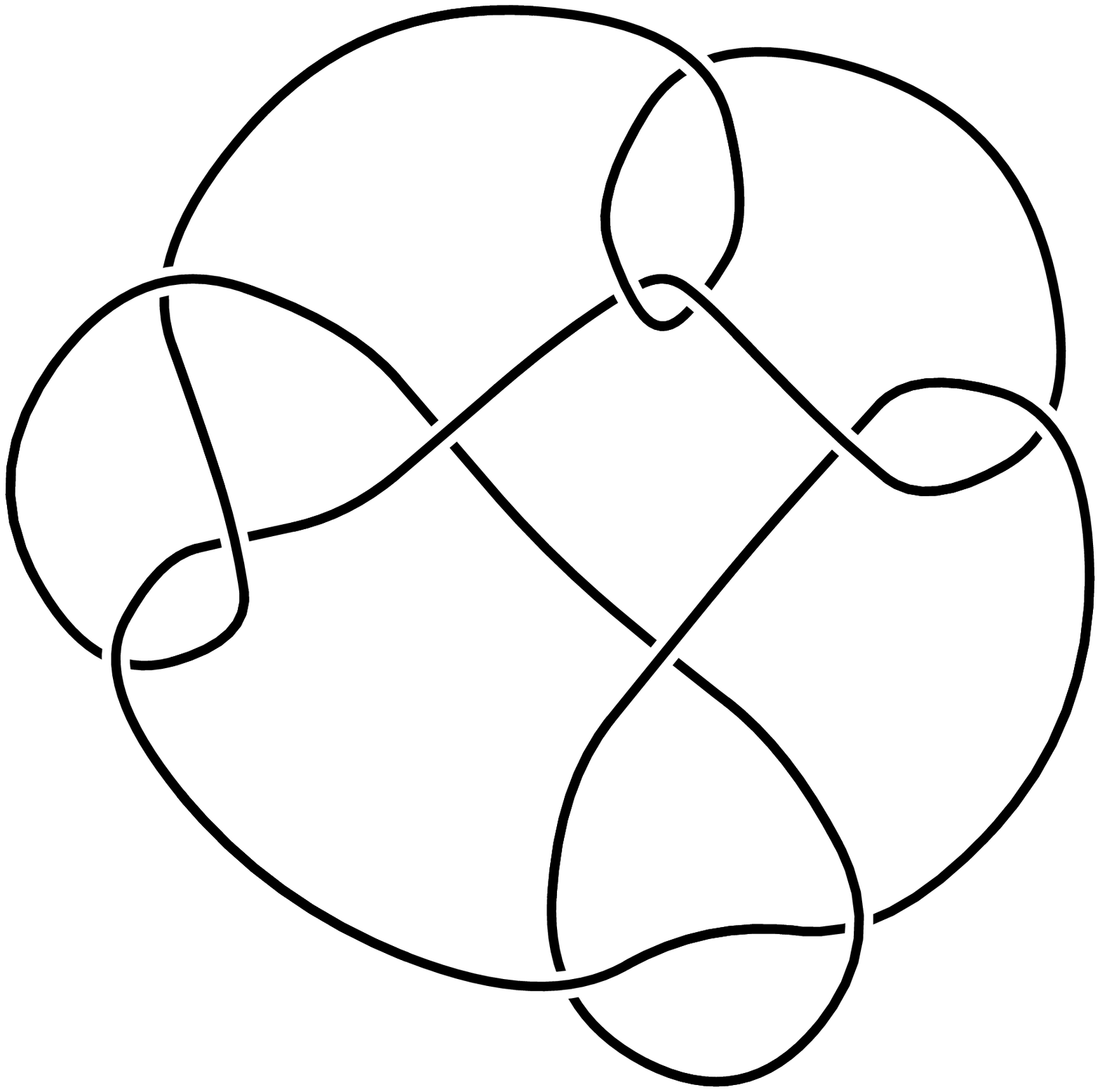}
    \\[-10pt]
    $12^N_{56}$ & $12^N_{57}$ & $12^N_{221}$
    \\[10pt]
    \hline
    &&\\[-10pt]
    \includegraphics[width=75pt]{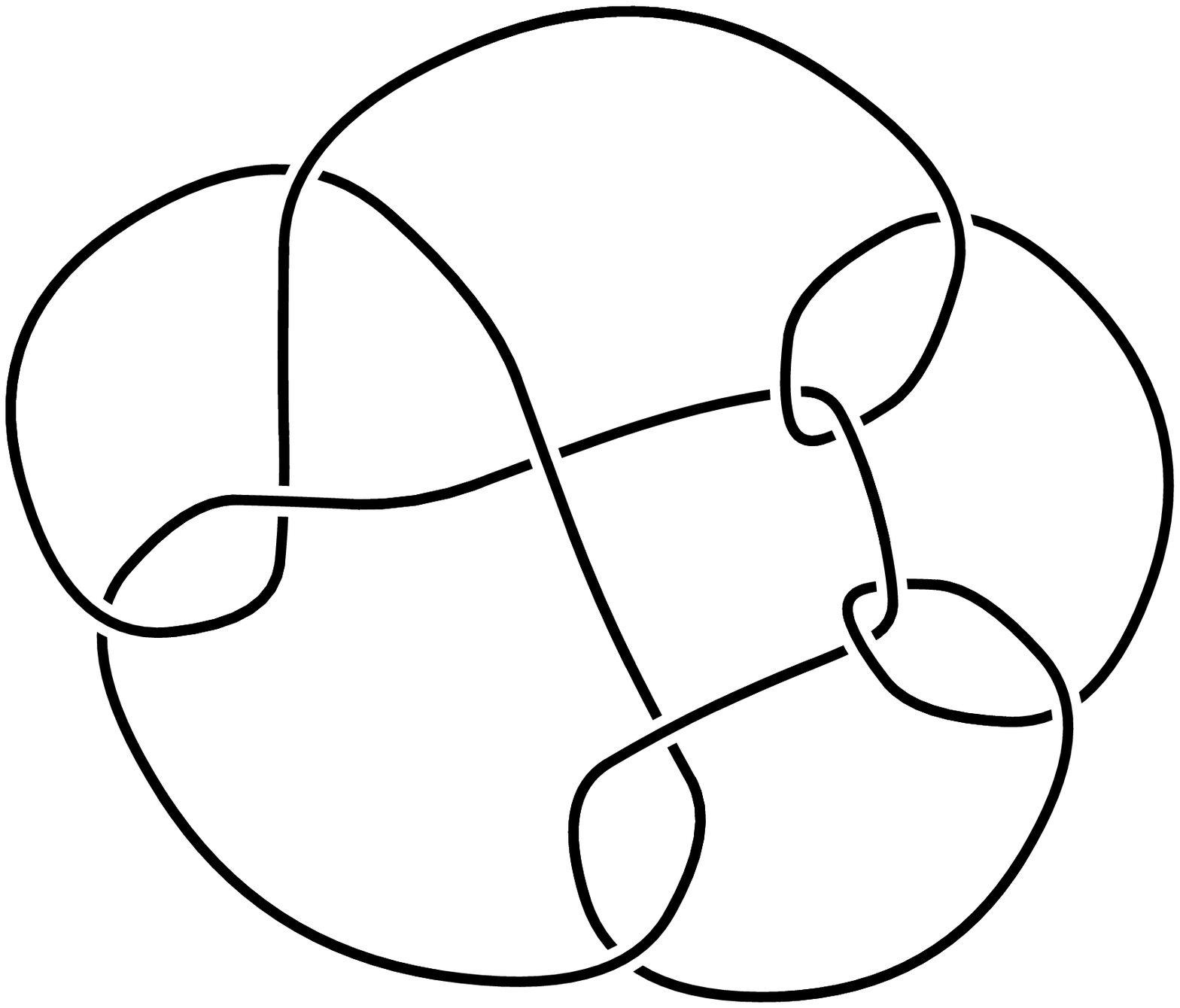}
    &
    \includegraphics[width=75pt]{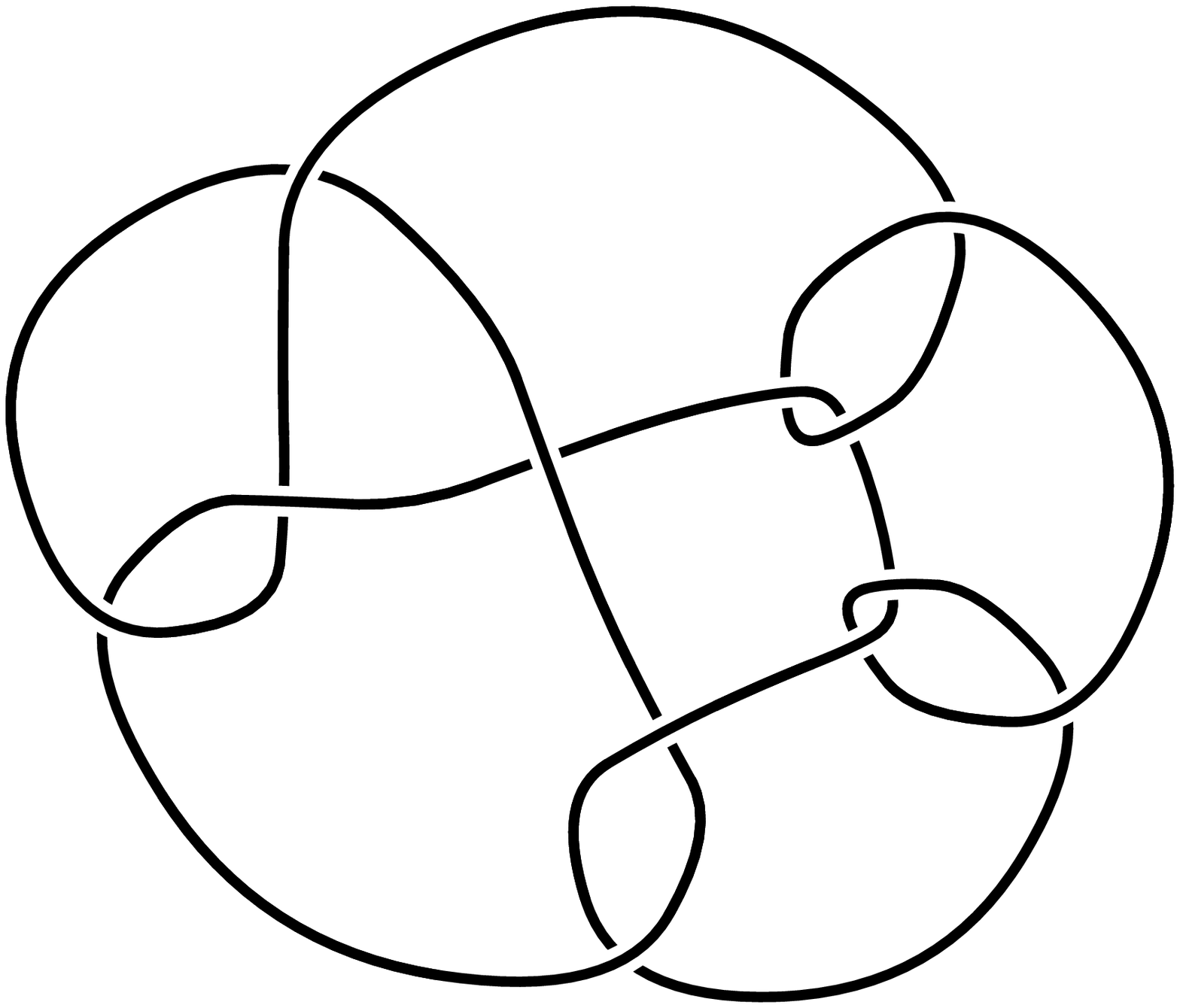}
    &
    \includegraphics[width=75pt]{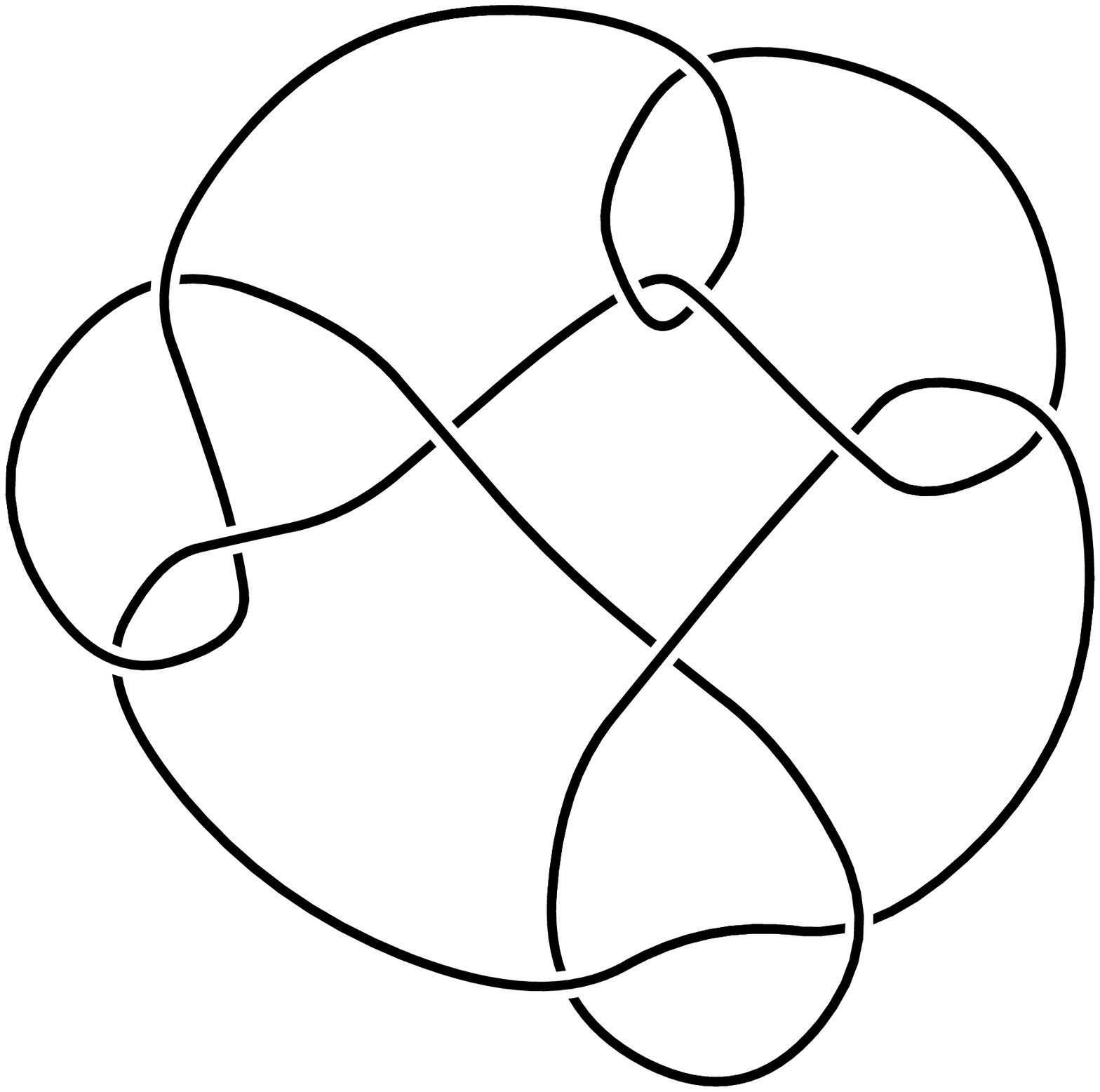}
    \\[-10pt]
    $12^N_{60}$ & $12^N_{61}$ & $12^N_{219}$
    \\[10pt]
    \hline
    &&\\[-10pt]
    \includegraphics[width=75pt]{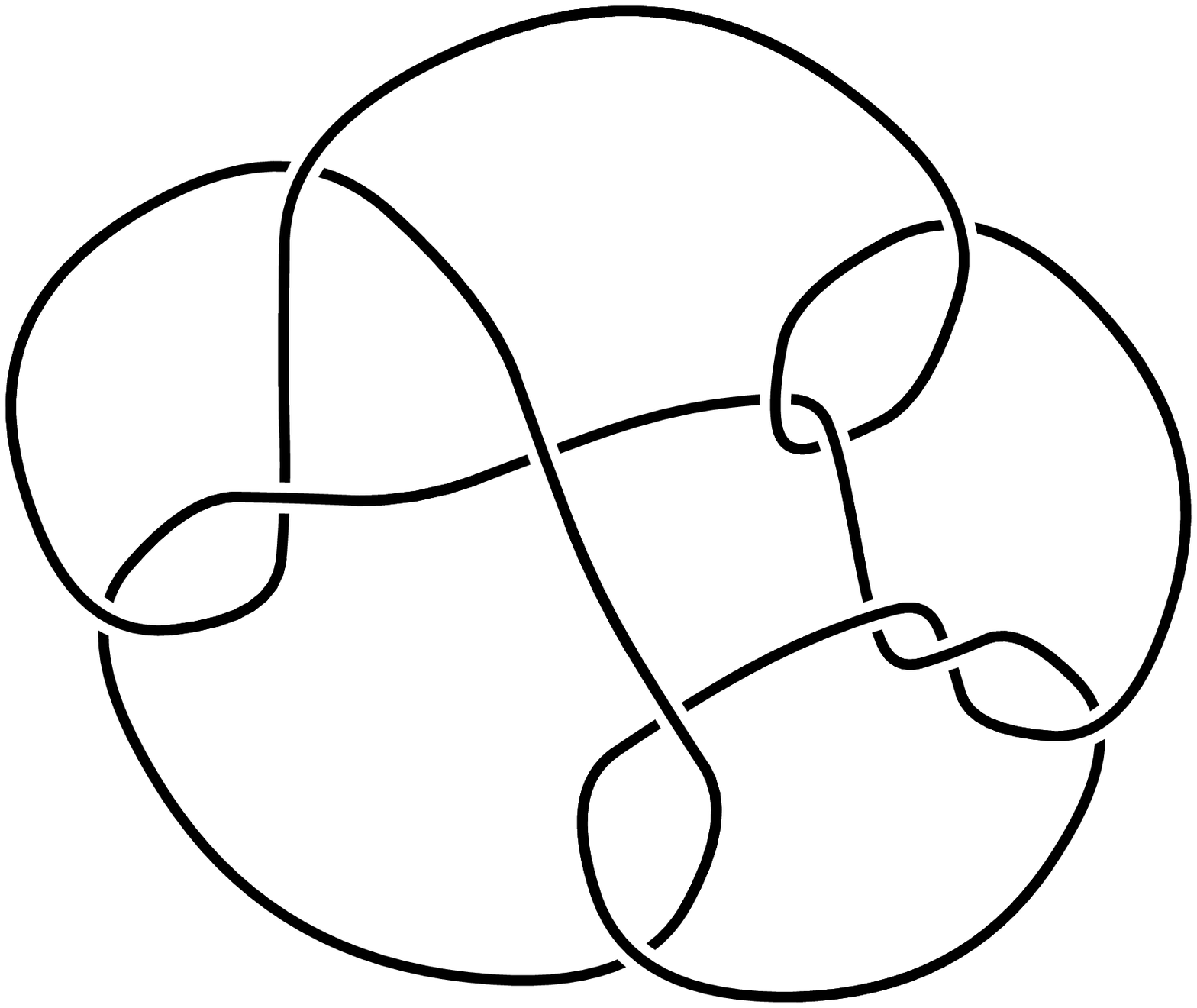}
    &
    \includegraphics[width=75pt]{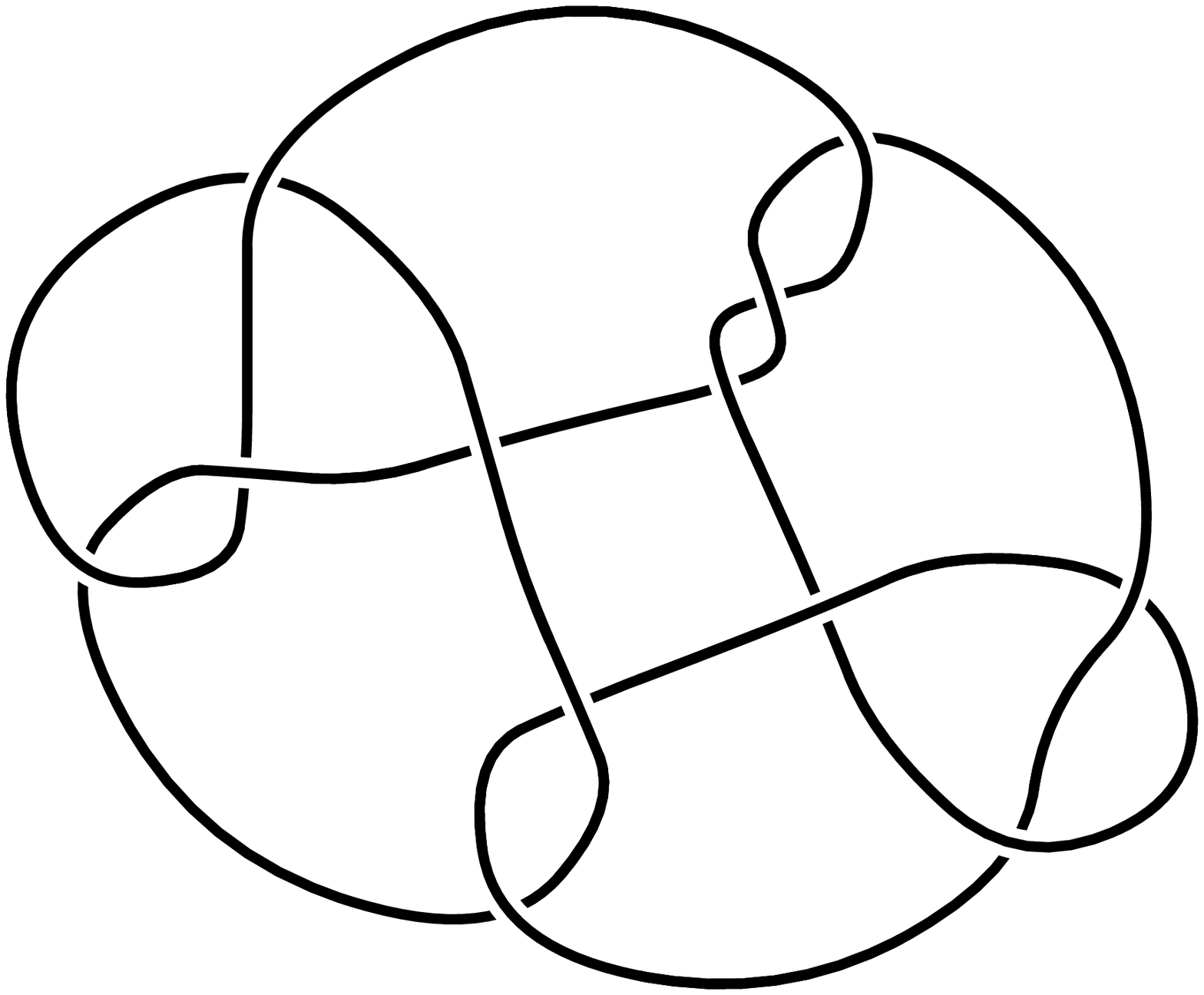}
    &
    \includegraphics[width=75pt]{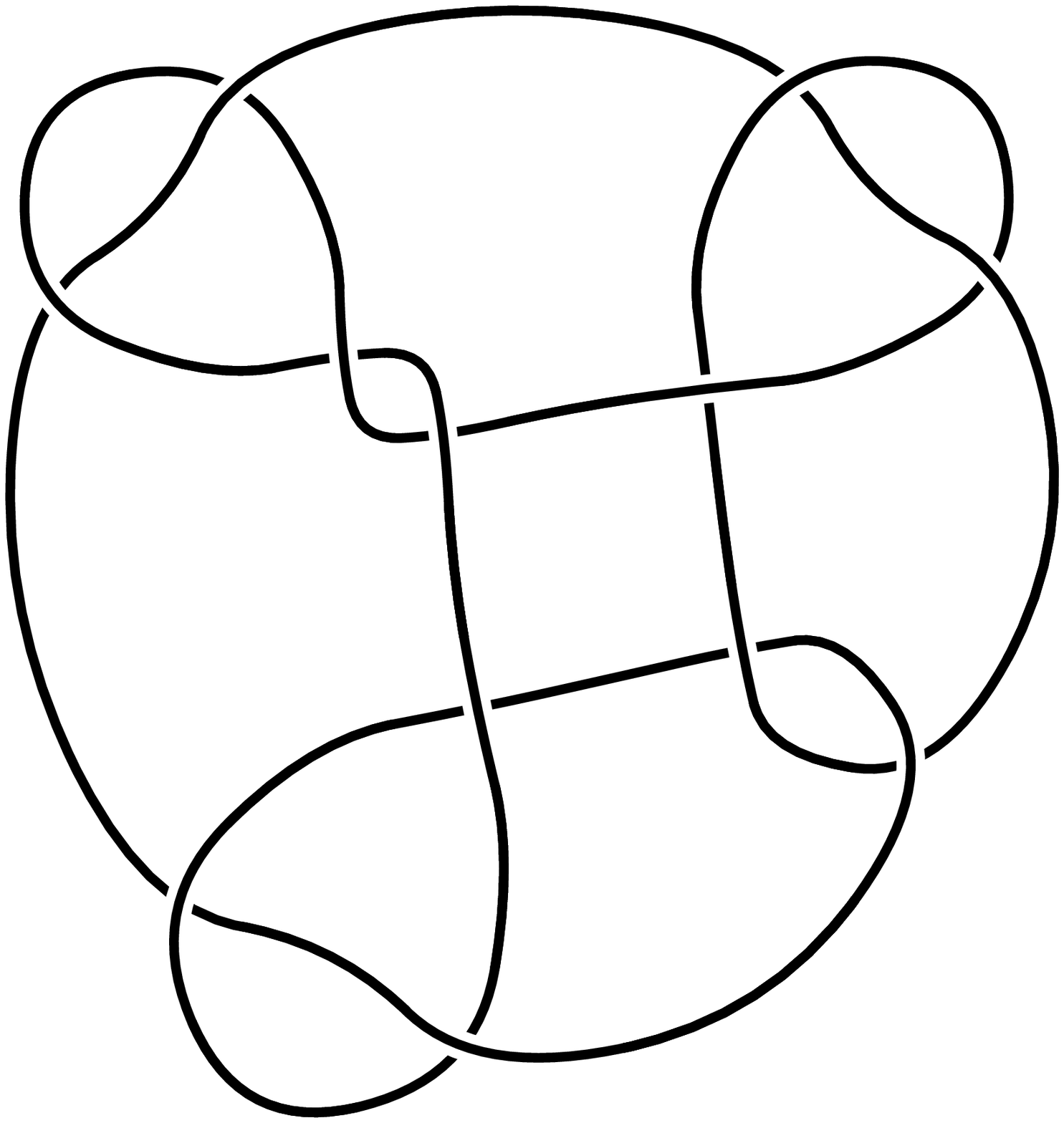}
    \\[-10pt]
    $12^N_{62}$ & $12^N_{66}$ & $12^N_{224}$
  \end{tabular}
  \caption{Nonalternating $12$-crossing mutant cliques 6/6}
  \label{figure:Nonalternating12crossingmutantcliques6of6}
  \end{centering}
\end{figure}

\begin{figure}[htbp]
  \begin{centering}
  \begin{tabular}{ccc}
    \includegraphics[width=75pt]{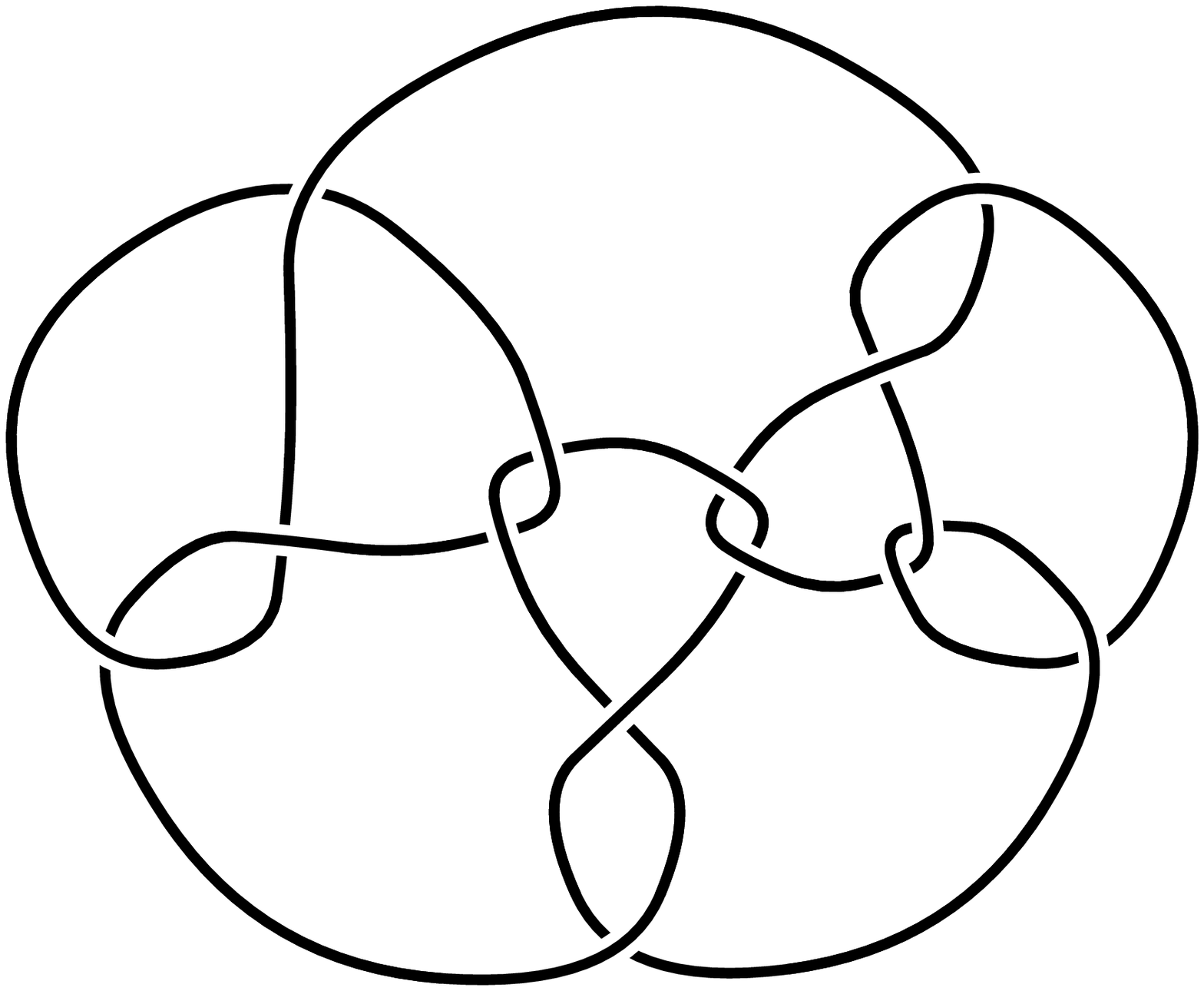}
    &
    \includegraphics[width=75pt]{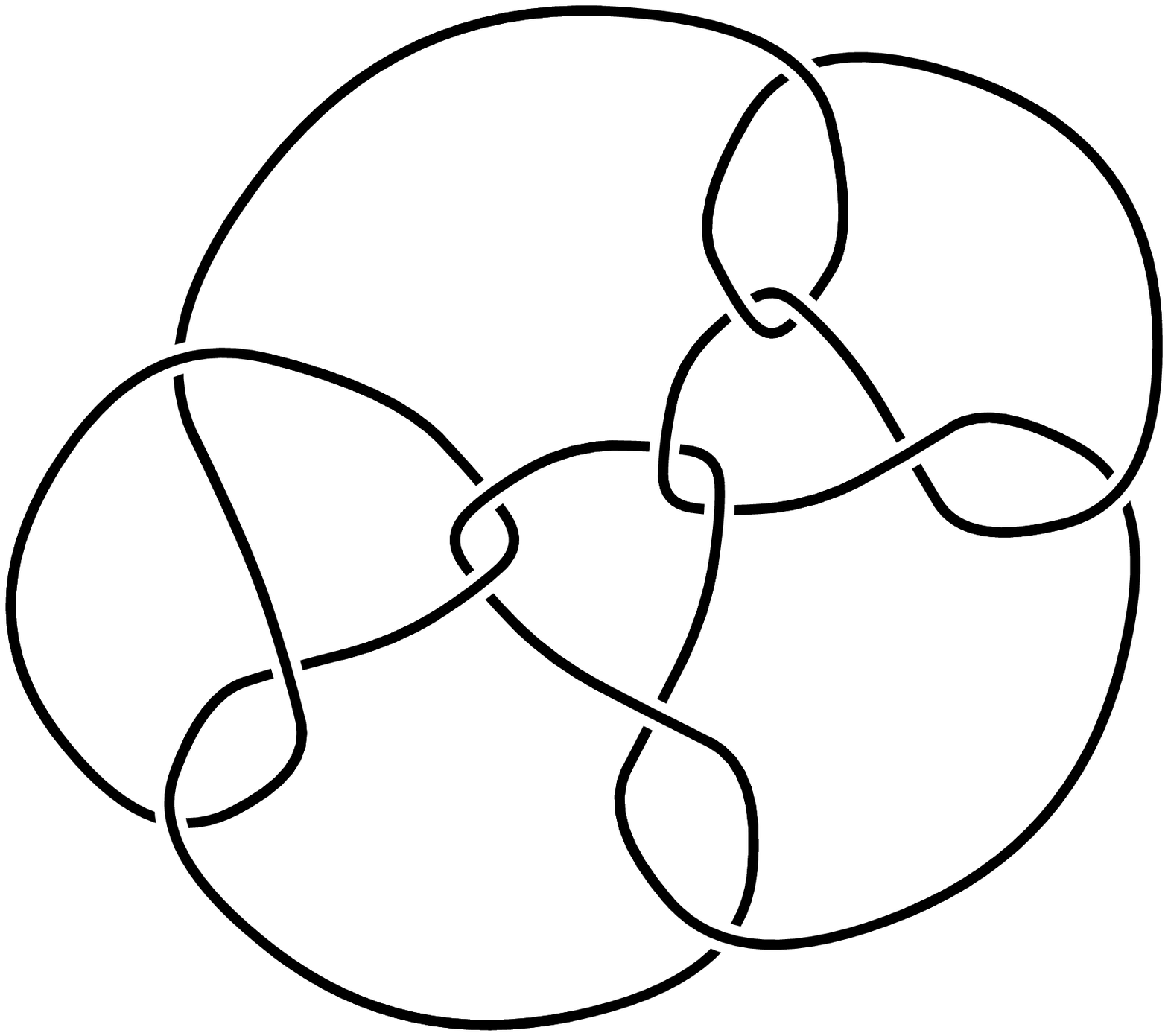}
    &
    \includegraphics[width=75pt]{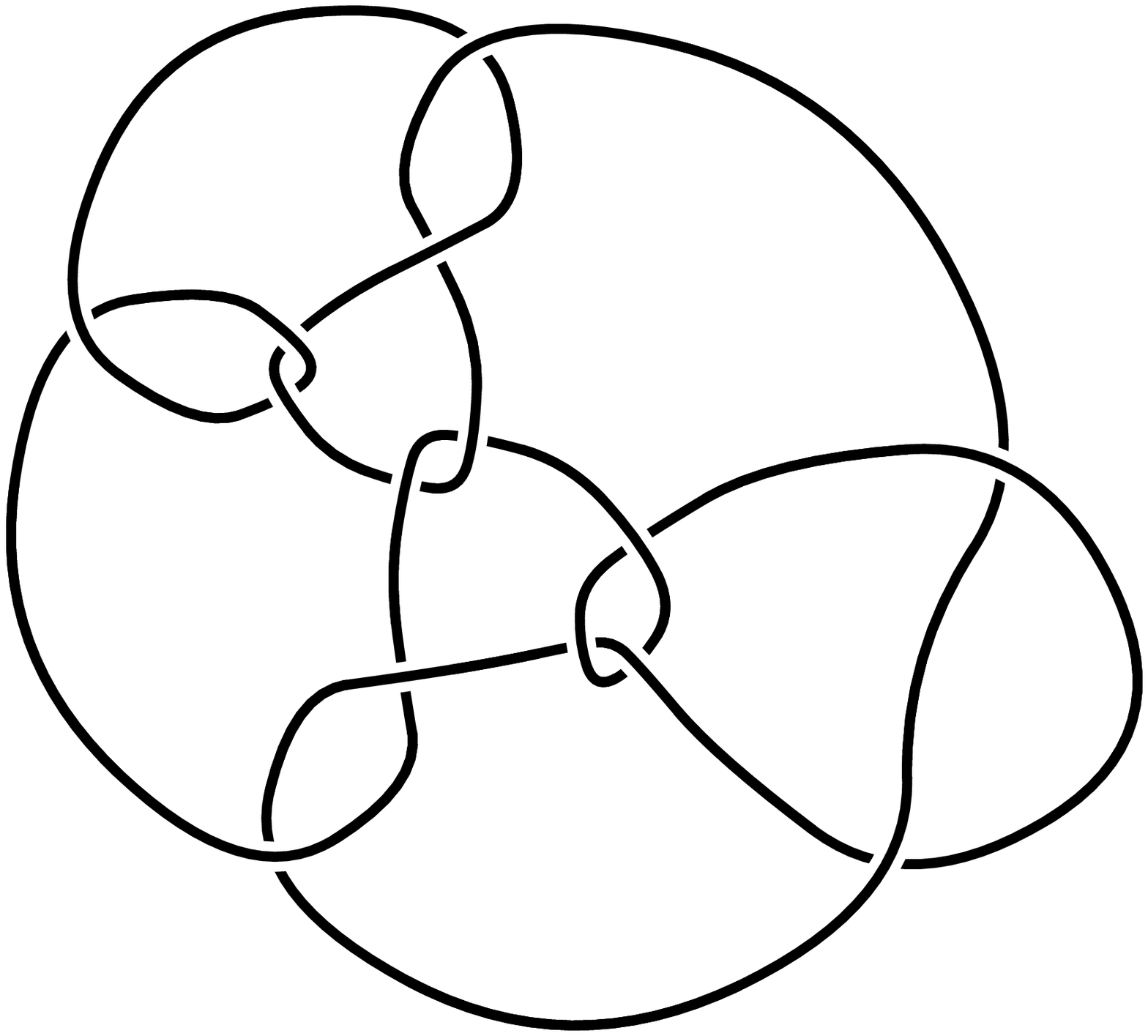}
    \\[-10pt]
    $14^A_{506}$ & $14^A_{486}$ & $14^A_{731}$
    \\[10pt]
    \hline
    &&\\[-10pt]
    \includegraphics[width=75pt]{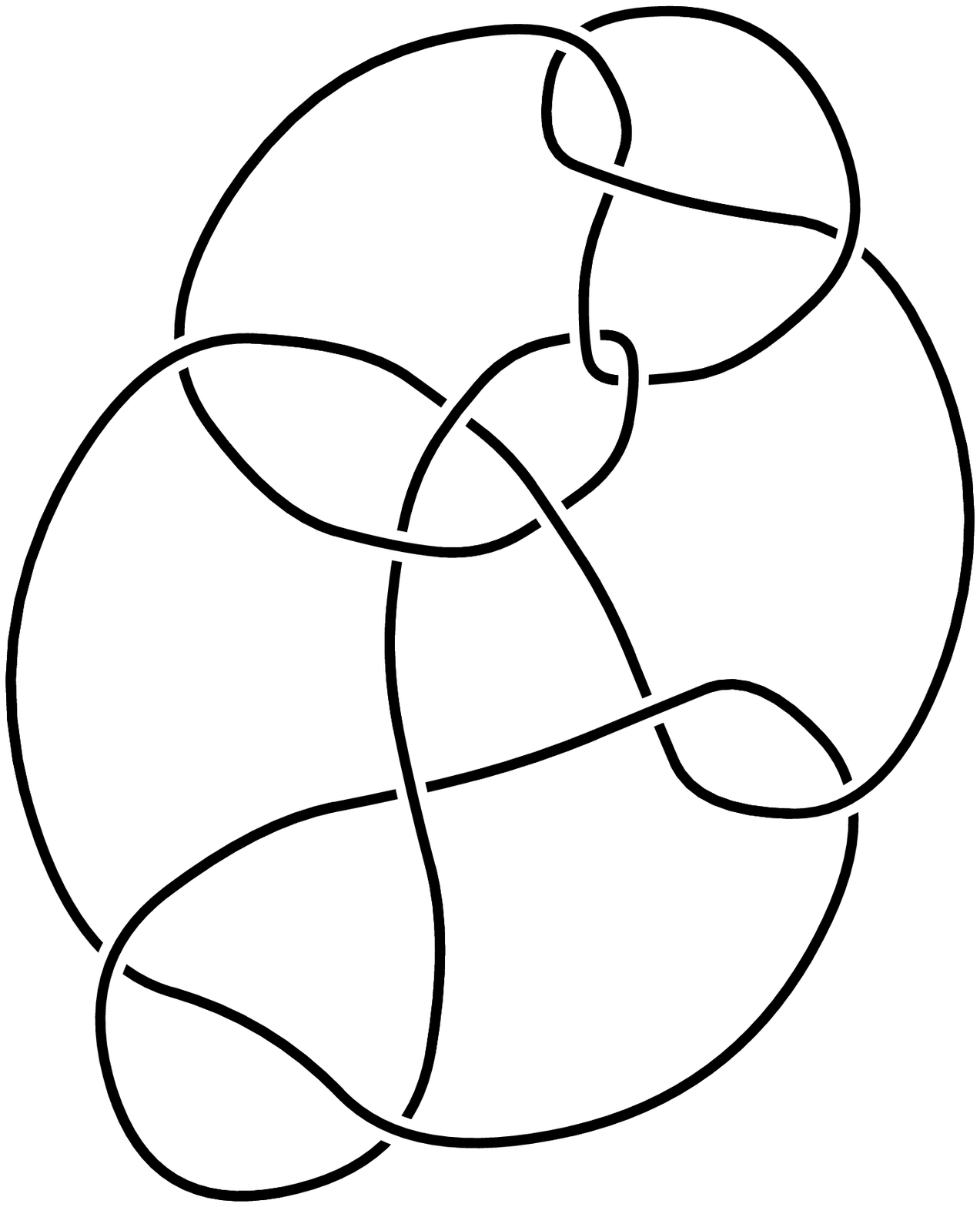}
    &
    \includegraphics[width=75pt]{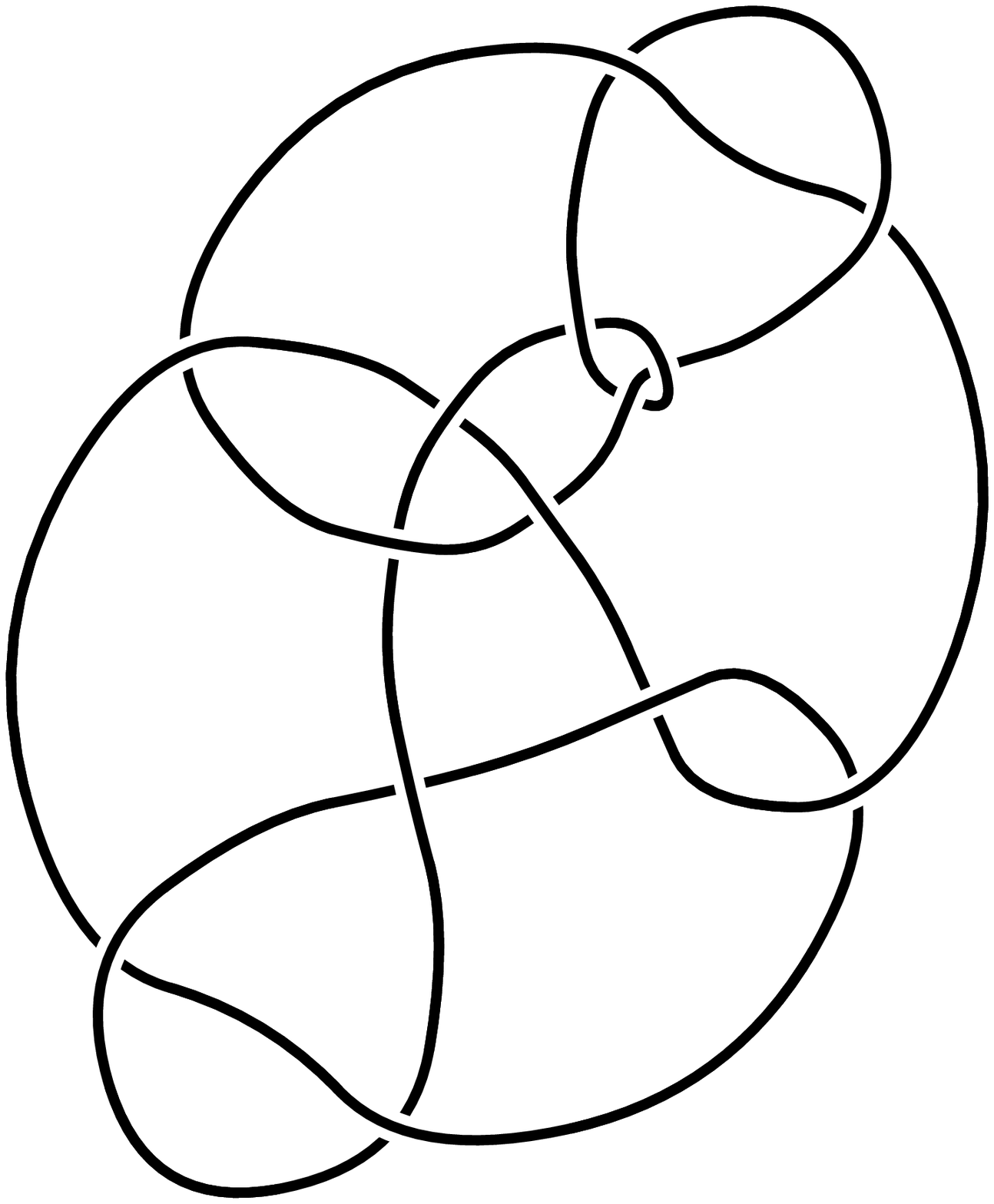}
    &
    \includegraphics[width=75pt]{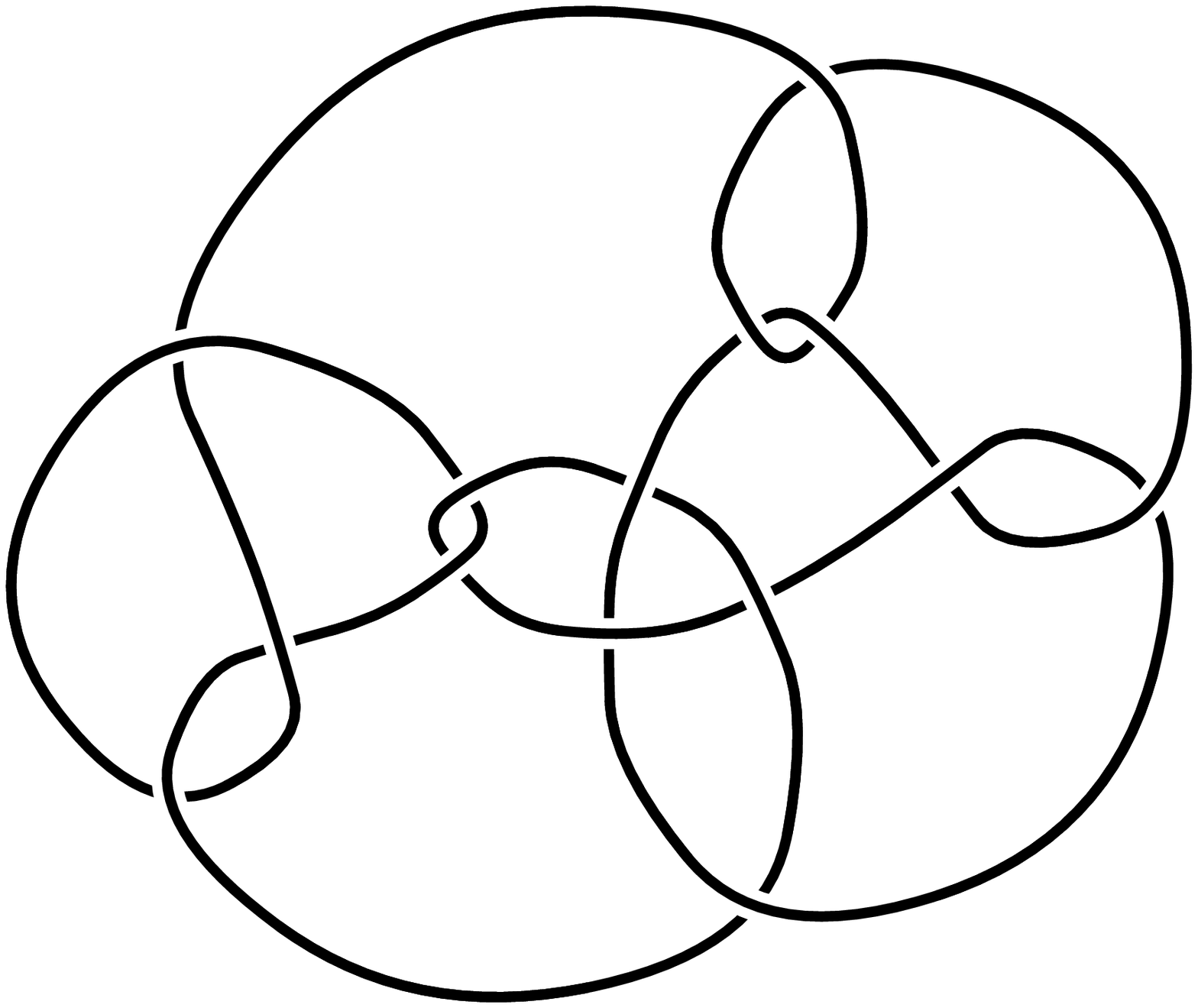}
    \\[-10pt]
    $14^A_{680}$ & $14^A_{509}$ & $14^A_{585}$
    \\[10pt]
    \hline
    &&\\[-10pt]
    \includegraphics[width=75pt]{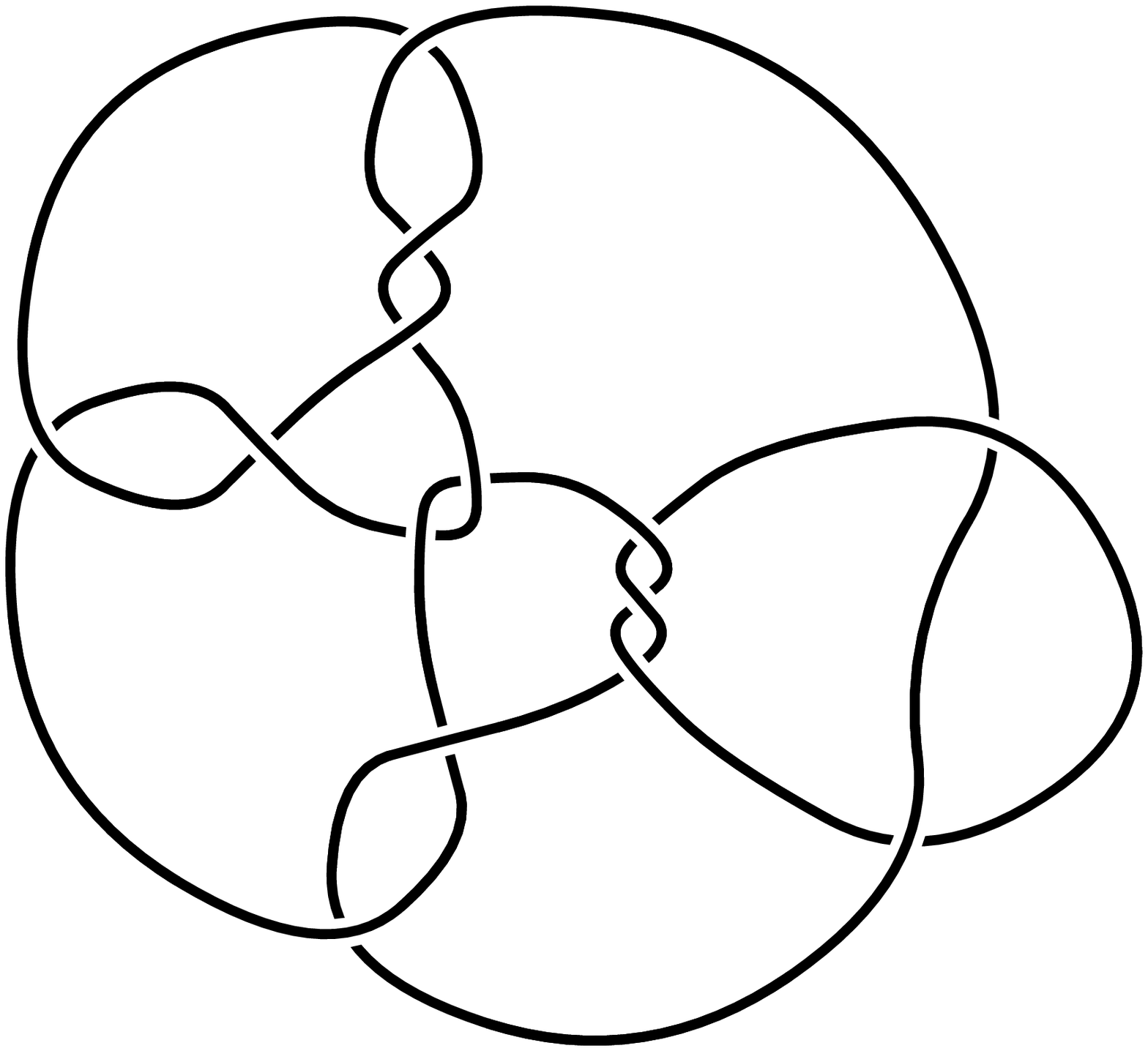}
    &
    \includegraphics[width=75pt]{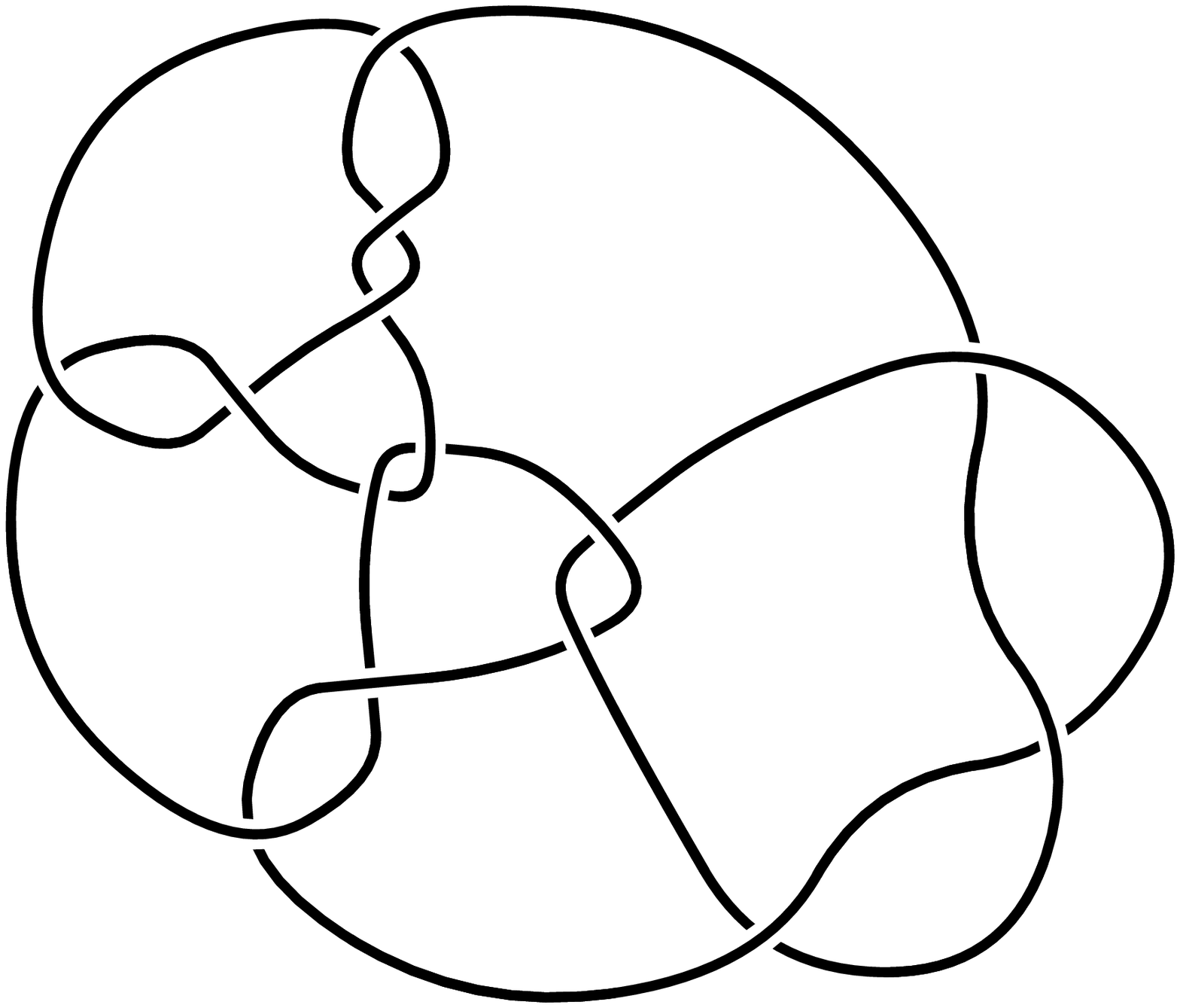}
    &
    \includegraphics[width=75pt,angle=90]{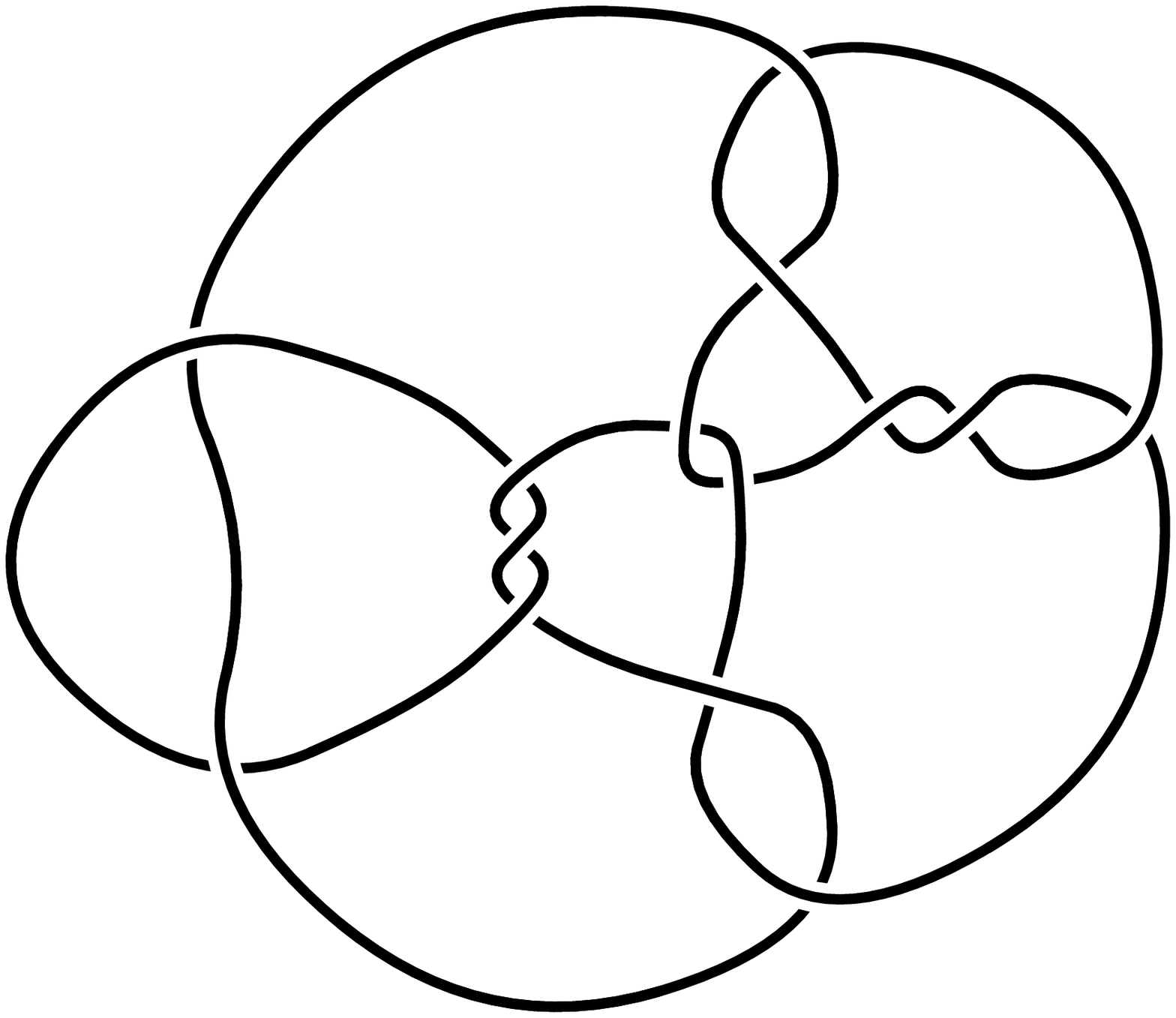}
    \\[-10pt]
    $14^A_{12813}$ & $14^A_{12807}$ & $14^A_{12875}$
    \\[10pt]
    \hline
    &&\\[-10pt]
    \includegraphics[width=75pt,angle=90]{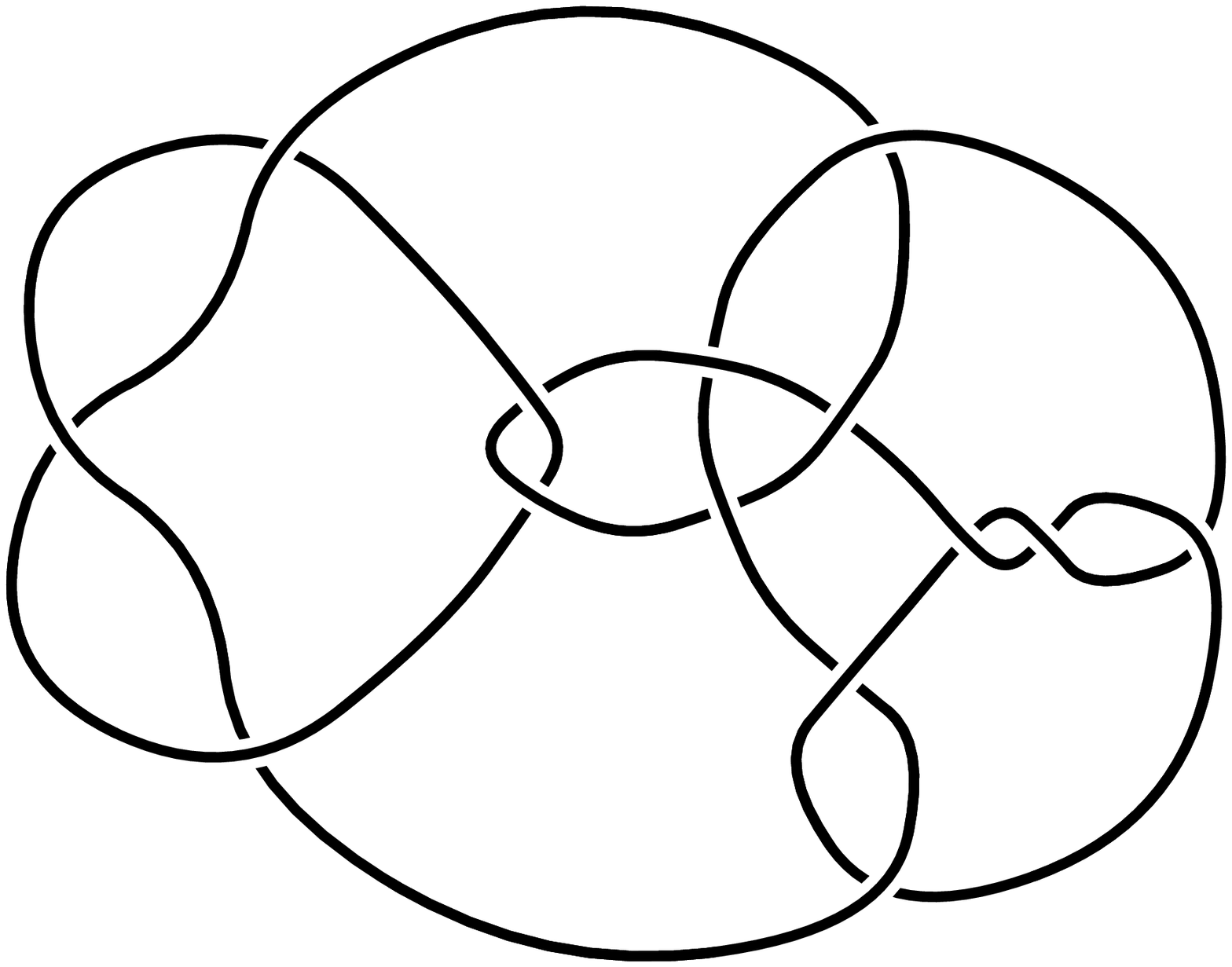}
    &
    \includegraphics[width=75pt]{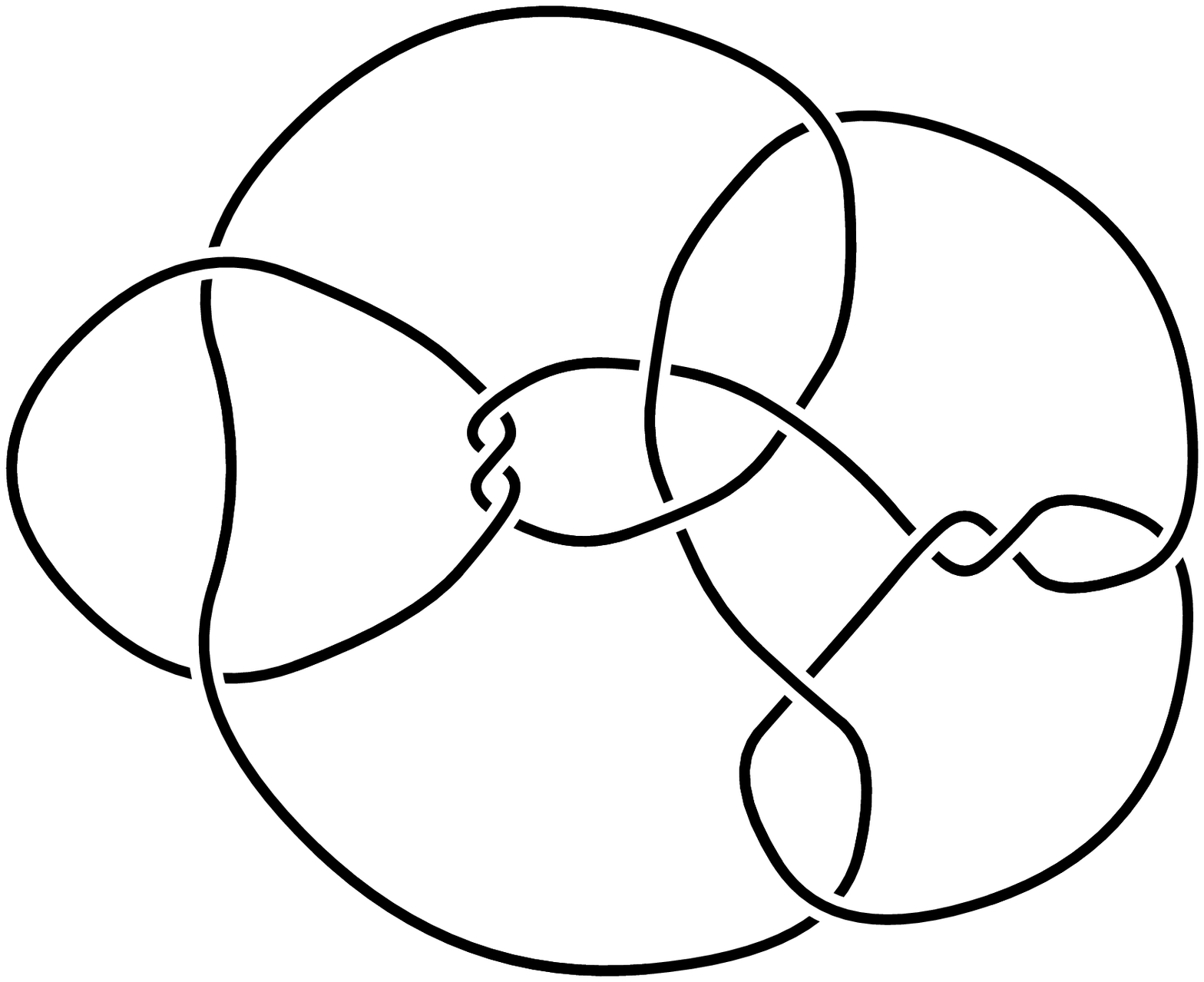}
    &
    \includegraphics[width=75pt]{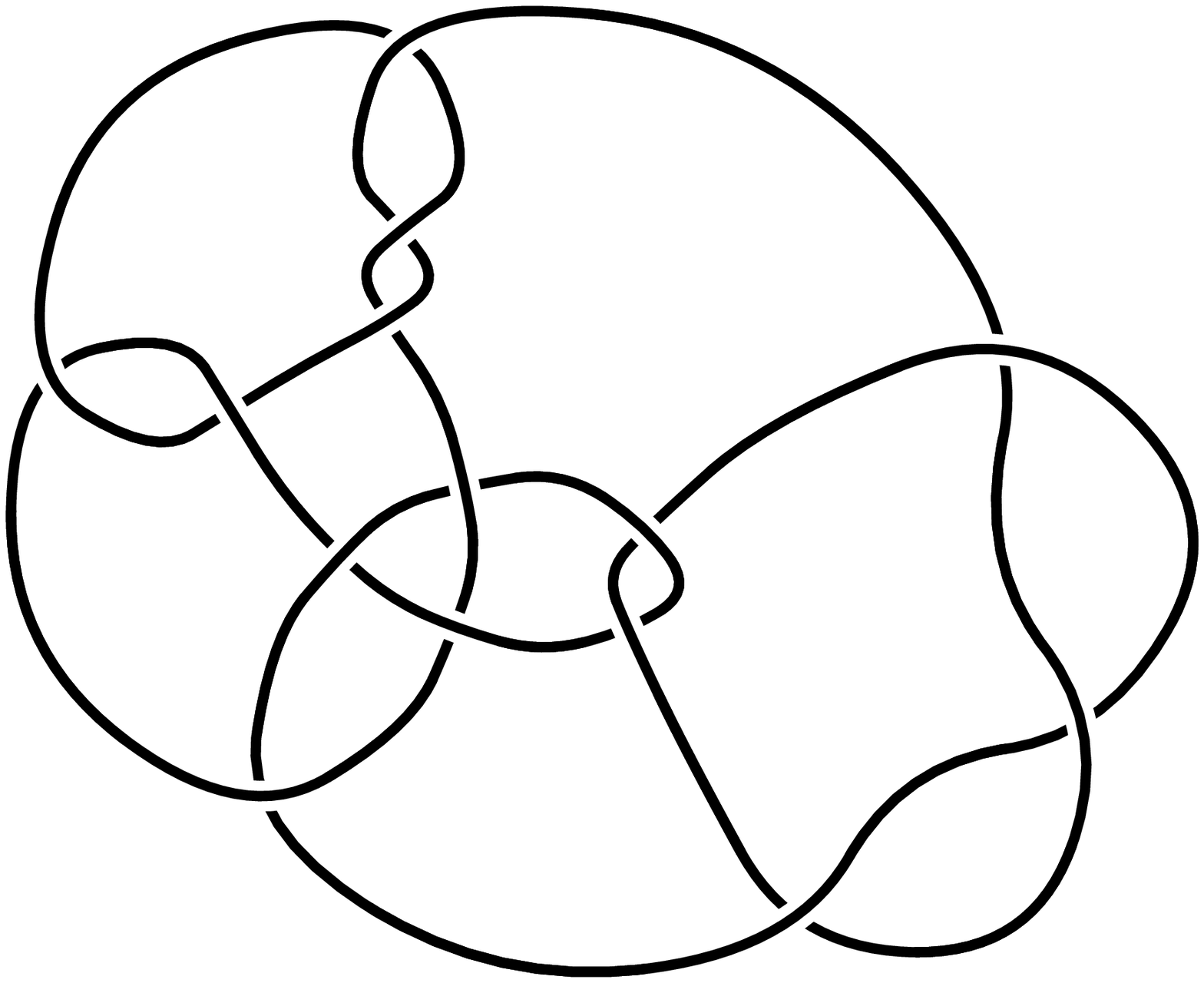}
    \\[-10pt]
    $14^A_{12858}$ & $14^A_{12815}$ & $14^A_{12830}$
    \\[10pt]
    \hline
    &&\\[-10pt]
    \includegraphics[width=75pt]{Graphics/Mutant_Cliques/14/A/14A_13107}
    &
    \includegraphics[width=75pt]{Graphics/Mutant_Cliques/14/A/14A_13109}
    &
    \includegraphics[width=75pt]{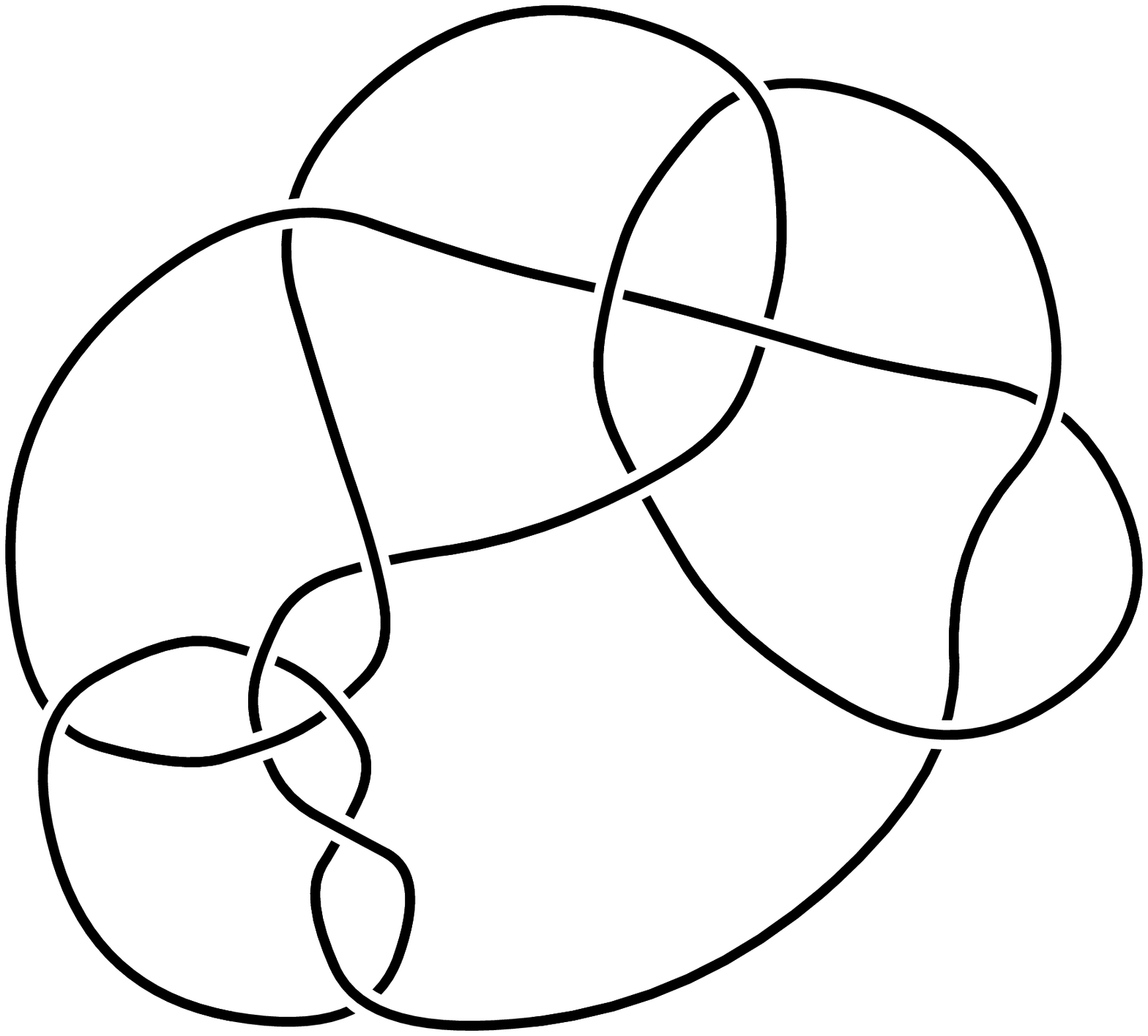}
    \\[-10pt]
    $14^A_{13107}$ & $14^A_{13109}$ & $14^A_{13489}$
  \end{tabular}
  \caption{Alternating $14$-crossing mutant cliques containing both chiral and achiral elements 1/2.}
  \label{figure:Alternating14crossingmutantcliques1of2}
  \end{centering}
\end{figure}

\begin{figure}[htbp]
  \begin{centering}
  \begin{tabular}{ccc}
    \includegraphics[width=75pt]{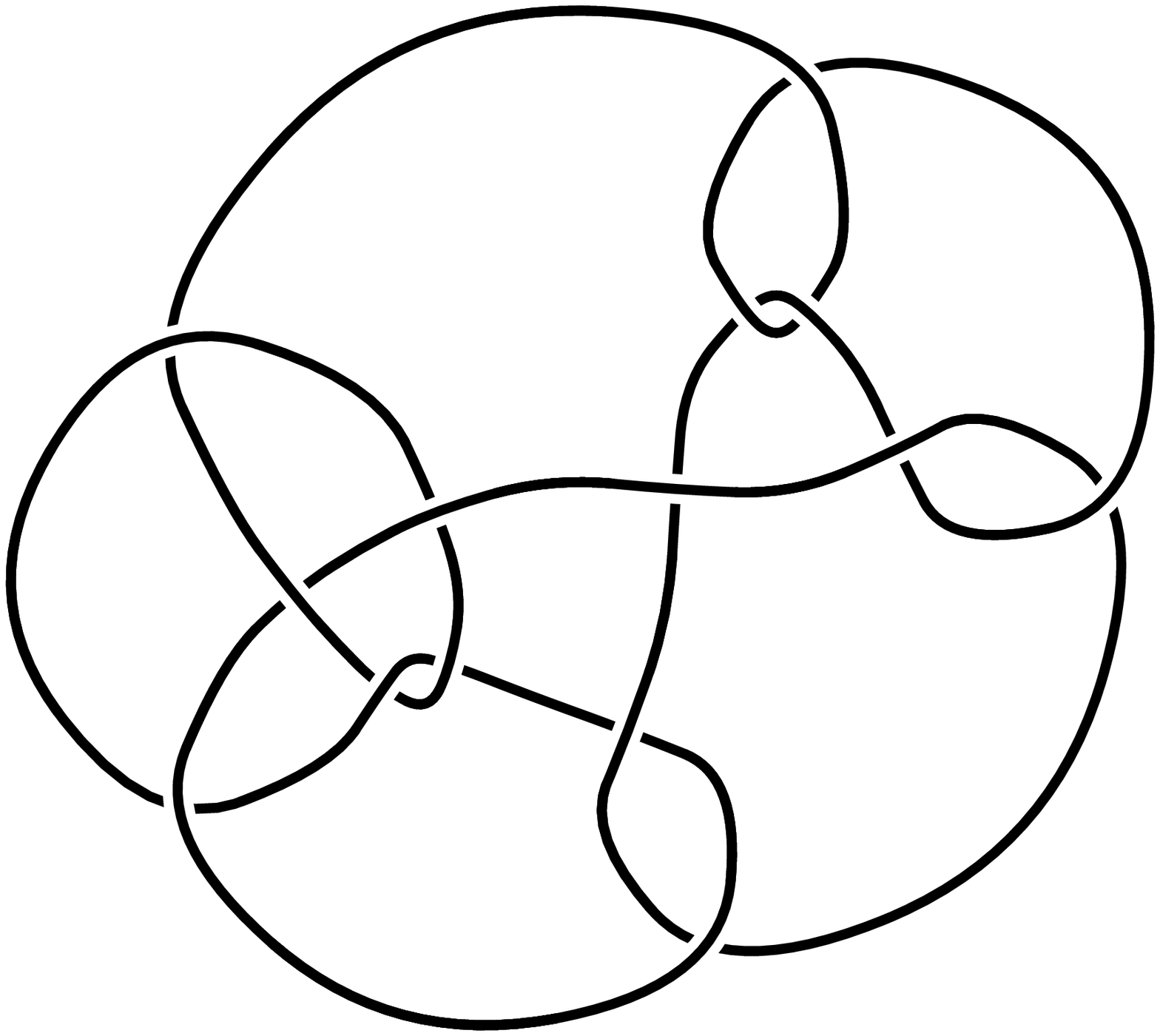}
    &
    \includegraphics[width=75pt]{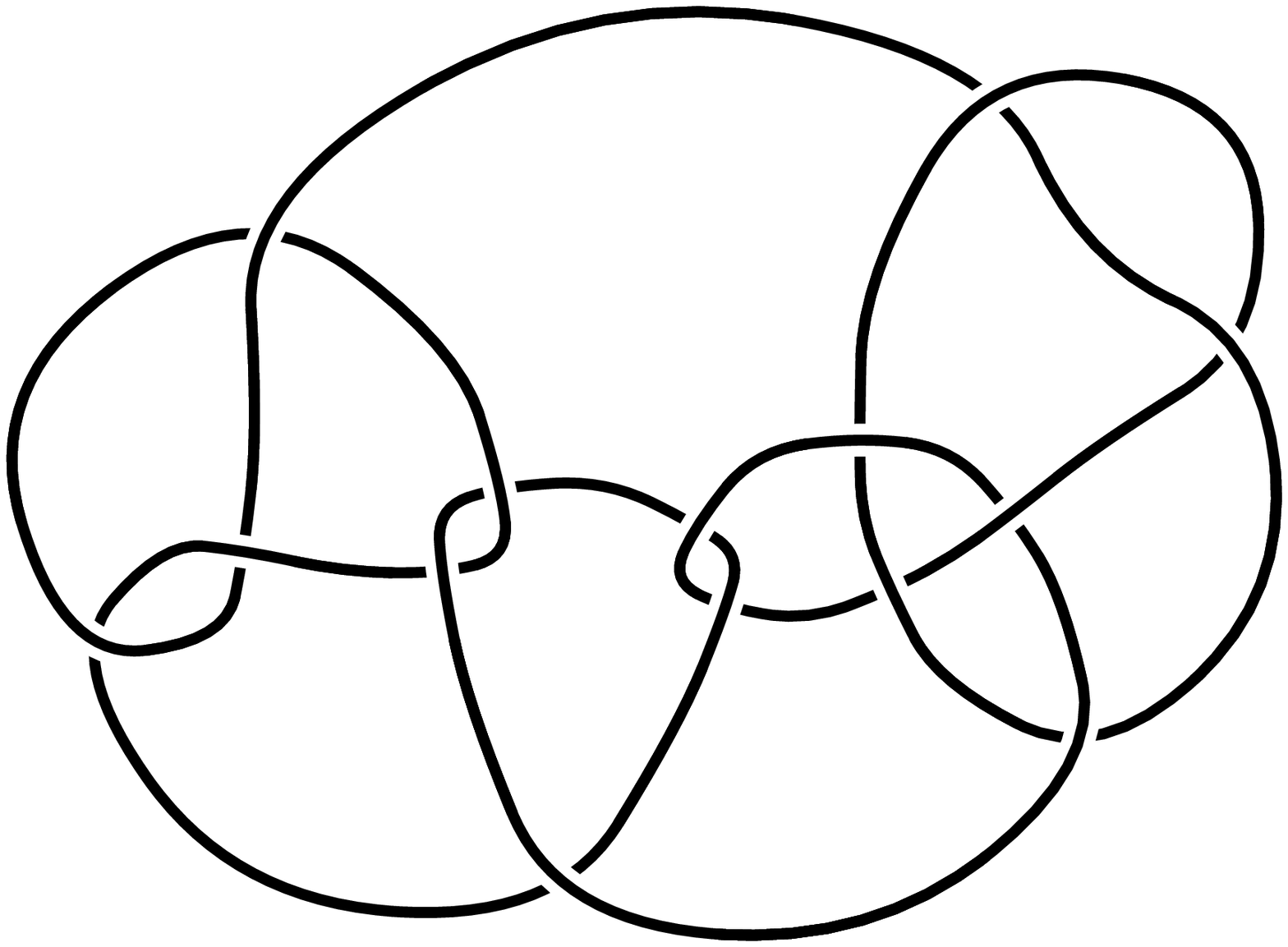}
    &
    \includegraphics[width=75pt]{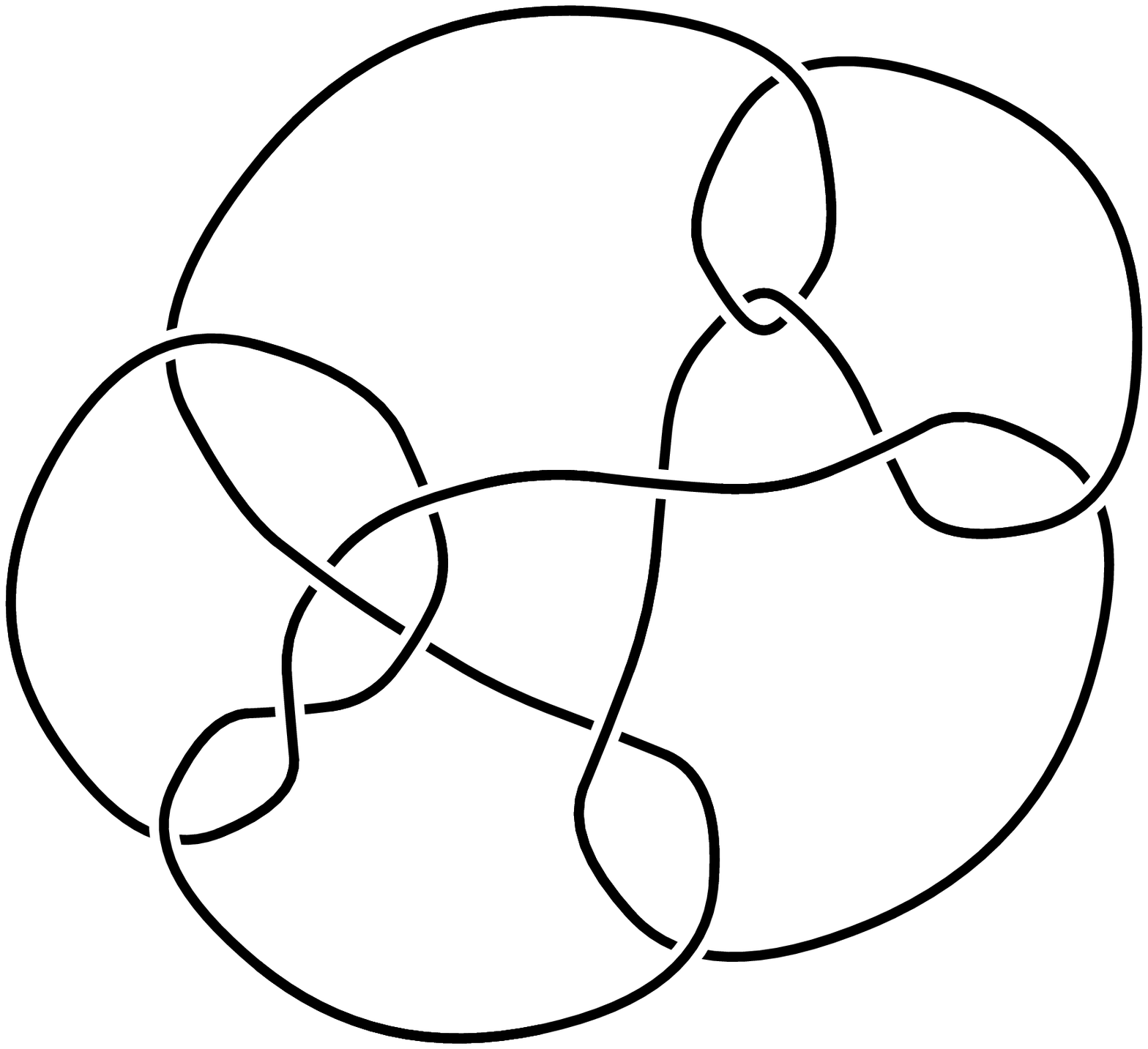}
    \\[-10pt]
    $14^N_{1309}$ & $14^N_{1327}$ & $14^N_{1497}$
    \\[10pt]
    \hline
    &&\\[-10pt]
    \includegraphics[width=75pt]{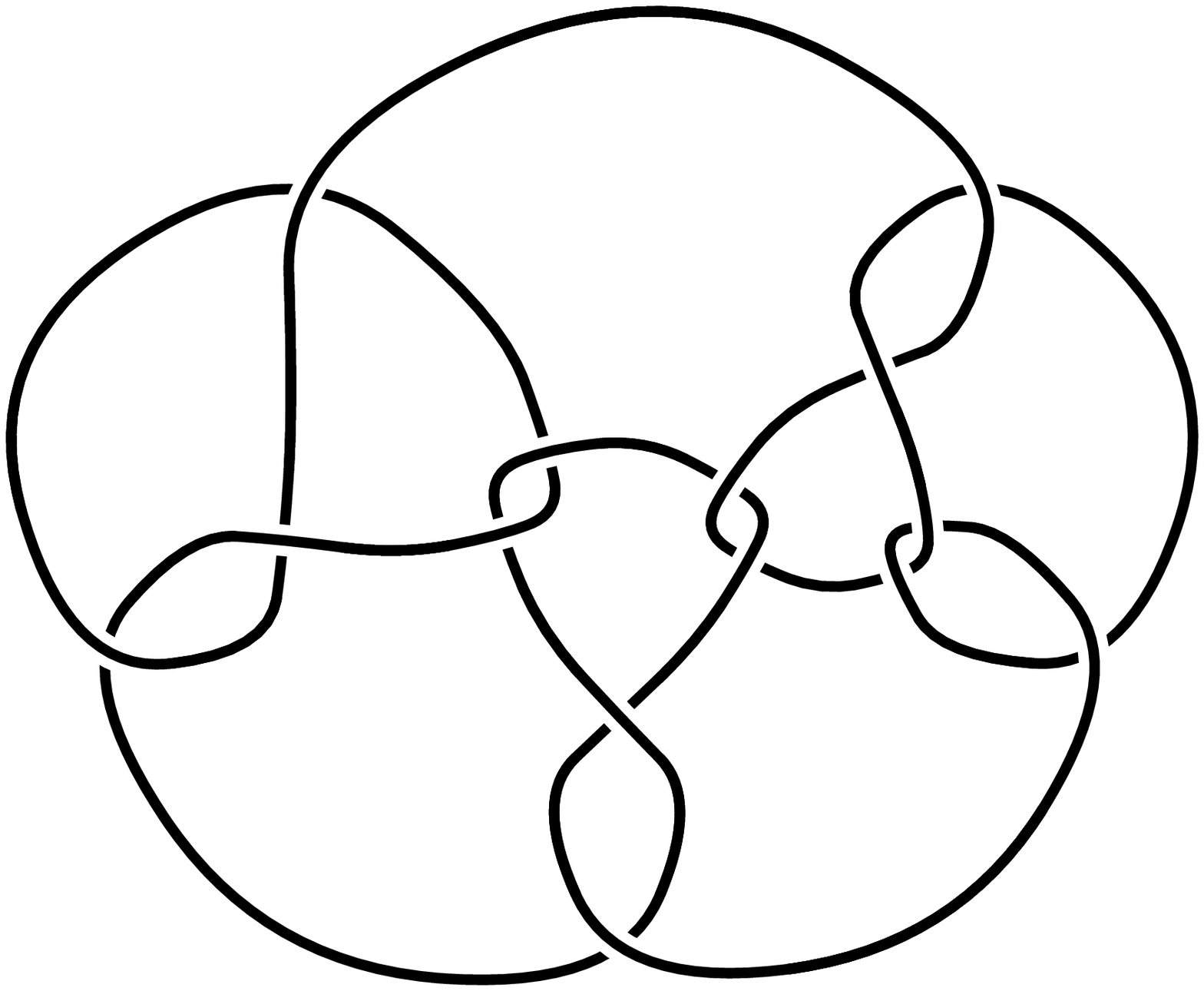}
    &
    \includegraphics[width=75pt]{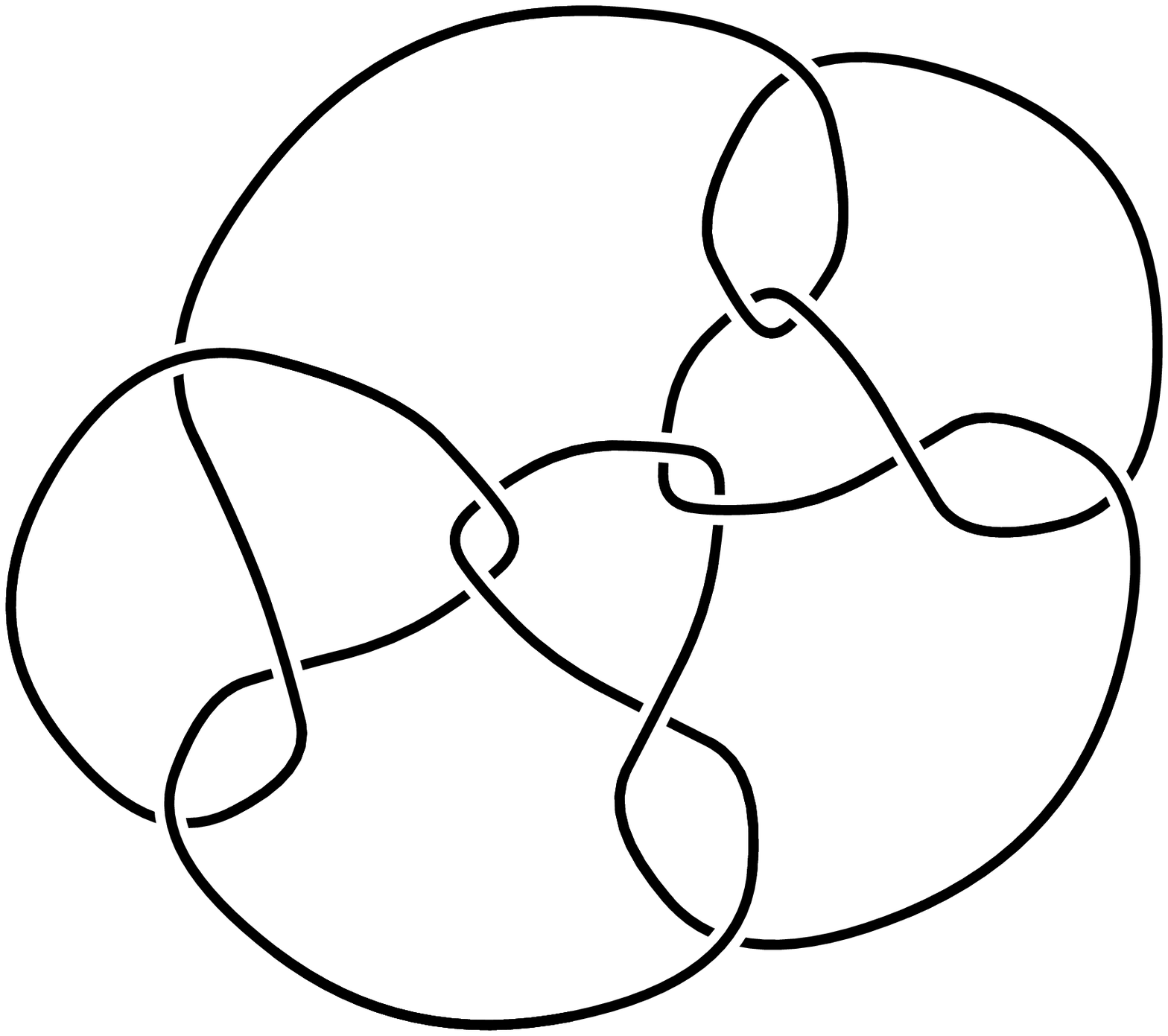}
    &
    \includegraphics[width=75pt]{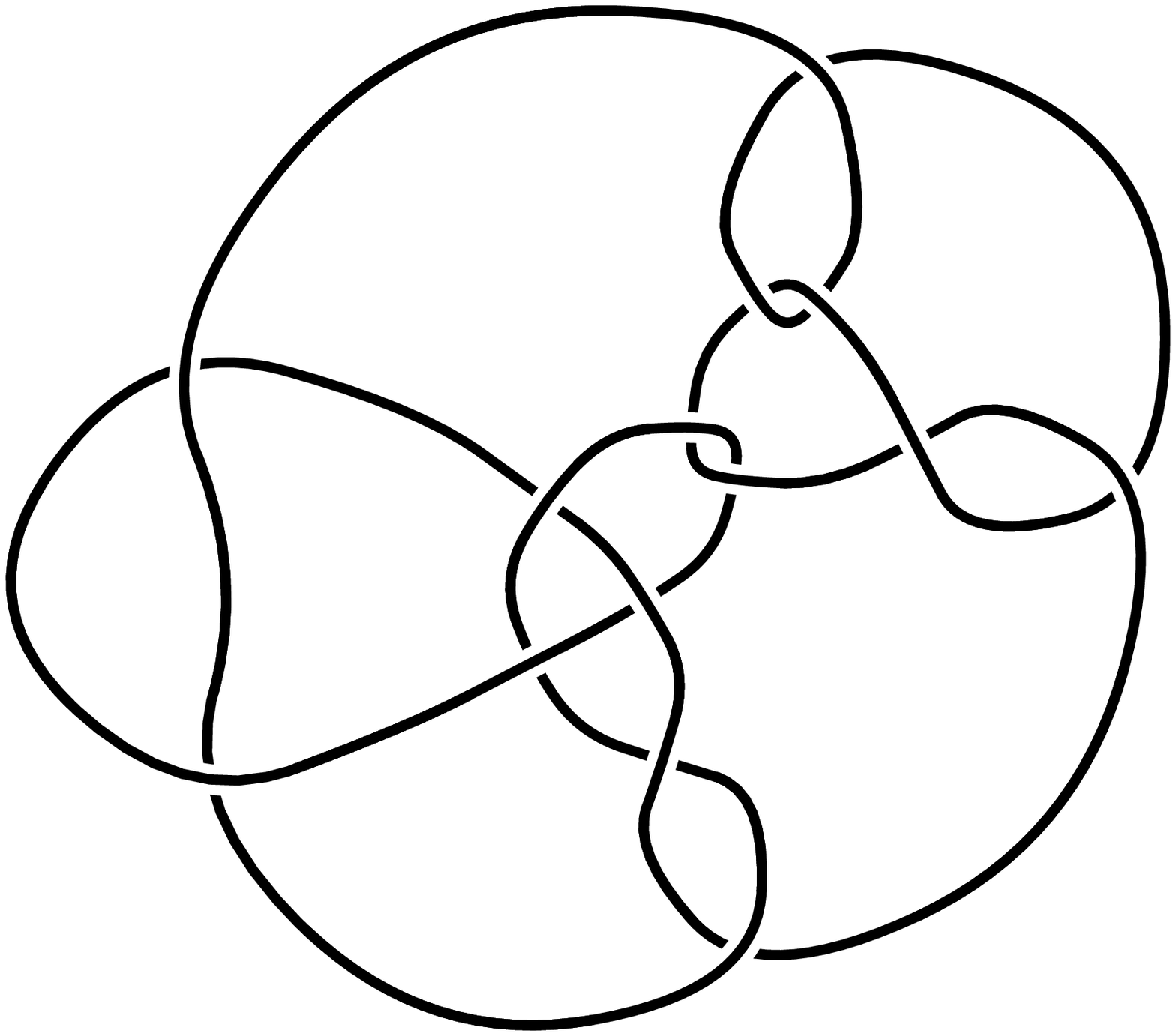}
    \\[-10pt]
    $14^N_{1641}$ & $14^N_{1552}$ & $14^N_{2132}$
    \\[10pt]
    \hline
    &&\\[-10pt]
    \includegraphics[width=75pt]{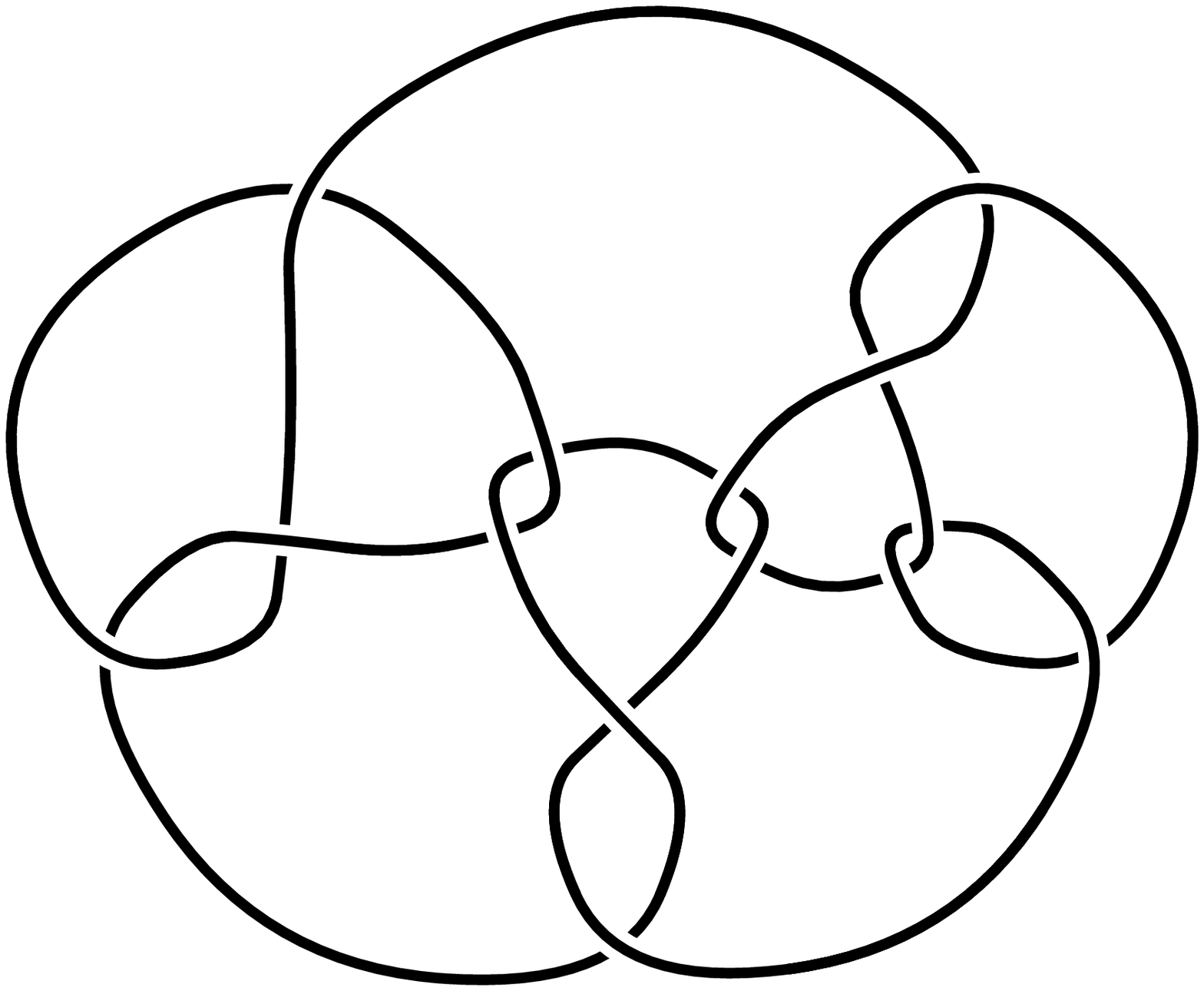}
    &
    \includegraphics[width=75pt]{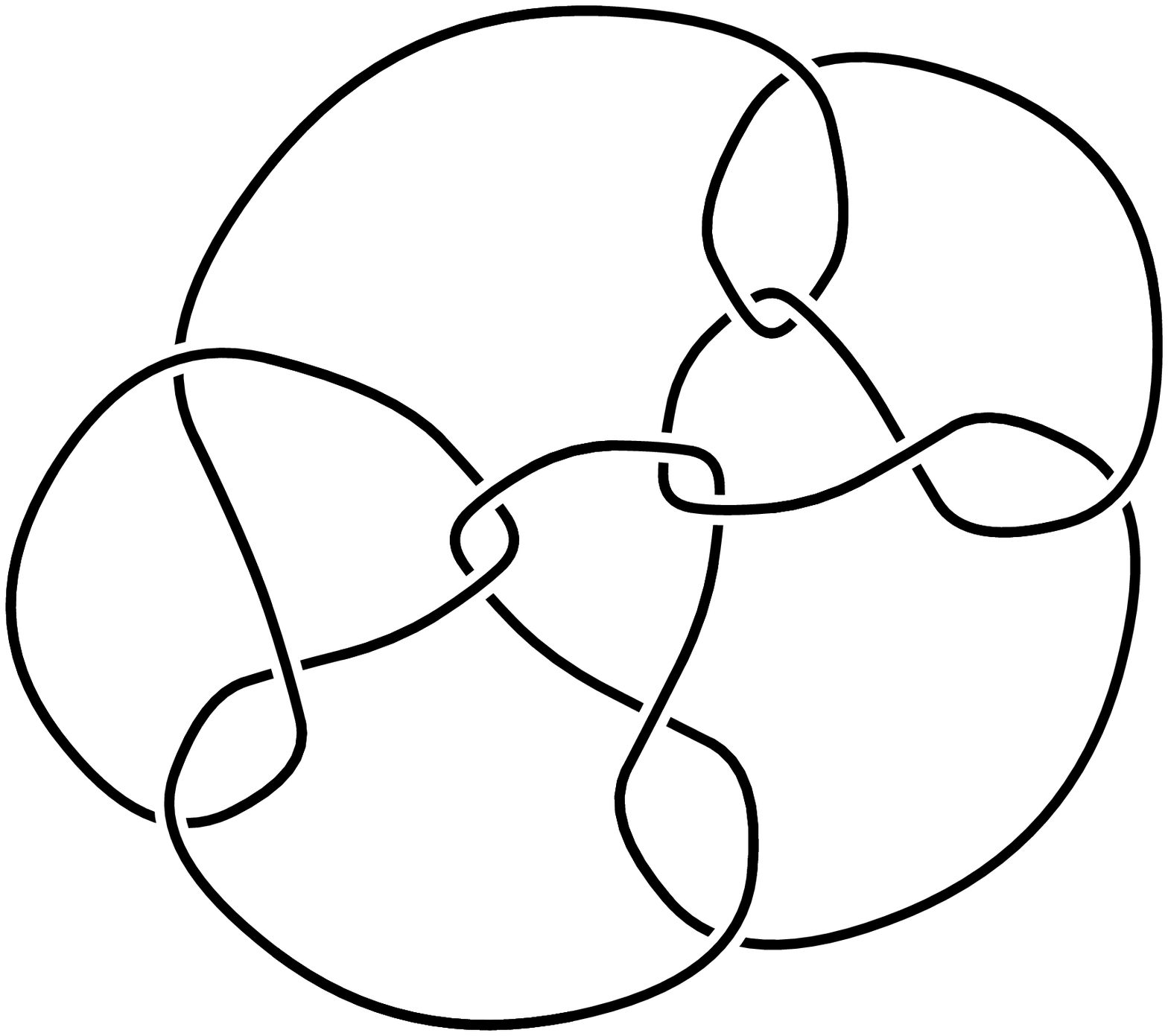}
    &
    \includegraphics[width=75pt]{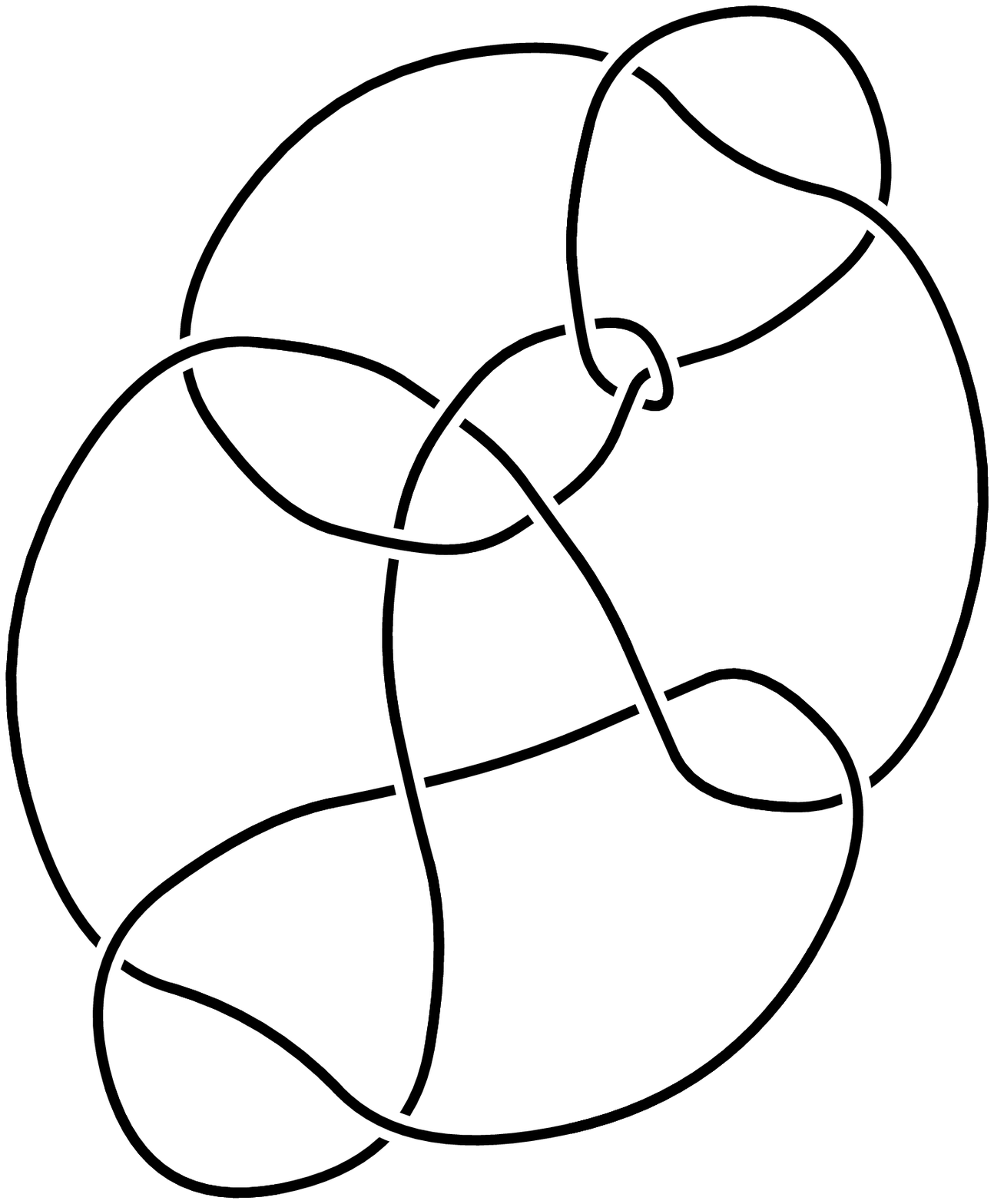}
    \\[-10pt]
    $14^N_{1644}$ & $14^N_{1555}$ & $14^N_{1671}$
    \\[10pt]
    \hline
    &&\\[-10pt]
    \includegraphics[width=75pt]{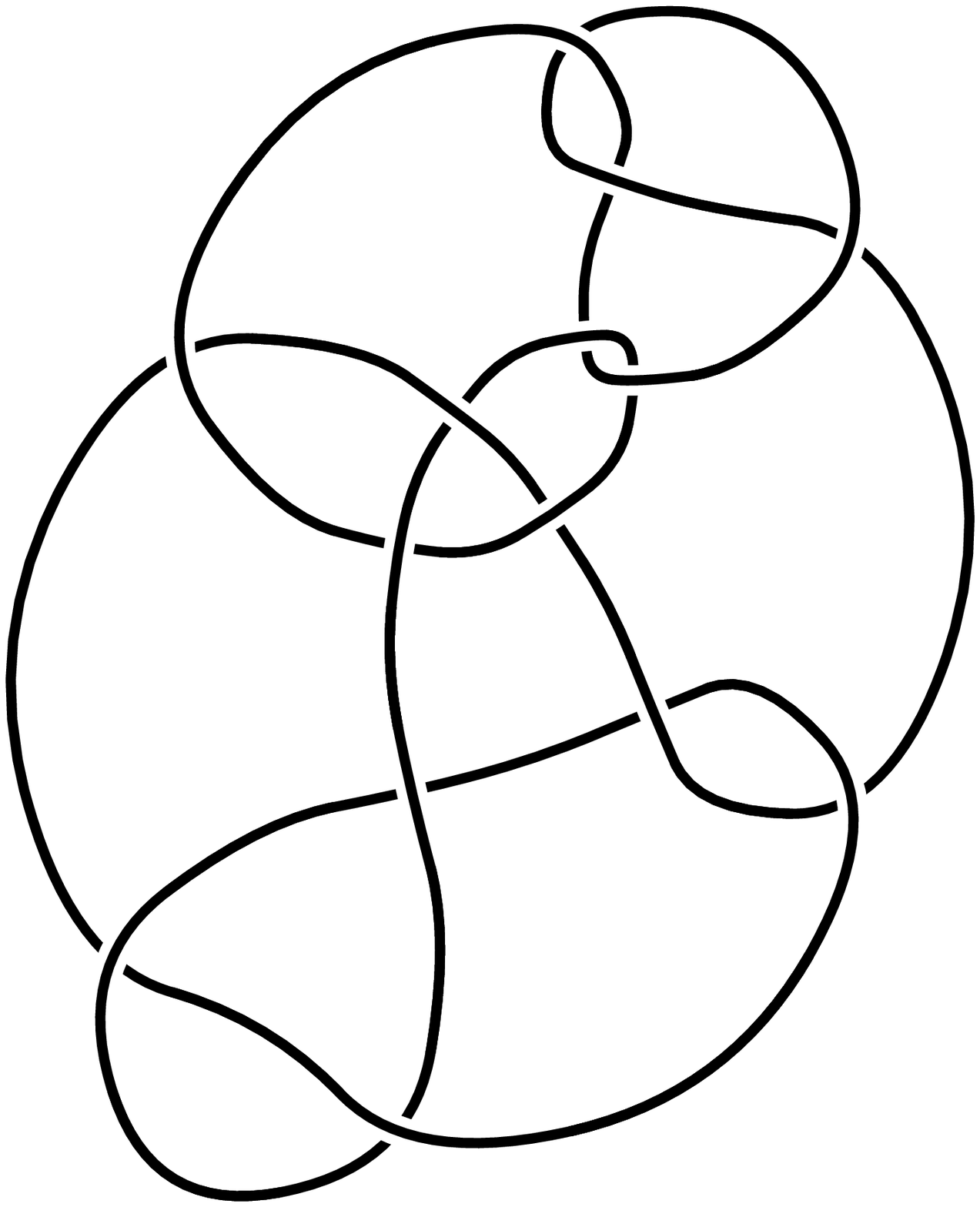}
    &
    \includegraphics[width=75pt]{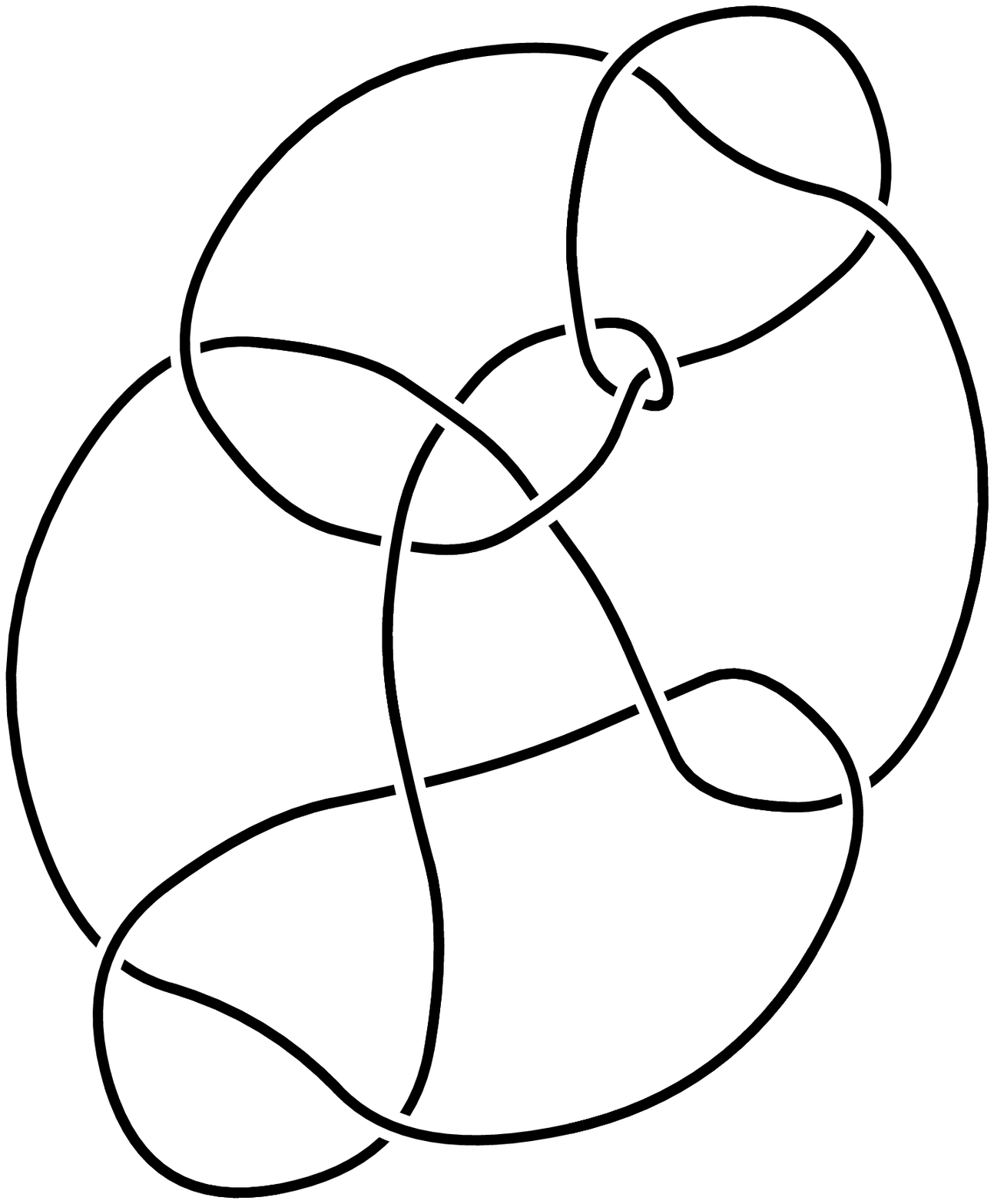}
    &
    \includegraphics[width=75pt]{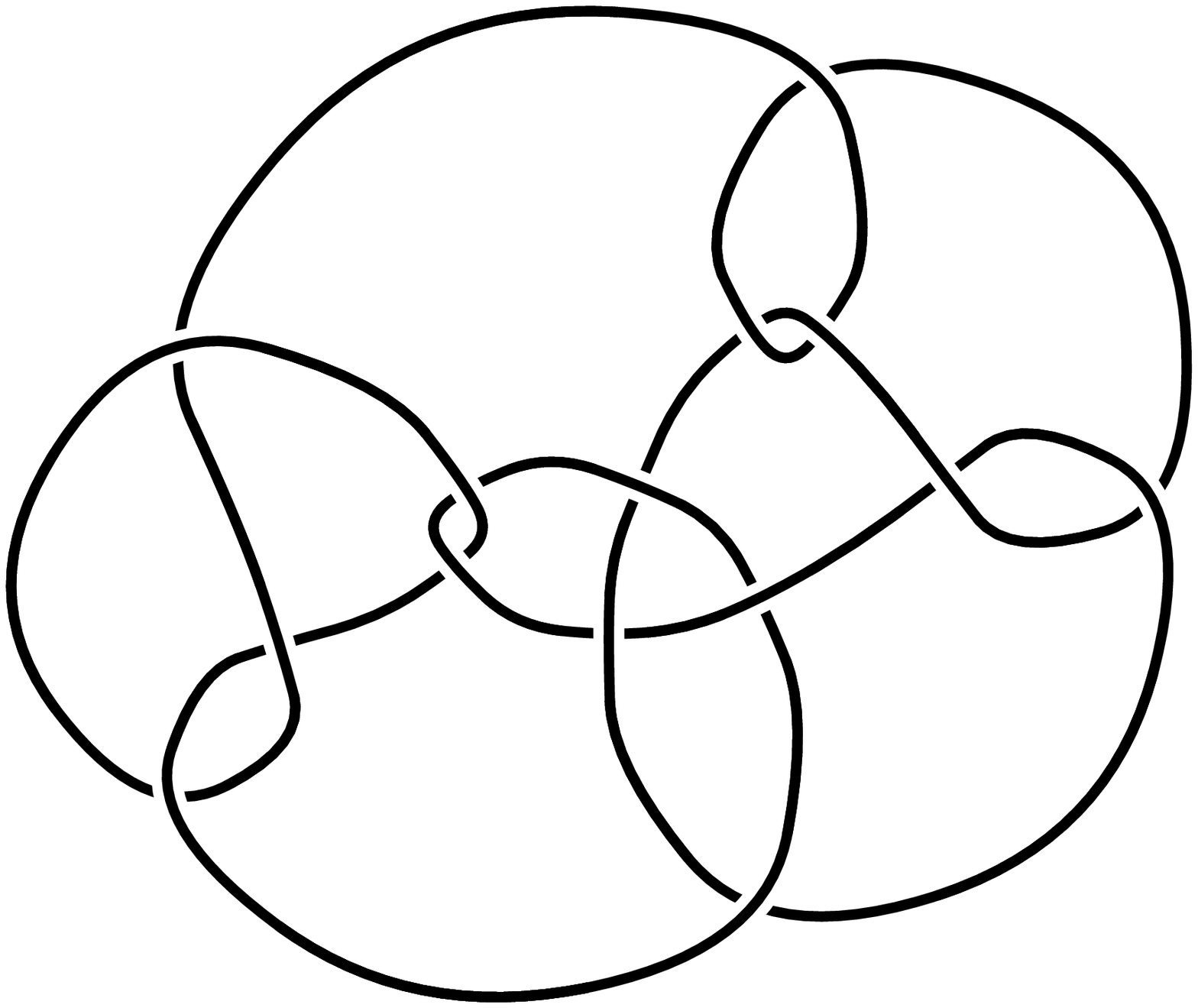}
    \\[-10pt]
    $14^N_{2164}$ & $14^N_{1669}$ & $14^N_{1925}$
  \end{tabular}
  \caption{Nonalternating $14$-crossing mutant cliques containing both chiral and achiral elements.}
  \label{figure:Nonalternating14crossingmutantcliques}
  \end{centering}
\end{figure}

On the following pages are illustrated some mutant cliques mentioned within the main text.
Figures~\ref{figure:Alternating12crossingmutantcliques1of4}--\ref{figure:Nonalternating12crossingmutantcliques6of6}
illustrate the mutant cliques within $\mathbb{K}_{12}$, and
Figures~\ref{figure:Alternating14crossingmutantcliques1of2}--\ref{figure:Nonalternating14crossingmutantcliques}
illustrate the mutant cliques containing both chiral and achiral knots within $\mathbb{K}_{14}$.

\begin{figure}[htbp]
  \begin{centering}
  \includegraphics[width=120pt]{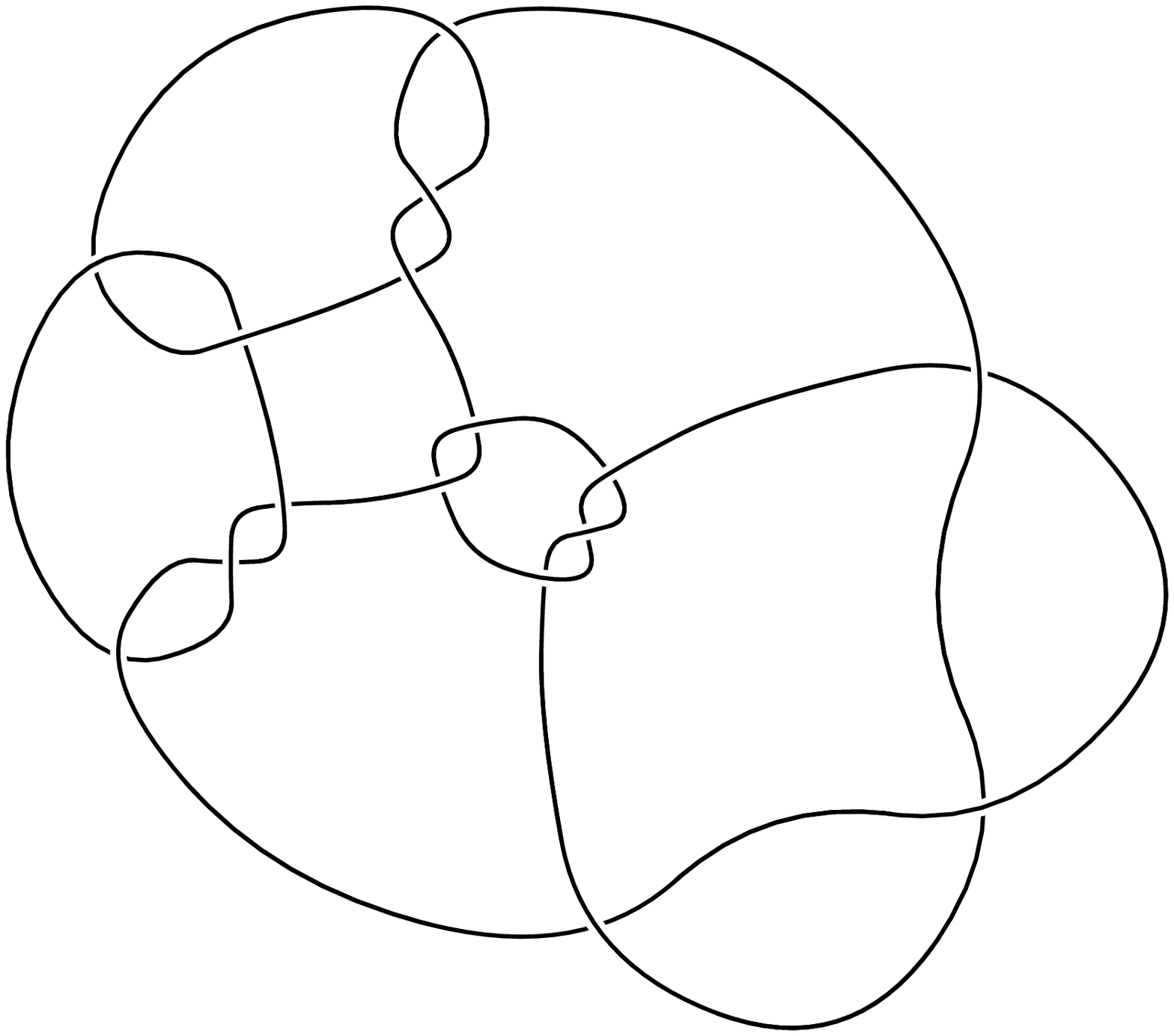}
  \hspace{20pt}
  \includegraphics[width=120pt]{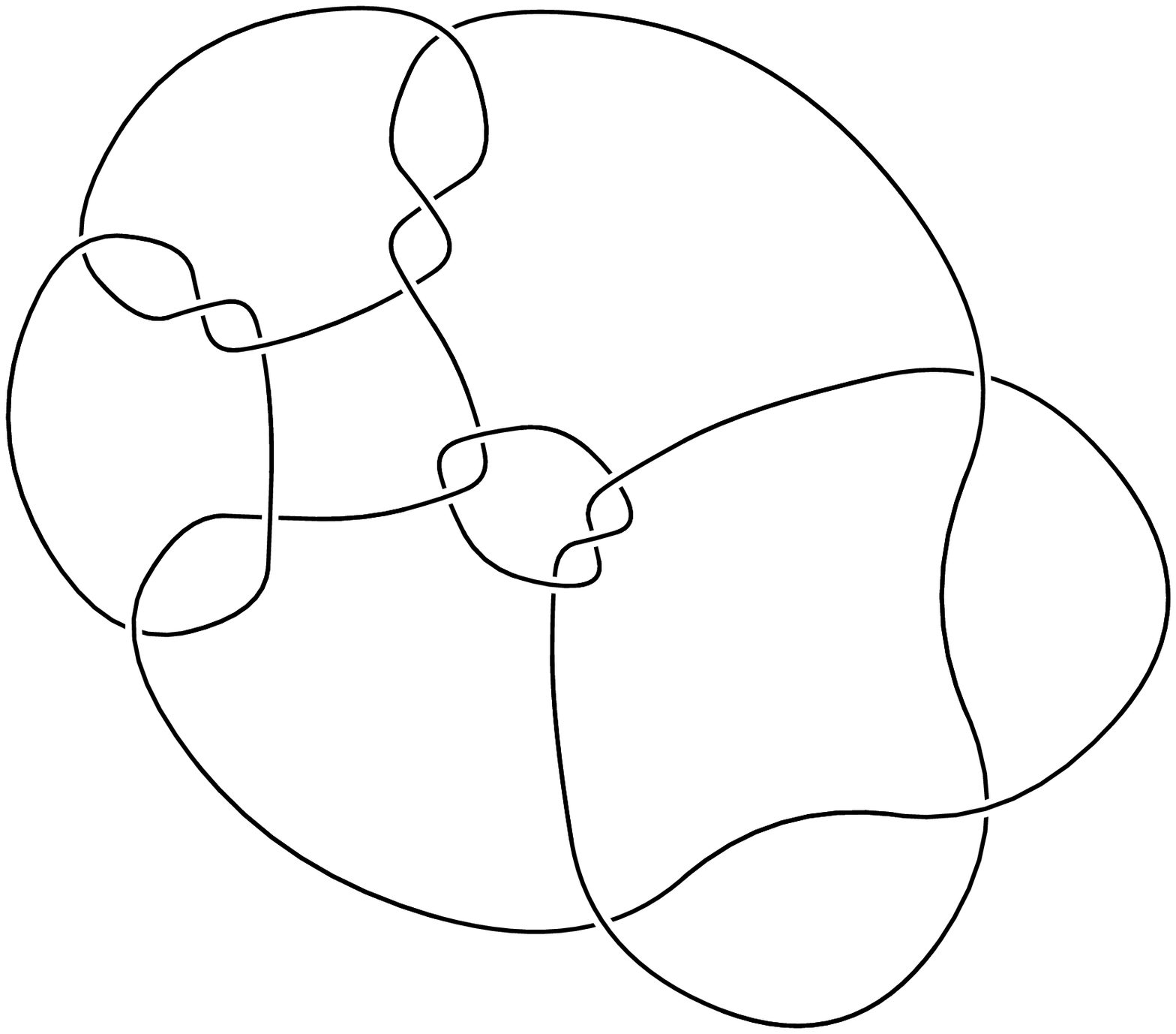}
  \caption{%
    The example of Kanenobu of a chiral-achiral pair of mutant prime knots:
    the chiral $16^A_{259088}$ and the achiral $16^A_{259984}$.
  }%
  \label{figure:16A259088and16A259984}
  \end{centering}
\end{figure}

\begin{figure}[htbp]
  \begin{centering}
  \begin{tabular}{cc@{\hspace{10pt}}|@{\hspace{10pt}}cc}
    \includegraphics[width=75pt]{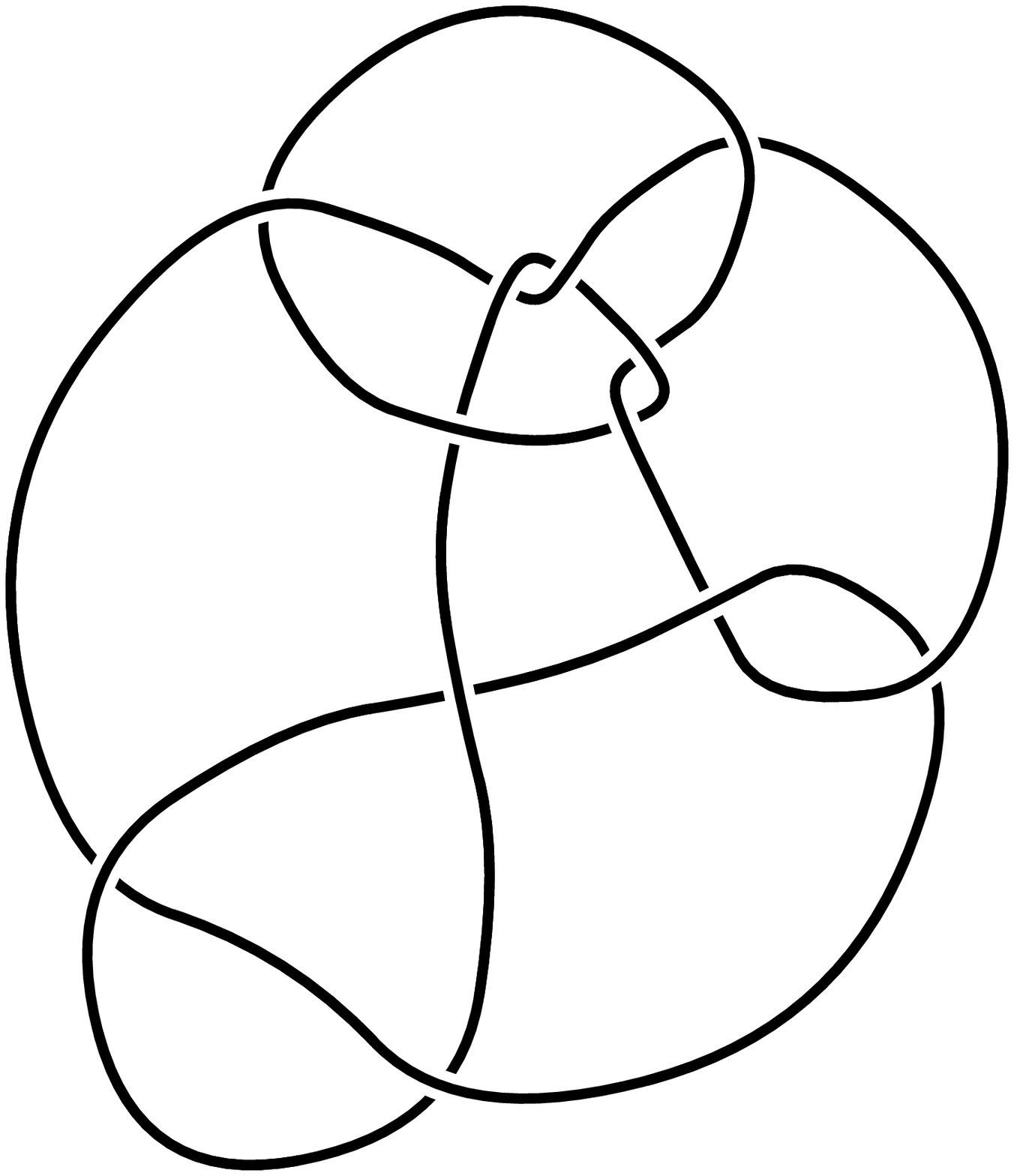}
    &
    \includegraphics[width=75pt]{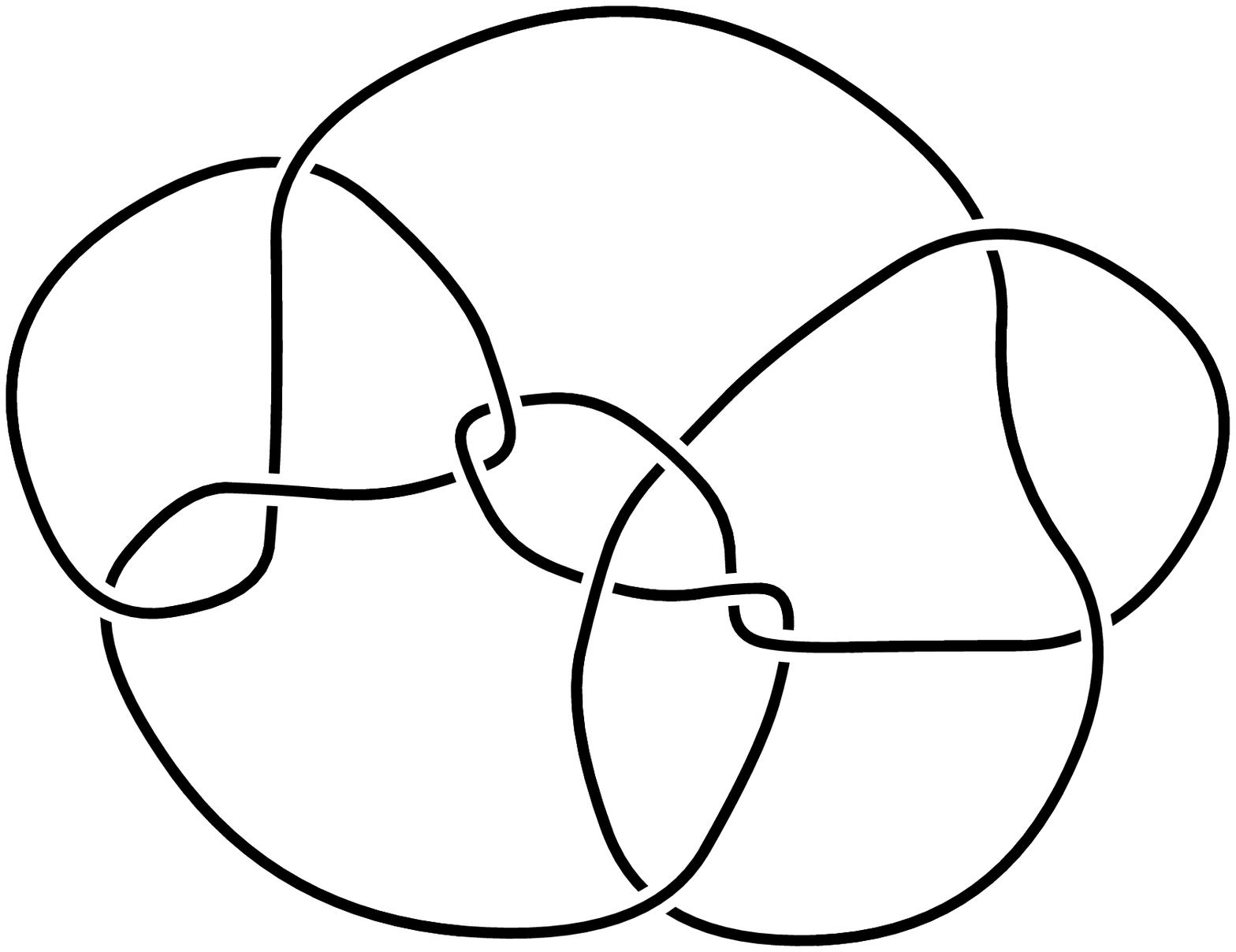}
    &
    \includegraphics[width=75pt]{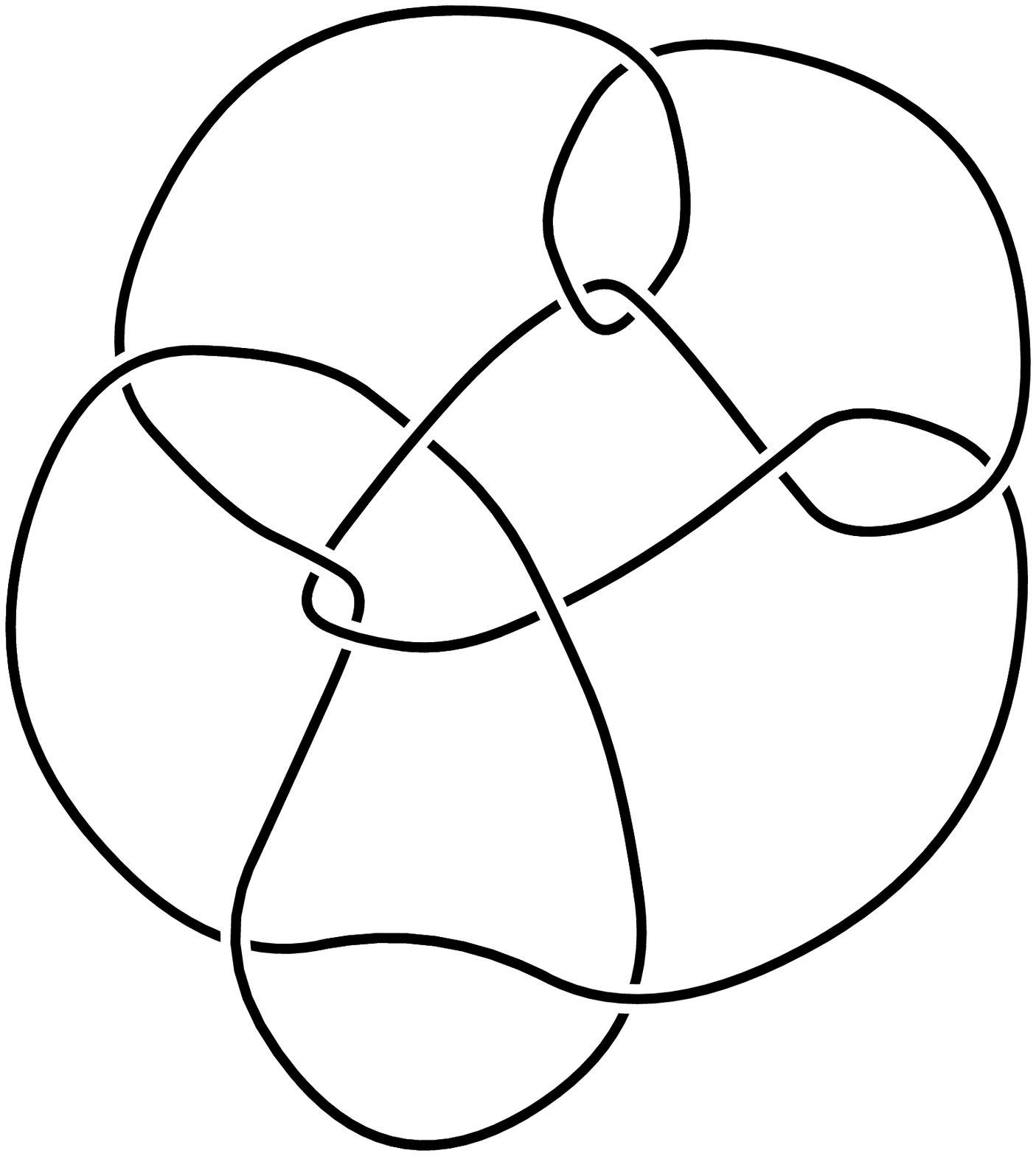}
    &
    \includegraphics[width=75pt]{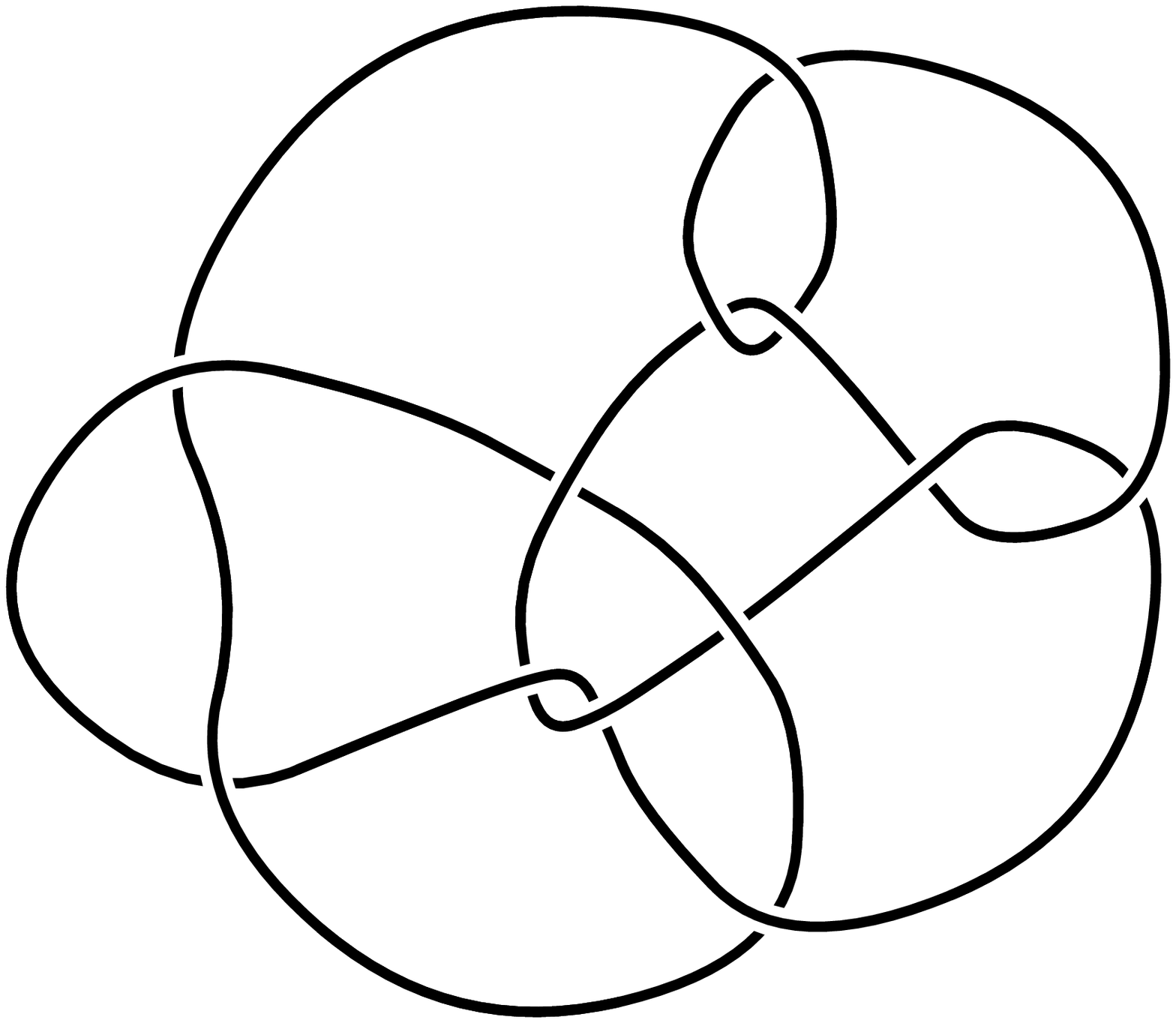}
    \\[-10pt]
    $12^A_{101}$ & $12^A_{115}$ & $12^A_{102}$ & $12^A_{107}$
    \\[10pt]
    \hline
    &&&\\[-10pt]
    \includegraphics[width=75pt]{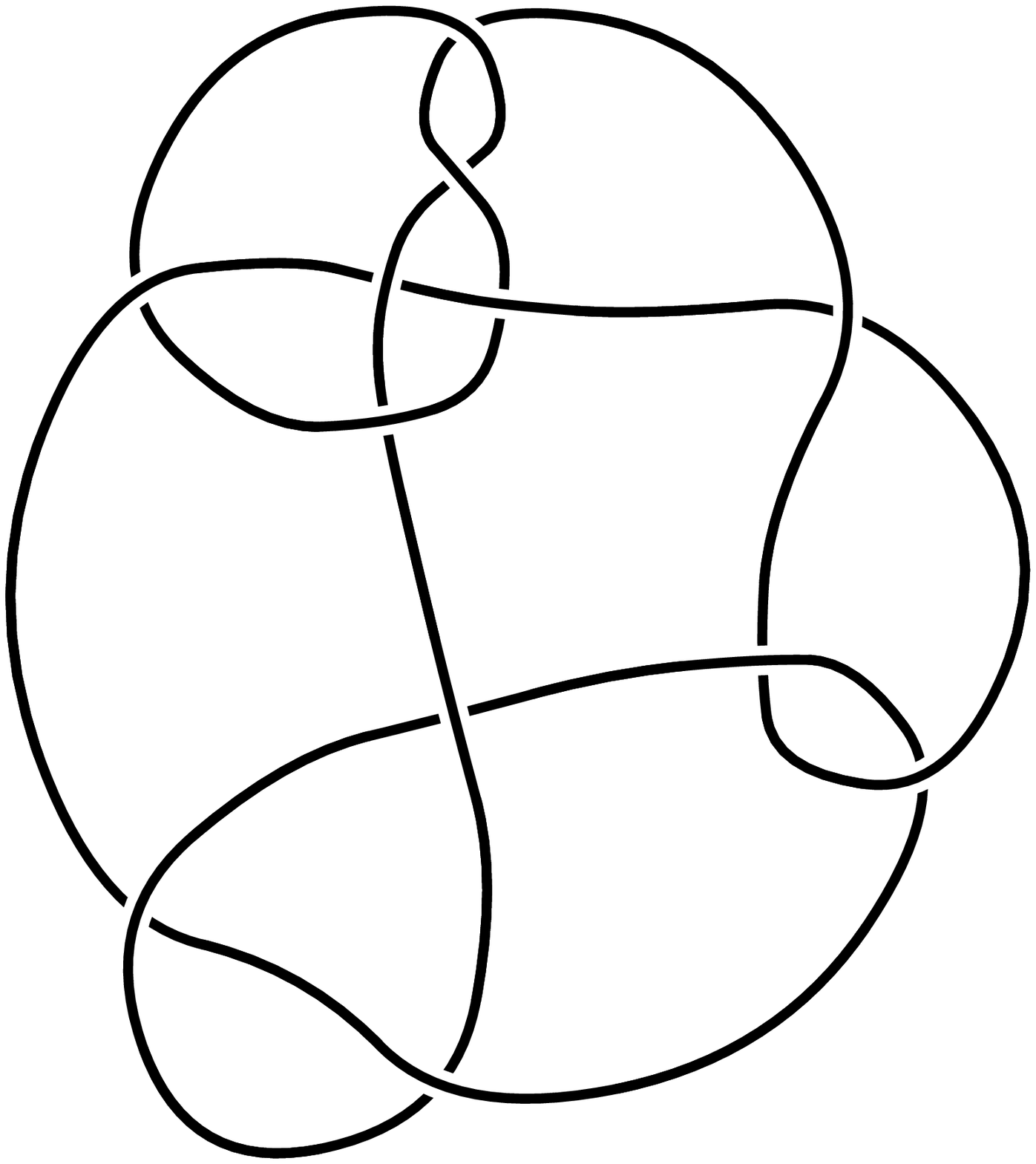}
    &
    \includegraphics[width=75pt]{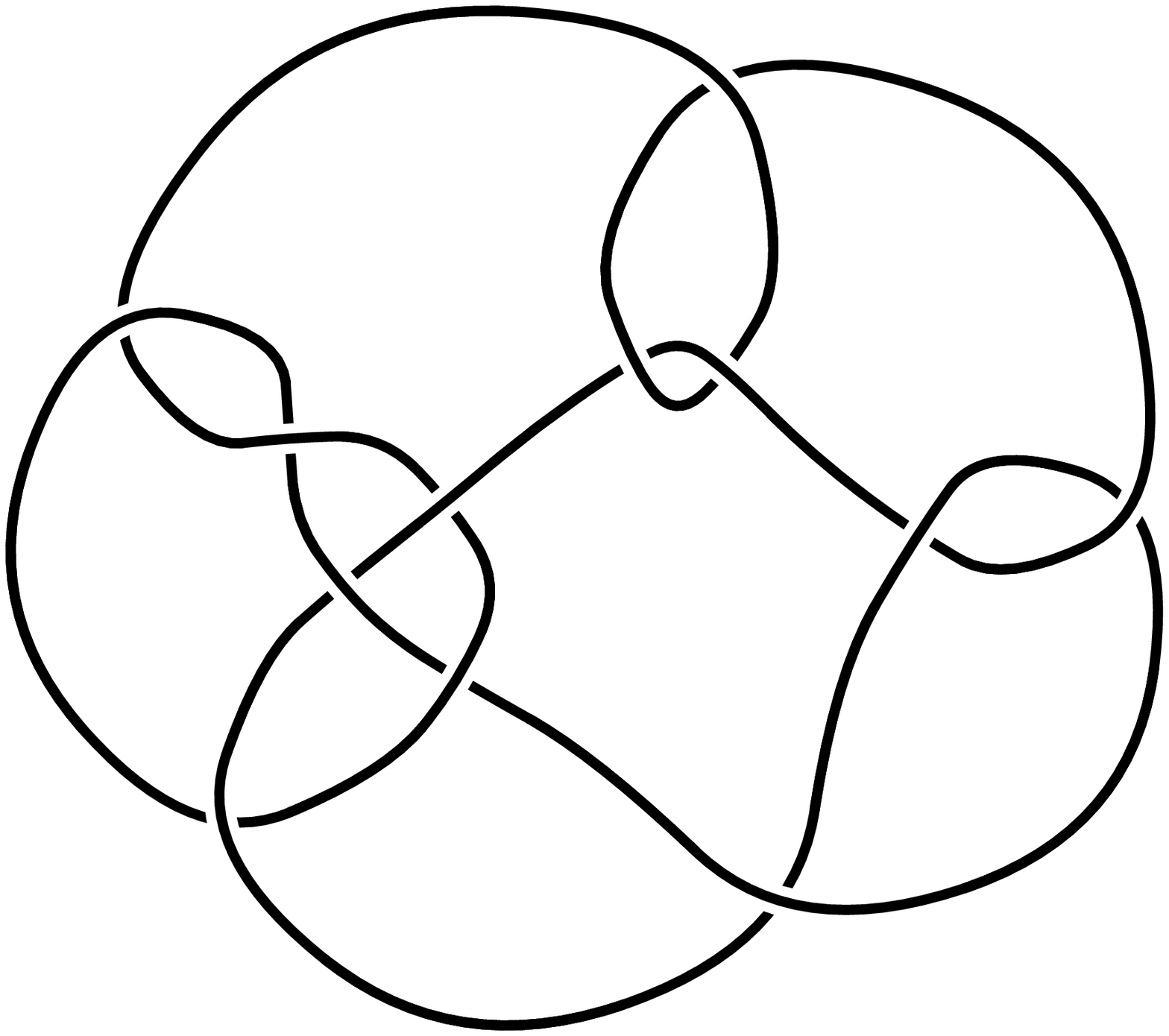}
    &
    \includegraphics[width=75pt]{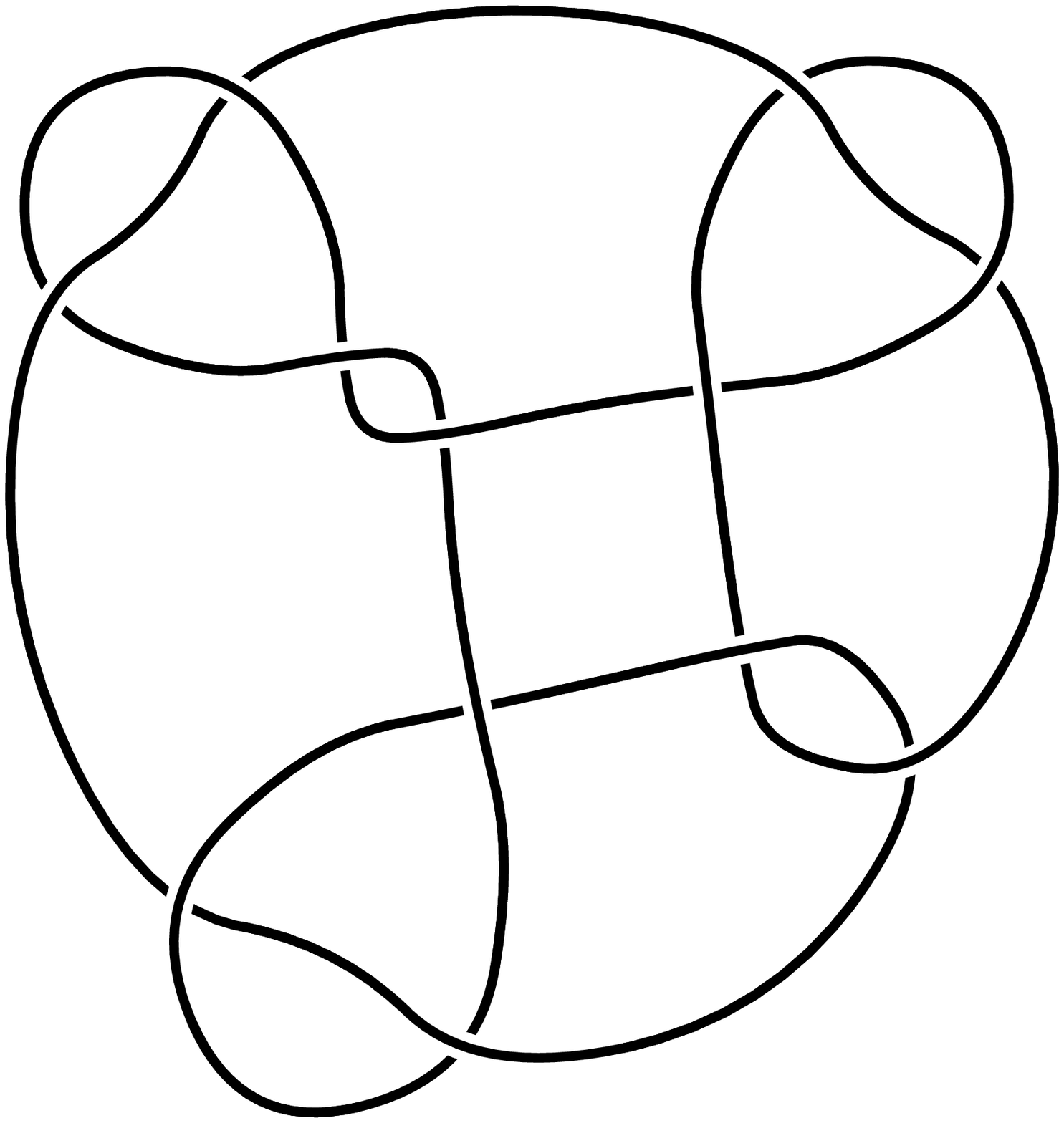}
    &
    \includegraphics[width=75pt]{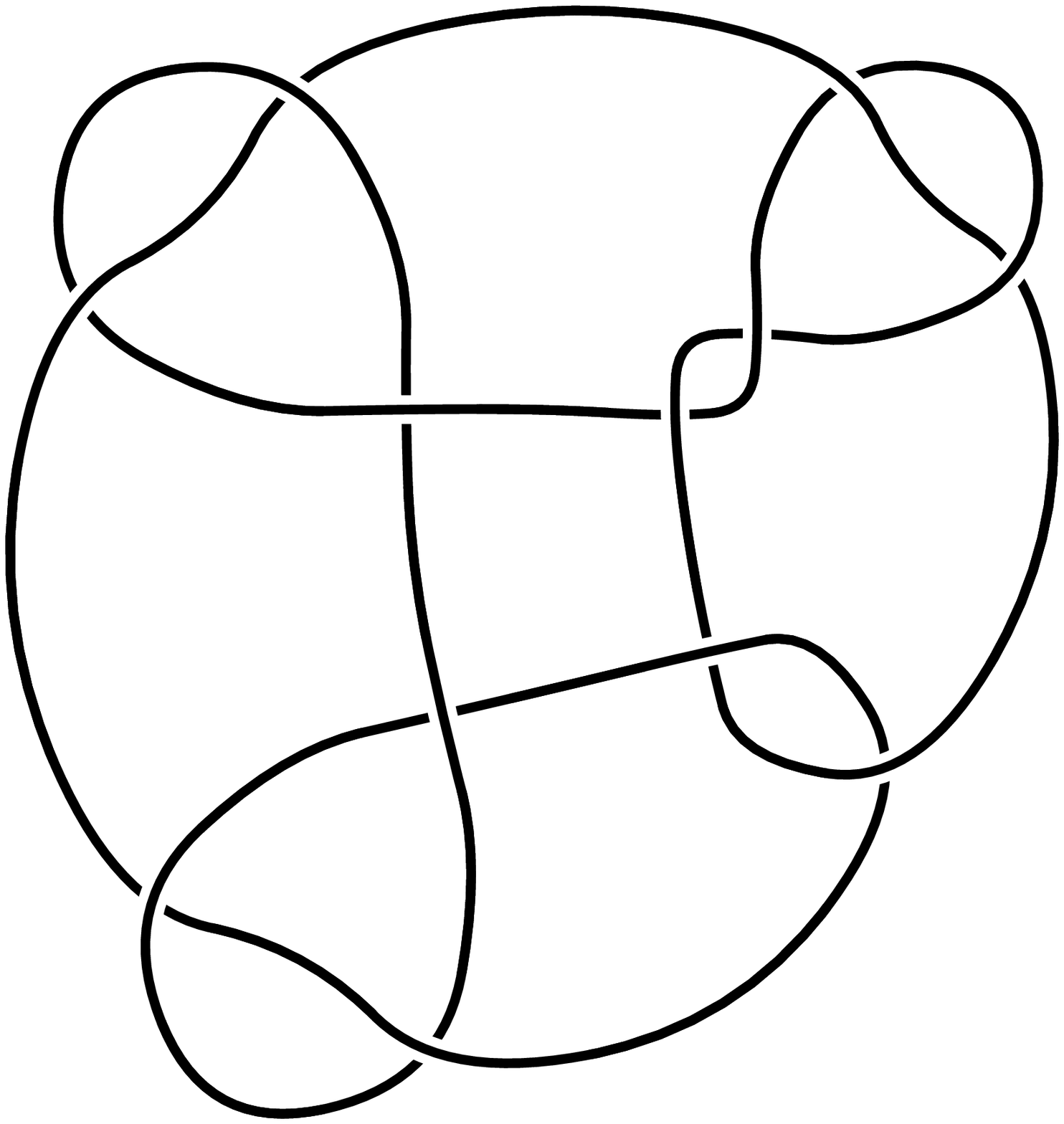}
    \\[-10pt]
    $12^A_{108}$ & $12^A_{120}$ & $12^A_{114}$ & $12^A_{117}$
    \\[10pt]
    \hline
    &&&\\[-10pt]
    \includegraphics[width=75pt]{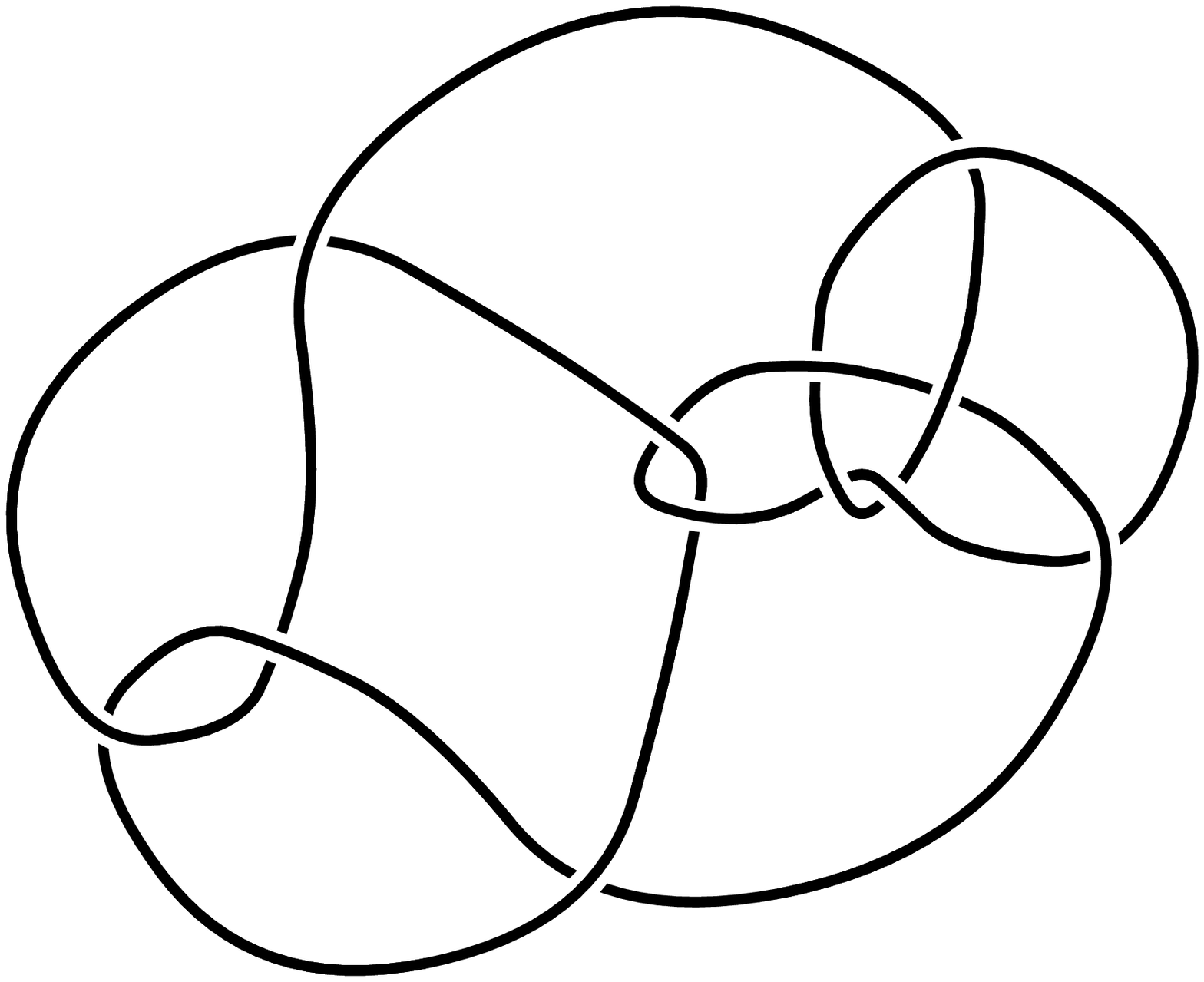}
    &
    \includegraphics[width=75pt]{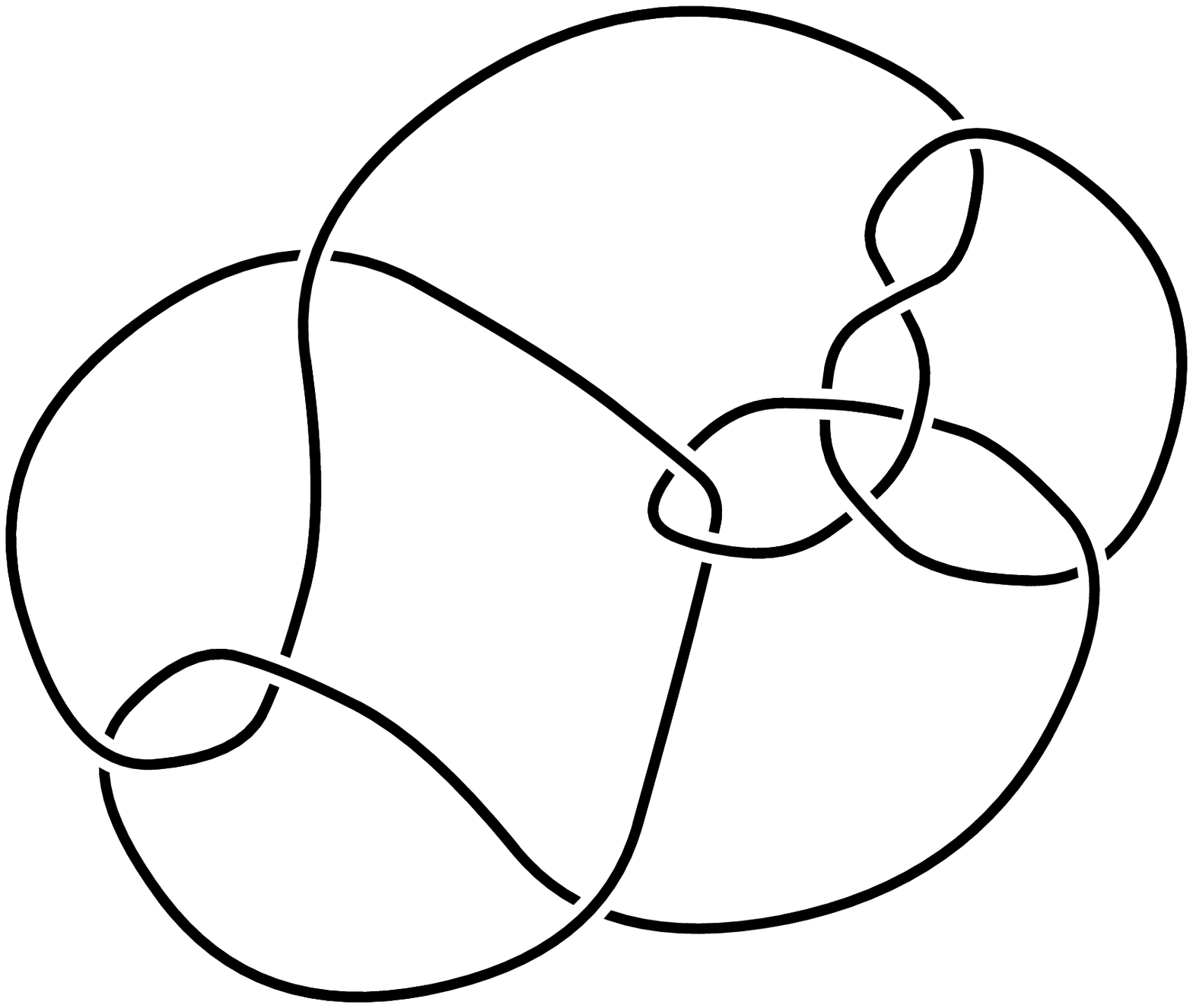}
    &
    \includegraphics[width=75pt]{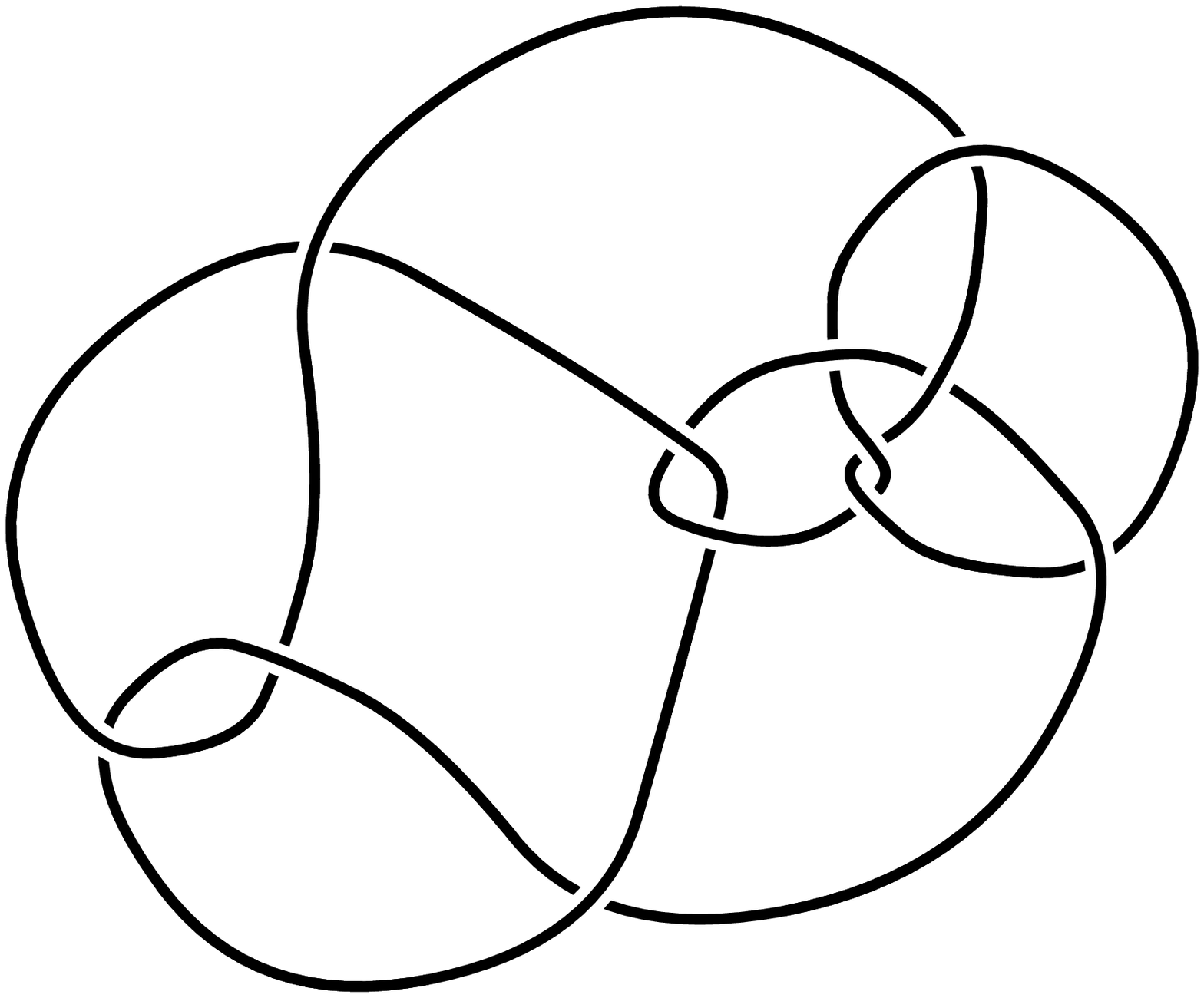}
    &
    \includegraphics[width=75pt]{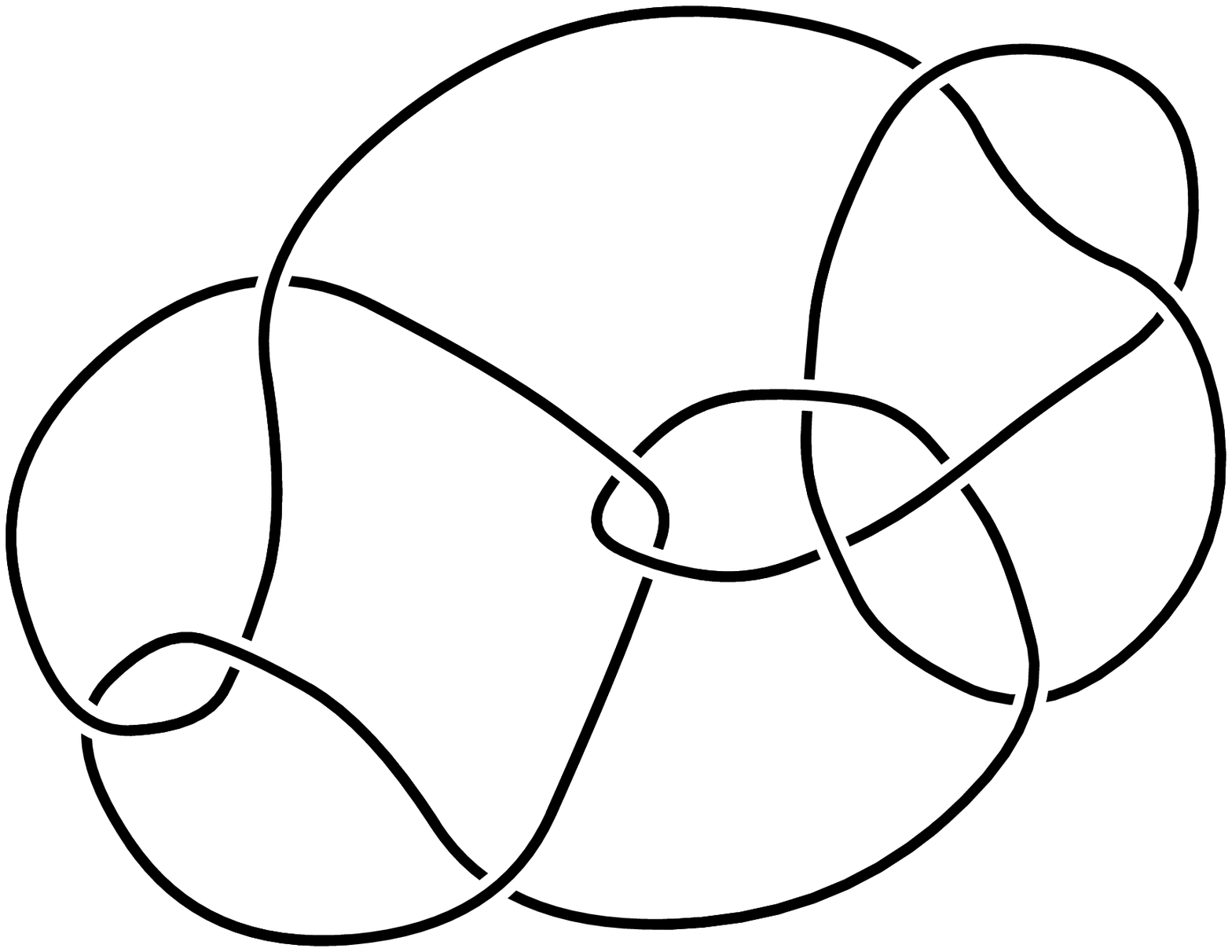}
    \\[-10pt]
    $12^A_{126}$ & $12^A_{132}$ & $12^A_{131}$ & $12^A_{133}$
    \\[10pt]
    \hline
    &&&\\[-10pt]
    \includegraphics[width=75pt]{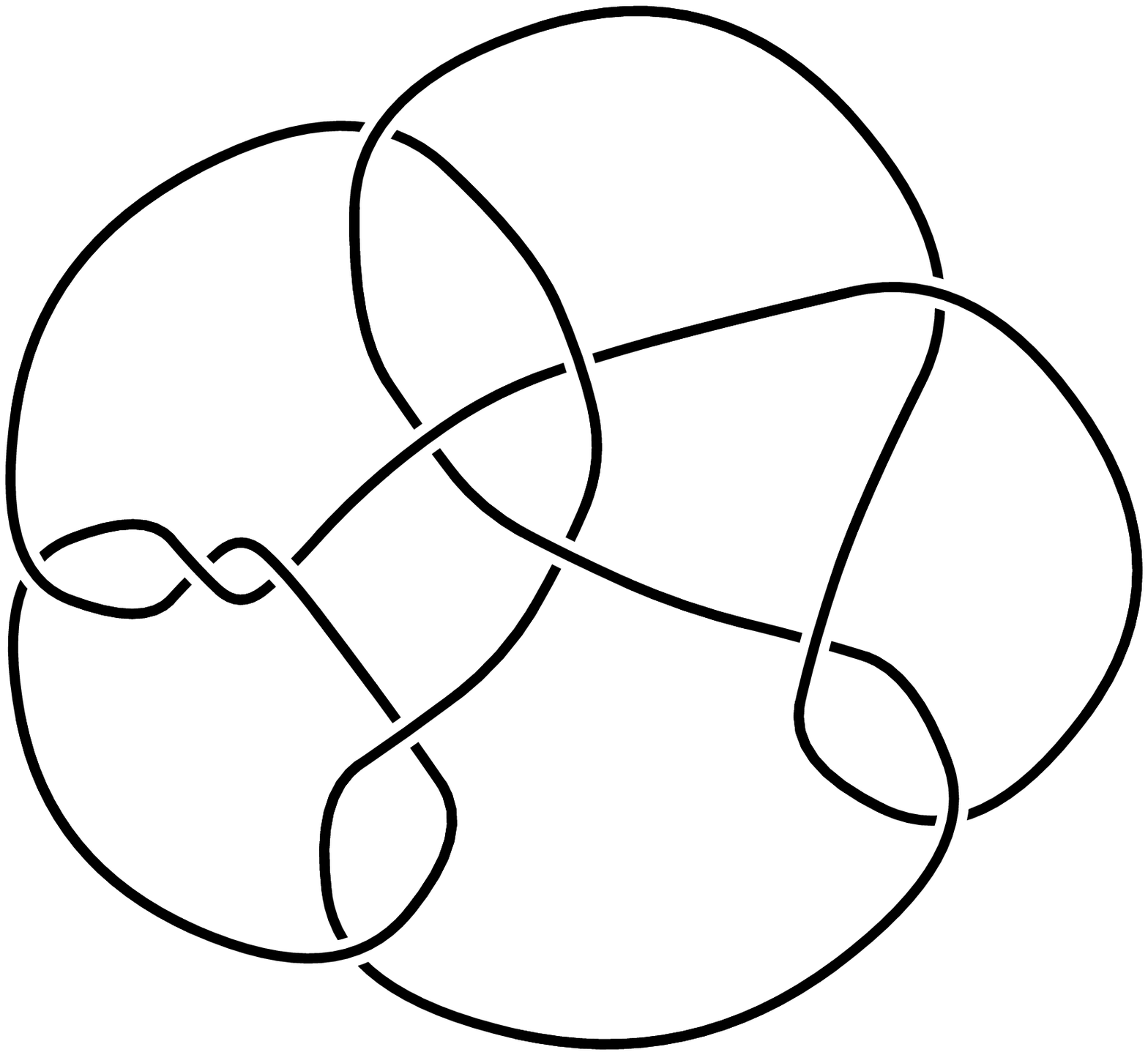}
    &
    \includegraphics[width=75pt]{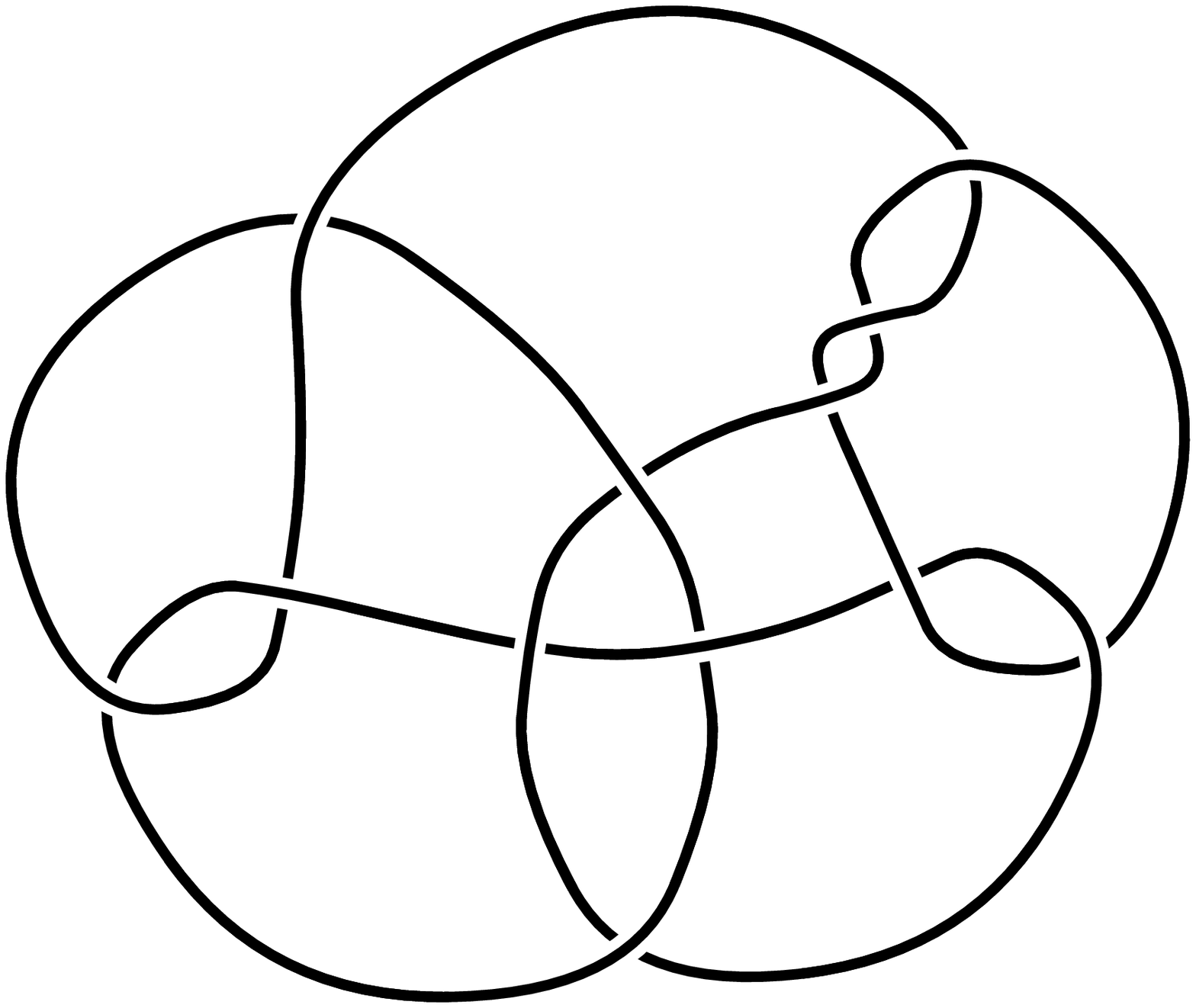}
    &
    \includegraphics[width=75pt]{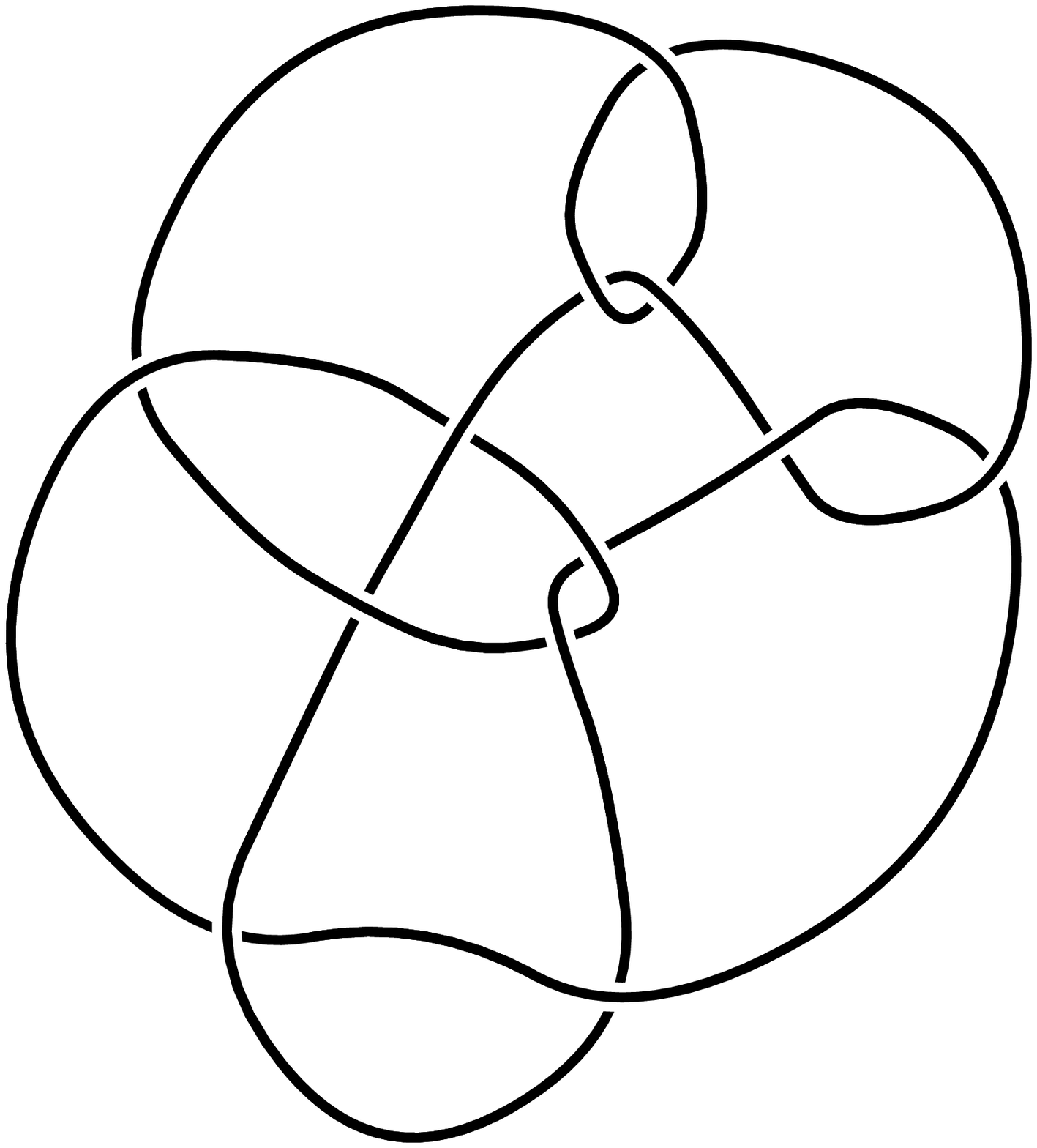}
    &
    \includegraphics[width=75pt]{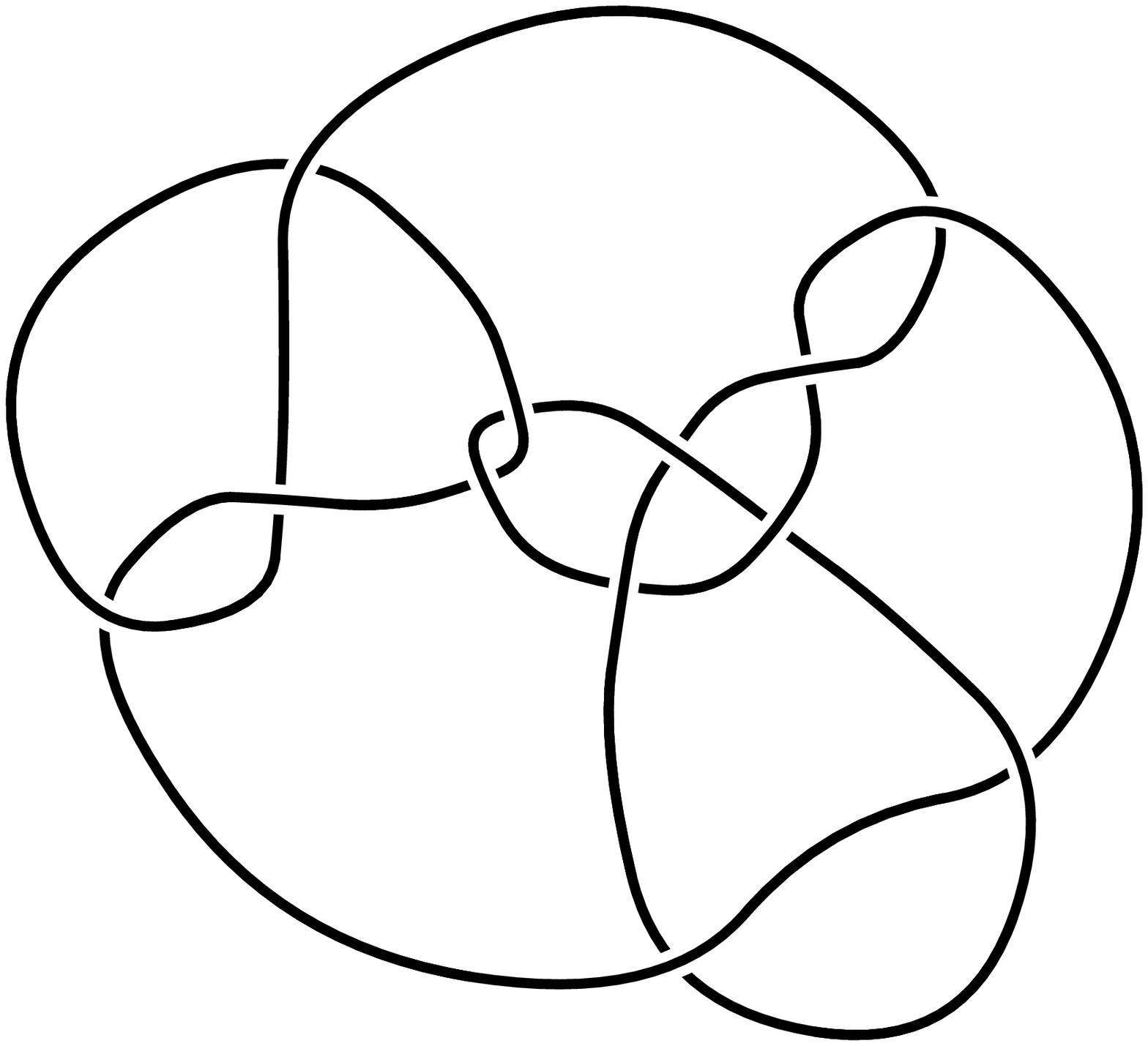}
    \\[-10pt]
    $12^A_{134}$ & $12^A_{188}$ & $12^A_{154}$ & $12^A_{162}$
    \\[10pt]
    \hline
    &&&\\[-10pt]
    \includegraphics[width=75pt]{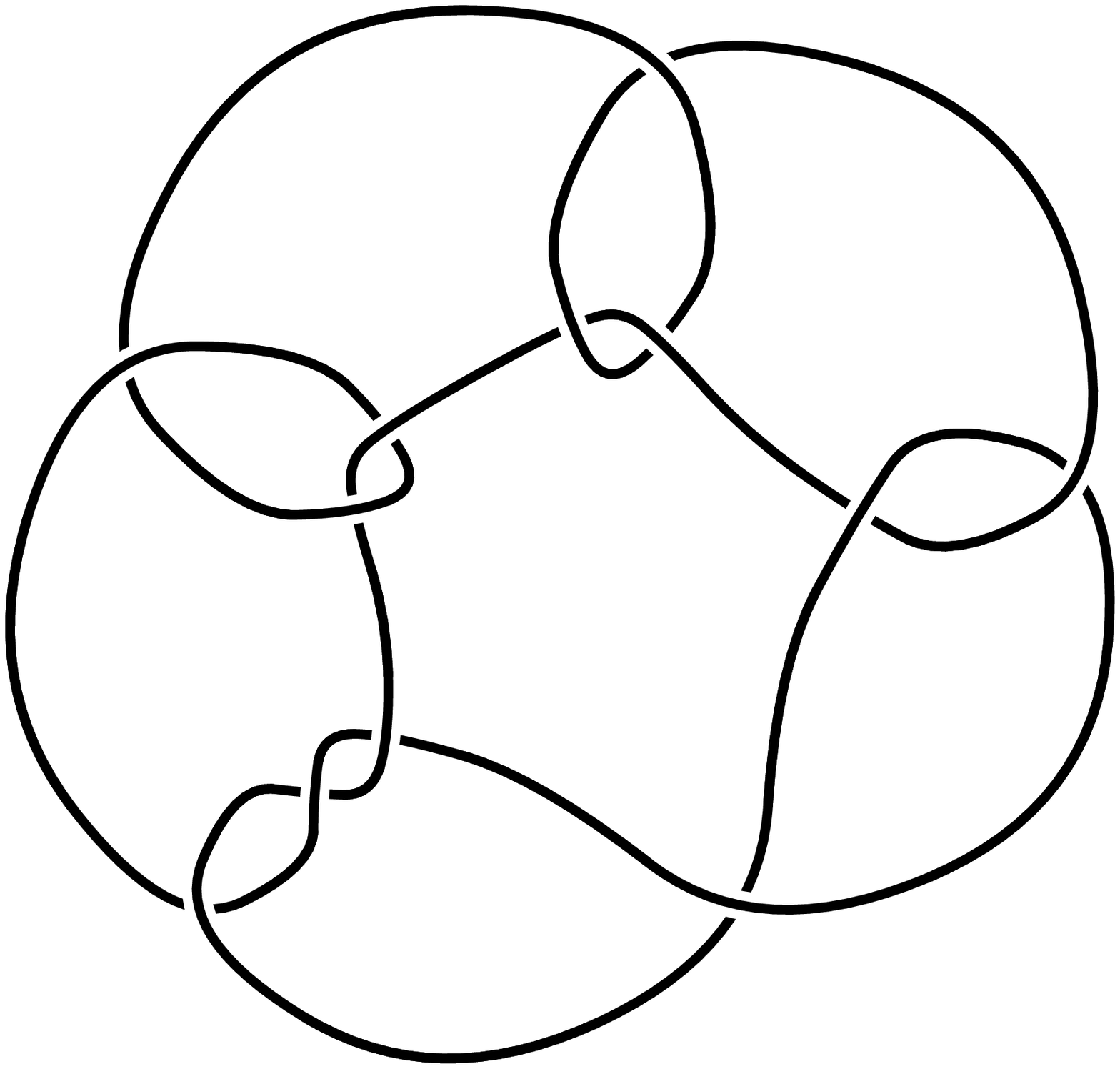}
    &
    \includegraphics[width=75pt]{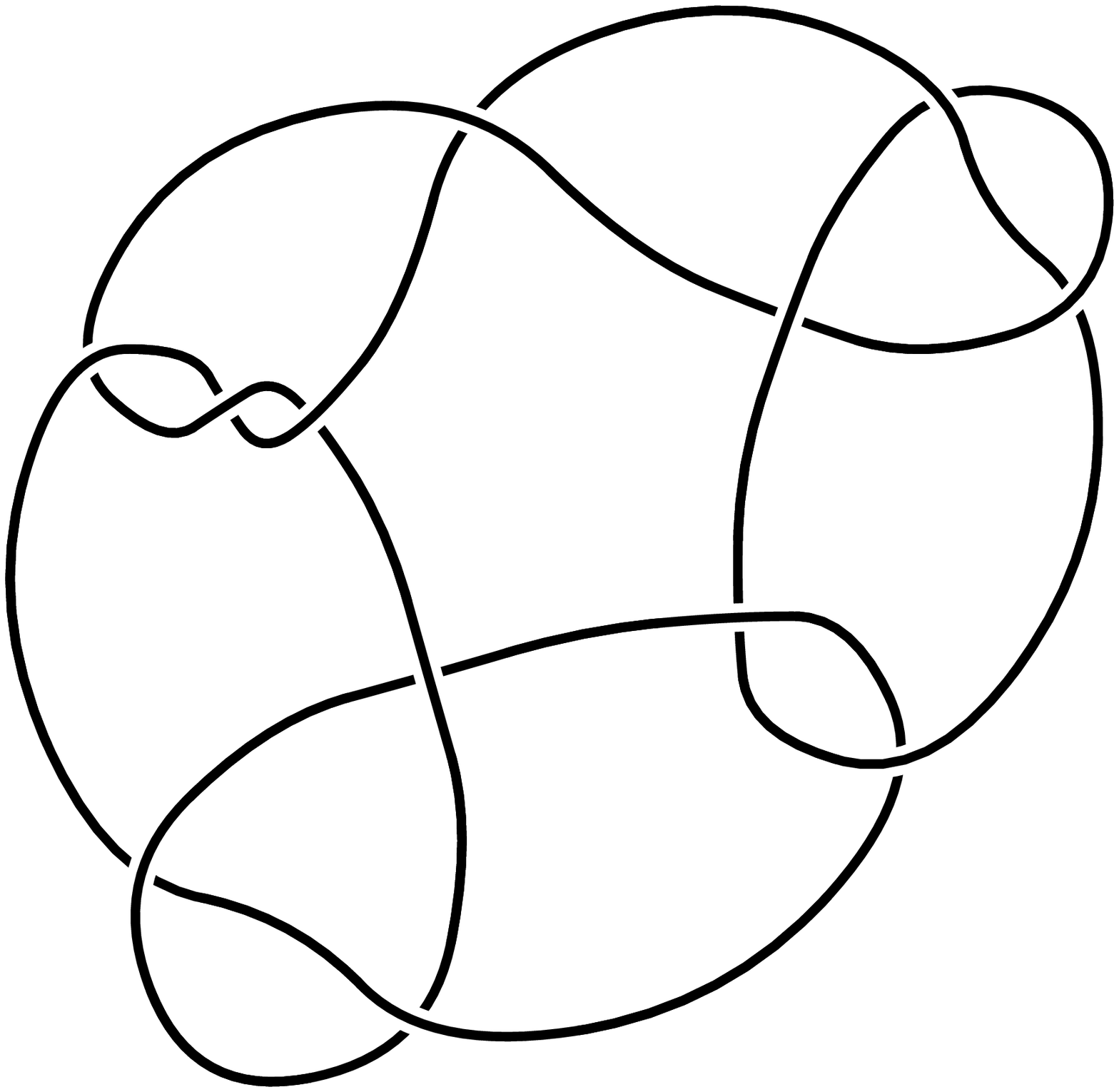}
    &
    \includegraphics[width=75pt]{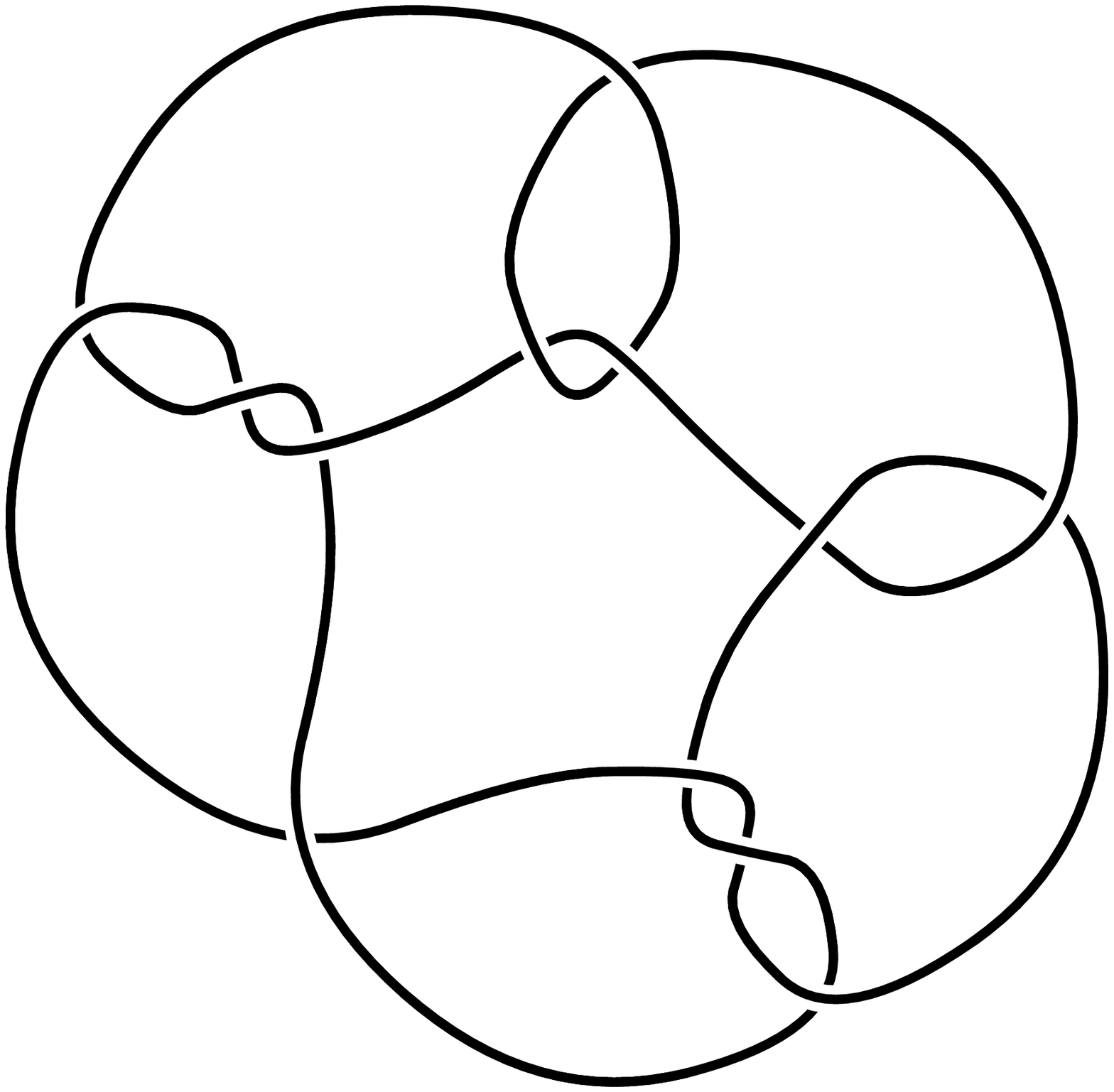}
    &
    \includegraphics[width=75pt]{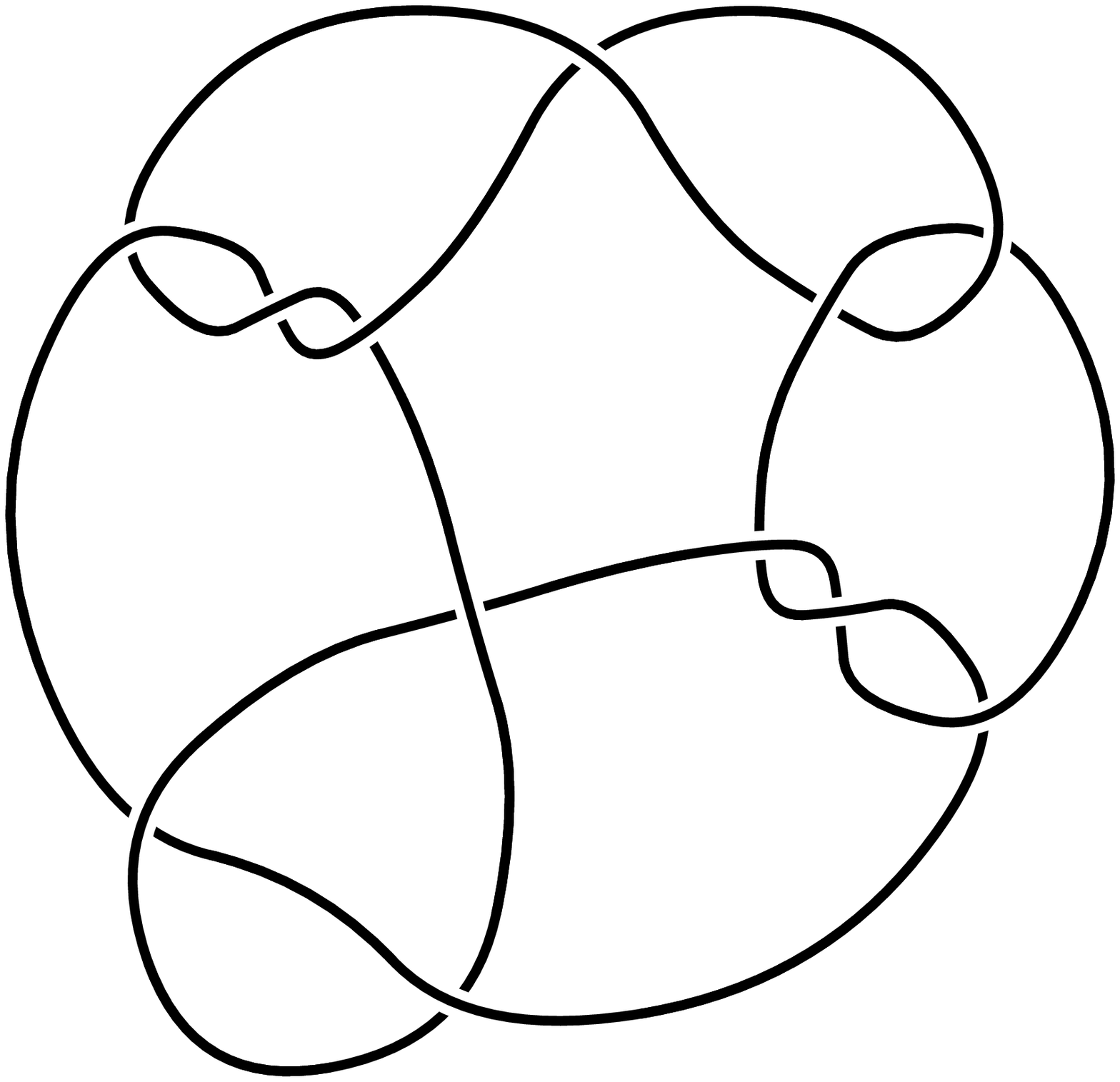}
    \\[-10pt]
    $12^A_{164}$ & $12^A_{166}$ & $12^A_{167}$ & $12^A_{692}$
  \end{tabular}
  \caption{Alternating $12$-crossing mutant cliques 2/4}
  \end{centering}
\end{figure}

\begin{figure}[htbp]
  \begin{centering}
  \begin{tabular}{cc@{\hspace{10pt}}|@{\hspace{10pt}}cc}
    \includegraphics[width=75pt]{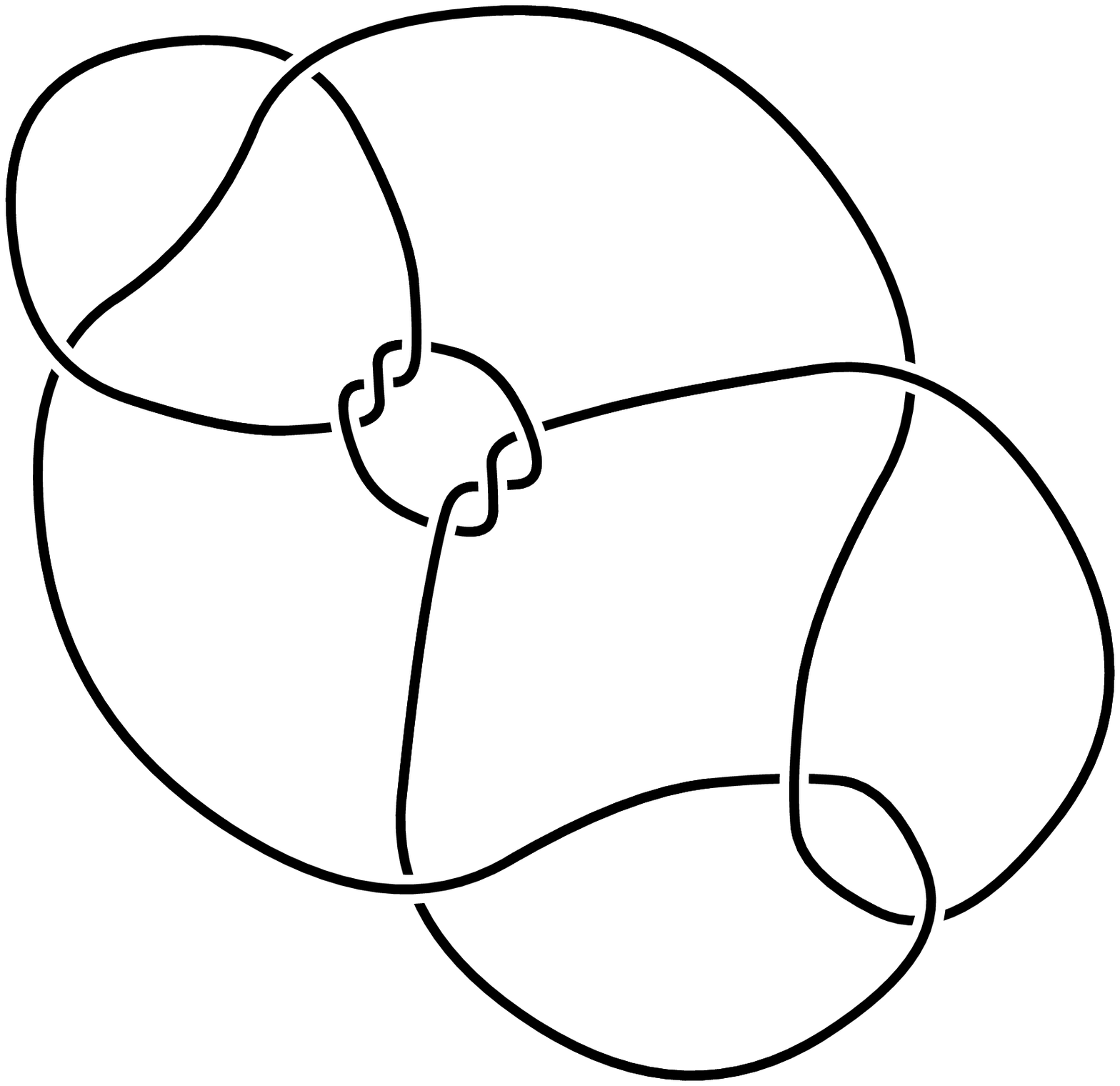}
    &
    \includegraphics[width=75pt]{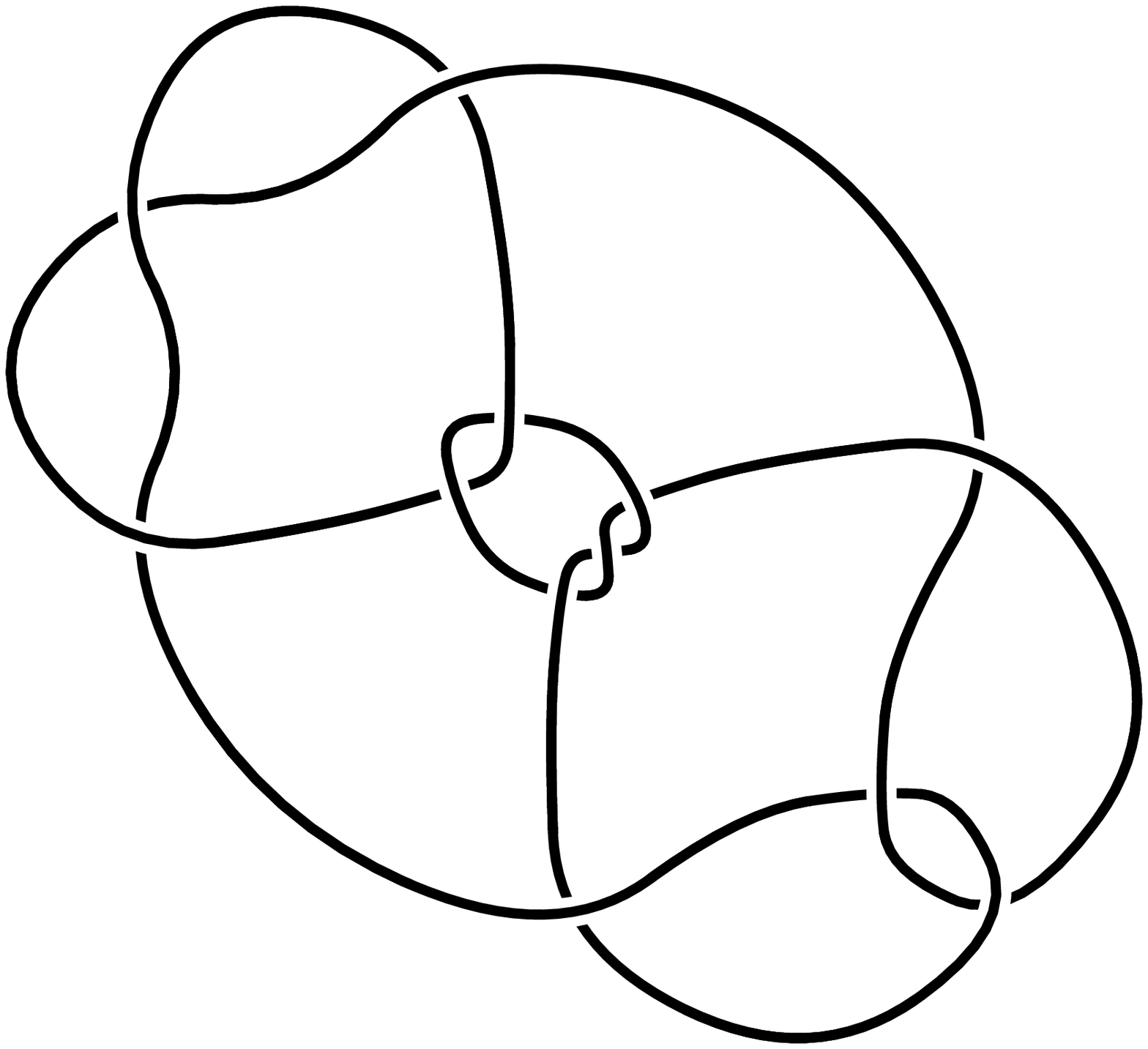}
    &
    \includegraphics[width=75pt]{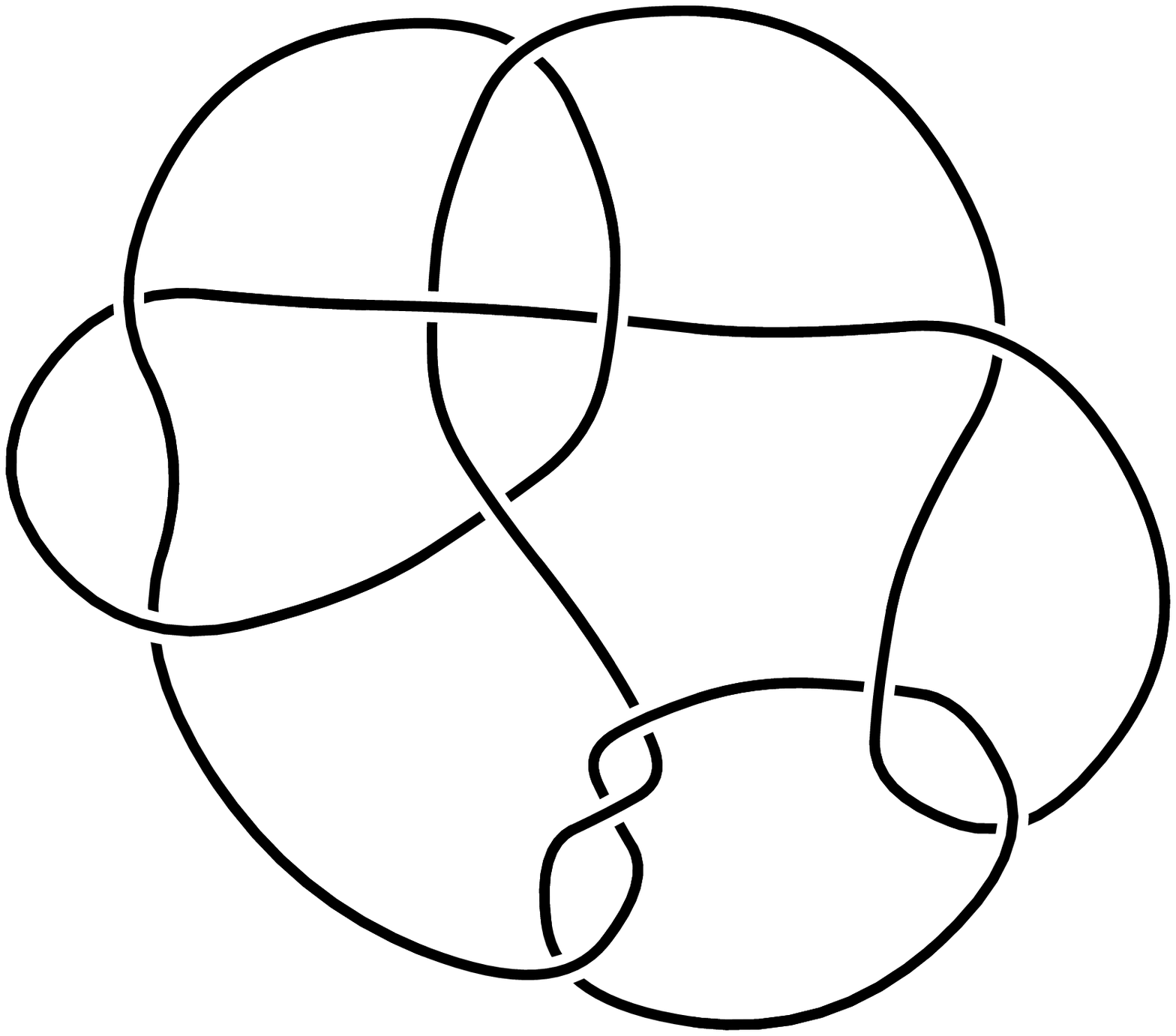}
    &
    \includegraphics[width=75pt]{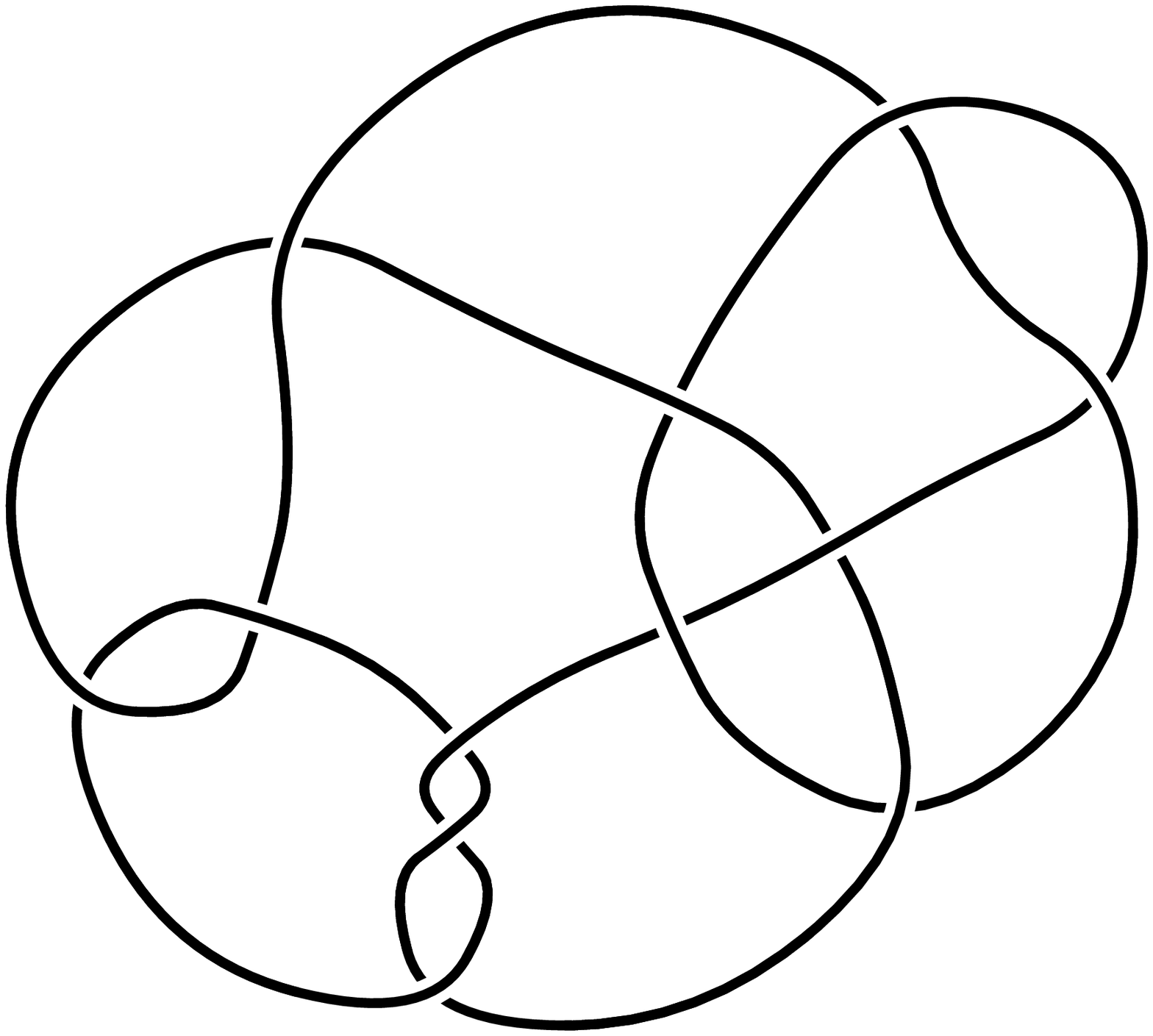}
    \\[-10pt]
    $12^A_{195}$ & $12^A_{693}$ & $12^A_{639}$ & $12^A_{680}$
    \\[10pt]
    \hline
    &&&\\[-10pt]
    \includegraphics[width=75pt]{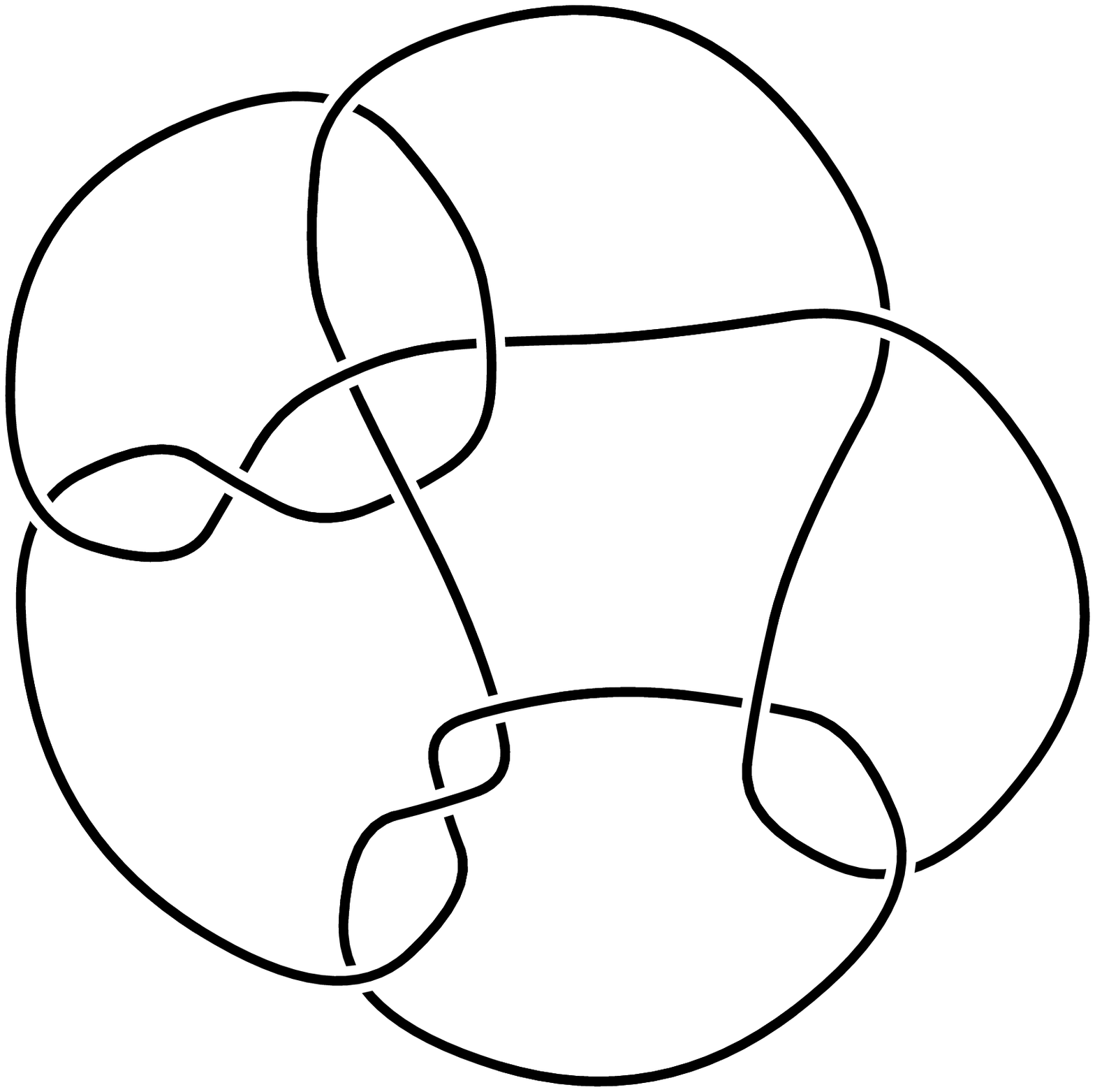}
    &
    \includegraphics[width=75pt]{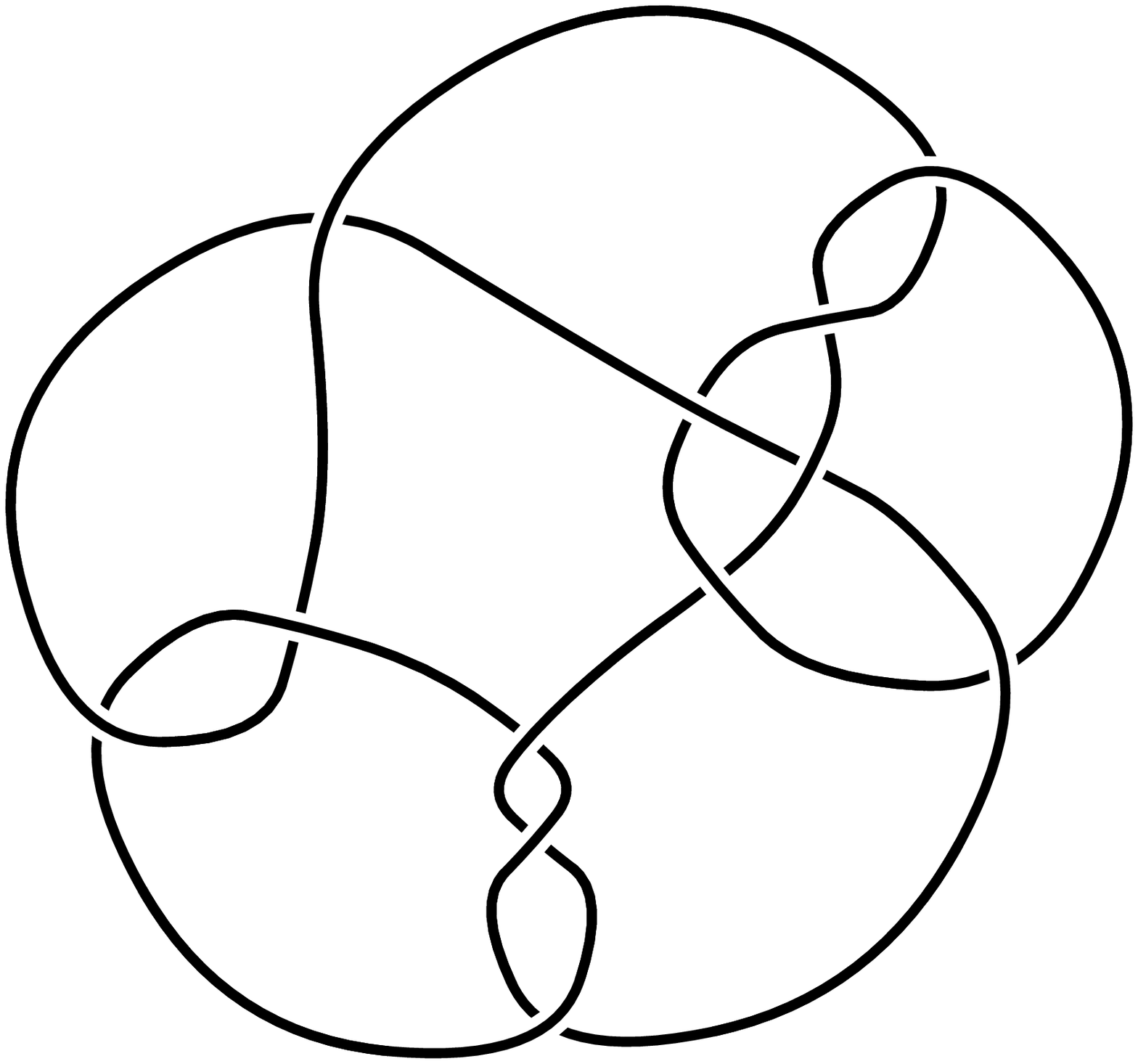}
    &
    \includegraphics[width=75pt]{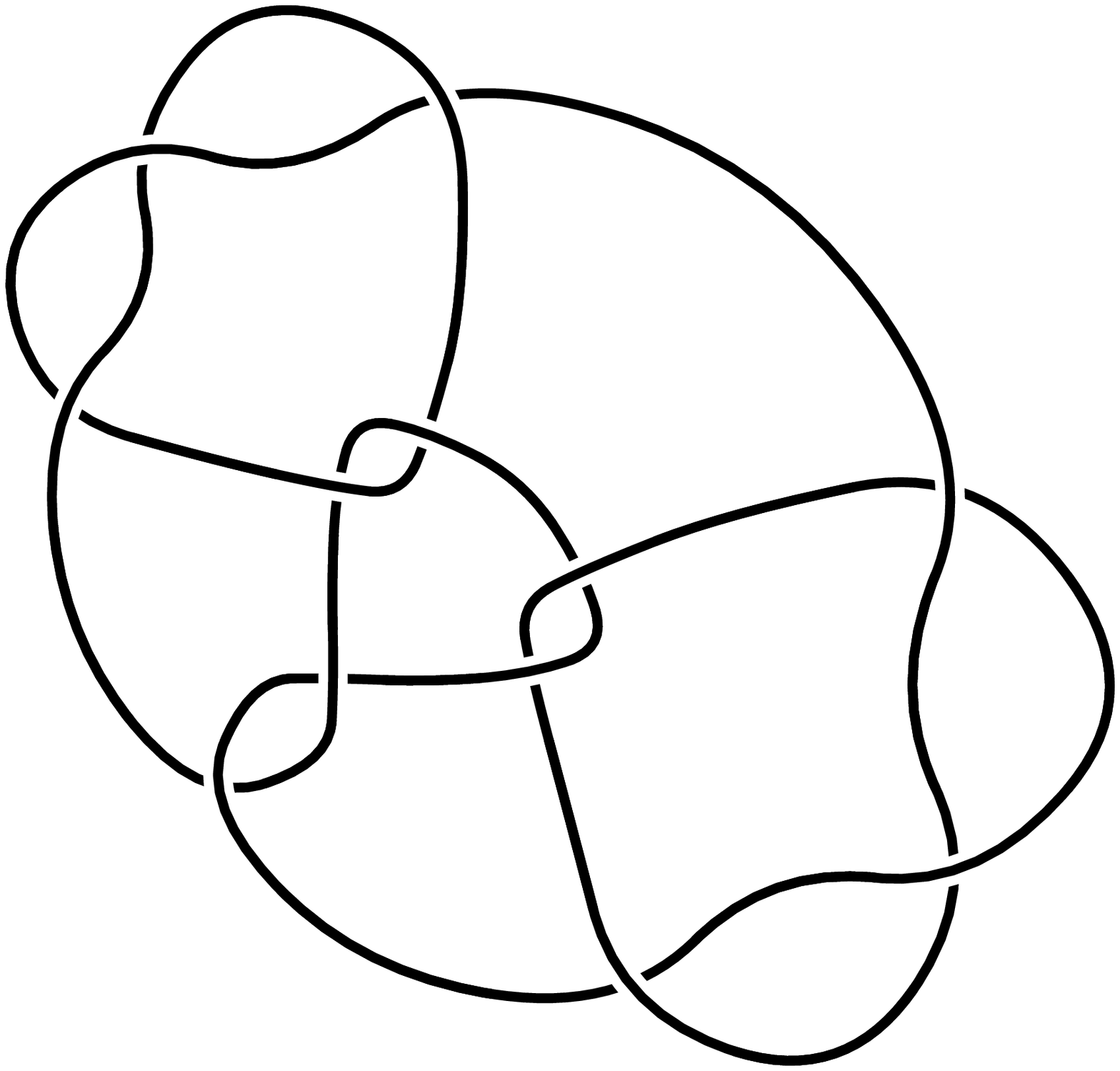}
    &
    \includegraphics[width=75pt]{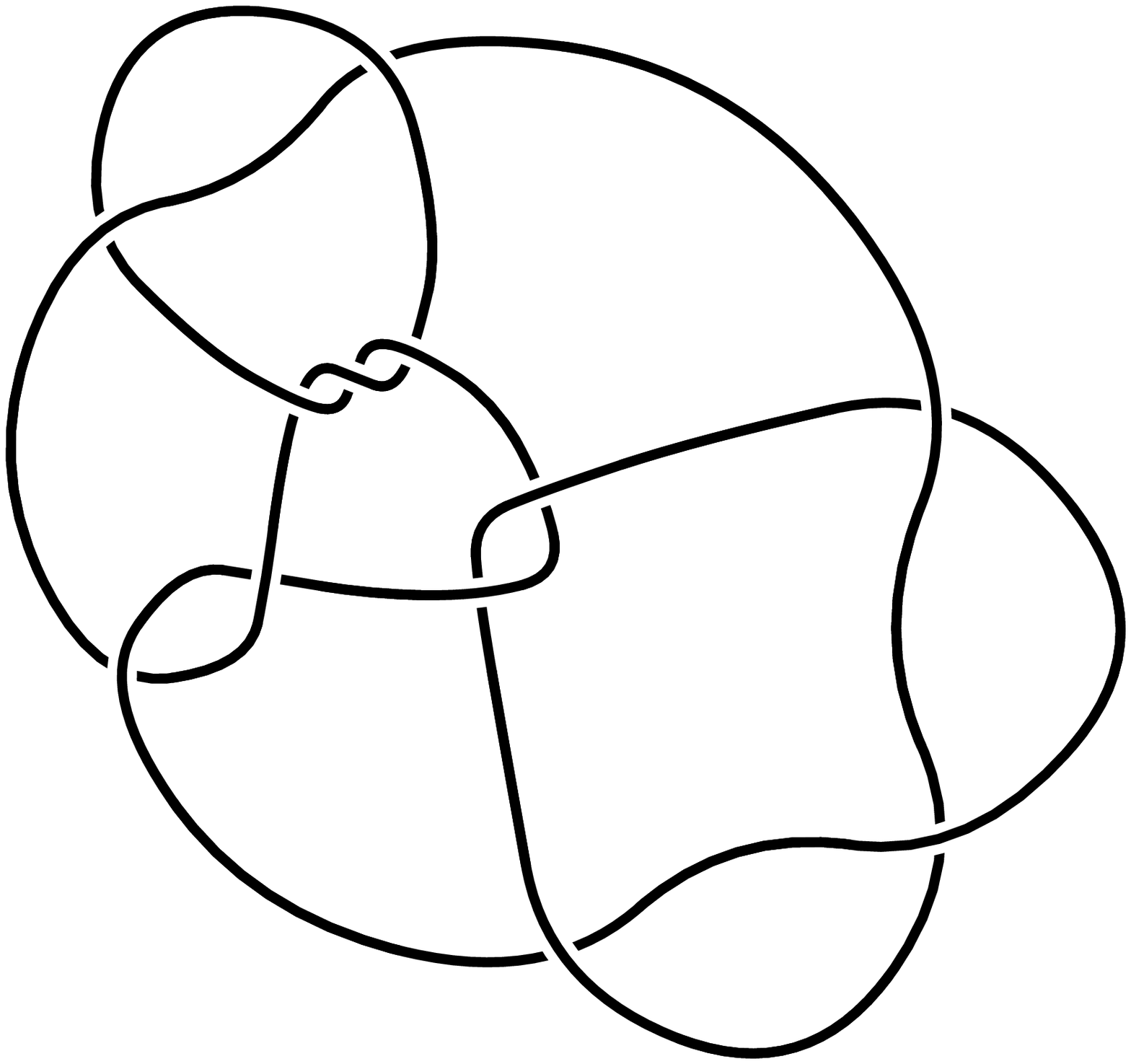}
    \\[-10pt]
    $12^A_{675}$ & $12^A_{688}$ & $12^A_{811}$ & $12^A_{817}$
    \\[10pt]
    \hline
    &&&\\[-10pt]
    \includegraphics[width=75pt]{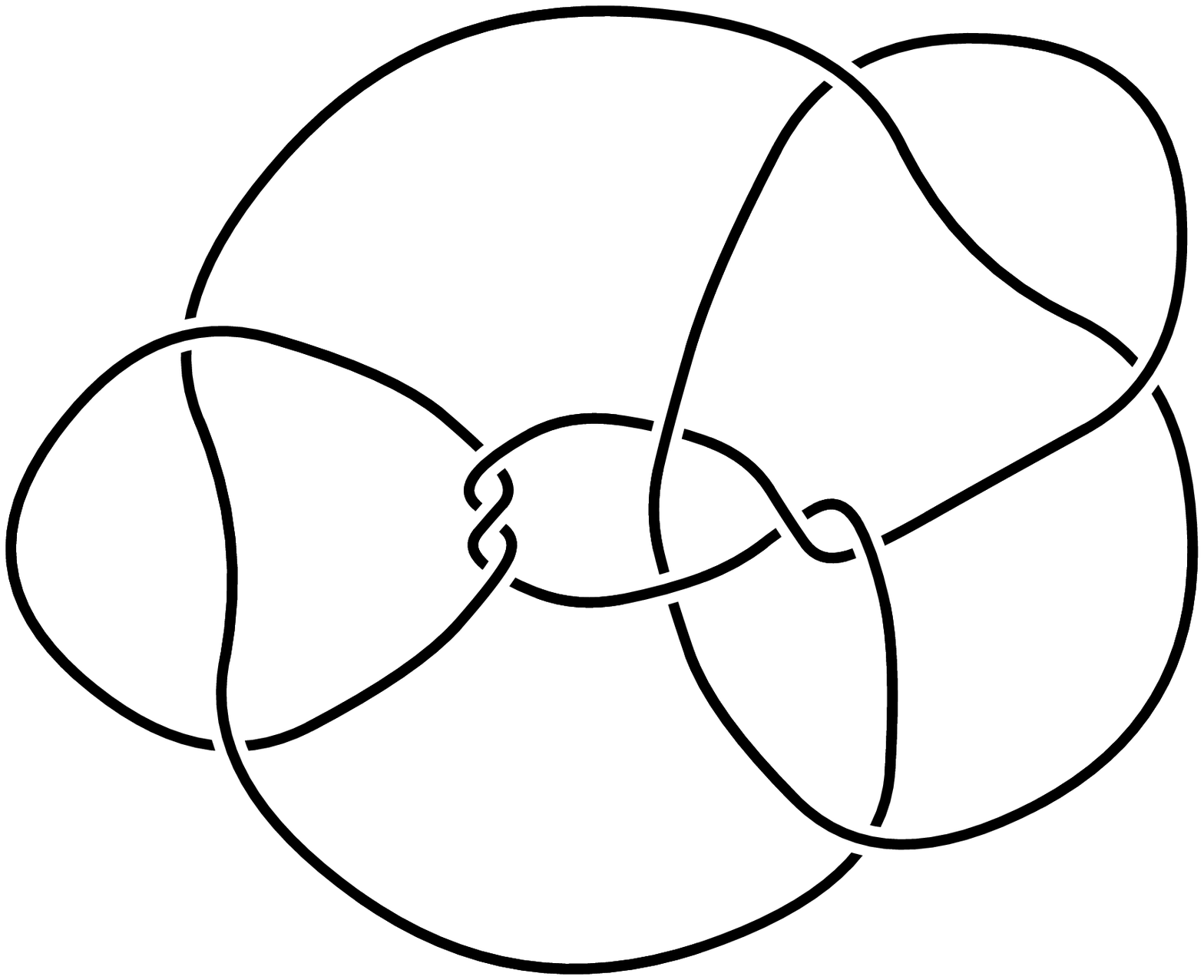}
    &
    \includegraphics[width=75pt]{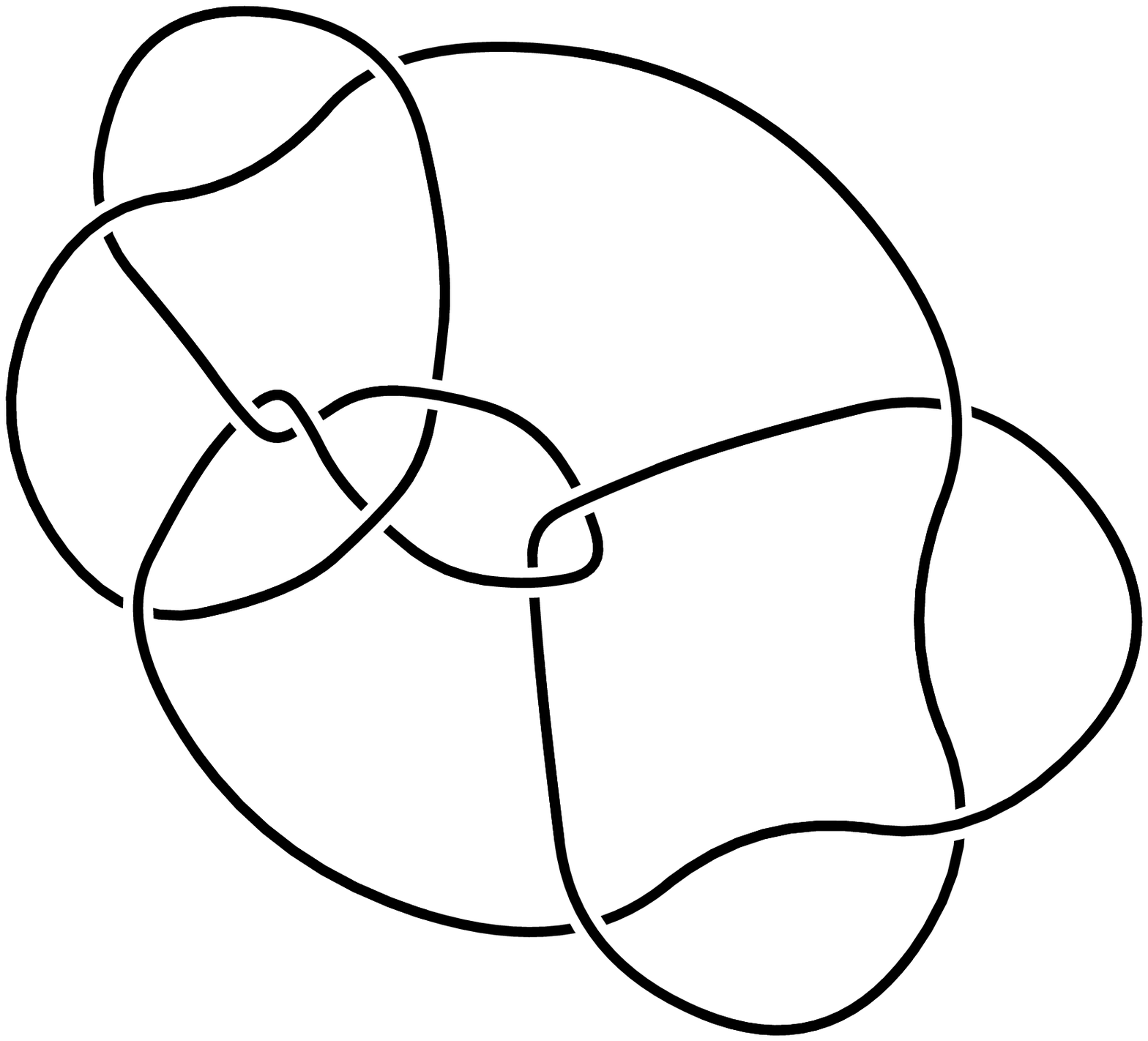}
    &
    \includegraphics[width=75pt]{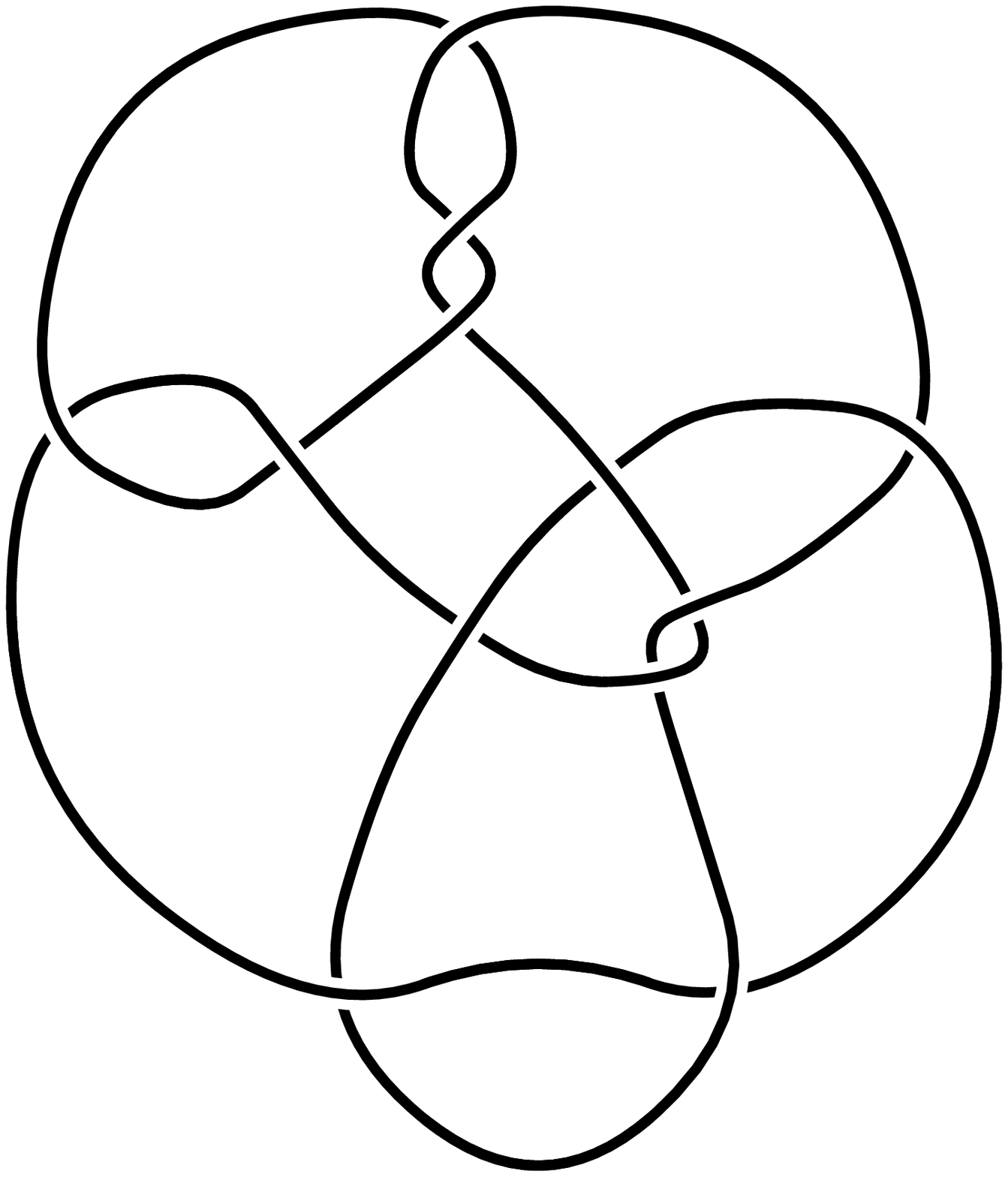}
    &
    \includegraphics[width=75pt]{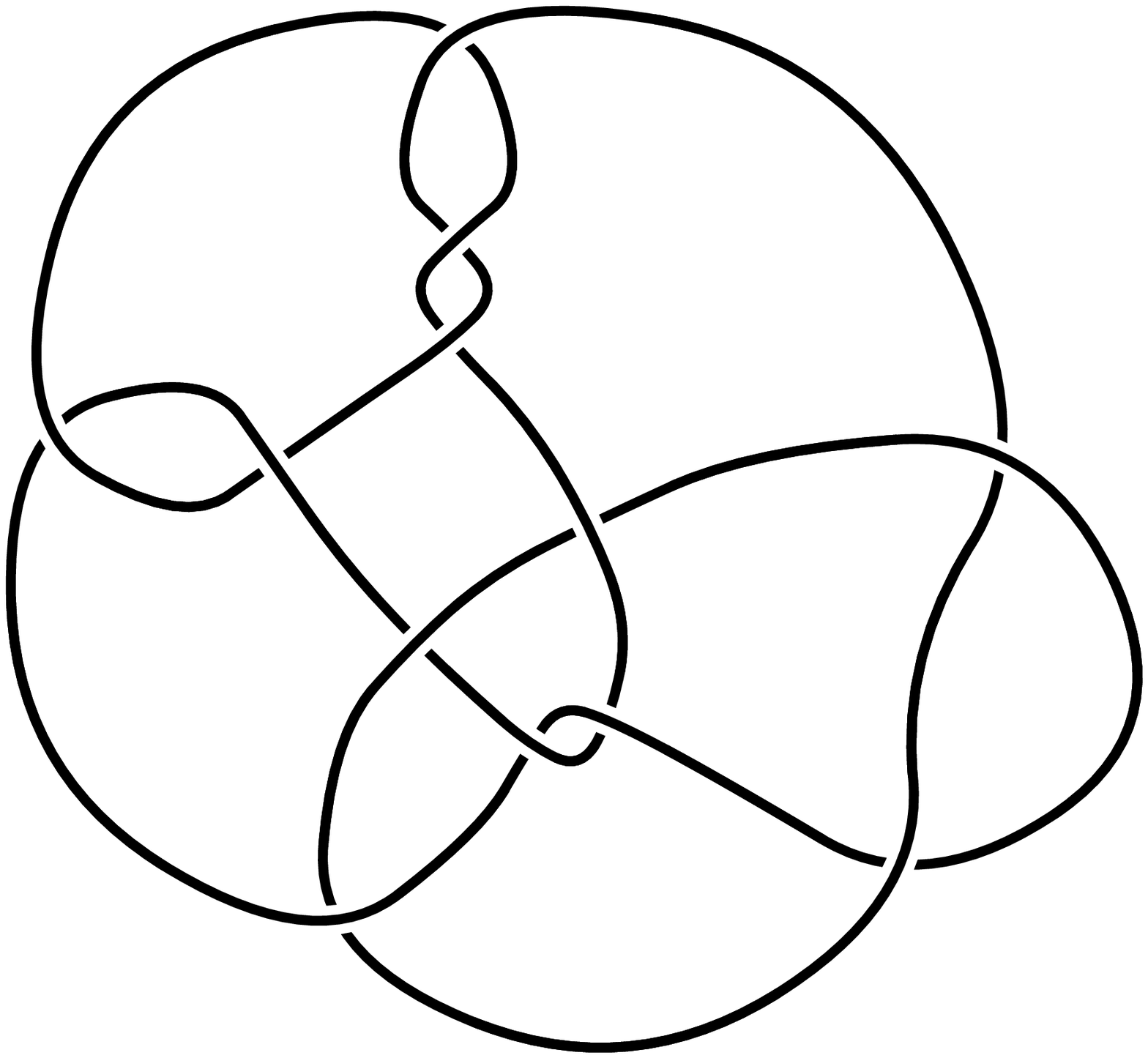}
    \\[-10pt]
    $12^A_{829}$ & $12^A_{832}$ & $12^A_{830}$ & $12^A_{831}$
    \\[10pt]
    \hline
    &&&\\[-10pt]
    \includegraphics[width=75pt]{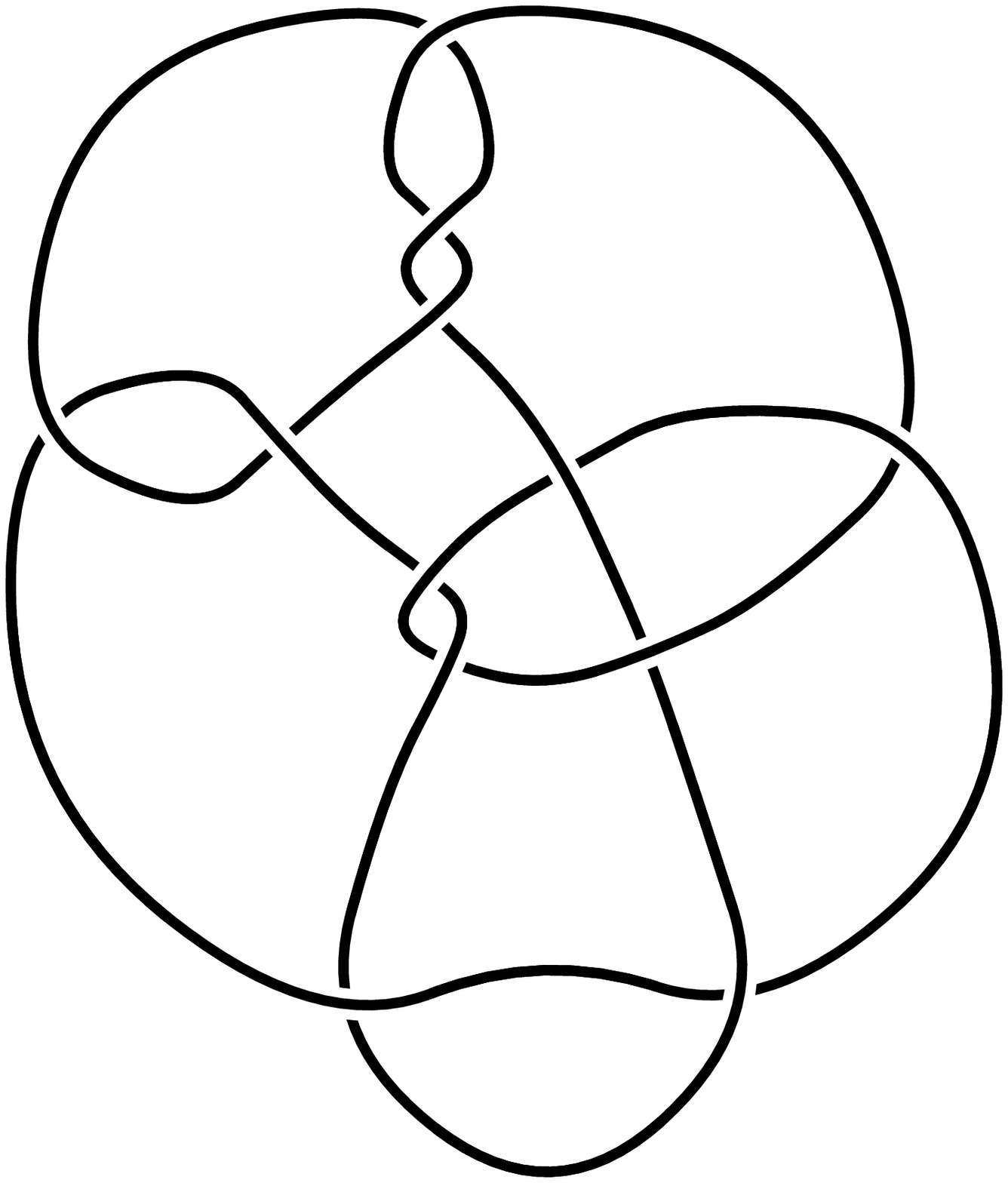}
    &
    \includegraphics[width=75pt]{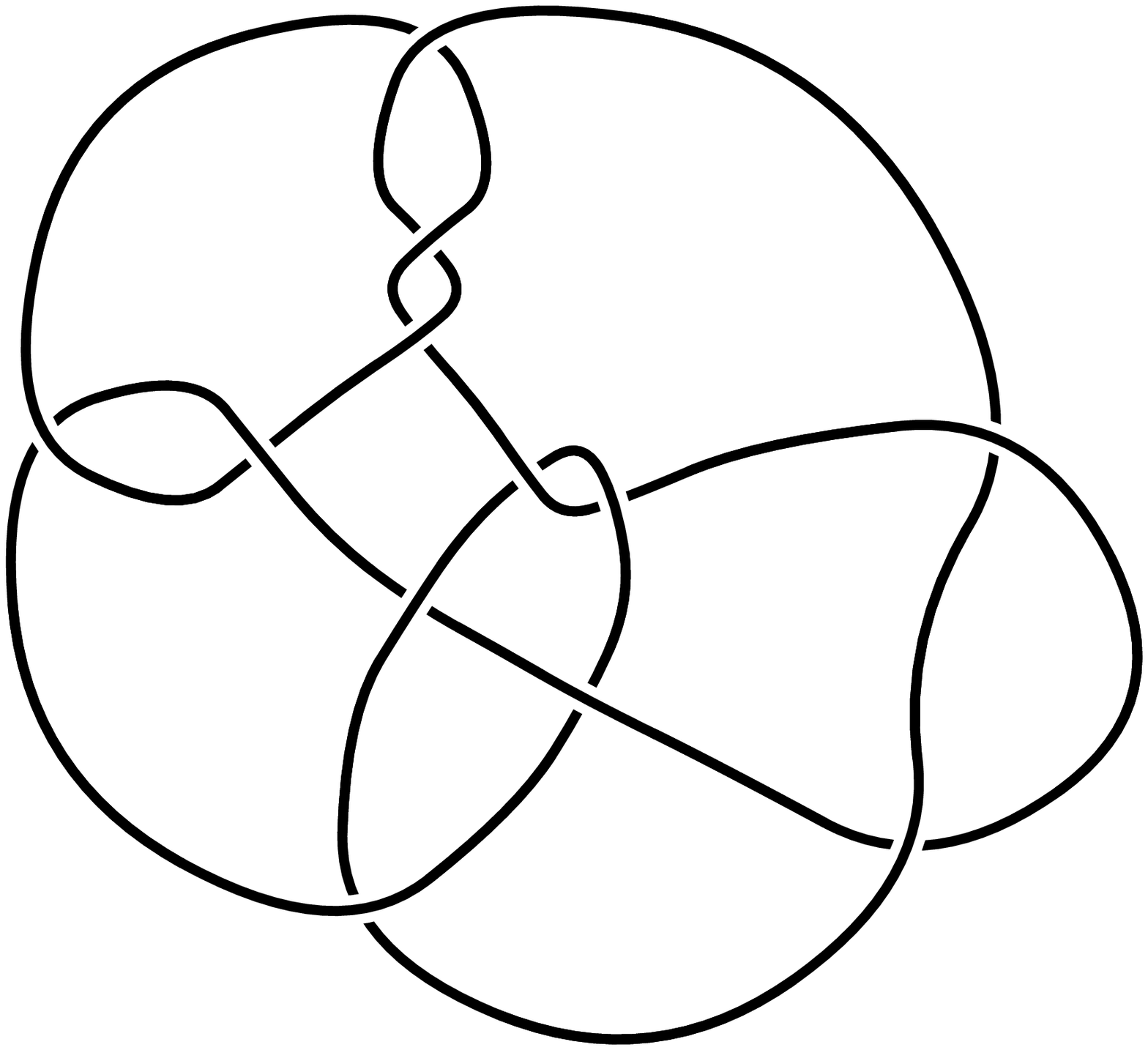}
    &
    &
    \\[-10pt]
    $12^A_{844}$ & $12^A_{846}$ & &
  \end{tabular}
  \caption{Alternating $12$-crossing mutant cliques 3/4}
  \end{centering}
\end{figure}

\begin{figure}[htbp]
  \begin{centering}
  \begin{tabular}{ccc}
    \includegraphics[width=75pt]{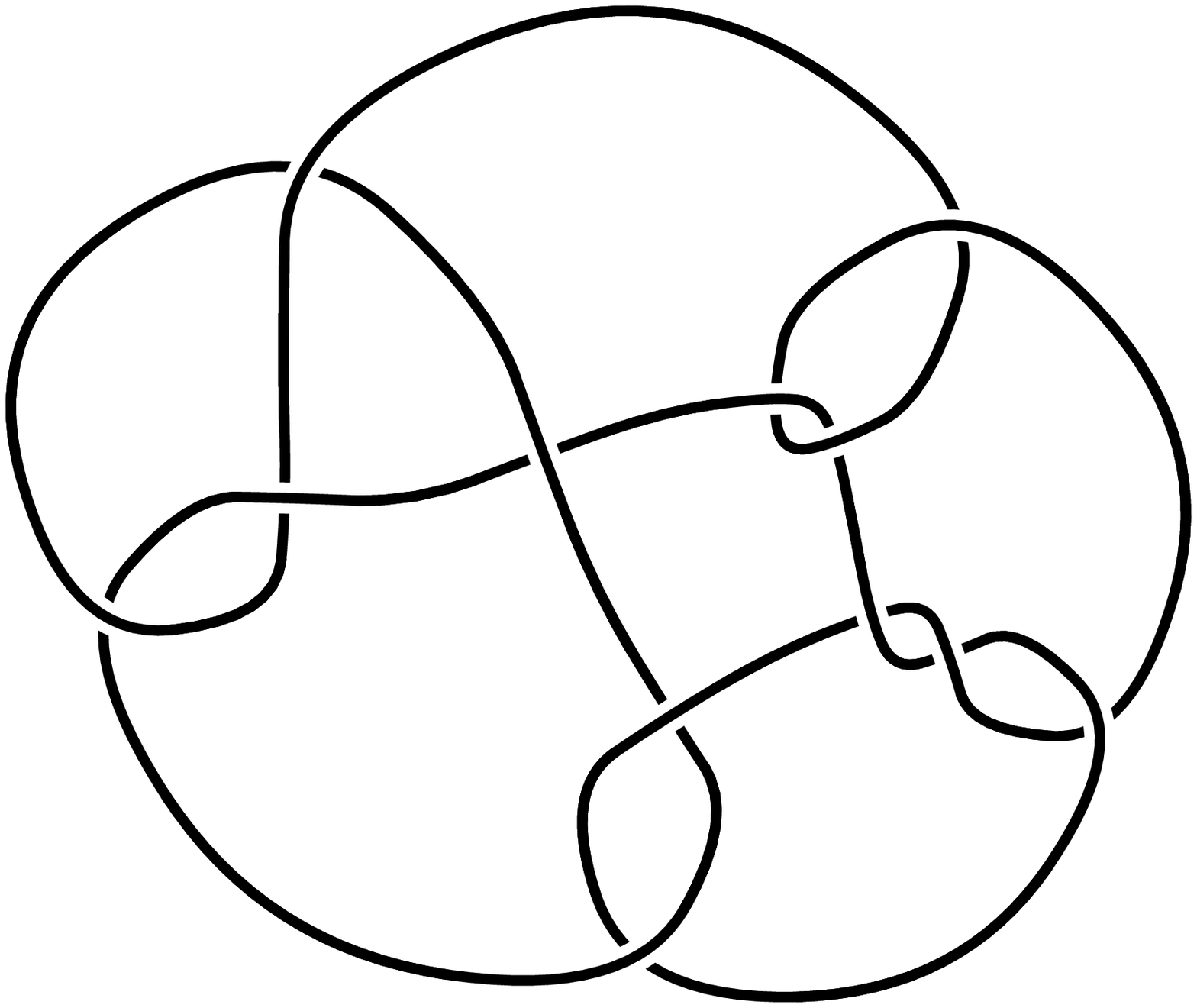}
    &
    \includegraphics[width=75pt]{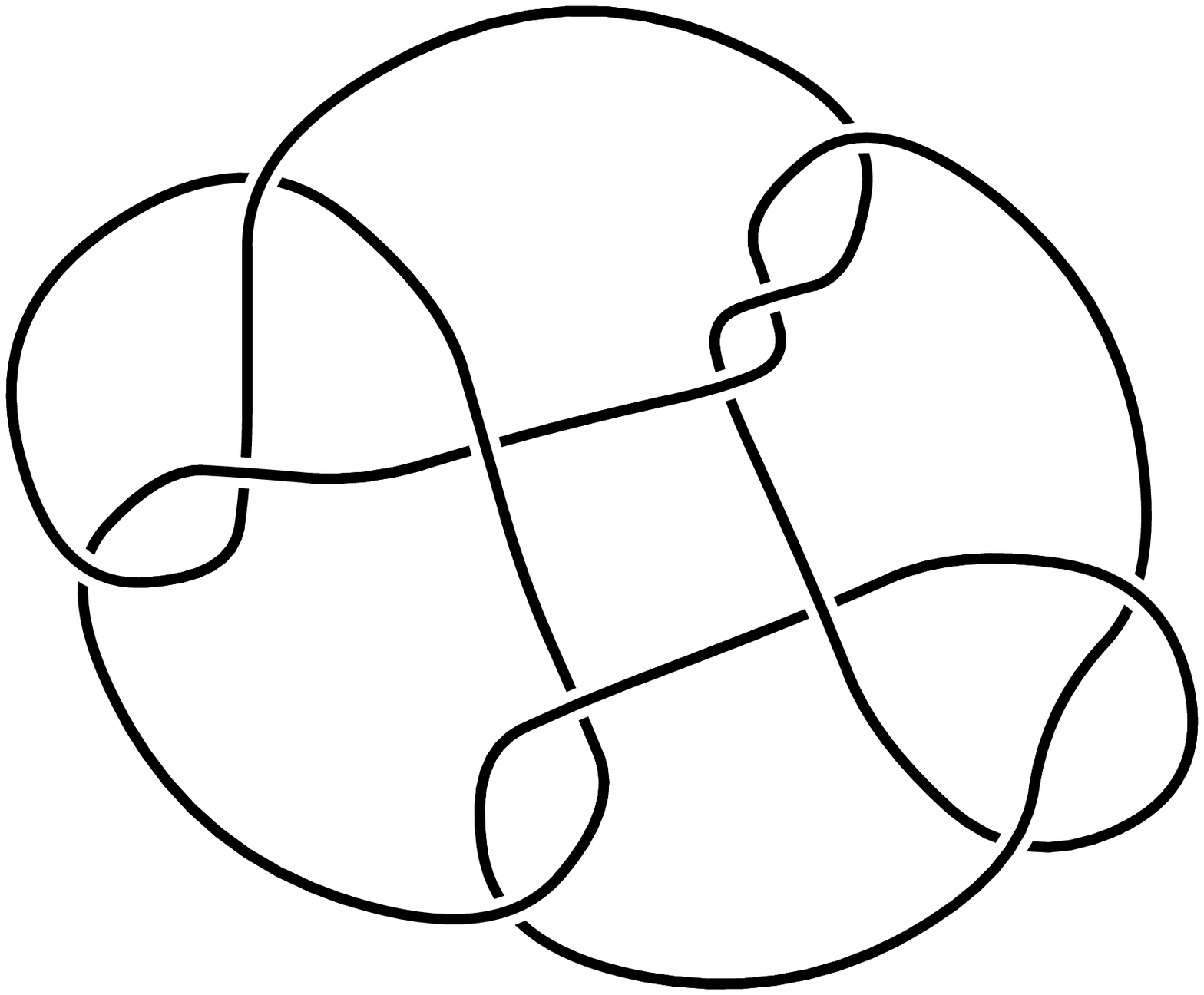}
    &
    \includegraphics[width=75pt]{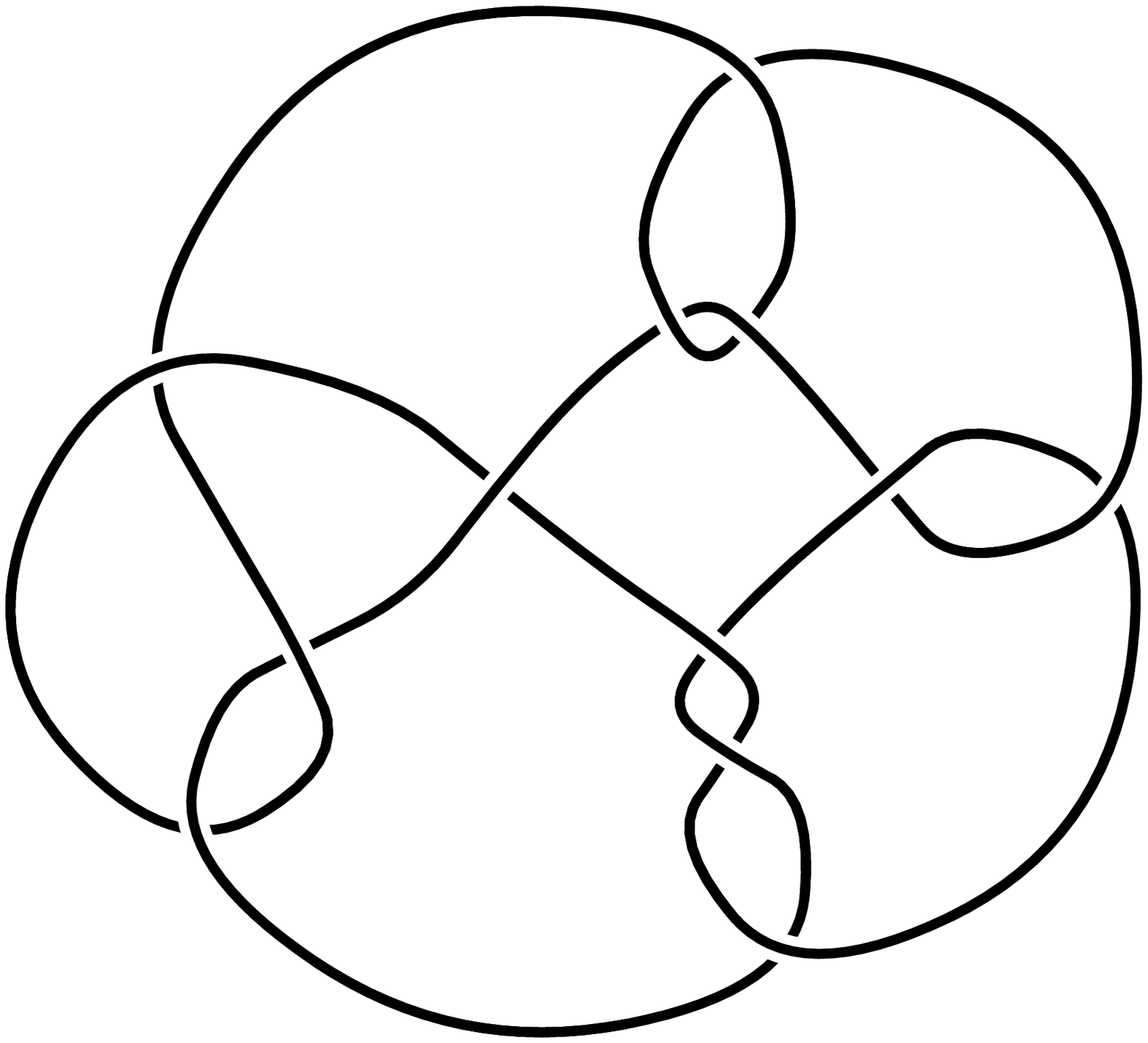}
    \\[-10pt]
    $12^A_{30}$ & $12^A_{33}$ & $12^A_{157}$
    \\[10pt]
    \hline
    &&\\[-10pt]
    \includegraphics[width=75pt]{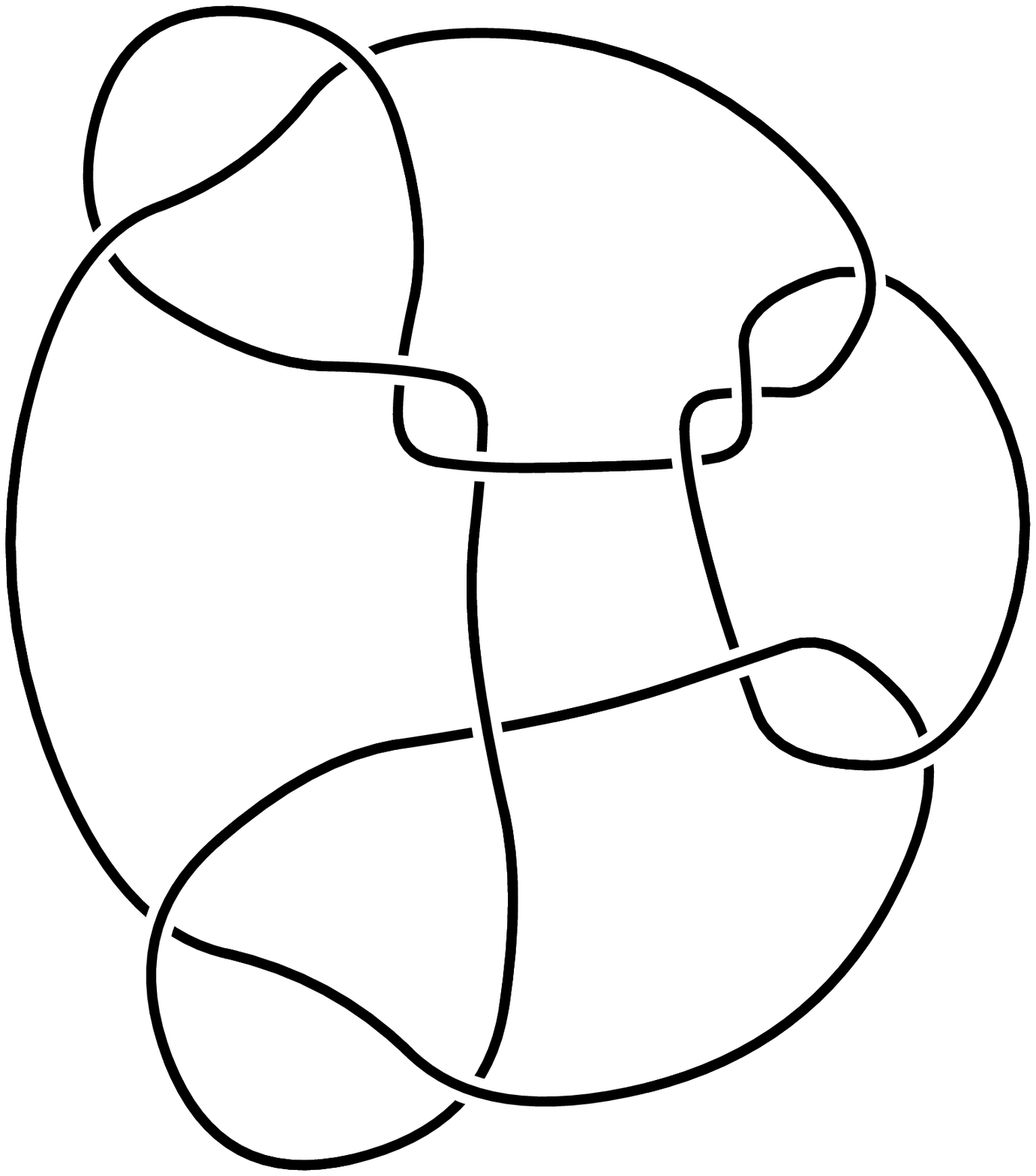}
    &
    \includegraphics[width=75pt]{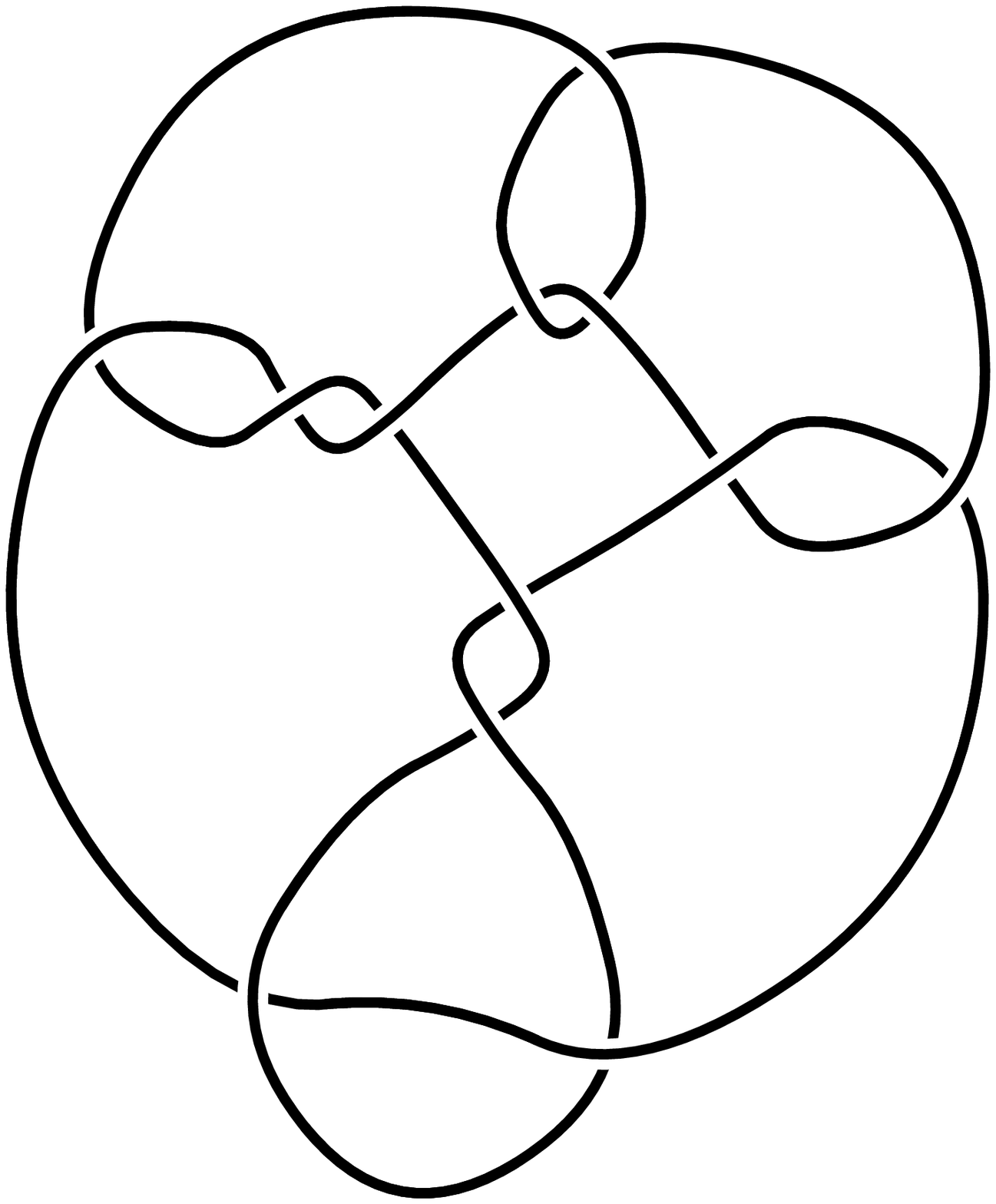}
    &
    \includegraphics[width=75pt]{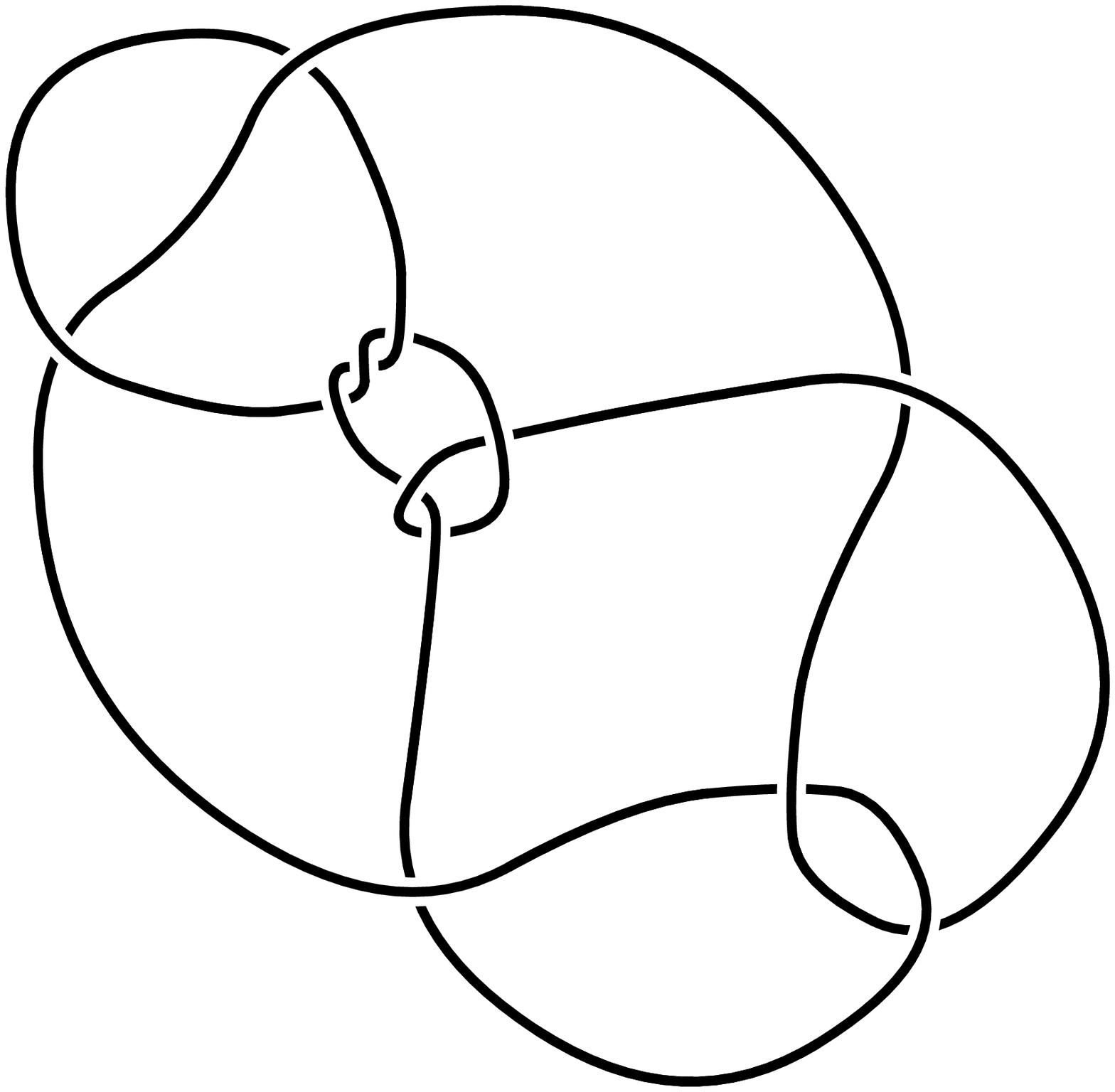}
    \\[-10pt]
    $12^A_{116}$ & $12^A_{122}$ & $12^A_{182}$
  \end{tabular}
  \caption{Alternating $12$-crossing mutant cliques 4/4}
  \end{centering}
\end{figure}


\begin{figure}[htbp]
  \begin{centering}
  \begin{tabular}{cc@{\hspace{10pt}}|@{\hspace{10pt}}cc}
    \includegraphics[width=75pt]{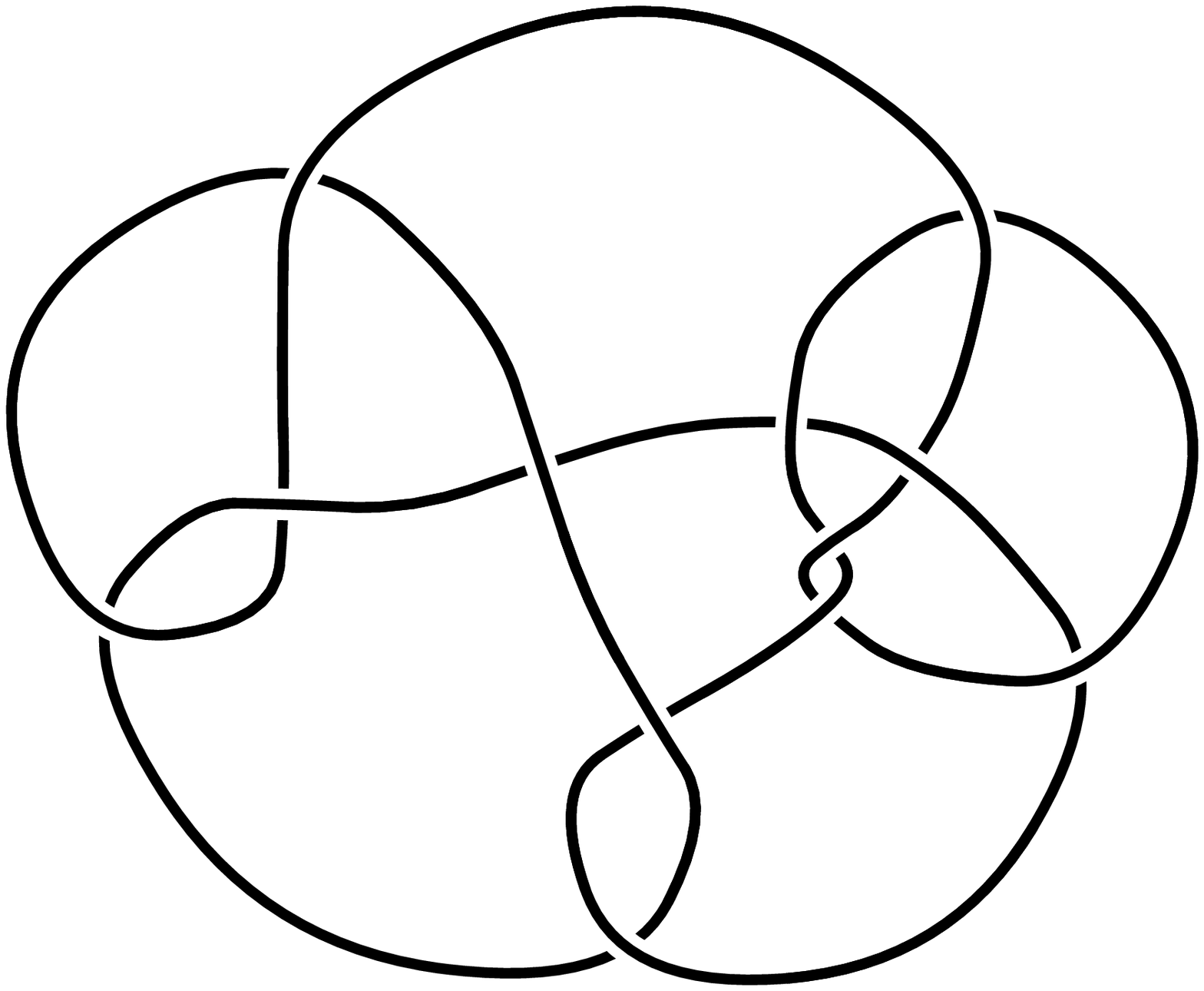}
    &
    \includegraphics[width=75pt]{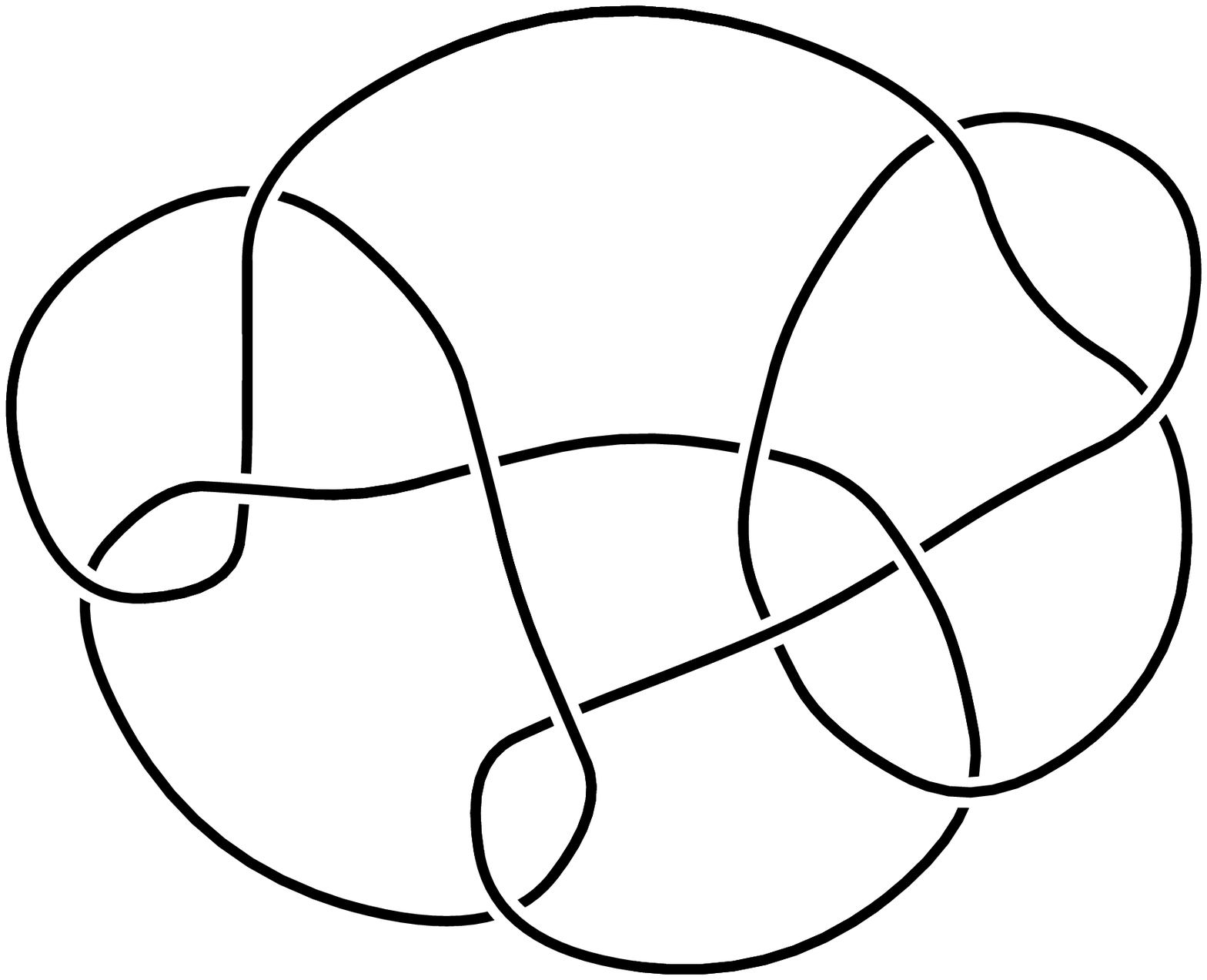}
    &
    \includegraphics[width=75pt]{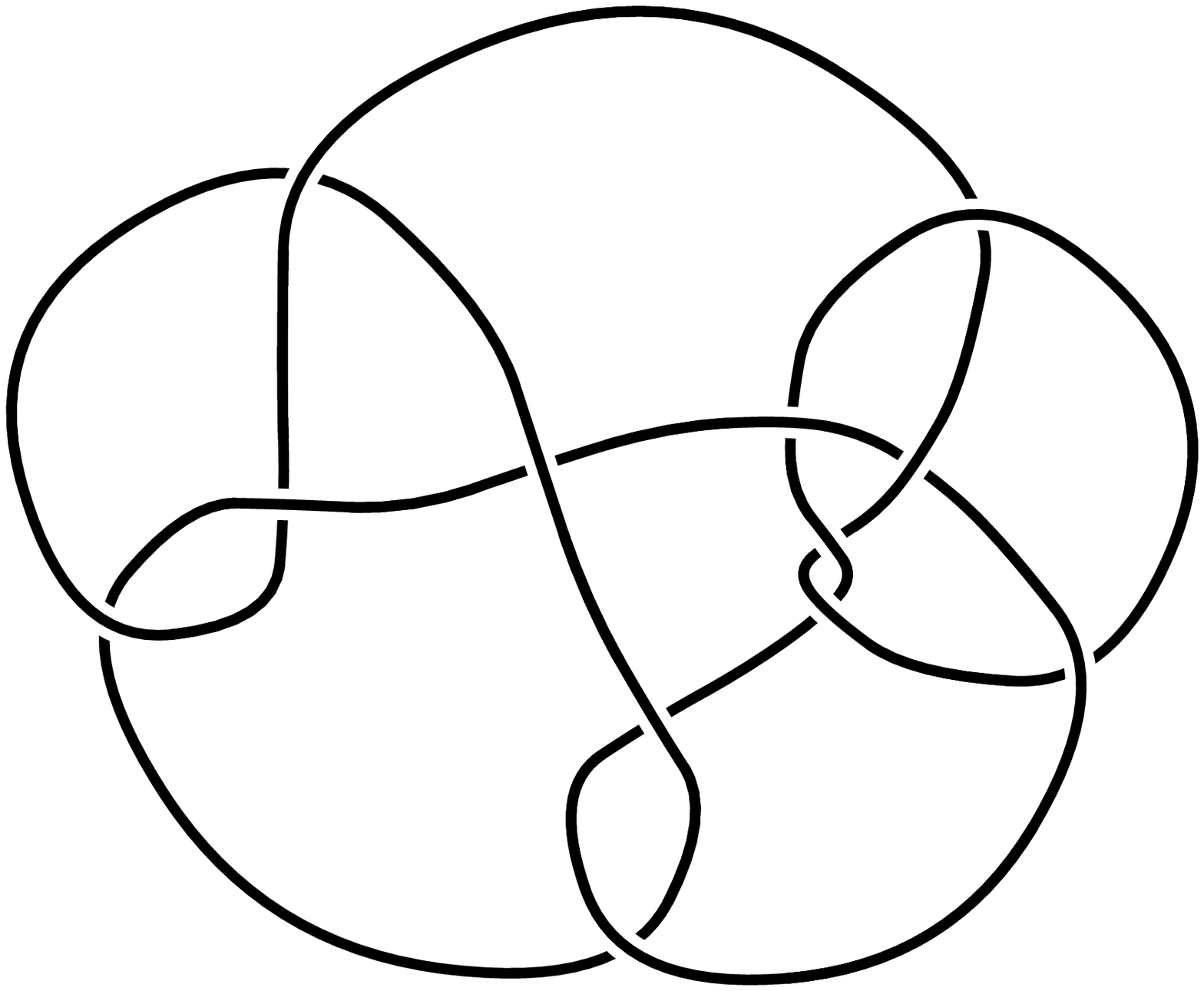}
    &
    \includegraphics[width=75pt]{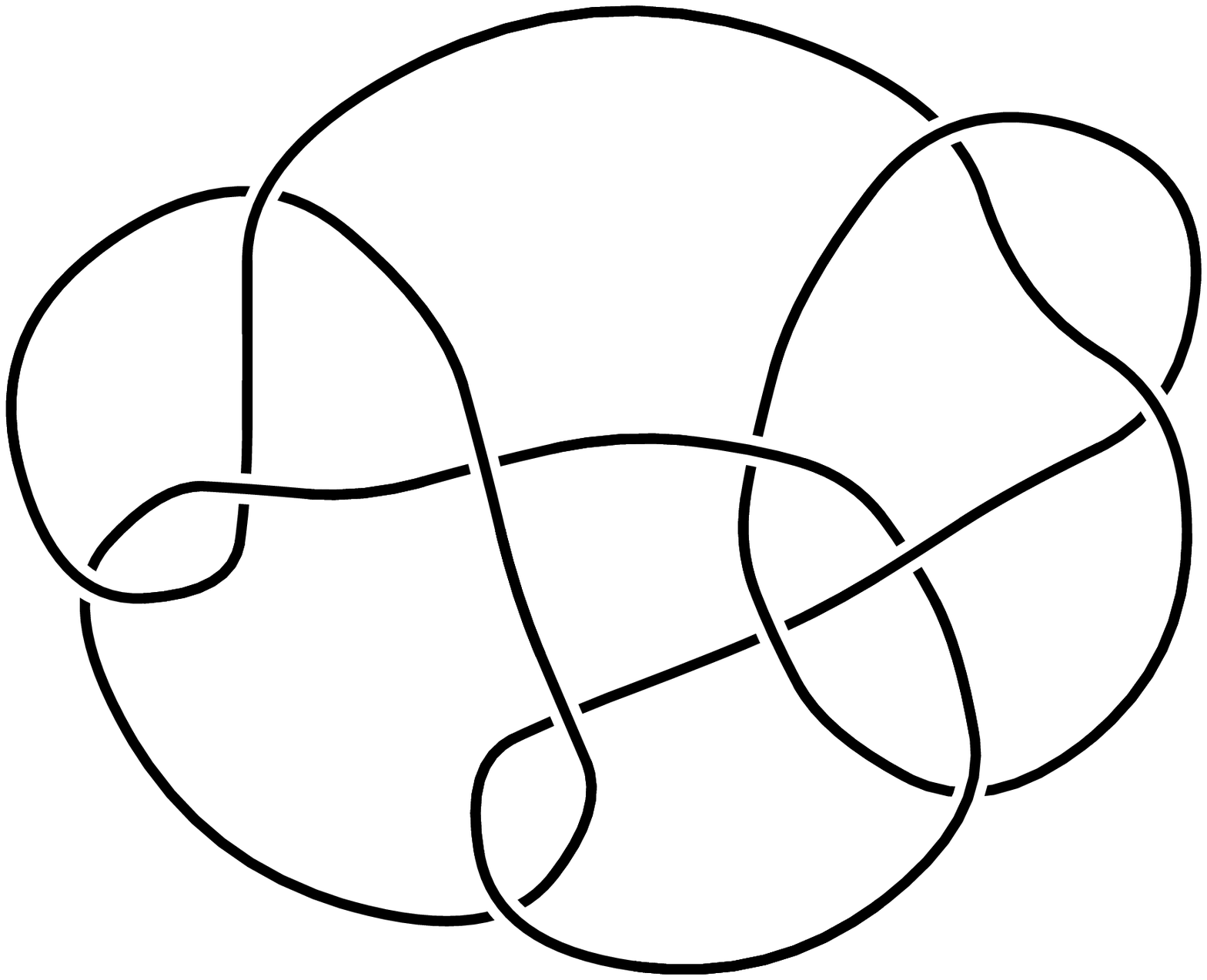}
    \\[-10pt]
    $12^N_{21}$ & $12^N_{29}$ & $12^N_{22}$ & $12^N_{30}$
    \\[10pt]
    \hline
    &&&\\[-10pt]
    \includegraphics[width=75pt]{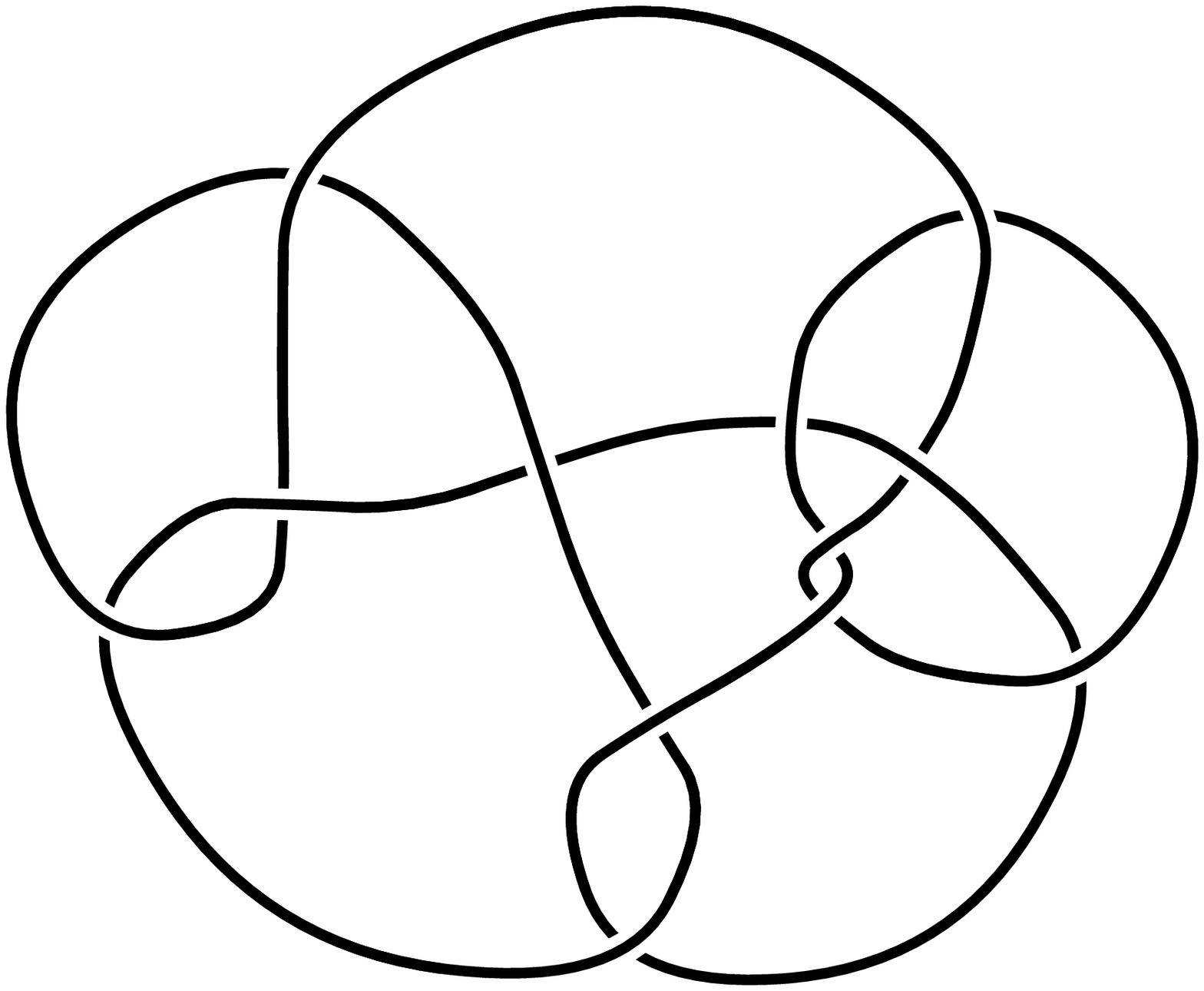}
    &
    \includegraphics[width=75pt]{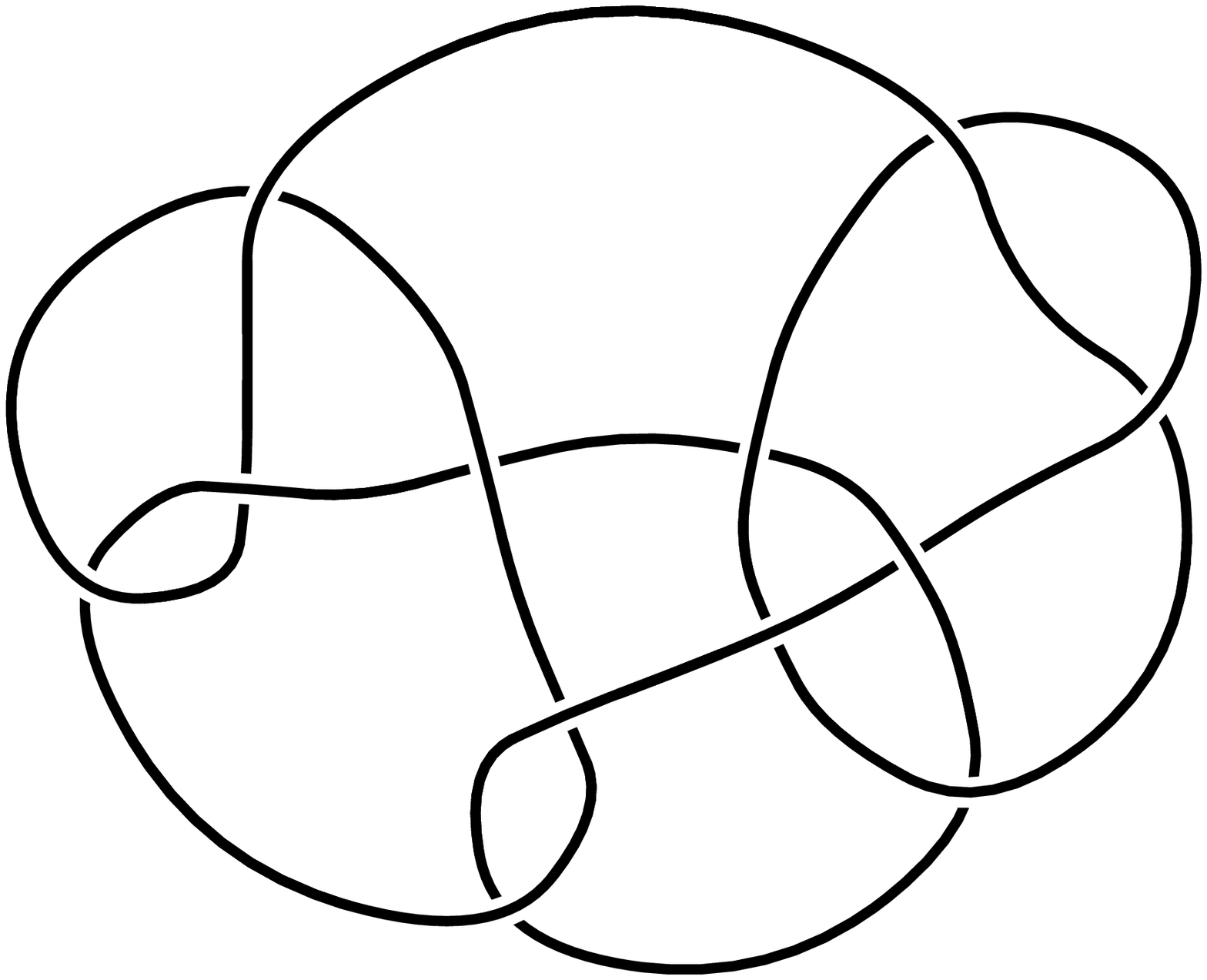}
    &
    \includegraphics[width=75pt]{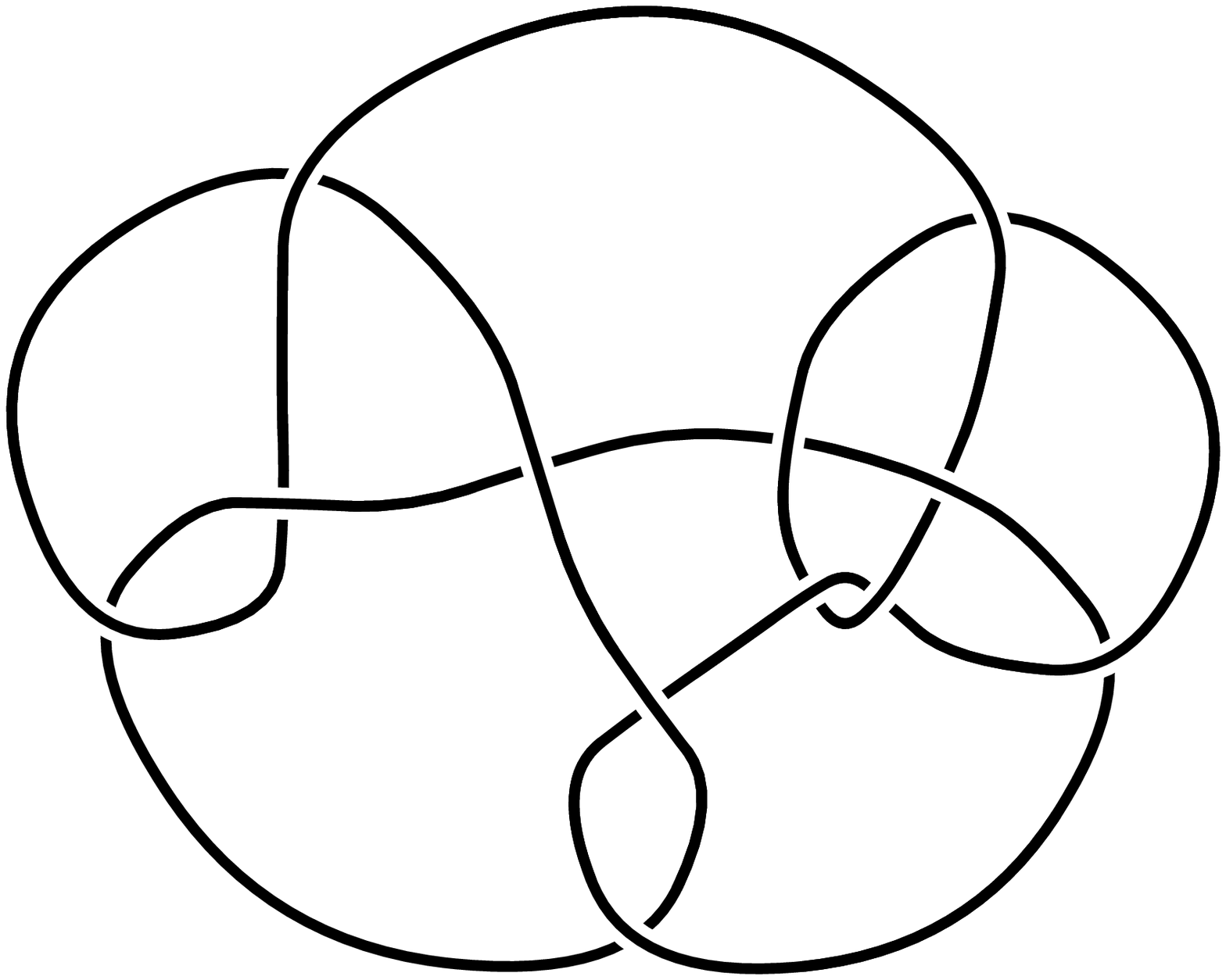}
    &
    \includegraphics[width=75pt]{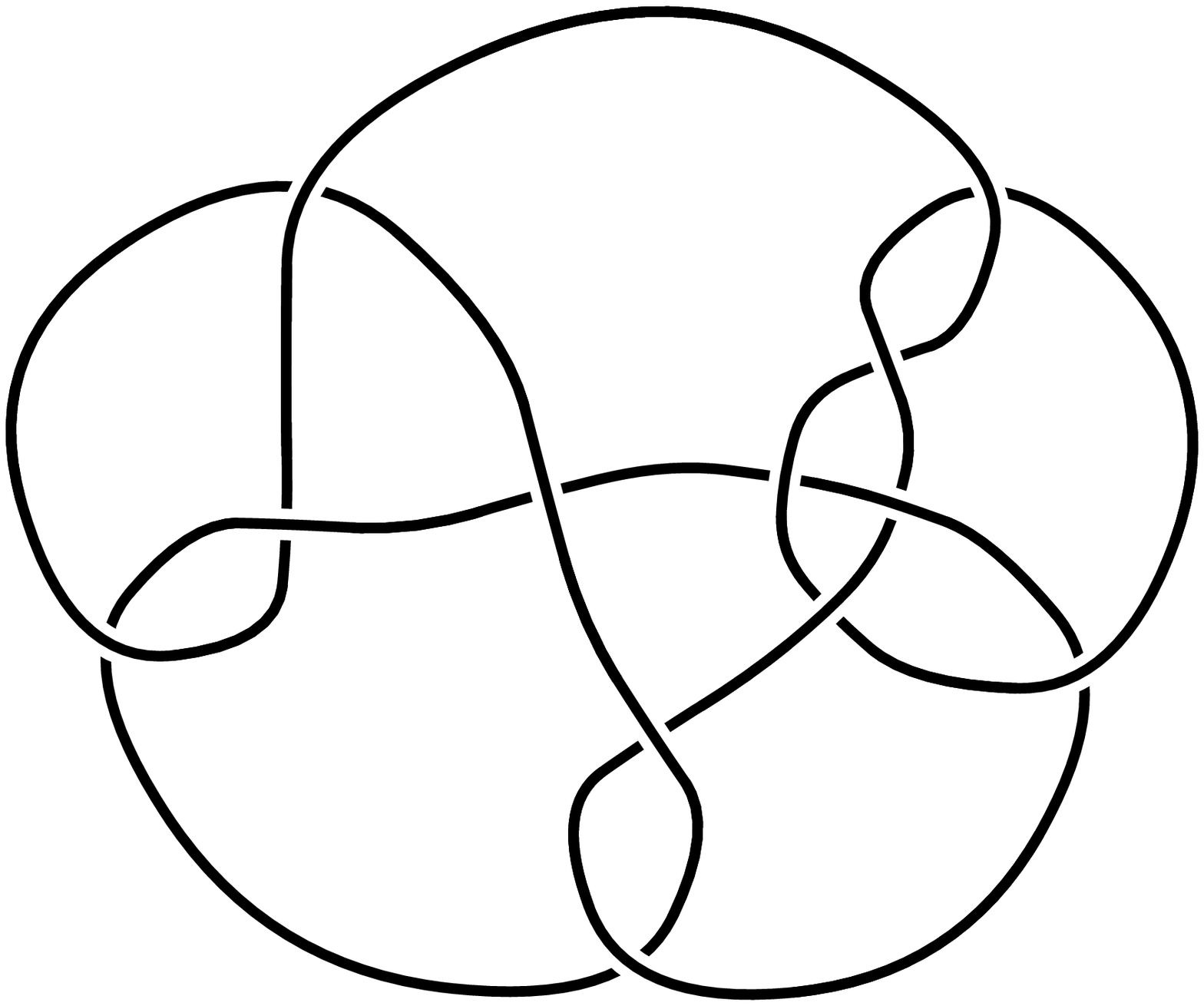}
    \\[-10pt]
    $12^N_{23}$ & $12^N_{31}$ & $12^N_{26}$ & $12^N_{32}$
    \\[10pt]
    \hline
    &&&\\[-10pt]
    \includegraphics[width=75pt]{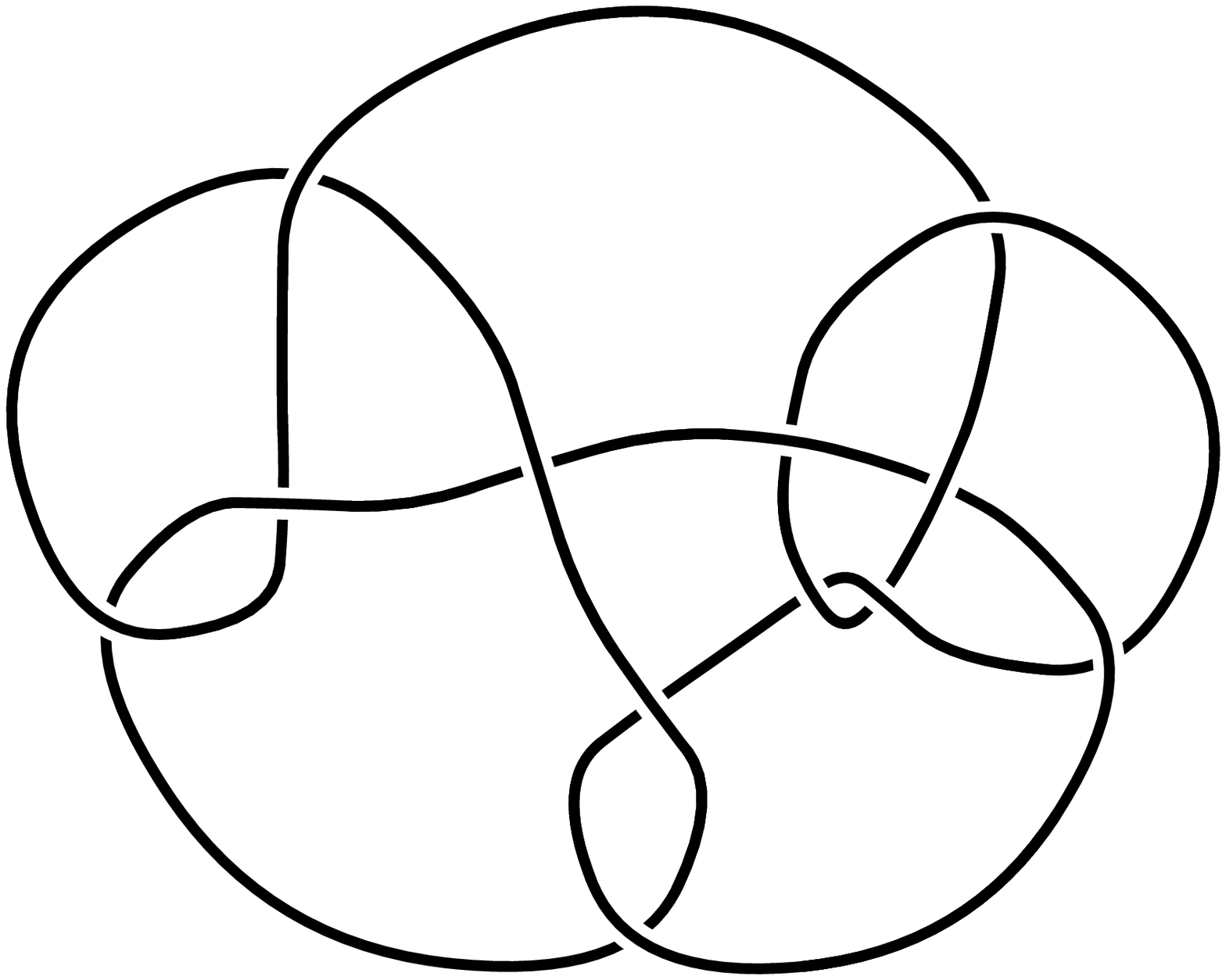}
    &
    \includegraphics[width=75pt]{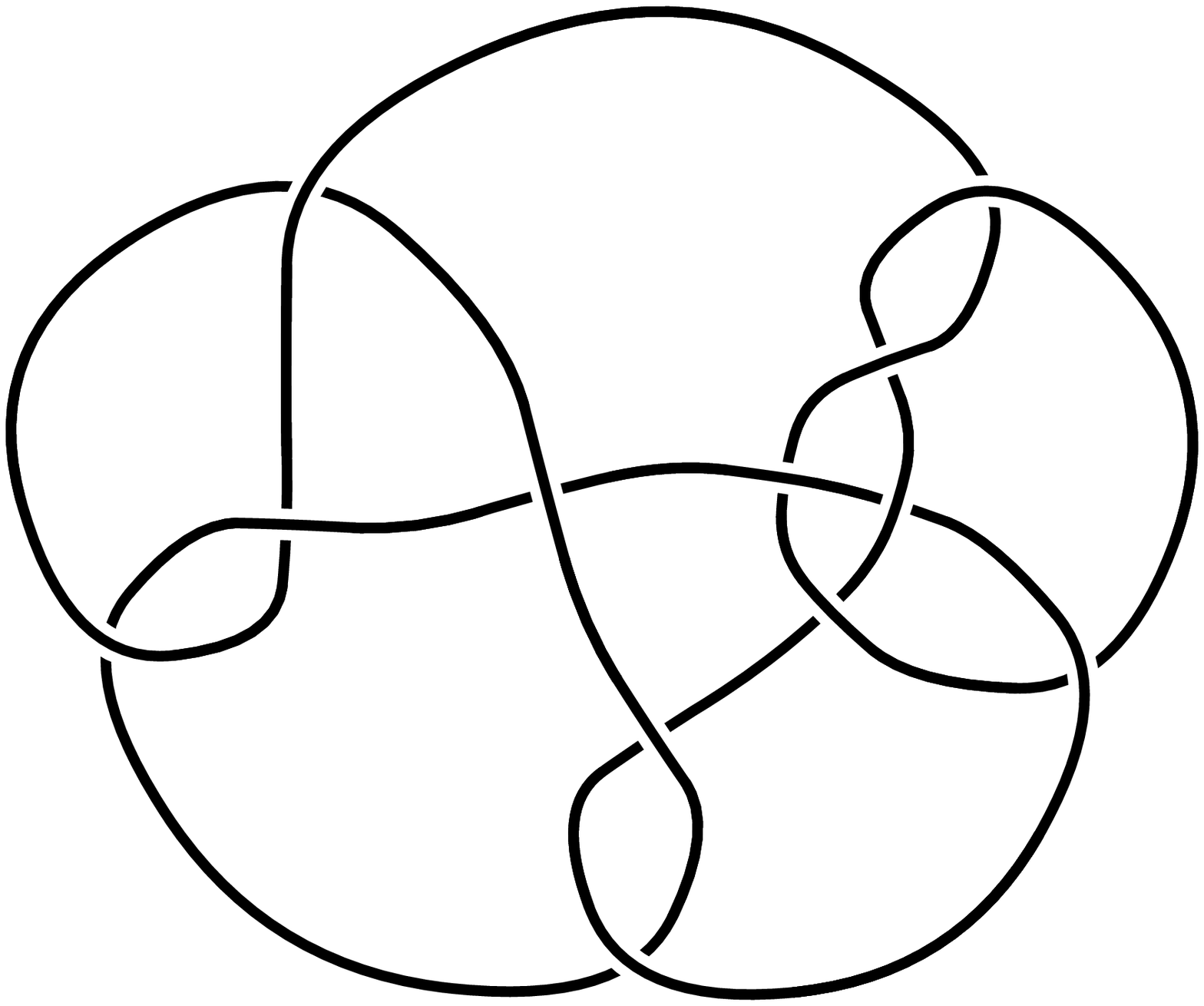}
    &
    \includegraphics[width=75pt]{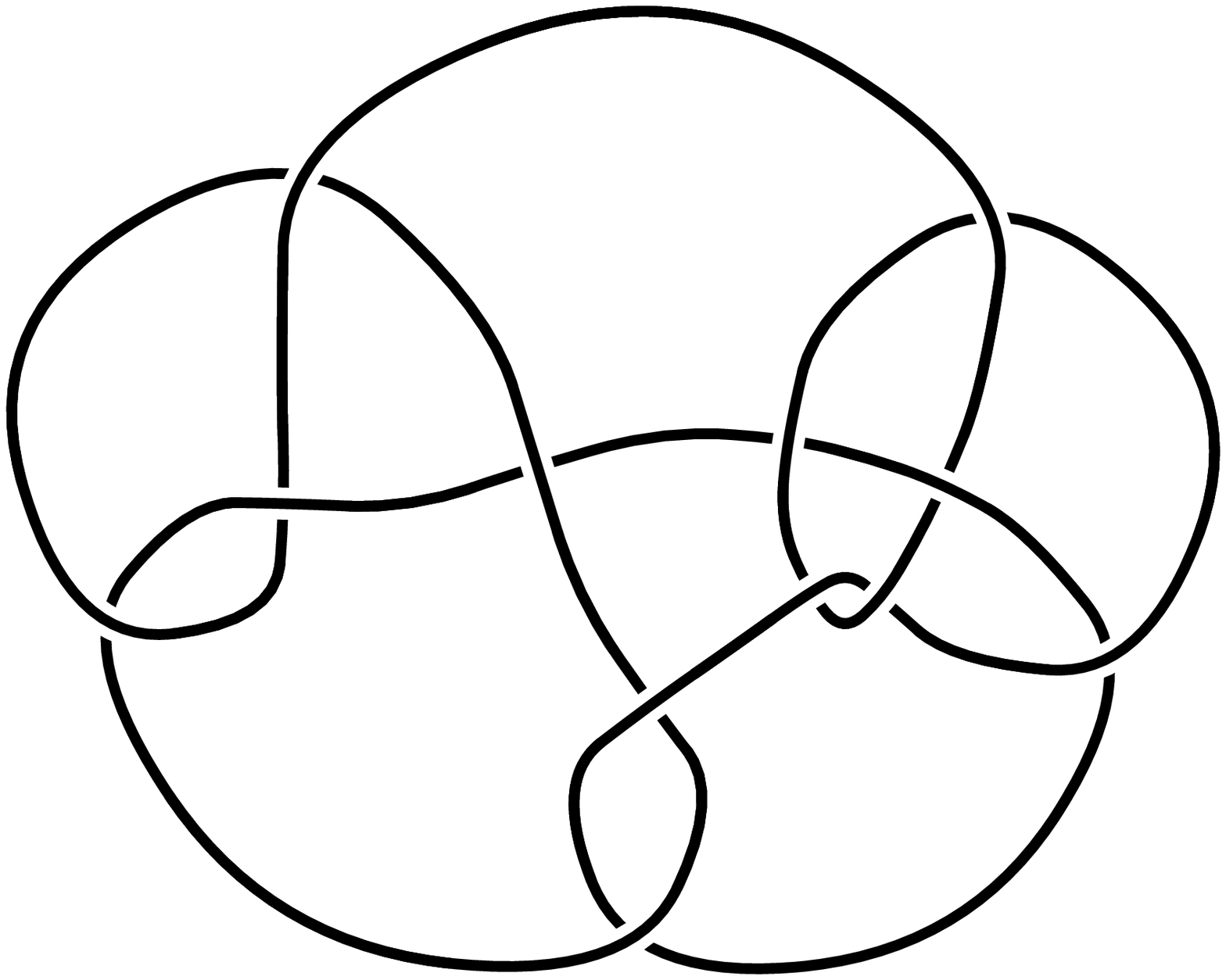}
    &
    \includegraphics[width=75pt]{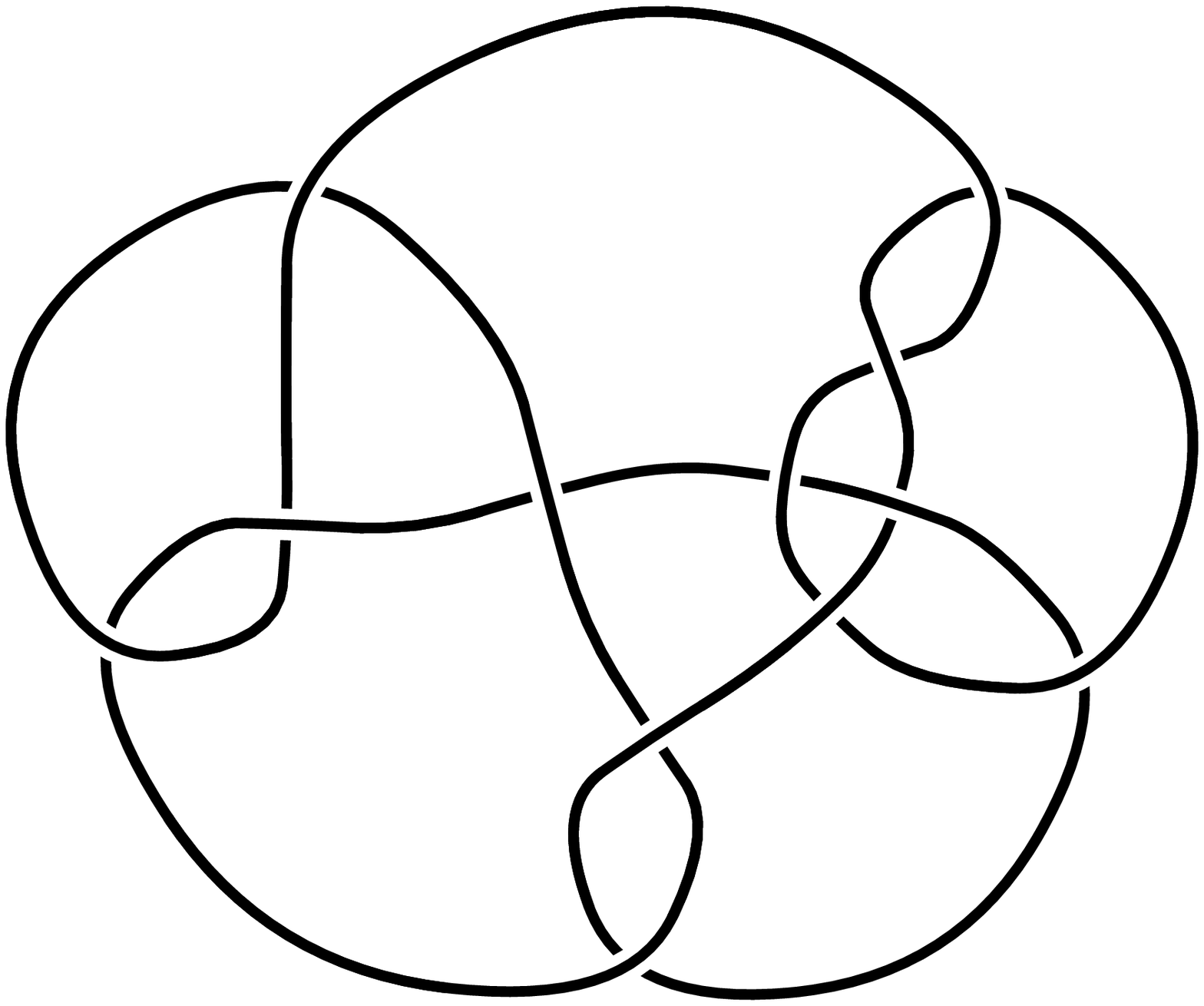}
    \\[-10pt]
    $12^N_{27}$ & $12^N_{33}$ & $12^N_{28}$ & $12^N_{34}$
    \\[10pt]
    \hline
    &&&\\[-10pt]
    \includegraphics[width=75pt]{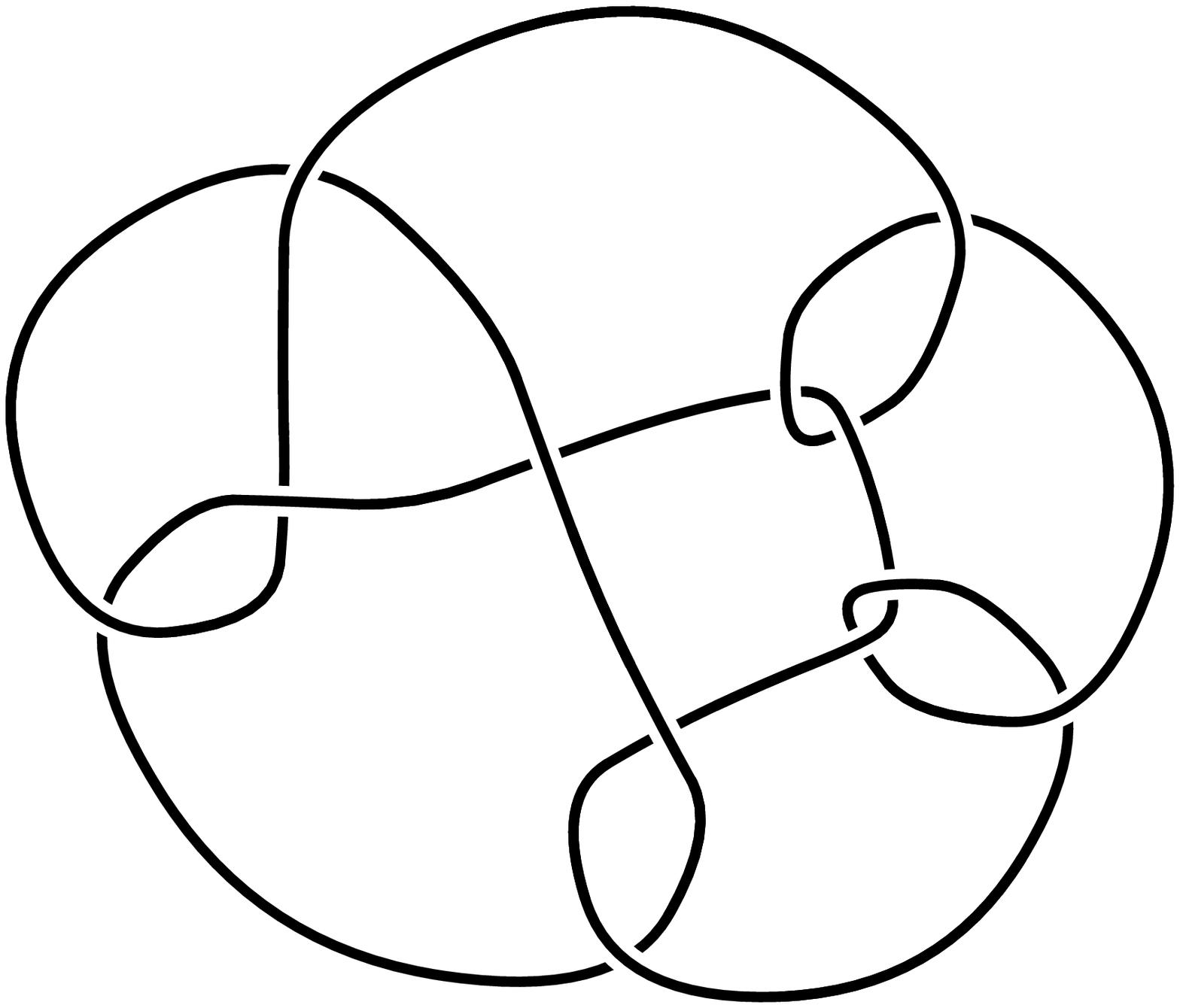}
    &
    \includegraphics[width=75pt]{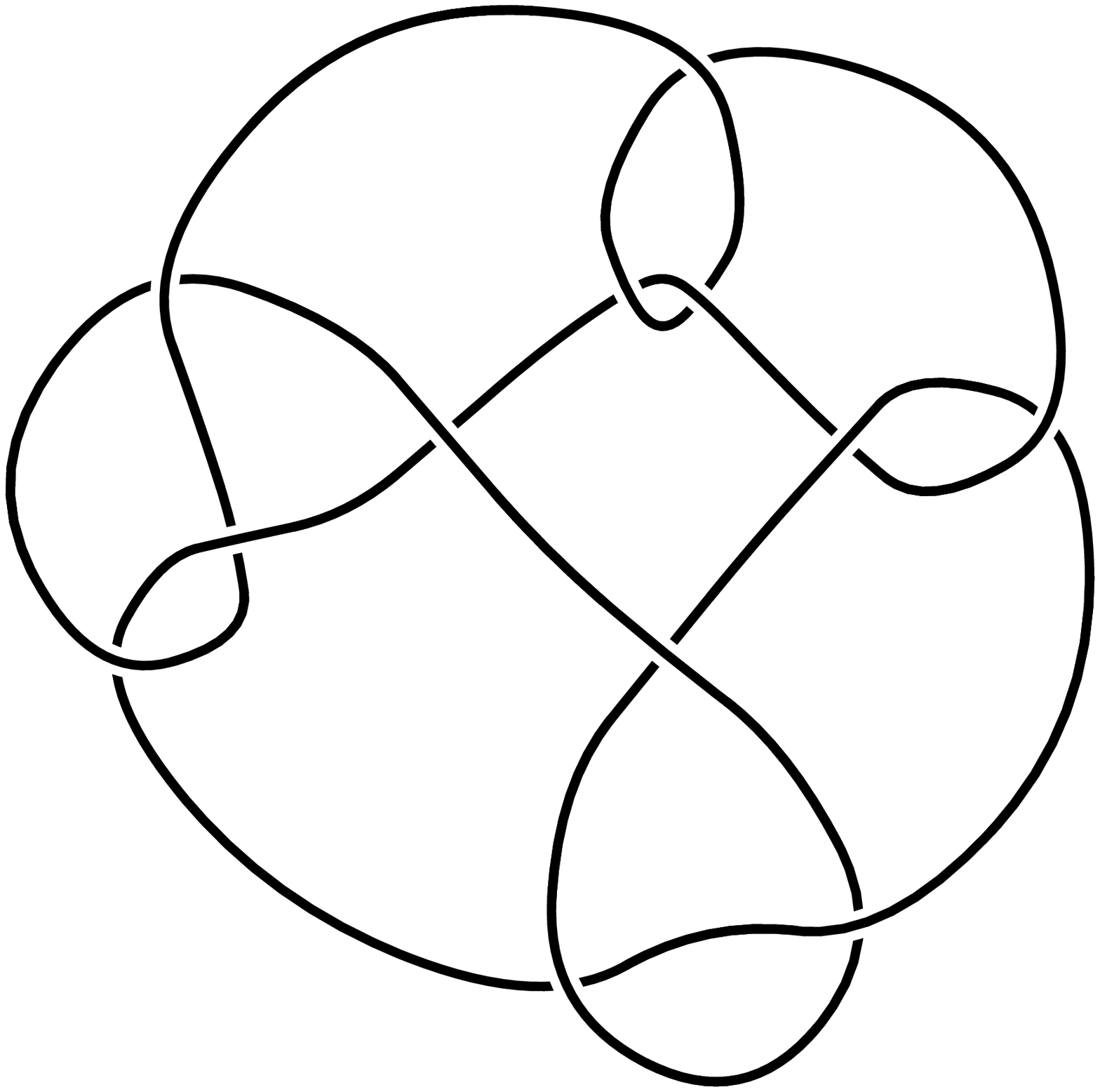}
    &
    \includegraphics[width=75pt]{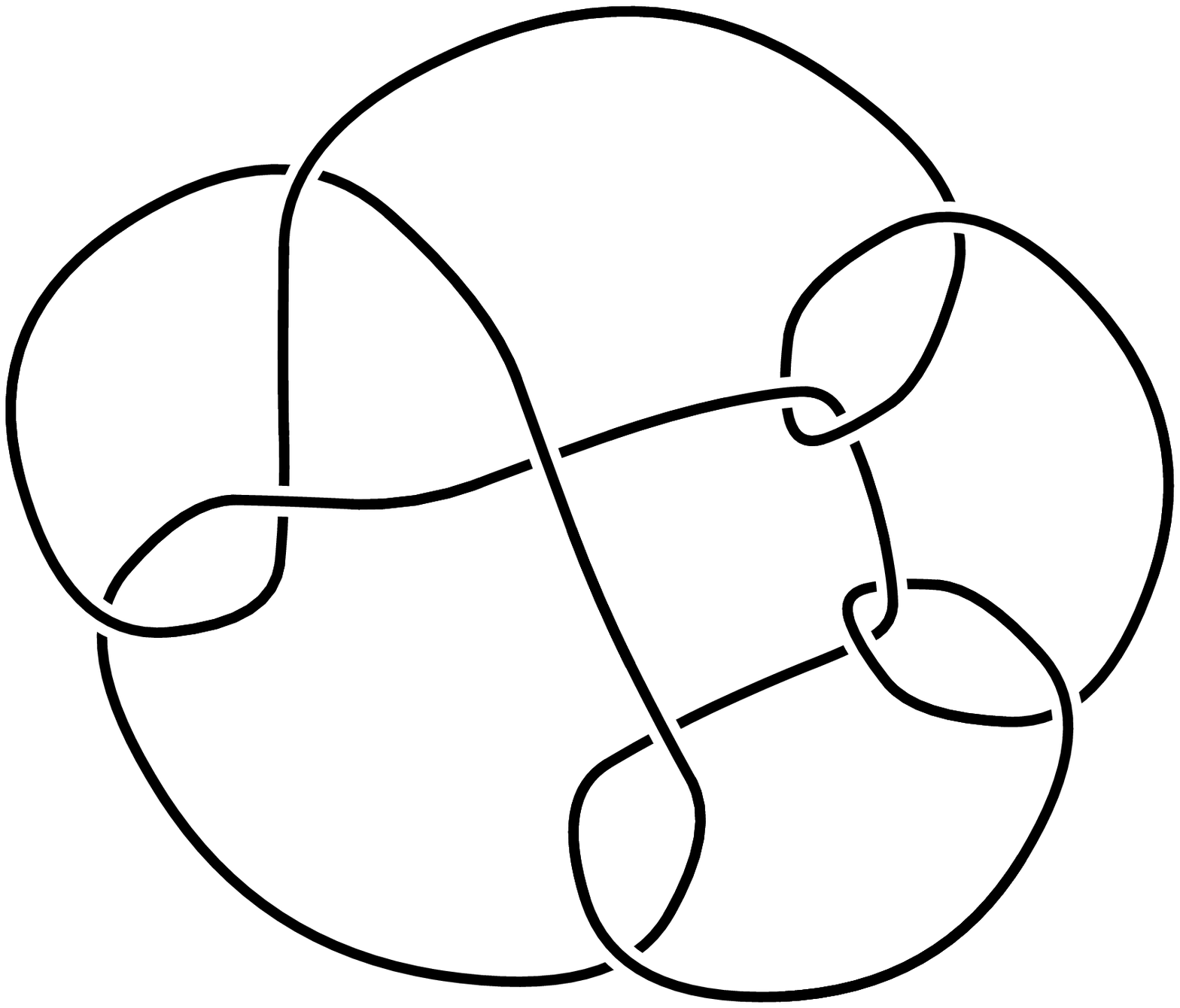}
    &
    \includegraphics[width=75pt]{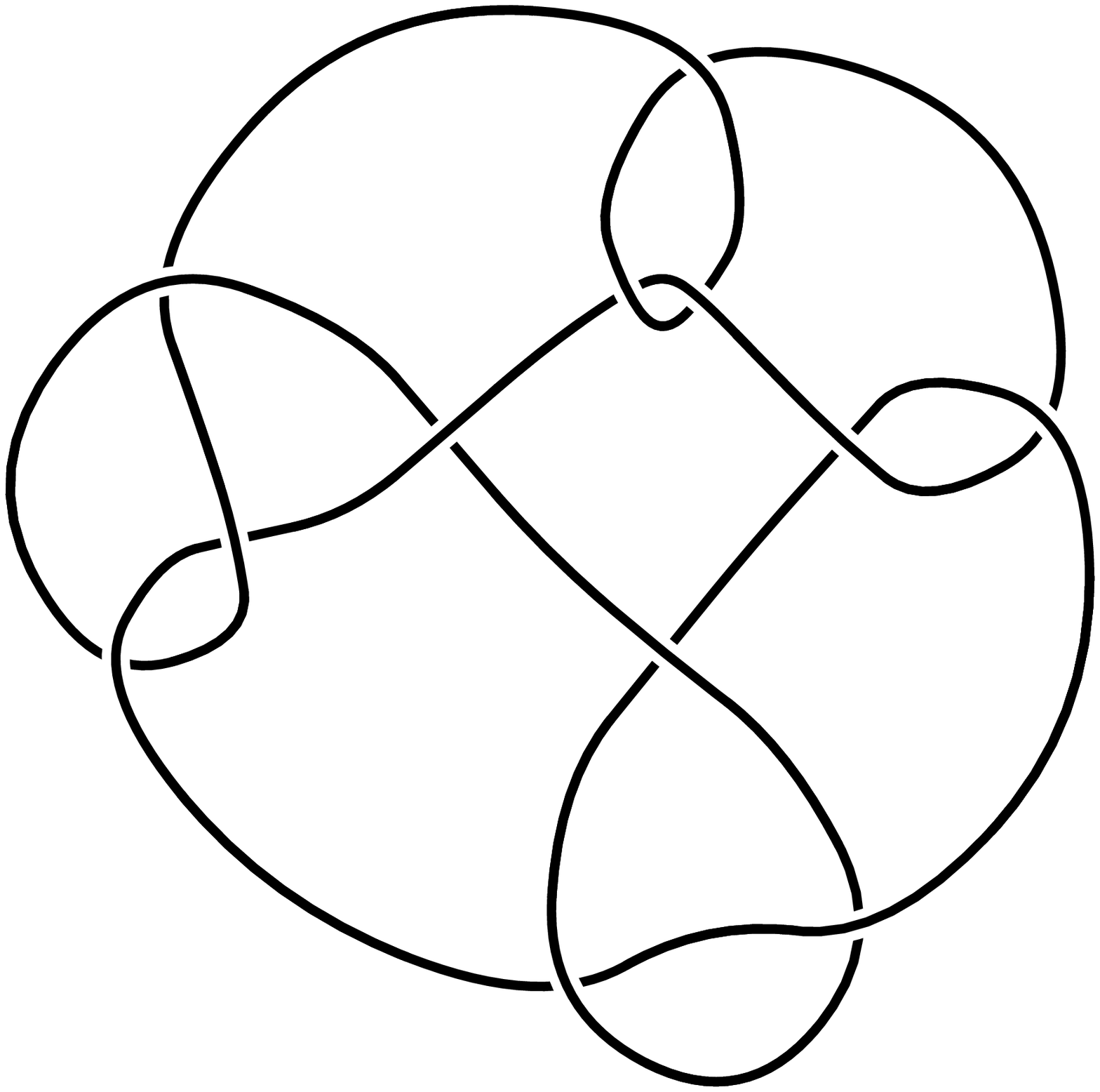}
    \\[-10pt]
    $12^N_{55}$ & $12^N_{223}$ & $12^N_{58}$ & $12^N_{222}$
    \\[10pt]
    \hline
    &&&\\[-10pt]
    \includegraphics[width=75pt]{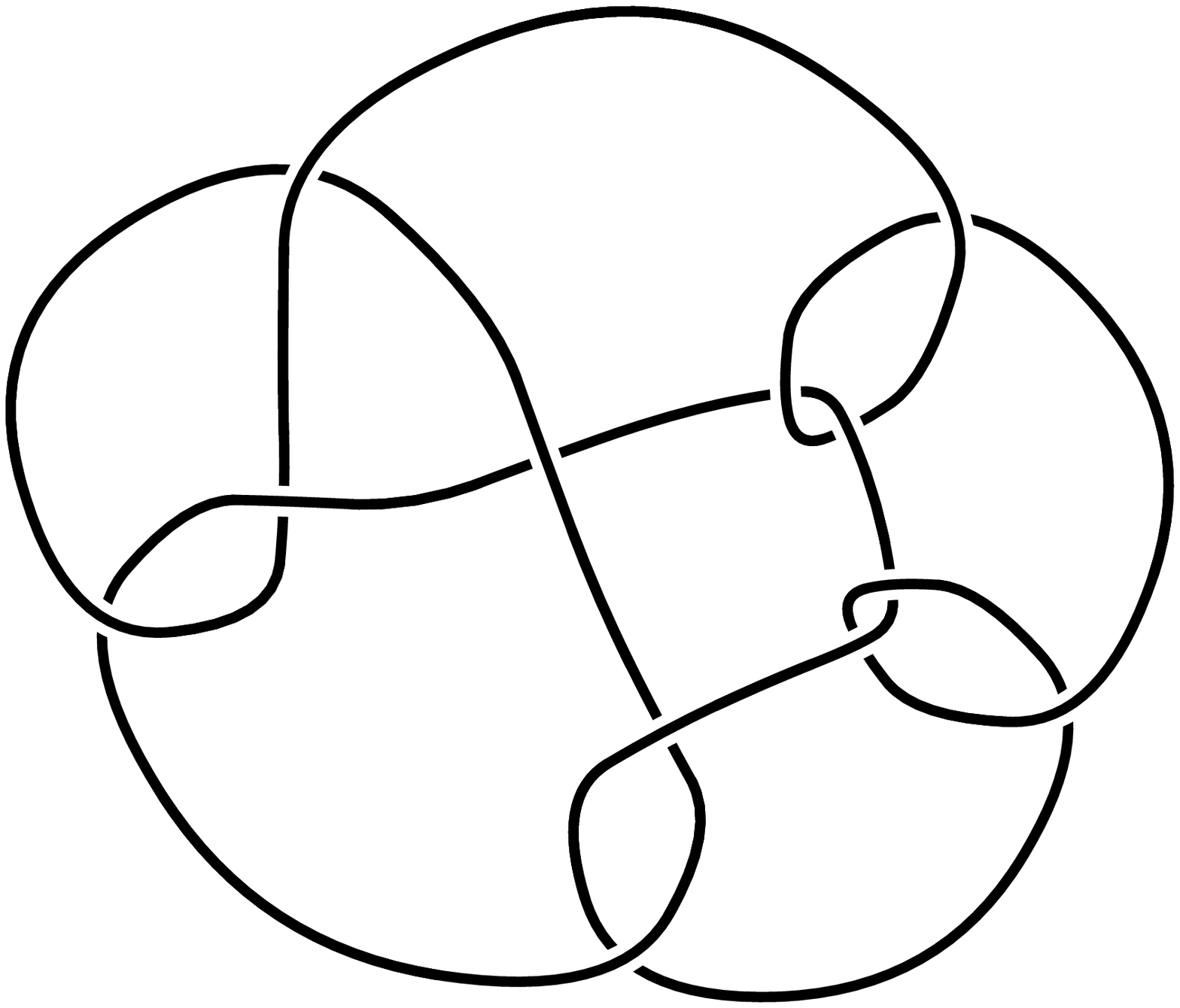}
    &
    \includegraphics[width=75pt]{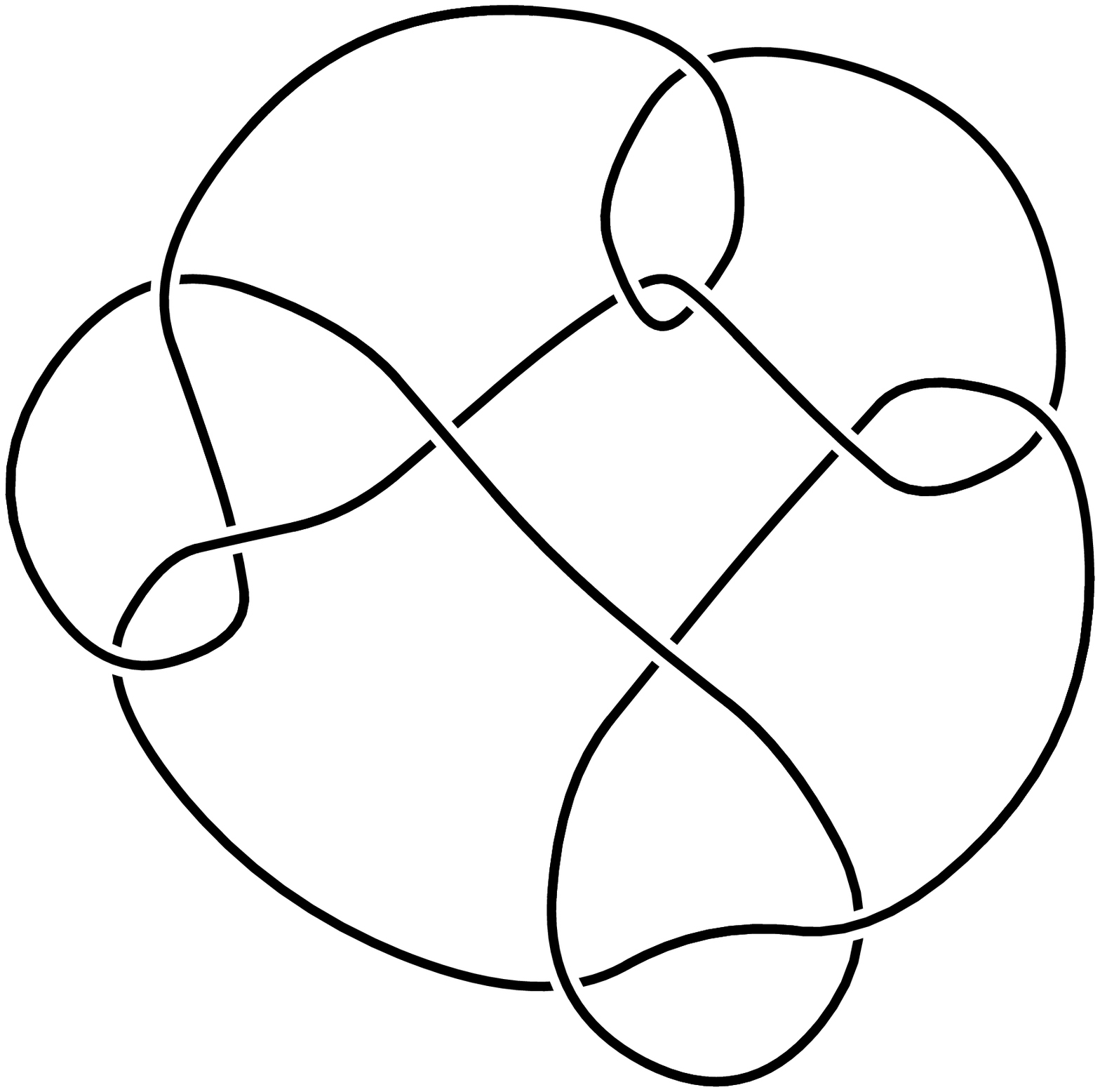}
    &
    \includegraphics[width=75pt]{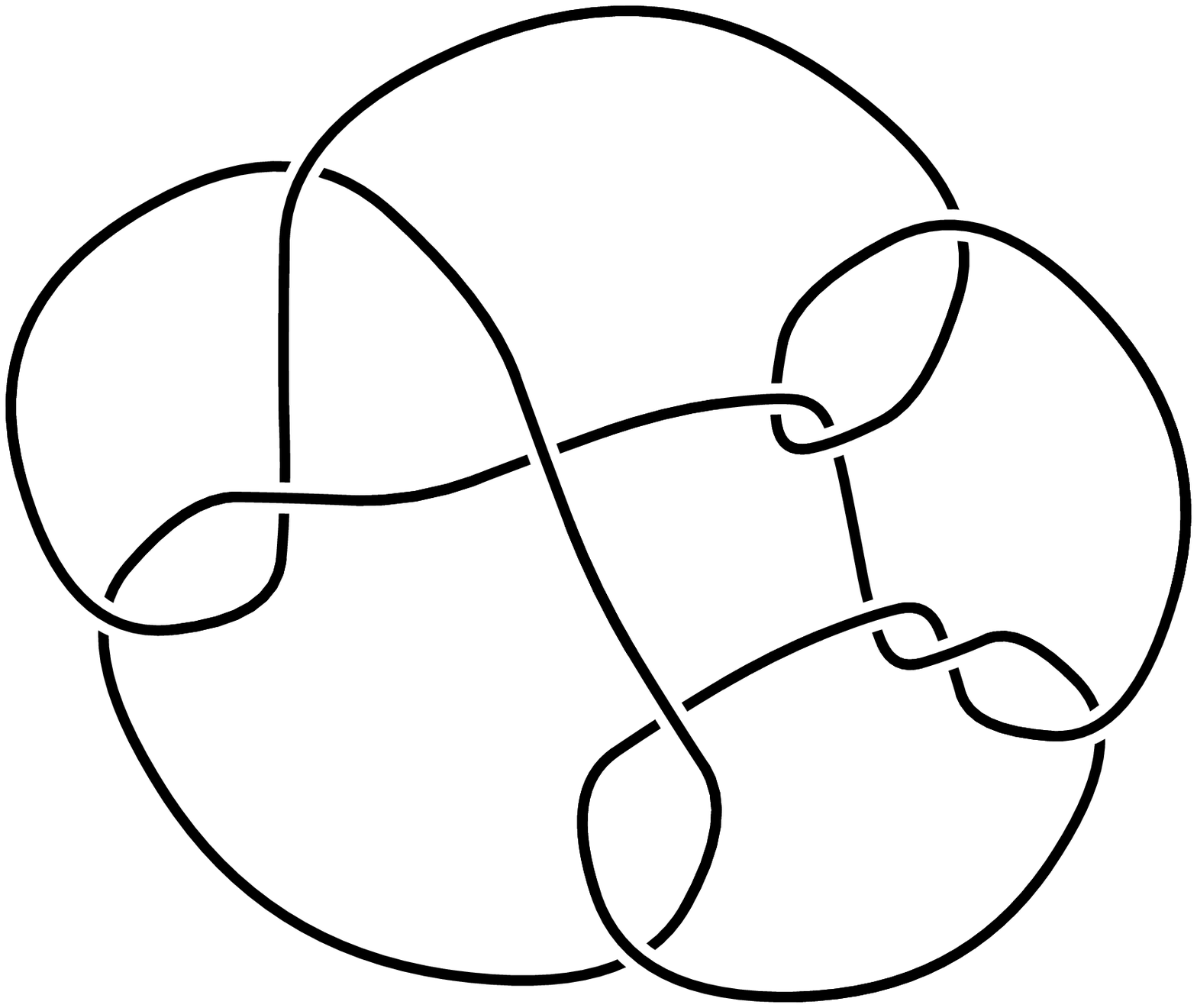}
    &
    \includegraphics[width=75pt]{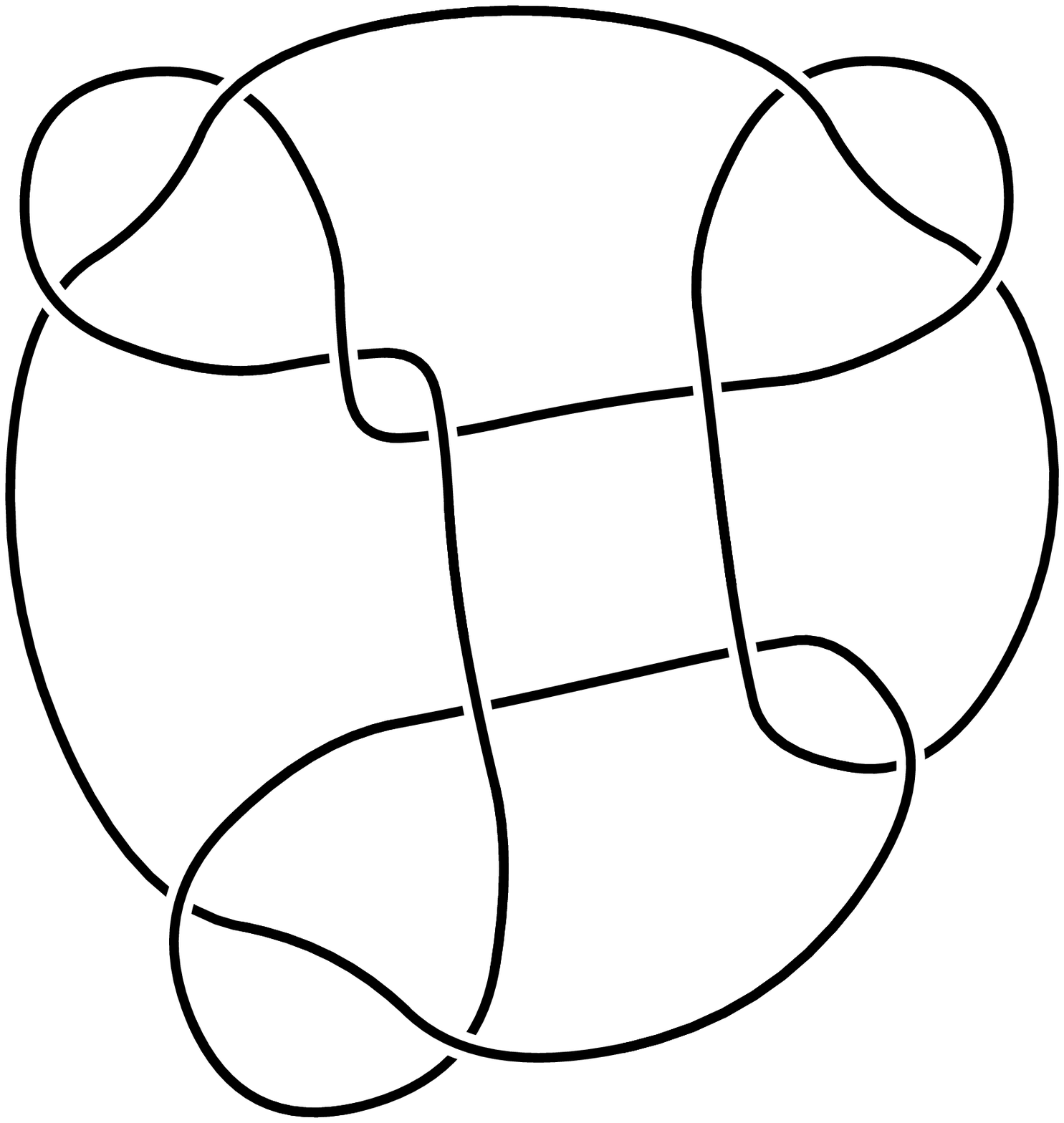}
    \\[-10pt]
    $12^N_{59}$ & $12^N_{220}$ & $12^N_{63}$ & $12^N_{225}$
  \end{tabular}
  \caption{Nonalternating $12$-crossing mutant cliques 1/6}
  \end{centering}
\end{figure}

\begin{figure}[htbp]
  \begin{centering}
  \begin{tabular}{cc@{\hspace{10pt}}|@{\hspace{10pt}}cc}
    \includegraphics[width=75pt]{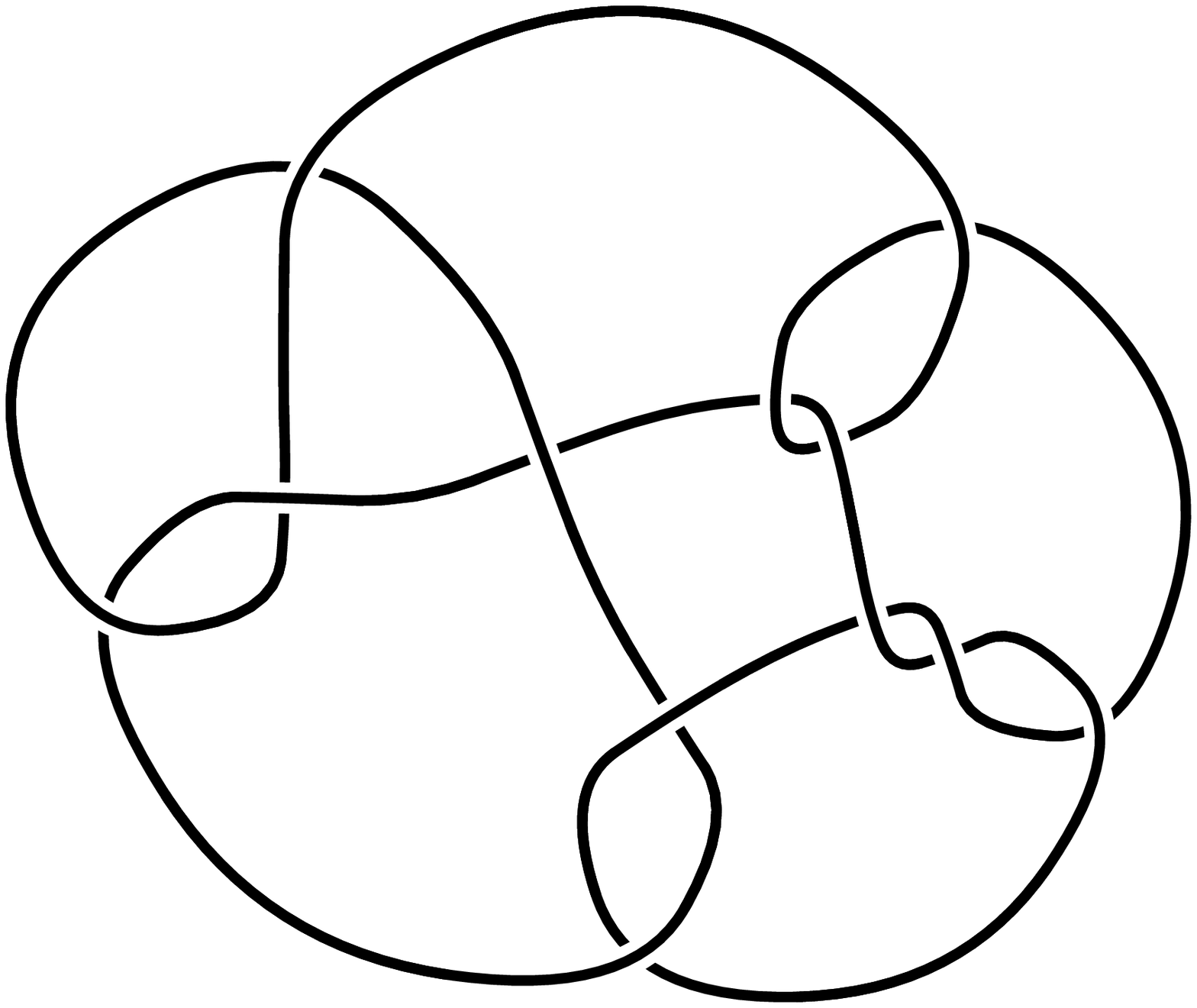}
    &
    \includegraphics[width=75pt]{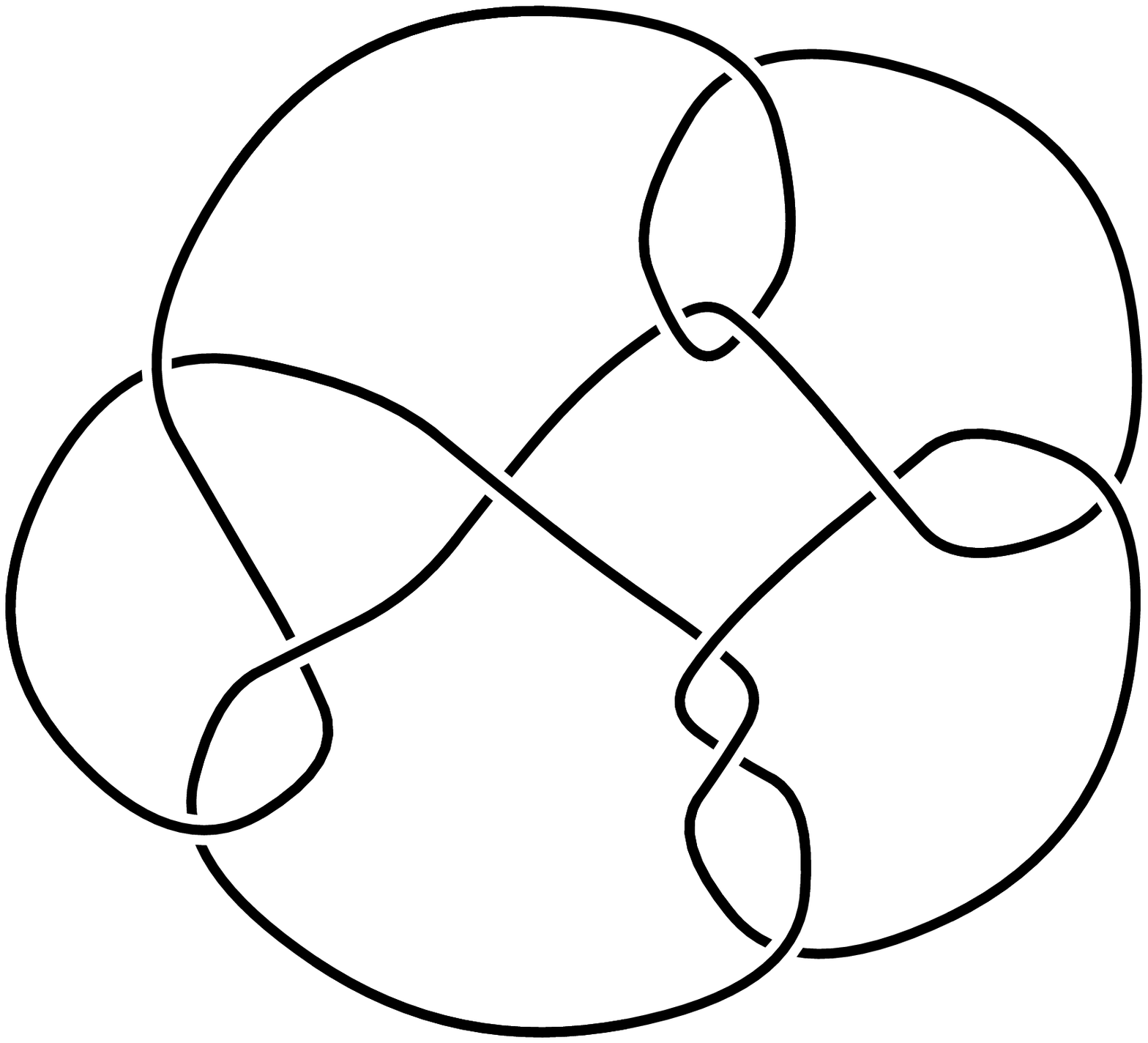}
    &
    \includegraphics[width=75pt]{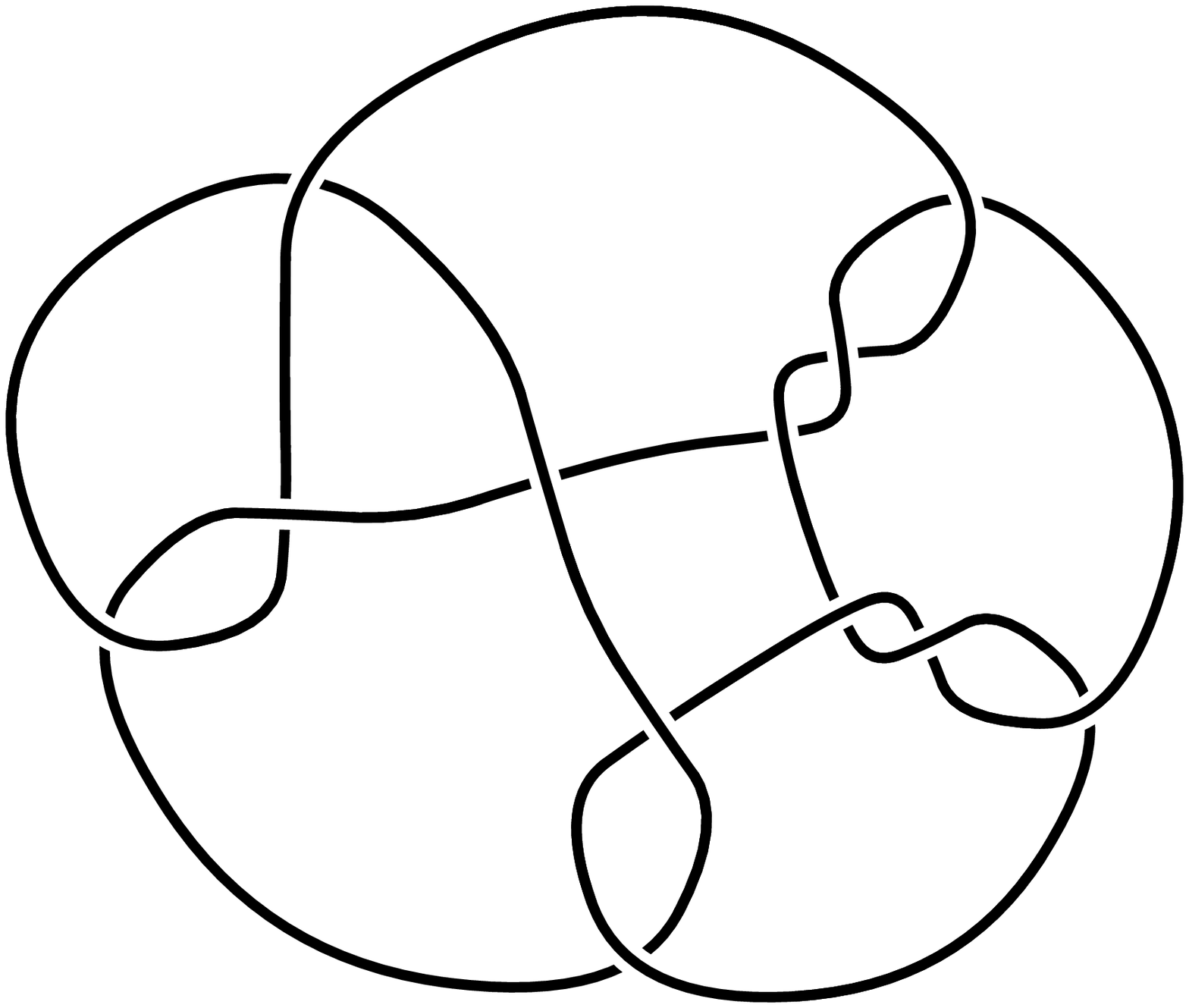}
    &
    \includegraphics[width=75pt]{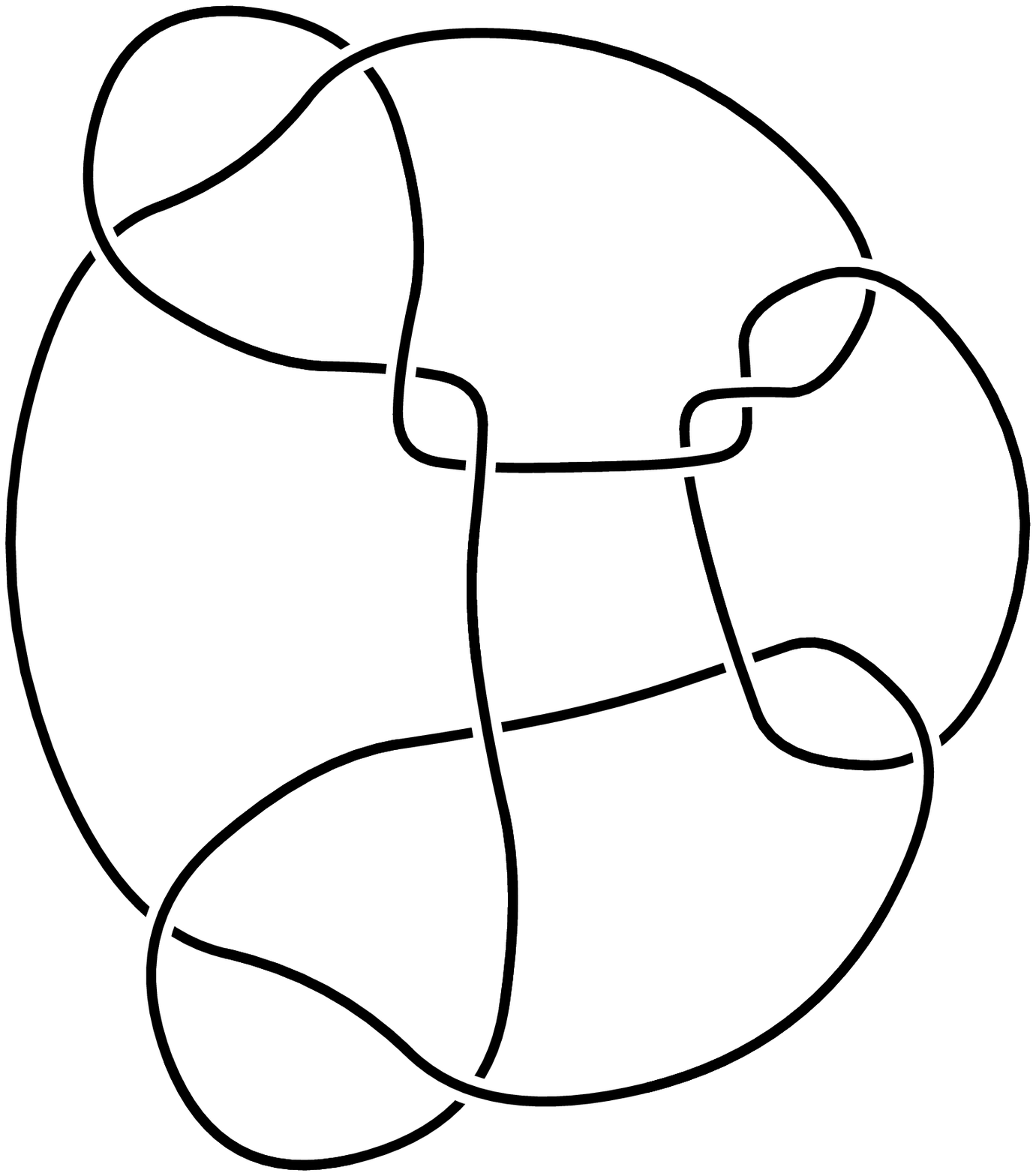}
    \\[-10pt]
    $12^N_{64}$ & $12^N_{261}$ & $12^N_{67}$ & $12^N_{229}$
    \\[10pt]
    \hline
    &&&\\[-10pt]
    \includegraphics[width=75pt]{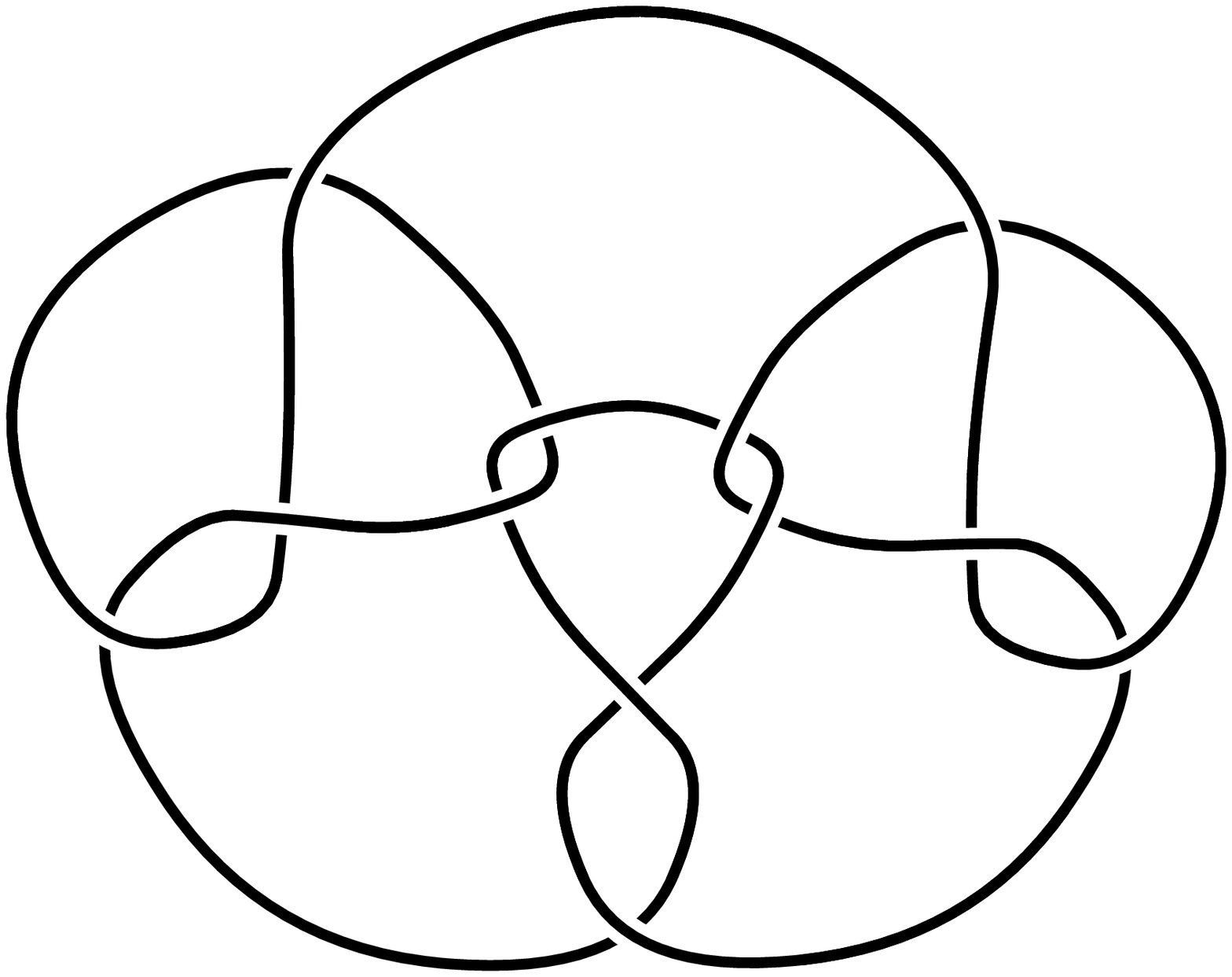}
    &
    \includegraphics[width=75pt]{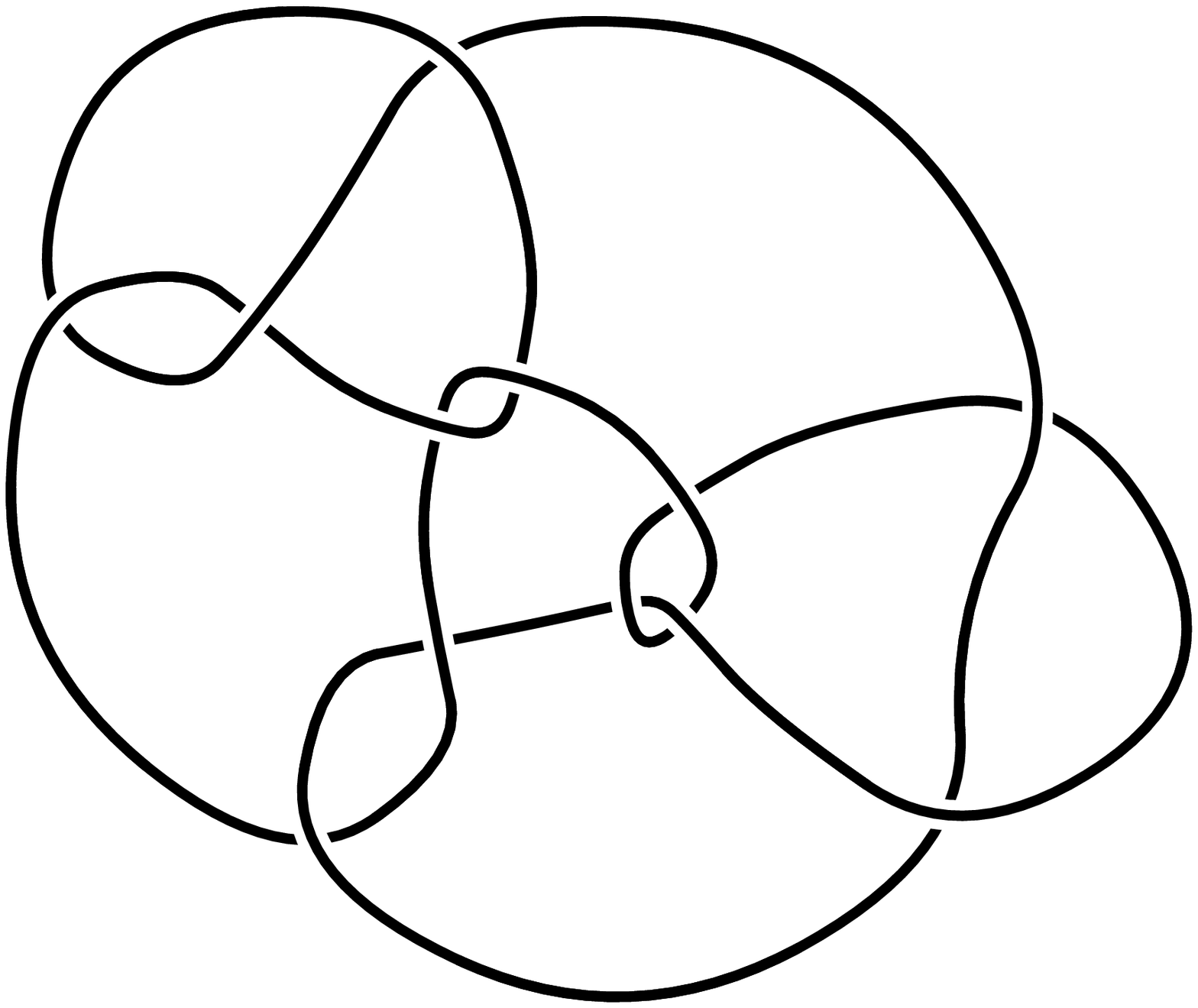}
    &
    \includegraphics[width=75pt]{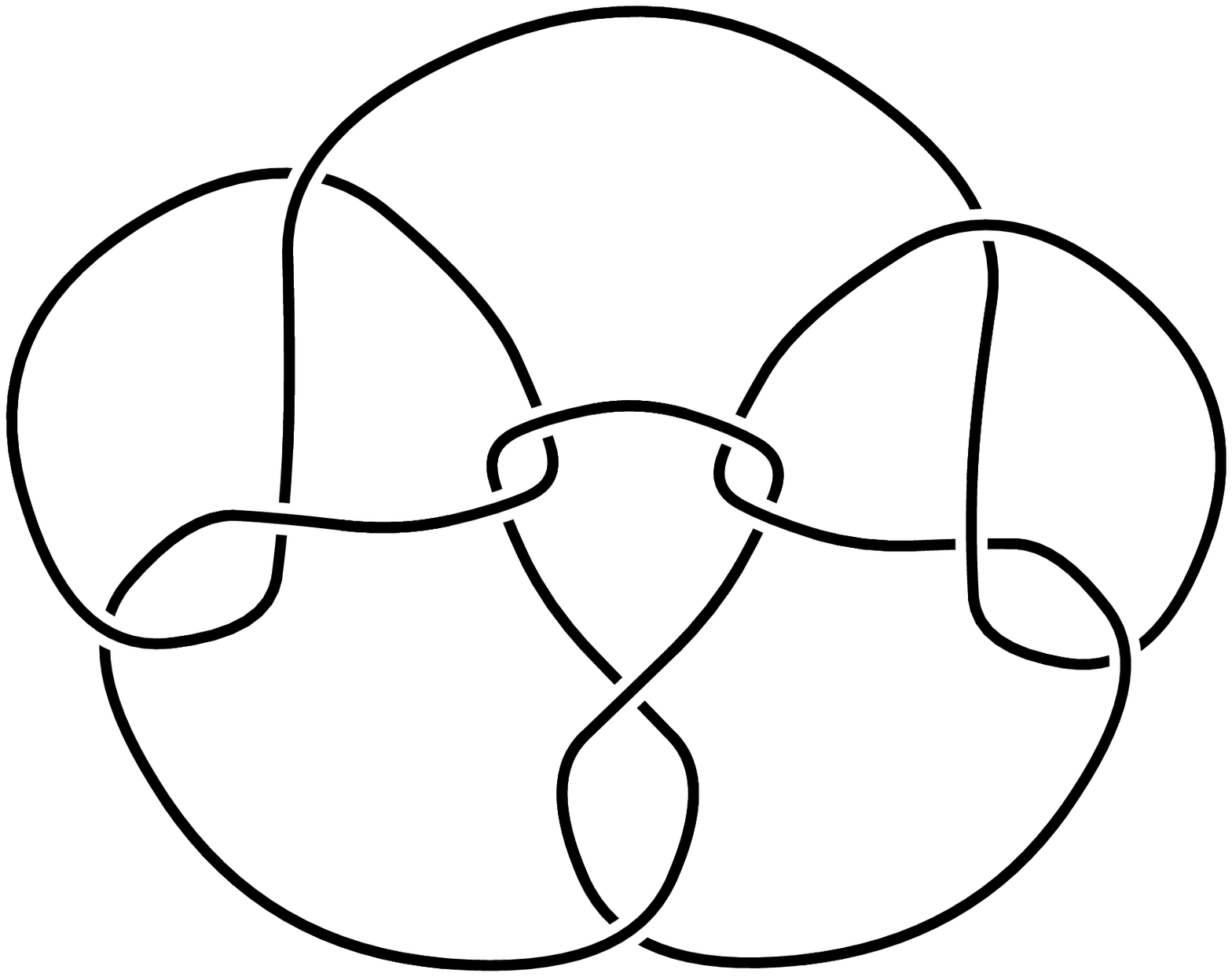}
    &
    \includegraphics[width=75pt]{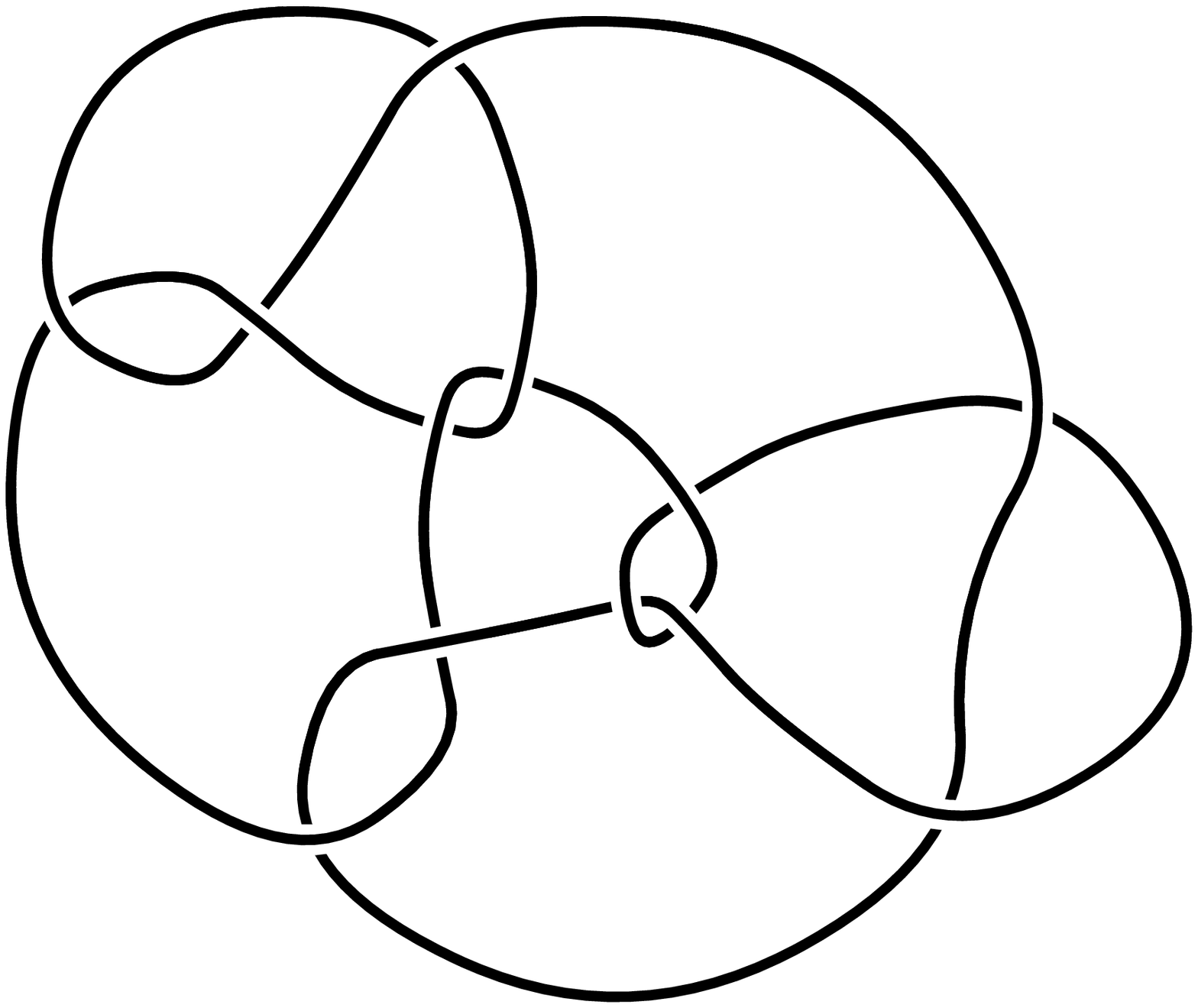}
    \\[-10pt]
    $12^N_{85}$ & $12^N_{130}$ & $12^N_{86}$ & $12^N_{131}$
    \\[10pt]
    \hline
    &&&\\[-10pt]
    \includegraphics[width=75pt]{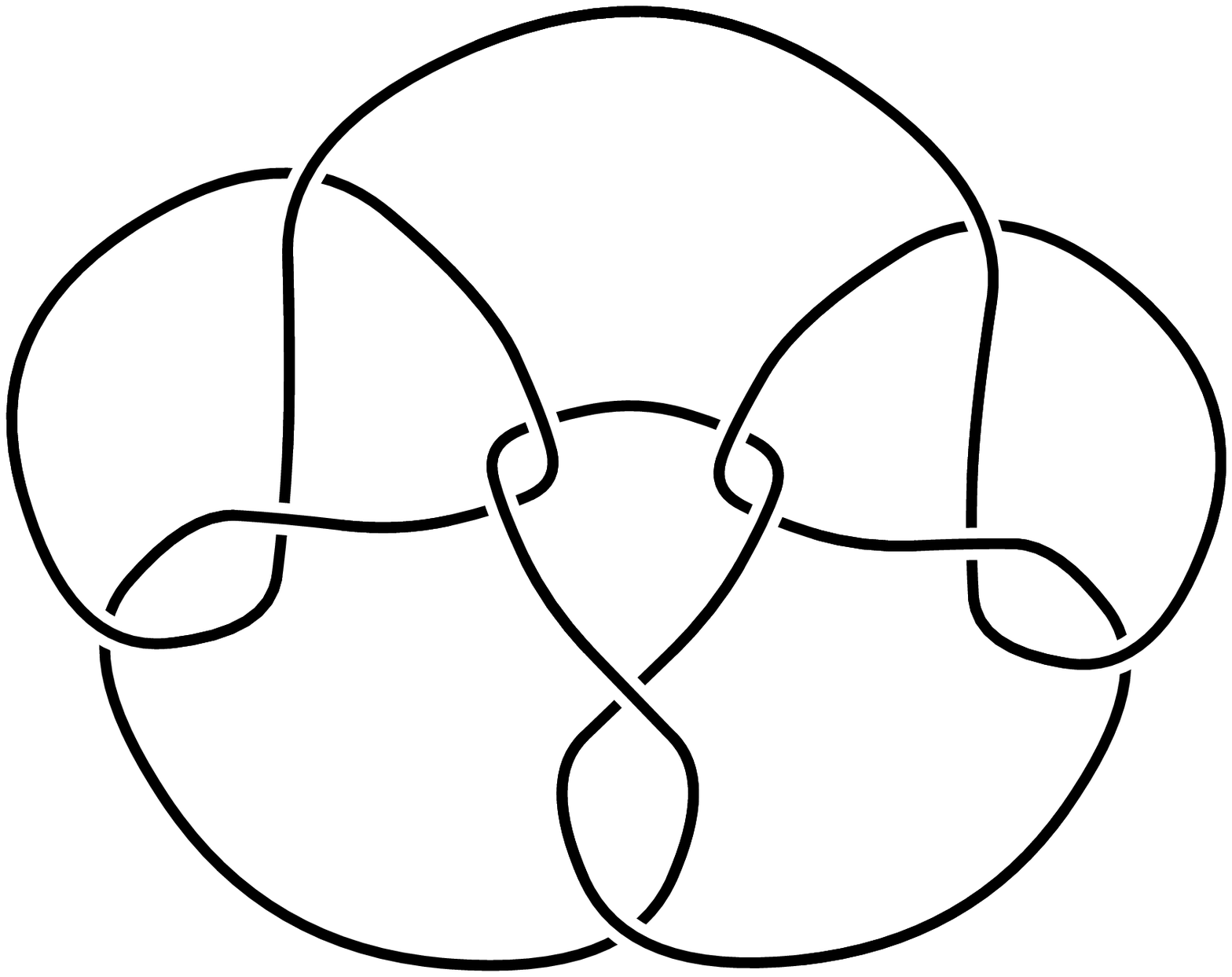}
    &
    \includegraphics[width=75pt]{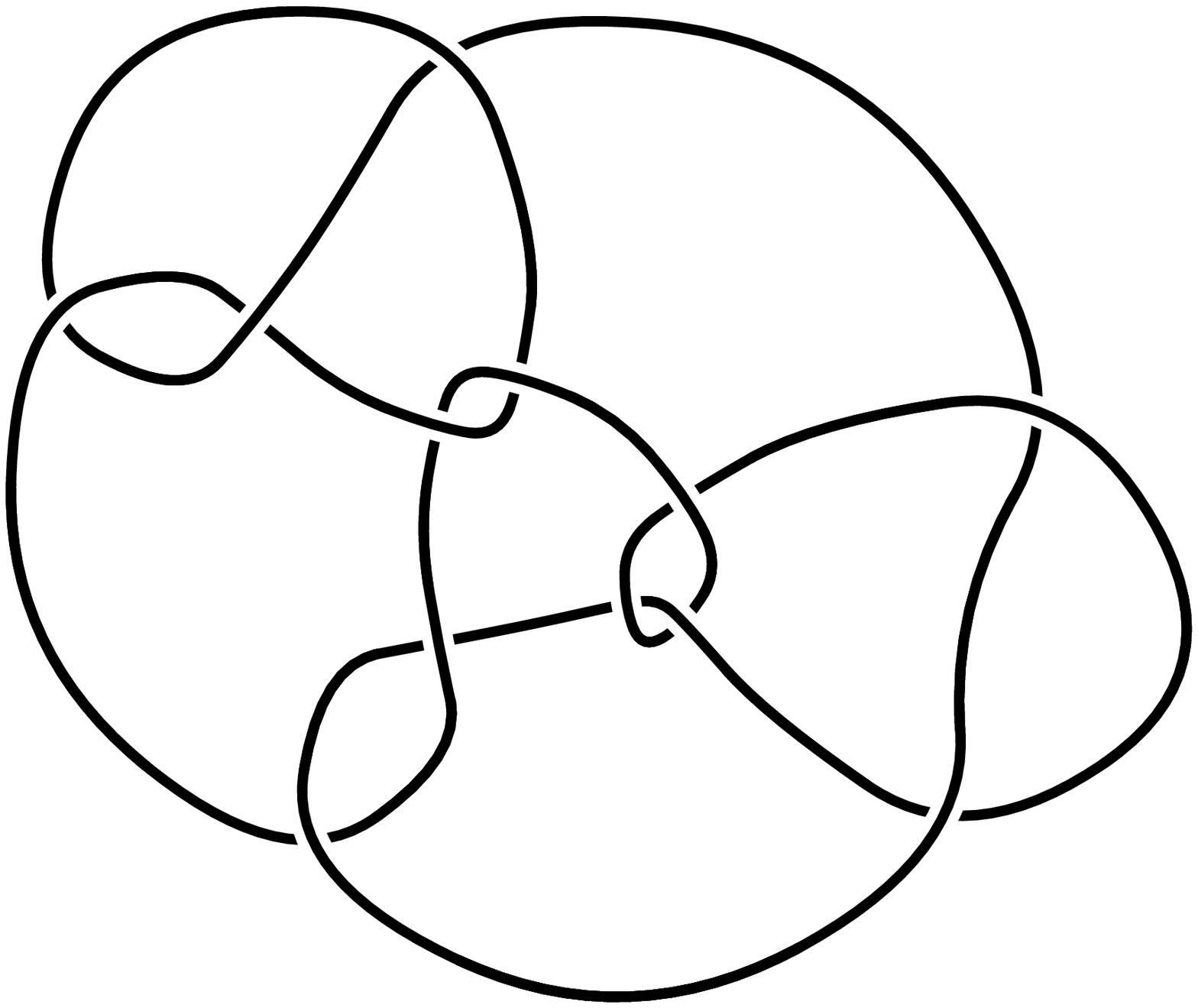}
    &
    \includegraphics[width=75pt]{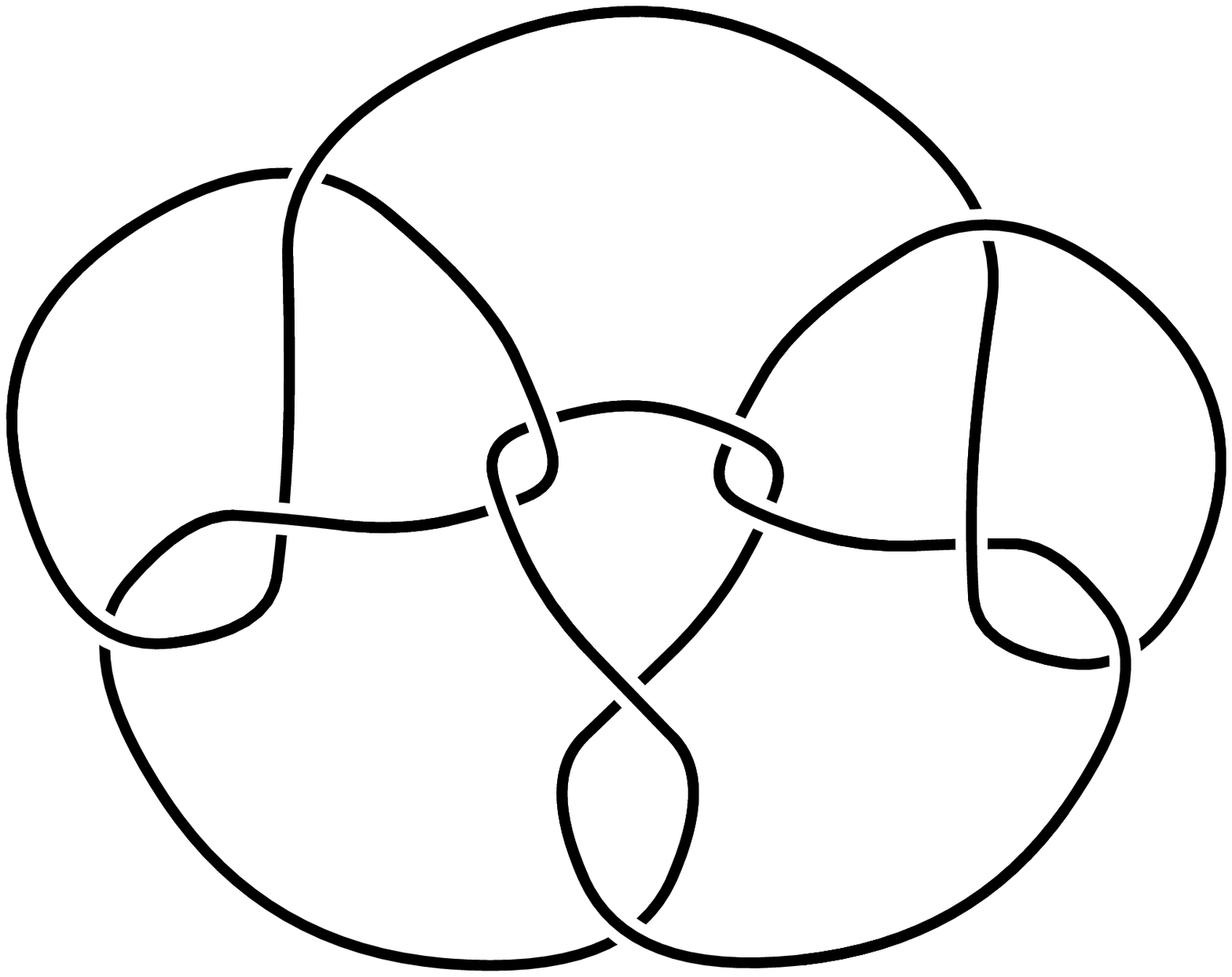}
    &
    \includegraphics[width=75pt]{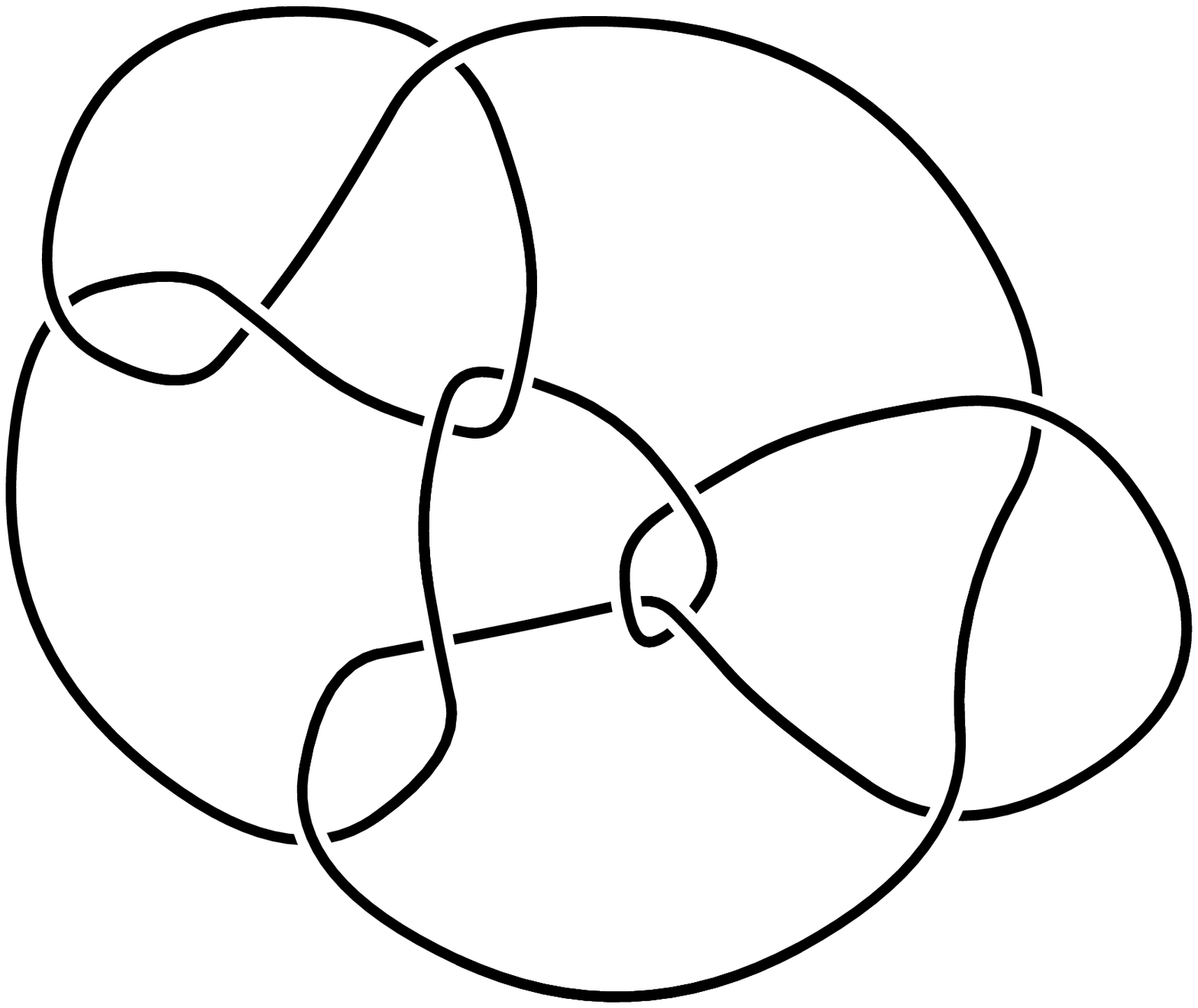}
    \\[-10pt]
    $12^N_{87}$ & $12^N_{132}$ & $12^N_{88}$ & $12^N_{133}$
    \\[10pt]
    \hline
    &&&\\[-10pt]
    \includegraphics[width=75pt]{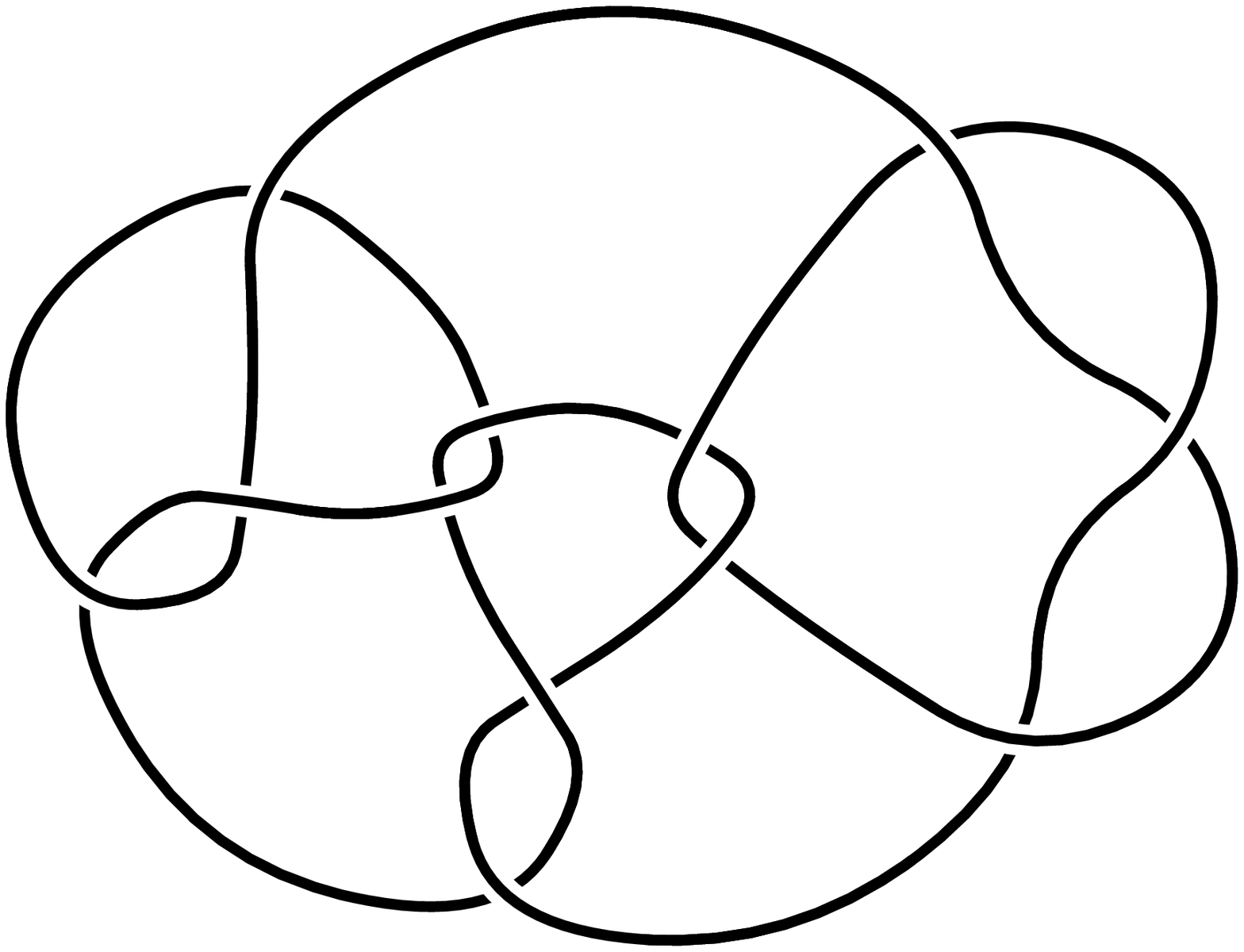}
    &
    \includegraphics[width=75pt]{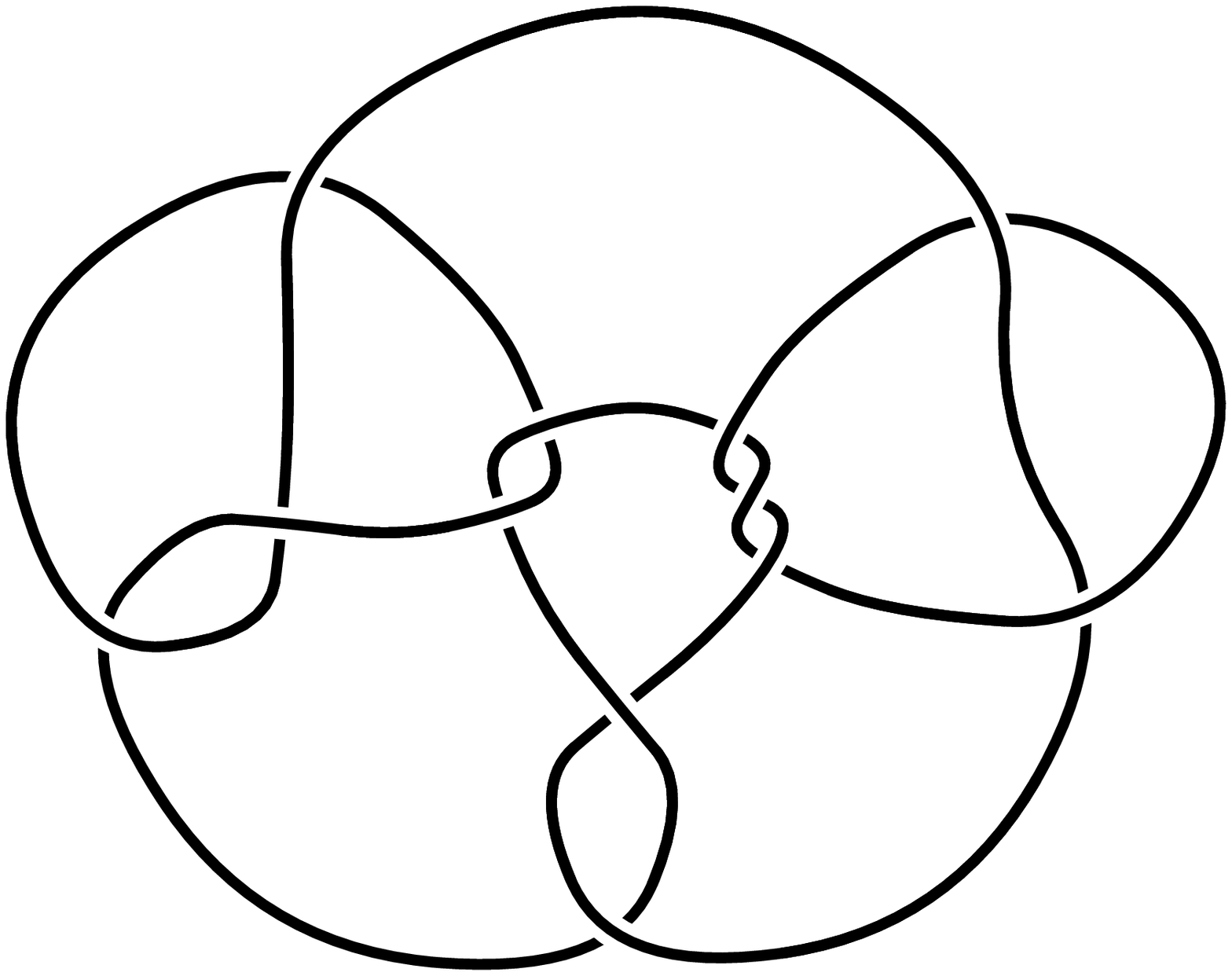}
    &
    \includegraphics[width=75pt]{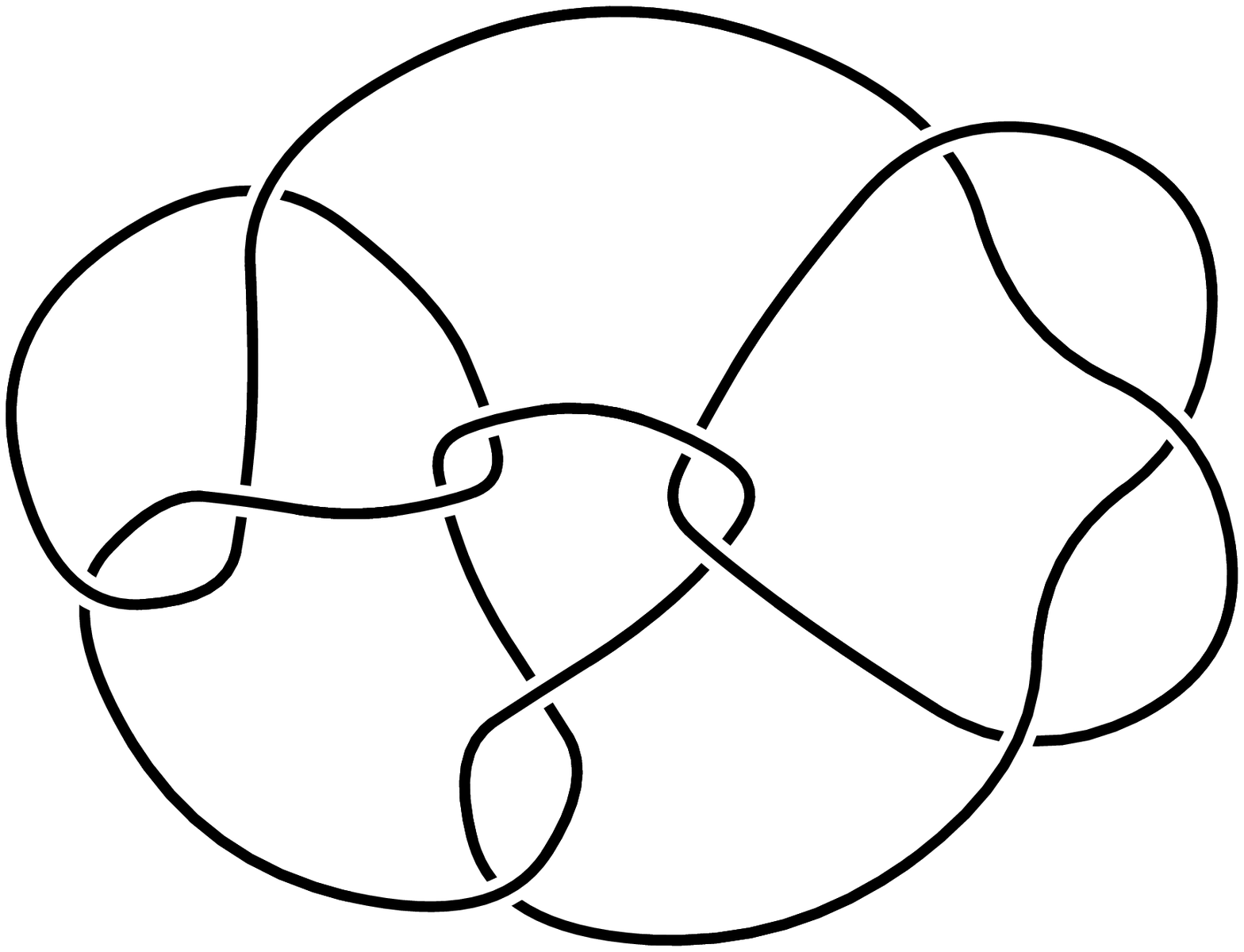}
    &
    \includegraphics[width=75pt]{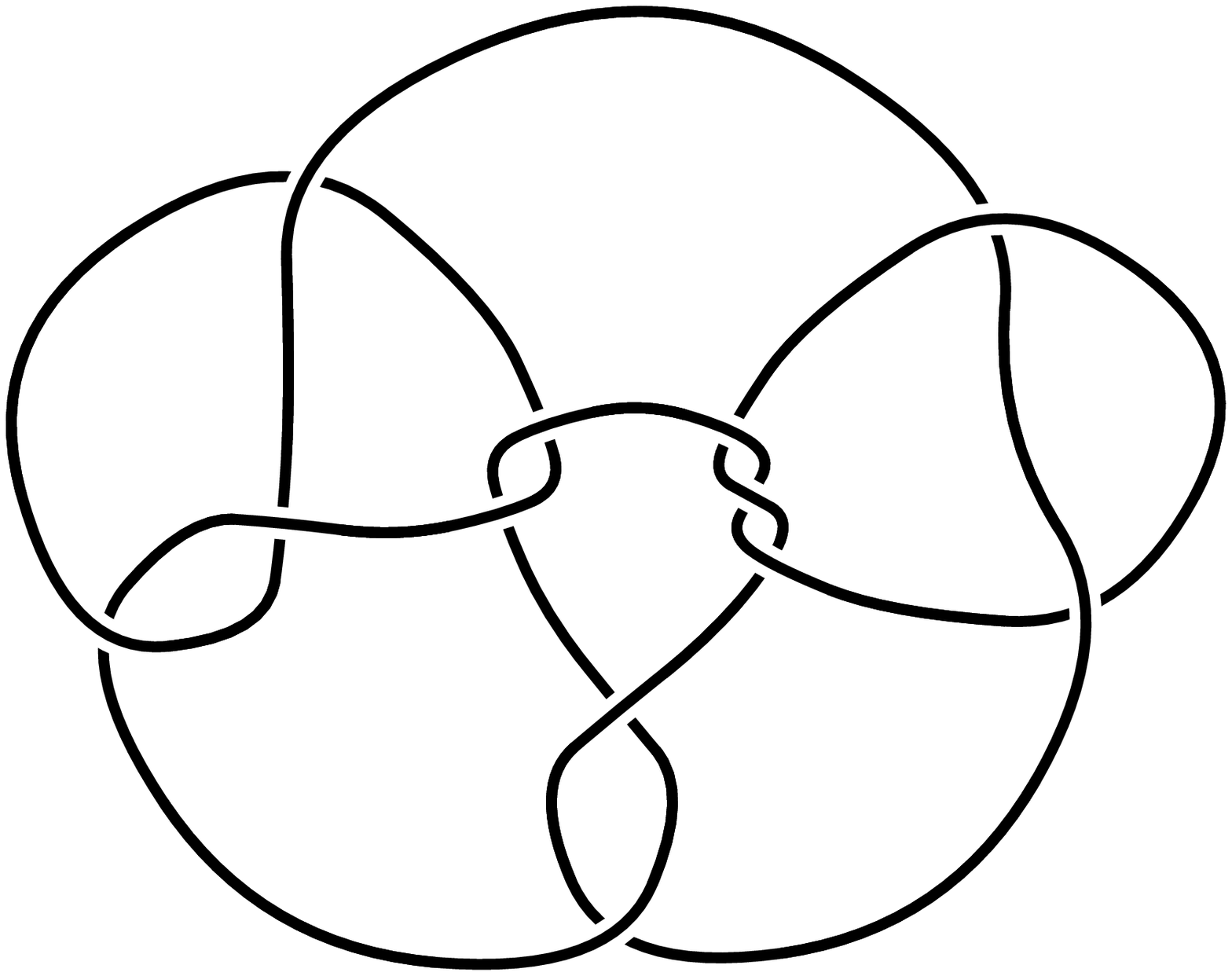}
    \\[-10pt]
    $12^N_{89}$ & $12^N_{134}$ & $12^N_{90}$ & $12^N_{135}$
    \\[10pt]
    \hline
    &&&\\[-10pt]
    \includegraphics[width=75pt]{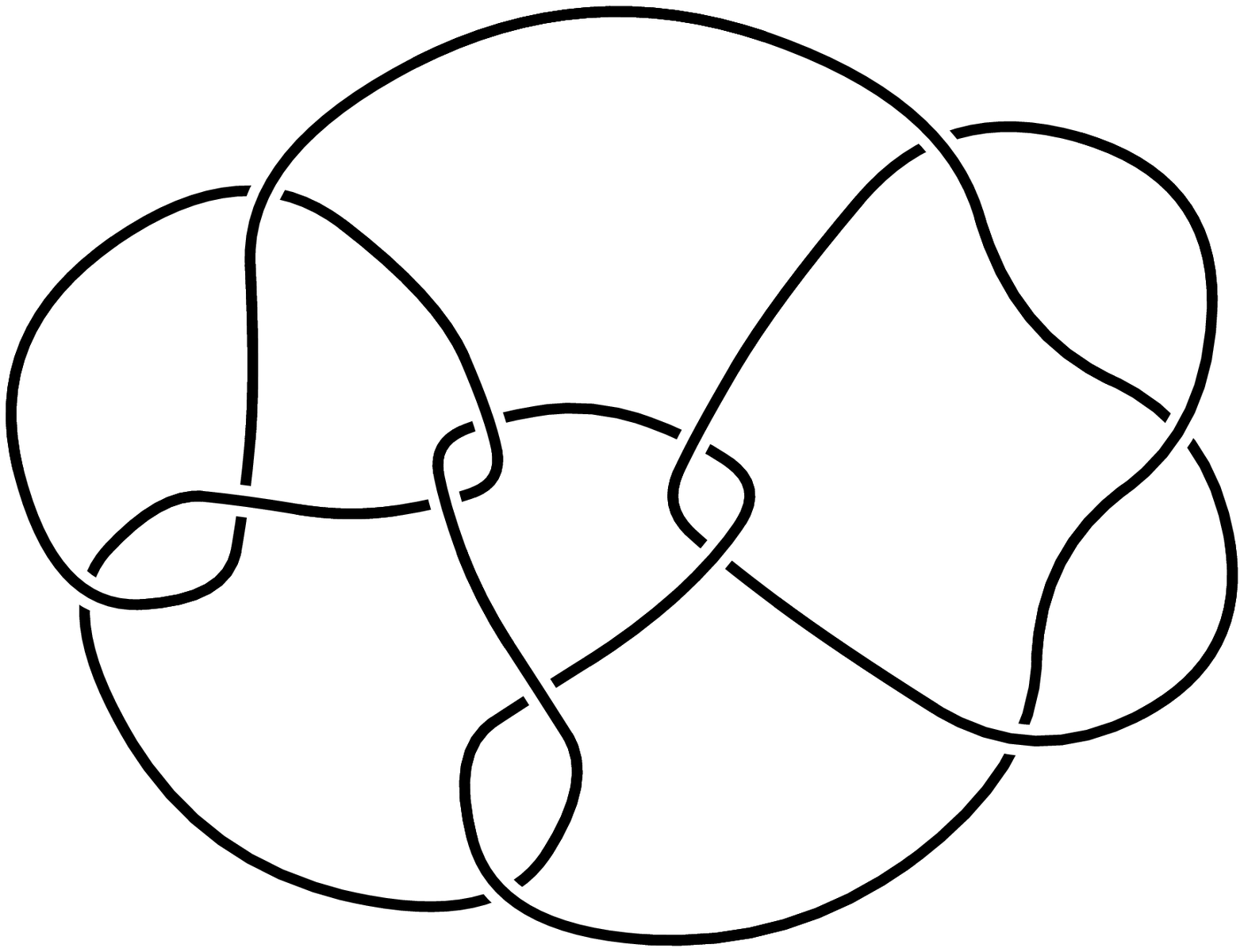}
    &
    \includegraphics[width=75pt]{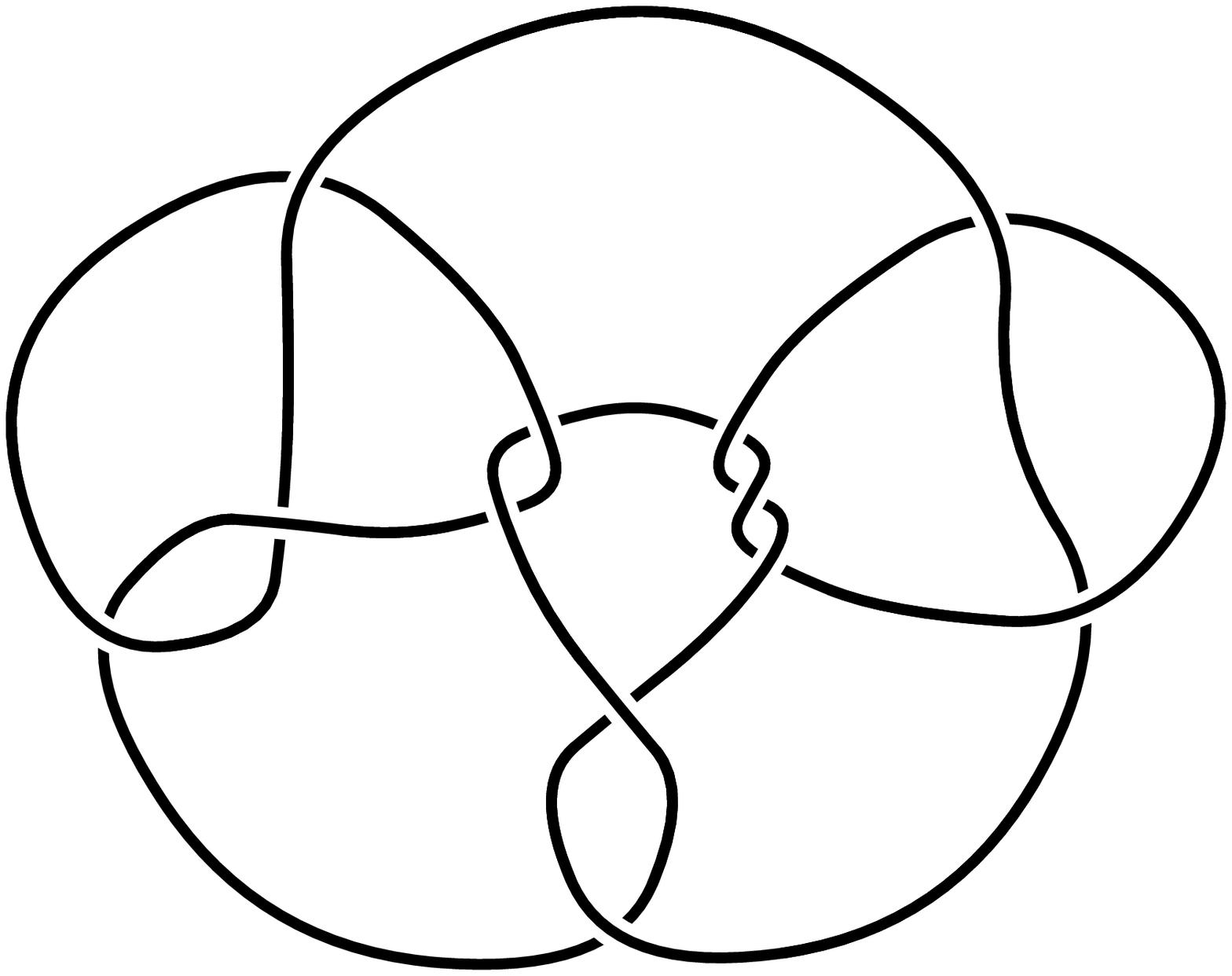}
    &
    \includegraphics[width=75pt]{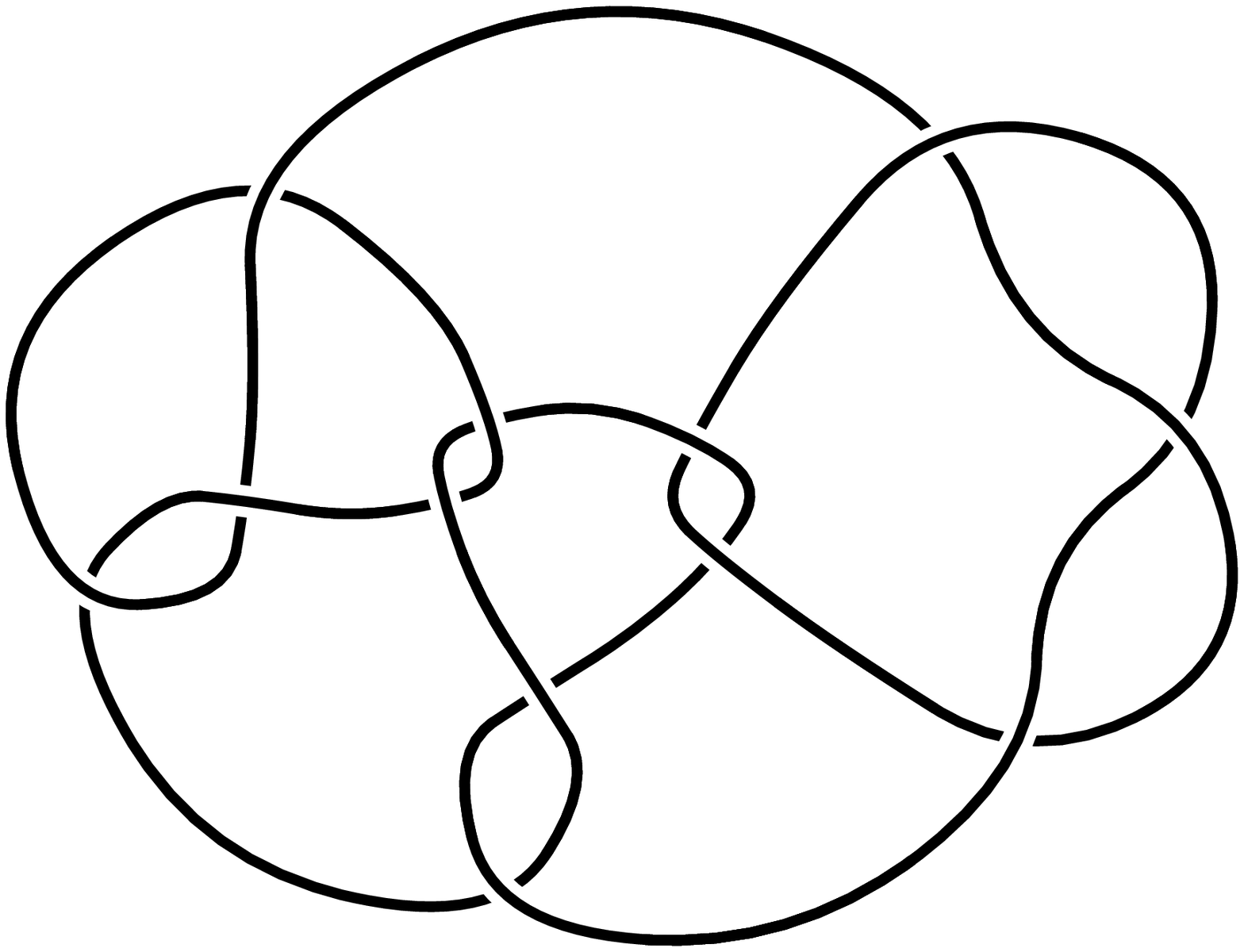}
    &
    \includegraphics[width=75pt]{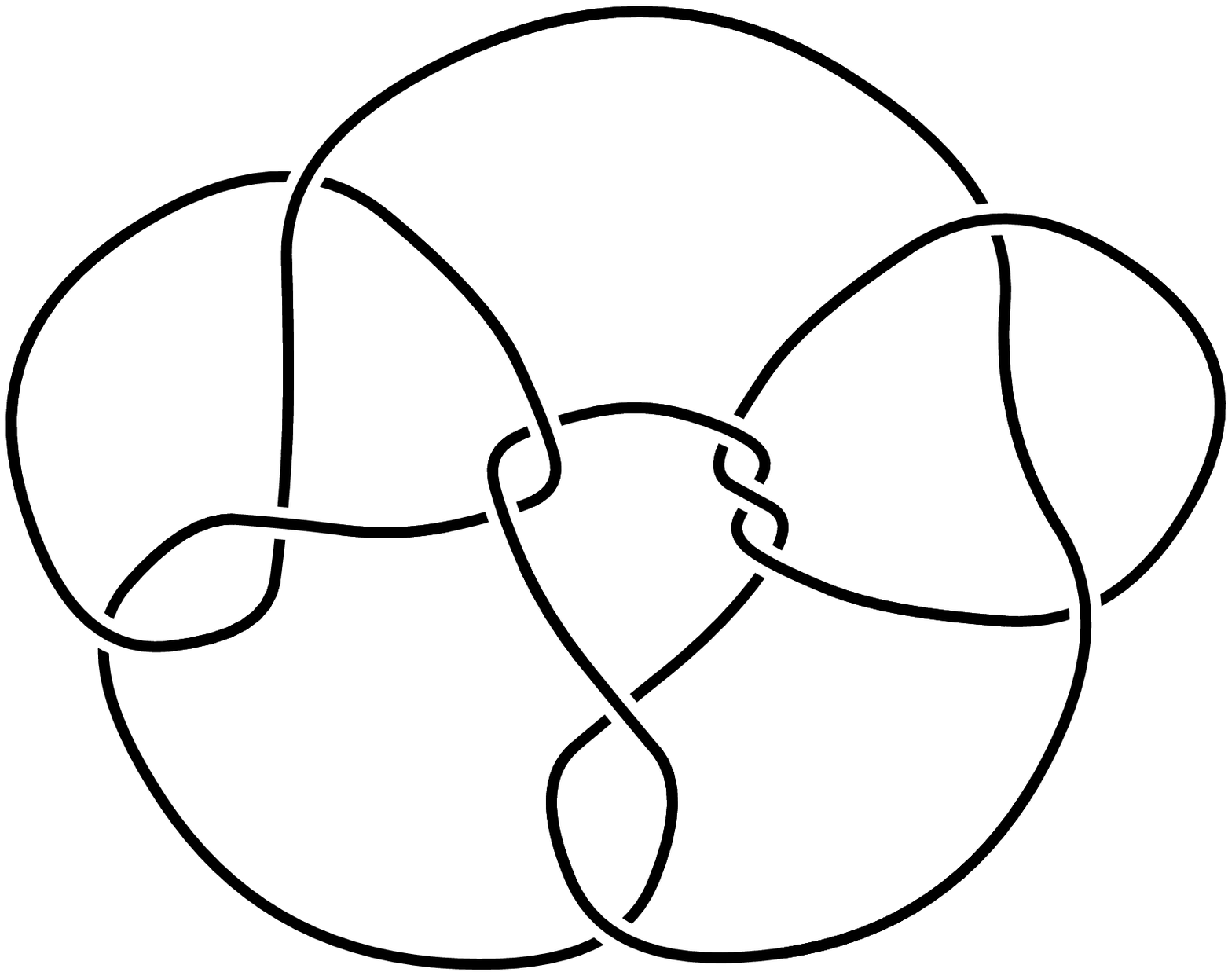}
    \\[-10pt]
    $12^N_{91}$ & $12^N_{136}$ & $12^N_{92}$ & $12^N_{137}$
  \end{tabular}
  \caption{Nonalternating $12$-crossing mutant cliques 2/6}
  \end{centering}
\end{figure}

\begin{figure}[htbp]
  \begin{centering}
  \begin{tabular}{cc@{\hspace{10pt}}|@{\hspace{10pt}}cc}
    \includegraphics[width=75pt]{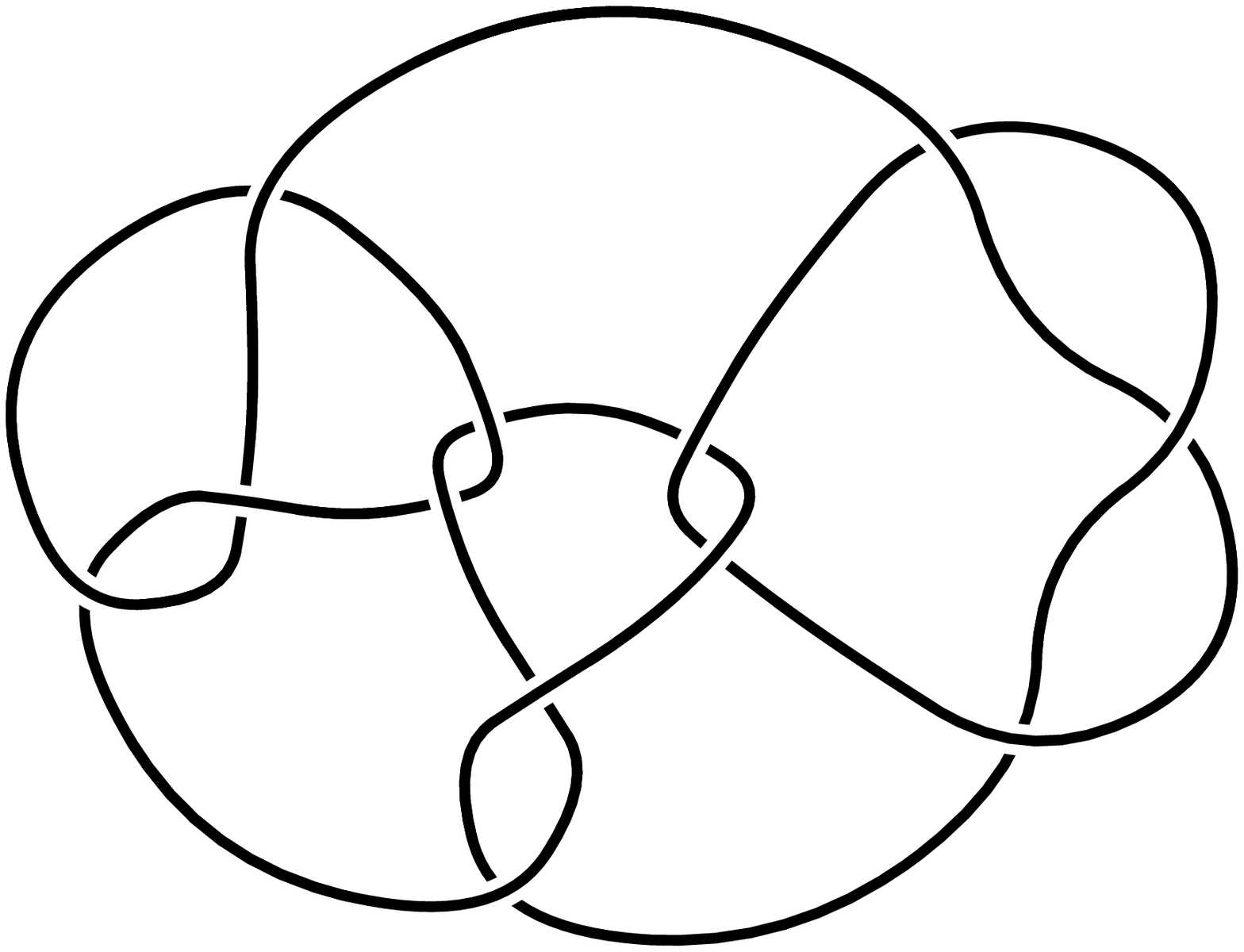}
    &
    \includegraphics[width=75pt]{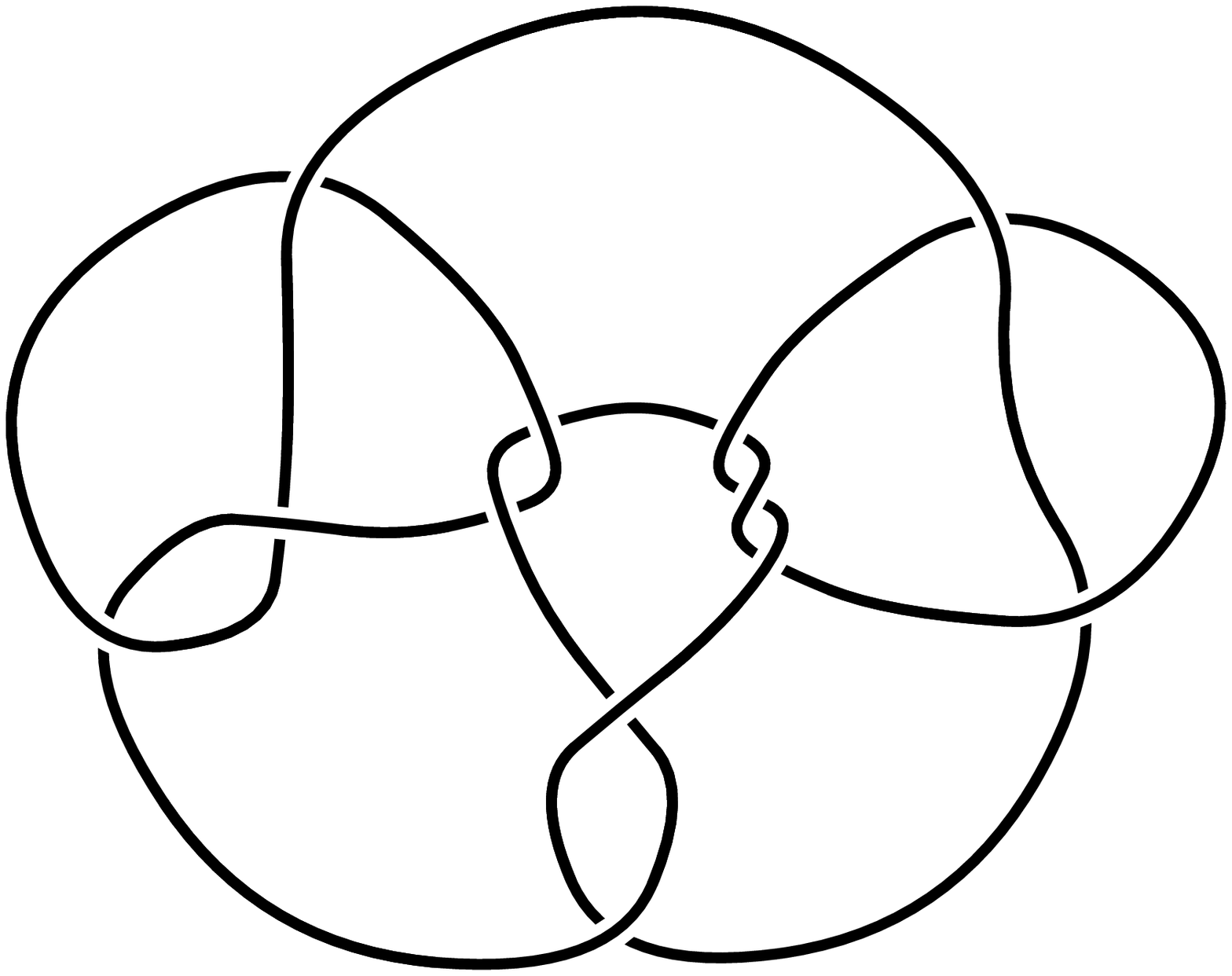}
    &
    \includegraphics[width=75pt]{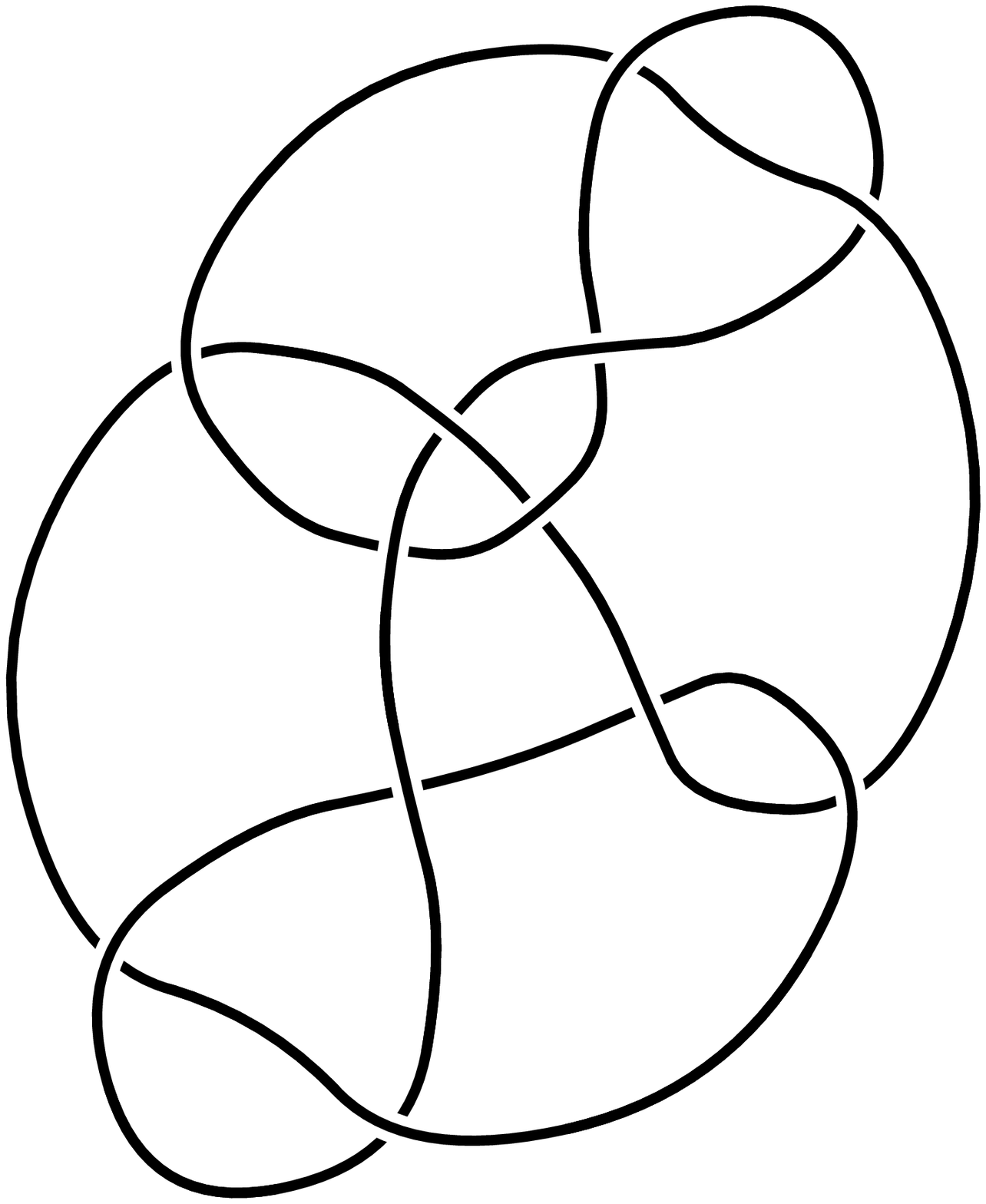}
    &
    \includegraphics[width=75pt]{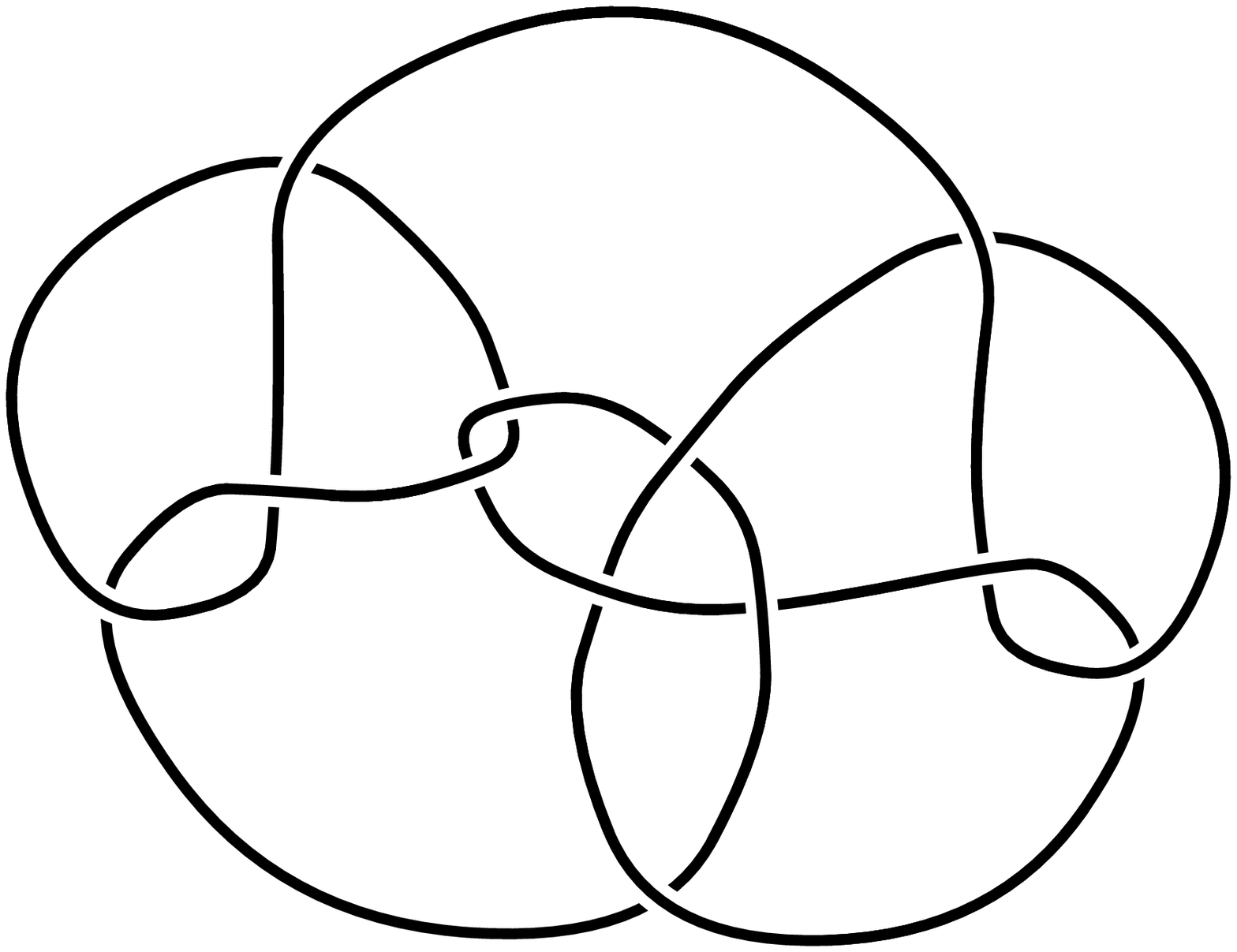}
    \\[-10pt]
    $12^N_{93}$ & $12^N_{138}$ & $12^N_{98}$ & $12^N_{125}$
    \\[10pt]
    \hline
    &&&\\[-10pt]
    \includegraphics[width=75pt]{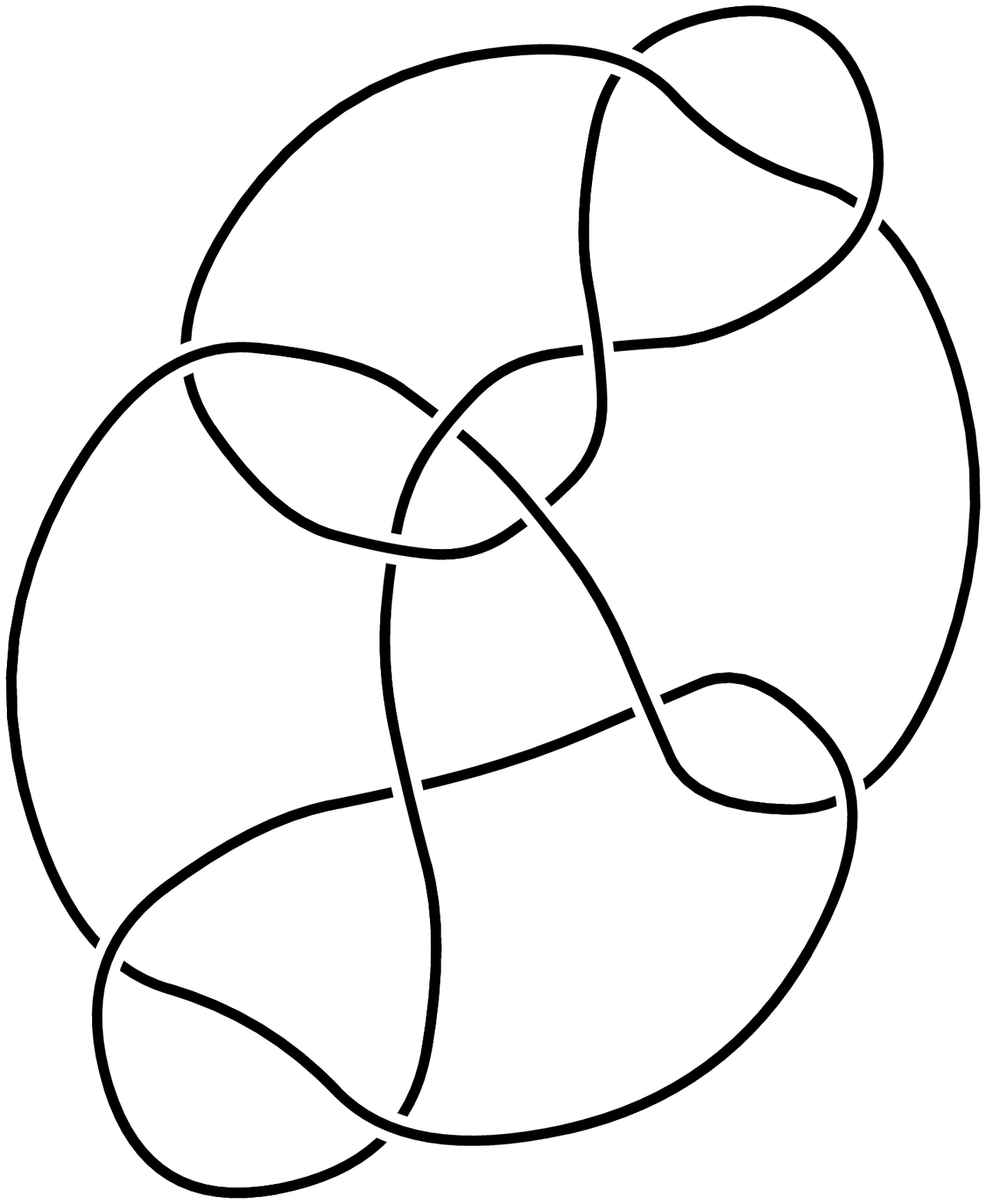}
    &
    \includegraphics[width=75pt]{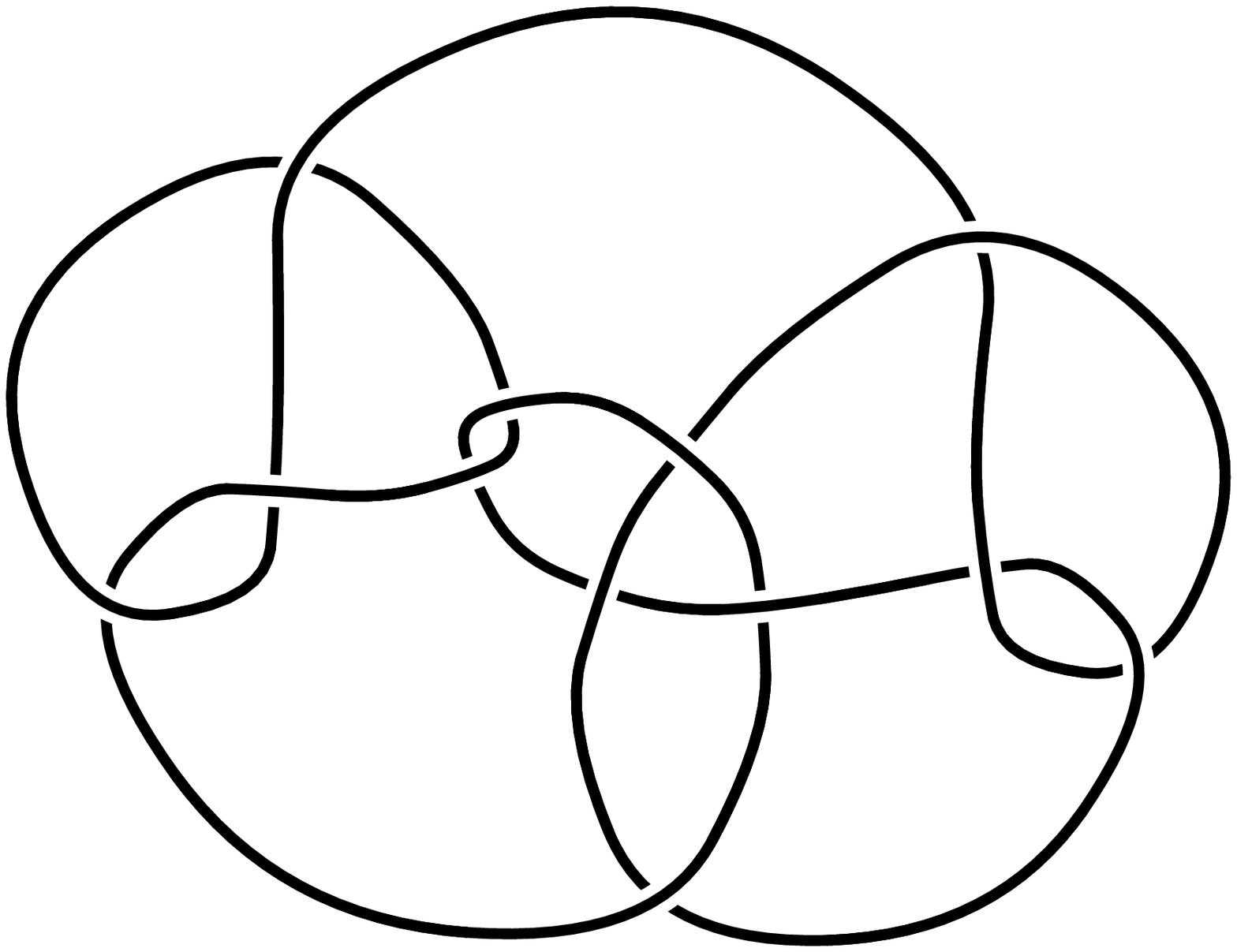}
    &
    \includegraphics[width=75pt]{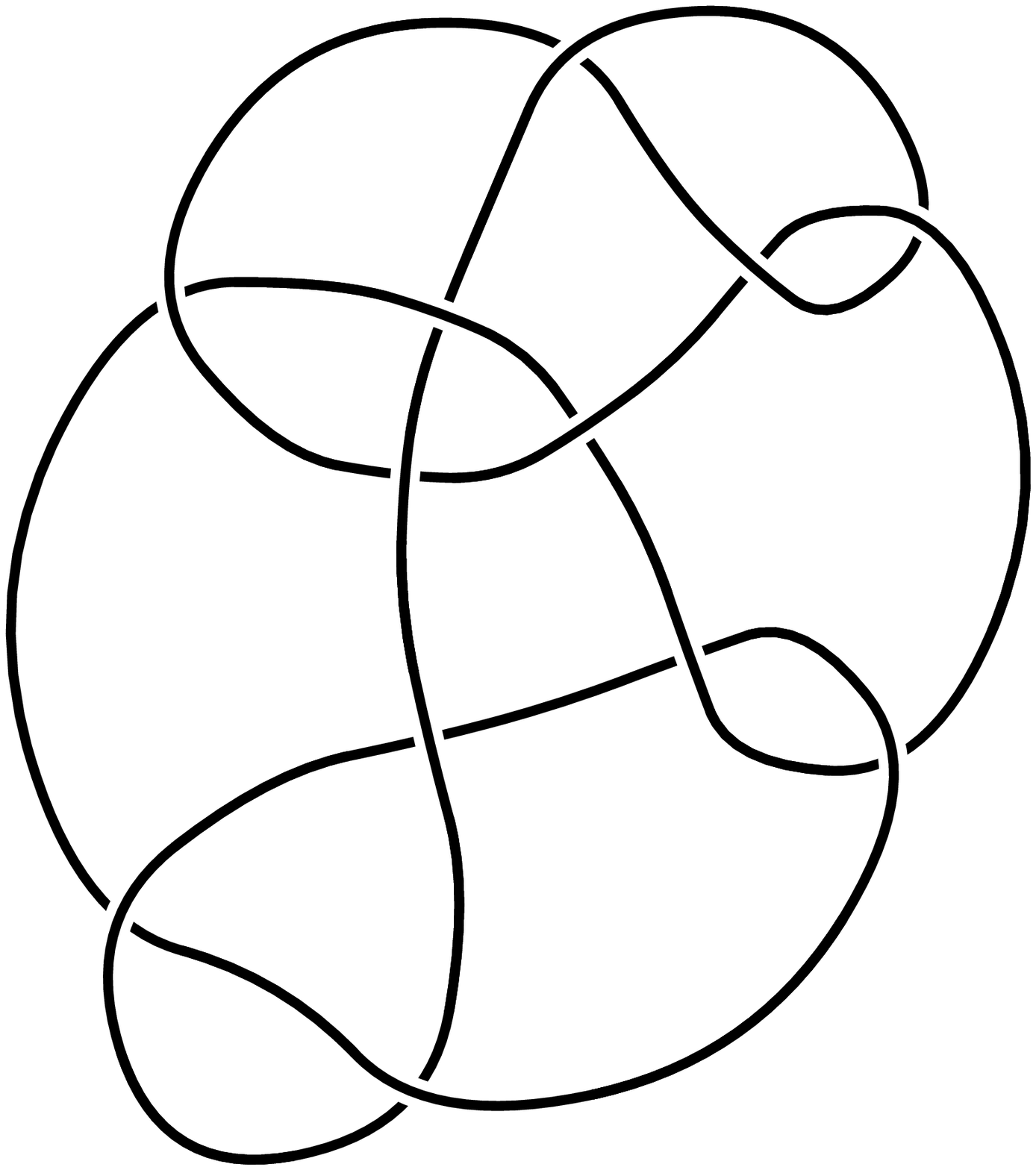}
    &
    \includegraphics[width=75pt]{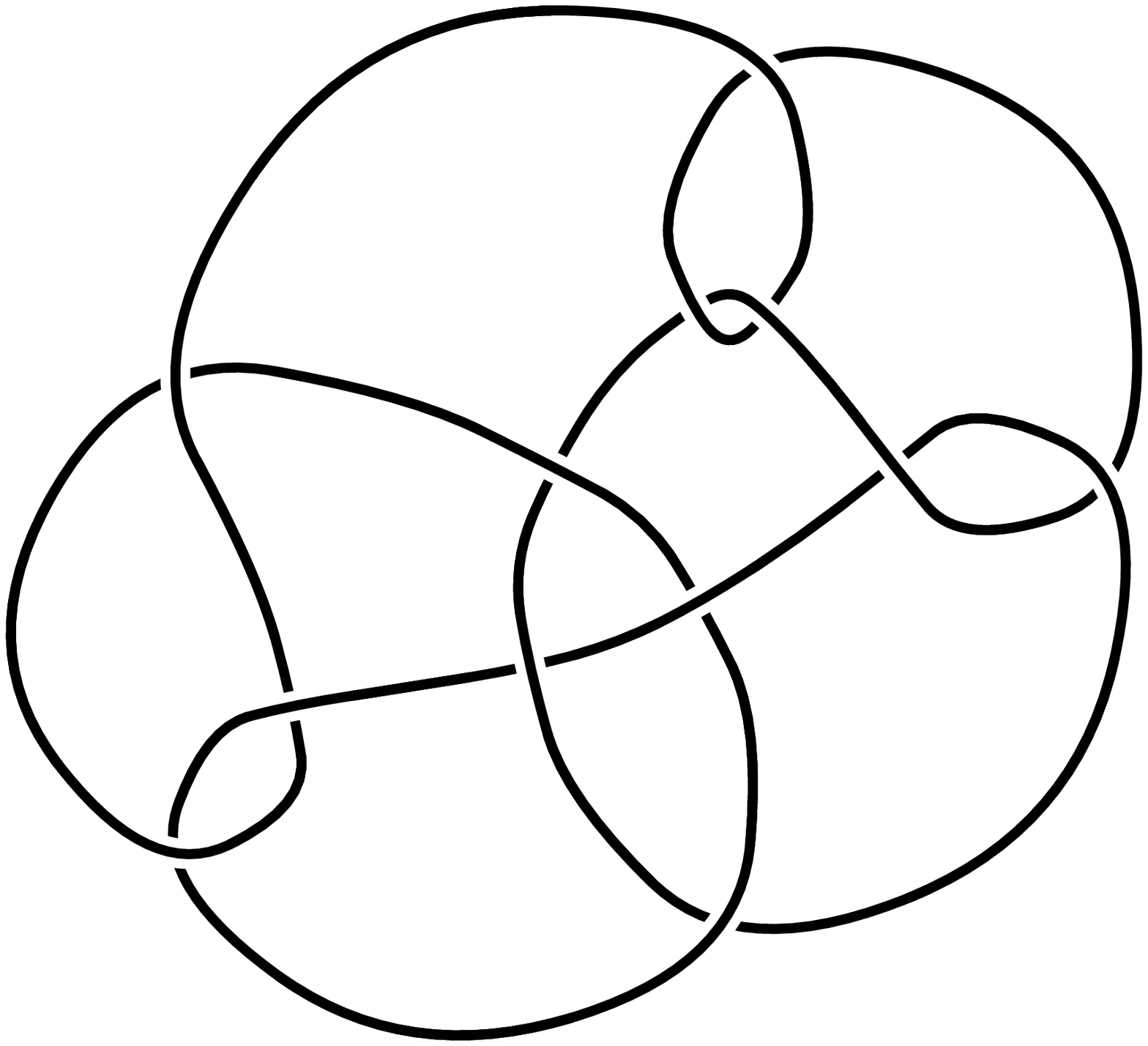}
    \\[-10pt]
    $12^N_{99}$ & $12^N_{126}$ & $12^N_{122}$ & $12^N_{127}$
    \\[10pt]
    \hline
    &&&\\[-10pt]
    \includegraphics[width=75pt]{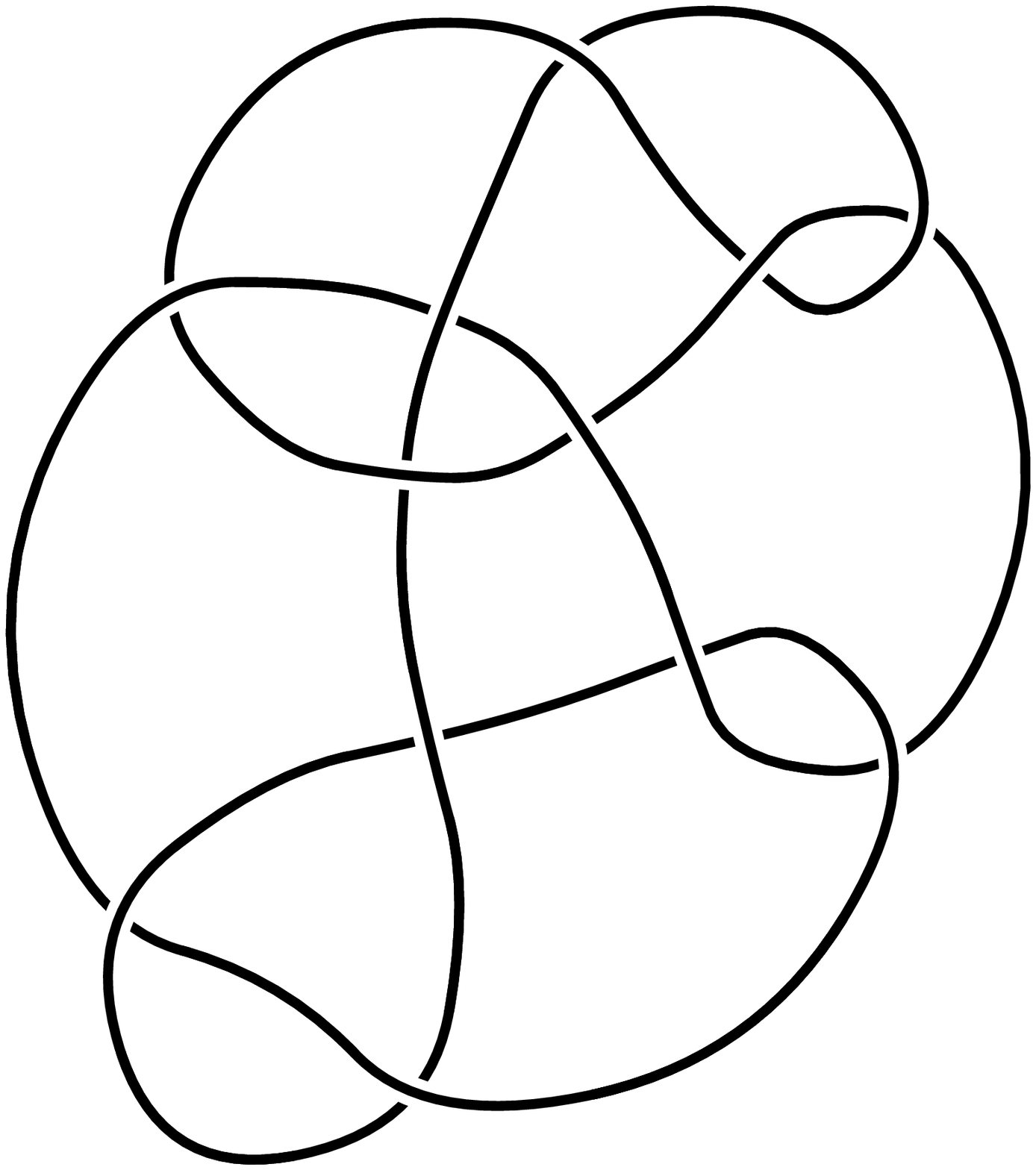}
    &
    \includegraphics[width=75pt]{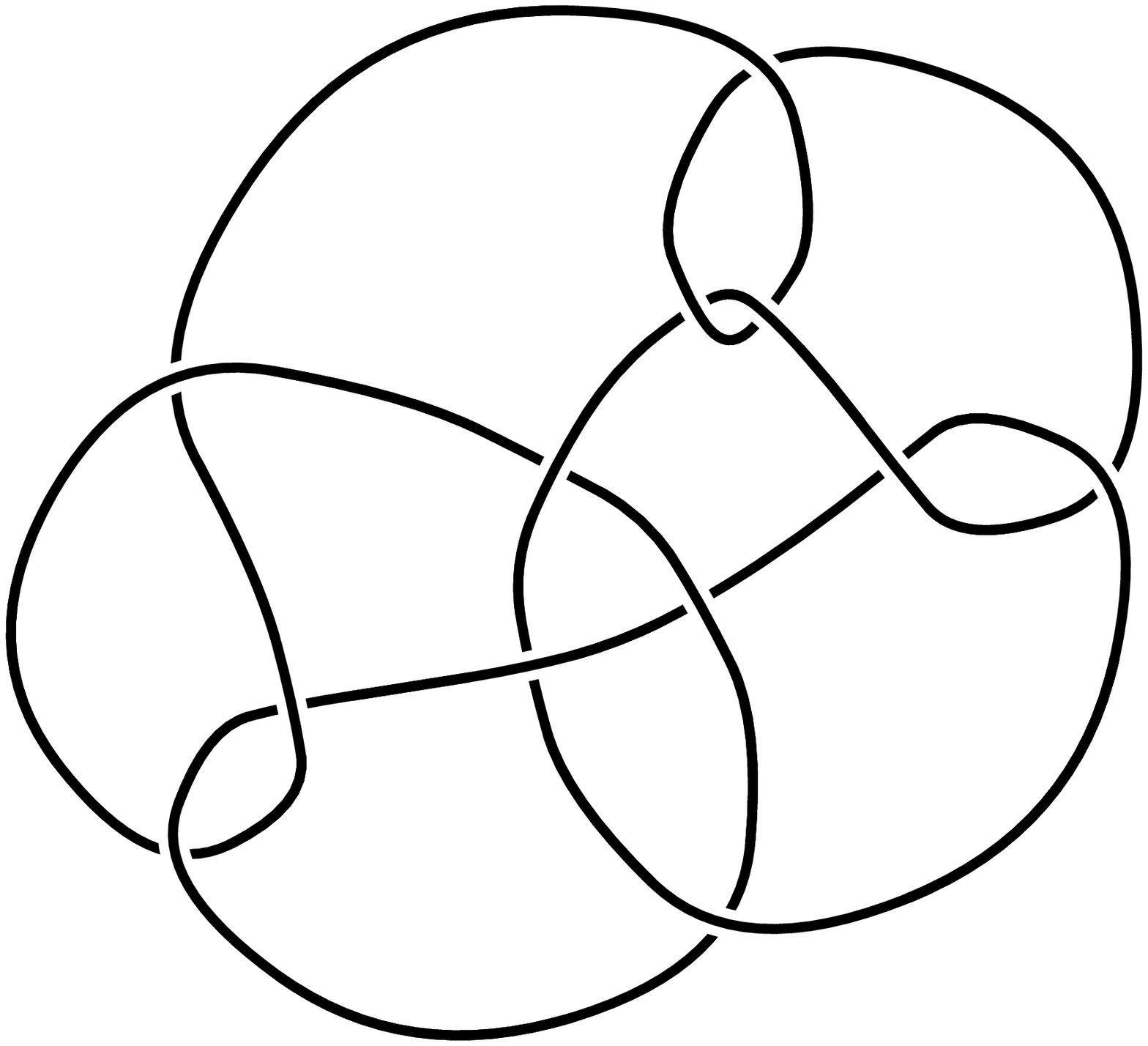}
    &
    \includegraphics[width=75pt]{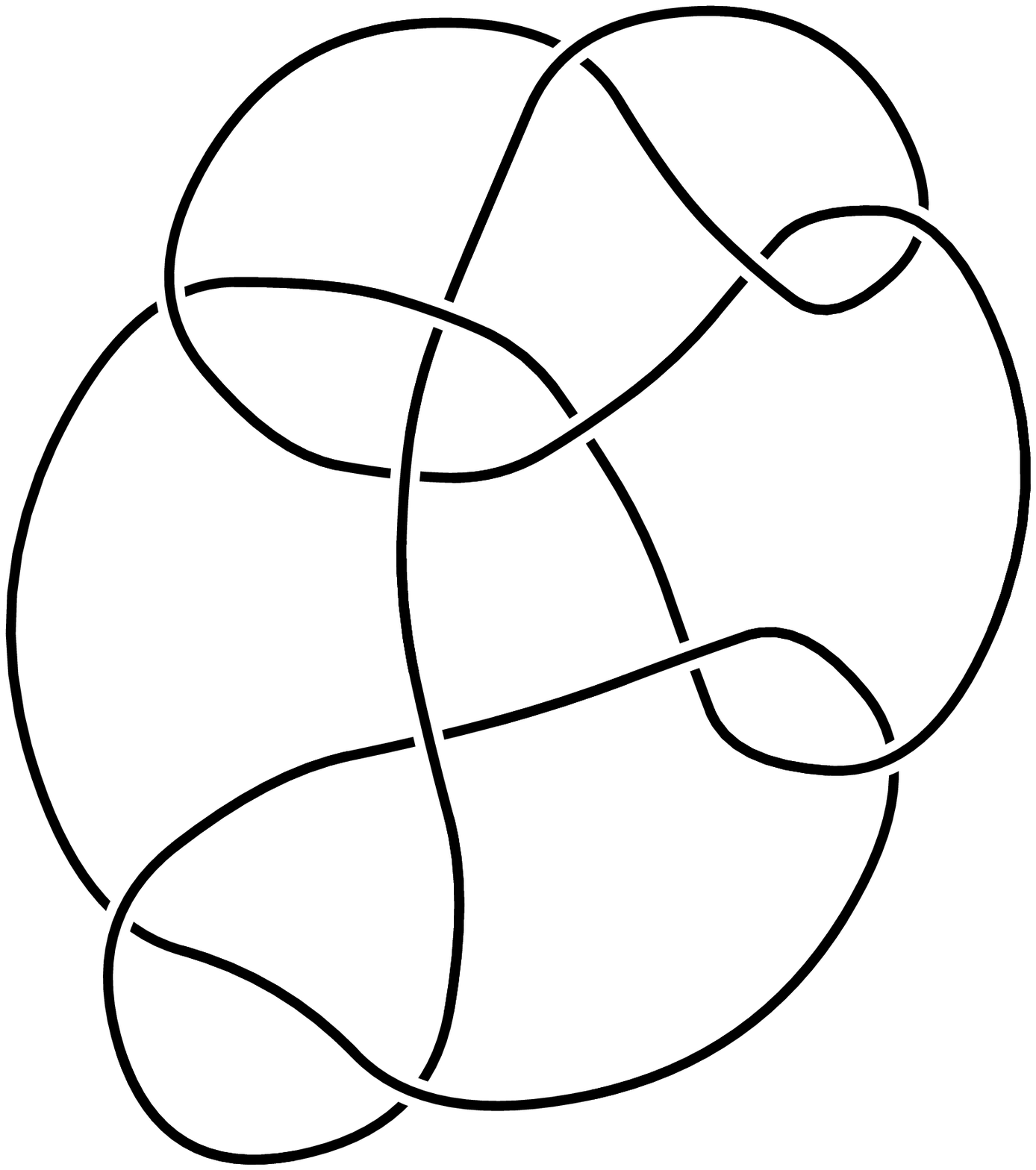}
    &
    \includegraphics[width=75pt]{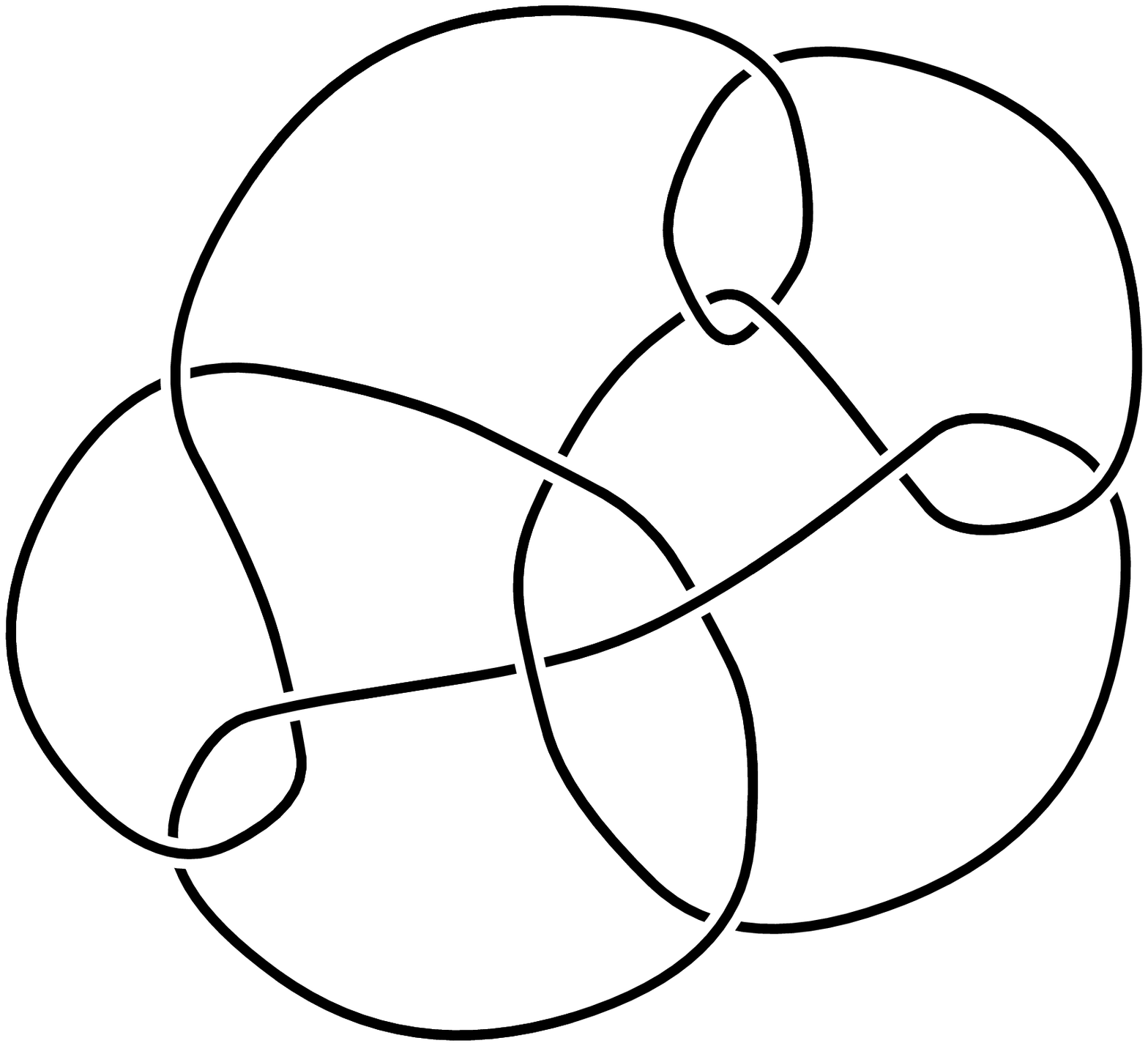}
    \\[-10pt]
    $12^N_{123}$ & $12^N_{128}$ & $12^N_{124}$ & $12^N_{129}$
    \\[10pt]
    \hline
    &&&\\[-10pt]
    \includegraphics[width=75pt]{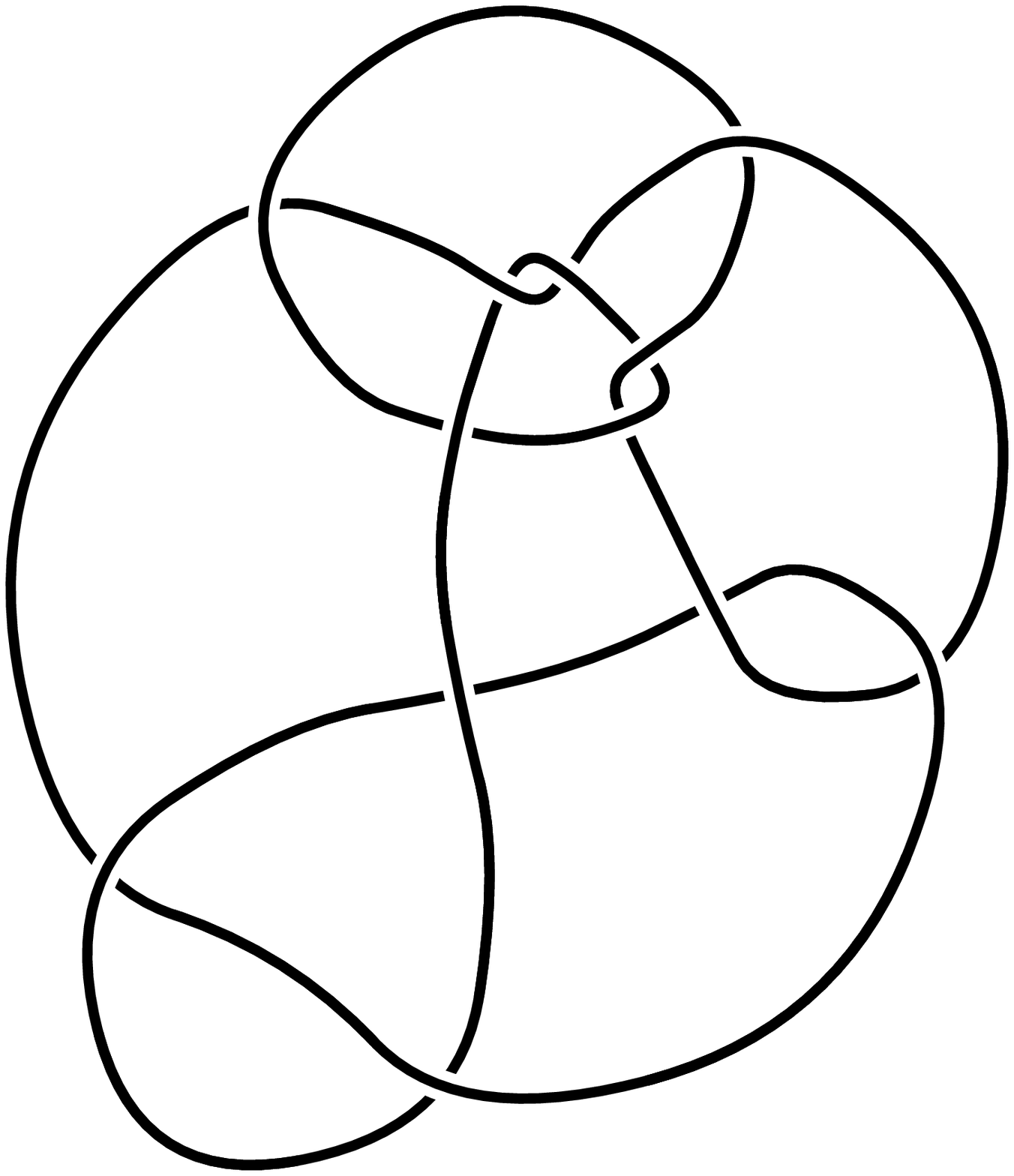}
    &
    \includegraphics[width=75pt]{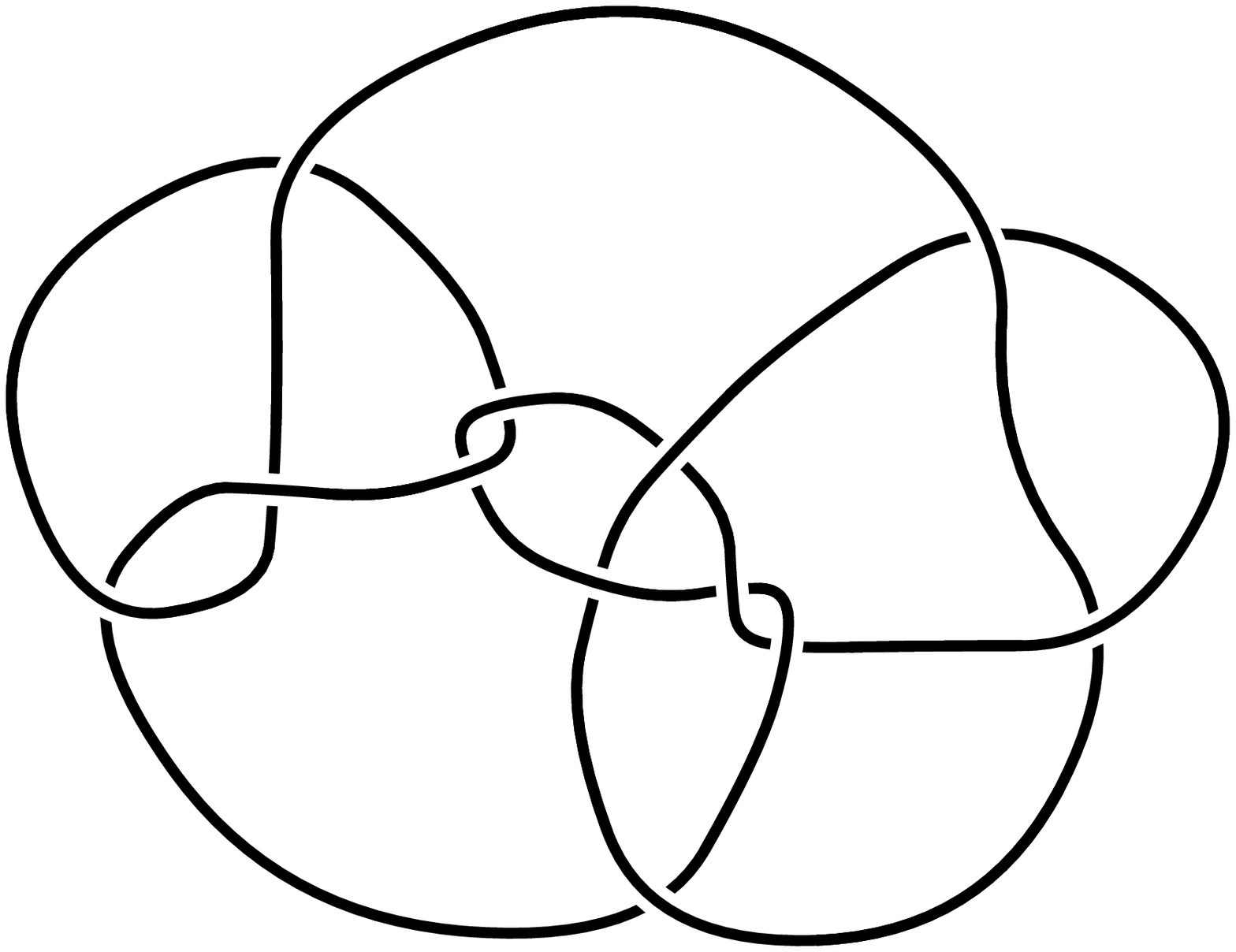}
    &
    \includegraphics[width=75pt]{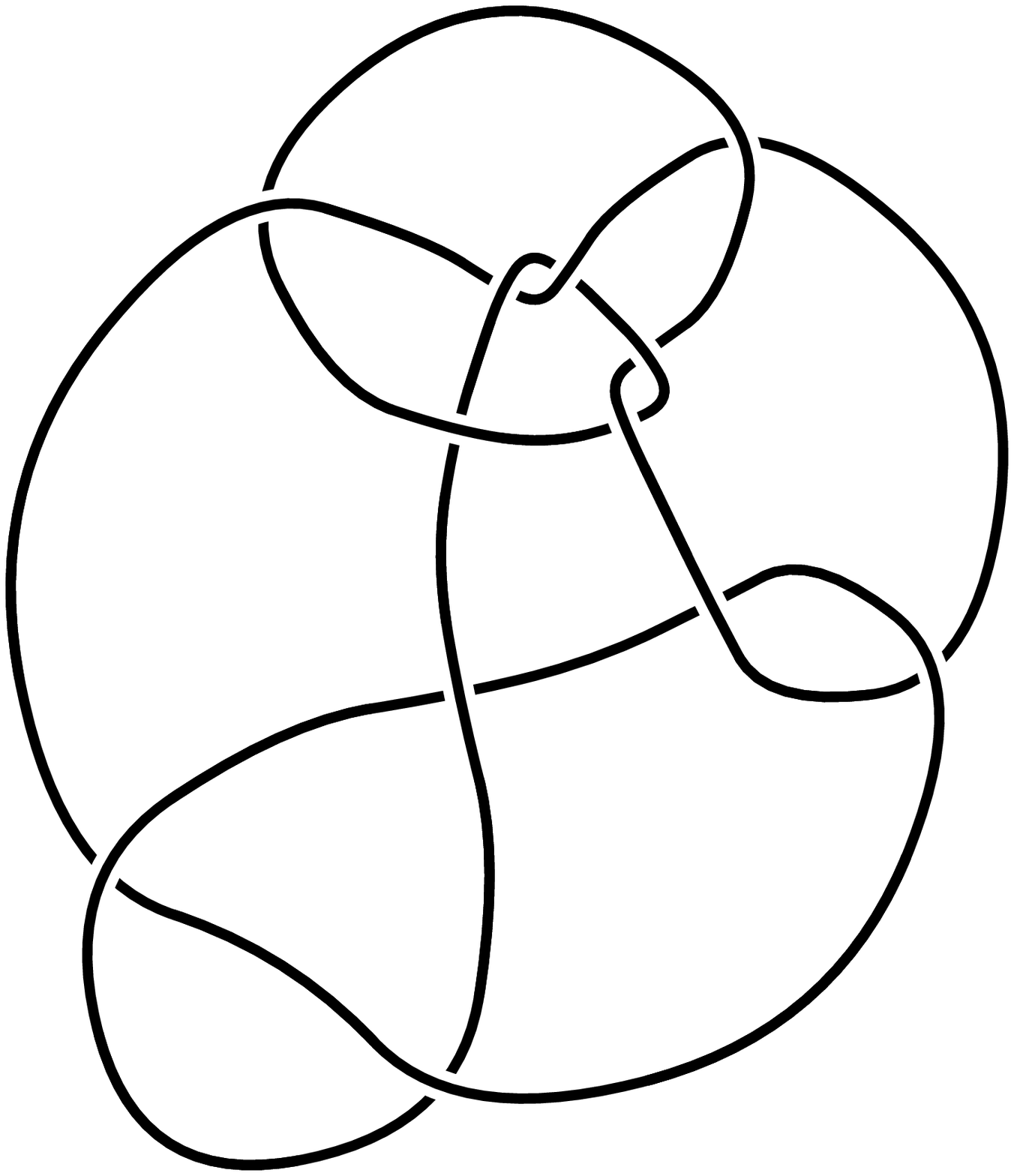}
    &
    \includegraphics[width=75pt]{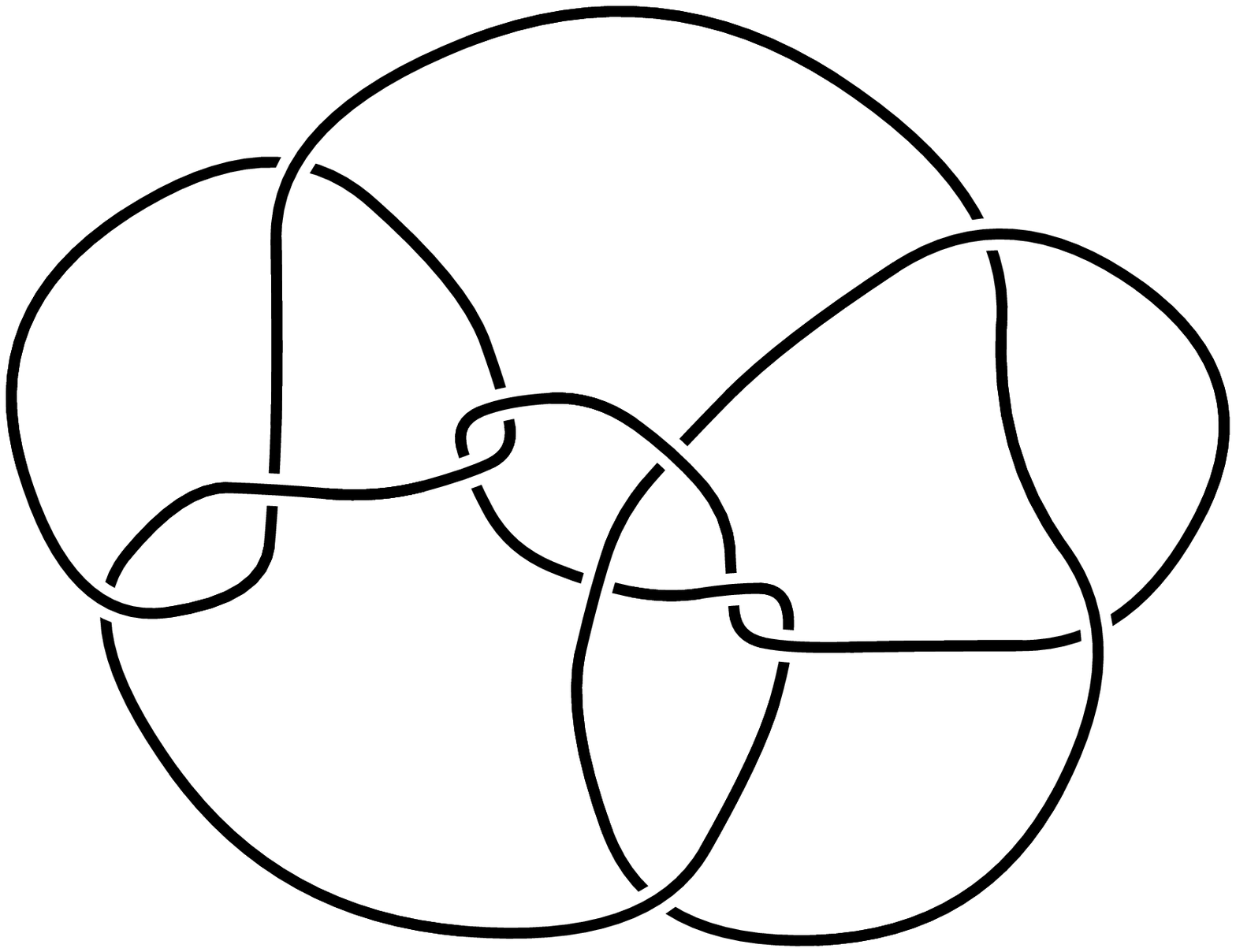}
    \\[-10pt]
    $12^N_{205}$ & $12^N_{226}$ & $12^N_{206}$ & $12^N_{227}$
    \\[10pt]
    \hline
    &&&\\[-10pt]
    \includegraphics[width=75pt]{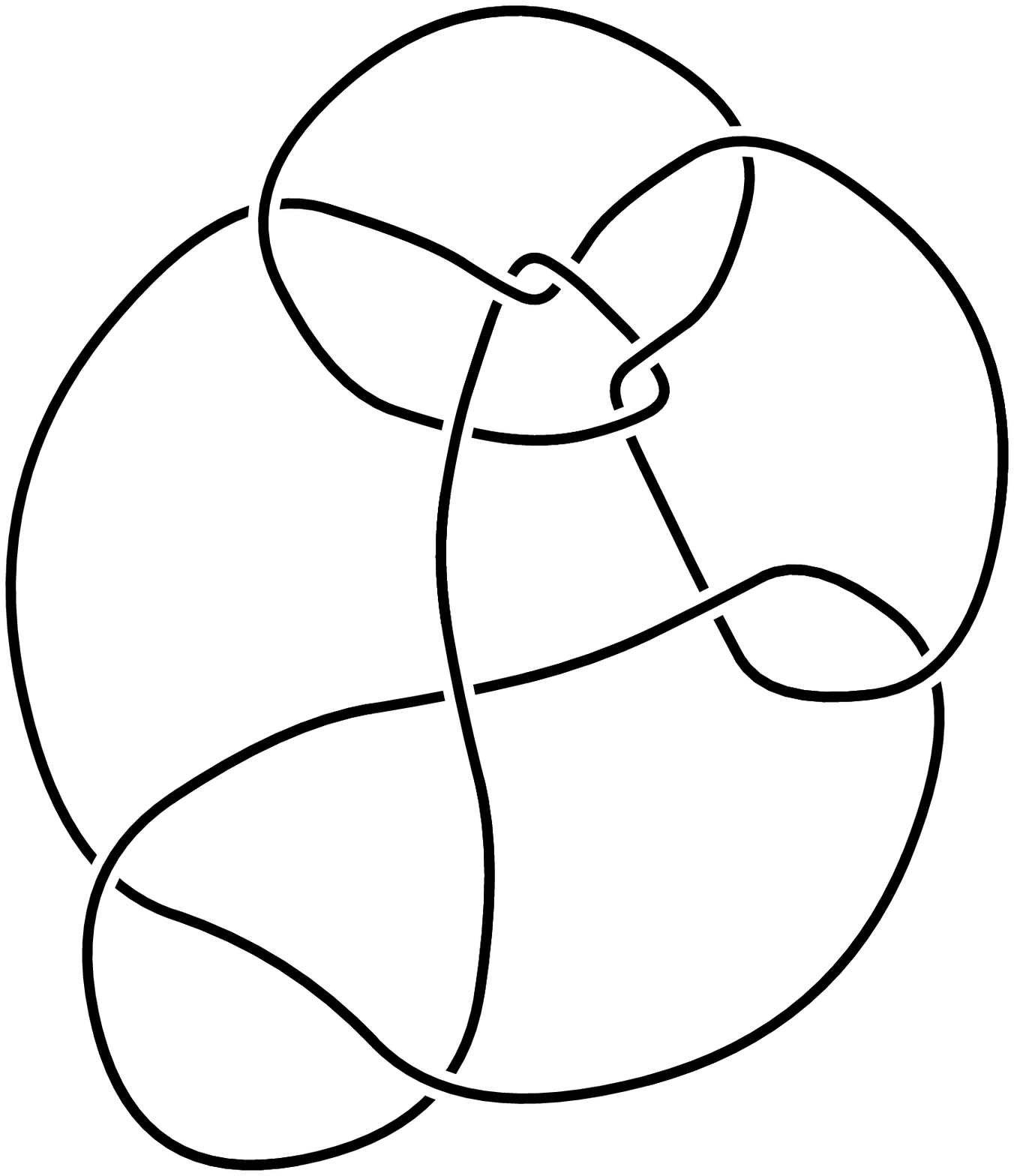}
    &
    \includegraphics[width=75pt]{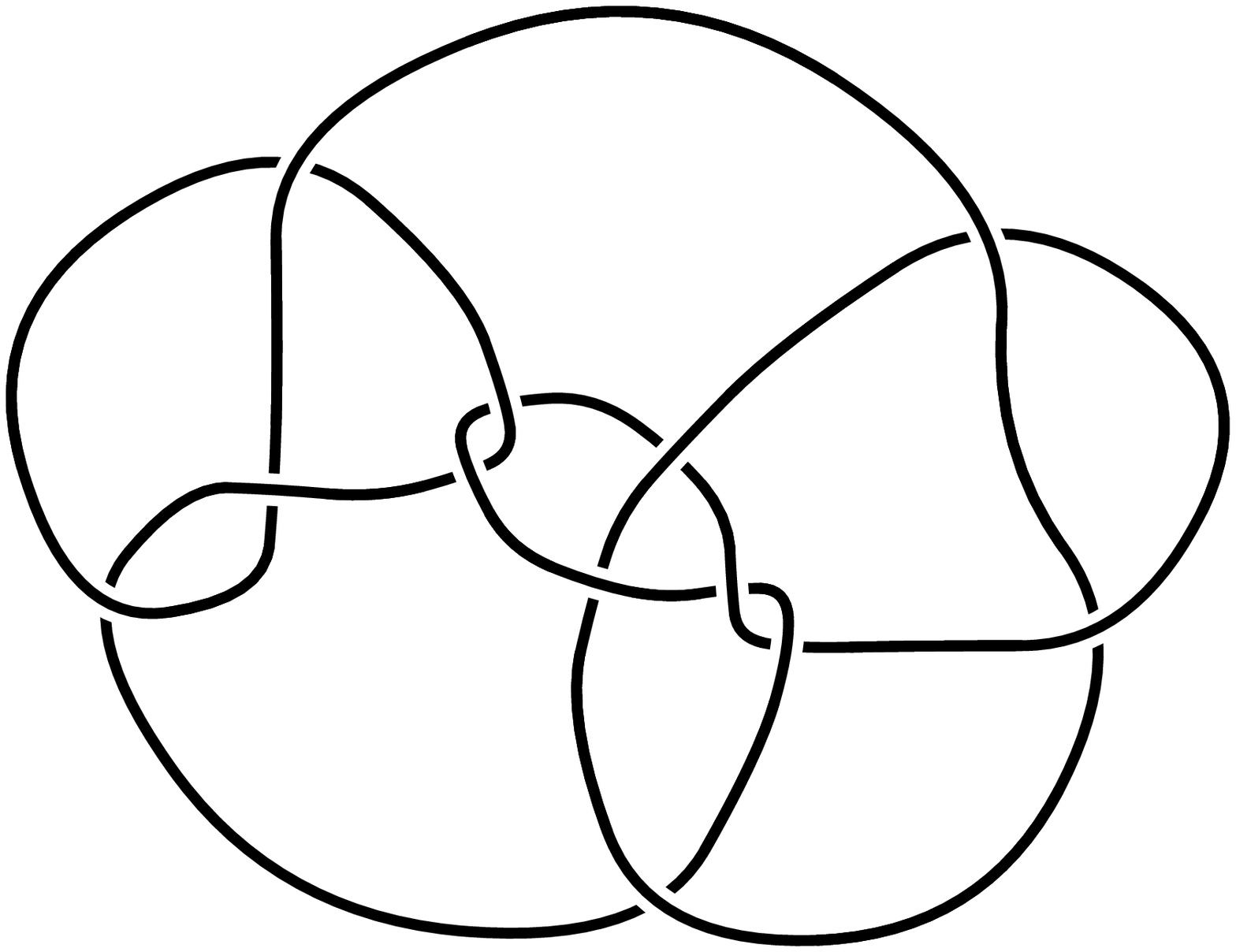}
    &
    \includegraphics[width=75pt]{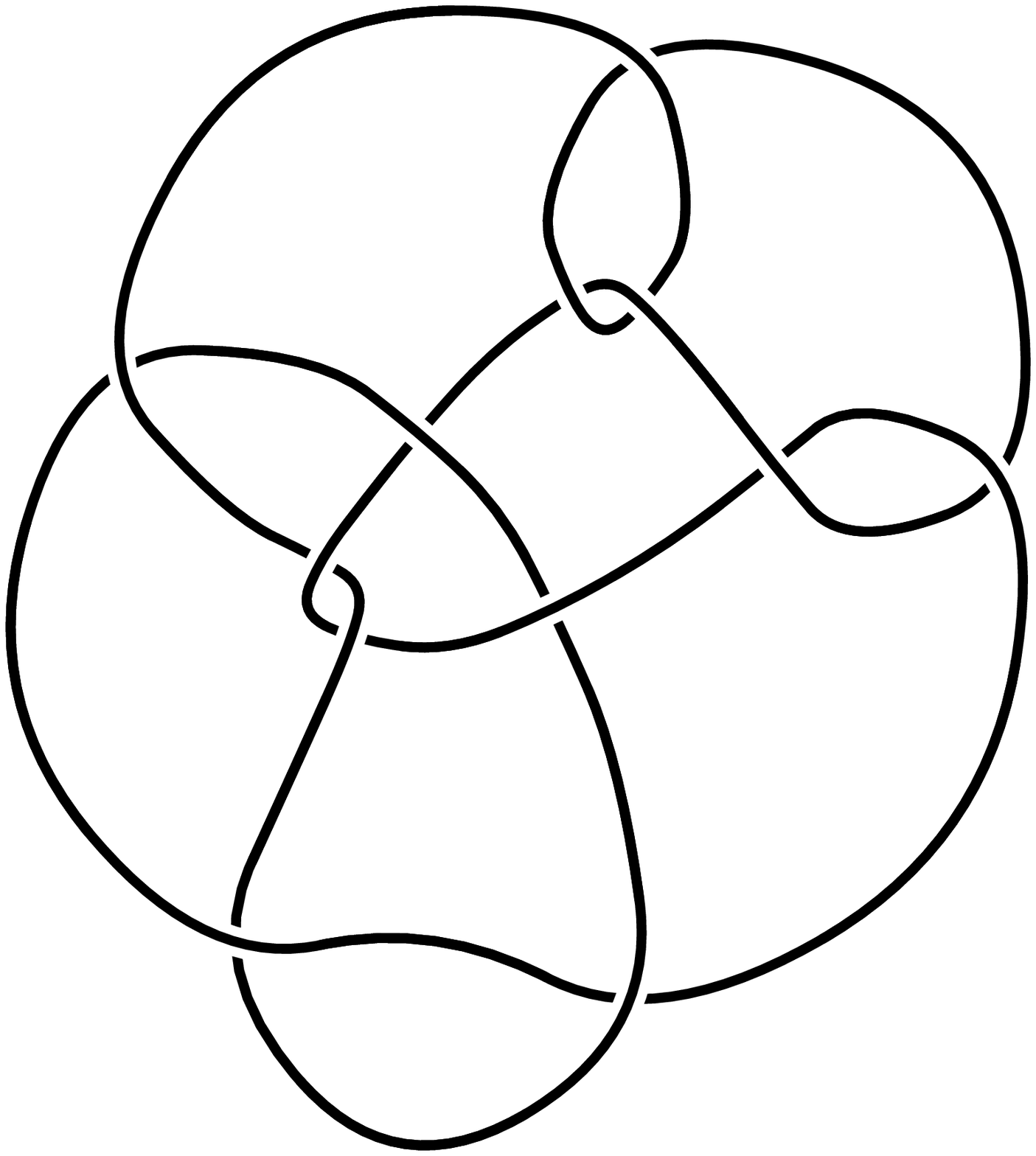}
    &
    \includegraphics[width=75pt]{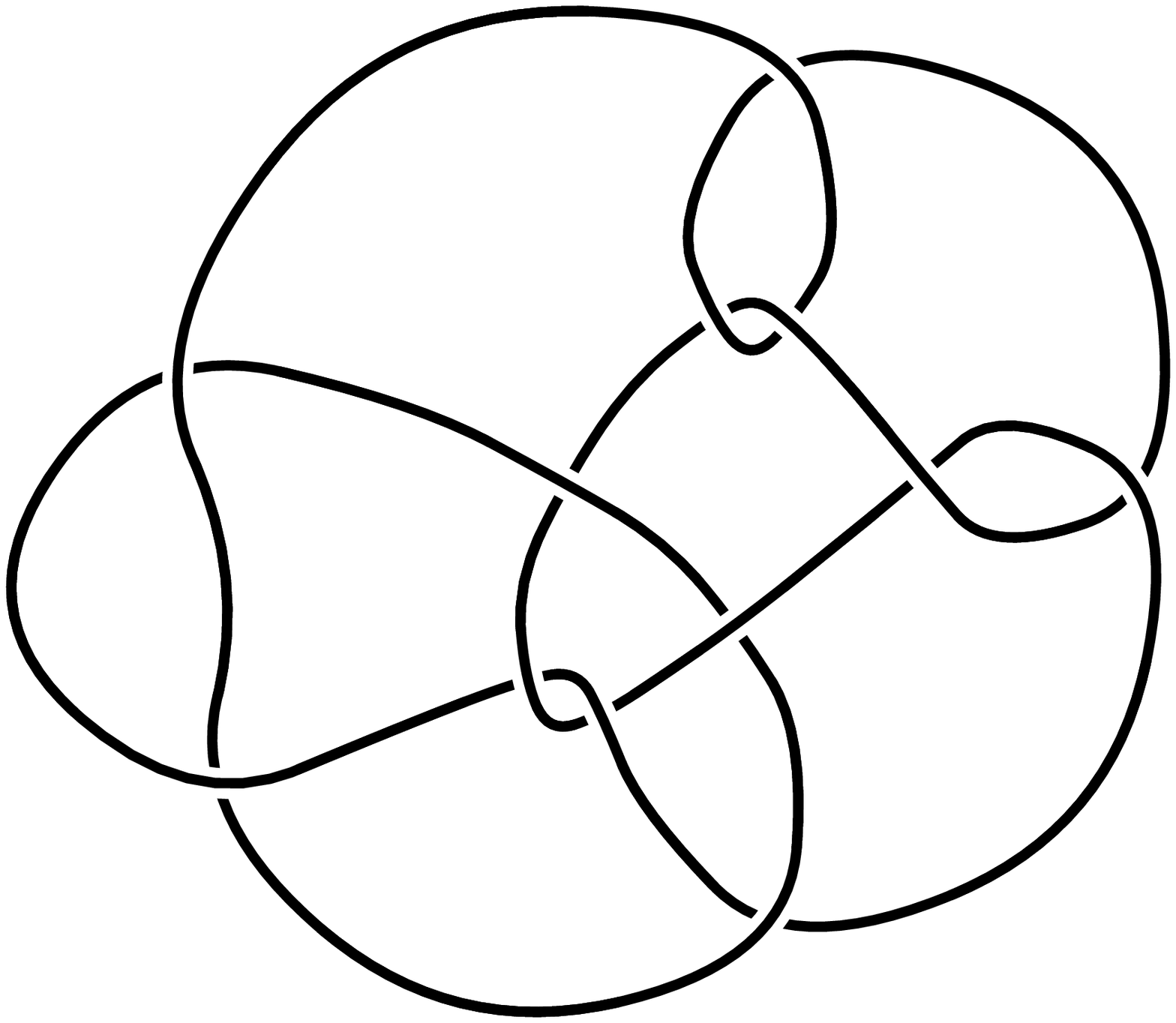}
    \\[-10pt]
    $12^N_{207}$ & $12^N_{228}$ & $12^N_{208}$ & $12^N_{212}$
  \end{tabular}
  \caption{Nonalternating $12$-crossing mutant cliques 3/6}
  \end{centering}
\end{figure}

\begin{figure}[htbp]
  \begin{centering}
  \begin{tabular}{cc@{\hspace{10pt}}|@{\hspace{10pt}}cc}
    \includegraphics[width=75pt]{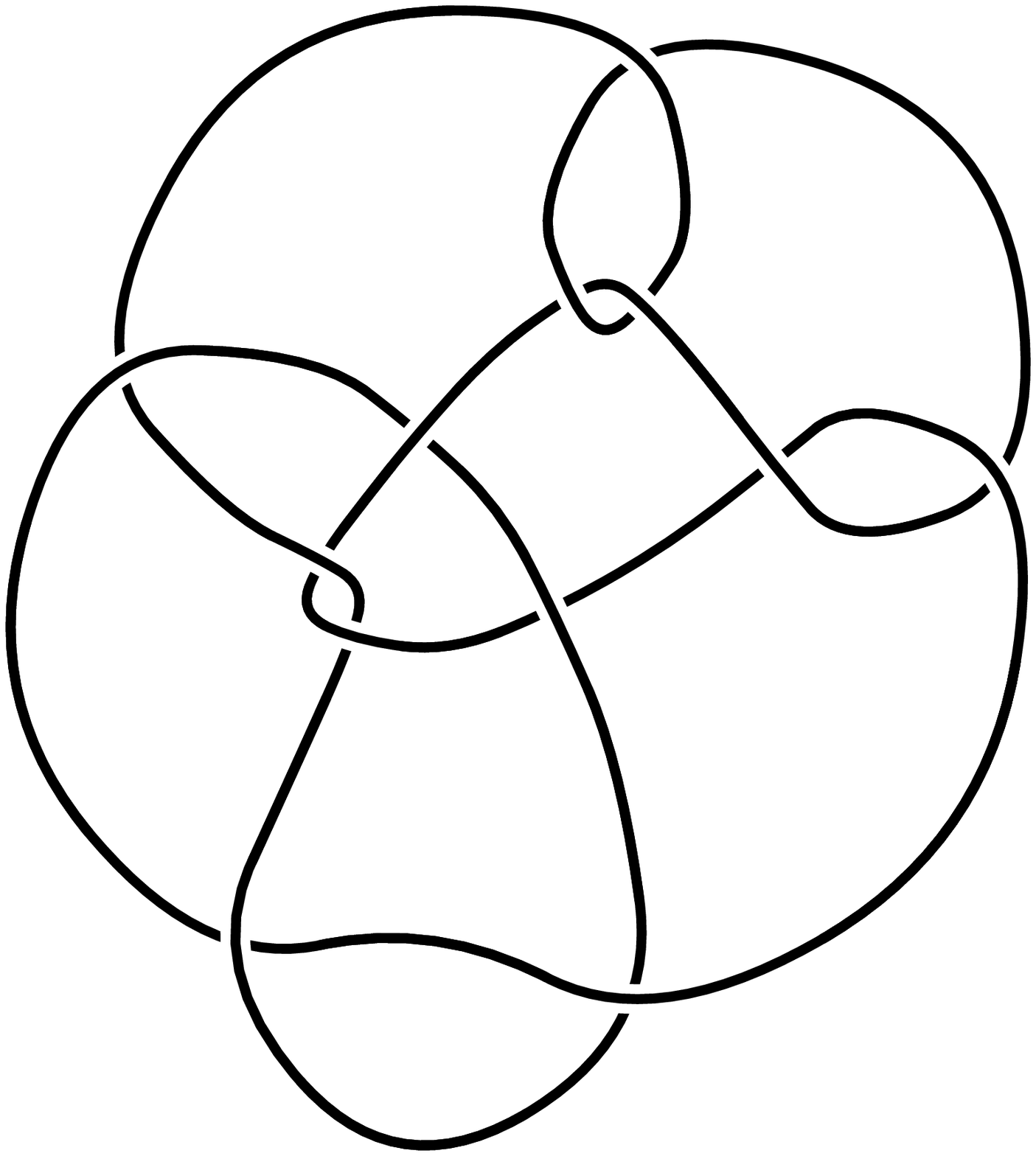}
    &
    \includegraphics[width=75pt]{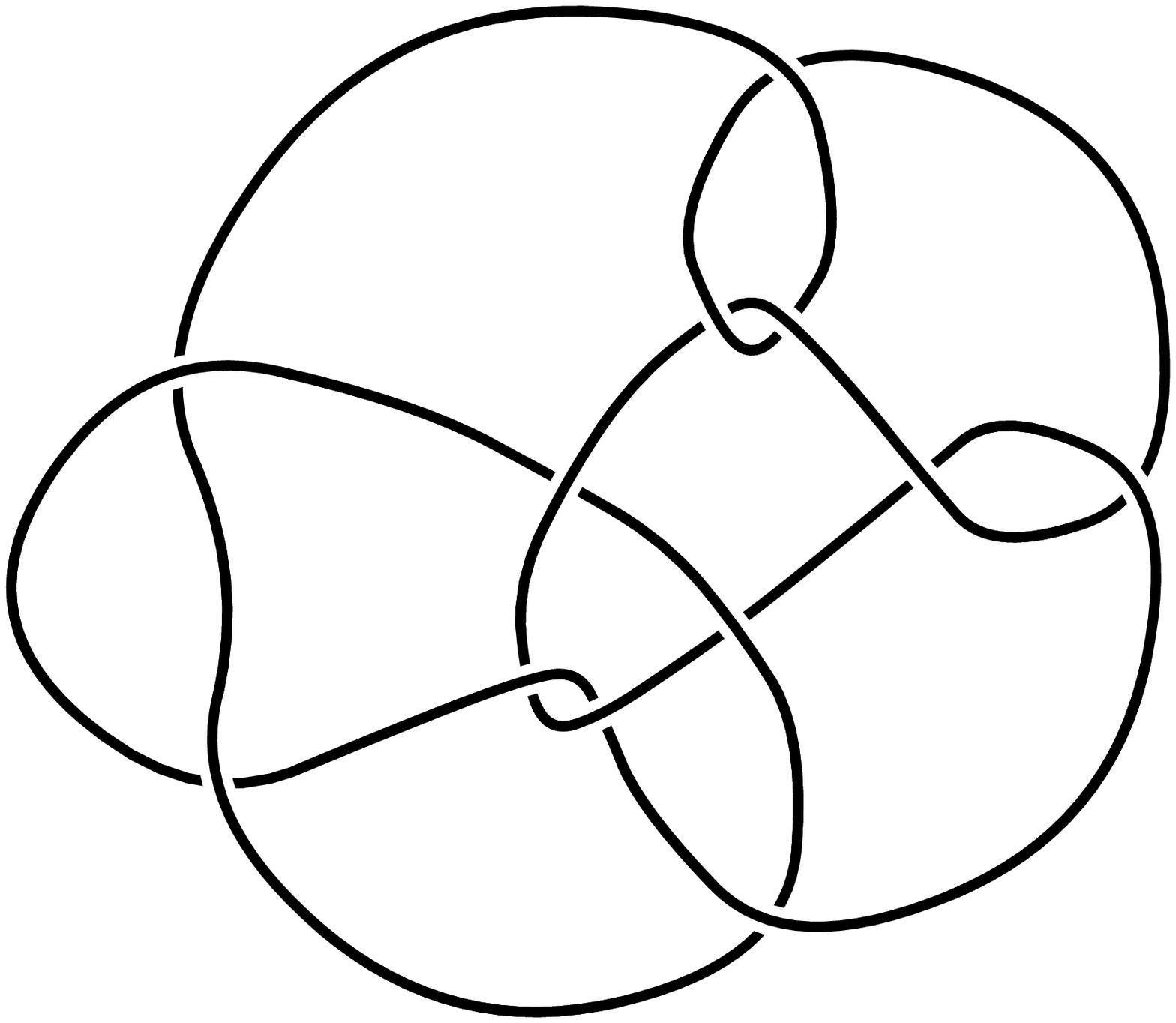}
    &
    \includegraphics[width=75pt]{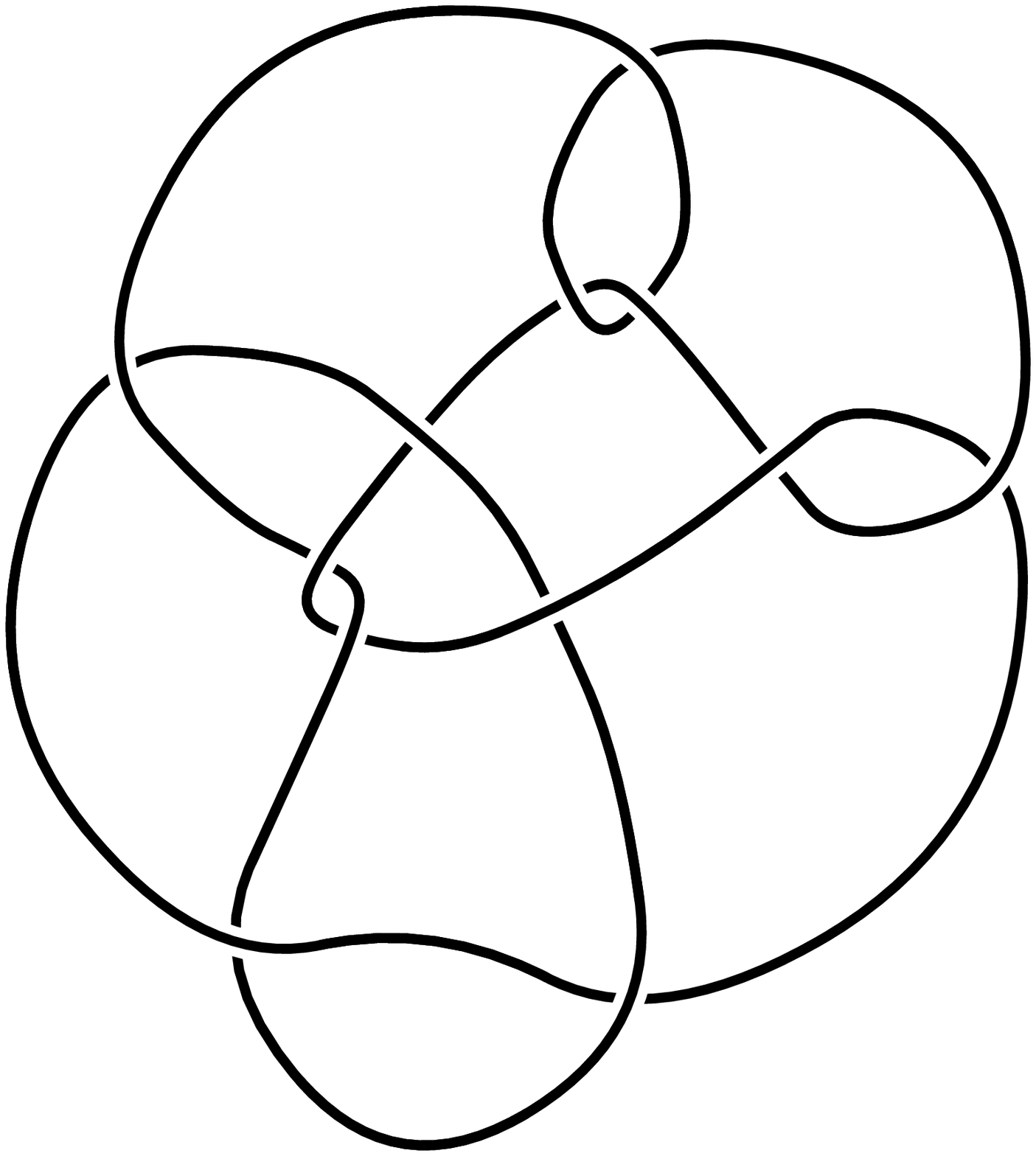}
    &
    \includegraphics[width=75pt]{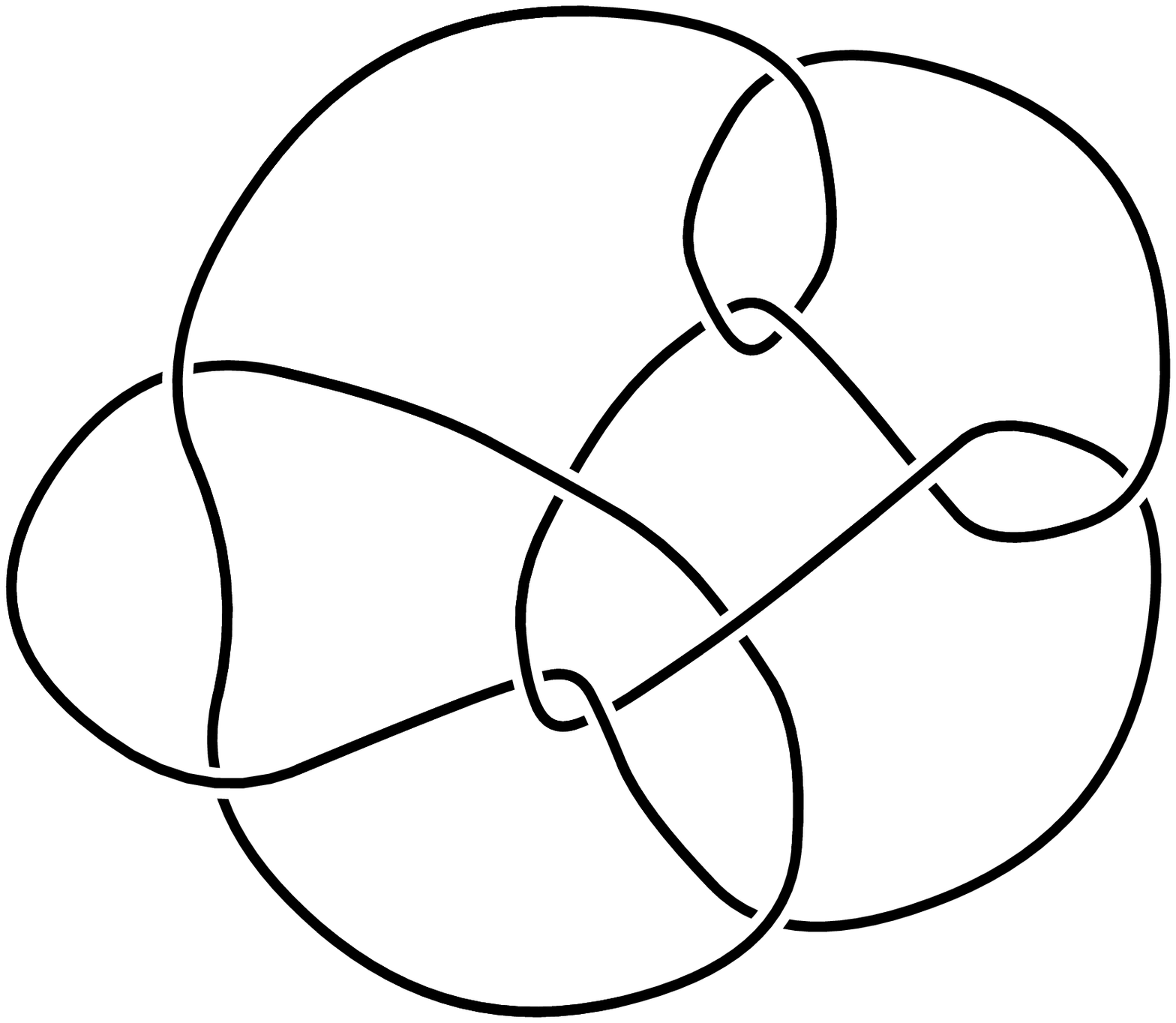}
    \\[-10pt]
    $12^N_{209}$ & $12^N_{213}$ & $12^N_{210}$ & $12^N_{214}$
    \\[10pt]
    \hline
    &&&\\[-10pt]
    \includegraphics[width=75pt]{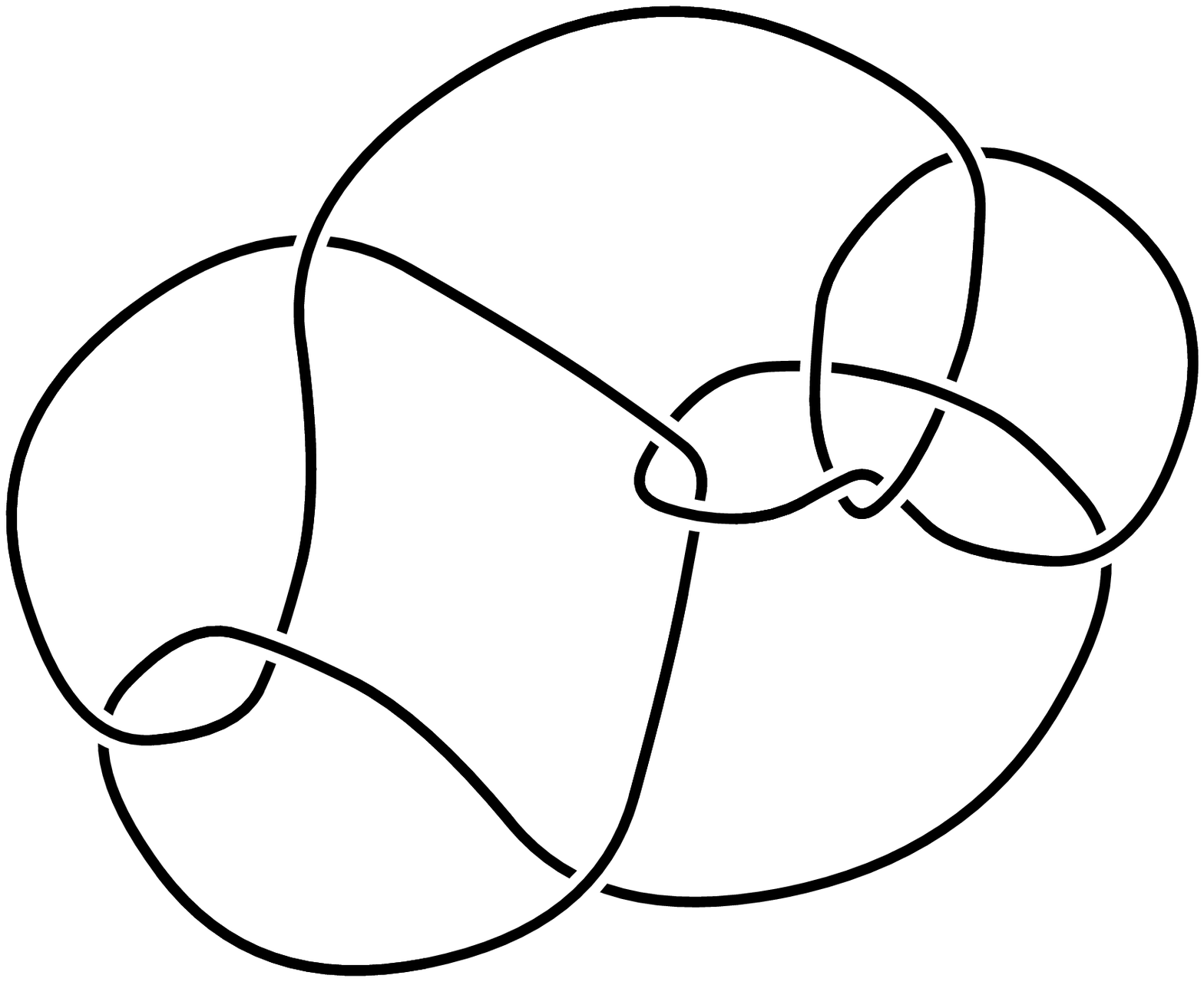}
    &
    \includegraphics[width=75pt]{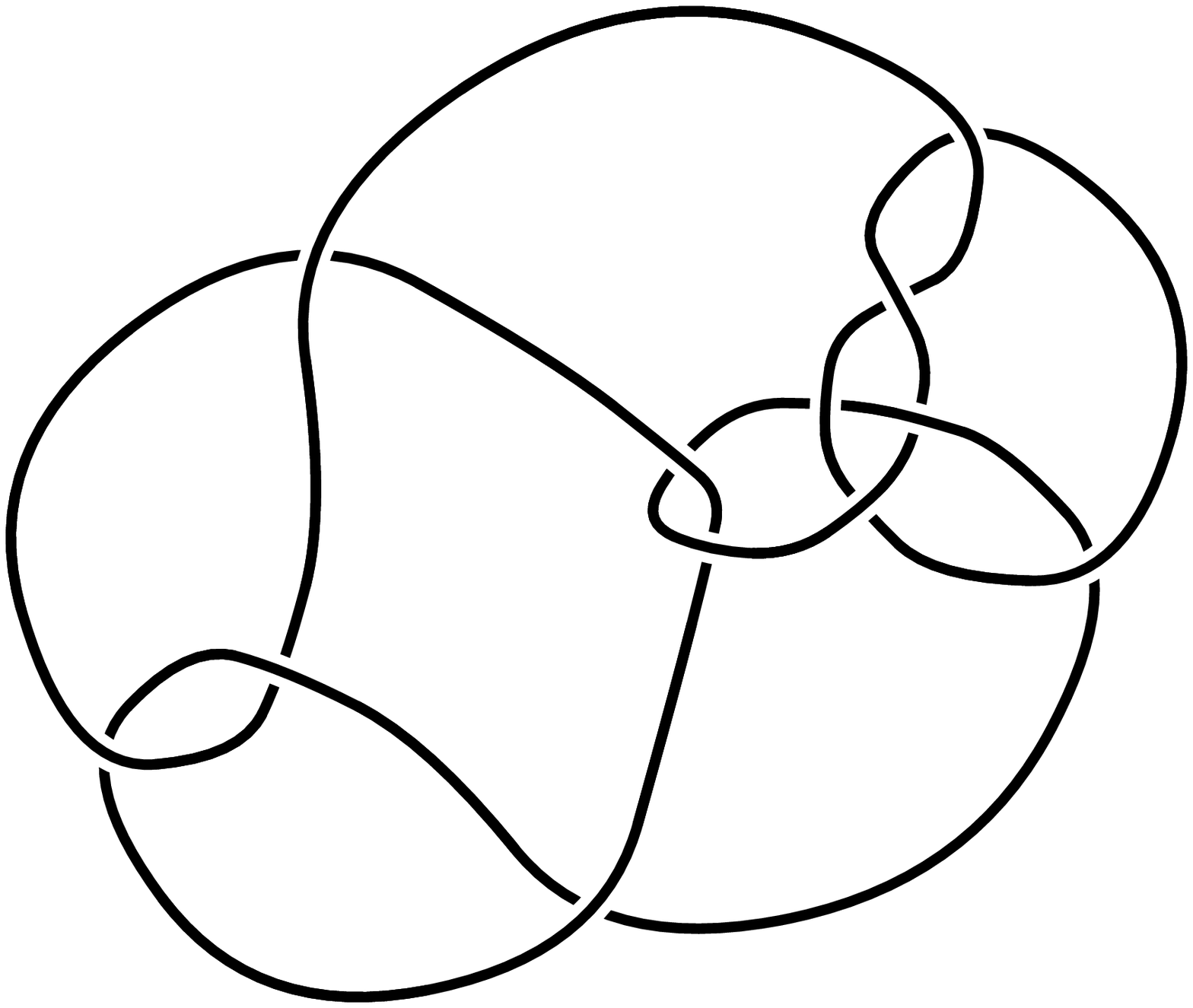}
    &
    \includegraphics[width=75pt]{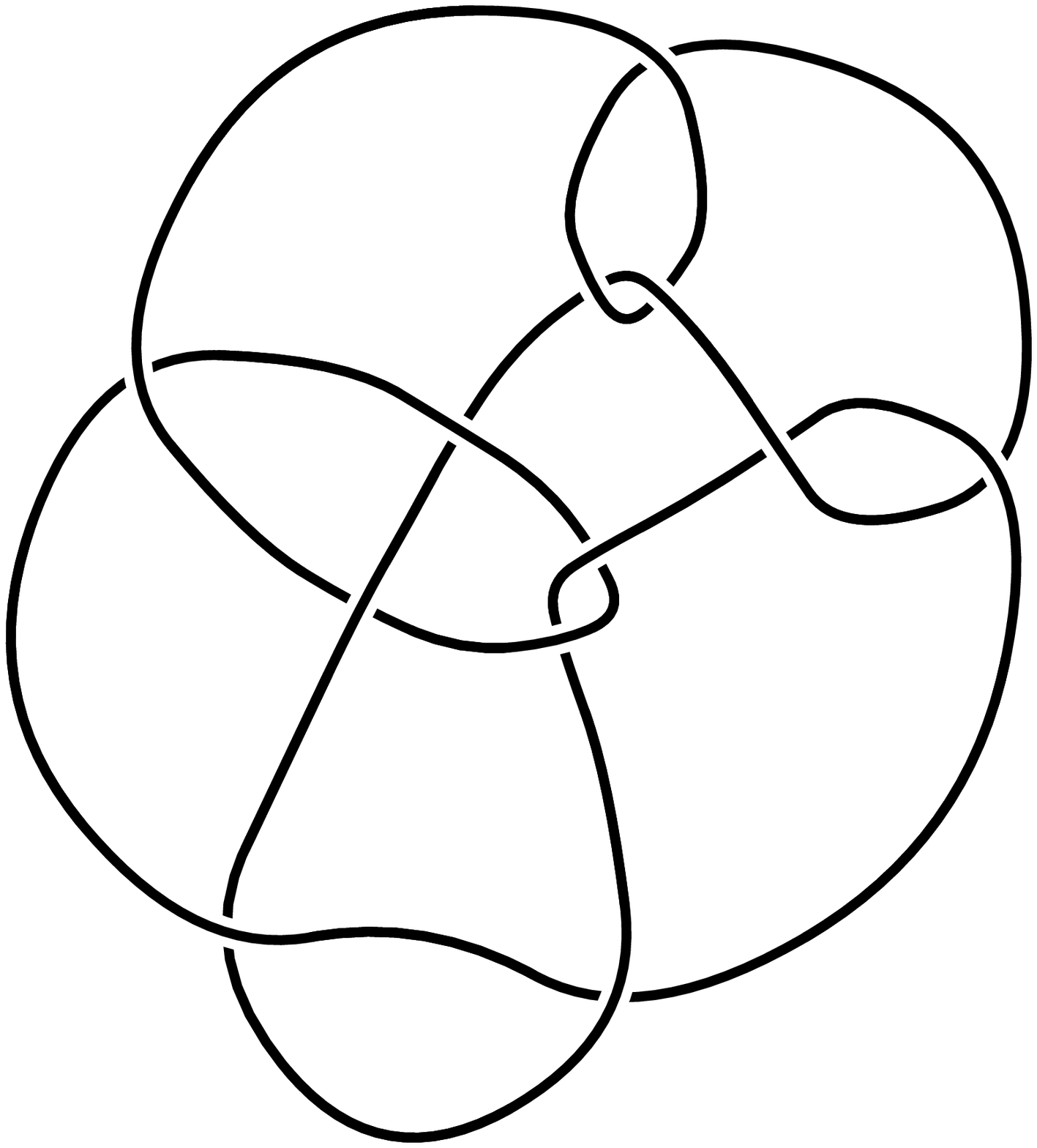}
    &
    \includegraphics[width=75pt]{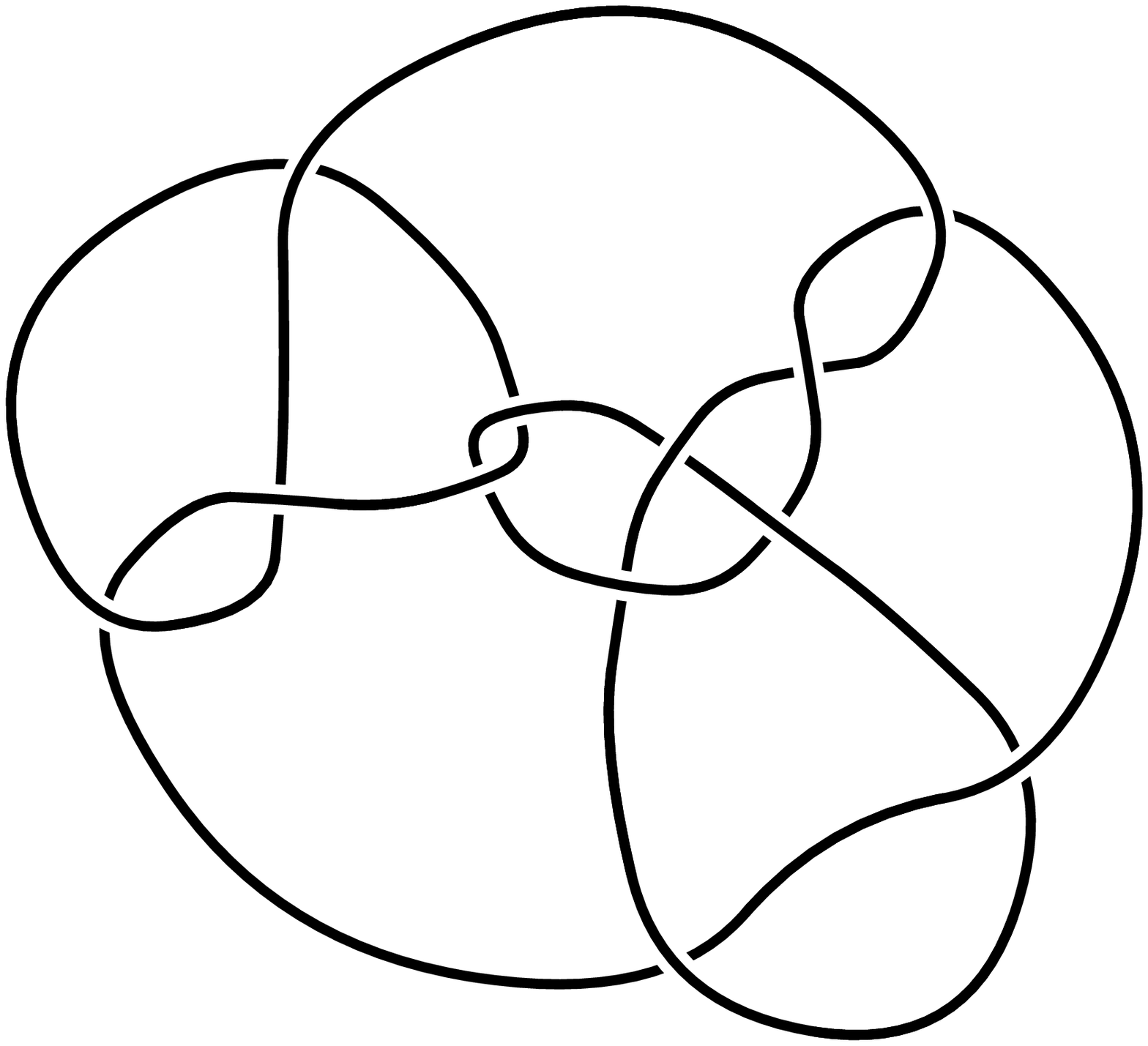}
    \\[-10pt]
    $12^N_{231}$ & $12^N_{232}$ & $12^N_{252}$ & $12^N_{262}$
    \\[10pt]
    \hline
    &&&\\[-10pt]
    \includegraphics[width=75pt]{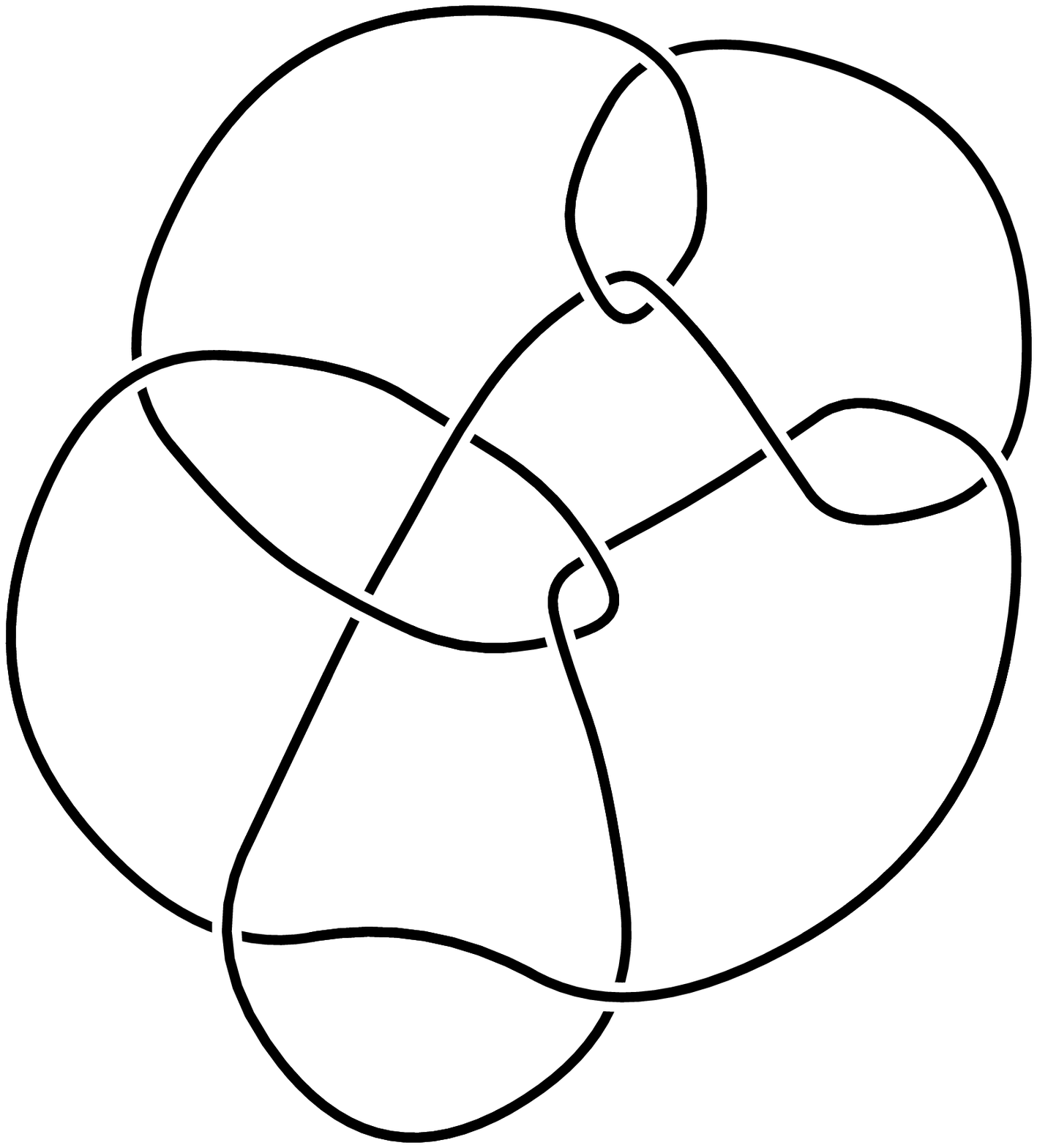}
    &
    \includegraphics[width=75pt]{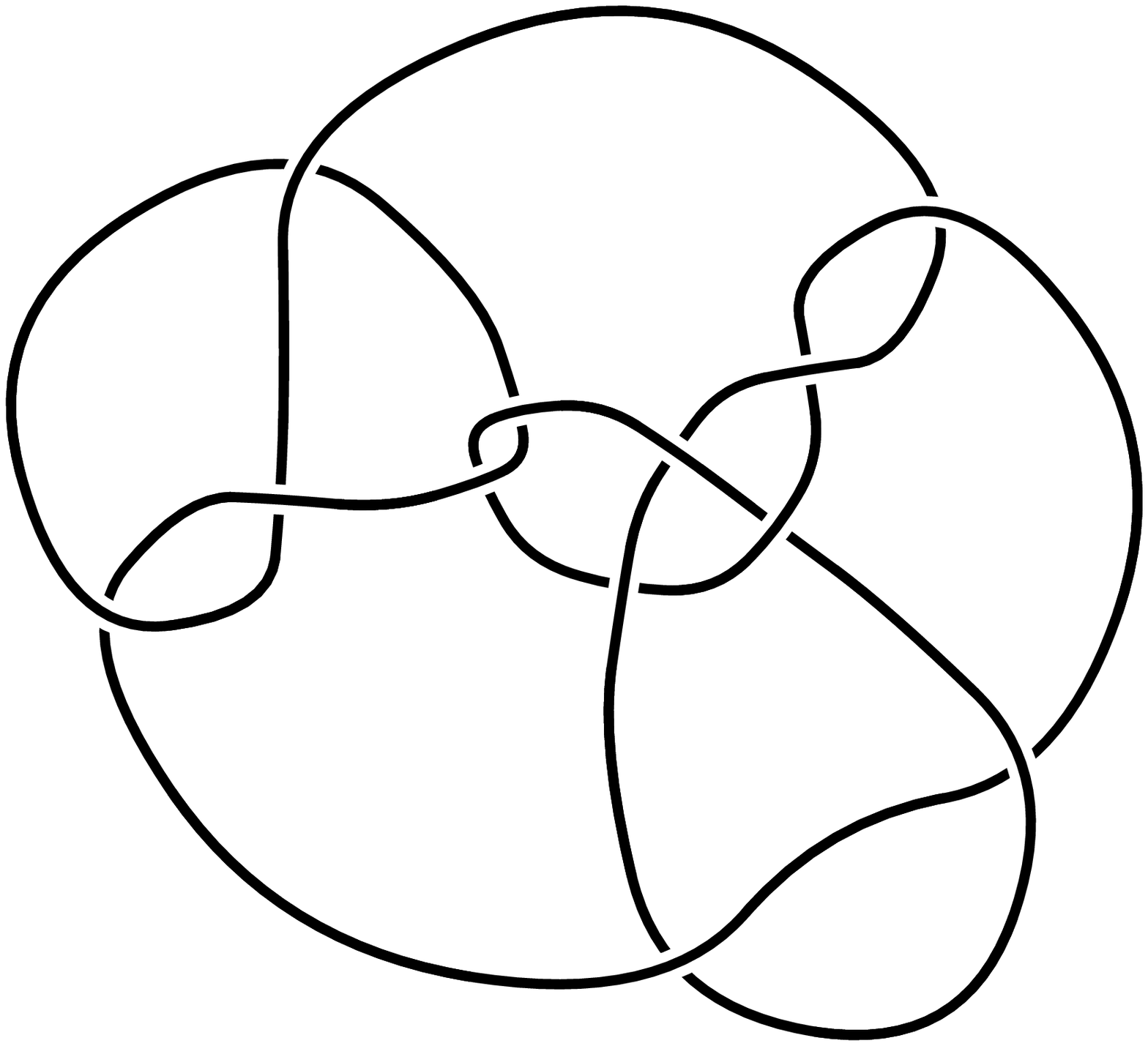}
    &
    \includegraphics[width=75pt]{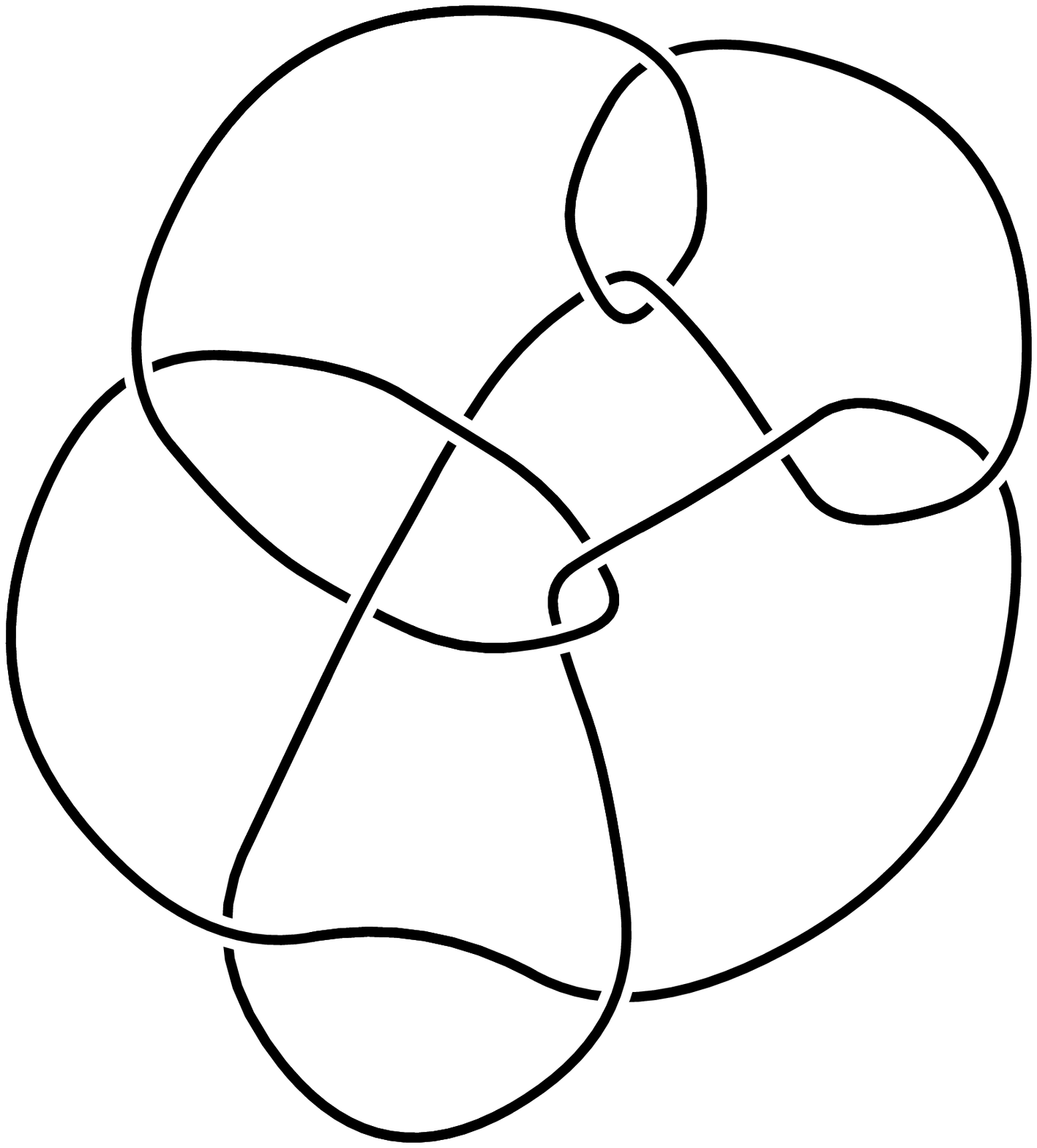}
    &
    \includegraphics[width=75pt]{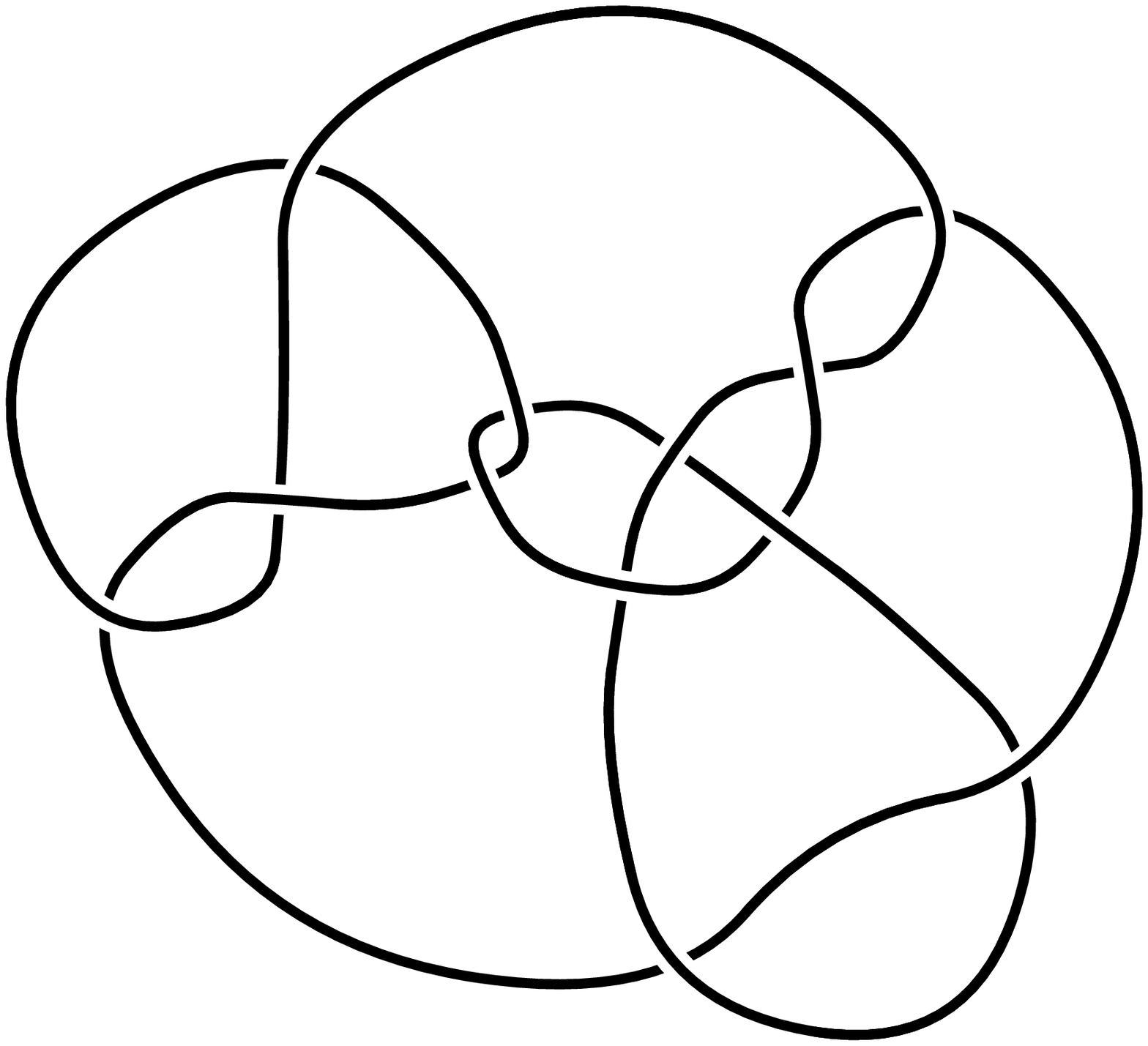}
    \\[-10pt]
    $12^N_{255}$ & $12^N_{263}$ & $12^N_{256}$ & $12^N_{264}$
    \\[10pt]
    \hline
    &&&\\[-10pt]
    \includegraphics[width=75pt]{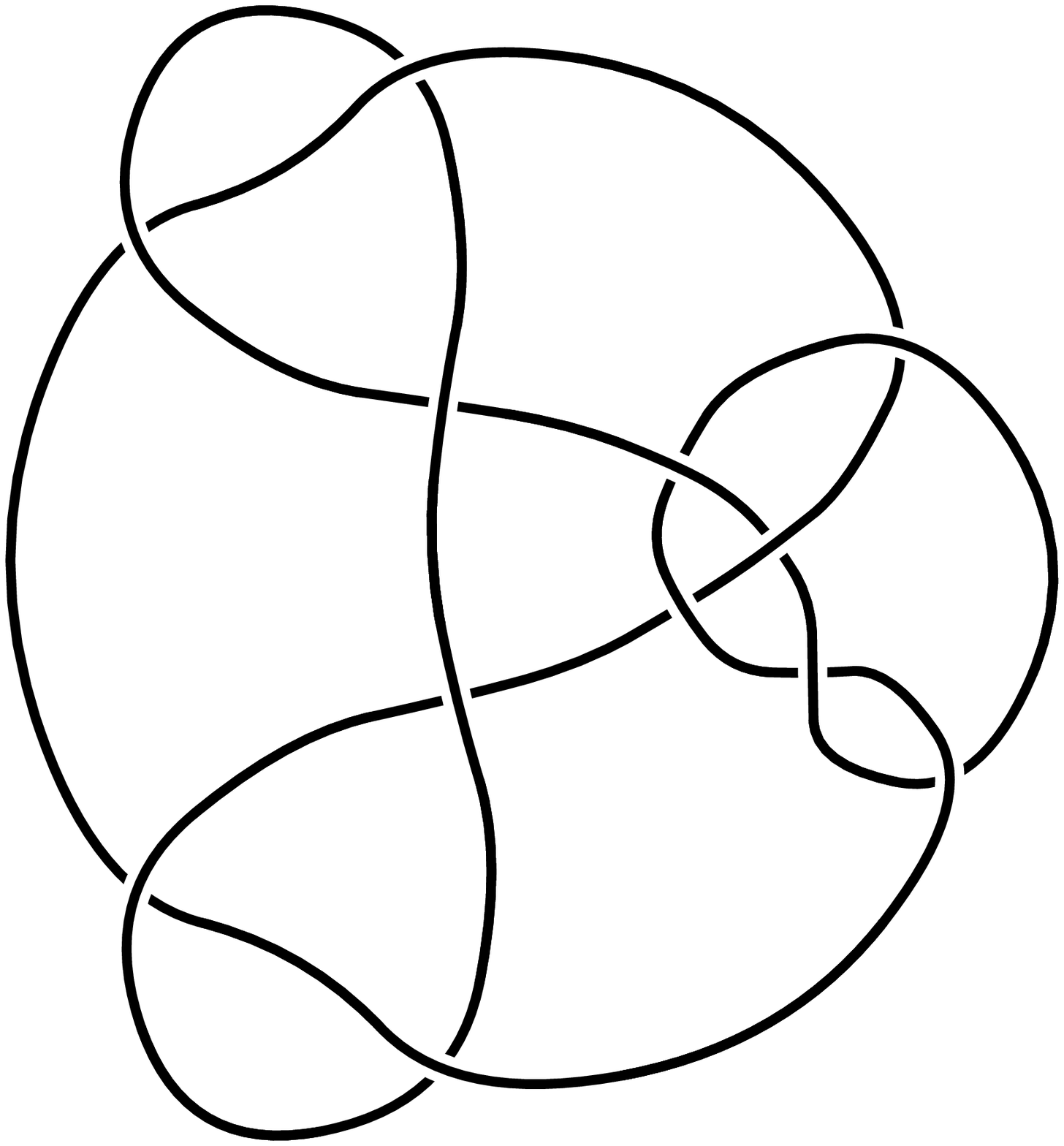}
    &
    \includegraphics[width=75pt]{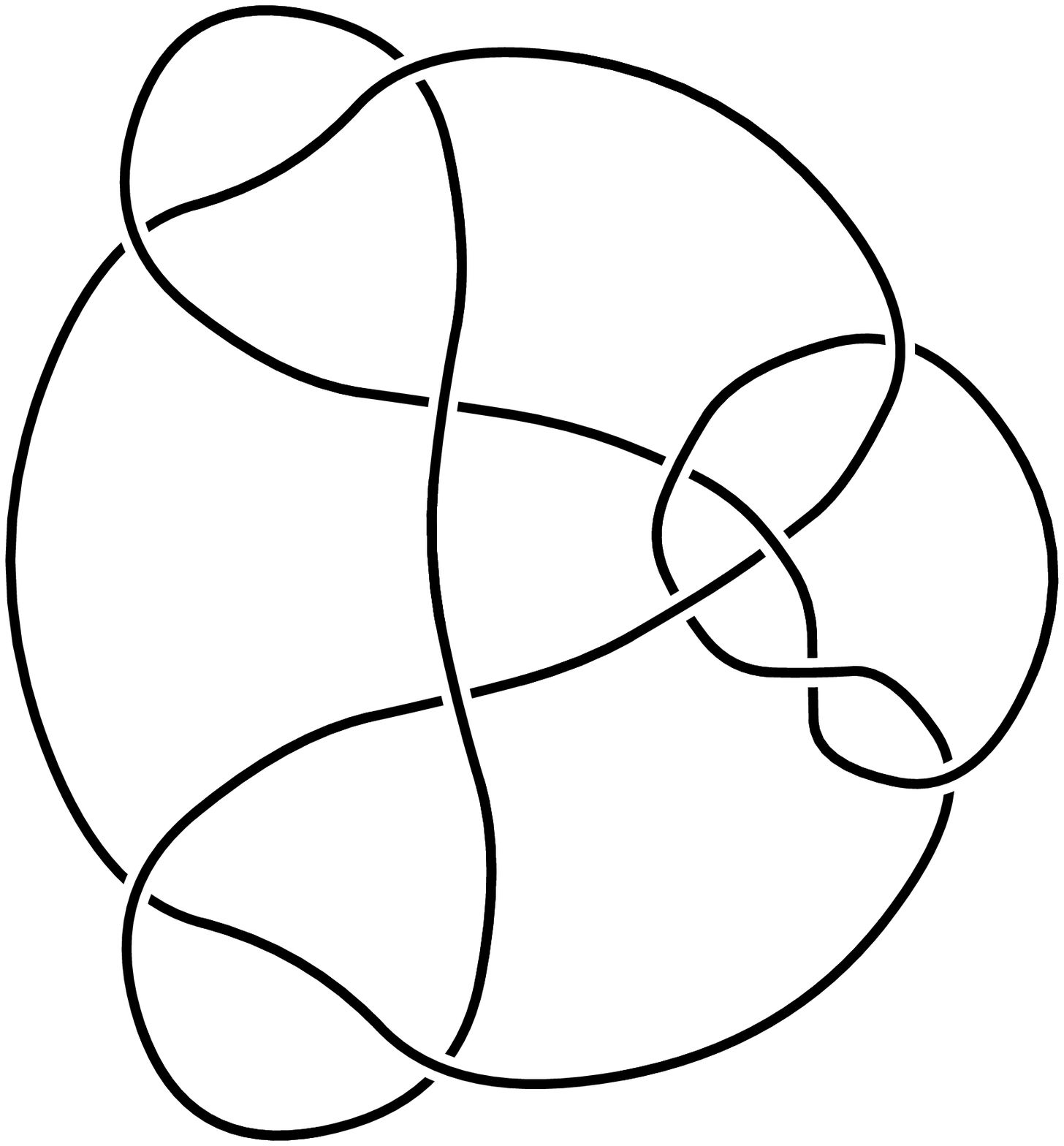}
    &
    \includegraphics[width=75pt]{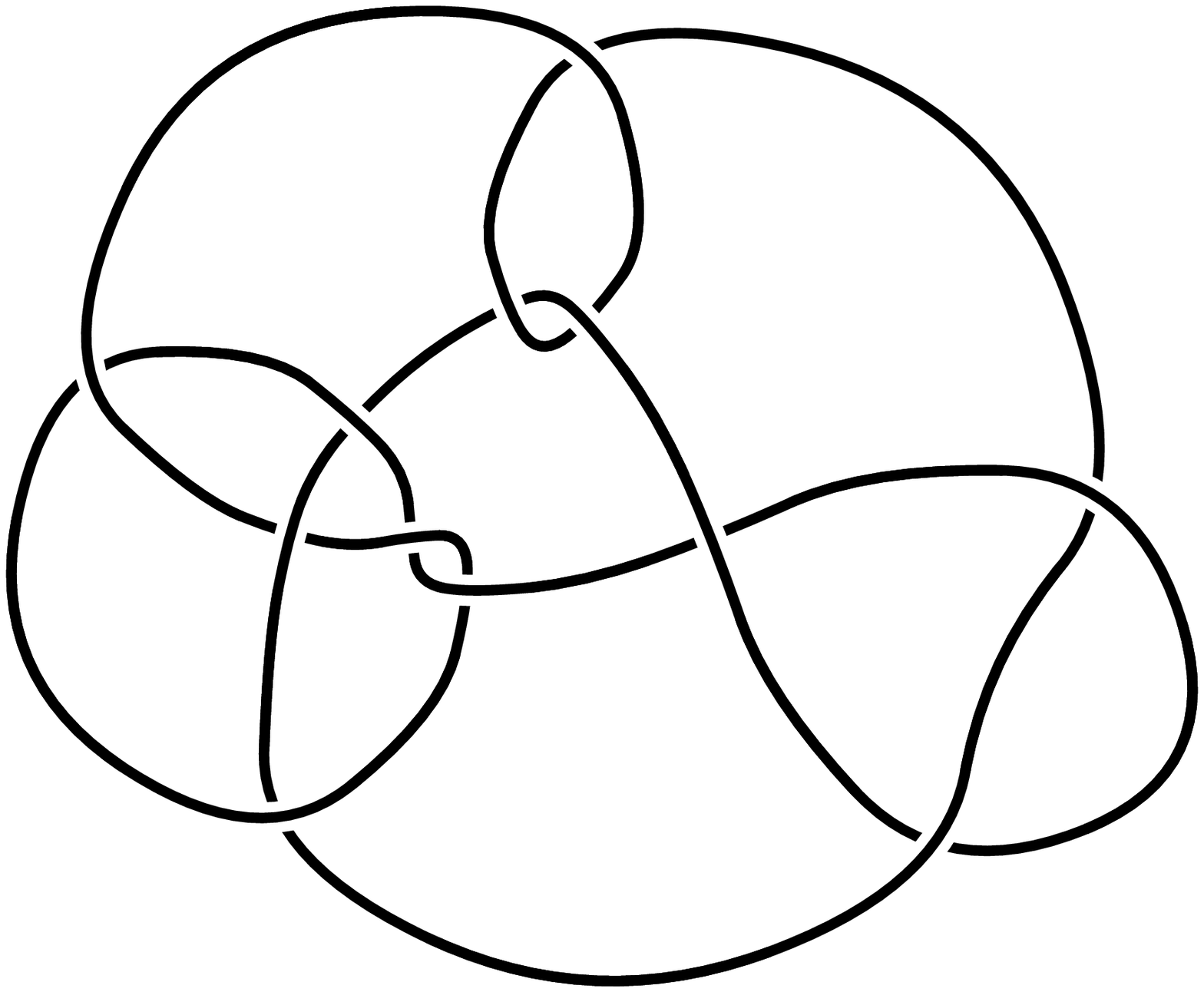}
    &
    \includegraphics[width=75pt]{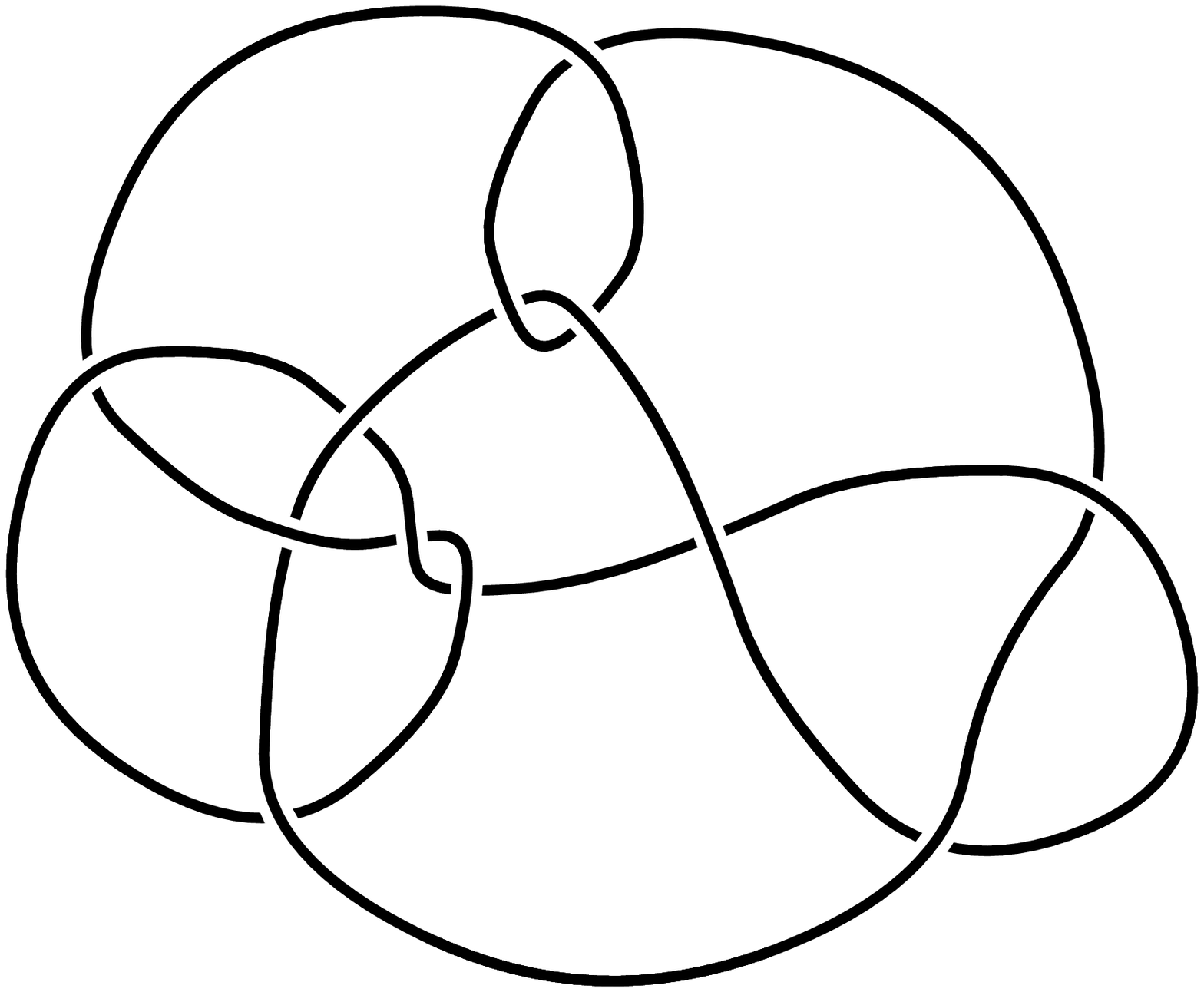}
    \\[-10pt]
    $12^N_{364}$ & $12^N_{365}$ & $12^N_{421}$ & $12^N_{422}$
    \\[10pt]
    \hline
    &&&\\[-10pt]
    \includegraphics[width=75pt]{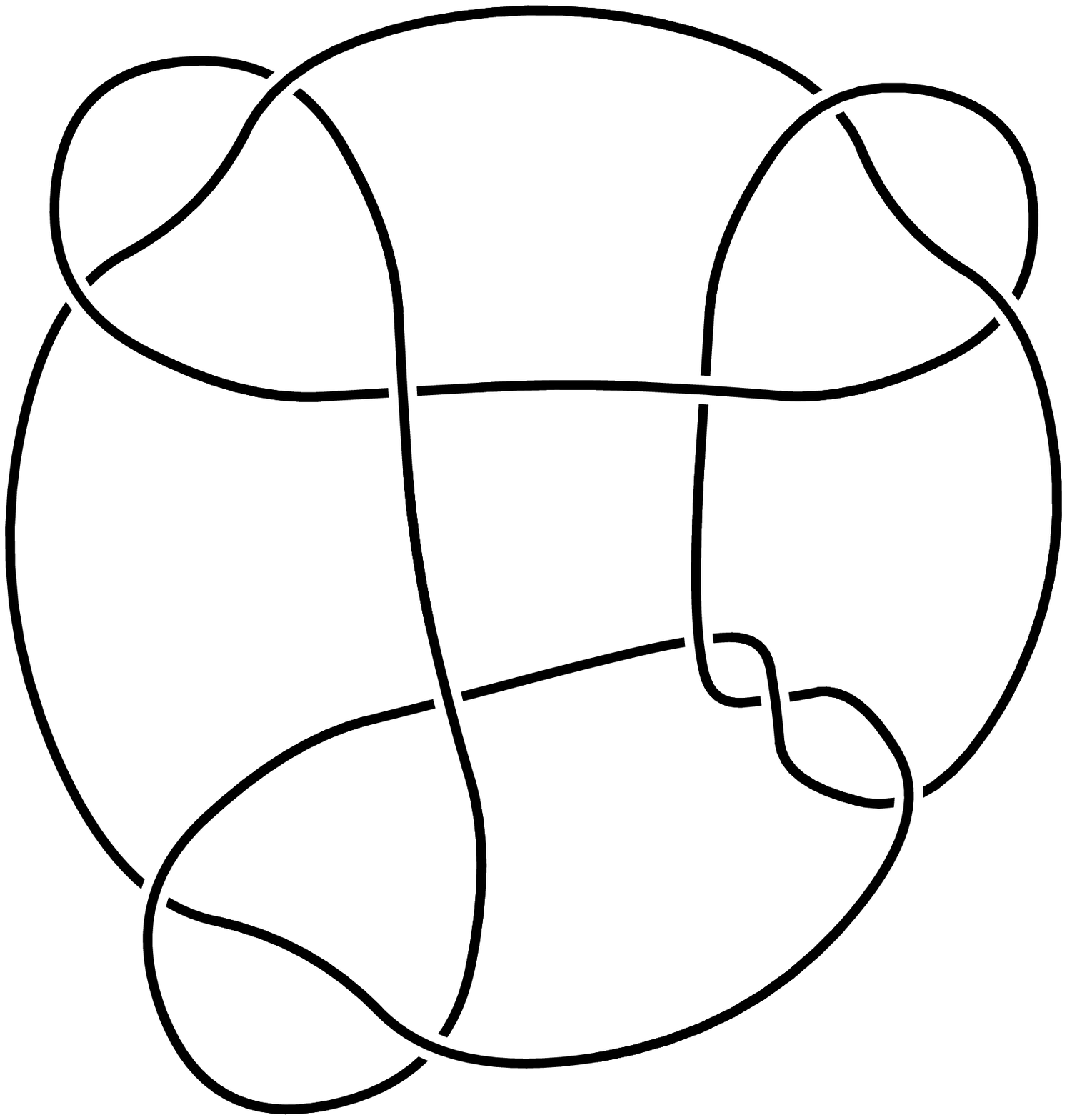}
    &
    \includegraphics[width=75pt]{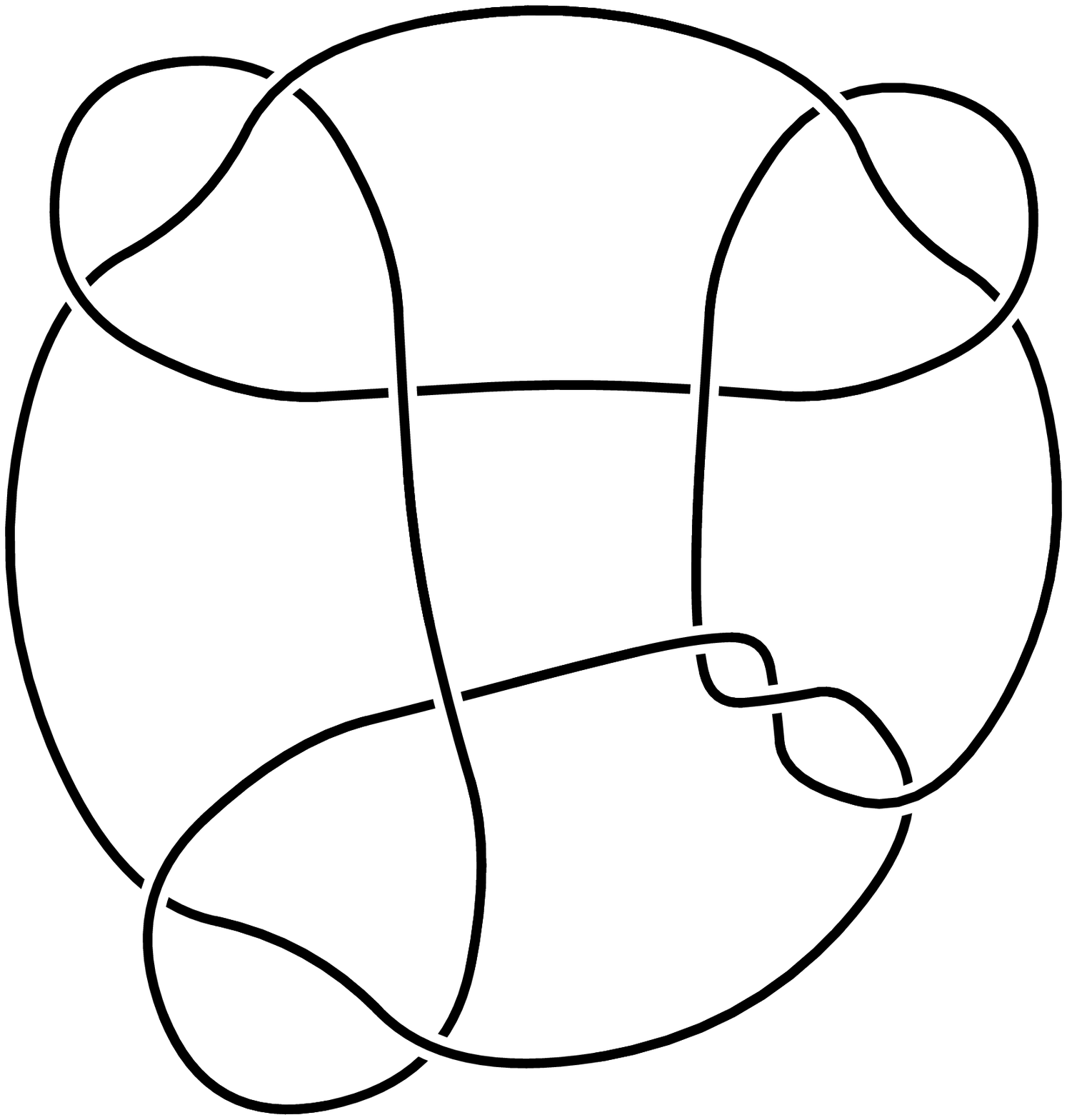}
    &
    \includegraphics[width=75pt]{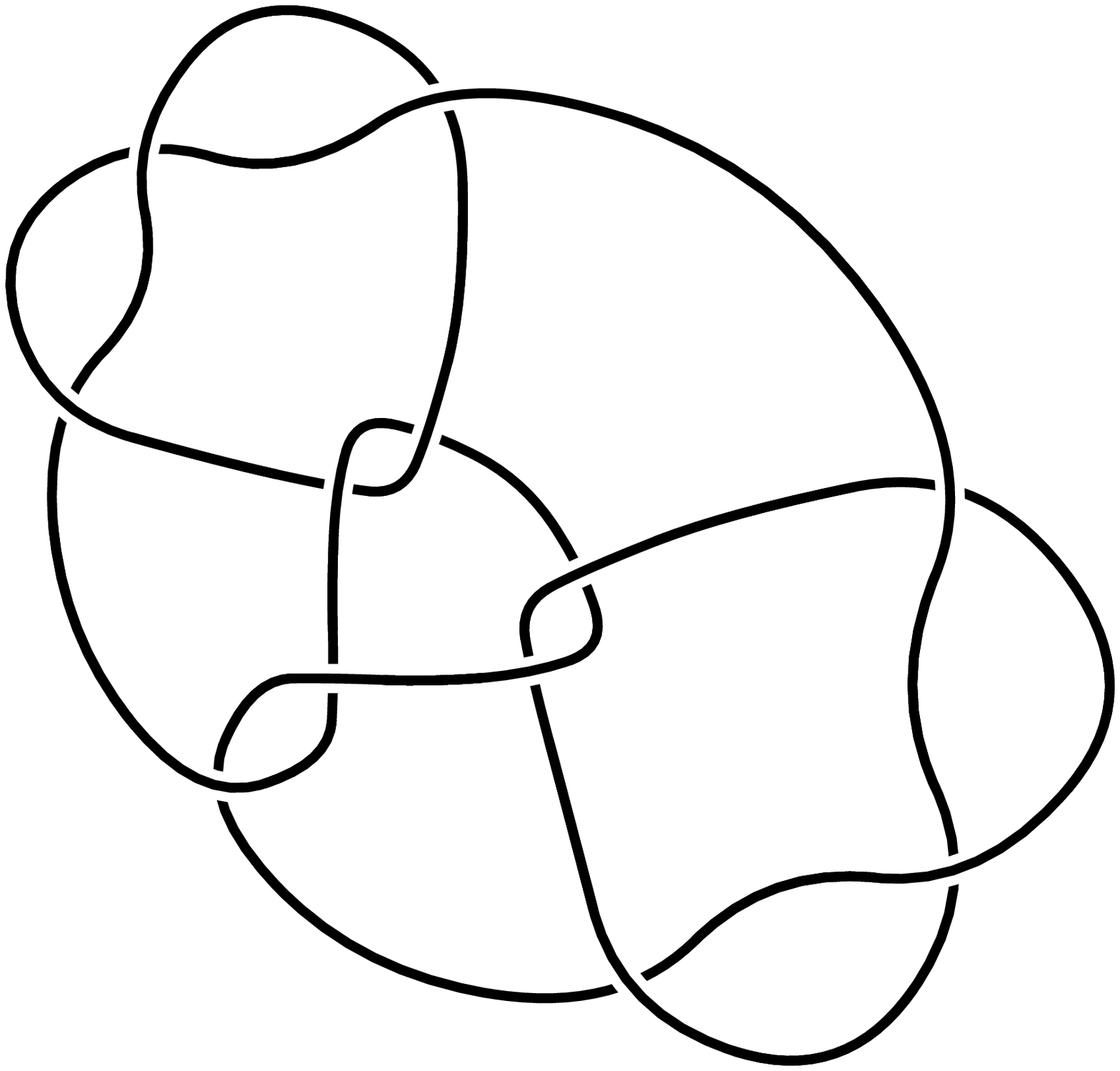}
    &
    \includegraphics[width=75pt]{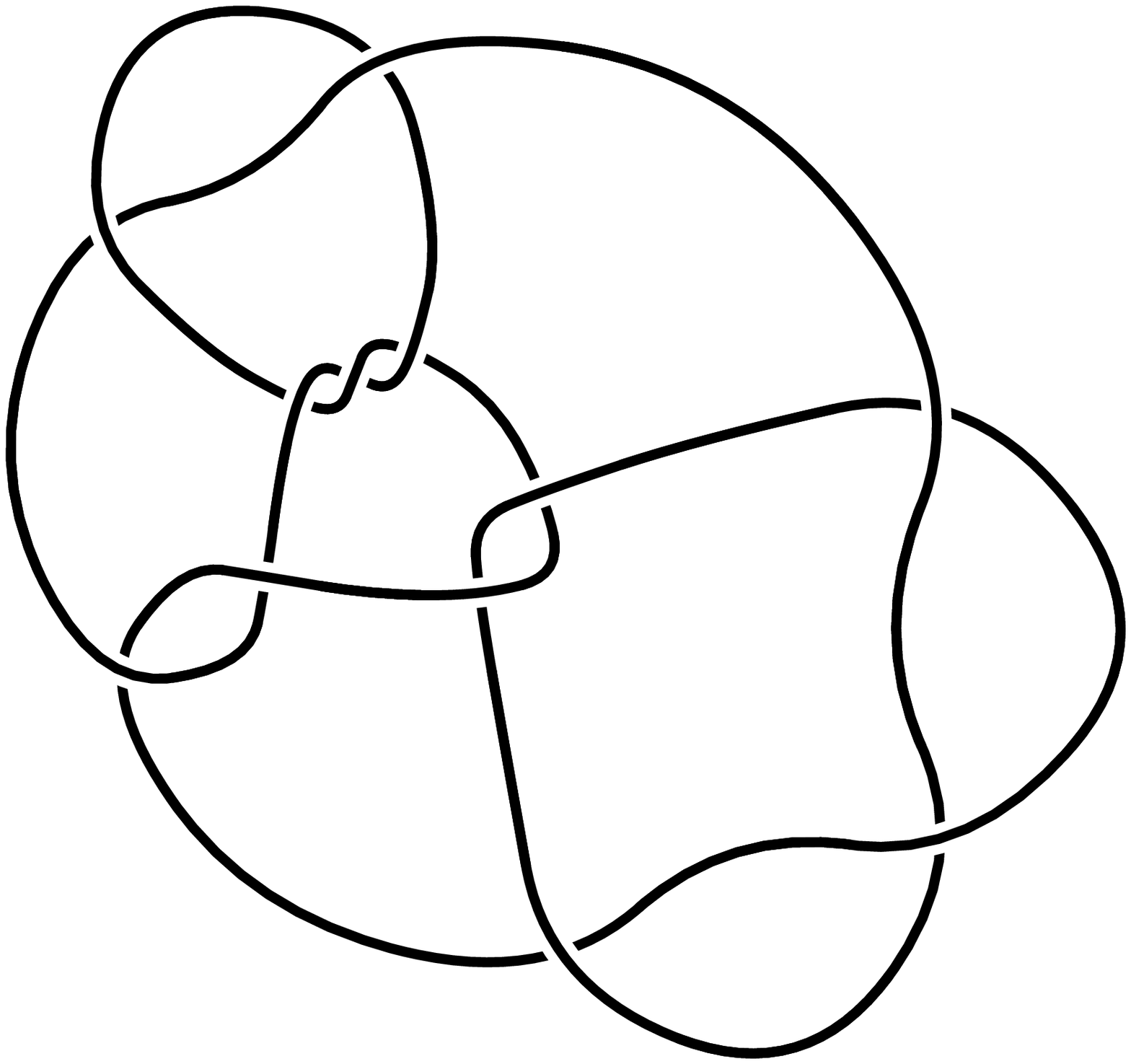}
    \\[-10pt]
    $12^N_{553}$ & $12^N_{556}$ & $12^N_{670}$ & $12^N_{681}$
  \end{tabular}
  \caption{Nonalternating $12$-crossing mutant cliques 4/6}
  \end{centering}
\end{figure}

\begin{figure}[htbp]
  \begin{centering}
  \begin{tabular}{cc@{\hspace{10pt}}|@{\hspace{10pt}}cc}
    \includegraphics[width=75pt]{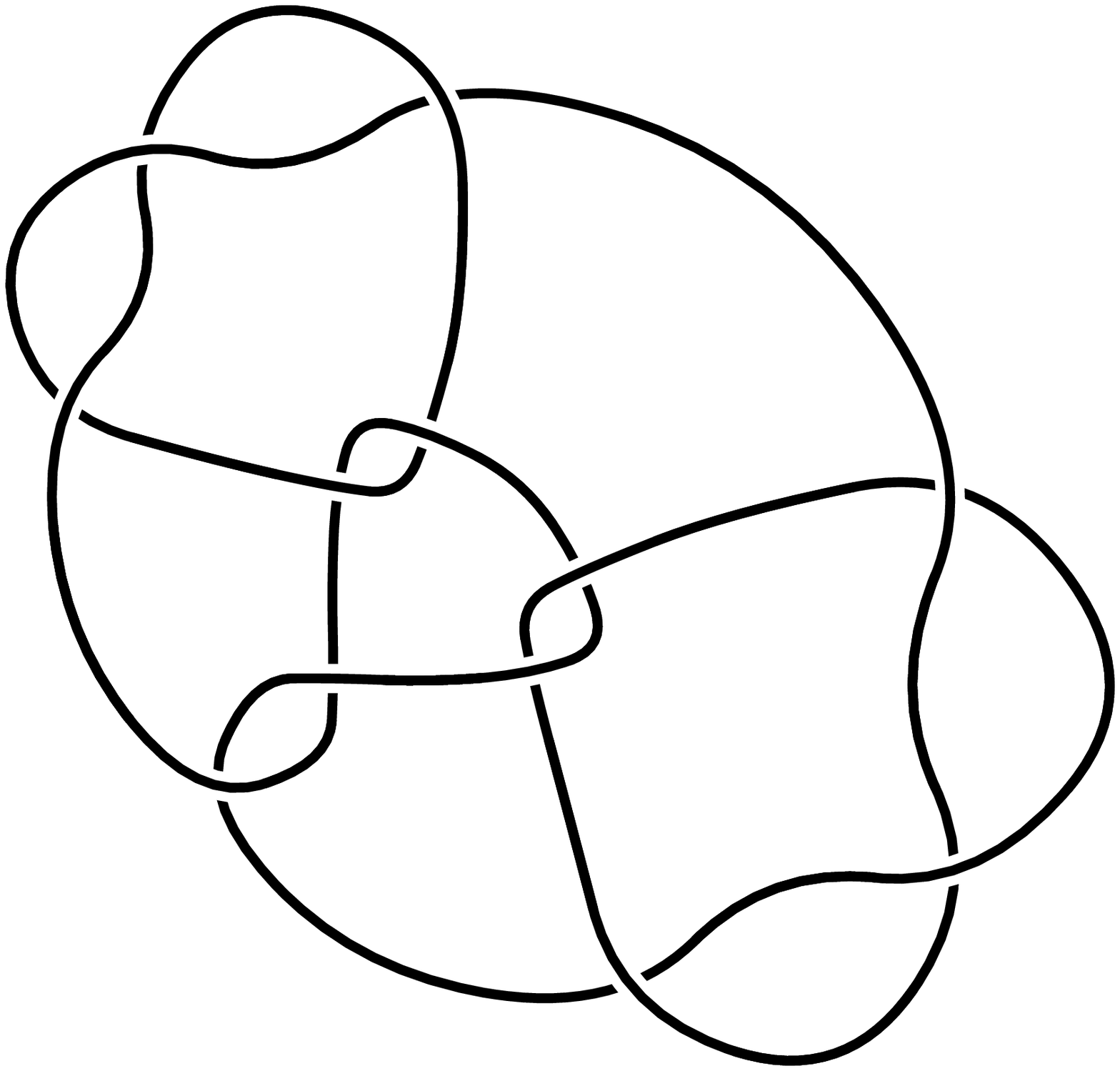}
    &
    \includegraphics[width=75pt]{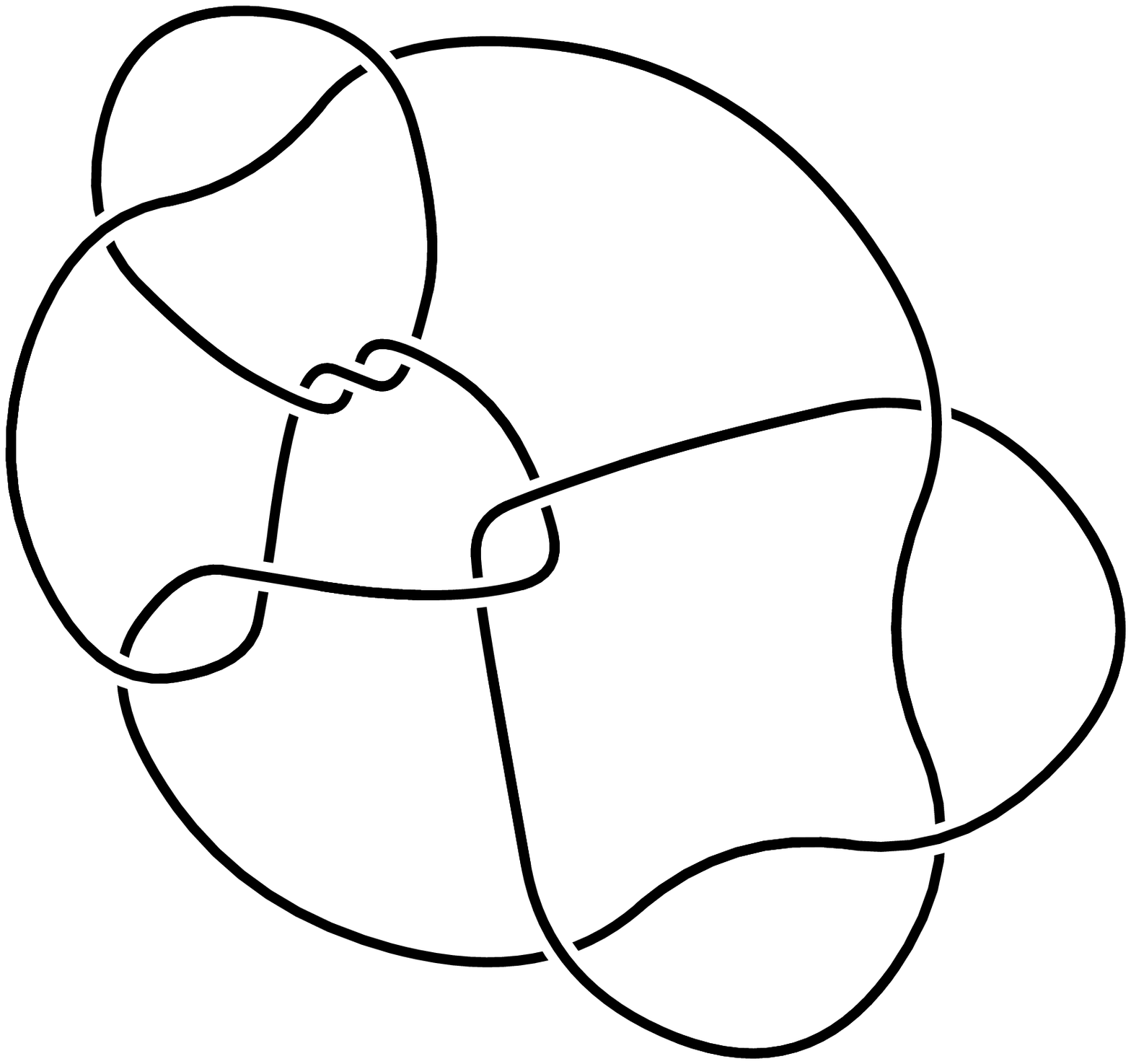}
    &
    \includegraphics[width=75pt]{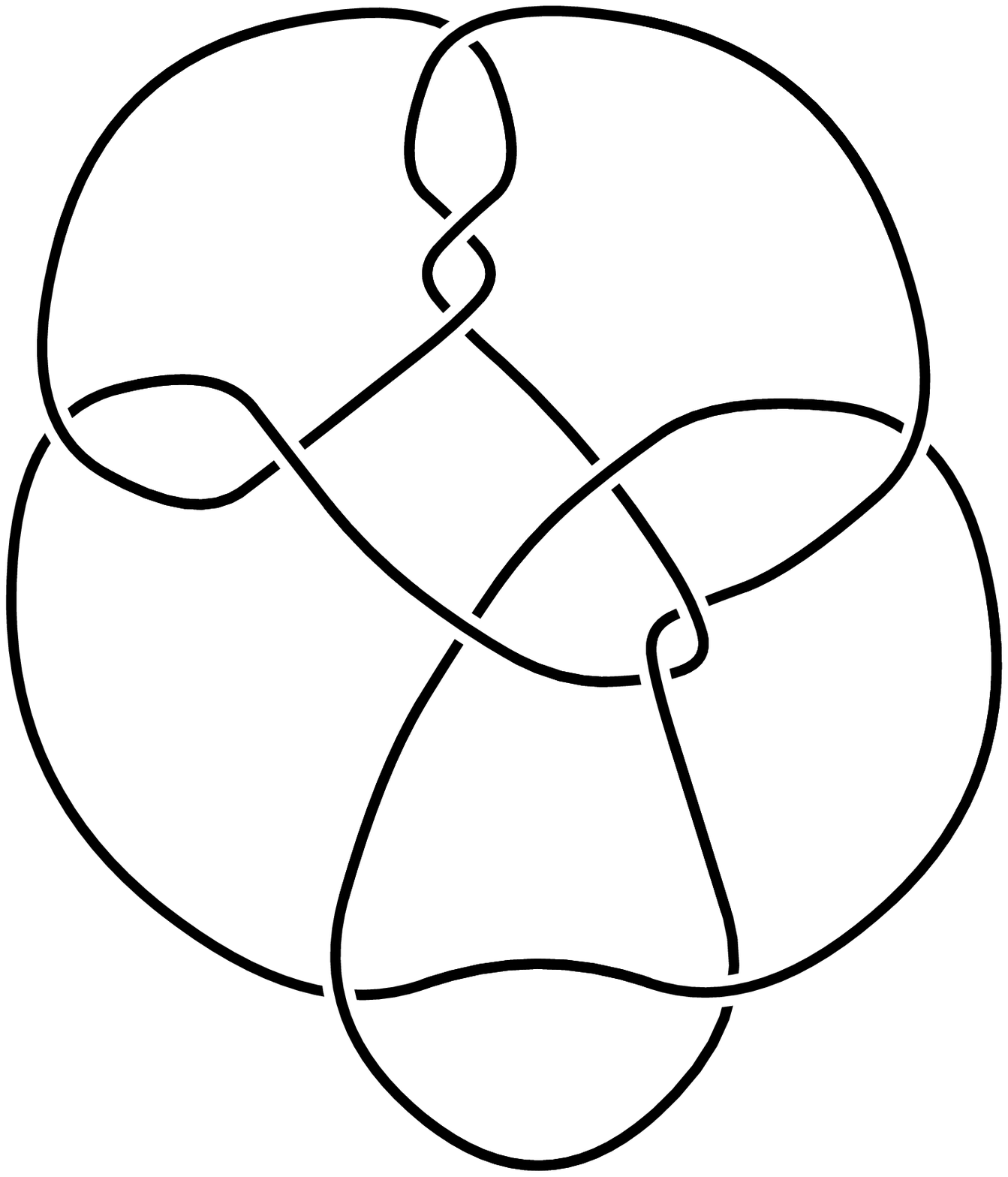}
    &
    \includegraphics[width=75pt]{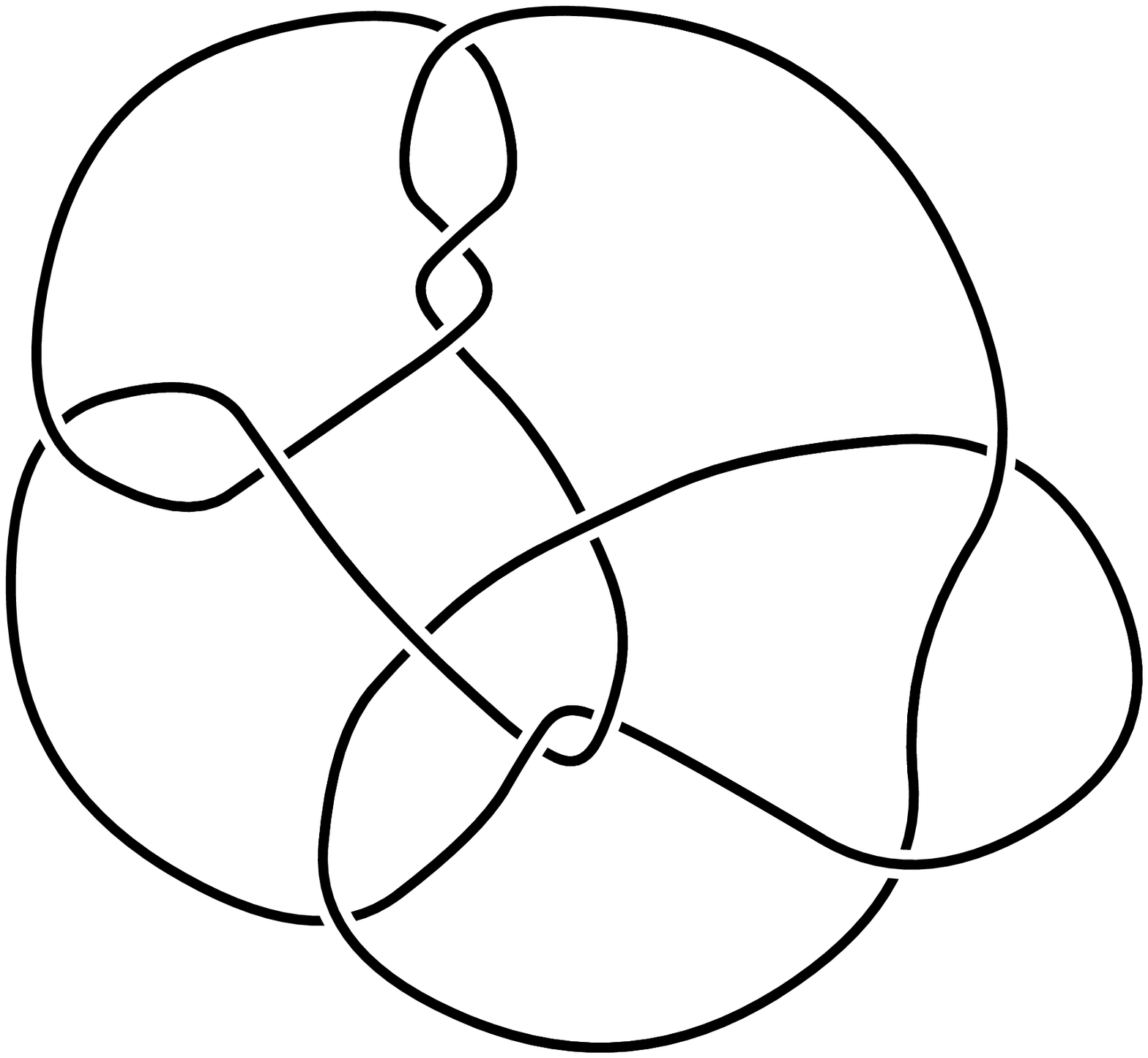}
    \\[-10pt]
    $12^N_{671}$ & $12^N_{682}$ & $12^N_{691}$ & $12^N_{692}$
    \\[10pt]
    \hline
    &&&\\[-10pt]
    \includegraphics[width=75pt]{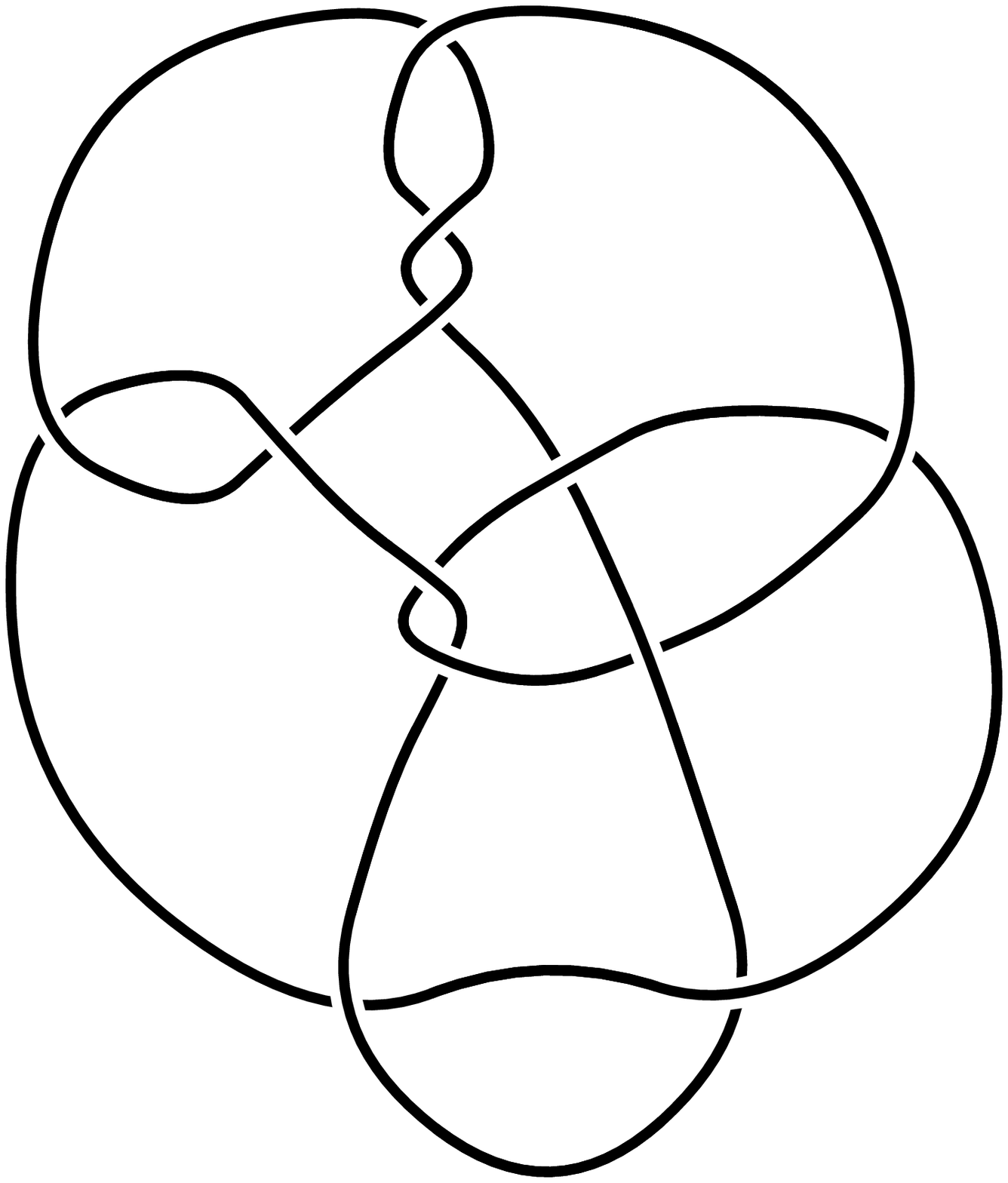}
    &
    \includegraphics[width=75pt]{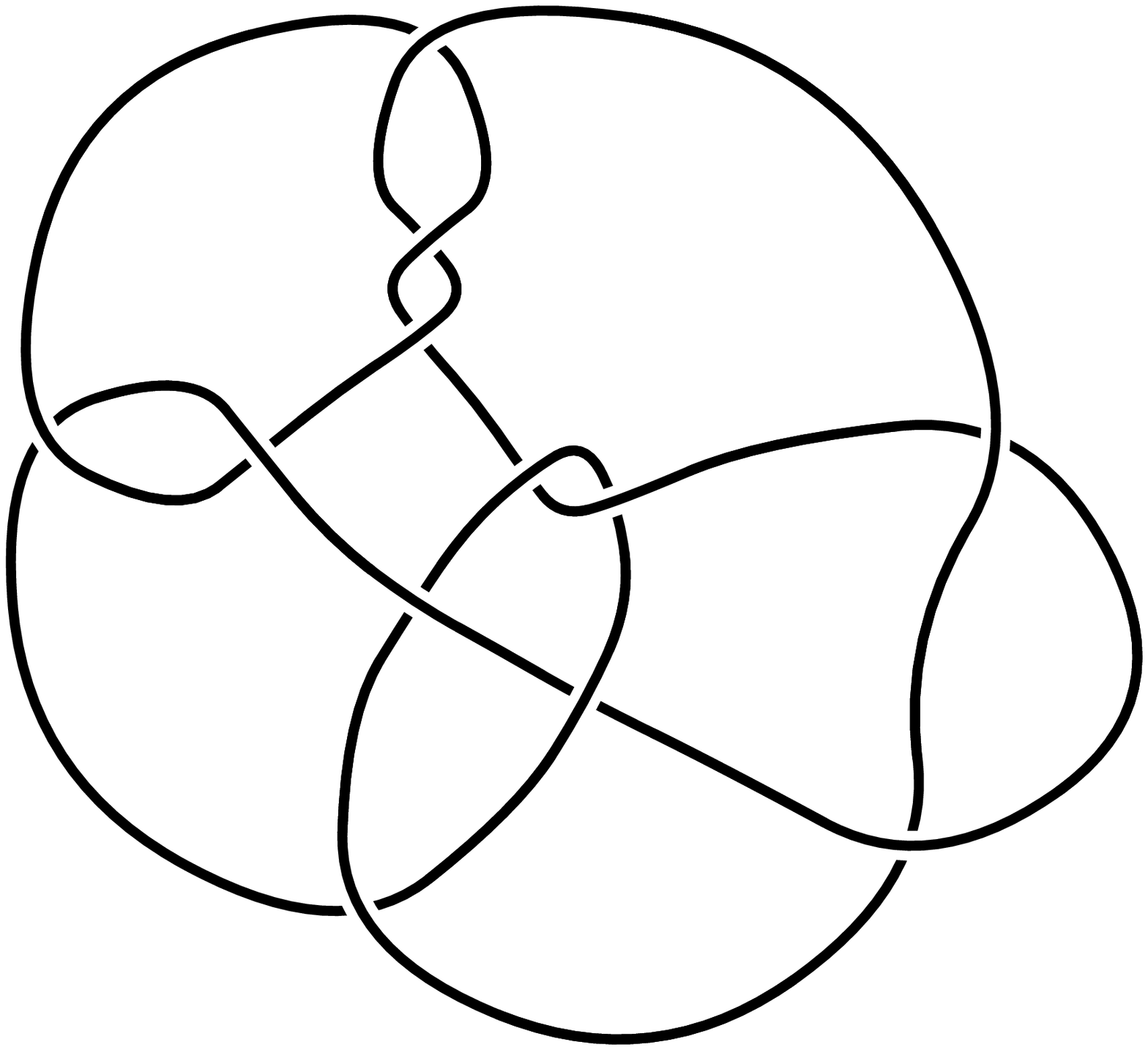}
    &
    &
    \\[-10pt]
    $12^N_{693}$ & $12^N_{696}$ &  &
  \end{tabular}
  \caption{Nonalternating $12$-crossing mutant cliques 5/6}
  \end{centering}
\end{figure}

\begin{figure}[htbp]
  \begin{centering}
  \begin{tabular}{ccc}
    \includegraphics[width=75pt]{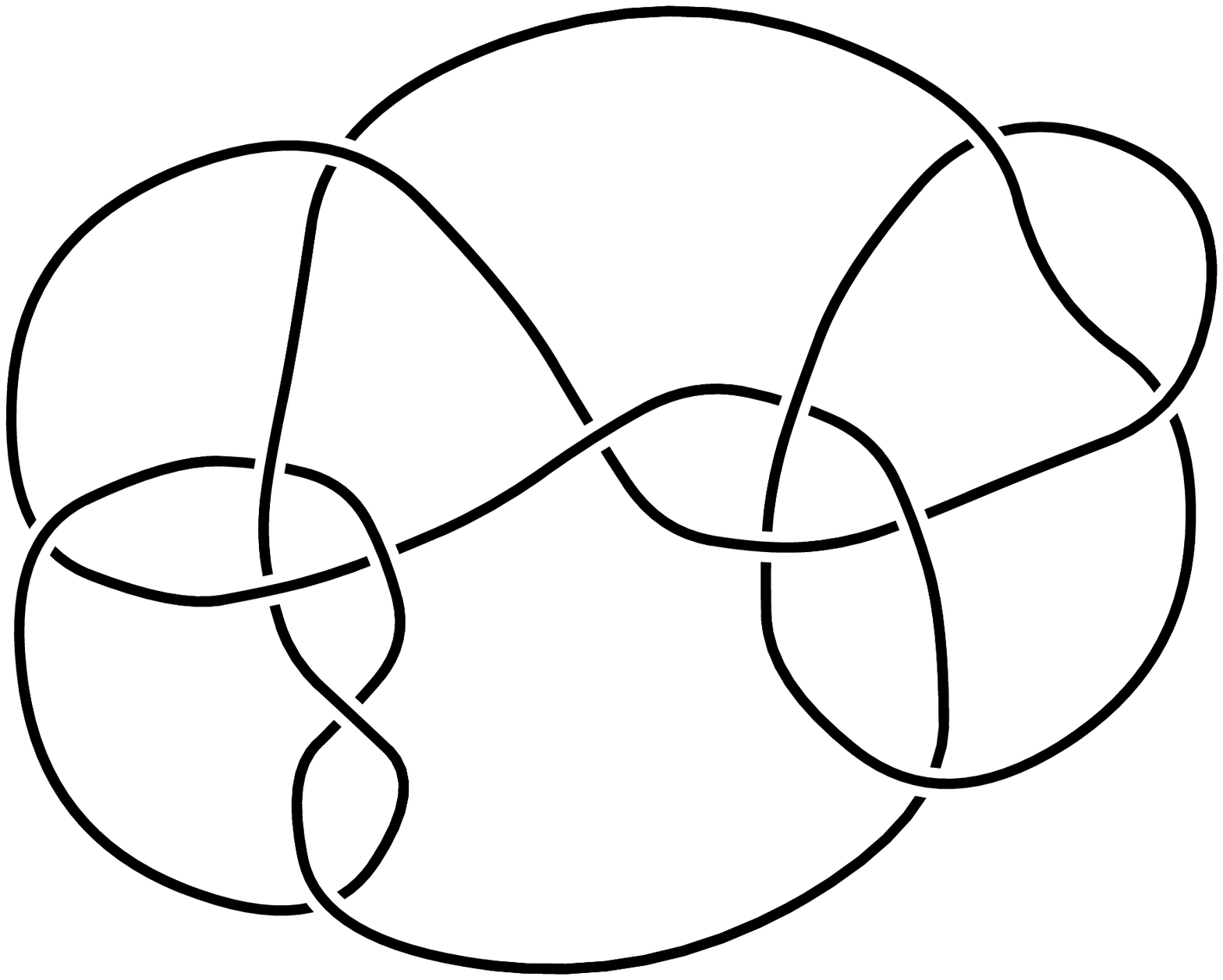}
    &
    \includegraphics[width=75pt]{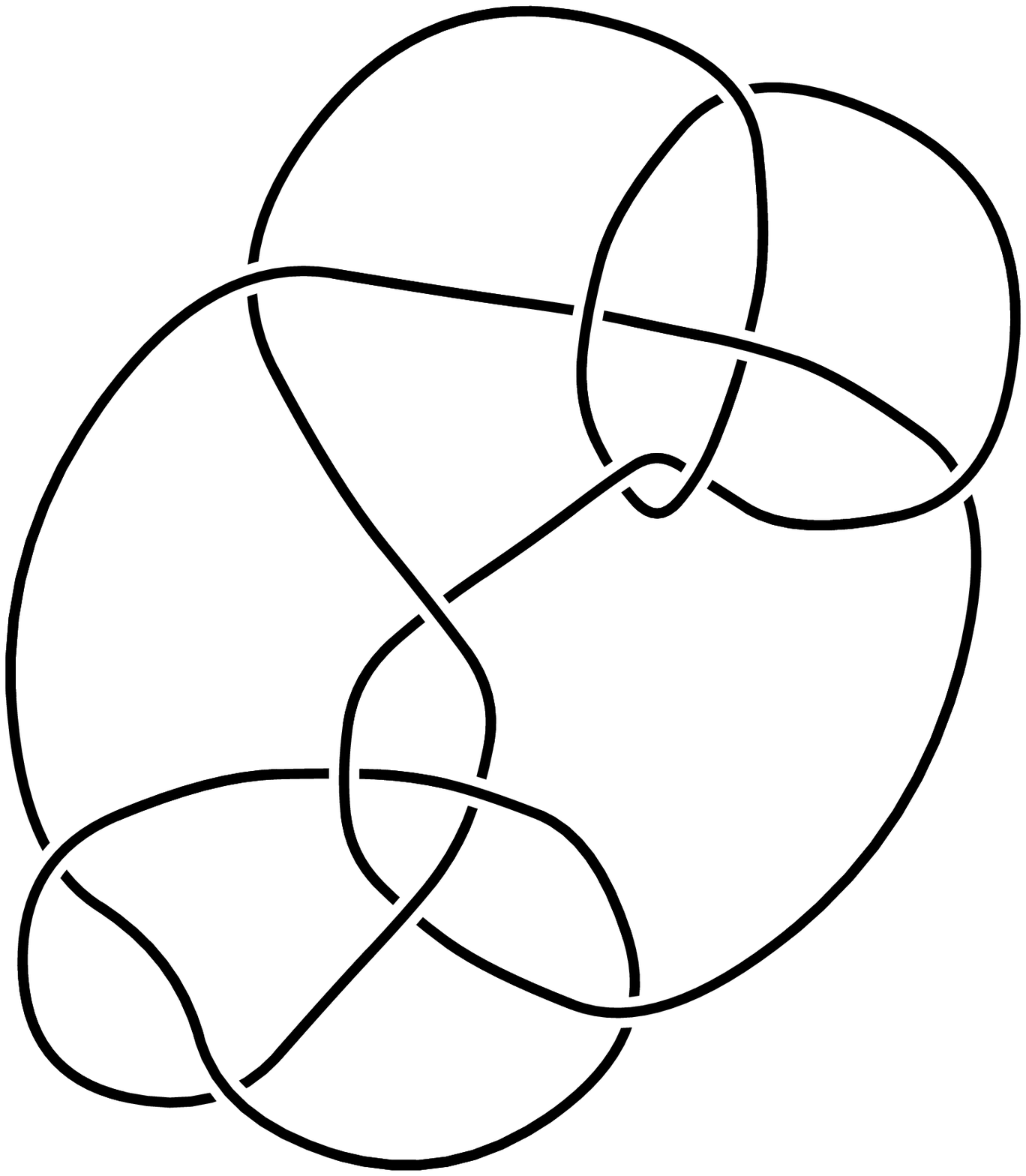}
    &
    \includegraphics[width=75pt]{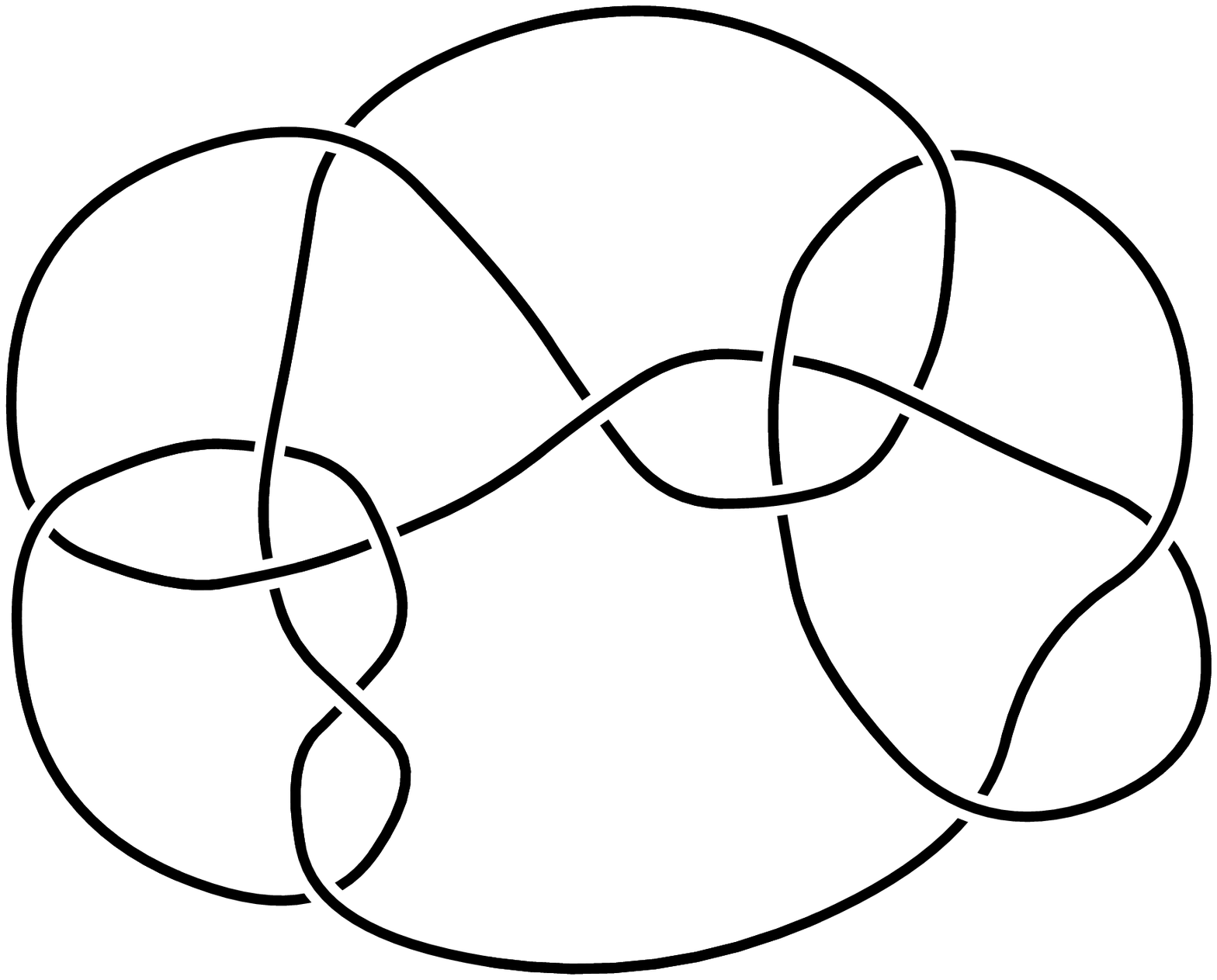}
    \\[-10pt]
    $14^A_{13262}$ & $14^A_{13269}$ & $14^A_{13506}$
    \\[10pt]
    \hline
    &&\\[-10pt]
    \includegraphics[width=75pt]{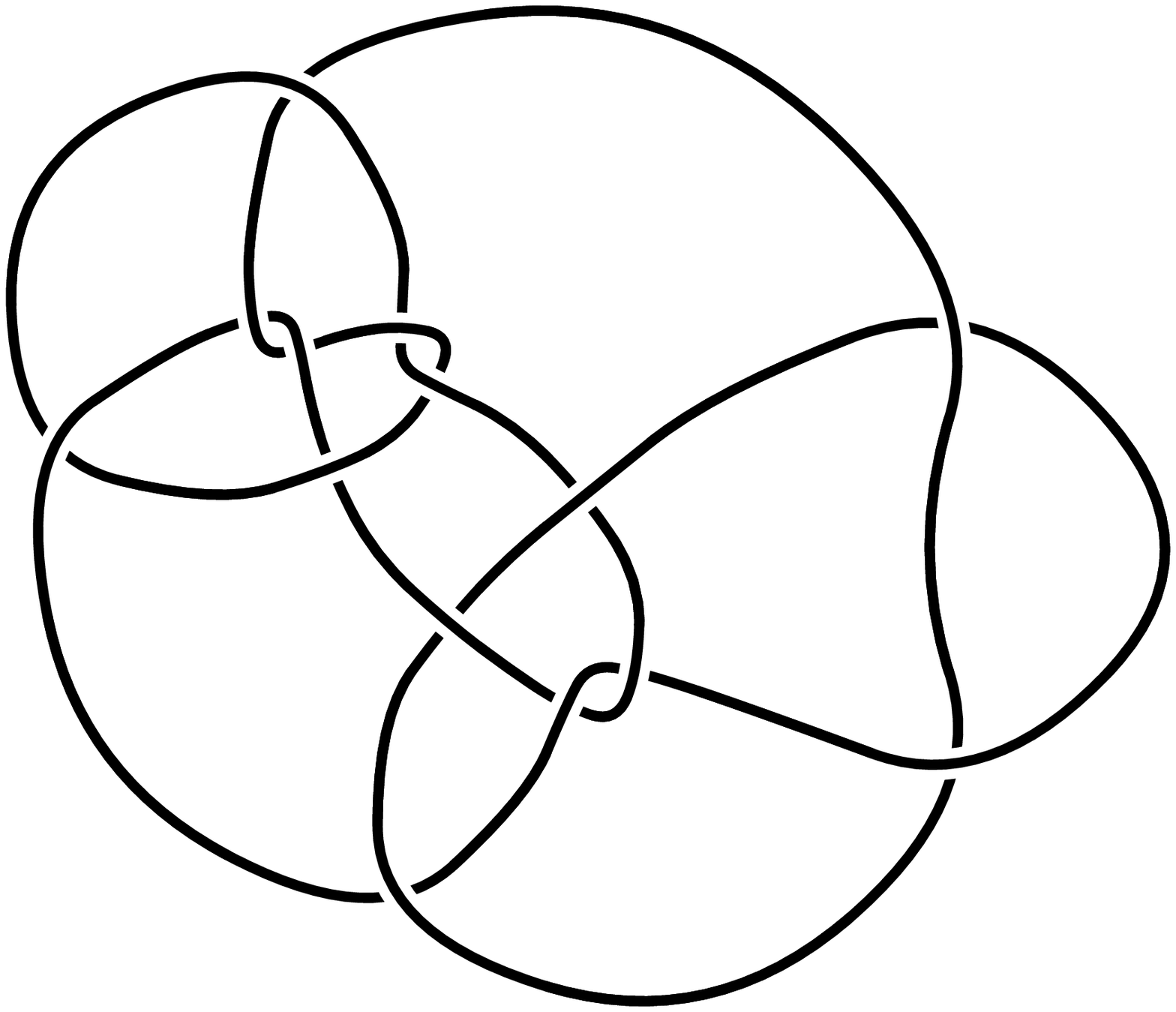}
    &
    \includegraphics[width=75pt]{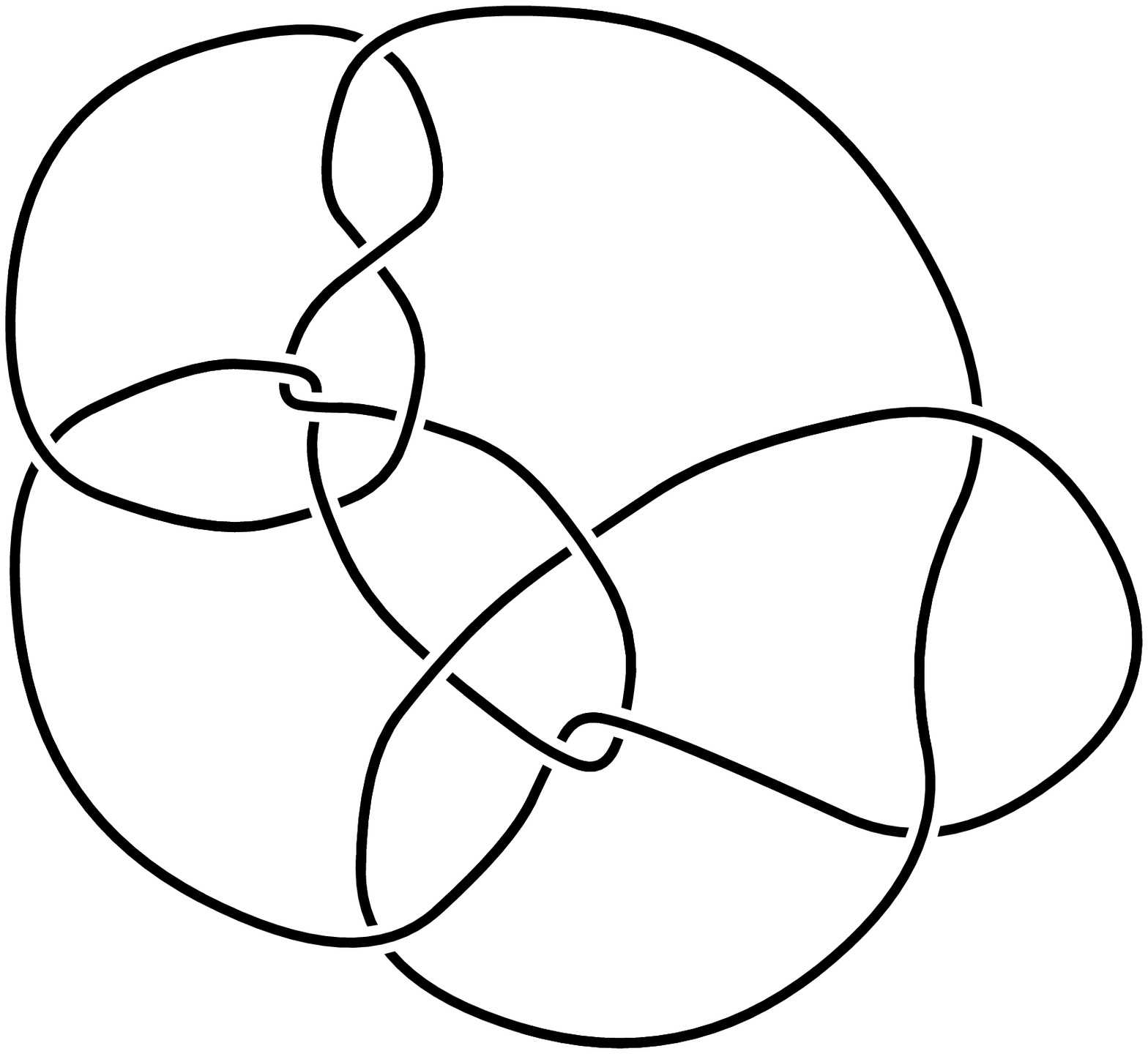}
    &
    \includegraphics[width=75pt]{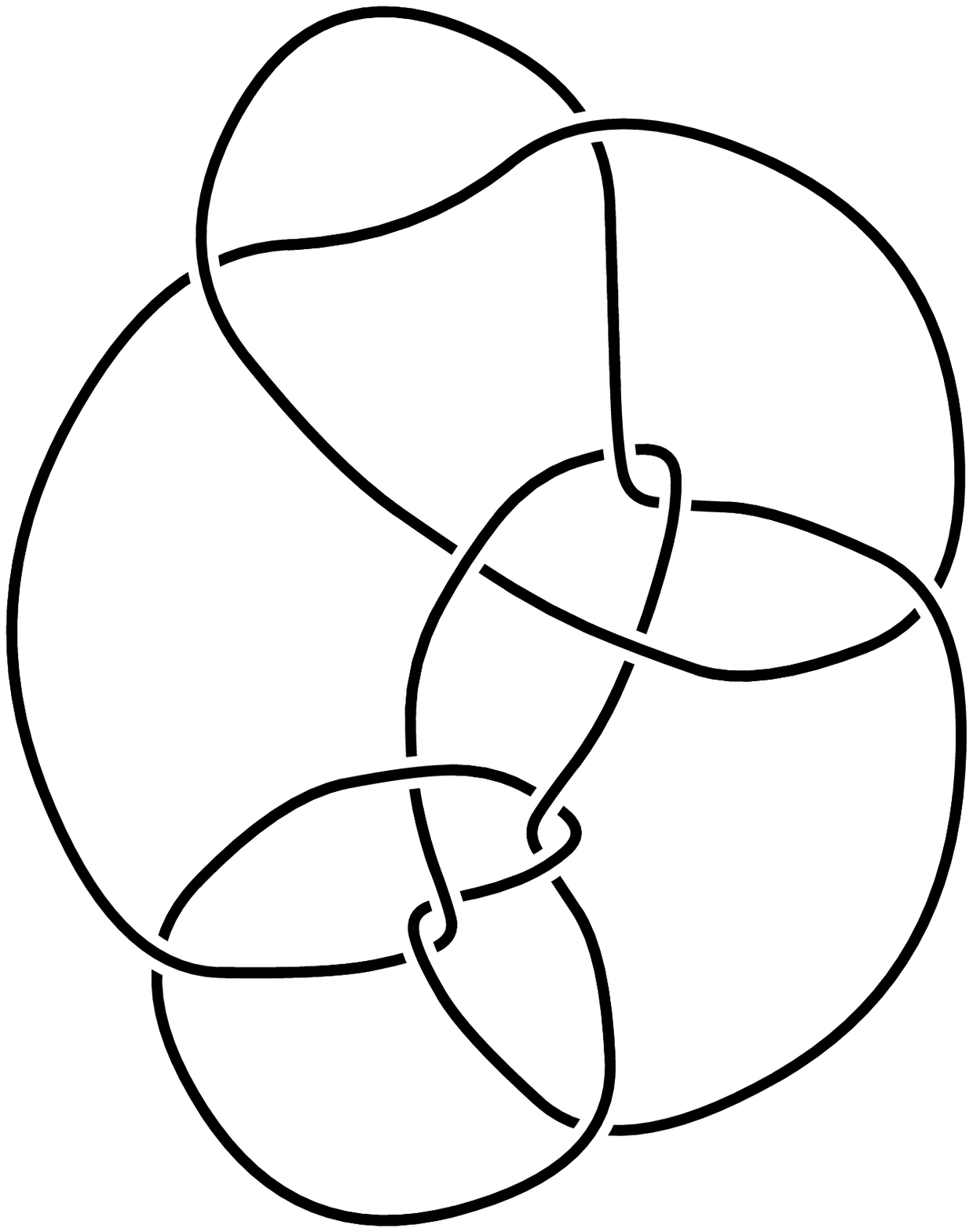}
    \\[-10pt]
    $14^A_{14042}$ & $14^A_{14043}$ & $14^A_{14671}$
    \\[10pt]
    \hline
    &&\\[-10pt]
    \includegraphics[width=75pt]{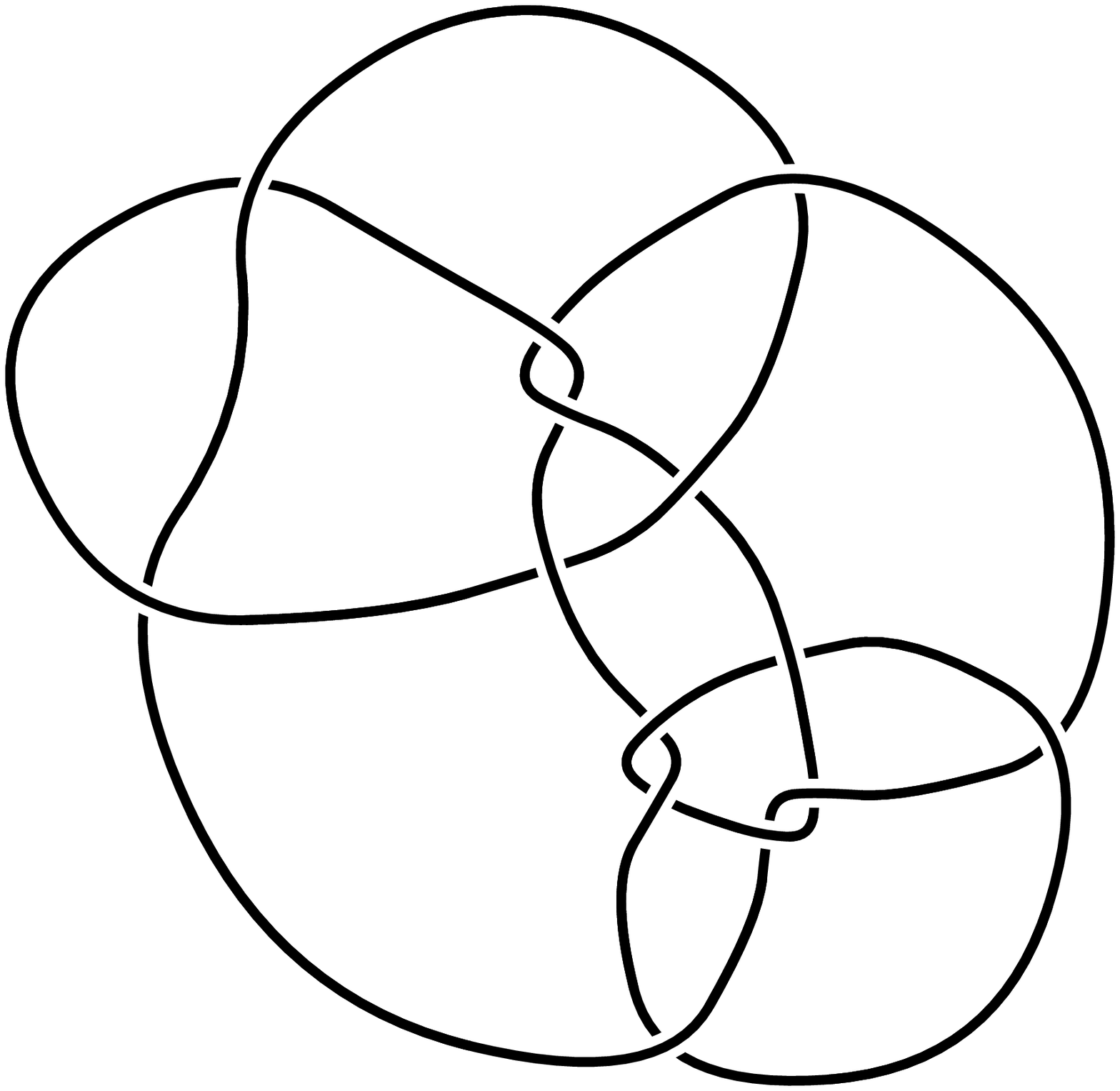}
    &
    \includegraphics[width=75pt]{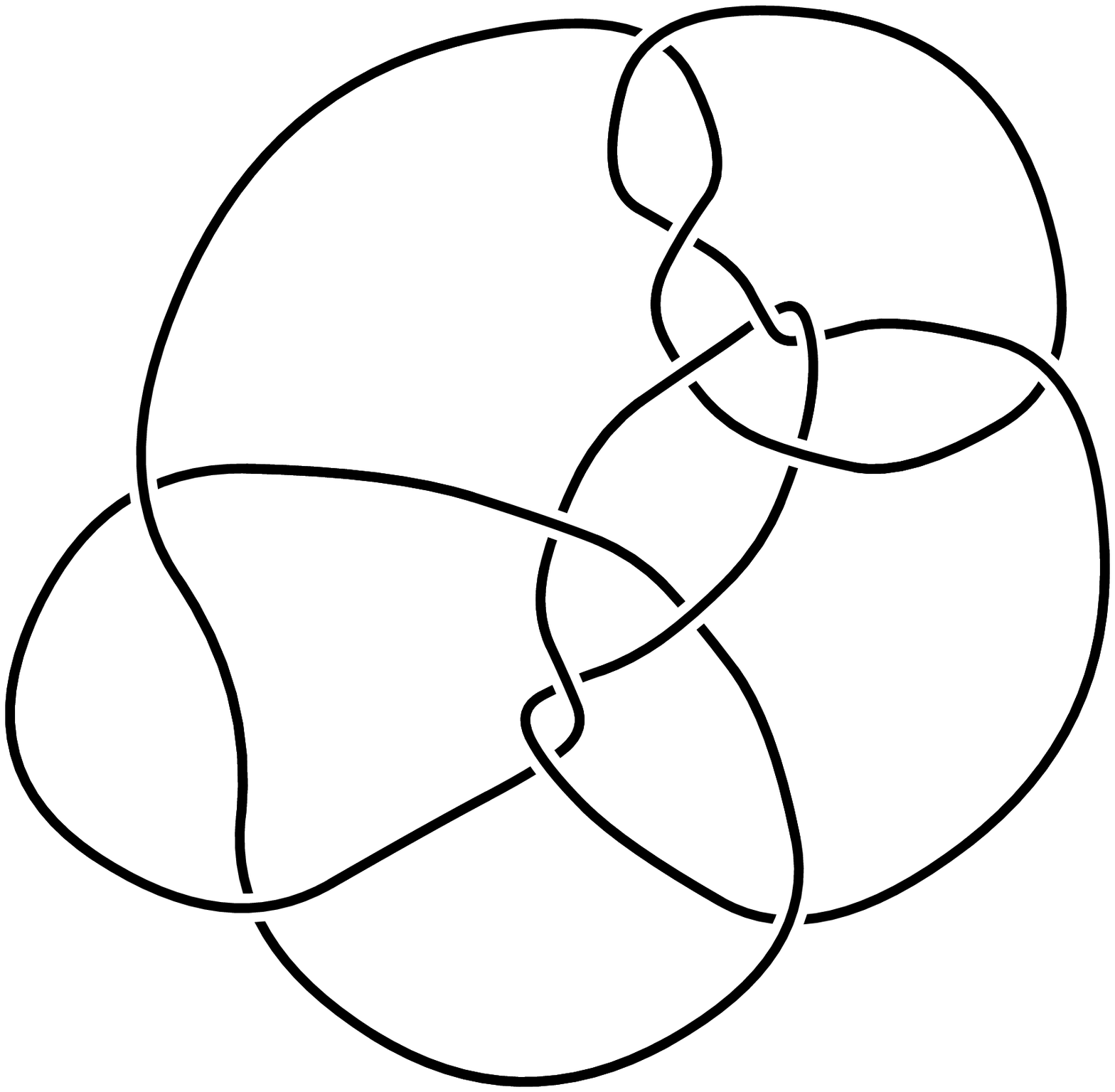}
    &
    \includegraphics[width=75pt]{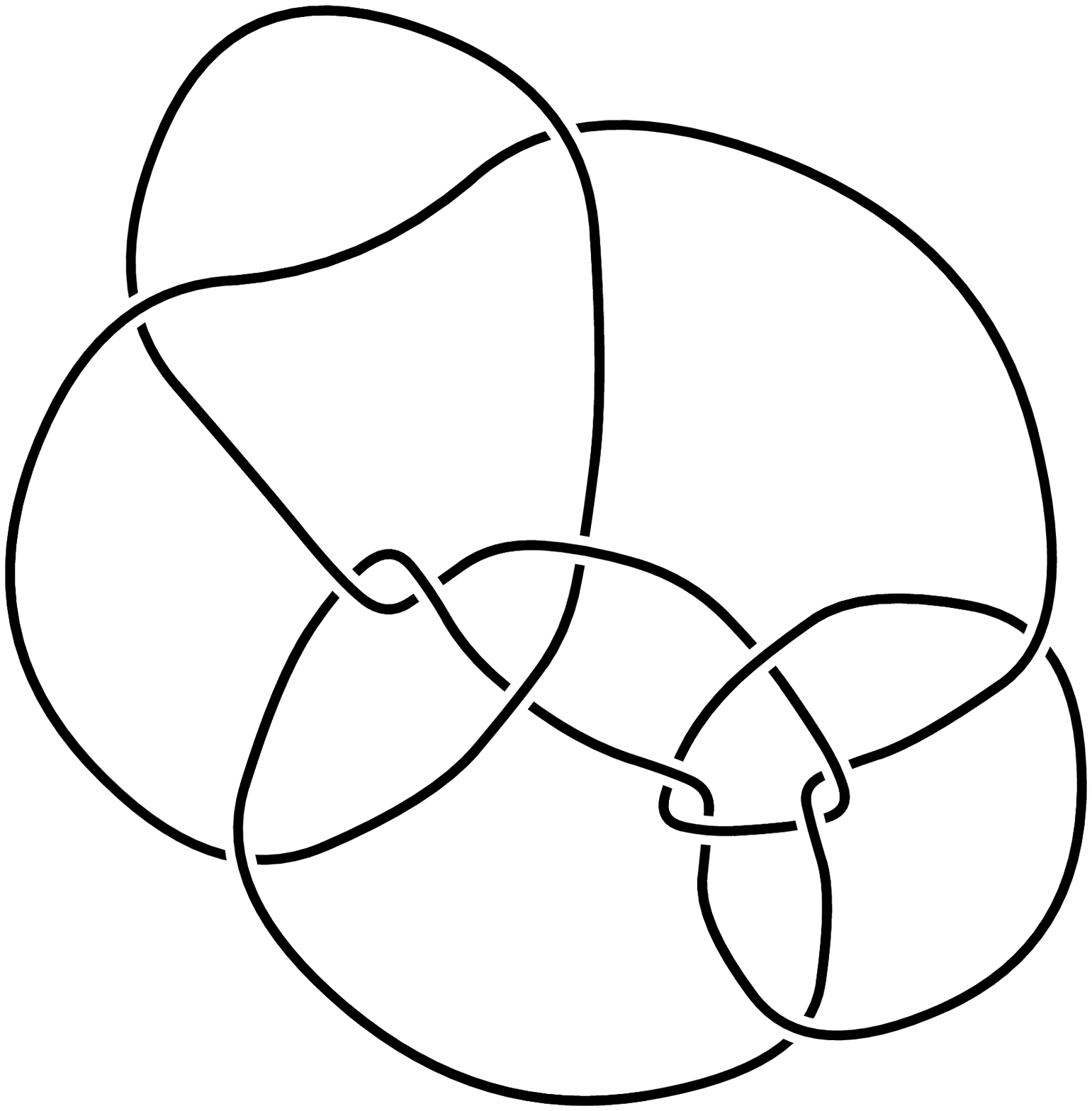}
    \\[-10pt]
    $14^A_{17268}$ & $14^A_{17265}$ & $14^A_{17275}$
    \\[10pt]
    \hline
    &&\\[-10pt]
    \includegraphics[width=75pt]{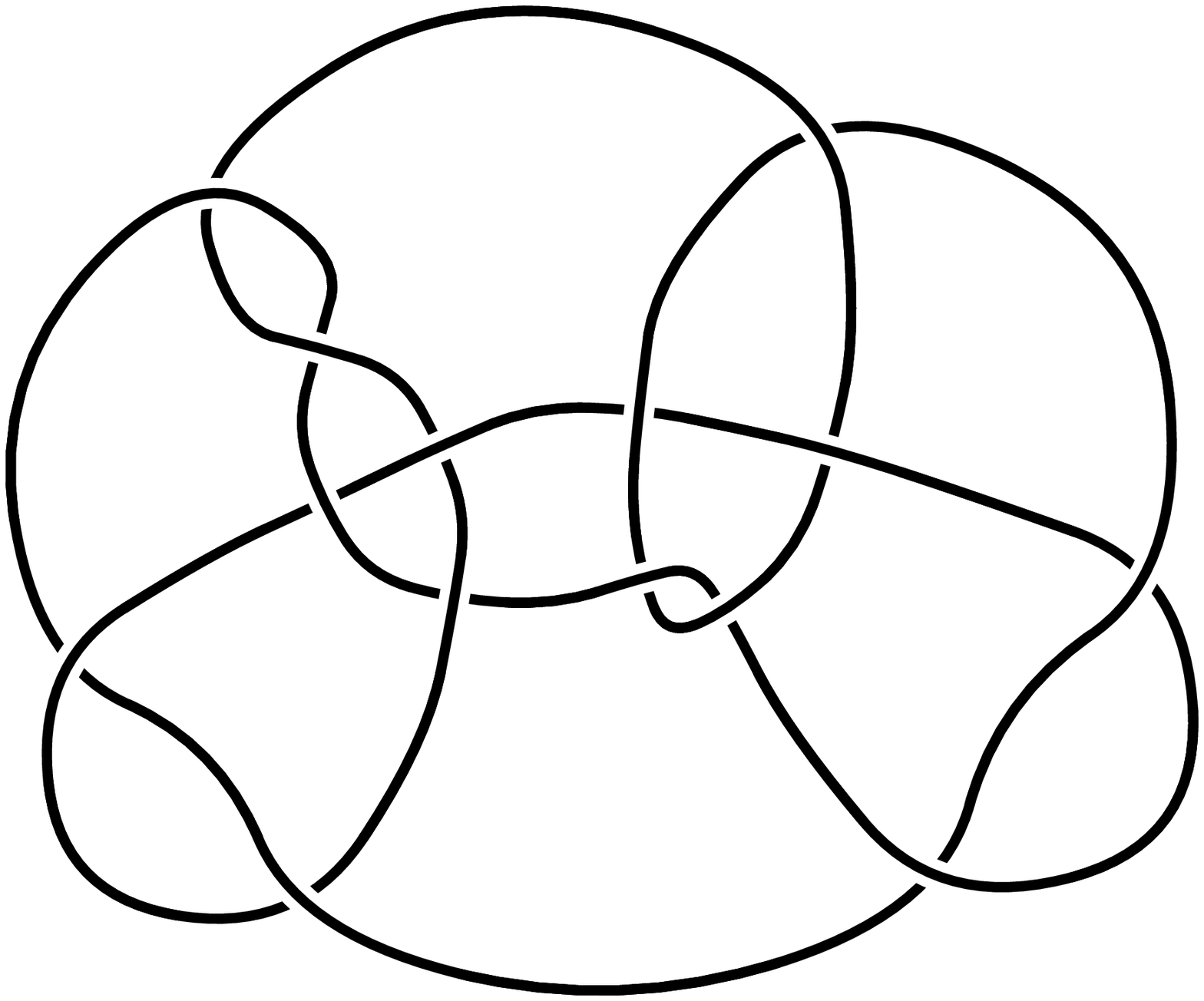}
    &
    \includegraphics[width=75pt,angle=90]{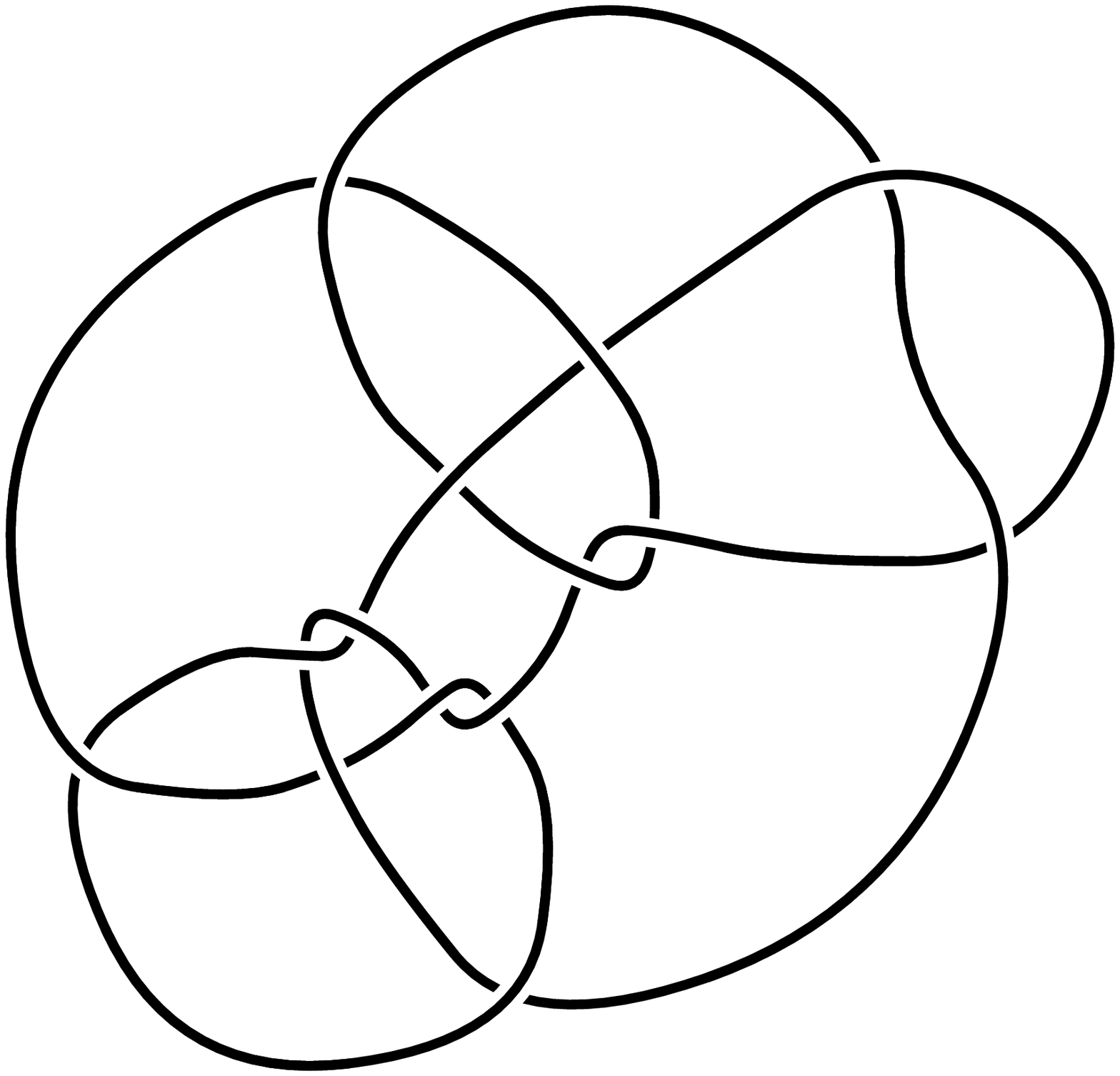}
    &
    \includegraphics[width=75pt]{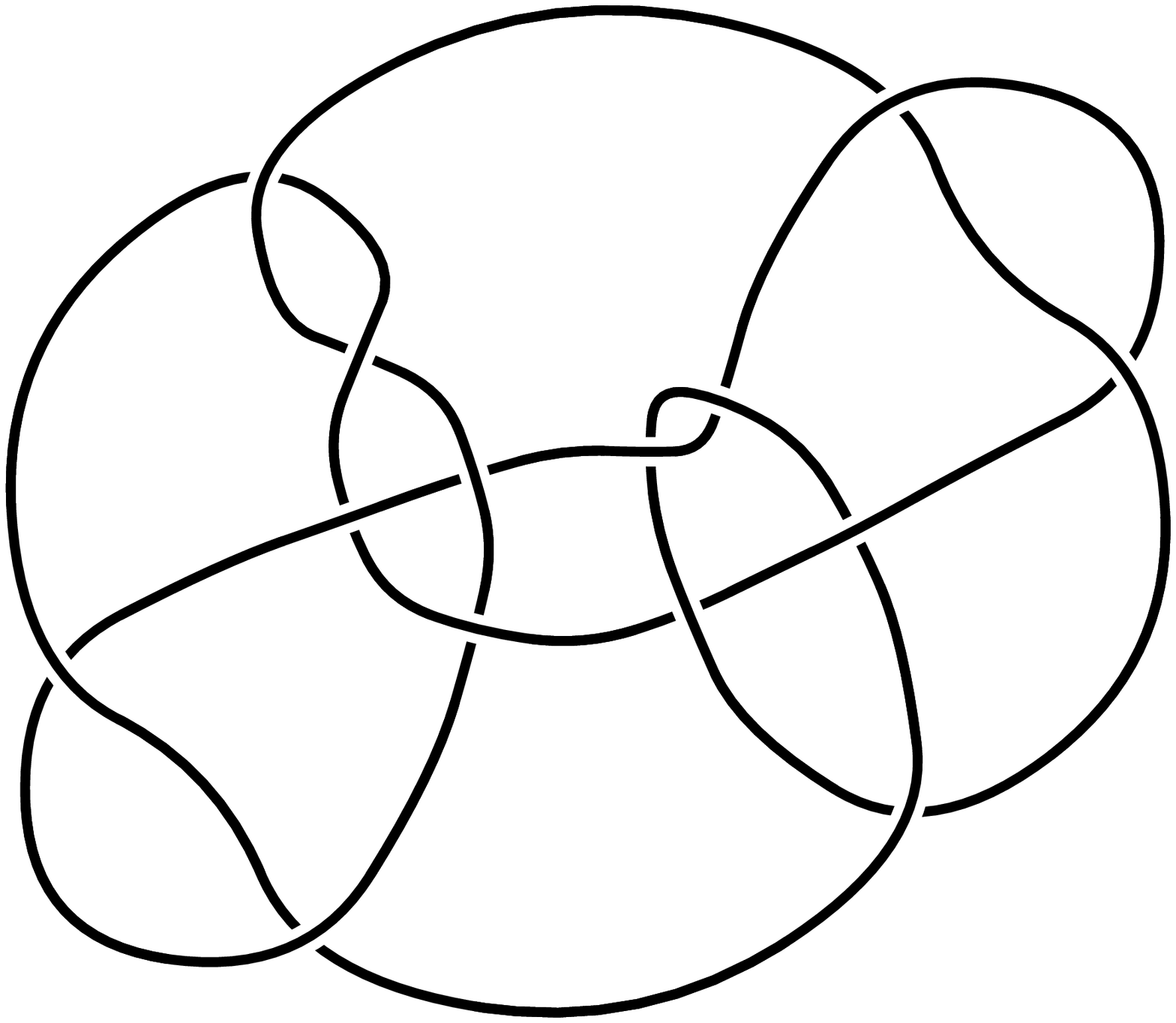}
    \\[-10pt]
    $14^A_{17533}$ & $14^A_{17531}$ & $14^A_{17680}$
  \end{tabular}
  \caption{Alternating $14$-crossing mutant cliques containing both chiral and achiral elements 2/2.}
  \end{centering}
\end{figure}

\end{document}